\documentclass[reqno]{amsart}
\usepackage{bm,amsmath,amsthm,amssymb,mathtools,verbatim,amsfonts,tikz-cd,mathrsfs}
\usepackage{enumitem}
\usepackage{float}

\usepackage[hidelinks,colorlinks=true,linkcolor=blue, citecolor=black,linktocpage=true]{hyperref}
\usepackage{graphicx}
\usepackage[alphabetic]{amsrefs}
\usepackage{caption}
\usepackage{morefloats}
\usepackage[top=1.3in, bottom=1.1in, left=1.3in, right=1.3in,marginpar=1in]{geometry}
\usepackage{subfigure}

\usetikzlibrary{matrix,arrows,decorations.pathmorphing}

\usepackage{chngcntr}
\counterwithin{figure}{section}

\newtheorem{thm}{Theorem}[section]
\newtheorem{prop}[thm]{Proposition}

\newtheorem{cor}[thm]{Corollary}
\newtheorem{lem}[thm]{Lemma}

\theoremstyle{definition}
\newtheorem{define}[thm]{Definition}

\theoremstyle{remark}
\newtheorem{rem}[thm]{Remark}
\newtheorem{example}[thm]{Example}

\newcommand{\ve}[1]{\boldsymbol{\mathbf{#1}}}
\newcommand{\R}{\mathbb{R}}

\newcommand{\ul}[1]{\underline{#1}}

\newcommand{\Z}{\mathbb{Z}}
\newcommand{\N}{\mathbb{N}}
\newcommand{\Q}{\mathbb{Q}}
\newcommand{\C}{\mathbb{C}}

\renewcommand{\d}{\partial}
\renewcommand{\subset}{\subseteq}

\renewcommand{\tilde}{\widetilde}
\renewcommand{\bar}{\overline}

\newcommand{\iso}{\cong}

\DeclareMathOperator{\Aut}{{Aut}}

\DeclareMathOperator{\can}{{can}}

\DeclareMathOperator{\codim}{{codim}}

\DeclareMathOperator{\Diff}{{Diff}}

\DeclareMathOperator{\ev}{{ev}}

\DeclareMathOperator{\gr}{{gr}}

\DeclareMathOperator{\Hom}{{Hom}}

\DeclareMathOperator{\id}{{id}}
\DeclareMathOperator{\ind}{{ind}}
\DeclareMathOperator{\Int}{{int}}

\DeclareMathOperator{\im}{{im}}
\DeclareMathOperator{\Jet}{{Jet}}

\DeclareMathOperator{\Map}{{Map}}

\DeclareMathOperator{\Mor}{{Mor}}

\DeclareMathOperator{\Pin}{{Pin}}

\DeclareMathOperator{\Spin}{{Spin}}

\DeclareMathOperator{\re}{{Re}}

\DeclareMathOperator{\Sing}{{Sing}}
\DeclareMathOperator{\Span}{{Span}}

\newcommand{\bA}{\mathbb{A}}
\newcommand{\bB}{\mathbb{B}}

\newcommand{\bD}{\mathbb{D}}
\newcommand{\bE}{\mathbb{E}}
\newcommand{\bF}{\mathbb{F}}

\newcommand{\bH}{\mathbb{H}}
\newcommand{\bI}{\mathbb{I}}

\newcommand{\bT}{\mathbb{T}}

\newcommand{\bX}{\mathbb{X}}

\newcommand{\bXI}{\mathbb{XI}}

\newcommand{\cA}{\mathcal{A}}
\newcommand{\cB}{\mathcal{B}}
\newcommand{\cC}{\mathcal{C}}
\newcommand{\cD}{\mathcal{D}}
\newcommand{\cE}{\mathcal{E}}
\newcommand{\cF}{\mathcal{F}}

\newcommand{\cH}{\mathcal{H}}
\newcommand{\cI}{\mathcal{I}}
\newcommand{\cJ}{\mathcal{J}}
\newcommand{\cK}{\mathcal{K}}
\newcommand{\cL}{\mathcal{L}}
\newcommand{\cM}{\mathcal{M}}
\newcommand{\cN}{\mathcal{N}}

\newcommand{\cP}{\mathcal{P}}
\newcommand{\cQ}{\mathcal{Q}}
\newcommand{\cR}{\mathcal{R}}

\newcommand{\cT}{\mathcal{T}}
\newcommand{\cU}{\mathcal{U}}
\newcommand{\cV}{\mathcal{V}}
\newcommand{\cW}{\mathcal{W}}

\newcommand{\cY}{\mathcal{Y}}

\newcommand{\frA}{\mathfrak{A}}
\newcommand{\frB}{\mathfrak{B}}

\newcommand{\frF}{\mathfrak{F}}

\newcommand{\frH}{\mathfrak{H}}
\newcommand{\frI}{\mathfrak{I}}

\newcommand{\frP}{\mathfrak{P}}

\newcommand{\frX}{\mathfrak{X}}
\newcommand{\frY}{\mathfrak{Y}}

\newcommand{\frf}{\mathfrak{f}}

\newcommand{\frh}{\mathfrak{h}}

\newcommand{\frm}{\mathfrak{m}}

\newcommand{\fro}{\mathfrak{o}}
\newcommand{\frp}{\mathfrak{p}}

\newcommand{\frs}{\mathfrak{s}}
\newcommand{\frt}{\mathfrak{t}}
\newcommand{\fru}{\mathfrak{u}}

\newcommand{\frw}{\mathfrak{w}}
\newcommand{\frx}{\mathfrak{x}}
\newcommand{\fry}{\mathfrak{y}}
\newcommand{\frz}{\mathfrak{z}}

\newcommand{\scC}{\mathscr{C}}
\newcommand{\scE}{\mathscr{E}}

\newcommand{\scS}{\mathscr{S}}

\newcommand{\cCFK}{\mathcal{C\hspace{-.5mm}F\hspace{-.3mm}K}}

\newcommand{\CF}{\mathit{CF}}
\newcommand{\uCF}{\ul{\mathit{CF}}}
\newcommand{\HF}{\mathit{HF}}
\newcommand{\HFI}{\mathit{HFI}}

\newcommand{\CFK}{\mathit{CFK}}
\newcommand{\HFK}{\mathit{HFK}}

\newcommand{\CFI}{\mathit{CFI}}

\newcommand\CFKi{\CFK^\infty}
\newcommand\HFIhat{\widehat{\HFI}}
\newcommand\HFhat{\widehat{\HF}}

\newcommand{\XI}{\mathbb{XI}}

\newcommand{\PD}{\mathit{PD}}

\newcommand{\xs}{\ve{x}}
\newcommand{\ys}{\ve{y}}
\newcommand{\zs}{\ve{z}}
\newcommand{\ws}{\ve{w}}

\newcommand{\ps}{\ve{p}}
\newcommand{\qs}{\ve{q}}
\newcommand{\as}{\ve{\alpha}}
\newcommand{\bs}{\ve{\beta}}
\newcommand{\gs}{\ve{\gamma}}
\newcommand{\ds}{\ve{\delta}}
\newcommand{\taus}{\ve{\tau}}
\newcommand{\sigmas}{\ve{\sigma}}
\newcommand{\Ds}{\ve{\Delta}}
\newcommand{\xis}{\ve{\xi}}
\newcommand{\zetas}{\ve{\zeta}}

\newcommand{\Dt}{\Delta}

\renewcommand{\a}{\alpha}
\renewcommand{\b}{\beta}
\newcommand{\g}{\gamma}
\newcommand{\dt}{\delta}
\newcommand{\veps}{\varepsilon}

\newcommand{\app}{\mathrm{app}}

\usepackage{leftidx}

\DeclareMathOperator{\emb}{{emb}}
\newcommand{\bCF}{\ve{\CF}}
\newcommand{\bHF}{\ve{\HF}}
\newcommand{\bHFI}{\ve{\HFI}}
\newcommand{\bCFI}{\ve{\CFI}}

\newcommand{\buCF}{\underline{\ve{\CF}}}
\newcommand{\buHF}{\underline{\ve{\HF}}}

\DeclareMathOperator{\Cone}{{Cone}}

\DeclareMathOperator{\St}{{St}}

\newcommand{\Ss}[1]{\scriptstyle{#1}}
\newcommand{\Sss}[1]{\scriptscriptstyle{#1}}

\numberwithin{equation}{section}

\newcommand{\exc}{E}
\newcommand{\ar}{\mathrm{a.r.}}

\newcommand{\llsquare}{[\hspace{-.5mm}[}
\newcommand{\rrsquare}{]\hspace{-.5mm}]}
\newcommand{\wind}{\Sss{\mathrm{W}}}
\newcommand{\twind}{\mathrm{W}}

\newcommand\co{\colon}
\newcommand{\lab}[1]{$\scriptstyle #1$}

\newcommand{\MS}[1]{\footnote{\color{cyan} {\bf MS:} #1}}
\newcommand{\IZ}[1]{\footnote{\color{red} {\bf IZ:} #1}}

\newcommand{\aux}{\mathrm{aux}}
\newcommand{\cen}{\mathrm{cen}}
\newcommand{\oneh}{1\text{-}\mathrm{h}}
\newcommand{\twoh}{2\text{-}\mathrm{h}}
\newcommand{\thrh}{3\text{-}\mathrm{h}}

\newcommand{\Bopp}{{\tilde{B}}}
\newcommand{\vopp}{{\tilde{v}}}

\newcommand{\bBopp}{\tilde{\bB}}

\newcommand{\knotU}{\mathscr{U}}
\newcommand{\knotV}{\mathscr{V}}

\newcommand{\Iota}{\mathsf{Iota}}

\newcommand{\SF}{\mathit{SF}}

\newcommand{\alg}{\mathrm{alg}}

\allowdisplaybreaks

\def\dl {\bunderline{d}}
\def\Vu {\overline{V}}
\def\Vl {\sunderline{V}}
\newcommand{\bunderline}[1]{\underline{#1\mkern-2mu}\mkern2mu }
\newcommand{\sunderline}[1]{\underline{#1\mkern-3mu}\mkern3mu }
\def\du {\bar{d}}

\renewcommand{\re}{\mathrm{re}}

\newcommand{\st}{\mathrm{st}}
\newcommand{\un}{\mathrm{un}}

\title{Surgery exact triangles in involutive Heegaard Floer homology}
\author{Kristen Hendricks}
\thanks{KH was partially supported by NSF grant DMS-2019396 and a Sloan Research Fellowship.}
\address{Department of Mathematics, Rutgers University, New Brunswick, NJ, USA}
\email{kristen.hendricks@rutgers.edu}
\author{Jennifer Hom}
\thanks{JH was partially supported by NSF grant DMS-1552285.}
\address{School of Mathematics, Georgia Institute of Technology, Atlanta, GA, USA}
\email{hom@math.gatech.edu}
\author{Matthew Stoffregen}
\address{Department of Mathematics, Michigan State University, East Lansing, MI, USA}
\email{stoffre1@msu.edu}
\thanks{MS was partially supported by NSF grant DMS-1702532.}
\author{Ian Zemke}
\address{Department of Mathematics\\Princeton University\\  Princeton, NJ, USA}
\email{izemke@math.princeton.edu}
\thanks{IZ was partially supported by NSF grant DMS-1703685.}

\begin{document}

\begin{abstract} 
We establish a surgery exact triangle for involutive Heegaard Floer homology by using a doubling model of the involution. We use this exact triangle to give an involutive version of Ozsv\'ath-Szab\'o's mapping cone formula for knot surgery. As an application, we use this surgery formula to give examples of integer homology spheres that are not homology cobordant to any linear combination of Seifert fibered spaces.
\end{abstract}

\maketitle

\tableofcontents

\section{Introduction}
Heegaard Floer homology, defined by Ozs\'vath-Szab\'o \cite{OSDisks}, is a powerful suite of invariants for studying 3-manifolds, knots inside of them, and 4-manifold cobordisms between them. The first author and Manolescu \cite{HMInvolutive} put additional structure on this package of invariants, in the form of a homotopy involution, leading to involutive Heegaard Floer homology. Their refinement was motivated by Manolescu's resolution of the triangulation conjecture. Galewski-Stern \cite{GalewskiStern} and Matumoto \cite{Matumoto} reduced the triangulation conjecture in dimension at least $5$ to the assertion that there exists an integral homology sphere with Rokhlin invariant one that is of order two in $\Theta^3_\Z$, the 3-dimensional integral homology cobordism group. Manolescu \cite{ManolescuPin2Triangulation} used a Pin$(2)$-equivariant version of Seiberg-Witten Floer homology to show that there exists no such homology sphere, proving that there exist manifolds in every dimension $\geq 5$ that cannot be triangulated.

The involutive Heegaard Floer package has been particularly useful in studying $\Theta^3_\Z$; see, for example, \cite{hendricks-hom-lidman}, \cite{DHSThomcob}. To a closed $3$-manifold $Y$, the first author and Manolescu \cite{HMInvolutive} associate a homotopy involution
\[ \iota \co \CF^-(Y) \to \CF^-(Y). \]
The pair $(\CF^-(Y), \iota)$ is called an \emph{$\iota$-complex} and, up to homotopy equivalence, is an invariant of $Y$. Hendricks, Manolescu, and Zemke \cite{HMZConnectedSum} show that up to a weaker notion of equivalence, called \emph{local equivalence}, the pair $(\CF^-(Y), \iota)$ is an invariant of the homology cobordism class of $Y$.

The goal of this paper is to extend various properties of Heegaard Floer homology, such as the surgery exact triangle and mapping cone formula, to the involutive setting. We show that up to local equivalence, the involutive mapping cone formula simplifies dramatically, paving the way for effective computations. As an application, we use this formula to prove:
\begin{thm}\label{thm:SFSintro}
Seifert fibered spaces do not generate the homology cobordism group $\Theta^3_\Z$.
\end{thm}

\subsection{Involutive surgery exact sequences}

One particularly useful feature of Heegaard Floer homology is the surgery exact triangle. Given a framed knot in a closed 3-manifold $Y$, Ozsv\'ath-Szab\'o \cite{OSProperties} prove that there are exact triangles 
\begin{equation}\label{eq:HFexacttriangle1}
\cdots \to \bHF^-(Y)\to \bHF^-(Y_0(K))\to \bHF^-(Y_1(K))\to \cdots
 \end{equation}
where $\bHF^-(Y)$ denotes the minus flavor of Heegaard Floer homology with completed coefficients.

Our paper centers on establishing the surgery exact triangle for involutive Heegaard Floer homology:
\begin{thm} \label{thm:exact-sequence-1}
 Suppose that $Y$ is a closed 3-manifold and $K$ is a framed knot. There is an exact triangle
 \[
\cdots \to \bHFI^-(Y)\to \bHFI^-(Y_0(K))\to \bHFI^-(Y_1(K))\to \cdots.
 \]
 The same holds for the plus version.
\end{thm}
The surgery exact triangle was previously been established by the first author and Lipshitz \cite{HLInvolutiveBordered}*{Theorem~7.1} for the hat version of involutive Heegaard Floer homology, but remained elusive for the minus theory. Our exact sequence mirrors the exact sequence of Ozsv\'{a}th and Szab\'{o} \cite{OSProperties}*{Theorem~9.1} in ordinary Heegaard Floer homology, as well as Lin's exact sequence for $\Pin(2)$-equivariant Monopole Floer homology \cite{LinExact}*{Theorem~1}.

Additionally, we prove the following analog of \cite{OSIntegerSurgeries}*{Theorem~3.1} for involutive Heegaard Floer homology:

\begin{thm} \label{thm:exact-sequence-2}
Suppose that $Y$ is an integer homology 3-sphere, and $K$ is a knot. For any $n\in \Z$ and $m\ge 1$, there is an exact triangle
\[
\cdots \to \bHFI^-(Y_n(K))\to \bHFI^-(Y_{n+m}(K))\to \underline{\bHFI}^-(Y)\to \cdots.
\]
\end{thm}

In Theorem~\ref{thm:exact-sequence-2}, $\underline{\bHFI}^-(Y)$ denotes a twisted version of involutive Heegaard Floer homology. Furthermore,
\begin{equation}
\underline{\bHFI}^-(Y)\iso \bigoplus_{i=1}^N \bHFI^-(Y)\oplus \bigoplus_{i=1}^M (\bHF^-(Y)[-1]\oplus \bHF^-(Y)), \label{eq:decomposition-twisted-HFI-intro}
\end{equation}
where 
\[
(N,M)=\begin{cases} (1,(m-1)/2)&\text{ if $m$ is odd,} \\
(0,m/2) &\text{ if $m$ is even and $n$ is odd,}\\
(2,(m-2)/2) &\text{ if $m$ and $n$ are even}.
\end{cases}
\]

\subsection{The involutive mapping cone formula}
The exact triangle in Theorem \ref{thm:exact-sequence-2} is a key ingredient in establishing Ozsv\'{a}th and Szab\'{o}'s mapping cone formula \cite{OSIntegerSurgeries} for the Heegaard Floer homology of surgery along a knot. 
The mapping cone formula gives a powerful relation between the knot Floer complex and the Heegaard Floer complex of surgeries. We first recall the formula, and then discuss the involutive counterparts we prove in our present work.

Given a knot $K$ in an integer homology 3-sphere $Y$, define $A_s^-\subset \CFK^\infty(Y,K)$ to be the subcomplex generated over $\bF[U]$ by elements $[\xs,i,j]$ with $i\le 0$ and $j\le s$, which satisfy $A(\xs)-j+i=0$, where $A$ denotes the Alexander grading. We write $\bm{A}^-_s$ for its completion over $\bF\llsquare U\rrsquare$. Define $B_s^-$ to be the subcomplex generated by elements $[\xs,i,j]$ satisfying the same Alexander grading restriction, and $i\le 0$, and write $\bm{B}_s^-$ for its completion. Define
\[
\bA^-:=\prod_{s\in \Z} \bm{A}_s^-\quad \text{and} \quad \bB^-:=\prod_{s\in \Z} \bm{B}_s^-.
\]
Ozsv\'{a}th and Szab\'{o} \cite{OSIntegerSurgeries} construct a map $D_n\colon \bA^-\to \bB^-$ so that
\begin{equation}
\bCF^-(Y_n(K))\simeq \bX_n^-:= \Cone(D_n\colon \bA^-\to \bB^-). \label{eq:mapping-cone-OS}
\end{equation}
We prove the following analog of~\eqref{eq:mapping-cone-OS} for involutive Heegaard Floer homology:

\begin{thm}\label{thm:cone-intro} Suppose $n\ge 1$, and $K$ is a knot in $Y$, an integer homology 3-sphere. Then $\bCFI^-(Y_n(K))$ is homotopy equivalent over the ring $\bF\llsquare U\rrsquare [Q]/Q^2$ to a complex of the form:
\[
\bXI^-_n(Y,K):= \begin{tikzcd}[column sep=2cm, row sep=2cm, labels=description]
\bA^-
	\arrow[d, "Q(\id +\iota_{\bA})",swap]
	\arrow[dr, dashed, "Q H"]
	\arrow[r, "D_n"]
& \bB^-
	\arrow[d, "Q(\id+\iota_{\bB})"]\\
Q \bA^-
	\arrow[r, "D_n"]
& Q \bB^-.
\end{tikzcd}
\]
Here $D_n$ is the map in Ozsv\'{a}th and Szab\'{o}'s mapping cone formula, $\iota_{\bA}$ is a map which is determined by the knot involution on $\CFK^\infty(K)$, $\iota_{\bB}$ is a map which is determined by the 3-manifold involution on $\CF^-(Y)$, and $H$ is a homogeneously graded map making the diagram into a chain complex. Equivalently, the iota-complex $(\CF^-(Y_n(K)), \iota)$ is homotopy equivalent (as an iota-complex) to $\bX^-_n$, equipped with the involution $\iota_{\bA}+H+\iota_{\bB}$.
\end{thm}

A more detailed version of Theorem~\ref{thm:cone-intro} is proven in Section~\ref{sec:involutive-cone}. Theorem~\ref{thm:cone-intro} also has a refinement over $\Spin^c$ structures, which we briefly summarize. There is a natural identification of $\Spin^c(Y_n(K))$ with $\Z_n$, where $\Spin^c$ conjugation is identified with the map $i\mapsto -i$. Of most interest to involutive Heegaard Floer homology are the \emph{self-conjugate} $\Spin^c$ structures, which are the ones satisfying $i\equiv -i \mod n$. If $n$ is odd, there is one self-conjugate $\Spin^c$ structure on $Y_n(K)$, and if $n$ is even, there are two. The complex $\bXI^-_n$ naturally splits over conjugacy classes of $\Spin^c$ structures on $Y_n(K)$, giving the involutive Heegaard Floer homology of $Y_n(K)$ in each conjugacy class. See Section~\ref{sec:involutive-cone} for more details.

\subsection{Knots in \texorpdfstring{$S^3$}{S3}}

One of the most useful properties of Ozsv\'{a}th and Szab\'{o}'s mapping cone formula is that for a knot $K$ in $S^3$, the map $D_n$ is computable up to chain homotopy from the chain complex $\CFK^\infty(K)$. For knots in other 3-manifolds, computation of  $D_n$ usually requires knowledge of  additional holomorphic curves, beyond those used to compute the knot Floer complex. We prove a similar result for the involutive mapping cone formula:

\begin{thm}\label{thm:computable}
If $K$ is a knot in $S^3$, and $r\in \Q$, then $\CFI^-(S^3_r(K))$ is computable from the pair $(\CFK^\infty(K),\iota_K)$.
\end{thm}

We briefly describe what goes into Theorem~\ref{thm:computable}. If $K$ is a knot in an integer homology 3-sphere, then the pair $(\CFK^\infty(Y,K),\iota_K)$ may naturally be viewed as an object in a category of algebraic \emph{$\iota_K$-complexes}; see Section~\ref{sec:iota-complexes} for more details. These consist of chain complexes $C_K$ over a two-variable polynomial ring $\bF[\knotU,\knotV]$, equipped with an endomorphism $\iota_K$, which together satisfy a list of algebraic properties satisfied by any knot Floer complex. We say that an $\iota_K$-complex is of \emph{L-space type} if it algebraically looks like the knot Floer complex of a knot in a 3-manifold $Y$ satisfying $\HF^-(Y)\iso \bF[U]$. Concretely, this amounts to requiring that $H_*(B_s)\iso \bF[U]$. In Section~\ref{sec:canonical-mapping-cone}, we construct an algebraically defined map
\[
\XI_n^{\alg}\colon \frac{\{\iota_K\text{-complexes of L-space type}\}}{\iota_K\text{-homotopy equivalence}}\to \frac{\{\iota\text{-complexes}\}}{\iota\text{-homotopy equivalence}},
\]
which gives the maps $D_n$, $\iota_{\bA},$ $\iota_{\bB}$, and $H$, appearing in Theorem~\ref{thm:cone-intro}, via a specified procedure. We call $\XI_n^{\alg}$ the \emph{algebraic} involutive mapping cone. The map $\XI_n^{\alg}$ enjoys several functoriality properties, with respect to morphisms of $\iota_K$-complexes.

  The key property that characterizes $\XI_n^{\alg}$ is that the maps $D_n$ and $H$ admit a factorization though the inclusions of $A_s^-$ into $B_s^-$ and $A_s^-$ into ${\Bopp}^-_s$, where ${\Bopp}^-_s$ denotes the subcomplex of $\CFK^\infty(K)$ generated over $\bF$ by elements $[\xs,i,j]$ with $j\le s$.

In Section~\ref{sec:strong-cone-formula} we prove that
\[
(\bCF^-(S^3_n(K),\iota))\simeq  \bXI^{\alg}_n(\CFK^\infty(K),\iota_K).
\]

Although Theorem~\ref{thm:computable} allows for algorithmic computation of $\CFI^-(S^3_{n}(K))$ whenever the $\iota_K$-complex $(\CFK^\infty(K),\iota_K)$ is known, it often produces very large chain complexes. For applications to homology cobordism, it is sufficient to work with the local equivalence class of $\CFI^-(S^3_n(K))$. We prove the following:
\begin{thm}\label{thm:local-class-intro}
 Suppose $K$ is a knot in $S^3$ and $n>0$ is an integer.
 \begin{enumerate}
 \item The $\iota$-complex $\CFI^-(S^3_n(K),[0])$ is locally equivalent to $(A_0^-(K),\iota_K)$. Here, $[0]$ is the self-conjugate $\Spin^c$ structure of $S^3_n(K)$ corresponding to $[0]\in \Z/n$ under the standard correspondence $\Spin^c(S^3_n(K))\iso \Z/n$.
 \item  The $\iota$-complex $(\CF^-(S^3_{2n}(K),[n]), \iota)$ is locally equivalent to the complex
  \[
 \begin{tikzcd}[column sep={1cm,between origins},labels=description] 
A_{n}^-
	\ar[dr, "v"]
& & A_n^- \ar[dl,"v"]\\
& B_n^-
\end{tikzcd}
\]
 with the involution which swaps the two copies of $A_n^-$, and fixes $B_n^-$.
 \end{enumerate}
\end{thm}

The fact that for $n>0$ the local equivalence class of $\CFI^-(S^3_n(K),[0])$ is determined by $(A_0^-(K),\iota_K)$ is perhaps unsurprising in light of \cite[Proposition 1.6]{NiWu}, where Ni--Wu show that $d(S^3_n(K), [0])$ (which can be thought of as the local equivalence class of $\CF^-(S^3_n(K), [0])$ in the non-involutive setting) is determined by $A_0^-(K)$. Theorem~\ref{thm:local-class-intro} also has an analog for rational surgeries. See Proposition~\ref{prop:local-class-surgeries-rational} for a precise statement. 

In \cite{HMInvolutive}, the first author and Manolescu use involutive Heegaard Floer homology to associate two new correction terms $\dl(Y, \frs)$ and $\du(Y,\frs)$ to a rational homology sphere $Y$ with a conjugation invariant $\Spin^c$ structure $\frs$, and two related knot concordance invariants $\Vl_0(K)$ and $\Vu_0(K)$. (The definitions of these invariants are reviewed in Sections \ref{sec:3-manifolds} and \ref{sec:involutive-knots}.) As a consequence of Theorem~\ref{thm:local-class-intro}, and its analog for rational surgeries, we compute the involutive correction terms of $p/q$ surgery in terms of $p/q$ and the involutive concordance invariants.

\begin{prop}\label{prop:correction-terms-intro}
 Suppose that $p,q>0$ are relatively prime integers and $K$ is a knot in $S^3$.
\begin{enumerate}
\item If $p$ and $q$ are odd, then there is one self-conjugate $\Spin^c$ structure  $[0]$ on $S^3_{p/q}(K)$, and furthermore
\[
\dl(S^3_{p/q}(K),[0])=d(L(p,q),[0])-2 \Vl_0(K), \qquad \du(S^3_{p/q}(K),[0])=d(L(p,q),[0])-2\Vu_0(K).
\]
\item If $p$ is even and $q$ is odd, then there are two self-conjugate $\Spin^c$ structures, which we denote $[0]$ and $[p/2q]$. The correction terms for $[0]$ are the same as the previous case, while
\[
\dl(S^3_{p/q}(K), [p/2q])=d(S_{p/q}^3(K), [p/2q]), \qquad \du(S_{p/q}^3(K), [p/2q])=d(L(p,q), [p/2q]).
\]
\item If $p$ is odd and $q$ is even, then there is just one self-conjugate $\Spin^c$ structure $[p/2q]$, and
\[
\dl(S^3_{p/q}(K), [p/2q])=d(S_{p/q}^3(K), [p/2q]), \qquad \du(S_{p/q}^3(K), [p/2q])=d(L(p,q), [p/2q]).
\]
\end{enumerate} 
\end{prop}

Note that Proposition \ref{prop:correction-terms-intro} is a generalization of \cite[Theorem 1.6]{HMInvolutive}, which is the case $q=1$ and $p \geq g(K)$. Proposition~\ref{prop:correction-terms-intro} has a particularly simple form for $+1$ surgeries:

\begin{cor}If $K$ is a knot in $S^3$, then
\[
\dl(S^3_{+1}(K))=-2\Vl_0(K)\quad \text{and} \quad \du(S^3_{+1}(K))=-2\Vu_0(K).
\]
\end{cor}

\subsection{An application to the homology cobordism group}\label{sec:applicationtohomcob}
In the 1950s, the Rokhlin homomorphism
\[ \mu \co \Theta^3_\Z \to \Z/2\Z, \]
showed that $\Theta^3_\Z$ is non-trivial. Since then significant progress has been made in understanding the structure of the 3-dimensional homology cobordism group, starting with Fintushel-Stern's proof \cite{FintushelStern85}, using techniques from Donaldson's proof of the diagonalization theorem, that the Poincar\'e homology sphere has infinite order in $\Theta^3_\Z$. Subsequently, Fintushel-Stern \cite{FintushelStern90} and Furuta \cite{furuta} established that $\Theta^3_\Z$ contains a $\Z^\infty$-subgroup, using Yang-Mills gauge theory; the third author \cite{stoffregen-sums} has given an alternate proof of this fact using $\mathrm{Pin}(2)$-equivariant Seiberg-Witten theory.

Dai, Hom, Stoffregen and Truong \cite{DHSThomcob} described a homomorphism
\[
\hat{h}\colon \Theta_{\Z}^3\to \hat{\frI},
\]
where $\hat{\frI}$ denotes the group of \emph{almost $\iota$-complexes}; roughly, $\hat{\frI}$ consists of $\iota$-complexes up to a slightly weaker version of local equivalence. They construct an infinite family of linearly independent surjective homomorophisms
\[ \phi_n \co \hat{\frI} \to \Z \]
to show that $\Theta^3_\Z$ contains a direct summand isomorphic to $\Z^\infty$.

Despite this progress in understanding the structure of $\Theta^3_\Z$, there has been relatively less progress on understanding exactly which kinds of manifolds can represent a given class $\xi \in \Theta^3_\Z$. 

The first results about this question were in the positive direction. Livingston \cite{livingston} showed that any class in $\Theta^3_\Z$ is represented by an irreducible integer homology sphere.  Later, Myers \cite{Myers} proved that in fact any class admits a hyperbolic representative.  
	
	In the negative direction, Fr{\o}yshov (in unpublished work), F. Lin \cite{LinExact}, and the third author \cite{StoffregenSeifertFibered} constructed classes $\xi \in \Theta^3_\Z$ that are not represented by any Seifert fibered spaces.  However, none of these proofs suffices to establish that a homology cobordism class is not represented by a \textit{connected sum} of Seifert spaces.  The classes constructed by Fr{\o}yshov and in \cite{StoffregenSeifertFibered} are both connected sums of Seifert spaces: the proofs proceeded by describing, respectively, the behavior of mod 2 instanton homology and $\mathrm{Pin}(2)$-equivariant Seiberg-Witten Floer homology of connected sums, as well as the Floer homologies of Seifert spaces.  The class $\xi$ constructed by Lin is surgery on a certain alternating knot, but has Floer homology consistent with it being representable as a connected sum of Seifert spaces.
All of these arguments depend on the fact that instanton Floer homology and Seiberg-Witten Floer homology of Seifert fibered spaces are relatively easy to understand, and then using that these theories can be understood under connected sum.

The group $\hat{\frI}$ admits $\hat{h}(\Theta_{\SF})$ as a summand \cite[Figure 14]{DHSThomcob}, where $\Theta_{\SF}$ denotes the span of Seifert fibered homology $3$-spheres in $\Theta^3_{\Z}$. Furthermore, building on \cite{dai-manolescu}, \cite{dai-stoffregen} and \cite{DaiConnected}, Dai-Hom-Stoffregen-Truong prove $\hat{h}(\Theta_{\SF})\iso \Z^{\infty}$, and that $\hat{h}(\Theta_{\SF})$ is a proper subgroup of $\hat{\frI}$. Due to the challenging nature of performing computations, they were unable to produce an integer homology 3-sphere $Y$ such that $\hat{h}(Y)\not \in \hat{h}(\Theta_{\SF}).$ Using Theorems~\ref{thm:computable} and~\ref{thm:local-class-intro}, we prove the following:

\begin{thm}\label{thm:nYThetaSF}
Let $Y=S_{+1}^3(-2T_{6,7}\# T_{6,13}\# T_{-2,3;2,5})$. Then $n [Y]$ is non-trivial in $\Theta_{\Z}^3/\Theta_{\SF}$ for any non-zero $n\in \Z$. 
\end{thm}

Note that Theorem \ref{thm:SFSintro} follows immediately from Theorem \ref{thm:nYThetaSF}. The proof of Theorem \ref{thm:nYThetaSF} relies on the fact from \cite[Theorem 8.1]{DHSThomcob} (cf. the proof of \cite[Theorem 1.1]{DaiConnected}) that the almost-local equivalence class of linear combinations of Seifert fibered spaces takes on a particularly simple form. Note that the manifold $Y$ in Theorem \ref{thm:nYThetaSF} is somewhat ad hoc, and does not immediately generalize to show that $\Theta_{\Z}^3/\Theta_{\SF}$ is infinitely generated. We plan to address this question in future work.

\begin{rem}
More generally, one can consider $\Theta_{\mathit{AR}}$, the subgroup of $\Theta^3_\Z$ generated by almost-rationally plumbed 3-manifolds. (For the precise definition of an almost-rational plumbing, see \cite{NemethiAR}.) Note that Seifert fibered spaces are almost-rationally plumbed. By \cite[Theorem 1.1]{dai-stoffregen}, $\hat{h}(\Theta_{\SF}) = \hat{h}(\Theta_{\mathit{AR}})$, so the proof of Theorem \ref{thm:nYThetaSF} actually shows that  $n [Y]$ is non-trivial in $\Theta_{\Z}^3/\Theta_{\mathit{AR}}$; that is, almost-rationally plumbed 3-manifolds do not generate $\Theta^3_\Z$.
\end{rem}

Theorem \ref{thm:nYThetaSF} also has the following corollary, which was suggested to us by Tye Lidman.

\begin{cor}\label{cor:Montesinos} The Montesinos knots do not generate the smooth concordance group $\mathcal C$.
\end{cor}

\begin{proof} Recall that the double branched cover of a Montesinos knot is a Seifert fibred space, and that furthermore the operation of taking branched double covers of knots sends connected sums of knots to connected sums of 3-manifolds and concordant knots to spin rational homology cobordant 3-manifolds. The local equivalence class and almost local equivalence class of the pair $(\CF^-(Y), \iota)$ for a 3-manifold $Y$ are invariant under not only integer homology cobordism but spin rational homology cobordism \cite{HMZConnectedSum, DHSThomcob}. Therefore it suffices to show that the manifold $Y$ of Theorem \ref{thm:nYThetaSF} is a double branched cover. Torus knots are strongly invertible, and connected sums of strongly invertible knots are strongly invertible \cite{Sakumainvertible}. Moreover, cables of strongly invertible knots are strongly invertible, see for example \cite{Watsonsurgical}, so $K=-2T_{6,7}\# T_{6,13}\# T_{-2,3;2,5}$ is a strongly invertible knot. Therefore $Y= S^3_{+1}(K)$ is a double branched cover $\Sigma(K')$ of some knot $K'$ \cite{Mdbc}. Hence $K'$ is not concordant to any connected sum of Montesinos knots.
\end{proof}

\begin{rem}
An alternate proof of Corollary \ref{cor:Montesinos} can be obtained as follows. Lowrance \cite[Theorem 1.1]{Lowrance} (see also \cite[Corollary 1.3]{StipsiczSzaboThickness}) proved that Montesinos knots have thickness at most one, from which it follows that if $K$ is a Montesinos knot, then $\varphi_n(K) = 0$ for $n>2$, where $\varphi_n$ is the concordance invariant from \cite{DHSTmore}. Since $\bigoplus_{n \in \N} \varphi_n \co \mathcal{C} \to \Z^\infty$ is a surjective homomorphism, it follows that Montesinos knots do not generate the concordance group.
\end{rem}

\begin{rem} A similar computation to the proof of Theorem \ref{thm:nYThetaSF} shows the following. We recall from \cite[Remark 9.7]{HLInvolutiveBordered} that a homology sphere $Y$ is \emph{$\widehat{\HFI}$-trivial} if $\widehat{\HFI}(Y)\iso \widehat{\HF}(Y)\oplus \bF$, where $Q$ vanishes on $\widehat{\HF}$ and is non-vanishing on the remaining generator.  Heretofore, no 3-manifolds have been shown to be $\widehat{\HFI}$-non-trivial. In Section~\ref{sec:surgeryalongK}, we compute that $+1$ surgery on $-2T_{4,5}\# T_{4,9}$ is $\widehat{\HFI}$-non-trivial. \end{rem}

\subsection{The involutive hypercube of a surgery triad}

Theorems~\ref{thm:exact-sequence-1}, \ref{thm:exact-sequence-2} and \ref{thm:cone-intro} all follow from the existence of a certain \emph{hypercube of chain complexes}. See Section~\ref{sec:hypercubes} for the definition of a  hypercube of chain complexes. The relevant version of the hypercube necessary for Theorems~\ref{thm:exact-sequence-2} and \ref{thm:cone-intro} is the following:

\begin{thm} \label{thm:main-hypercube}  Suppose that $K$ is a knot in an integer homology 3-sphere $Y$, $m,n\in \Z$ and $m\ge 1$. Then there is a hypercube of chain complexes of the following form
\[
\begin{tikzcd}[column sep=1.5cm, row sep=1.5cm, labels=description]
\bCF^-(Y_n(K))
	\arrow[drr, dashed, "h"]
	\arrow[dr,"f_1"]
	\arrow[dd, "\iota_1"]
	\arrow[dddrr,dotted]
\\[-.7cm]
&[-.8cm] \bCF^-(Y_{n+m}(K))
	\arrow[r,"f_2"]
	\arrow[ddr,dashed]
&[.8cm] \bigoplus_{i=1}^m \bCF^-(Y)
	\arrow[dd, "\iota_3"]
\\ 
\bCF^-(Y_n(K))
	\arrow[dr,"f_1"]
	\arrow[drr, dashed, "h"]
\\[-.7cm]
& \bCF^-(Y_{n+m}(K))
	\arrow[r,"f_2"]
	\arrow[from=uu, crossing over, "\iota_2"]
	\arrow[from=uuul,crossing over, dashed]
&
\bigoplus_{i=1}^m \bCF^-(Y)
\end{tikzcd}
\]
Furthermore, the maps $f_1$, $f_2$ and $h$ are the ones appearing in Ozsv\'{a}th and Szab\'{o}'s exact sequence \cite{OSIntegerSurgeries}.
\end{thm}

In Theorem~\ref{thm:main-hypercube}, the complexes $\bCF^-(Y_n(K))$ and $\bCF^-(Y_{n+m}(K))$ are the ordinary Heegaard Floer complexes, summed over all $\Spin^c$ structures, and with coefficients in the power series ring $\bF\llsquare U\rrsquare$. The maps $\iota_1$ and $\iota_2$ are chain homotopic to the ordinary involutions on Heegaaard Floer homology. The map $\iota_3$ interchanges the summands of $\bigoplus_{i=1}^m \bCF^-(Y)$, but can otherwise be identified with the involution on $\bCF^-(Y)$.

\subsection{The doubling model of the involution}

\begin{figure}[H]
	\centering
\begingroup%
  \makeatletter%
  \providecommand\color[2][]{%
    \errmessage{(Inkscape) Color is used for the text in Inkscape, but the package 'color.sty' is not loaded}%
    \renewcommand\color[2][]{}%
  }%
  \providecommand\transparent[1]{%
    \errmessage{(Inkscape) Transparency is used (non-zero) for the text in Inkscape, but the package 'transparent.sty' is not loaded}%
    \renewcommand\transparent[1]{}%
  }%
  \providecommand\rotatebox[2]{#2}%
  \newcommand*\fsize{\dimexpr\f@size pt\relax}%
  \newcommand*\lineheight[1]{\fontsize{\fsize}{#1\fsize}\selectfont}%
  \ifx\svgwidth\undefined%
    \setlength{\unitlength}{181.50578471bp}%
    \ifx\svgscale\undefined%
      \relax%
    \else%
      \setlength{\unitlength}{\unitlength * \real{\svgscale}}%
    \fi%
  \else%
    \setlength{\unitlength}{\svgwidth}%
  \fi%
  \global\let\svgwidth\undefined%
  \global\let\svgscale\undefined%
  \makeatother%
  \begin{picture}(1,0.52196068)%
    \lineheight{1}%
    \setlength\tabcolsep{0pt}%
    \put(0,0){\includegraphics[width=\unitlength,page=1]{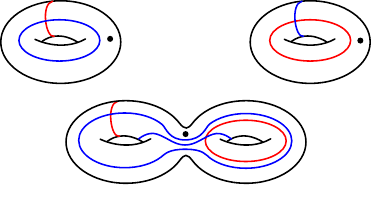}}%
    \put(0.03419609,0.27992072){\color[rgb]{0,0,0}\makebox(0,0)[t]{\lineheight{1.25}\smash{\begin{tabular}[t]{c}$\cH$\end{tabular}}}}%
    \put(0.93126085,0.2569635){\color[rgb]{0,0,0}\makebox(0,0)[t]{\lineheight{1.25}\smash{\begin{tabular}[t]{c}$\bar{\cH}$\end{tabular}}}}%
    \put(0.49044565,0.0044334){\color[rgb]{0,0,0}\makebox(0,0)[t]{\lineheight{1.25}\smash{\begin{tabular}[t]{c}$D(\cH)$\end{tabular}}}}%
    \put(0,0){\includegraphics[width=\unitlength,page=2]{fig194.pdf}}%
  \end{picture}%
\endgroup%

	\caption{Realizing the involution on $S^3$ by doubling.
	 }\label{fig:194}
\end{figure}

 The involution in \cite{HMInvolutive} is given using an arbitrary sequence of Heegaard moves from $\cH$ to $\bar{\cH}$. In this paper, we make use of a special sequence of moves connecting $\cH$ and $\bar{\cH}$, which we call the \emph{doubling} model of the involution. The key observation is that if $\Sigma$ is a Heegaard splitting of $Y$, containing a basepoint $w$, then we can construct another Heegaard splitting $D(\Sigma)$ of $Y$, with  $D(\Sigma)\iso \Sigma\# \bar{\Sigma}$, which is embedded as the boundary of a regular neighborhood of $\Sigma\setminus N(w)$. We think of $D(\Sigma)$ as lying half-way between $\Sigma$ and $\bar{\Sigma}$. Attaching curves for doubled diagrams are described in Section~\ref{sec:doubled-diagrams}.  A schematic example is shown in Figure~\ref{fig:194}.  One helpful property of the doubling operation is that if $\cH$ is a diagram for $Y$, and $D(\cH)$ is a double, then the maps from naturality relating $\CF^-(\cH)$ and  $\CF^-(D(\cH))$ have a conceptually simple form, and similarly for the maps from $\CF^-(D(\cH))$ to $\CF^-(\bar{\cH})$; see Proposition~\ref{prop:doubling-transition-formula}.

In the context of Heegaard Floer homology, doubled diagrams have been studied in \cite{ZemDuality}, \cite{JZContactHandles} and \cite{JZStabilizationDistance}. We note that a doubled Heegaard diagram can also be thought of as a diagram obtained by gluing together two bordered Heegaard diagrams \cite{LOTBordered} for genus $g$ handlebodies. 

Finally, we note that the doubling model  has an interesting topological relation to work of Hass, Thompson and Thurston; see  \cite{HTT-Heegaard-splitting}*{p. 2032}. If $\Sigma$ is a Heegaard surface in $Y$, it follows from the Reidemeister--Singer theorem that $\Sigma$ and $\bar{\Sigma}$ become isotopic as oriented surfaces after some number of stabilizations. The doubling model shows that $g=g(\Sigma)$ stabilizations is always sufficient. According to \cite{HTT-Heegaard-splitting}*{Theorem~1.1}, for each $g>1$, there is a closed, oriented 3-manifold with a genus $g$ Heegaard surface $\Sigma$ which does not become isotopic to $\bar{\Sigma}$ until at least $g$ stabilizations.
 
 \subsection{Organization}
 The paper is organized as follows. We begin by recalling some background on involutive Heegaard Floer homology and Heegaard Floer homology with twisted coefficients in Section \ref{sec:InvolutiveHF}. In Section \ref{sec:iota-complexes}, we describe the algebraic framework of $\iota$-complexes and $\iota_K$-complexes, and define an algebraic operation which takes an algebraic $\iota_K$-complex $\scC$, and gives an algebraically defined model for the involutive mapping cone $\XI_n^{\alg}(\scC)$. In Section \ref{sec:knotcomputations}, we use the involutive mapping cone formula to prove Theorem~\ref{thm:nYThetaSF}. In Section \ref{sec:hypercubes}, we present some background on homological algebra, in particular hyperboxes of chain complexes. Next, in Section \ref{sec:doubled-diagrams}, we describe the operation of doubling a Heegaard diagram, use this operation to give a conceptually simple formula for the involution, and sketch the construction of the main hypercube used to prove Theorem \ref{thm:main-hypercube}. In Section \ref{sec:moduli-spaces}, we deal with some of the analytic input, and in Section \ref{sec:degeneratingconnectedsumtubes}, we describe a neck-stretching technique that will be particularly helpful to us. In Section \ref{sec:stabilization}, we describe several  stabilization operations for Heegaard diagrams, triples and quadruples, and in Section~\ref{sec:stabilizations-holomorphic}, we describe their interaction with the counts of holomorphic curves. In Section \ref{sec:holotriangles}, we consider Heegaard triples where one set of attaching curves is obtained via a small isotopy of another and show that in favorable circumstances, the counts of holomorphic triangles are appropriately simple. In Section \ref{sec:previous-results},  we review some background on Ozsv\'{a}th and Szab\'{o}'s surgery exact sequence \cite{OSIntegerSurgeries} and describe the maps which appear along the top and bottom faces of our main hypercube. Sections \ref{sec:hypercubes-almost-complex-structure}--\ref{sec:flip-hypercube} deal with the construction of different pieces of the main hypercube, and Section \ref{sec:N-null-homotopy} deals with a technical detail from Section~\ref{sec:constructing-C-cen}. In Section \ref{sec:involutive-cone}, we use the involutive hypercube to prove Theorem \ref{thm:cone-intro} and Theorem \ref{thm:exact-sequence-2}. In Section \ref{sec:strong-cone-formula}, we prove that for a knot $K$ in L-space integer homology sphere $Y$, the algebraic involutive mapping cone from Section \ref{sec:iota-complexes} computes $\HFI^-(Y_n^3(K))$. Finally, in Section \ref{sec:rational}, we consider rational surgery as well as 0-surgery.

 \subsection{Acknowledgements}
 We would like to thank Sherry Gong, Tye Lidman, Robert Lipshitz, Ciprian Manolescu, and JungHwan Park for helpful conversations. We would like to thank Danny Ruberman for helpful comments on an earlier draft. We are grateful to the Frontiers in Floer Homology Workshop at Boston College in August 2019, the University of Oregon Floer Homotopy Workshop in August 2019, and the Geometry and Topology Workshop at UCLA in January 2020, where parts of this work were completed. Finally, we are grateful to the anonymous referees for their thoughtful reading and helpful comments.

\section{Involutive Heegaard Floer homology}\label{sec:InvolutiveHF}

In this section, we recall background on \emph{involutive Heegaard Floer homology}, due to the first author and Manolescu \cite{HMInvolutive}, as well as some variations which are relevant to the exact sequences and knot surgery formula.

\subsection{3-manifolds} \label{sec:3-manifolds}

If $\cH=(\Sigma,\as,\bs,w)$ is a Heegaard diagram of a 3-manifold $Y$, then $\bar{\cH}:=(\bar{\Sigma}, \bs, \as,w)$ is also a Heegaard diagram for $Y$, where $\bar{\Sigma}$ denotes $\Sigma$, with its orientation reversed. We say $\cH$ and $\bar{\cH}$ are \emph{conjugate}. If $\varpi\subset \Spin^c(Y)$ denotes a conjugacy class of $\Spin^c$ structures, and $\cH$ is a Heegaard diagram, we write
\[
\CF^-(\cH,\varpi):=\bigoplus_{\frs\in \varpi} \CF^-(\cH,\frs).
\]
Note that $\varpi$ either consists  of one $\Spin^c$ structure $\frs$ satisfying $\frs=\bar\frs$, or is a 2-element set of the form $\{\frs,\bar\frs\}$.

There is a tautological chain isomorphism
\[
\eta\colon \CF^-\left(\bar{\cH},\bar{\frs}\right)\to \CF^-\left(\cH, \frs\right). 
\]
There is additionally a chain homotopy equivalence 
\[
\Psi_{\cH\to \bar{\cH}}\colon \CF^-\left(\cH,\frs\right)\to \CF^-\left(\bar\cH,\frs\right),
\]
obtained by picking a set of Heegaard moves connecting $\bar{\cH}$ and $\cH$, and using the induced map from naturality.

Hendricks and Manolescu study the involution $\iota\colon \CF^-(\cH,\varpi)\to \CF^-(\cH,\varpi)$ defined as the composition $\iota:=\eta\circ \Psi_{\cH\to \bar{\cH}}$.  They define the involutive Heegaard Floer chain complex to be
\[
\CFI^-(Y,\varpi):=\Cone\left(Q\cdot (\iota+\id)\colon \CF^-(\cH,\varpi)\to Q\cdot\CF^-(\cH,\varpi)\right),
\]
where $Q$ is a formal variable, which we view as satisfying $Q^2=0$.

Note that if $\varpi=\{\frs,\bar\frs\}$ where $\frs\neq \bar\frs$, then 
\[
\HFI^-(Y,\varpi)\iso \HF^-(Y,\frs)[-1]\oplus \HF^-(Y,\bar\frs).
\]
See \cite{HMInvolutive}*{Proposition~4.5}. For a conjugation invariant $\Spin^c$ structure, Hendricks and Manolescu extract two involutive correction terms $\dl(Y, \frs)$ and $\du(Y, \frs)$ defined via
\[ \dl(Y,\frs) = \max \{r \mid \exists \ x \in \HFI^-_r(Y, \frs), \forall \ n, \ U^nx\neq 0 \ \text{and} \ U^nx \notin \operatorname{Im}(Q)\}-1 \]
and
\[ \du(Y,\frs) = \max \{r \mid \exists \ x \in \HFI^-_r(Y,\frs),\forall \ n, U^nx\neq 0; \exists \ m\geq 0 \ \operatorname{s.t.} \ U^m x \in \operatorname{Im}(Q)\}. \]

Involutive Heegaard Floer homology also has a connected sum formula \cite{HMZConnectedSum}, extending Ozsv\'{a}th and Szab\'{o}'s connected sum formula for ordinary Heegaard Floer homology \cite{OSProperties}*{Theorem~1.5}. If $Y_1$ and $Y_2$ are two 3-manifolds, equipped with self-conjugate $\Spin^c$ structures $\frs_1$ and $\frs_2$, then 
\begin{equation}
(\CF^-(Y_1\# Y_2, \frs_1\# \frs_2),\iota_{Y_1\# Y_2})\simeq (\CF^-(Y_1)\otimes \CF^-(Y_2), \iota_{Y_1}\otimes \iota_{Y_2}),
\label{eq:connected-sums-3-manifolds}
\end{equation} where $\simeq$ denotes homotopy equivalence of iota-complexes.

\subsection{Knots}
\label{sec:involutive-knots}
If $K$ is a null-homologous knot in an integer homology 3-sphere $Y$, Hendricks and Manolescu also describe a skew-filtered homotopy automorphism $\iota_K\colon \CFK^\infty(Y,K)\to \CFK^\infty(Y,K)$ \cite{HMInvolutive}*{Section~6.1}, which we describe presently.

We recall that a Heegaard diagram for $(Y,K)$ consists of a tuple $\cH=(\Sigma,\as,\bs,w,z)$ such that the following hold:
\begin{enumerate}[label=(\textit{KD}-\arabic*), ref=\textit{KD}-\arabic*, leftmargin=*]
\item\label{knot-diagram-1} $(\Sigma,\as,\bs)$ is a Heegaard diagram for $Y$. In particular $Y\setminus \Sigma$ is the disjoint union of two handlebodies $U_{\a}$ and $U_{\b}$.
\item\label{knot-diagram-2} $K\cap \Sigma=\{w,z\}$, and furthermore $K$ intersects $\Sigma$ positively at $z$ and negatively at $w$.
\item\label{knot-diagram-3} The tangles $K\cap U_{\a}$ and $K\cap U_{\b}$ are boundary parallel.
\end{enumerate}

If $\cH=(\Sigma,\as,\bs,w,z)$ is a diagram for $(Y,K)$, then $\bar{\cH}:=(\bar{\Sigma}, \bs,\as,z,w)$ is also a diagram for $(Y,K)$. There is a tautological isomorphism
\[
\eta_K\colon \CFK^\infty(\bar{\cH})\to \CFK^\infty(\cH)
\]
given by $\eta_K([\xs,i,j])=[\xs,j,i]$.

The two diagrams $\bar{\cH}$ and $\cH$ both represent $(Y,K)$. However, there is an additional subtlety that the roles of the basepoints have been reversed, and the naturality theorem of \cite{JTNaturality} holds only for \emph{based} knots and links. Nonetheless, we may pick an automorphism $\psi$ of $(Y,K)$, which is supported in a neighborhood of $K$, and which switches $w$ and $z$ via a half twist along $K$. Then $\psi_*(\cH)$ and $\bar{\cH}$ are both diagrams for $(Y,K)$, with the same basepoints. Hence there is a naturality map
\begin{equation}
\Psi_{\psi_*(\cH)\to \bar{\cH}}\colon \CFK^\infty(\psi_*(\cH))\to \CFK^\infty(\bar{\cH}).\label{eq:naturality-map-knots}
\end{equation}
Hendricks and Manolescu define the knot involution $\iota_K$ as
\[
\iota_K:=\eta_K \circ \Psi_{\psi_*(\cH)\to \bar{\cH}}\circ \psi_*.
\] 
Note that technically $\iota_K$ is not an involution, since $\iota_K^2$ is homotopic to  Sarkar's map $\id+\Phi\Psi$ (cf. \cite{SarkarMovingBasepoints}); however, it is a convention to refer to $\iota_K$ as the knot involution.

One key property established by Hendricks and Manolescu is the relation of the involutive knot Floer complex to the involutive Heegaard Floer complex of surgeries. If $K$ is a knot in $S^3$, then according to \cite{HMInvolutive}*{Theorem~1.5}, $(\CF^-(S^3_n(K),[0]),\iota_{S^3_n(K)})$ is equivariantly homotopy equivalent to the iota-complex $(A_0(K),\iota_K)$ for $n\gg 0$ (in fact, $n\ge g_3(K)$ is sufficient). Here $[0]$ represents the conjugation invariant $\Spin^c$ structure identified with $[0]$ in the usual identification of $\Spin^c$ structures on $S^3_n(K)$ with $\Z/n\Z$.

From involutive knot Floer homology, Hendricks and Manolescu extract two concordance invariants $\Vl_0(K)$ and $\Vu_0(K)$, which may be defined via the formula
$$\dl(S^3_n(K), [0]) = \frac{n-1}{4} -2\Vl_0(K), \ \ \ \ \du(S^3_n(K), [0])=\frac{n-1}{4} -2\Vu_0(K).$$
Finally, involutive knot Floer homology enjoys a connected sum formula, similar to ~\eqref{eq:connected-sums-3-manifolds}. If $K_1$ and $K_2$ are two knots in $S^3$, then the knot involution $\iota_{K_1\# K_2}$ on $\CFK^\infty(K_1\# K_2)\iso \CFK^\infty(K_1)\otimes \CFK^\infty(K_2)$ is homotopic to 
\[
\iota_1\otimes \iota_2+\Phi\iota_1\otimes \Psi\iota_2,
\]
where $\Phi$ and $\Psi$ are the two \emph{formal derivative} maps, which are endomorphisms of knot Floer homology. See \cite{ZemConnectedSums}*{Theorem~1.1} for more details.

\subsection{Twisted 3-manifold complexes}
\label{sec:twisted-complexes-def}

If $K$ is a knot in a 3-manifold $Y$, then there is a twisted version of the complex $\CF^-(Y)$, which naturally appears in the proof of the mapping cone formula \cite{OSIntegerSurgeries}. In this section, we describe its definition, as well as define an involution. To simplify the discussion, we focus on the case that $Y$ is an integer homology sphere. We will revisit the construction for rational homology 3-spheres in Section~\ref{sec:rational}, when we consider rational surgeries.

\begin{define} 
Suppose $(K,w,z)$ is a doubly pointed knot. An \emph{enhanced orientation} $\hat{\fro}=(\fro,p)$ on $K$ consists of an orientation $\fro$ on $K$, together with a choice of special basepoint $p\in \{w,z\}$.
\end{define}
We write $\cH_{\hat{\fro}}$ or $\cH_{\fro,p}$ for a doubly pointed Heegaard diagram, equipped with an enhanced orientation of the underlying knot.

Suppose $\cH=(\Sigma,\as,\bs,w,z)$ is a diagram for $K$ in $Y$, and further, suppose $K$ is equipped with an enhanced orientation $\hat{\fro}$. We now describe the complex $\buCF^-(\cH_{\hat{\fro}})$. We define $\buCF^-(\cH_{\hat{\fro}})$ to be the free $\bF\llsquare U\rrsquare\otimes \bF[\Z/m]$-module generated by intersection points $\ve{x}\in \bT_{\a}\cap \bT_{\b}$. If $a\in \Z/m$, write $T^a$ for the corresponding element of the group ring $\bF[\Z/m]$. We equip $\buCF^-(\cH_{\hat{\fro}})$ with the differential
\begin{equation}
\d \xs=\sum_{\ys\in \bT_{\a}\cap \bT_{\b}} \sum_{\substack{\phi\in \pi_2(\xs,\ys)\\ \mu(\phi)=1}} \# (\cM(\phi)/\R)\cdot  U^{n_p(\phi)} T^{n_{z}(\phi)-n_w(\phi)}\cdot \ys, \label{eq:twisted-CF-def}
\end{equation}
where $p\in \{w,z\}$ is the special point of $\hat{\fro}$. Note that the ordering of the basepoints in the exponent of $T$ is induced by the orientation of $K$. More precisely, $z$ is the basepoint of $K$ corresponding to a positive intersection of $K$ and $\Sigma$, and $w$ is the basepoint of $K$ corresponding to the negative intersection.

It is convenient also to have some additional notation. We will also write $(\Sigma,\as,\bs,\hat{w},z)$ or $(\Sigma,\as,\bs,w,\hat{z})$ for a diagram of a knot with an enhanced orientation. The order of the basepoints indicates the orientation, through the convention in \eqref{knot-diagram-2}, and the hat indicates the special basepoint.

There is a trivialization
\begin{equation}
\theta_p\colon \buCF^-(\cH_{\fro,p})\to \bCF^-(\Sigma,\as,\bs,p)\otimes \bF[\Z/m]
\label{eq:theta-w-def}
\end{equation}
given by
\[
\theta_p(T^s \cdot \xs)=T^{s+A_{w,z}(\xs)}\cdot  \xs.
\]
extended equivariantly over the action of $U$. Here, $A_{w,z}(\xs)$ denotes the absolute Alexander grading of $\xs$, as an element of $\widehat{\CFK}(\Sigma,\as,\bs,w,z)$. (Recall that the relative Alexander grading is given by $A_{w,z}(\xs,\ys)=n_z(\phi)-n_w(\phi)$ where $\phi$ is any class in $\pi_2(\xs,\ys)$.)

\begin{lem}
The map $\theta_p$ is a chain isomorphism.
\end{lem}
\begin{proof} The map $\theta_p$ is obviously an isomorphism of modules, so it is sufficient to show that it is a chain map. Suppose that $\phi\in \pi_2(\xs,\ys)$ is a class with Maslov index 1. Then $\phi$ contributes a summand of
\begin{equation}
\# (\cM(\phi)/\R)\cdot U^{n_p(\phi)}T^{A_{w,z}(\ys)+n_z(\phi)-n_{w}(\phi)}\cdot \ys \label{eq:theta-w-chain-iso-1}
\end{equation}
to $(\theta_p\circ \d)(\xs)$, and a summand of
\begin{equation}
\# (\cM(\phi)/\R)\cdot U^{n_p(\phi)}T^{A_{w,z}(\xs)}\cdot \ys\label{eq:theta-w-chain-iso-2}
\end{equation}
to $(\d \circ \theta_p)(\xs)$. Equations~\eqref{eq:theta-w-chain-iso-1} and~\eqref{eq:theta-w-chain-iso-2} are clearly equal, so $\theta_p\d+\d\theta_p=0$.
\end{proof}

Writing $\bar{\cH}=(\bar{\Sigma},\bs,\as,z,w)$, there is a tautological map
\begin{equation}
\underline{\eta}\colon \buCF^-(\bar{\cH}_{\hat{\fro}})\to \buCF^-(\cH_{\hat{\fro}}), \label{eq:tautoloical-eta-def-twisted}
\end{equation}
given by the formula
\[
\underline{\eta}(U^i T^j\cdot \xs)=U^i T^{-j} \xs.
\]
It is straightforward to see that $\underline{\eta}$ is a chain isomorphism. Note that in~\eqref{eq:tautoloical-eta-def-twisted}, the diagram $\bar{\cH}_{\hat{\fro}}$ denotes $\bar{\cH}$ decorated with the enhanced orientation $\hat{\fro}$. In particular, the enhanced orientation is the same in the domain and codomain of $\underline{\eta}$.

The last ingredient of the involution is the \emph{flip-maps}, which change the special point of an enhanced orientation. These take the form of maps
\[
\begin{split}
\frF_{w\to z}&\colon \buCF^-(\cH_{\fro,w})\to \buCF^-(\cH_{\fro,z})\quad \text{and}\\
\frF_{z\to w}&\colon \buCF^-(\cH_{\fro,z})\to \buCF^-(\cH_{\fro,w}).
\end{split}
\]
which we describe presently.

There is diffeomorphism map 
\[
\phi_*\colon \bCF^-(\Sigma,\as,\bs,z)\to \bCF^-(\Sigma,\as,\bs,w),
\]
obtained by moving $z$ to $w$ along a subarc of $K$, and then composing with the change of diagram maps from naturality. In principle, there is a dependence on the choice of subarc, however in all cases we will make a specific choice. The flip map $\frF_{z\to w}$ is defined via the formula
\begin{equation}
\frF_{z\to w}:=T^n\cdot\theta_w^{-1} \circ \left(\phi_*\otimes  \id_{\bF[\Z/m]}\right)\circ  \theta_z.
\label{eq:flipmap-def-z->w}
\end{equation}
(Here $n$ is the surgery coefficient that appears in the exact sequence). The other flip map $\frF_{w\to z}$ is defined via the formula
\begin{equation}
\frF_{w\to z}:=T^{-n}\cdot\theta_z^{-1} \circ \left((\phi^{-1})_*\otimes  \id_{\bF[\Z/m]}\right)\circ  \theta_{w}.
\label{eq:flipmap-def-w->z}
\end{equation}

Finally,  we define the involution $\underline{\iota}$ of $\buCF^-(\Sigma,\as,\bs,w,z)$ via the formula
\begin{equation}
\underline{\iota}:=\frF_{z\to w}\circ   \underline{\eta}\circ \Psi_{\psi(\cH)\to \bar{\cH}}\circ \psi_*.
 \label{eq:def-iota-3}
\end{equation}

We note that in Theorem~\ref{thm:main-hypercube}, we are identifying $\bigoplus_{i=1}^m \bCF^-(Y)$ with $\buCF^-(\cH_{\fro,w})$, for a Heegaard diagram $\cH$ of $(Y,K)$, as well as a fixed orientation  $\fro$ on $K$. Furthermore, the map labeled $\iota_3$ therein is the map $\underline{\iota}$ in this notation.

\section{Iota complexes and mapping cones}
\label{sec:iota-complexes}
In this section, we describe the algebraic framework of $\iota$-complexes and $\iota_K$-complexes, and define an algebraic operation which takes an algebraic $\iota_K$-complex $\scC$, and gives an algebraically defined model for the involutive mapping cone $\XI_n^{\alg}(\scC)$. In this section, we focus on the case that $n$ is a non-zero integer. The construction has natural analogs for rational surgeries and 0-surgeries. See Sections~\ref{subsec:rational-surgeries} and~\ref{subsec:0-surgery}.

\subsection{\texorpdfstring{$\iota$}{iota}-complexes}

We begin by defining $\iota$-complexes and $\iota_K$-complexes. Compare the following to \cite{HMZConnectedSum}*{Definition~8.1}:

\begin{define}
\label{def:iota-complex}
 An \emph{$\iota$-complex} is a chain complex $(C,\d)$, which is free and finitely generated over $\bF[U]$, equipped with an endomorphism $\iota$. Furthermore, the following hold:
\begin{enumerate}
\item\label{iota-1} $C$ is equipped with a $\Z$-grading, such that $U$ has grading $-2$. We call this grading the \emph{Maslov} or \emph{homological} grading.
\item\label{iota-2} There is a grading preserving isomorphism $U^{-1}H_*(C)\iso \bF[U,U^{-1}]$.
\item \label{iota-3} $\iota$ is a grading preserving chain map and $\iota^2\simeq \id$.
\end{enumerate}
\end{define}

\begin{rem}
We will usually abuse notation and refer to the involutive mapping cone $\XI_{n}^{\alg}(\scC)=(\bX_n(\scC),\iota_{\bX})$, constructed in Section~\ref{sec:canonical-mapping-cone}, as an $\iota$-complex. However, the reader should recall that $\bX_n$ is defined over $\bF\llsquare U\rrsquare$, is infinitely generated, and also satisfies $U^{-1} H_*( \bX_n)\iso \bigoplus^n \bF\llsquare U, U^{-1}]$.
\end{rem}

\begin{define}\label{def:iota-homomorphism}
If $(C,\iota)$ and $(C',\iota')$ are two $\iota$-complexes, an \emph{$\iota$-homomorphism} is a homogeneously graded, $\bF[U]$-equivariant chain map $F\colon C\to C'$ such that $F\iota+\iota' F\simeq 0$. Two $\iota$-complexes $(C,\iota)$ and $(C',\iota')$ are \emph{$\iota$-homotopy equivalent} if there exist  $\iota$-homomorphisms $F \colon (C, \iota) \to (C', \iota')$ and $G \colon (C', \iota') \to (C, \iota)$ such that $FG \simeq \id$ and $GF \simeq \id$.
\end{define}

We write $\Iota$ for the category whose objects are $\iota$-complexes, and whose morphisms are $\iota$-homomorphisms. Note that $\Iota$ is naturally a preadditive category, but is not naturally a dg-category. There is a related category $\underline{\Iota}$, which we call the category of $\iota$-complexes with \emph{enhanced morphisms}, as follows:

\begin{define}\label{def:enhanced-iota-homomorphisms}
If $(C,\iota)$ and $(C',\iota')$ are $\iota$-complexes, we define the group of \emph{enhanced $\iota$-morphisms} to be
\[
\underline{\Mor}((C,\iota),(C',\iota')):=\Hom_{\bF[U]}(C,C')\oplus \Hom_{\bF[U]}(C,C')[1],
\]
where $[1]$ denotes a grading shift. We define
\[
\d_{\underline{\Mor}}(F,h)=(\d' F+F\d, F\iota+\iota'F+\d' h+h\d),
\]
which makes $\underline{\Iota}$ into a dg-category.
 An enhanced $\iota$-\emph{homo}morphism is an enhanced $\iota$-morphism $(F,h)$ satisfying $\d_{\underline{\Mor}}(F,h)=0$. An \emph{enhanced $\iota$-homotopy equivalence} is an enhanced $\iota$-homomorphism $(F,h)$, such that there exists an enhanced $\iota$-homomorphism $(G,k)$ such that $(G,k)\circ (F,h)$ and $(F,h)\circ (G,k)$ both differ from $(\id,0)$, the identity enhanced $\iota$-morphism, by a boundary in $\underline{\Mor}$.
\end{define}

\begin{rem}
\begin{enumerate}
\item Note that an enhanced $\iota$-homomorphism $(F,h)$ consists of a chain map $F \colon C \to C'$ such that $F \iota + \iota' F \simeq 0$ via the chain homotopy $h$.
\item We will see later that an enhanced $\iota$-homomorphism $(F,h)$ is an enhanced $\iota$-homotopy equivalence if and only if $F$ is an ordinary homotopy equivalence of chain complexes over $\bF[U]$. See Lemma~\ref{lem:homotopy-equivalence-iota-complexes}.
\end{enumerate}
\end{rem}

Composition of enhanced $\iota$-morphisms is given by
\[
(F,h)\circ (G,k)=(F G, hG+Fk).
\]
The group of enhanced $\iota$-morphisms $(F,h)$ is canonically isomorphic to
\begin{equation}
 \Hom_{\bF[U,Q]/Q^2}\left(\Cone(C\xrightarrow{Q(1+\iota)} Q C),\Cone(C'\xrightarrow{Q(1+\iota')} Q C')\right)
 \label{eq:enhanced-morphism-equivalent}
\end{equation}
via $(F,h)\mapsto F+ Qh$.

\subsection{Local equivalence} \label{subsec:local}

To study the homology cobordism group, we consider a special type of $\iota$-morphism, called a \emph{local map}:
\begin{define} Suppose $(C,\iota)$ and $(C',\iota')$ are two $\iota$-complexes.
\begin{enumerate}
\item A \emph{local map} from $(C,\iota)$ to $(C',\iota')$ is a grading preserving $\iota$-homomorphism $F\colon C\to C'$, which induces an isomorphism from $U^{-1} H_*(C)$ to $U^{-1} H_*(C')$.
\item We say that $(C,\iota)$ are $(C',\iota')$ are \emph{locally equivalent} if there is a local map from $(C,\iota)$ to $(C',\iota')$, as well as a local map from $(C',\iota')$ to $(C,\iota)$.
\end{enumerate}
\end{define}
The set of local equivalence classes forms an abelian group, denoted $\frI$, with group structure given by
\begin{equation}
(C,\iota)*(C',\iota')=(C\otimes C',\iota\otimes \iota'). \label{eq:group-relation-I}
\end{equation}
See \cite{HMZConnectedSum}*{Section~8}. According to \cite{HMZConnectedSum}*{Theorem~1.8}, the map
\[
Y\mapsto [(\CF^-(Y), \iota)]
\]
 determines a homomorphism from $\Theta_{\Z}^3$ to $\frI$. 

Given an $\iota$-complex $(C,\iota)$, we may define correction terms $d(C)$, $\dl(C)$, and $\du(C)$ as in Section \ref{sec:3-manifolds}, all of which are preserved by local equivalence. 

\subsection{Almost \texorpdfstring{$\iota$}{iota}-complexes and almost local equivalence}

There is an additional, weaker, equivalence relation between iota complexes, as follows.

\begin{define}\cite{DHSThomcob}*{Definition 3.1}
Let $C_1$ and $C_2$ be free, finitely generated chain complexes over $\bF[U]$,  such that each $C_i$ has an absolute $\Q$-grading and a relative $\Z$-grading with respect to which $U$ has grading $-2$. Two grading-preserving $\bF[U]$-module homomorphisms
\[ 
f, g 
\colon
 C_1 \rightarrow C_2
 \]
are \emph{homotopic mod $U$}, denoted $f \simeq g \mod U$, if there exists an $\bF[U]$-module homomorphism $H \colon C_1 \rightarrow C_2$ such that $H$ increases grading by one and
\[
 f + g + H \circ \d + \d \circ H \in \im U.
 \]
\end{define}

\begin{define}[{\cite{DHSThomcob}*{Definition 3.2}}]
\label{def:almost-iota-complex}
 An \emph{almost $\iota$-complex} is a chain complex $(C,\d)$, which is free and finitely generated over $\bF[U]$, equipped with an endomorphism $\bar{\iota}$. Furthermore, the following hold:
\begin{enumerate}
\item\label{almost-iota-1} $C$ has an absolute $\Z$-grading, such that $U$ has grading $-2$.
\item\label{almost-iota-2} There is a grading preserving isomorphism $U^{-1}H_*(C)\iso \bF[U,U^{-1}]$.
\item \label{almost-iota-3} $\bar{\iota}$ is a grading-preserving, $\bF[U]$-module homomorphism $\bar{\iota} \colon C \rightarrow C$ such that
	\[
	 \bar{\iota} \circ \d + \d \circ \bar{\iota} \in \im U \qquad \text{ and } \qquad \bar{\iota}^2 \simeq \id \mod U. \]
\end{enumerate}
It will often be convenient to consider the map
\[
\omega = 1 + \bar{\iota}.
\]
\end{define}

\noindent Of course, every $\iota$-complex has an associated almost $\iota$-complex. Note that the definition of an almost $\iota$-complex is obtained from the definition of an $\iota$-complex by relaxing every statement about $\iota$ to merely be true mod $U$.

\begin{define}[\cite{DHSThomcob}*{Definition 3.4}]
\label{def:aim}
An \emph{almost $\iota$-homomorphism} from $\cC_1 = (C_1, \bar{\iota}_1)$ to $\cC_2 = (C_2, \bar{\iota}_2)$ is a grading-preserving, $\bF[U]$-equivariant chain map
\[ f \colon C_1 \rightarrow C_2 \]
such that
\[ f \circ \bar{\iota} \simeq \bar{\iota} \circ f \mod U. \]
\end{define}

\begin{define}[\cite{DHSThomcob}*{Definition~3.5}]
\begin{enumerate}
\item If $(C,\bar{\iota})$ and $(C',\bar{\iota}')$ are almost $\iota$-complexes, an \emph{almost local map} is a degree preserving almost $\iota$-homomorphism $F\colon C\to C'$, which induces an isomorphism from $U^{-1} H_*(C)$ to $U^{-1} H_*(C')$.
\item We say that $(C,\bar{\iota})$ are $(C',\bar{\iota}')$ are \emph{almost locally equivalent} if there is an almost local map from $(C,\bar{\iota})$ to $(C',\bar{\iota}')$, as well as an almost local map from $(C',\bar{\iota}')$ to $(C,\bar{\iota})$.
\end{enumerate}
\end{define}

According to \cite{DHSThomcob}*{Theorem~3.25}, the set of almost local equivalence classes forms a totally ordered abelian group, denoted $\hat{\frI}$. As with the local equivalence group, the group structure is given by
\begin{equation}
(C,\bar{\iota})*(C',\bar{\iota}')=(C\otimes C',\bar{\iota}\otimes \bar{\iota}'). \label{eq:group-relation-Ihat}
\end{equation}
There are homomorphisms
\[
\Theta_{\Z}^3 \rightarrow \frI \rightarrow \widehat{\frI}
\] 
The structure of the group $\widehat{\frI}$ is recalled in further detail in Section \ref{sec:knotcomputations}.

\subsection{\texorpdfstring{$\iota_K$}{iota-K}-complexes}

There is a refinement of the notion of an $\iota$-complex for knots, called an $\iota_K$-complex. Before stating the definition, we state some algebraic background. Suppose that $(C_K,\d)$ is a free, finitely generated complex over the ring $\bF[\knotU,\knotV]$. There are two naturally associated maps
\[
\Phi,\Psi\colon C_K\to C_K,
\]
as follows. We write $\d$ as a matrix with respect to a free $\bF[\knotU,\knotV]$-basis of $C_K$. We define $\Phi$ to be the endomorphism obtained differentiating each entry of this matrix with respect to $\knotU$. We define $\Psi$ to be the endomorphism obtained by differentiating each entry with respect to $\knotV$. These maps naturally appear in the context of knot Floer homology, see \cites{SarkarMaslov, ZemQuasi, ZemCFLTQFT}. The maps $\Phi$ and $\Psi$ are independent of the choice of basis, up to $\bF[\knotU,\knotV]$-equivariant chain homotopy \cite{ZemConnectedSums}*{Corollary~2.9}.

We say an $\bF$-linear map $F\colon C_K\to C_K'$ is \emph{skew-$\bF[\knotU,\knotV]$-equivariant} if 
\[
F\circ \knotV=\knotU\circ F\quad \text{and} \quad F\circ \knotU=\knotV\circ F.
\]

We may view a free complex over $\bF[\knotU,\knotV]$ also as an infinitely generated complex over $\bF[U]$, where $U$ acts diagonally via $U=\knotU \knotV$. Concretely, if $B=\{\xs_1,\dots, \xs_n\}$ is an $\bF[\knotU,\knotV]$-basis, then an $\bF[U]$-basis is given by the elements $\knotU^i\cdot  \xs_k$ and $\knotV^j\cdot  \xs_k$, ranging over all $i\ge 0$, $j\ge 0$ and $k\in \{1,\dots, n\}$.

Compare the following to \cite{HMInvolutive}*{Definition~6.2} and \cite{ZemConnectedSums}*{Definition~2.2}.

\begin{define}\label{def:iota-K-complex}
\begin{enumerate}
\item  An \emph{$\iota_K$-complex} $(C_K,\d, \iota_K)$ is a finitely generated, free chain complex $(C_K,\d)$ over $\bF[\knotU,\knotV]$, equipped with a skew-equivariant endomorphism $\iota_K$ satisfying 
 \[
\iota_K^2\simeq \id+\Phi\Psi. 
 \]
\item We say an $\iota_K$-complex $(C_K,\d,\iota_K)$ is of \emph{$\Z H S^3$-type} if there are two $\Z$ valued gradings, $\gr_{\ws}$ and $\gr_{\zs}$, such that $\knotU$  and $\knotV$ have $(\gr_{\ws},\gr_{\zs})$-bigrading $(-2,0)$ and $(0,-2)$, respectively. We assume $\d$ has $(\gr_{\ws},\gr_{\zs})$-bigrading $(-1,-1)$, and that $\iota_K$ switches $\gr_{\ws}$ and $\gr_{\zs}$. Furthermore, we assume that $A:=\tfrac{1}{2}(\gr_{\ws}-\gr_{\zs})$ is integer valued. We call $A$ the \emph{Alexander} grading, and we call $\gr_{\ws}$ and $\gr_{\zs}$ the \emph{Maslov} gradings.  Writing $\cA_s\subset C_K$ for the subspace in Alexander grading $s$, we assume that there is a grading preserving isomorphism $U^{-1} H_*(\cA_s)\iso \bF[U,U^{-1}]$ for all $s\in \Z$. 
\item We say an $\iota_K$-complex $(C_K,\d,\iota_K)$ is of \emph{L-space type} if it is of $\Z H S^3$-type and additionally $H_*(\cB_s)\iso \bF[U]$, where $\cB_s$ is the subset of $\knotV^{-1}C_K$ in Alexander grading $s$.
\end{enumerate}
\end{define}

In Definition~\ref{def:iota-K-complex}, an $\iota_K$-complex of $\Z H S^3$-type is equipped with two Maslov gradings, $\gr_{\ws}$ and $\gr_{\zs}$. We note that in the literature, usually one considers just the $\gr_{\ws}$-grading, which is referred to as the homological grading.

\begin{rem}\label{rem:iotaKU}
If $(C_K,\d,\iota_K)$ is an $\iota_K$-complex, then $\iota_K$ commutes with $U=\knotU\knotV$, and hence we can view $C_K$ as an (infinitely generated) complex with an $\bF[U]$-equivariant endomorphism $\iota_K$. 
\end{rem}

If $K$ is a knot in an integer homology 3-sphere $Y$, then we obtain an $\iota_K$-complex in the sense of Definition~\ref{def:iota-K-complex} by considering the version of the knot Floer complex generated by monomials $\knotU^i \knotV^j \cdot \xs$, for $i,j\ge 0$, where $\knotU$ counts the basepoint $w$ and $\knotV$ counts the basepoint $z$. This contains equivalent information to the more standard version of the knot Floer complex generated by tuples satisfying $[\xs,i,j]$ with $A(\xs)-j+i=0$. The correspondence is given by
\begin{equation}
\knotU^i \knotV^j \cdot \xs\longleftrightarrow [\xs,-i,-j]. \label{eq:correspondence-variable-brackets}
\end{equation}

Analogous to Definitions~\ref{def:iota-homomorphism} and~\ref{def:enhanced-iota-homomorphisms}, we have the following notions of morphism of $\iota_K$-complexes:

\begin{define} Suppose that  $\mathscr{C}=(C_K,\d, \iota_K)$ and $\mathscr{C}'=(C'_K,\d',\iota_K')$ are $\iota_K$-complexes.
\begin{enumerate}
\item  An \emph{$\iota_K$-homomorphism} from $\mathscr{C}$ to $\mathscr{C}'$ consists of an $\bF[\knotU,\knotV]$-equivariant chain map $F\colon C_K\to C_K'$, which satisfies $\iota_K'F+F\iota_K\eqsim 0$ (where $\eqsim$ denotes skew-equivariantly chain homotopy equivalence).
\item The group of \emph{enhanced $\iota_K$-morphisms} is
\[
\underline{\Mor}(\scC,\scC'):=\Hom_{\bF[\knotU,\knotV]}(C_K,C_K')\oplus \bar{\Hom}_{\bF[\knotU,\knotV]}(C_K,C_K')[1,1],
\]
where $\bar{\Hom}_{\bF[\knotU,\knotV]}(C_K,C_K')$ denotes the group of $\bF[\knotU,\knotV]$-skew-equivariant maps. The differential on $\underline{\Mor}(\scC_K,\scC'_K)$ is given by 
\[
\d_{\underline{\Mor}}(F,g)=(F\d+\d'F, F\iota_K+\iota_K' F+\d' g+g\d).
\]
We say $(F,g)$ is an enhanced $\iota_K$-\emph{homo}morphism if $\d_{\underline{\Mor}}(F,g)=0$. Two enhanced $\iota_K$-morphisms are $\iota_K$-\emph{homotopic} if their sum is a boundary in $\underline{\Mor}(\scC_K,\scC_K')$.
\end{enumerate}
\end{define}

We write $\Iota_K$ for the category of $\iota_K$-complexes, with morphism set equal to the set of $\iota_K$-homomorphisms. We write $\underline{\Iota}_K$ for the category of $\iota_K$-complexes with enhanced $\iota_K$-morphisms.

There is also a natural notion of local equivalence of $\iota_K$-complexes:
\begin{define} Suppose that $\scC=(C_K,\d,\iota_K)$ and $\scC'=(C_K',\d',\iota_K')$ are $\iota_K$-complexes of $\Z HS^3$-type.
\begin{enumerate}
\item A local map from $\scC$ to $\scC'$ is a $(\gr_{\ws},\gr_{\zs})$-preserving, $\bF[\knotU,\knotV]$-linear chain map $F\colon C_K\to C_K'$, which satisfies $F\iota_K+\iota_K' F\eqsim 0$, and which induces an isomorphism from $H_*( U^{-1} \cA_s)$ to $H_*(U^{-1} \cA'_s)$ for all $s$.
\item We say $\scC$ and $\scC'$ are \emph{locally equivalent} if there are local maps in both directions.
\end{enumerate}
\end{define}

According to \cite{ZemConnectedSums}*{Theorem~1.5}, the set of local equivalence classes of $\iota_K$-complexes forms a group, denoted $\frI_K$. The group operation is given by
\[
(C_K,\d,\iota_K)* (C_K',\d',\iota_K')\mapsto (C_K\otimes_{\bF[\knotU,\knotV]} C'_K, \d\otimes 1+1\otimes \d', \iota_K\otimes \iota_K'+\iota_K\Phi\otimes \iota_K'\Psi).
\]
Furthermore,  the map
\[
K\mapsto [(\cCFK^-(K), \iota_K)]
\]
determines a homomorphism from the concordance group to $\frI_K$, where $\cCFK^-(K)$ denotes the version of knot Floer homology generated by monomials $\knotU^i \knotV^j \cdot \xs$, where $i,j\ge 0$.

\subsection{The algebraic involutive mapping cone}
\label{sec:canonical-mapping-cone}

In this section, we define the algebraic involutive mapping cone 
\[
\XI^{\alg}_n(\scC)=(\bX_n(\scC),\iota_{\bX})
\]
of an $\iota_K$-complex $\scC$ of L-space type. The $\bF\llsquare U\rrsquare$ complex $\bX_n(\scC)$ coincides with the minus version of the mapping cone complex constructed by Ozsv\'{a}th and Szab\'{o} \cite{OSIntegerSurgeries}. If $(F,g)\colon \scC\to \scC'$ is an enhanced $\iota_K$-homomorphism, then we will also construct an enhanced $\iota$-homomorphism,
\[
\XI_n^{\alg}(F,g)\colon \XI_n^{\alg}(\scC)\to \XI_n^{\alg}(\scC'),
\]
which is well defined up to homotopy of enhanced $\iota$-morphisms. We now summarize the important properties of $\XI_n^{\alg}$:

\begin{prop} \label{prop:canonical-cone}
Let $\scC$ and $\scC'$ be $\iota_K$-complexes of L-space type.
\begin{enumerate}
\item\label{prop:canonical-cone-1}  The pair $\XI^{\alg}_n(\scC)=(\bX(\scC),\iota_{\bX})$ is an $\iota$-complex, which is well defined up to $\iota$-homotopy equivalence.
\item\label{prop:canonical-cone-2}   If $(F,g)$ is an enhanced $\iota_K$-homomorphism from $\scC$ to $\scC'$, then $\XI_n^{\alg}(F,g)$ is an enhanced $\iota$-homomorphism from $\bXI_n^{\alg}(\scC)$ to $\bXI_n^{\alg}(\scC')$. If $(F_1,g_1)$ and $(F_2,g_2)$ are two enhanced $\iota_K$-morphisms which are $\iota_K$-homotopic, then $\XI_n^{\alg}(F_1,g_1)$ and $\XI_n^{\alg}(F_2,g_2)$ are $\iota$-homotopic.
\item\label{prop:canonical-cone-3}  If $(F,g)$ is an enhanced $\iota_K$-homomorphism, such that $F$ is an $\iota_K$-homotopy equivalence, then $\XI_n^{\alg}(F,g)$ is an enhanced $\iota$-homotopy equivalence.
\item\label{prop:canonical-cone-4}  If $(F,g)$ is an enhanced $\iota$-homomorphism, such that $F$ is a local map, then $\XI_n^{\alg}(F,g)$ is also a local map.
\end{enumerate}
\end{prop}

In this section, we describe the constructions of $\XI_n^{\alg}(\scC)=(\bX_n(\scC),\iota_{\bX})$ and $\XI_n^{\alg}(F,g)$, and prove parts part~\eqref{prop:canonical-cone-2}--\eqref{prop:canonical-cone-4} of Proposition~\ref{prop:canonical-cone}. In the subsequent Section~\ref{sec:canonical-cone-is-iota-cx}, we prove part~\eqref{prop:canonical-cone-1}.

Recall that $A_s$ denotes the subset of $(\knotU,\knotV)^{-1}C_K$ generated by monomials $\knotU^{i} \knotV^{j}\cdot \xs$  satisyfing $A(\xs)+j-i=0$, $i\ge 0$ and $j\ge -s$. Recall also that $B_s$ denotes the subset of $(\knotU,\knotV)^{-1}C_K$ generated by monomials $\knotU^i \knotV^j \cdot \xs$ satisfying $A(\xs)+j-i$, and $i\ge 0$, but with no restriction on $j$. We define $\Bopp_s$ to be the subcomplex generated by $\knotU^i \knotV^j\cdot \xs$ satisfying $A(\xs)+j-i=0$, and $j\ge -s$, but with no restriction on $i$.

Note that in the literature, $A_s$ is often given the equivalent definition of being generated by $[\xs,i,j]$ with $A(\xs)-j+i=0$, $i\le 0$ and $j\le s$, and similarly for $B_s$.  The correspondence between the two definitions is given in~\eqref{eq:correspondence-variable-brackets}.
  
We define $\bm{A}_s$ and $\bm{B}_s$ to be the completions over $\bF\llsquare U\rrsquare$ of $A_s$ and $B_s$. We set
\[
\bA:=\prod_{s\in \Z} \bm{A}_s\quad \text{and} \quad \bB:=\prod_{s\in \Z} \bm{B}_s.
\]  
\begin{rem} Here and throughout, we slightly extend the notion of a graded ring to allow for infinite sums of homogeneously graded elements, so that $\mathbb F \llsquare U \rrsquare$ is a graded ring over which $A_s$ and $B_s$ are graded modules, and likewise for analogous cases. \end{rem}

  Next, we let $\frA$ be the $\bF\llsquare \knotU,\knotV\rrsquare$-module obtained by completing $C_K$ with respect to the actions of $\knotU$ and $\knotV$. We let $\frB$ be defined as the $\bF\llsquare \knotU,\knotV\rrsquare$-completion of the module $\knotV^{-1} \cdot C_K$, and we let $\tilde{\frB}$ be the completion of the module $\knotU^{-1}\cdot C_K$. There is a canonical \emph{shift map}
\begin{equation}
\scS\colon \frA\to \bA, \label{eq:shift-map-def}
\end{equation}
which is an $\bF\llsquare U\rrsquare$-equivariant chain isomorphism, as we now describe. The map $\scS$ sends $\cA_s$ to $A_s$, via multiplication by $V^{-s}$. This uniquely determines a map on the completion. There are similar isomorphisms between $\frB$ and $\bB$, and $\tilde{\frB}$ and $\tilde{\bB}$. Abusing notation slightly, we denote these isomorphisms also by $\scS$.

Next, we pick a homogeneously graded chain homotopy equivalence
\[
\frF\colon \Bopp_s\to B_{s+n},
\]
which exists since $\scC$ is of L-space type. We set $\bX_n(\scC)$ to be  the mapping cone
\[
\bX_n(\scC):=\Cone( D_n\colon \bA\to \bB),
\]
where $D_n:=v+\frF \vopp.$ Here, $v$ denotes the inclusion of $A_s$ into $B_s$, and $\vopp$ denotes the inclusion of $A_s$ into $\Bopp_s$. The map $\frF\vopp$ coincides with the map $h$ in Ozsv\'{a}th and Szab\'{o}'s notation \cite{OSIntegerSurgeries}.

 We now construct the involution $\iota_{\bX}$, which will have the form shown in~\eqref{eq:iota-X-def}.
\begin{equation}
\iota_{\bX}=\begin{tikzcd}[column sep=3cm, row sep=1.5cm,labels=description]
\bA
	\arrow[d,"\iota_{\bA}"]
	\arrow[dr,dashed, "H\tilde{v}"]
&
\bB
	\arrow[d, "\iota_{\bB}"]\\
\bA
&\bB.
\end{tikzcd}
\label{eq:iota-X-def}
\end{equation}

The map $\iota_{\bA}$ is induced by $\iota_K$ and the shift isomorphism $\scS$, as we now describe. By definition, $\iota_K$ sends $\cA_s$ to $\cA_{-s}$. Hence, with respect to the canonical shift isomorphisms $A_s\iso \cA_s$ and $A_{-s}\iso \cA_{-s}$, given by multiplication by $\knotV^{-s}$ and $\knotV^{s}$, respectively, the map $\iota_{\bA}$ is defined on $A_s$ to be 
\[ \iota_{\bA} := U^{s} \iota_K. \]
The map $\iota_{\bB}$ is defined on $B_s$ to be
\[ \iota_{\bB} := \frF U^s \iota_K. \]
Note, $\iota_K$ sends $B_s$ to $\Bopp_{0}$, $U^s$ sends $\Bopp_0$ to $\Bopp_{-s}$, and $\frF$ sends $\Bopp_{-s}$ to $B_{-s+n}$.

We now define the map $H$ appearing in~\eqref{eq:iota-X-def}. We will construct it to send $\Bopp_s$ to $B_{-s}$. Consider the map
\[
N:=D_n \iota_{\bA}+\iota_{\bB} D_n,
\]
which sends $A_s$ to $B_{-s}\oplus B_{-s+n}$. The component of $N$ which maps $A_s$ to $B_{-s+n}$ consists of 
\begin{equation}
\Pi_{-s+n}N_s=\frF \vopp U^s\iota_K+\frF U^s \iota_K  v,\label{eq:L-s->-s+n}
\end{equation}
where  $\Pi_{-s+n}$ denotes projection of $\bB$ onto $B_{-s+n}$.
 Since $\iota_K$ is skew-equivariant, we have 
 \begin{equation}
 U^s \iota_K  v=\vopp U^s \iota_K, \label{eq:iota-k-skew-filtered}
 \end{equation}
  as maps from $A_s$ to $\Bopp_{-s}$, so \eqref{eq:L-s->-s+n} vanishes. We now consider the component of $N$ which maps $A_s$ to $B_{-s}$. This
takes the form
\begin{equation}
\Pi_{-s} N_s=vU^s\iota_K+\frF U^s\iota_K\frF\vopp. \label{eq:L-s->-s}
\end{equation}
Using~\eqref{eq:iota-k-skew-filtered}, we see that~\eqref{eq:L-s->-s} is equal to $(U^s\iota_K+\frF U^s\iota_K\frF)\vopp$. The maps $U^s\iota_K$ and $\frF U^s\iota_K\frF$ are both homogeneously graded, $\bF[U]$-equivariant chain homotopy equivalences from $\Bopp_s$ to $B_{-s}$. Since $C_K$ is of L-space type, $U^s\iota_K$ and $\frF U^s\iota_K\frF$ must be chain homotopic over $\bF[U]$, as maps from $\bBopp$ to $\bB$. We define the map $H$ appearing in \eqref{eq:iota-X-def} to be any $+1$ graded chain homotopy between $U^s\iota_K$ and $\frF U^s\iota_K\frF$, which sends $\Bopp_s$ to $B_{-s}$.

The construction above gives a map $\iota_{\bX}$ which is clearly a chain map from $\bX_n(\scC)$ to itself. Hence, we may view $\XI_n^{\alg}(\scC)$ as inducing a chain complex over $\bF\llsquare U\rrsquare [Q]/Q^2$ (though we have not established that $\iota_{\bX}^2\simeq \id$).

Before addressing the well-definedness of $\XI_n^{\alg}(\scC)$, it is helpful to first address functoriality of the construction, from which well-definedness will follow easily. The description of enhanced $\iota$-morphisms in ~\eqref{eq:enhanced-morphism-equivalent} implies that the notion of an enhanced $\iota$-morphism is independent from whether $\iota_{\bX}^2\simeq \id$, and hence we can tackle parts~\eqref{prop:canonical-cone-2}--\eqref{prop:canonical-cone-4} of Proposition~\ref{prop:canonical-cone} before considering part~\eqref{prop:canonical-cone-1}.

 Suppose that $\scC$ and $\scC'$ are $\iota_K$-complexes, and $(F,g)\colon \scC\to \scC'$ is an enhanced $\iota_K$-homomorphism. Let $\frF$ and $H$ be choices of maps used in the construction of $\XI_n^{\alg}(\scC)$, and let $\frF'$ and $H'$ be analogous choices for $\XI_n^{\alg}(\scC')$.

Our goal will be to construct length 2 and 3 maps which turn the diagram in Figure~\ref{fig:cube-for-F-XI} into a chain complex. (Such a diagram is called a \emph{hypercube of chain complexes}; see Section~\ref{sec:hypercubes}). Along the left face of Figure~\ref{fig:cube-for-F-XI}, we use the maps $F$ and $g$ (technically, we are implicitly conjugating with the shift map in~\eqref{eq:shift-map-def}, though we suppress this from the notation).

\begin{figure}[ht!]
\begin{tikzcd}[column sep={3.2cm,between origins},row sep={1.5cm,between origins},labels=description]
\bA
	\ar[dd, swap,"Q(1+\iota_{\bA})"]
	\ar[dr,  "F"]
	\ar[rr, " D_n"]
	\ar[ddrr,dashed,near end,"Q H\vopp"]
	\ar[dddrrr,dotted,sloped,"Q(Kv+J\vopp)"]
&&[-.8cm]
\bB
	\ar[dd, "Q(1+\iota_{\bB})"]
	\ar[dr,"F"]
	\ar[dddr,dashed, "QM"]
&
\\
&\bA'
&&
\bB'
	\ar[dd, "Q(1+\iota_{\bB}')"]
	\ar[from=ulll,dashed,crossing over, "L\vopp"]
	\ar[from=ll,crossing over, "D_n'"]
\\[2cm]
Q\bA
	\ar[rr, "D_n"]
	\ar[dr,"F"]
	\ar[drrr,dashed, "L\vopp"]
&&Q\bB
	\ar[dr, " F"]	
&
\\
&
Q\bA'
	\ar[rr, "D_n'"]
	\ar[from =uu, crossing over,"Q(1+\iota_{\bA}')"]
	\ar[from=uuul,dashed, crossing over, "Qg"]
	&&
Q\bB'
	\ar[from=uull,crossing over, dashed, "Q H'\vopp' "]
\end{tikzcd}
\caption{A diagram encoding $\XI_n^{\alg}(F,g)$.}
\label{fig:cube-for-F-XI}
\end{figure}

We first describe the map $L$, appearing in Figure~\ref{fig:cube-for-F-XI}.  The sum of length 1 compositions along the top face (as well as the bottom face) is $F D_n+D_n' F$, which expands to 
 \begin{equation}
F v+F \frF \vopp+v'F+\frF'\vopp' F. \label{eq:commutators-top+bottom-canonical}
 \end{equation}
Here, $v$ and $\vopp$ denote the inclusion maps for $\scC$, while $v'$ and $\vopp'$ denote the inclusion maps for $\scC'$. Since $F$ is induced by an $\bF[\knotU,\knotV]$-equivariant map, we have $Fv=v'F$  and $\vopp' F=F\vopp$. Hence,~\eqref{eq:commutators-top+bottom-canonical} simplifies to 
\[
(F\frF+\frF'F)\vopp.
\]
As maps from $\bBopp$ to $\bB'$, the maps $F\frF$ and $\frF'F$ have the same grading, and map $\Bopp_s$ to $B'_{s+n}$. Since $\Bopp_s$ and $B'_{s+n}$ are homotopy equivalent to $\bF[U]$, and $\frF$ and $\frF'$ are homotopy equivalences, it follows that the maps $F\frF$ and $\frF' F$ are homotopic to each other. We let $L$ be a $+1$ graded map such that 
\[
[\d, L]=F \frF +\frF'F.
\]

We now construct the map $M$, appearing on the right side of the cube in Figure~\ref{fig:cube-for-F-XI}. We let $M$ be any $\bF[U]$-equivariant, homogeneously $+1$ graded map which sends $B_s$ to $B'_{-s+n}$, and satisfies
\[
[\d,M]=F\iota_{\bB}+\iota_{\bB}'F.
\]
The existence of $M$ follows from the fact that $B_s$ and $B_{s}'$ are homotopy equivalent to $\bF[U]$, using the same logic that gave the existence of $L$.

We now construct the maps $K$ and $J$ appearing along the length 3 arrow (dotted) in Figure~\ref{fig:cube-for-F-XI}. We wish to choose these maps so that the total diagram is a chain complex (or equivalently, so that the maps make the diagram into a hypercube of chain complexes; see Section~\ref{sec:hypercubes}).

Let $N$ denote the composition of the length 1 and 2 maps which feature in the length 3 hypercube relation, i.e.:
 \begin{equation}
 N:=M D_n+\iota_{\bB}' L \vopp+H'\vopp' F+F H\vopp+D_n'g+L \vopp \iota_{\bA}.\label{eq:length-3-relation-hypercube-for-X}
 \end{equation}
First, an easy computation shows that $[\d, N] = 0$ (in fact, this is automatic from the length 1 and 2 relations of the diagram in Figure~\ref{fig:cube-for-F-XI}).

We expand $N$, and use the fact that $F$ is induced by an $\bF[\knotU,\knotV]$-equivariant map, and $\iota_K$ and $g$ are induced by skew-equivariant maps, to obtain that
\[
N=q v+r\vopp,
\]
where
\[
q=M+\frF'g+L\iota_{\bA},\quad\text{and} \quad r=M\frF+\iota_{\bB}'L+H'F+F H+g.
\]
We note that $qv$ maps $A_s$ to $B_{-s+n}'$ while $r\vopp$ maps $A_s$ to $B_{-s}'$. It is easy to check that both $q$ and $r$ are $+1$ graded chain maps. Since $\scC$ and $\scC'$ are of L-space type, every $+1$ graded map from $B_s$ to $B_{-s+n}'$, or from $\Bopp_s$ to $B_{-s}'$ is a boundary, since $\Hom_{\bF[U]}(\bF[U],\bF[U])$ is $\{0\}$ in gradings greater than 0. Hence, $q$ and $r$ are both chain homotopic to zero. We pick $K$ and $J$ to be any $+1$ graded null-homotopies of $q$ and $r$ which send $B_s$ to $B_{-s+n}'$ and $\Bopp_s$ to $B_{-s}'$, respectively.

We now address the well-definedness of $\XI_n^{\alg}(F,g)$ from the choices made in the construction:

\begin{proof}[Proof of Part~\eqref{prop:canonical-cone-2} of Proposition~\ref{prop:canonical-cone}] Suppose that $(F_1,g_1)$ and $(F_2,g_2)$ are two enhanced homomorphisms from $\scC$ to $\scC'$, which are $\iota_K$-homotopic as enhanced morphisms. We wish to show that the maps $\XI_n^{\alg}(F_1,g_1)$ and $\XI_n^{\alg}(F_2,g_2)$ are enhanced $\iota$-chain homotopic. (Our argument when $(F_1,g_1)=(F_2,g_2)$ will prove that $\XI_n^{\alg}(F,g)$ is independent of the choices in the construction, up to enhanced $\iota$-homotopy).

Let $\frF$ and $H$ be the maps used to construct $\XI_n^{\alg}(\scC)$, and let $\frF'$ and $H'$ be the maps used to construct $\XI_n^{\alg}(\scC')$. Let $L_1$, $M_1$, $J_1$ and $K_1$ be the maps used in the construction of $\XI_n^{\alg}(F_1,g_1)$, and let $L_2$, $M_2$, $J_2$ and $K_2$ be the analogs for $\XI_n^{\alg}(F_2,g_2)$.

 We now wish to construct an enhanced $\iota$-morphism $\cE$ from $\XI_n^{\alg}(\scC)$ to $\XI_n^{\alg}(\scC')$ such that
\[
\d_{\ul{\Mor}}(\cE)=\XI^{\alg}_{n}(F_1,g_1)+\XI^{\alg}_n(F_2,g_2).
\]  
  First, we pick $(\phi,\eta)\in \ul{\Mor}(\scC_1,\scC_2)$ such that $\d_{\ul{\Mor}}(\phi,\eta)=(F_1+F_2,g_1+g_2)$, which exists by hypothesis.

\begin{figure}[H]
\begin{tikzcd}[column sep={2cm,between origins},row sep={1cm,between origins},labels=description]
\bA
	\ar[dr,  "F_1+F_2"]
	\ar[dddrrr,dotted,sloped,"Q(K_1+K_2)v+Q(J_1+J_2)\vopp"]
&&[.5cm]
\bB
	\ar[dr,"F_1+F_2"]
	\ar[dddr,dashed, "Q(M_1+M_2)"]
&
\\
&\bA'
&&
\bB'
	\ar[from=ulll,dashed,crossing over, "(L_1+L_2)\vopp "]
\\[2cm]
Q\bA
	\ar[dr,"F_1+F_2"]
	\ar[drrr,dashed, "(L_1+L_2)\vopp"]
&&Q\bB
	\ar[dr, " F_1+F_2"]	
&
\\
&
Q\bA'
	\ar[from=uuul,dashed, crossing over, "Q(g_1+g_2)"]
	&&
Q\bB'
\end{tikzcd}
\begin{tikzcd}[column sep={2cm,between origins},row sep={1cm,between origins},labels=description]
\bA
	\ar[dr,  "\phi"]
	\ar[dddrrr,dotted,"Q\kappa v+Q \zeta \vopp",sloped]
&&[.5cm]
\bB
	\ar[dr,"\phi"]
	\ar[dddr,dashed, "Q \mu"]
&
\\
&\bA'
&&
\bB'
	\ar[from=ulll,dashed,crossing over, "\lambda \vopp"]
\\[2cm]
Q\bA
	\ar[dr,"\phi"]
	\ar[drrr,dashed, "\lambda\vopp"]
&&Q\bB
	\ar[dr, "\phi"]	
&
\\
&
Q\bA'
	\ar[from=uuul,dashed, crossing over, "Q\eta"]
	&&
Q\bB'
\end{tikzcd}
\caption{The enhanced $\iota$-morphisms $\XI_n^{\alg}(F_1,g_1)+\XI_n^{\alg}(F_2,g_2)$ (left)  and $\cE$ (right).}
\label{fig:cube-for-E-XI}
\end{figure}

The relation 
\begin{equation}
\d_{\ul{\Mor}}(\cE)=\XI^{\alg}_{n}(F_1,g_1)+\XI^{\alg}_n(F_2,g_2)\label{eq:dMor(E)=XI+XI}
\end{equation}
can be compactly restated as follows. Consider the $\bF\llsquare U \rrsquare$-module $M$ appearing on the right side of Figure~\ref{fig:cube-for-E-XI} (consisting of the direct sum of $\bA$, $\bB$, $\bA'$, $\bB'$ and so forth). We build an endomorphism $\hat{\cE}$ of $M$, as follows. We begin with $\cE$ as show on the right side of Figure~\ref{fig:cube-for-E-XI}, and add $H\vopp$, $D_n$, $Q(1+\iota_{\bA})$, $Q(1+\iota_{\bB})$ to the back face of the cube, and also add the corresponding maps to the front face. Then, we add the internal differentials of $\bA$, $\bB$, and so forth. Equation~\eqref{eq:dMor(E)=XI+XI} is equivalent to the relation 
\begin{equation}
\hat{\cE}^2=\XI^{\alg}_{n}(F_1,g_1)+\XI^{\alg}_n(F_2,g_2).\label{eq:hat-E-squared}
\end{equation}

We now construct $\lambda$, $\mu$, $\kappa$ and $\zeta$. Our construction will be similar to our construction of the maps $\XI_n^{\alg}(F_i,g_i)$. 

We begin by noting that~\eqref{eq:hat-E-squared} is already satisfied along the left face of the cube, since the desired relations are
\[
F_1+F_2=[\d, \phi]\quad \text{and} \quad g_1+g_2+\iota_{\bA'} \phi+\phi \iota_{\bA}=[\d, \eta],
\]
which is equivalent to $\d_{\ul{\Mor}}(\phi,\eta)=(F_1,g_1)+(F_2,g_2)$.

 Next, we wish $\lambda$ to satisfy
\begin{equation}
D_n'\phi+\phi D_n+(L_1+L_2)\vopp=[\d, \lambda\vopp].
\label{eq:construct-lambda}
\end{equation}
Rearranging, and using the fact that $\phi$ is induced by an $\bF[\knotU,\knotV]$-equivariant map, we see ~\eqref{eq:construct-lambda} is equivalent to
\[
(\phi \frF+\frF' \phi+L_1+L_2)\vopp=[\d,\lambda]\vopp.
\]
The map $\phi\frF+\frF'\phi+L_1+L_2$ is a $+1$ graded chain map from $\bBopp$ to $\bB'$, which sends $\Bopp_s$ to $B_{s+n}$, and hence is null-homotopic. We let $\lambda$ be any $+2$ graded homotopy, which sends $\Bopp_s$ to $B_{s+n}$.

Next, we construct $\mu$. The desired length 2 relation for the right face of the cube is
\[
\phi\iota_{\bB}+\iota_{\bB'} \phi+M_1+M_2=[\d,\mu].
\]
It is straightforward to check that $\phi\iota_{\bB}+\iota_{\bB'} \phi+M_1+M_2$ is a $+1$ graded chain map which sends $B_s$ to $B_{-s+n}$, and hence is null-homotopic. We let $\mu$ be any $+2$ graded null-homotopy, which also sends $B_s$ to $B_{-s+n}$.

Finally, we construct $\kappa$ and $\zeta$. The desired equality is
\begin{equation}
H'\vopp' \phi+D_n' \eta+\mu D_n+\phi H \vopp+\iota_{\bB'}\lambda\vopp +\lambda\vopp\iota_{\bA}+(K_1+K_2)v+(J_1+J_2)\vopp=[\d, \kappa v+\zeta \vopp].
\label{eq:length-3-map-kappe-zeta}
\end{equation}
It is easy to check that the left hand side of~\eqref{eq:length-3-map-kappe-zeta} decomposes as a sum $qv+r\vopp$, where $q$ sends $B_s$ to $B_{-s+n}'$ and $r$ sends $\Bopp_s$ to $B'_{-s}$. Furthermore both $q$ and $r$ are $+2$ graded chain maps, and hence are null-homotopic, since $\Hom_{\bF[U]}(\bF[U],\bF[U])$ is trivial in grading +2. We let $\kappa$ and $\zeta$ be $+3$ graded null-homotopies.

The above relations show exactly that $\d_{\ul{\Mor}}(\cE)=\XI_n^{\alg}(F_1,g_1)+\XI_n^{\alg}(F_2,g_2)$, so the proof is complete.
\end{proof}

We now consider the behavior of $\XI_n^{\alg}$ with respect to $\iota_K$-homotopy equivalences, and local maps:

\begin{proof}[Proof of parts~\eqref{prop:canonical-cone-3} and~\eqref{prop:canonical-cone-4} of Proposition~\ref{prop:canonical-cone}]
We begin with part~\eqref{prop:canonical-cone-3}. Suppose $(F,g)$ is an enhanced $\iota_K$-homomorphism from $\scC$ to $\scC'$, such that $F$ is a homotopy equivalence of chain complexes over $\bF[\knotU,\knotV]$. Clearly, the map $F$ induces a homotopy equivalence of chain complexes over $\bF\llsquare U\rrsquare$ from $\bA$ to $\bA'$. Similarly, $F$ also induces a homotopy equivalence between $\bB$ and $\bB'$. We pick a map $L\colon \bB\to \bB'$, as in the construction of $\XI_n^{\alg}(F,g)$. The same argument as in Lemma~\ref{lem:homotopy-equivalence-iota-complexes} implies that the map induced by $F$ and $L$ is a homotopy equivalence from $\bX_n(\scC)$ to $\bX_n(\scC')$, for any choice of  $L$. Lemma~\ref{lem:homotopy-equivalence-iota-complexes} now implies that the map from $\bXI_n^{\alg}(\scC)$ to $\bXI_n^{\alg}(\scC')$ is an $\iota$-homotopy equivalence, completing the proof of part~\eqref{prop:canonical-cone-3}.

Next, we consider part~\eqref{prop:canonical-cone-4}. For concreteness, we focus on the case that $n>0$, as the case that $n<0$ is a minor modification. Suppose that $(F,g)$ is an enhanced $\iota_K$-homomorphism, such that $F$ is a local map.  We want to show that the map from $\bX_n(\scC)$ to $\bX_n(\scC')$ induced by $F$ and $L$ becomes an isomorphism after inverting $U$. One should note that there is some subtlety in that localization at $U$ does not naturally commute with infinite direct products, so we need to argue somewhat carefully.  First, we truncate the mapping cone and consider the quotient complex $\bX_n(\scC)\langle b\rangle$ generated by $A_i$ for $-b\le i\le b$ and $B_i$ for $-b+n\le i\le b$, for some large $b$. Analogously to \cite{OSIntegerSurgeries}*{Lemma~4.3}, the projection maps onto the quotient complexes $\bX_n(\scC)\langle b\rangle$
 and $\bX_n(\scC')\langle b\rangle$ are quasi-isomorphisms. Since these are finite direct sums, we may invert $U$ on each summand. Each summand is homotopy equivalent to $\bF\llsquare U, U^{-1}]$. Furthermore, both $v_s$ and $\frF \vopp_s$ become quasi-isomorphisms after inverting $U$. From this description, we see that projection of $U^{-1}\bX_n(\scC)\langle b\rangle$ onto $U^{-1}A_0\oplus \cdots \oplus U^{-1}A_{n-1}$ is a quasi-isomorphism. Furthermore, under the projection map onto these summands, the map $F$ is intertwined with the identity map from $\bF\llsquare U, U^{-1}]$ to itself. In particular, the map induced by $F$ is an isomorphism from $U^{-1} H_*(\bX_n(\scC))$ to $U^{-1} H_*(\bX_n(\scC'))$, completing the proof.
\end{proof}

\subsection{The algebraic cone is an iota-complex}
\label{sec:canonical-cone-is-iota-cx}

In this section, we prove part~\eqref{prop:canonical-cone-1} of Proposition~\ref{prop:canonical-cone}, i.e. that the pair $(\bX_n(\scC), \iota_{\bX})$ is an iota complex whenever $\scC$ is an $\iota_K$-complex of L-space type. We begin by proving that the algebraic analogs of the large surgery complexes are iota-complexes:

\begin{lem}\label{lem:large-surgeries-iota-complex}
 Suppose $\cC_K=(C_K,\d,\iota_K)$ is an $\iota_K$-complex. Then $\cC_K$ is also an $\iota$-complex, if we view $C_K$ instead as a chain complex over $\bF[U]$. In particular $\iota_K^2\simeq \id$, $\bF[U]$-equivariantly, on each $\cA_s$.
\end{lem}
\begin{proof} We will write $\cC=(C,\d,\iota)$, for $\cC_K$, viewed as a complex over $\bF[U]$. We will write down a $+1$ graded endomorphism of $\bF[U]$-modules $H\colon C\to C$ which satisfies 
\[
\iota^2+\id=[\d, H].
\]
The map $H$ we write down will not in general be a map of $\bF[\knotU,\knotV]$-modules.

First, it is sufficient to write down a homotopy between $\Phi\Psi$ and $0$, since by definition, $\iota_K^2\simeq \Phi\Psi+\id$.  There is an $\bF[U]$-equivariant map
\[
\Omega\colon C\to C,
\]
defined similarly to $\Phi$ and $\Psi$, except keeping track of the changes of the $U$ variable. More concretely, if $\xs$ and $\ys$ are basis elements of an $\bF[U]$-basis of $C$, and $\d(\xs)$ contains a summand of $U^k\cdot \ys$, then $\Omega(\xs)$ contains a summand of $k U^{k-1} \cdot \ys$. In terms of our original $\bF[\knotU,\knotV]$-basis of $C_K$, $\Omega$ has the following description. If $\langle \d(\knotU^i \knotV^j \cdot \xs),\ys\rangle=\knotU^{i+m}\knotV^{j+n}$, then
\[
\langle  \Omega(\knotU^i\knotV^j\cdot \xs),\ys\rangle:=\left( \min(i+m, j+n)-\min(i,j)\right) \knotU^{i+m-1} \knotV^{j+n-1}.
\]

There is an $\bF[U]$-equivariant map $H_{\Omega^2}\colon C\to C$, of homogeneous grading $+1$, which satisfies
\begin{equation}
\Omega^2=[\d, H_{\Omega^2}].
\label{eq:Phi-hat-squares-to-zero}
\end{equation}
If  $\langle \d(\xs),\ys\rangle=\knotU^m \knotV^n$, 
\[
\langle H_{\Omega^2}(\knotU^i \knotV^j \cdot \xs),\ys\rangle= \frac{\ell(\ell-1)}{2}\cdot \knotU^{i+m-2}\knotV^{j+n-2},
\]
where  $\ell=\min(i+m,j+n)-\min(i,j)$. Equation~\eqref{eq:Phi-hat-squares-to-zero} follows from the same argument as in \cite[Lemma 14.18]{ZemGraphTQFT}.

We make the following claims:
\begin{equation}
\begin{split}
\Phi+ \knotV \Omega&=[\d, H_{\Phi}]\\
\Psi+ \knotU \Omega&=[\d, H_{\Psi}],
\end{split}
\label{eq:Phi-Psi-hat=Phi}
\end{equation}
where 
\[
\begin{split}
H_\Phi(\knotU^i \knotV^j \cdot \xs)&=\max (0,i-j) \knotU^{i-1} \knotV^j\cdot  \xs, \quad \text{and} 
\\
 H_\Psi(\knotU^i \knotV^j \cdot \xs)&=\max(0,j-i) \knotU^i \knotV^{j-1} \cdot \xs.
 \end{split}
\]
It is easy to check that $H_\Phi$ and $H_{\Psi}$ give well defined endomorphisms of $C_K$
 (in particular, there are no negative exponents when applied to elements $\knotU^i \knotV^j\cdot \xs$ with $i,j\ge 0$).

We now prove the first equation of~\eqref{eq:Phi-Psi-hat=Phi}. Suppose that $\langle \d(\xs), \ys\rangle  =\knotU^m \knotV^n$. We compute directly from the definitions that
\[
\begin{split}
\langle \Phi(\knotU^i \knotV^j\cdot \xs),\ys\rangle&=m\cdot \knotU^{i+m-1} \knotV^{j+n}, \\
\langle \knotV \Omega(\knotU^i \knotV^j \cdot \xs),\ys\rangle &=(\min(i+m,j+n)-\min(i,j))\cdot \knotU^{i+m-1}\knotV^{j+n},\\
\langle \d H_\Phi(\knotU^i \knotV^j\cdot \xs) ,\ys \rangle&=\max(0,i-j)\cdot \knotU^{i+m-1}\knotV^{j+n},\\
\langle H_\Phi \d (\knotU^i \knotV^j\cdot \xs),\ys\rangle& = \max (0,i+m-j-n)\cdot \knotU^{i+m-1}\knotV^{j+n}.
\end{split}
\]
Hence, to prove the first equation of ~\eqref{eq:Phi-Psi-hat=Phi}, it is sufficient to show that if $i,j,m,n$ are integers, then
\begin{equation}
\max(0,i-j)+\min(i,j)+\max(0,i+m-j-n)+ \min(m+i,n+j)+m\equiv 0 \pmod{2}. \label{eq:homotopy-Phi-Phi-hat-key-comp}
\end{equation}
Note that $\max(i+m-j-n,0)=\max(i+m,j+n)-j-n$ and $\max(0,i-j)=\max(i,j)-j$. Furthermore, $\max(a,b)+\min(a,b)=a+b$, for any $a$ and $b$. Hence, we manipulate the left hand side of~\eqref{eq:homotopy-Phi-Phi-hat-key-comp}, as follows:
\[
\begin{split}
&\max(i,j)+\min(i,j)+\max(i+m,j+n)+\min(i+m,j+n)-j+m-n-j\\
&\qquad \qquad =(i+j)+(i+m+j+n)-j+m-n-j\\
& \qquad \qquad \equiv 0\pmod{2},
\end{split}
\]
establishing~\eqref{eq:homotopy-Phi-Phi-hat-key-comp}, and hence the first relation of~\eqref{eq:Phi-Psi-hat=Phi}. The second relation of~\eqref{eq:Phi-Psi-hat=Phi} follows by an easy modification.

Combining \eqref{eq:Phi-hat-squares-to-zero}, \eqref{eq:Phi-Psi-hat=Phi} and the $\bF[\knotU,\knotV]$-equivariance of $\Phi$, we conclude
\begin{equation}
\begin{split}
\Phi\Psi&\simeq \Phi \knotU\Omega\\
&=\knotU\Phi\Omega\\
& \simeq \knotU\knotV \Omega \Omega\\
&\simeq 0.
\end{split}
\end{equation}
\end{proof}

The maps $\iota_{\bA}$ and $\iota_K$ are conjugate under $\scS$. Hence, we obtain the following corollary:

\begin{cor}
 The map $\iota_{\bA}\colon \bA\to \bA$ satisfies $\iota_{\bA}^2\simeq \id$.
\end{cor}

We now investigate the homotopy constructed in Lemma~\ref{lem:large-surgeries-iota-complex} in more detail. We write
\[
H_{\frA}\colon \frA\to \frA
\]
for the null-homotopy of $\Phi\Psi$ constructed in Lemma~\ref{lem:large-surgeries-iota-complex}.
Concretely,
\begin{equation}
H_{\frA}:=\Phi H_\Psi+\knotU H_\Phi \Omega+\knotU\knotV H_{\Omega^2}.\label{eq:H-A-def}
\end{equation}

It is helpful now to abuse notation slightly, and write $v$ for the inclusion of $\frA$ into $\frB$, and $\vopp$ for the inclusion of $\frA$ into $\tilde{\frB}$.

In our proof of part~\eqref{prop:canonical-cone-1} of Proposition~\ref{prop:canonical-cone}, we will encounter the expressions $v H_{\frA}$ and $\vopp H_{\frA}$. The general algebraic strategy is to factor the maps through an initial factor of $v$ or $\vopp$. Hence, one might optimistically hope that
\begin{equation}
v H_{\frA}=H_{\frB} v,  \quad \text{and} \quad  \vopp H_{\frA}=H_{\tilde{\frB}} \vopp,\label{eq:unsatisfiable-equality}
\end{equation}
 for some $\bF[U]$-equivariant maps $H_{\frB}\colon \frB\to \frB$ and $H_{\tilde{\frB}}\colon \tilde{\frB}\to \tilde{\frB}$. Such maps do not exist in general, as we illustrate in Example~\ref{ex:fig-8-factorization}.

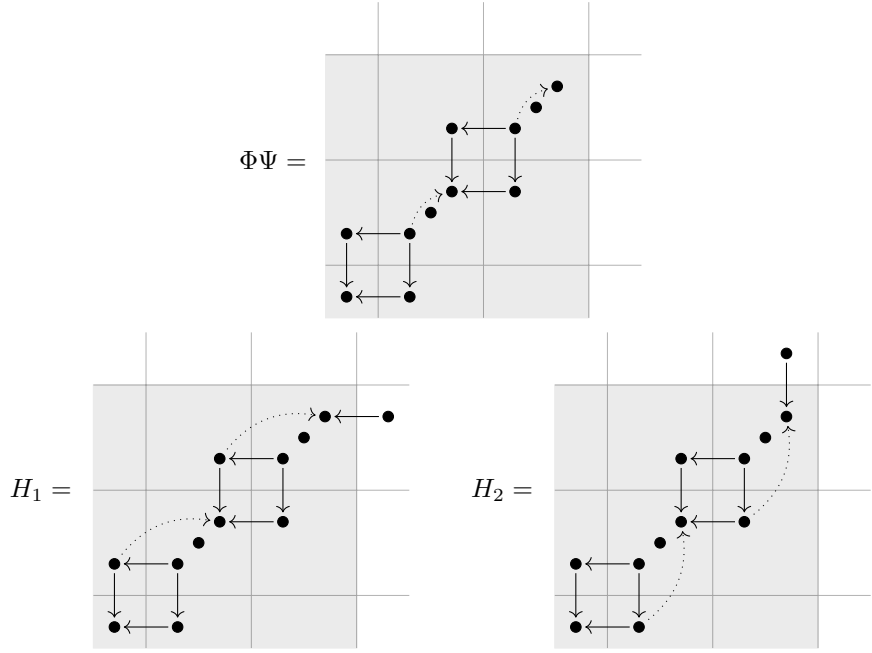
\begin{figure}[ht!]
\[
\begin{tikzpicture}[scale=1.4]
	\draw[step=1, black!30!white, very thin] (-0.5, -0.5) grid (2.5, 2.5);
	\draw[fill] (-.3,-.3) circle [radius=.05] node[auto] (x_0) {};
	\draw[fill] (-.3,.3) circle [radius=.05] node[auto] (x_1) {};
	\draw[fill] (.3,-.3) circle [radius=.05] node[auto] (x_2) {};
	\draw[fill] (.3,.3) circle [radius=.05] node[auto] (x_3) {};
	\draw[fill] (.5,.5) circle [radius=.05] node[auto] (x_4) {};
	\draw[fill] (.7,.7) circle [radius=.05] node[auto] (x_5) {};
	\draw[fill] (.7,1.3) circle [radius=.05] node[auto] (x_6) {};
	\draw[fill] (1.3, 1.3) circle [radius=.05] node[auto] (x_7) {};
	\draw[fill] (1.3, .7) circle [radius=.05] node[auto] (x_8) {};
	\draw[fill] (1.5,1.5) circle [radius=.05] node[auto] (x_9) {};
	\draw[fill] (1.7,1.7) circle [radius=.05] node[auto] (x_10) {};
	\draw (-1,1) node[auto] {$\Phi\Psi=$};
	\draw[->] (x_3) to (x_2);
	\draw[->] (x_3) to (x_1);
	\draw[->] (x_2) to (x_0);
	\draw[->] (x_1) to (x_0);
	\draw[->] (x_7) to (x_6);
	\draw[->] (x_7) to (x_8);
	\draw[->] (x_6) to (x_5);
	\draw[->] (x_8) to (x_5);
	\draw[->,dotted, bend left] (x_7) to (x_10);
	\draw[->, dotted, bend left] (x_3) to (x_5);
	\draw[fill, fill opacity=.08,draw opacity=.08] (-.5,-.5) rectangle (2,2);
\end{tikzpicture}
\]
\[
\begin{tikzpicture}[scale=1.4]
	\draw[step=1, black!30!white, very thin] (-0.5, -0.5) grid (2.5, 2.5);
	\draw[fill] (-.3,-.3) circle [radius=.05] node[auto] (x_0) {};
	\draw[fill] (-.3,.3) circle [radius=.05] node[auto] (x_1) {};
	\draw[fill] (.3,-.3) circle [radius=.05] node[auto] (x_2) {};
	\draw[fill] (.3,.3) circle [radius=.05] node[auto] (x_3) {};
	\draw[fill] (.5,.5) circle [radius=.05] node[auto] (x_4) {};
	\draw[fill] (.7,.7) circle [radius=.05] node[auto] (x_5) {};
	\draw[fill] (.7,1.3) circle [radius=.05] node[auto] (x_6) {};
	\draw[fill] (1.3, 1.3) circle [radius=.05] node[auto] (x_7) {};
	\draw[fill] (1.3, .7) circle [radius=.05] node[auto] (x_8) {};
	\draw[fill] (1.5,1.5) circle [radius=.05] node[auto] (x_9) {};
	\draw[fill] (1.7,1.7) circle [radius=.05] node[auto] (x_10) {};
	\draw[fill] (2.3,1.7) circle [radius=.05] node[auto] (x_11) {};
	\draw (-1,1) node[auto] {$H_1=$};
	\draw[->] (x_3) to (x_2);
	\draw[->] (x_3) to (x_1);
	\draw[->] (x_2) to (x_0);
	\draw[->] (x_1) to (x_0);
	\draw[->] (x_7) to (x_6);
	\draw[->] (x_7) to (x_8);
	\draw[->] (x_6) to (x_5);
	\draw[->] (x_8) to (x_5);
	\draw[->] (x_11) to (x_10);
	\draw[->, dotted, bend left] (x_6) to (x_10);
	\draw[->, dotted, bend left] (x_1) to (x_5);
	\draw[fill, fill opacity=.08,draw opacity=.08] (-.5,-.5) rectangle (2,2);
\end{tikzpicture}
\qquad
\begin{tikzpicture}[scale=1.4]
	\draw[step=1, black!30!white, very thin] (-0.5, -0.5) grid (2.5, 2.5);
	\draw[fill] (-.3,-.3) circle [radius=.05] node[auto] (x_0) {};
	\draw[fill] (-.3,.3) circle [radius=.05] node[auto] (x_1) {};
	\draw[fill] (.3,-.3) circle [radius=.05] node[auto] (x_2) {};
	\draw[fill] (.3,.3) circle [radius=.05] node[auto] (x_3) {};
	\draw[fill] (.5,.5) circle [radius=.05] node[auto] (x_4) {};
	\draw[fill] (.7,.7) circle [radius=.05] node[auto] (x_5) {};
	\draw[fill] (.7,1.3) circle [radius=.05] node[auto] (x_6) {};
	\draw[fill] (1.3, 1.3) circle [radius=.05] node[auto] (x_7) {};
	\draw[fill] (1.3, .7) circle [radius=.05] node[auto] (x_8) {};
	\draw[fill] (1.5,1.5) circle [radius=.05] node[auto] (x_9) {};
	\draw[fill] (1.7,1.7) circle [radius=.05] node[auto] (x_10) {};
	\draw[fill] (1.7, 2.3) circle [radius=.05] node[auto] (x_11) {};
	\draw (-1,1) node[auto] {$H_2=$};
	\draw[->] (x_3) to (x_2);
	\draw[->] (x_3) to (x_1);
	\draw[->] (x_2) to (x_0);
	\draw[->] (x_1) to (x_0);
	\draw[->] (x_7) to (x_6);
	\draw[->] (x_7) to (x_8);
	\draw[->] (x_6) to (x_5);
	\draw[->] (x_8) to (x_5);
	\draw[->] (x_11) to (x_10);
	\draw[->, dotted, bend right] (x_8) to (x_10);
	\draw[->, dotted, bend right] (x_2) to (x_5);
	\draw[fill, fill opacity=.08,draw opacity=.08] (-.5,-.5) rectangle (2,2);
\end{tikzpicture}
\]
\caption{Two null-homotopies of $\Phi\Psi$ for the figure-8. The solid arrows are differentials of $\CFK^\infty$, and the dotted arrows denote $\Phi\Psi$, $H_1$ or $H_2$. The shaded region is $A_0$. }
\label{fig:figure-8-tikz}
\end{figure}

\begin{example}\label{ex:fig-8-factorization}
Consider the figure-8 knot. The complex $A_0$, and the map $\Phi\Psi$ are shown in Figure~\ref{fig:figure-8-tikz}, as well as  two natural choices of null-homotopies, $H_1$ and $H_2$, of $\Phi\Psi$. The map $H_2$ extends to a $\bF[U]$-equivariant map over $B_0$, and the map $H_1$ extends to a map over $\Bopp_0$. There is no $\bF[U]$-equivariant null-homotopy of $\Phi\Psi$ on $A_0$ which extends over both $B_0$ and $\Bopp_0$.
\end{example}

Despite Example~\ref{ex:fig-8-factorization}, it turns out that we may construct maps $H_{\frB}$ and $H_{\tilde{\frB}}$ which satisfy ~\eqref{eq:unsatisfiable-equality} up to a further chain homotopy, as we now describe.

We define $H_{\frB}\colon \frB\to \frB$ via the formula
\begin{equation}
H_{\frB}:= \Phi H_{\Psi}+\knotV^{-1} \knotU  \Phi H_\Phi+\knotV^{-1} \knotU  H_{\Phi^2}.
\label{def:H-B}
\end{equation}
Note that $\knotV$ is invertible on $\frB$. The map $H_{\Phi^2}$ in~\eqref{def:H-B} is defined similarly to $H_{\Omega^2}$, as follows. If $\langle \d (\xs), \ys\rangle=\knotU^m \knotV^n$, then
\[
\langle H_{\Phi^2}(\knotU^i \knotV^j\cdot \xs),\ys\rangle =\frac{m(m-1)}{2}\cdot \knotU^{i+m-2} \knotV^{j+n}.
\]
Thinking 2-categorically, we think of the map $H_{\frB}$ as being constructed from the following sequence of homotopies:
\[
\Phi\Psi\xrightarrow{\Phi H_{\Psi}}   \Phi \knotU \Omega=\knotV^{-1}\knotU  \Phi\knotV \Omega\xrightarrow{\knotV^{-1}\knotU  \Phi H_{\Phi}} \knotV^{-1}  \knotU  \Phi^2\xrightarrow{\knotV^{-1}  \knotU   H_{\Phi^2}}0.
\]
(Here, an arrow denotes a homotopy between the maps at the two ends). It is straightforward to check that the expression in equation~\eqref{def:H-B} gives a well-defined map from $\frB$ to $\frB$.

Next, we define $H_{\tilde{\frB}}\colon \frA\to \tilde{\frB}$ via the formula
\begin{equation}
H_{\tilde{\frB}}:= H_{\Phi} \Psi+\knotU^{-1} \knotV  H_\Psi \Psi+\knotU^{-1} \knotV  H_{\Psi^2}.
\label{eq:H-B-tilde-def}
\end{equation}
Here, $H_{\Psi^2}$ is the map which satisfies
\[
\langle H_{\Psi^2}(\knotU^i \knotV^j \cdot \xs), \ys\rangle =\frac{n(n-1)}{2} \cdot \knotU^{i+m} \knotV^{j+n-2},
\] 
if $\langle \d(\xs),\ys\rangle=\knotU^m \knotV^n$.
We think of $H_{\tilde{\frB}}$ as being associated to the following sequence of homotopies:
\[
\Phi \Psi\xrightarrow{  H_\Phi \Psi}  \knotV \Omega \Psi= \knotU^{-1} \knotV \knotU\Omega \Psi\xrightarrow{\knotU^{-1} \knotV  H_{\Psi} \Psi} \knotU^{-1} \knotV  \Psi \Psi\xrightarrow{\knotU^{-1} \knotV  H_{\Psi^2}} 0. 
\]

The maps $H_{\frB}$ and $H_{\tilde{\frB}}$ satisfy
\[
[\d, H_{\frB}]= \Phi \Psi\quad \text{and} \quad [\d, H_{\tilde{\frB}}]= \Phi \Psi.
\]

The following technical lemma relates $H_{\frA}$, $H_{\frB}$ and $H_{\tilde{\frB}}$:
\begin{lem}\label{lem:devious-homotopies} Suppose that $(C_K,\d,\iota_K)$ is an $\iota_K$-complex of $L$-space type.
\begin{enumerate}
\item\label{lem:homotopies-pt3} There is an $\bF[U]$-equivariant map $h_{\frA,\frB}\colon \frA\to \frB$ satisfying
\[
[\d, h_{\frA,\frB}]=v H_{\frA}+H_{\frB}v.
\]
Furthermore, $h_{\frA,\frB}$ may be taken to have homogeneous grading $+2$, and send $\cA_s$ to $\cB_s$.
\item\label{lem:homotopies-pt4} There is an $\bF[U]$-equivariant map $h_{\frA,\tilde{\frB}}\colon \frA\to \tilde{\frB}$ satisfying
\[
[\d, h_{\frA,\tilde{\frB}}] =\vopp H_{\frA}+H_{\tilde{\frB}}\tilde{v}.
\]
Furthermore, $h_{\frA,\tilde{\frB}}$ may be taken to have homogeneous grading $+2$ and send $\cA_s$ to $\tilde{\cB}_s$.
\end{enumerate}
\end{lem}
\begin{proof} 

We begin by investigating part~\eqref{lem:homotopies-pt3}. We consider the sum $vH_{\frA}+H_{\frB}v$. Using the formulas in~\eqref{eq:H-A-def} and~\eqref{def:H-B}, we compute that if $\langle \d(\xs),\ys \rangle=\knotU^m \knotV^n$, then
\begin{equation}
\begin{split}
&\langle (v H_{\frA}+H_{\frB}v)(\knotU^{i}\knotV^j \cdot \xs), \ys \rangle\\
=&\bigg( \max(0,m+i-n-j)\cdot \left(\min(m+i,n+j)-\min(i,j)\right)\\
+&\frac{(\min(m+i,n+j)-\min(i,j))(\min(m+i,n+j)-\min(i,j)-1)}{2}\\
+&\max(0,i-j)m+\frac{m(m-1)}{2}\bigg) \knotU^{i+m-1} \knotV^{j+n-1}.
\end{split}
\label{eq:horrible-H-AB-1}
\end{equation}
To simplify~\eqref{eq:horrible-H-AB-1}, we introduce some notation. Write 
\begin{equation}
\ve{X}:=\knotU^i \knotV^j\cdot \xs\quad \text{and}\quad  \ve{Y}:=\knotU^{i+m}\knotV^{j+n}\cdot \ys.
\label{eq:boldXboldY-def}
\end{equation}
 Define
\[
\epsilon(\ve{X}):=\max(i-j,0)\quad \text{and} \quad \epsilon(\ve{Y}):=\max(m+i-n-j,0),
\]
which we view as elements of $\Z$.
Note that
\[
\begin{split}
\min(m+i,n+j)-\min(i,j)&=\min(0,n+j-m-i)+m+i-\min(0,j-i)-i\\
&=-\epsilon(\ve{Y})+\epsilon(\ve{X})+m,
\end{split}
\]
since $\min(0,a)=-\max(0,-a)$ for all $a\in \Z$. 
Hence, the coefficient of $\knotU^{i+m-1}\knotV^{j+n-1}$ in \eqref{eq:horrible-H-AB-1} becomes
\[
\begin{split}
&\epsilon(\ve{Y})(-\epsilon(\ve{Y})+\epsilon(\ve{X})+m)+\frac{(-\epsilon(\ve{Y})+\epsilon(\ve{X})+m)(-\epsilon(\ve{Y})+\epsilon(\ve{X})+m-1)}{2}\\
&+m\epsilon(\ve{X})+\frac{m(m-1)}{2}\\
&\equiv \frac{\epsilon(\ve{Y})(\epsilon(\ve{Y})-1)}{2}+\frac{\epsilon(\ve{X})(\epsilon(\ve{X})-1)}{2}\pmod 2.
\end{split}
\]
In particular, $vH_{\frA}+H_{\frB}v=[\d, h_{\frA,\frB}]$, where
\[
h_{\frA,\frB}(\ve{X})=U^{-1}\cdot \frac{\epsilon(\ve{X})(\epsilon(\ve{X})-1)}{2}\cdot \ve{X}.
\]
It is easy to check that $h_{\frA,\frB}$ is $\bF[U]$-equivariant and maps $\cA_s$ to $\cB_s$.

We now consider part~\eqref{lem:homotopies-pt4} of the main claim. The argument is similar to part~\eqref{lem:homotopies-pt3}, as we now describe. Similar to the previous case, using ~\eqref{eq:H-A-def} and~\eqref{eq:H-B-tilde-def}, we compute that if $\langle \d(\xs),\ys\rangle=\knotU^m \knotV^n$, then
\begin{equation}
\begin{split}
&\langle (\vopp H_{\frA}+H_{\tilde{\frB}}\vopp)(\knotU^{i}\knotV^j \cdot \xs), \ys \rangle\\
=&\bigg(m \max(0,j-i)+\max(0,m+i-n-j)(\min(m+i,n+j)-\min(i,j))\\
&+\frac{(\min(m+i,n+j)-\min(i,j))(\min(m+i,n+j)-\min(i,j)-1)}{2}\\
&+ n\max(0,i+m-j-n+1)+n \max(0,j+n-1-m-i)\\
&+\frac{n(n-1)}{2}\bigg)\cdot \knotU^{i+m-1} \knotV^{j+n-1}.
\end{split}
\label{eq:horrible-H-AB-2}
\end{equation}
We let $\ve{X}$ and $\ve{Y}$ be as in ~\eqref{eq:boldXboldY-def}, and define
\[
\nu(\ve{X}):=\max(0,j-i)\quad \text{and} \quad \nu(\ve{Y}):=\max(0,n+j-m-i).
\]
Noting that $\min(m+i,n+j)-\min(i,j)=-\nu(\ve{Y})+\nu(\ve{X})+n$, and also that $\max(0,a)+\max(0,-a)\equiv a\mod 2$, for any $a$, modulo 2 the coefficient of the right hand side of~\eqref{eq:horrible-H-AB-2} becomes
\begin{equation}
\begin{split}
&m\nu(\ve{X})+(m+i-n-j)(-\nu(\ve{Y})+\nu(\ve{X})+n)+\nu(\ve{Y})(-\nu(\ve{Y})+\nu(\ve{X})+n)\\
&+\frac{(-\nu(\ve{Y})+\nu(\ve{X})+n)(-\nu(\ve{Y})+\nu(\ve{X})+n-1)}{2}\\
&+n(i+m-j-n)+n+\frac{n(n-1)}{2}.
\end{split}
\label{eq:nu-homotopy-2}
\end{equation}
Rearranging ~\eqref{eq:nu-homotopy-2}, we obtain
\begin{equation}
\begin{split}
&-\frac{\nu(\ve{Y})(\nu(\ve{Y})-1)}{2}+\frac{\nu(\ve{X})(\nu(\ve{X})-1)}{2}\\
&+(n+m)\nu(\ve{X})+(m+i-n-j)(-\nu(\ve{Y})+\nu(\ve{X}))+n\\
&+n(n-1)+2n(m+i-n-j).
\end{split}
\label{eq:nu-homotopy-3}
\end{equation}
Noting that, modulo 2, $(i-j)\nu(\ve{X})=\nu(\ve{X})$ and $(m+i-n-j)\nu(\ve{Y})=\nu(\ve{Y})$, we see that, modulo 2, \eqref{eq:nu-homotopy-3} is equal to
\begin{equation}
\frac{\nu(\ve{Y})(\nu(\ve{Y})+1)}{2}+\frac{\nu(\ve{X})(\nu(\ve{X})+1)}{2}+n.
\end{equation}
We define 
\[
h(\ve{X}):=U^{-1}\cdot \frac{\nu(\ve{X})(\nu(\ve{X})+1)}{2}\cdot \ve{X}\quad \text{and} \quad k(\knotU^i \knotV^j\cdot \ve{x}):=j \cdot \knotU^{i-1}\knotV^{j-1}\cdot \ve{x}.
\]
Note that $h$ is $\bF[U]$-equivariant, while $k$ is not. Also note that $h$ maps $\cA_s$ to $\tilde{\cB}_s$.

The above manipulation shows that
\begin{equation}
\tilde{v}H_{\frA}+H_{\tilde{\frB}}\tilde{v}=[\d, h]+[\d,k].\label{eq:almost-bbop-homotopy}
\end{equation}
We note that $[\d,k]=\knotU^{-1} \Psi \vopp$.  The map $\knotU^{-1} \Psi$ is $\bF[U]$-equivariant, maps $\tilde{\cB}_s$ to $\tilde{\cB}_s$, and shifts both $\gr_{\ws}$ and $\gr_{\zs}$ by $+1$ (recall that $\gr_{\ws}=\gr_{\zs}+2s$ on $\tilde{\cB}_s$).  Since $C_K$ is of $L$-space type, there is an $\bF[U]$-equivariant map $h_0\colon \tilde{\cB}_s\to \tilde{\cB}_s$ such that  $\knotU^{-1} \Psi=[\d,h_0]$, as endomorphisms of $\tilde{\cB}_s$. We extend $h_0$ to an endomorphism of all of $\tilde{\frB}$, by declaring $h_0$ to be equivariant with respect to $\knotU$ and $\knotV$.  We set $h_{\frA,\tilde{\frB}}$ to be
\[
h_{\frA,\tilde{\frB}}:=h+h_0\vopp,
\]
which manifestly satisfies $[\d, h_{\frA,\tilde{\frB}}]=\vopp H_{\frA}+H_{\tilde{\frB}} \vopp$, completing the proof.
\end{proof}

\begin{prop} 
If $\scC=(C_K,\d,\iota_K)$ is an $\iota_K$-complex of L-space type, then $\bXI^{\alg}_{n}(\scC)$ is an $\iota$-complex.
\end{prop}
\begin{proof}  Let $H$ and $\frF$ denote the maps chosen in the construction of $\iota_{\bX}$. We compute that
\begin{equation}
\iota_{\bX}^2+\id=\begin{tikzcd}[column sep=3cm, row sep=1.5cm,labels=description]
\bA
	\arrow[d,"\iota_{\bA}^2+\id"]
	\arrow[dr,dashed,sloped, "H\iota_{\bA}v+\iota_{\bB}H\tilde{v}"]
&
\bB
	\arrow[d, "\iota_{\bB}^2+\id"]\\
\bA
&\bB.
\end{tikzcd}
\label{eq:iota^2-is-simple}
\end{equation}
Here we have implicitly used the relation $U^s\iota_K v=\vopp U^s\iota_K$, which follows since $\iota_K$ is $\bF[\knotU,\knotV]$-equivariant. Our proof of the main claim will be to show that
\begin{equation}
\iota_{\bX}^2+\id\simeq \begin{tikzcd}[column sep=3cm, row sep=1.5cm,labels=description]
\bA
	\arrow[d,"0"]
	\arrow[dr,dashed,sloped, "qv+r\tilde{v}"]
&
\bB
	\arrow[d, "0"]\\
\bA
&\bB
\end{tikzcd}
\label{eq:iota-complex-goal}
\end{equation}
where $q\colon \bB\to \bB$ and $r\colon \tilde{\bB}\to \bB$ are $\bF\llsquare U\rrsquare$-equivariant maps such that $q$ sends $B_s$ to $B_{s}$ and $r$ sends $\tilde{B}_s$ to $B_{s+n}$, such that $q$ and $r$ are $+1$ graded. Note that if this is the case, then $q$ and $r$ must be chain maps. Furthermore, $q\simeq 0$ and $r\simeq 0$, since $\scC$ is of L-space type. We obtain a homotopy giving the relation $\iota_{\bX}^2+\id\simeq 0$, by defining a map from $\bX_n(\scC)$ to $\bX_n(\scC)$ which puts this null-homotopy of $qv+r\tilde{v}$ along the diagonal.

More generally, suppose that $F\colon \bX_n(\scC)\to \bX_n(\scC)$ is a grading preserving map which is filtered with respect to the mapping cone filtration, i.e., $F$ takes the following form:
  \[
F=\begin{tikzcd}[column sep=3cm, row sep=1.5cm,labels=description]
\bA
	\arrow[d,"J"]
	\arrow[dr,dashed, "K"]
&
\bB
	\arrow[d, "L"]\\
\bA
&\bB.
\end{tikzcd}
\]
We say that $F$  is \emph{simple} if the diagonal map $K$ decomposes as 
\[
K=qv+r\vopp,
\]
where $q\colon \bB\to \bB$ and $r\colon \tilde{\bB}\to \bB$ are homogeneously graded, $\bF\llsquare U\rrsquare$-equivariant maps, such that $q$ sends $B_s$ to $B_s$, and $r$ sends $\Bopp_s$ to $B_{s+n}$. Equation~\eqref{eq:iota^2-is-simple} implies that $\iota_{\bX}^2+\id$ is simple.

Given an $\bF\llsquare U\rrsquare$-equivariant chain map $F\colon \bX_n(\scC)\to \bX_n(\scC)$, which is filtered with respect to the mapping cone filtration, it is natural to change the maps $J$, $K$ and $L$ individually. We make the following observations: 
\begin{enumerate}[label=($\kappa$-\arabic*), ref=$\kappa$-\arabic*,leftmargin=*]
\item\label{item-homotopy-1} If $F'$ is obtained by replacing $L$ with $[\d_{\bB},\kappa ] +L$ and replacing $K$ with $K+\kappa D_n$, for a $+1$ graded map $\kappa \colon \bB\to \bB$, then $F'\simeq F$. If $F$ is simple and $\kappa $ sends $B_s$ to $B_s$, then $F'$ is also simple.
\item\label{item-homotopy-2} If $F'$ is obtained by replacing $J$ with $J+[\d_{\bA},\kappa ]$, where $\kappa \colon \bA\to \bA$ is a $+1$ graded map, and also replacing $K$ with $K+D_n \kappa $, then
\[
F'\simeq F.
\] 
Furthermore, $F'$ is simple if $F$ is simple and $D_n \kappa $ decomposes as $jv+k\tilde{v}$, for homogeneously graded maps $j\colon  \bB\to \bB$ and $k\colon \tilde{\bB}\to \bB$ such that $j$ maps $B_s$ to $B_s$ and $g$ maps $\tilde{B}_s$ to $B_{s+n}$, and $q$ and $r$ are grading preserving (with respect to the mapping cone grading).
\item\label{item-homotopy-3} If $F'$ is obtained from $F$ by adding $[\d, \kappa]$ to the diagonal map $K$, where $\kappa\colon \bA\to \bB$ is any $+1$ graded map of $\bF\llsquare U\rrsquare$-modules, then $F'\simeq F$.
\end{enumerate}
All of the above observations are straightforward to prove.

We return to the expression for $\iota^2_{\bX}+\id$ in equation~\eqref{eq:iota^2-is-simple}. Firstly, $\iota_{\bB}^2\simeq \id$ on $\bB$ (automatically, since $\scC$ is L-space type), so using observation \eqref{item-homotopy-1}, $\iota_{\bX}^2+\id$ is chain homotopic to a simple map which maps $\bA$ to $\bA$ via $\iota_{\bA}^2+\id$, and maps $\bB$ to $\bB$ via $0$.

Next, we take an $\bF[\knotU,\knotV]$-equivariant map $\kappa $ satisfying $[\d,\kappa ]=\iota_K^2+\id+\Phi\Psi$ on $C_K$. Conjugating with the shift isomorphism, we obtain a map from $\bA$ to $\bB$, which we also denote by $\kappa $ (abusing notation slightly). Since $\kappa $ is induced by an $\bF[\knotU,\knotV]$-equivariant map, we have
\[
\begin{split}
D_n \kappa &= (v+\frF \tilde{v})\kappa =\kappa v +\frF \kappa \tilde{v}.
\end{split}
\]
By \eqref{item-homotopy-2}, we may replace the component of $\iota_{\bX}^2+\id$ mapping $\bA$ to $\bA$ (which takes the form $\iota_{\bA}^2+\id$) with $\Phi\Psi$, while retaining simplicity of the map.

By Lemma~\ref{lem:large-surgeries-iota-complex}, there is an $\bF\llsquare U\rrsquare$-equivariant map $H_{\bA}\colon \bA\to \bA$ which satisfies
$\Phi\Psi=[\d_{\bA},H_{\bA}].$
By \eqref{item-homotopy-2}
\[
\iota_{\bX}^2+\id\simeq 
\begin{tikzcd}[column sep=3cm, row sep=1.5cm,labels=description]
\bA
	\arrow[d,"0"]
	\arrow[dr,dashed,sloped, "qv+r\tilde{v}+D_n H_{\bA}"]
&
\bB
	\arrow[d, "0"]\\
\bA
&\bB.
\end{tikzcd}
\]
Note that
\[
D_n H_{\bA}=vH_{\bA}+\frF\vopp H_{\bA}. 
\]
Conjugating the maps constructed in Lemma~\ref{lem:devious-homotopies} with the shift map $\scS$, we obtain maps $H_{\bB}\colon \bB\to \bB$ and $H_{\bBopp}\colon \bBopp\to \bBopp$ such that
\[
vH_{\bA}\simeq H_{\bB}v\quad \text{and} \quad \vopp H_{\bA}\simeq H_{\bBopp} \vopp.
\]
Hence, by \eqref{item-homotopy-3}
\[
\iota_{\bX}^2+\id\simeq 
\begin{tikzcd}[column sep=3cm, row sep=1.5cm,labels=description]
\bA
	\arrow[d,"0"]
	\arrow[dr,dashed,sloped, "(q+H_{\bB})v+(r+H_{\bBopp})\tilde{v}"]
&
\bB
	\arrow[d, "0"]\\
\bA
&\bB,
\end{tikzcd}
\]
which completes the proof, since we have reduced to the form described in \eqref{eq:iota-complex-goal}.
\end{proof}

\subsection{The local equivalence class of surgeries}
\label{sec:local-class-surgeries}

In this section, we prove that the local equivalence class of $\XI_{n}^{\alg}(\scC)$ has a simple form:

\begin{prop}\label{prop:local-equivalence-class} 
Suppose $\scC$ is an $\iota_K$-complex of L-space type and $n>0$.
\begin{enumerate}
\item\label{prop:local-equivalence-class-0} Then $\bXI_{n}^{\alg}(\scC,[0])$ is locally equivalent to $A_0(\scC)$.
\item\label{prop:local-equivalence-class-even} Furthermore, $\XI_{2n}^{\alg}(\scC, [n])$ is locally equivalent to the complex
 \begin{equation}
\begin{tikzcd} [column sep={1cm,between origins},labels=description] 
A_{n}
\ar[dr, "v"]
& & A_n \ar[dl,"v"]\\
& B_n
\end{tikzcd}
\label{eq:local-class-even-2}
\end{equation}
 with the involution which swaps the two copies of $A_n$, and fixes $B_n$.
\end{enumerate}
\end{prop}

We split the proof of Proposition~\ref{prop:local-equivalence-class} into two parts:

\begin{proof}[Proof of part~\eqref{prop:local-equivalence-class-0} of Proposition~\ref{prop:local-equivalence-class}]  We assume $n=1$, to simplify the notation. Projection gives a local map from $\bXI_{+1}^{\alg}(\scC)$ to $A_0(\scC)$, so it suffices to construct a local map
\[
F\colon \bXI_{+1}^{\alg}(\scC)\langle k\rangle\to \bXI_{+1}^{\alg}(\scC)\langle k+1\rangle.
\]
 On the summands $A_{-k+1},\dots, A_{k-1},B_{-k+1},\dots, B_k$, we define $F$ to be the identity. On the remaining summands, we define $F$ as in Figure~\ref{fig:local-map-F0-schematic}, where $\alpha$, $\beta$, and $\gamma$ are maps which we define shortly.
 
 \begin{figure}[H]
 \begin{tikzcd}[column sep={1cm,between origins},labels=description] 
 &&A_{-k}
 	\ar[dr,dashed, "h"]
 	\ar[dd,"\id"]
 	\ar[dddl, bend right=5, "\beta "]
 	\ar[ddll,bend right=10, "\alpha"]
 &
 &\cdots & & 
 A_{k}
 	\ar[dl,dashed,"v" ]
 	\ar[dd,"\id"]
 	\ar[ddrr,bend left=10, "\gamma"]
&\, &
 \\ 
 &&&
B_{-k+1}
 	\ar[dd,"\id"]
&\cdots 
& B_{k}
	\ar[dd,"\id"]
&\, 
\\[1cm]
A_{-k-1}
	\ar[dr,"h",dashed]&
&A_{-k}
  	\ar[dl, "v",dashed]
  	\ar[dr,dashed, "h"]
 &
 &\cdots & & 
A_{k}
	\ar[dl,dashed,"v"]
	\ar[dr,dashed, "h"]&&
    A_{k+1}\ar[dl,dashed, "v"]
 \\ 
& B_{-k}&& B_{-k+1} &\cdots & B_{k}&&B_{k+1} 
 \end{tikzcd}
 \caption{The map $F\colon \bXI_{+1}^{\alg}(\scC)\langle k\rangle\to \bXI_{+1}^{\alg}(\scC)\langle k+1\rangle$. Solid arrows denote $F$, while dashed arrows denote internal differentials of $\XI_{+1}^{\alg}(\scC)\langle k\rangle$ and $\XI_{+1}^{\alg}(\scC)\langle k+1\rangle$.}
 \label{fig:local-map-F0-schematic}
 \end{figure}
 
We recall that $h_s\colon A_s\to B_{s+1}$ is the composition of $\frF_s \vopp_s$, where $\frF_s\colon \tilde{B}_s\to B_{s+1}$ is a choice of homogeneously graded homotopy equivalence, and $\vopp_s$ is the inclusion of $A_s$ into $\tilde{B}_s$. The choice of $\frF_s$ does not affect the homotopy type of the mapping cone complex, so we make the following convenient choice of $\frF_s$:
\begin{equation}
\frF_s:= U^s \iota_K. \label{eq:convenient-F=iota}
\end{equation}
In~\eqref{eq:convenient-F=iota}, $\iota_K$ gives a homotopy equivalence between $\tilde{B}_s=C(j\le s)$ and $C(i\le s)$, and $U^s$ gives a chain isomorphism between $C(i\le s)$ and $B=C(i\le 0)$. With these definitions, the involutions take the form
\[
\iota_{\bB}=\iota_K^2\quad \text{and} \quad \iota_{\bA,s}:=U^s \iota_K.
\]
(The formula for $\iota_{\bB}$ looks odd, but follows from the definition $\iota_{\bB,s}:=\frF_s U^{-s}\iota_K$).

 With this choice of $\frF$, the property of $F$ being a chain map is equivalent to the relations
\begin{equation}
\begin{split}
[\d,\alpha]&=0,\\
h \alpha+v&=[\d, \beta],\\
[\d, \gamma]&=0,\\
 v\gamma+h&=0.
\end{split}
\label{eq:relations-for-F-chain-map}
\end{equation}
In~\eqref{eq:relations-for-F-chain-map}, the maps are to be applied only to $A_{-k}$ or $A_k$, as appropriate, and also $\d$ denotes the internal differential of the summands $A_s$ and $B_s$ (in particular, not the mapping cone differential, which includes  $v$ and $h$).
We set
\[
\alpha=U^{k+1} \iota_K\quad \text{and} \quad \gamma=U^{k}\iota_K.
\]
These are clearly chain maps, so the first and third relations of~\eqref{eq:relations-for-F-chain-map} are satisfied. The fourth is also satisfied, since $h_{k}=U^k \iota_K\vopp_k$.
The second relation of~\eqref{eq:relations-for-F-chain-map}, which involves $\alpha$ and $\beta$, is equivalent to
\[
(\iota_K^2 +\id)v=[\d,\beta].
\]
We pick any $\bF[U]$-equivariant, homogeneously $+1$ graded map $J\colon B_s\to B_s$ satisfying $[\d, J]=\iota_K^2+\id$, which exists because $B_s\simeq \bF[U]$. We define
\begin{equation}
\beta=Jv.\label{eq:beta-def-Jv}
\end{equation}

With these definitions, $F$ becomes a chain map. Furthermore, $F$ is homogeneously graded with respect to the grading on the mapping cone, and clearly becomes an isomorphism on homology after inverting $U$.

It remains to show that 
\begin{equation}
F\iota_{\bX}+\iota_{\bX} F\simeq 0.\label{eq:F-involutive-local}
\end{equation}
Note that $F\iota_{\bX}+\iota_{\bX} F=0$ on all of the summands of $\bXI_{+1}(\scC)\langle k\rangle$ except for $A_{-k}$ and $A_k$. On $A_k$ and $A_{-k}$, we expand out $[F,\iota_{\bX}]$ schematically in Figure~\ref{fig:expanded-computation-[F,iota]}.
 
 \begin{figure}[H]
$[F,\iota_{\bX}]|_{A_{-k}}=$
 \begin{tikzcd}[column sep={1cm,between origins},labels=description] 
 A_{-k}
 \ar[d, "\iota_{\bA}"]
 \ar[drr, bend left=30, "H_{-k}\vopp"]
&&&
\\
A_k
 \ar[d, "\id"]
\ar[dr, bend left=20, "\gamma"]
& &B_k \ar[d, "\id"]
\\
A_k &A_{k+1} &B_k
 \end{tikzcd}
 $ +$ 
 \begin{tikzcd}[column sep={1cm,between origins},labels=description] 
&&& A_{-k}
	\ar[d, "\id"]
	\ar[dll, "\alpha", bend right=30]
	\ar[dl, "\beta", bend right =20]
&\\
& A_{-k-1}
	\ar[dl, "H_{-k-1}\vopp", bend right =25]
	\ar[d, "\iota_{\bA}"]
& B_{-k}
	\ar[d, "\iota_{\bB}"]
& A_{-k}
	\ar[d, "\iota_{\bA}"]
	\ar[dr, bend left=25, "H_{-k}\vopp"]
&\\
B_{k+1}& A_{k+1}& B_{k+1}& A_k& B_k
 \end{tikzcd}
 
\vspace{5mm}
$[F,\iota_{\bX}]|_{A_k}=$
 \begin{tikzcd}[column sep={1cm,between origins},labels=description] 
 && A_k
	\ar[d, "\iota_{\bA}"]
\\
 &&A_{-k}
	\ar[dll, "\alpha", bend right =30]
	\ar[dl,"\beta", bend right=20]
	\ar[d, "\id"] 
 \\
 A_{-k-1}& B_{-k}& A_{-k}
 \end{tikzcd} 
\quad +\quad 
 \begin{tikzcd}[column sep={1cm,between origins},labels=description] 
& A_k
	\ar[d,"\id"]
	\ar[dr,"\gamma", bend left=20]
&&\\
&A_k
	\ar[dl, "H_{k}\vopp", bend right=20]
	\ar[d, "\iota_{\bA}"]
& A_{k+1}
	\ar[d, "\iota_{\bA}"]
\\
B_{-k}& A_{-k}& A_{-k-1}
 \end{tikzcd}  
 \caption{Computing $[F,\iota_{\bX}]$ on $A_{-k}$ (top row) and $A_{k}$ (bottom row).}
 \label{fig:expanded-computation-[F,iota]}
 \end{figure}

 After removing and canceling terms, we obtain the description of $[F,\iota_{\bX}]$ shown in Figure~\ref{fig:commutator-condensed-F,iota}.

 \begin{figure}[H]
 \begin{tikzcd}[column sep={1cm,between origins},labels=description] 
 &&A_{-k}
 	\ar[dr,dashed, "h"]
 	\ar[dddl, bend right=5, "\iota_{\bB}\beta +H_{-k-1} \vopp\alpha"]
 &
 &\cdots & & 
 A_{k}
 	\ar[dl,dashed,"v" ]
 	\ar[dddr, bend left=5, "\beta\iota_{\bA}+H_k\vopp"]
&\, &
 \\ 
 &&&
B_{-k+1}
&\cdots 
& B_{k}
&\, 
\\[1cm]
A_{k+1}
	\ar[dr,"v",dashed]&
&A_{k}
  	\ar[dl, "h",dashed]
  	\ar[dr,dashed, "v"]
 &
 &\cdots & & 
A_{-k}
	\ar[dl,dashed,"h"]
	\ar[dr,dashed, "v"]&&
    A_{-k-1}\ar[dl,dashed, "h"]
 \\ 
& B_{k+1}&& B_{k} &\cdots & B_{-k+1}&&B_{-k} 
 \end{tikzcd}
 \caption{The commutator $[F,\iota_{\bX}]$ (solid arrows). Note that the order of the indices is reversed along the bottom row. Dashed arrows denote internal differentials of $\bX_n(\scC)$.}
 \label{fig:commutator-condensed-F,iota}
 \end{figure}

Hence, it is sufficient to show that the maps 
\[
\iota_{\bB} \beta+H_{-k-1}\vopp \alpha\colon A_{-k}\to B_{k+1}\quad \text{and} \quad \beta \iota_{\bA}+H_k\vopp\colon A_k\to B_{-k}
\]
are null-homotopic. (Given such null-homotopies, we form a null-homotopy of $[F,\iota_{\bX}]$ by putting these null-homotopies in the same positions as in Figure~\ref{fig:local-map-F0-schematic}). Expanding out the definitions, we obtain
\[
\iota_{\bB}\beta+H_{-k-1}\vopp\alpha=(\iota_K^2 J +H_{-k-1} U^{k+1} \iota_K)v.
\]
The map $(\iota_K^2 J+H_{-k-1} U^{k+1} \iota_K)$ is a $+1$ graded chain map from $B_{-k}$ to $B_{k+1}$, and is hence null-homotopic, since both $B_{-k}$ and $B_{k+1}$ are homotopy equivalent to $\bF[U]$, with $1$ supported in the same grading.

Similarly, 
\[
\beta \iota_{\bA}+H_k \vopp= J U^k\iota_K\vopp+H_k \vopp,
\]
which is null-homotopic by the same reasoning. This concludes the claim about the $\Spin^c$ structure $[0]$.
\end{proof}

We now consider the local class of even surgeries:

\begin{proof}[Proof of part~\eqref{prop:local-equivalence-class-even} of Proposition~\ref{prop:local-equivalence-class}]

The proof of part~\eqref{prop:local-equivalence-class-0} carries over to show that $\XI_{2n}^{\alg}(\scC,[n])$ is locally equivalent to the truncation of $\XI_{2n}^{\alg}(\scC,[n])$ shown in ~\eqref{eq:local-class-even-1}, with the involution induced by $\XI^{\alg}_{2n}(\scC)$.
\begin{equation}
\begin{tikzcd} [column sep={1cm,between origins},labels=description] 
A_{-n}
\ar[dr, "h"]
& & A_n \ar[dl,"v"]\\
& B_n
\end{tikzcd}
\label{eq:local-class-even-1}
\end{equation}
As before, we can choose $\frF_{-n}=U^{-n} \iota_K$, so that $h=U^{-n} \iota_K \vopp$, in ~\eqref{eq:local-class-even-1}.

We define chain maps $F$ and $G$ between~\eqref{eq:local-class-even-2} and~\eqref{eq:local-class-even-1} in Figure~\ref{fig:maps-between-truncations}.
\begin{figure}[H]
\[
F=
\begin{tikzcd} [column sep={1.5cm,between origins},labels=description] 
A_{-n}
	\ar[dr, "h",dashed]
	\ar[dd,"U^{-n} \iota_K"]
& & A_n \ar[dl,"v",dashed]
	\ar[dd, "\id"]	
\\
& B_n
	\ar[dd, "\id"]
\\
A_{n}
\ar[dr, "v",dashed]
& & A_n \ar[dl,"v",dashed]\\
& B_n
\end{tikzcd} \qquad G=
\begin{tikzcd} [column sep={1.5cm,between origins},labels=description] 
A_{n}
	\ar[dr, "v",dashed]
	\ar[dd,"U^{n} \iota_K"]
	\ar[dddr, "Jv"]
& & A_n \ar[dl,"v",dashed]
	\ar[dd, "\id"]	
\\
& B_n
	\ar[dd, "\id"]
\\
A_{-n}
\ar[dr, "h",dashed]
& & A_n \ar[dl,"v",dashed]\\
& B_n
\end{tikzcd}
\]
 \caption{The maps $F$ and $G$ (solid arrows).}
 \label{fig:maps-between-truncations}
 \end{figure}
Here, $Jv$ is the same map as in~\eqref{eq:beta-def-Jv}. We compute $[F,\iota]$ and $[G,\iota]$ in Figure~\ref{fig:maps-between-truncations-commutators}.

\begin{figure}[H]
\[
[F,\iota]=
\begin{tikzcd} [column sep={1.5cm,between origins},labels=description] 
A_{-n}
	\ar[dr, "h",dashed]
	\ar[dddr,"H_{-n}\vopp"]
& & A_n \ar[dl,"v",dashed]
	\ar[dd, "\id+\iota_K^2"]	
\\
& B_n
	\ar[dd, "\id+\iota_K^2"]
\\
A_{n}
	\ar[dr, "v",dashed]
& & A_n \ar[dl,"v",dashed]\\
& B_n
\end{tikzcd} \qquad [G,\iota]=
\begin{tikzcd} [column sep={1.5cm,between origins},labels=description] 
A_{n}
	\ar[dr, "v",dashed]
	\ar[dd,"\id+\iota_K^2"]
	\ar[dddr,"Lv"]
& & A_n \ar[dl,"v",dashed]
	\ar[dddl,"Jv"]
\\
& B_n
\ar[dd, "\id+\iota_K^2"]
\\
A_{n}
\ar[dr, "v",dashed]
& & A_{-n} \ar[dl,"h",dashed]\\
& B_n
\end{tikzcd}
\]
 \caption{The commutators $[F,\iota]$ and $[G,\iota]$.}
 \label{fig:maps-between-truncations-commutators}
 \end{figure}
In Figure~\ref{fig:maps-between-truncations-commutators}, 
\[
L=H_{-n} U^n \iota_K +\iota_K^2 J .
\]

We claim that $[F,\iota]$ and $[G,\iota]$ are both null-homotopic. We construct null-homotopies of $[F,\iota]$ and $[G,\iota]$ in ~\eqref{eq:homotopies-F,G-iota}.
\begin{equation}
H_F=
\begin{tikzcd} [column sep={1.5cm,between origins},labels=description] 
A_{-n}
	\ar[dr, "h",dashed]
	\ar[dddr,"\alpha"]
& & A_n \ar[dl,"v",dashed]
	\ar[dd, "H_{\frA}"]	
	\ar[dddl, "\beta"]
\\
& B_n
	\ar[dd, "H_{\frB}"]
\\
A_{n}
	\ar[dr, "v",dashed]
& & A_n \ar[dl,"v",dashed]\\
& B_n
\end{tikzcd}
\qquad H_G=
\begin{tikzcd} [column sep={1.5cm,between origins},labels=description] 
A_{n}
	\ar[dr, "v",dashed]
	\ar[dd,"H_{\frA}"]
	\ar[dddr,"\gamma"]
& & A_n 
	\ar[dl,"v",dashed]
	\ar[dddl, "\delta"]
\\
& B_n
	\ar[dd,"H_{\frB}"]
\\
A_{n}
	\ar[dr, "v",dashed]
& & A_{-n} \ar[dl,"h",dashed]\\
& B_n
\end{tikzcd}
\label{eq:homotopies-F,G-iota}
\end{equation}
In~\eqref{eq:homotopies-F,G-iota} we write $H_{\frB}$ and $H_{\frA}$ for the homotopies constructed in Section~\ref{sec:canonical-cone-is-iota-cx} (plus an $F[\knotU,\knotV]$-equivariant homotopy between $\id+\iota_K^2$ and $\Phi\Psi$, which we suppress from the notation). The relation $[F,\iota]=[\d_{\bX},H_F]$ is equivalent to the relations
\begin{equation}
\begin{split}
\id+\iota_K^2&=[\d,H_{\frA}],\\
\id+\iota_K^2&=[\d, H_{\frB}],\\
H_{\frB}h+H_{-n}\vopp&=[\d,\alpha],\\
H_{\frB} v+v H_{\frA}&=[\d, \beta],
\end{split}
\label{eq:F-iota-homotopy-relations}
\end{equation}
(applied only to $A_n$ or $A_{-n}$, as appropriate). The first and second relations of ~\eqref{eq:F-iota-homotopy-relations} follow from the constructions of $H_{\frA}$ and $H_{\frB}$. The existence of a $\beta$ satisfying the fourth relation of~\eqref{eq:F-iota-homotopy-relations} follows from Lemma~\ref{lem:devious-homotopies}. The existence of $\alpha$ which satisfies the third relation follows since $H_{\frB}h+H_{-k}\vopp$ is a $+1$ graded chain map which factors through $\vopp$.

Similarly, $[G,\iota]=[\d_{\bX},H_G]$ is equivalent to the relations
 \begin{equation}
\begin{split}
\id+\iota_K^2&=[\d, H_{\frA}],\\
\id+\iota_K^2&=[\d, H_{\frB}],\\
H_{\frB}v+Jv&=[\d,\delta],\\
Lv+vH_{\frA}+H_{\frB}v&=[\d, \g].
\end{split}
\label{eq:G-iota-homotopy-relations}
\end{equation}
The first two relations are immediate. The existence of $\delta$ which satisfies the third relation follows since $H_{\frB}v+Jv$ is a $+1$ graded chain map which factors through $v$. The existence of $\g$ satisfying the fourth relation is as follows. By Lemma~\ref{lem:devious-homotopies}, $vH_{\frA}\simeq H_{\frB} v$. Furthermore, $L$ is a $+1$ graded chain map from $B_n$ to $B_n$, so is null-homotopic. The proof is complete.
\end{proof}

For negative surgery coefficients, we have the following analog of Proposition~\ref{prop:local-equivalence-class}:

\begin{prop}\label{prop:negative-surgeries-local-class} 
The local equivalence class of $\XI_{-n}^{\alg}(\scC,[0])$ coincides with the local class of $(A_0(-\scC),\iota_K)^\vee$, where $\vee$ denotes the dual $\iota$-complex, and $-\scC$ denotes the dual $\iota_K$-complex.
\end{prop}
\begin{proof} The proof is very similar to the proof of Proposition~\ref{prop:local-equivalence-class}, so we only sketch the important details. As before, we focus on the case that $n=1$. Inclusion gives a local map from $\XI_{-1}^{\alg}(\scC,[0])\langle k\rangle$ to $\XI_{-1}^{\alg}(\scC,[0])\langle k+1\rangle$, so it suffices to construct a local map $F$ in the opposite direction, which we do in Figure~\ref{fig:negative-surgeries-local-equivalence}.
\begin{figure}[H]
 \begin{tikzcd}[column sep={1cm,between origins},labels=description] 
&A_{-k-1}
	\ar[dr,"v",dashed]
	\ar[dl,"h",dashed]
	\ar[ddrr,bend right=20, "\gamma"]
	\ar[dddr, in=135,out=-90, "\beta"]
&&
A_{-k}
  	\ar[dl, "h",dashed]
  	\ar[dr,dashed, "v"]
  	\ar[dd, "\id"]
 &
 &\cdots & & 
A_{k}
	\ar[dl,dashed,"h"]
	\ar[dr,dashed, "v"]
	\ar[dd,"\id"]
&&
A_{k+1}
	\ar[dl,dashed, "h"]
	\ar[dr,dashed,"v"]
	\ar[ddll,bend left=20, "\delta"]
	\ar[dddlll,in=-15, out=-80, looseness=2, "\epsilon"]
&
 \\ 
B_{-k-2}
	\ar[ddrrrr, out=240,in=200,looseness=2, "\alpha"]
&& 
B_{-k-1}
	\ar[dd,"\id",crossing over]
&& 
B_{-k}
	\ar[dd,"\id"]
 &\cdots & 
B_{k-1}
	\ar[dd, "\id"]
&&
B_{k} 
	\ar[dd, "\id",crossing over]
&& 
B_{k+1}
	\ar[ddllll, out=-60, in=-20, looseness=2, "\zeta"]
\\[1cm]
& &&A_{-k}
 	\ar[dr,dashed, "v"]
 	\ar[dl, dashed, "h"]
 &
 &\cdots & & 
 A_{k}
 	\ar[dl,dashed,"h" ]
 	\ar[dr,dashed, "v"]
&\, 
&
 \\ 
 &&
B_{-k-1}
&&
B_{-k}
&\cdots 
& B_{k-1}
&&
B_k
\end{tikzcd}
\caption{The local equivalence $F\colon \XI_{-1}^{\alg}(\scC)\langle k+1\rangle\to \XI_{-1}^{\alg}(\scC)\langle k\rangle$.}
\label{fig:negative-surgeries-local-equivalence}
\end{figure}

The necessary relations for $F$ to be a chain map are
\begin{equation}
\begin{split}
[\d,\alpha]&=0,\\
[\d, \gamma]&=0,\\
\alpha h+v\gamma&=0,\\
v+h\gamma&=[\d,\beta],\\
[\d,\delta]&=0,\\
[\d,\zeta]&=0,\\
h\delta+\zeta v&=[\d,\epsilon],\\
h+v\delta&=0.
\end{split}
\label{eq:relations-for-negative-local-equivalence}
\end{equation}
In~\eqref{eq:relations-for-negative-local-equivalence}, the third and fourth relations have domain $A_{-k-1}$, while the last two have domain $A_{k+1}$. Toward~\eqref{eq:relations-for-negative-local-equivalence}, we make the following definitions:
\begin{equation}
\alpha:=U^{2k+1} \id,\quad \gamma:=U^k \iota_K, \quad \delta= U^{k+1} \iota_K,\quad \text{and} \quad \zeta=U^{2k+1} \id. \label{eq:defs-for-negative-local-eq}
\end{equation}
It is easily checked that with the definitions in~\eqref{eq:defs-for-negative-local-eq}, maps $\beta$ and $\epsilon$ may be found so that ~\eqref{eq:relations-for-negative-local-equivalence} is satisfied. Thus, $F$ is a homogeneously graded chain map, and clearly sends $\bF[U]$-non-torsion cycles to $\bF[U]$-non-torsion cycles.

The proof of Proposition~\ref{prop:local-equivalence-class} extends easily to show that $[F,\iota_{\bX}]$ is null-homotopic. We conclude that $\XI_{-1}^{\alg}(\scC,[0])$ is locally-equivalent to the complex
\begin{equation}
\begin{tikzcd} [column sep={1cm,between origins},labels=description] 
 & A_0
\ar[dl, "h"]
\ar[dr, "v"]&\\
B_{-1}&& B_0
\end{tikzcd}
\label{eq:negative-surgery-local-preliminary-form}
\end{equation}
with involution induced by $\XI^{\alg}_{-1}$. An argument analogous to the proof of part~\eqref{prop:local-equivalence-class-even} of Proposition~\ref{prop:local-equivalence-class} shows that~\eqref{eq:negative-surgery-local-preliminary-form} is locally equivalent to the complex
\begin{equation}
\begin{tikzcd} [column sep={1cm,between origins},labels=description] 
 & A_0
\ar[dl, "\vopp"]
\ar[dr, "v"]&\\
\Bopp_{0}&& B_0
\end{tikzcd}
\label{eq:negative-surgery-local-less-preliminary-form}
\end{equation}
with the involution which is $\iota_K$ on $A_0$, and is also $\iota_K$ on the $B$-summands, switching $B_0$ and $\Bopp_0$. We now claim that the complex in~\eqref{eq:negative-surgery-local-less-preliminary-form} is homotopy equivalent to the subcomplex $A_0^{\mathrm{big}}$ of $\CFK^\infty(K)$ generated by $[\xs,i,j]$ with $\min(i,j)\le 0$ and $A(\xs)+i-j=0$. To see this, note that the complex in~\eqref{eq:negative-surgery-local-less-preliminary-form} is $\Cone(\vopp\oplus v\colon A_0\to \Bopp_0\oplus B_0)$. Note $A_0=\Bopp_0\cap B_0$ and $A_0^{\mathrm{big}}=\Bopp_0+B_0$. We obtain an exact sequence
\begin{equation}
0\to A_0\xrightarrow{\vopp\oplus v} \Bopp_0\oplus B_0\xrightarrow{\Sigma} A_0^{\mathrm{big}}\to 0, \label{eq:A-0-big-exact-sequence}
\end{equation}
which manifestly induces a chain map $\Phi$ from $\Cone(\vopp\oplus v)\to A_0^{\mathrm{big}}$. A 
chain map $\Psi$ in the opposite direction is given by basic homological algebra.  Namely, we may construct $\bF[U]$-equivariant, splittings $f\colon A_0^{\mathrm{big}}\to \Bopp_0\oplus B_0$ and $g\colon \Bopp_0\oplus B_0\to A_0$ of the exact sequence in~\eqref{eq:A-0-big-exact-sequence} because $A_0^{\mathrm{big}}$ is a free module over $\bF[U]$. In general, $f$ and $g$ will not be chain maps. There is an induced map $\Psi\colon A_0^{\mathrm{big}}\to A_0\oplus \Bopp_0\oplus B_0=\Cone(\vopp\oplus v)$  given by $\Psi(x)=( g \d f x ,f x)$. Basic homological algebra shows that $\Phi$ and $\Psi$ are chain maps, and furthermore, homotopy inverses. The map $\Phi$ commutes with the involution, and hence it follows automatically that $\Psi$ commutes with the involution up to homotopy.

Furthermore, it is not hard to see that $(A_0^{\mathrm{big}}(\scC),\iota_K(\scC))$ is canonically isomorphic to $(A_0(-\scC),\iota_K(-\scC))^\vee$, where $\vee$ denotes the dual $\iota$-complex, completing the proof.
\end{proof}

\begin{rem}
 For $\iota_K$-complexes induced by knots in $S^3$, Proposition~\ref{prop:negative-surgeries-local-class} follows from Proposition~\ref{prop:local-equivalence-class}, together with the observation that $S_{-n}^3(K)=-S_{+n}^3(m(K))$. A similar statement holds more generally for null-homologous knots in L-spaces.
\end{rem}

\subsection{Correction terms of integer surgeries}

Proposition~\ref{prop:local-equivalence-class} has the following implication for the involutive correction terms of surgeries, which is the restriction of Proposition~\ref{prop:correction-terms-intro} to integer surgeries:

\begin{prop}\label{prop:correction-terms-main-computation}
 Suppose $K$ is a knot in $S^3$ and $n>0$. Then
\begin{equation}
\dl(S_n^3(K),[0])=d(L(n,1),[0])-2\Vl_0(K)\quad \text{and} \quad \du(S^3_{n}(K),[0])=d(L(n,1),[0])-2 \Vu_0(K),
\end{equation}
and
\begin{equation}
\dl(S_{2n}^3(K), [n])=d(S_{2n}^3(K),[n]), \quad \text{and} \quad \du(S^3_{2n}(K),[n])=d(L(2n,1),[n]).
\end{equation}
\end{prop}

We begin with a lemma which will allow us to pass to a smaller model of $\bX$, where $B_s=\bF[U]$:

\begin{lem}\label{lem:iotahomotopy}
Let $(C, \iota)$ be an $\iota$-complex. Then for any chain complex $C'$ homotopy equivalent to $C$, there exists a homotopy involution $\iota' \co C' \to C'$ such that the $\iota$-complex $(C', \iota')$ is homotopy equivalent to $(C, \iota)$. 
\end{lem}

\begin{proof}
Suppose that $f \co C \to C'$ and $g \co C' \to C$ such that $g \circ f \simeq \id$ via a homotopy $H$. Then it is straightforward to verify that $\iota' = f \circ \iota \circ g$ gives a homotopy involution on $C'$ and that $f \circ \iota \simeq \iota' \circ f$ via $f \circ \iota \circ H$ and $g \circ \iota' \simeq \iota \circ g$ via $H \circ \iota \circ g$.
\end{proof}

\begin{rem}\label{rem:precise}
Note that if $\bA \simeq \bA'$ via $f_s \co A_s \to A'_s$ and $f'_s \co A'_s \to A_s$ where $f'_s f_s \simeq \id$ via $H_s$, and $\bB \simeq \bB'$ via $g_s \co B_s \to B'_s$ and $g'_s \co B'_s \to B_s$ where $g'_s g_s \simeq \id$ via $J_s$, then $\bX :=\Cone( D \co \bA \to \bB) \simeq \Cone (gDf' \co \bA \to \bB) =: \bX'$ via the homotopy equivalence
\begin{align*}
	F \co \bX &\to \bX'  \\
	(a, b) &\mapsto (f(a), gDH(a) + g(b))
\end{align*} 
with homotopy inverse
\begin{align*}
	G \co \bX' &\to \bX  \\
	(a', b') &\mapsto (f'(a'), JDf'(a') + g'(b')).
\end{align*} 
(It's straightforward to verify that $G \circ F \simeq \id_{\bX}$. Furthermore, one can see that $F \circ G$ is homotopic to a chain isomorphism, i.e., an upper triangular matrix with ones along the diagonal, which we will denote $\phi$. Hence $F \circ G \circ \phi^{-1} \simeq \id_{\bX'}$. Then $G \simeq G \circ F \circ G \circ \phi^{-1} \simeq G \circ \phi^{-1}$.)
In particular, if we replace each $A_s$ and $B_s$ with homotopy equivalent complexes $A'_s$ and $B'_s$ respectively, then the induced mapping cone $\bX'$ still has the form as in Theorem \ref{thm:mapping-cone}.
\end{rem}

We now prove Proposition~\ref{prop:correction-terms-main-computation}, our computation of the correction terms of surgeries:

\begin{proof}[Proof of Proposition~\ref{prop:correction-terms-main-computation}] The first claim about the correction terms in the $\Spin^c$ structure
corresponding to $[0]$ is an immediate consequence of Proposition~\ref{prop:local-equivalence-class}, together with the definition of the absolute grading on the mapping cone complex (see \cite{OSIntegerSurgeries}*{Section~4}), the large surgery formula for involutive Heegaard Floer homology, and the definition of the involutive concordance invariants.

We now investigate the second claim, concerning even surgeries. Let $C$ be the complex
\[
\begin{tikzcd} [column sep={1cm,between origins},labels=description] 
A_{n}
\ar[dr, "v"]
& & A_n \ar[dl,"v"]\\
& B_n
\end{tikzcd}
\]

\noindent with involution $\iota$ which swaps the two copies of $A_n$ and is the identity on $B_n$. This represents the local equivalence class of $\XI_{2n}^{\alg}(\scC, [n]) \simeq \CF^-(S^3_{2n}(Y),[n])$.


For any element $x \in A_n$, let $x^{\ell}$ denote the copy of $x$ in the lefthand copy of $A_n$ and $x^r$ the copy in the righthand copy of $A_n$. By Lemma \ref{lem:iotahomotopy}, we may replace $B_n$ with an appropriately-graded copy of $\bF[U]$. More precisely, using Remark \ref{rem:precise} we choose some chain maps $g\colon B_n \rightarrow \bF[U]$ and $g' \colon \bF[U] \rightarrow B_n$, and a homotopy $J$ between $g'g$ and the identity map on $B_n$. Then taking the maps $f_s$ and $f_s'$ in the statement of the remark to be identity maps and the homotopy $H$ to be identically zero, we see that replacing $B_n$ by $\bF[U]$ replaces $\iota$ by $F \circ \iota \circ G$, where $F$ and $G$ are chain homotopy equivalences in the remark. We see that
\begin{align*}
(F \circ \iota \circ G)(x^{\ell}) &= x^r + Jv(x^r) \\
(F \circ \iota \circ G)(x^r) &= x^{\ell} + Jv(x^{\ell})	\\
(F \circ \iota \circ G)(b) &=b				 
\end{align*}

\noindent where the last term is forced by the fact that $b \in \bF[U]$. Let $M = Jv$, so that $\iota$ is replaced by $\iota+M$. Notice, in particular, that the element $x^{\ell}+x^r$ is always fixed by the replacement involution.

Let $a_n$ be a tower generator of $A_n$. That is, $a_n$ is a cycle of maximal grading such that $U^m[a_n]\neq 0$ for $m\geq 0$ in $H_*(A_n)$. Let $b_n$ be a generator of $B_n$. Then $v(a_n) = U^{V_n}b_n$. The element $[a_n^{\ell}+a_n^r]$ is then a tower generator in $H_*(C)$. By \cite{OSIntegerSurgeries}*{Section~4} and \cite{NiWu}*{Proposition~1.6}
\[
\gr(b_n)=d(L(2n,1),[n])-1 \qquad \qquad \gr(a_n^{\ell})=\gr(a_n^r)=d(S^3_{2n}(K), [n]) = d(L(2n,1),[n])-2V_n(K). \]

We now turn our attention to the computation of the involutive correction terms. Let $\partial_A$ denote the differential in the complex $A_n$ and $\partial_C$ denote the differential in the complex $C$. Consider the involutive complex
\[
CI = (C[-1]\otimes \bF[Q]/(Q^2), \partial^{\iota}=\partial_C + Q(1+\iota)).\]

\noindent Note that, as elements of $CI$, we have
\[
\gr(b_n)=d(L(2n,1),[n]) \qquad \qquad \gr(a_n^{\ell})=\gr(a_n^r)=d(S^3_{2n}(K), [n])+1 \]

\noindent because of the upward grading shift. Observe that $U^m(a_n^{\ell}+a_n^r)$ is $\iota$-invariant and is nonzero in $H_*(CI)$ for all $m\geq 0$. Therefore $[a_n^{\ell}+ a_n^r]$ generates a copy of $\bF[U]$ in $H_*(CI)$ in gradings of the same parity $d(S^3_{2n}(K),[n])+1$ modulo two, which must therefore be contained in the first tower in $H_*(CI)$. Hence $\dl(C)$ must be no less than 
\[
\gr(a_n^{\ell}+a_n^r)-1 = d(S^3_{2n}(K),[n])+1-1 = d(S^3_{2n}(K),[n]).\]

\noindent Since $d(S^3_{2n}(K),[n])=d(C)$ is also an upper bound on $\dl(C)$, we conclude that $[a_n^{\ell}+a_n^r]$ generates the first tower in $H_*(CI)$ and
\[
\dl(C)=d(S^3_{2n}(K),[n]).
\]

We now locate the second tower in $H_*(CI)$. We claim that the class $[b_n]$ generates the tower. First we check that $U^m[b_n]\neq 0$ for $m\geq 0$. For if not, then $U^mb_n$ lies in the image of $\partial^{\iota}$ for some $m$. Since $\partial^{\iota}(B_n[-1]\otimes\bF[Q]/(Q^2))$ lies entirely in $\mathrm{Im}(Q)$, we see that there must exist some $x_1,x_2, x_3, x_4 \in A_n$ such that 

\[U^mb_n = \partial^{\iota}(x_1^{\ell} + (x_2^{\ell}+x_2^{r}) + Qx_3^{\ell}+ Q(x_4^{\ell}+x_4^r)).\]

Then 

\begin{align*}
U^mb_n &= v(x_1^{\ell}) + \partial_A x_1^{\ell} + Q(x_1^{\ell}+x_1^{r}) +QM(x_1^{\ell}) + (\partial_A x_2^{\ell} + \partial_A x_2^r) + Qv(x_3^{\ell}) + Q\partial_Ax_3^{\ell} + Q(\partial_A x_4^{\ell} + \partial_A x_4^r) \\
		&= v(x_1^{\ell})+ \partial_A x_1^{\ell} + (\partial_A x_2^{\ell} + \partial_A x_2^r) + Q(v(x_3^{\ell}) + \partial_A x_3^{\ell}+ (x_1^{\ell}+x_1^{r}) + (\partial_A x_4^{\ell} + \partial_A x_4^r)+Mx_1^{\ell})
\end{align*}

We see that we must have $v(x_1^{\ell})=U^mb_n$ and $\partial_A x_4^{r}= x_1^{r}$. But then

\begin{align*}
\partial_C^2(x_4^{r}) &= \partial_C(x_1^{r}) + \partial_C(v(x_4^{r})) \\
						&= \partial_A(x_1^{r})+v(x_1^{r})+0 \\
						&=\partial_A(x_1^{r}) +U^mb_n \\
						&\neq 0.
\end{align*}

\noindent This is a contradiction, so $U^m[b_n] \neq 0$ for all $m\geq 0$. It remains to be shown that $[b_n] \neq U[\eta]$ for some $\eta \in CI$. Since $b_n$ is an element of maximal grading in $B_n[-1]\otimes \bF[Q]/(Q^2)$, the element $\eta$ must lie in $(A_n\oplus A_n)[-1]\otimes \bF[Q]/(Q^2)$. Suppose that 

\begin{align*}
b_n + U\eta &= \partial^{\iota}(x_1^{\ell} + (x_2^{\ell}+x_2^r) + Qx_3^{\ell} + Q(x_4^{\ell}+x_4^r)) \\
			&= v(x_1^{\ell}) + \partial_A(x_1^{\ell}) + Q(x_1^{\ell}+x_1^r) + QM(x_1^{\ell})+ (\partial_A x_2^{\ell} + \partial_A x_2^r) + Qv(x_3^{\ell}) + Q\partial_A x_3^{\ell} + Q(\partial_A x_4^{\ell} + \partial_A x_4^r).
\end{align*}

\noindent We observe that $v(x_1^{\ell})=b_n$ and $v(x_3^{\ell})=0$. We then have
\[
U\eta = \partial_A(x_1^{\ell}) + (\partial_A x_2^{\ell} + \partial_A x_2^r) + Q\partial_A x_3^{\ell} + Q((x^{\ell}_1+\partial_A x_4^{\ell}) + (x_1^r+ \partial_A x_4^r)) +QM(x_1^{\ell})
\]

Since $U\eta$ is in the image of $U$, in particular $x_1^{r}+\partial_Ax_4^{r} = Ux^{r}$ for some $x \in A_n$. However, then 

\begin{align*}
Uv(x^r)&=v(Ux^r)\\
			&=v(x_1^r+\partial_Ax_4^r)\\
			&=v(x_1^r) + \partial_Bv(x_4^r)\\
			&= b_n+0 \\
			&=b_n.
\end{align*}

This is a contradiction since $b_n$ is not in the image of $U$. We conclude that $[b_n]$ generates the second tower in $H_*(CI)$, and that
\[
\du(C) = \gr(b_n) = (d(L(2n,1),[n])-1)+1 = d(L(2n,1),[n]).
\]
\end{proof}
\section{An application to the homology cobordism group} \label{sec:knotcomputations}
The goal of this section is to use the involutive mapping cone formula to prove Theorem~\ref{thm:nYThetaSF}, that is, to find examples of homology spheres which are not homology cobordant to linear combinations of Seifert fibered spaces.  Our proof begins by computing the $\iota_K$-complexes of certain linear combinations of torus knots and their cables, in order to compute the almost-local equivalence class of $Y=S^3_{+1}(-2T_{6,7} \# T_{6,13} \# -T_{2,3;2,5})$. We then use these calculations to show that $n[Y]$ is non-trivial in $\Theta_{\Z}^3/\Theta_{\mathit{SF}}$ for any non-zero integer $n$. The proof relies on the fact from \cite[Theorem 8.1]{DHSThomcob} (cf. the proof of \cite[Theorem 1.1]{DaiConnected}) that the almost-local equivalence class of linear combinations of Seifert fibered spaces takes on a particularly simple form.

Let $T_{p,q;r,s}$ denote the $(r, s)$-cable of the $(p,q)$-torus knot, where $r$ denotes the longitudinal winding. The goal of this section is to compute the local equivalence class of the $\iota$-complex associated to $S^3_{+1}(K)$, where $K = -2T_{6,7} \# T_{6,13} \# -T_{2,3;2,5}$, in order to prove that $S^3_{+1}(K)$ is not homology cobordant to a linear combination of Seifert fibered spaces.

In Section \ref{sec:tensorprods}, we compute some tensor products of abstract $\iota_K$-complexes. In Section \ref{sec:iotaKexample}, we identify some of these abstract $\iota_K$-complexes with the $\iota_K$-complexes of specific knots in $S^3$. In Section \ref{sec:surgeryalongK}, we use Proposition \ref{prop:local-equivalence-class} to compute the local equivalence class of $(\CF^-(S^3_{+1}(K)), \iota)$. In Section \ref{sec:alec}, we use these calculations to prove Theorem \ref{thm:nYThetaSF}.

\subsection{Some tensor products}\label{sec:tensorprods}
We begin by computing the local equivalence classes of several tensor products of $\iota_K$-complexes. Throughout, we describe $\iota_K$-complexes as finitely generated chain complexes over $\bF[U, U^{-1}]$; that is, following Remark \ref{rem:iotaKU}, we consider $\cA_0$ as a chain complex over $\bF[U]$ with an $\bF[U]$-equivariant endomorphism $\iota_K$ and localize at $U$.

\begin{figure}
\subfigure[]{
\begin{tikzpicture}[scale=0.8]
	\draw[step=1, black!30!white, very thin] (-0.5, -0.5) grid (4.5, 4.5);

	\filldraw (3.7, 3.7) circle (2pt) node[label=right:{\lab{a}}] (a) {};
	\filldraw (3.3, 3.3) circle (2pt) node[label=below left:{\lab{b}}] (b) {};
	\filldraw (0.3, 3.3) circle (2pt) node[label=left:{\lab{c}}] (c) {};
	\filldraw (3.3, 0.3) circle (2pt) node[label=below:{\lab{d}}] (d) {};
	\filldraw (0.3, 0.3) circle (2pt) node[label=below:{\lab{e}}] (e) {};
	
	\draw[->, bend left=0] (b) to node[above,sloped]{\lab{}} (d);	
	\draw[->] (b) to (c);	
	\draw[->] (c) to (e);	
	\draw[->] (d) to (e);		
	
\end{tikzpicture}
\label{fig:C1}
}
\hspace{20pt}
\subfigure[]{
\begin{tikzpicture}[scale=0.8]
	\draw[step=1, black!30!white, very thin] (-0.5, -0.5) grid (5.5, 5.5);

	\filldraw (4.5, 0.5) circle (2pt) node[label=right:{\lab{x_0}}] (x_0) {};
	\filldraw (3.5, 0.5) circle (2pt) node[label=below:{\lab{x_1}}] (x_1) {};
	\filldraw (3.5, 3.5) circle (2pt) node[label=right:{\lab{x_2}}] (x_2) {};
	\filldraw (0.5, 3.5) circle (2pt) node[label=below:{\lab{x_3}}] (x_3) {};
	\filldraw (0.5, 4.5) circle (2pt) node[label=right:{\lab{x_4}}] (x_4) {};
	
	\draw[->](x_0) to (x_1);	
	\draw[->] (x_2) to (x_1);	
	\draw[->] (x_2) to (x_3);	
	\draw[->] (x_4) to (x_3);		
	
\end{tikzpicture}
\label{fig:C2}
}
\caption{Left, the complex $C_1$. Right, the complex $C_2$.}
\label{fig:C1C2}
\end{figure}
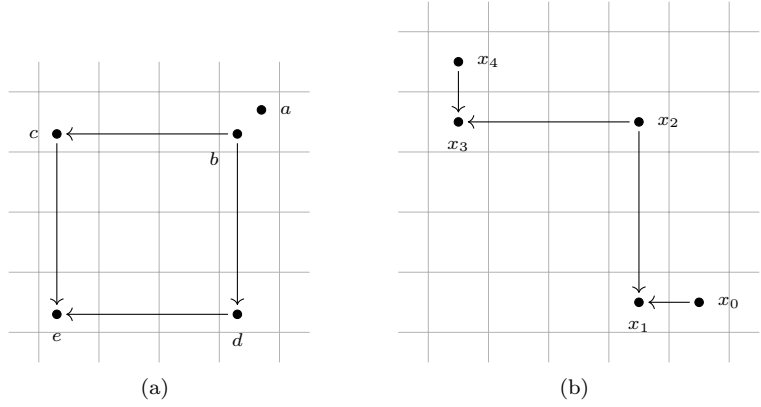

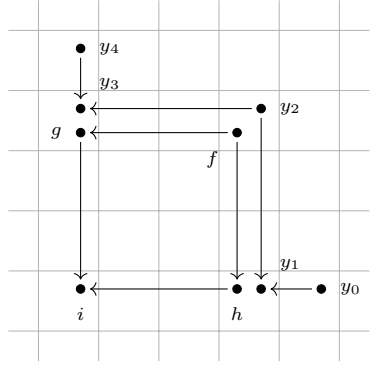
\begin{figure}
\[
\begin{tikzpicture}[scale=0.8]
	\draw[step=1, black!30!white, very thin] (-0.5, -0.5) grid (5.5, 5.5);

	\filldraw (4.7, 0.7) circle (2pt) node[label=right:{\lab{y_0}}] (y_0) {};
	\filldraw (3.7, 0.7) circle (2pt) node[label=above right:{\lab{y_1}}] (y_1) {};
	\filldraw (3.7, 3.7) circle (2pt) node[label=right:{\lab{y_2}}] (y_2) {};
	\filldraw (0.7, 3.7) circle (2pt) node[label=above right:{\lab{y_3}}] (y_3) {};
	\filldraw (0.7, 4.7) circle (2pt) node[label=right:{\lab{y_4}}] (y_4) {};
	\draw[->](y_0) to (y_1);	
	\draw[->] (y_2) to (y_1);	
	\draw[->] (y_2) to (y_3);	
	\draw[->] (y_4) to (y_3);	

	\filldraw (3.3, 3.3) circle (2pt) node[label=below left:{\lab{f}}] (f) {};
	\filldraw (0.7, 3.3) circle (2pt) node[label= left : {\lab{g}}] (g) {};
	\filldraw (3.3, 0.7) circle (2pt) node[label= below : {\lab{h}}] (h) {};
	\filldraw (0.7, 0.7) circle (2pt) node[label= below : {\lab{i}}] (i) {};
	\draw[->] (f) to (h);	
	\draw[->] (f) to (g);	
	\draw[->] (g) to (i);	
	\draw[->] (h) to (i);
	
\end{tikzpicture}
\]
\caption{The complex $C_3$, which is locally equivalent to $C_1 \otimes C_2$.}
\label{fig:C3}
\end{figure}

\begin{lem}\label{lem:C1xC2}
Let $C_1$ be the complex in Figure \ref{fig:C1C2} generated by $a, b, c, d, e$ with Maslov gradings
\[ M(a)=M(b) = 0, \quad  M(c)=M(d)=-1, \quad M(e)=-2,\]  
differential
\[ 
	\d a = 0, \quad
	\d b = c + d, \quad
	\d c = e, \quad
	\d d = e, \quad
	\d e = 0, 
\]
and $\iota_K$ (which we denote by $\iota_1$)
\[ 
	\iota_1(a) = a + U^{-1}e, \quad 
	\iota_1(b) = b+a, \quad 
	\iota_1(c) = d, \quad 
	\iota_1(d) = c, \quad 
	\iota_1(e) = e.
\]
Let $C_2$ be the complex in Figure \ref{fig:C1C2}, together with the endomorphism $\iota_2$
\[ \iota_2(x_i) = x_{4-i}, \quad 0 \leq i \leq 4.
\]
The tensor product $C_1 \otimes C_2$ is locally equivalent to the complex $C_3$ in Figure \ref{fig:C3}, where $\iota_3$ is given by
\begin{align*}
	\iota_3(y_i) &= y_{4-i} , \quad i \neq 2 \\
	\iota_3(y_2) &= y_2 + U^{-1}i \\
	\iota_3(f) &= f+y_2 \\ 
	\iota_3(g) &= h+y_1 \\ 
	\iota_3(h) &= g+y_3 \\ 
	\iota_3(i) &= i.
\end{align*}
\end{lem}

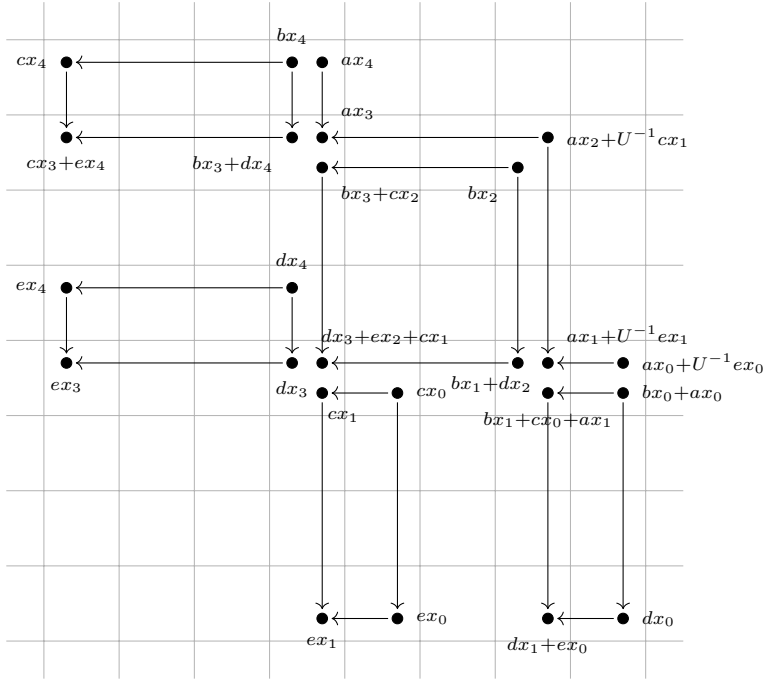
\begin{figure}
\[
\begin{tikzpicture}[scale=1]
	\draw[step=1, black!30!white, very thin] (-3.5, -3.5) grid (5.5, 5.5);

	\filldraw (4.7, 0.7) circle (2pt) node[label=right:{\lab{ax_0+U^{-1}ex_0}}] (y_0) {};
	\filldraw (3.7, 0.7) circle (2pt) node[label=above right:{\lab{ax_1+U^{-1}ex_1}}] (y_1) {};
	\filldraw (3.7, 3.7) circle (2pt) node[label=right:{\lab{ax_2+U^{-1}cx_1}}] (y_2) {};
	\filldraw (0.7, 3.7) circle (2pt) node[label=above right:{\lab{ax_3}}] (y_3) {};
	\filldraw (0.7, 4.7) circle (2pt) node[label=right:{\lab{ax_4}}] (y_4) {};
	\draw[->](y_0) to (y_1);	
	\draw[->] (y_2) to (y_1);	
	\draw[->] (y_2) to (y_3);	
	\draw[->] (y_4) to (y_3);	

	\filldraw (3.3, 3.3) circle (2pt) node[label=below left:{\lab{bx_2}}] (f) {};
	\filldraw (0.7, 3.3) circle (2pt) node[label=below right:{\lab{bx_3+cx_2}}] (g) {};
	\filldraw (3.3, 0.7) circle (2pt) node[label= {[xshift=-10pt, yshift=-18pt] \lab{bx_1+dx_2}}] (h) {};
	\filldraw (0.7, 0.7) circle (2pt) node[label= {[xshift=24pt, yshift=0pt] \lab{dx_3+ex_2+cx_1}}] (i) {};
	\draw[->] (f) to (h);	
	\draw[->] (f) to (g);	
	\draw[->] (g) to (i);	
	\draw[->] (h) to (i);

	\filldraw (4.7, 0.3) circle (2pt) node[label=right:{\lab{bx_0+ax_0}}] (j) {};
	\filldraw (3.7, 0.3) circle (2pt) node[label=below :{\lab{bx_1+cx_0+ax_1}}] (k) {};
	\filldraw (4.7, -2.7) circle (2pt) node[label= right: {\lab{dx_0}}] (l) {};
	\filldraw (3.7, -2.7) circle (2pt) node[label= below : {\lab{dx_1+ex_0}}] (m) {};
	\draw[->] (j) to (l);	
	\draw[->] (j) to (k);	
	\draw[->] (k) to (m);	
	\draw[->] (l) to (m);
	
	\filldraw (0.3, 4.7) circle (2pt) node[label= above:{\lab{bx_4}}] (j') {};
	\filldraw (0.3, 3.7) circle (2pt) node[label= below left :{\lab{bx_3+dx_4}}] (k') {};
	\filldraw (-2.7, 4.7) circle (2pt) node[label= left: {\lab{cx_4}}] (l') {};
	\filldraw (-2.7, 3.7) circle (2pt) node[label= below : {\lab{cx_3+ex_4}}] (m') {};
	\draw[->] (j') to (l');	
	\draw[->] (j') to (k');	
	\draw[->] (k') to (m');	
	\draw[->] (l') to (m');
	
	\filldraw (1.7, 0.3) circle (2pt) node[label=right:{\lab{cx_0}}] (n) {};
	\filldraw (0.7, 0.3) circle (2pt) node[label= {[xshift=8pt, yshift=-18pt] \lab{cx_1}}] (o) {};
	\filldraw (1.7, -2.7) circle (2pt) node[label= right : {\lab{ex_0}}] (p) {};
	\filldraw (0.7, -2.7) circle (2pt) node[label= below : {\lab{ex_1}}] (q) {};
	\draw[->] (n) to (p);	
	\draw[->] (n) to (o);	
	\draw[->] (o) to (q);	
	\draw[->] (p) to (q);
	
	\filldraw (0.3, 1.7) circle (2pt) node[label= above : {\lab{dx_4}}] (n') {};
	\filldraw (0.3, 0.7) circle (2pt) node[label= below : {\lab{dx_3}}] (o') {};
	\filldraw (-2.7, 1.7) circle (2pt) node[label= left : {\lab{ex_4}}] (p') {};
	\filldraw (-2.7, 0.7) circle (2pt) node[label= below : {\lab{ex_3}}] (q') {};
	\draw[->] (n') to (p');	
	\draw[->] (n') to (o');	
	\draw[->] (o') to (q');	
	\draw[->] (p') to (q');

\end{tikzpicture}
\]
\caption{The complex $C_1 \otimes C_2$.}
\label{fig:C1xC2}
\end{figure}

\begin{proof}
Let $(C_1, \iota_1)$ and $(C_2, \iota_2)$ be as above. Recall from \cite[Theorem 1.1]{ZemConnectedSums} that 
\[ \iota_{C_1\otimes C_2} \co C_1\otimes C_2 \to C_1\otimes C_2 \] 
is given by
\[ \iota_1 \otimes \iota_2 + (\Phi_1 \otimes \Psi_2) \circ (\iota_1 \otimes \iota_2). \]
where $\Phi_1$ and $\Psi_2$ are the formal derivatives of $\d$. (In practice, after fixing a basis, $(\Phi_1 \otimes \Psi_2)(x \otimes y)$ counts arrows leaving $x$ that point an odd number of units to the left and arrows leaving $y$ that point an odd number of units down, together with an overall factor of $U^{-1}$.)
For simplicity, we will denote $\iota_{C_1\otimes C_2}$ by $\iota_{12}$. We find that $\iota_{12}$ is given by
\begin{align*}
	\iota_{12} (ax_0) &= ax_4 + U^{-1}ex_4 \\
	\iota_{12} (ax_1) &= ax_3 +U^{-1}ex_3 \\
	\iota_{12} (ax_2) &=  ax_2 + U^{-1}ex_2 \\
	\iota_{12} (ax_3) &=  ax_1 + U^{-1}ex_1 \\
	\iota_{12} (ax_4) &=  ax_0 + U^{-1}ex_0 \\
	\iota_{12} (bx_0) &=  bx_4 + ax_4 + U^{-1} cx_3 \\
	\iota_{12} (bx_1) &=  bx_3 + ax_3 \\
	\iota_{12} (bx_2) &=  bx_2 + ax_2 + U^{-1}cx_1 \\
	\iota_{12} (bx_3) &=  bx_1 + ax_1 \\
	\iota_{12} (bx_4) &=  bx_0 + ax_0 \\
	\iota_{12} (cx_0) &=  dx_4 + U^{-1}ex_3 \\
	\iota_{12} (cx_1) &=  dx_3 \\
	\iota_{12} (cx_2) &=  dx_2 + U^{-1}ex_1 \\
	\iota_{12} (cx_3) &=  dx_1 \\
	\iota_{12} (cx_4) &=  dx_0 \\
	\iota_{12} (dx_0) &=  cx_4 \\
	\iota_{12} (dx_1) &=  cx_3 \\
	\iota_{12} (dx_2) &=  cx_2 \\
	\iota_{12} (dx_3) &=  cx_1 \\
	\iota_{12} (dx_4) &=  cx_0 \\
	\iota_{12} (ex_0) &=  ex_4 \\
	\iota_{12} (ex_1) &=  ex_3 \\
	\iota_{12} (ex_2) &=  ex_2 \\
	\iota_{12} (ex_3) &=  ex_1 \\
	\iota_{12} (ex_4) &=  ex_0.
\end{align*}
Consider the change of basis shown in Figure \ref{fig:C1xC2}. It is straightforward to verify that we may $\iota_{12}$-equivariantly decompose $C_1 \otimes C_2$ as 
\[ C_1 \otimes C_2 = X_1 \oplus X_2 \oplus X_3, \]
where $X_1$ consists of the ``staircase'' together with the three-by-three box (namely, the span of $ax_0+U^{-1}ex_0, ax_1+U^{-1}ex_1, ax_2+U^{-1}cx_1, ax_3, ax_4, bx_2, bx_3+cx_2, bx_1+dx_2, dx_3+ex_2+cx_1$), $X_2$ consists of the upper left and lower right rectangle (the span of $bx_4, cx_4, bx_3+dx_4, cx_3+ex_4, bx_0+ax_0, bx_1+cx_0+ax_1, dx_0, dx_1+ex_0$), and $X_3$ consists of the remaining two rectangles (the span of $dx_4, ex_4, dx_3, ex_3, cx_0, cx_1, ex_0, ex_1$).

Since $X_1$ supports $\HF^\infty$, it follows that $C_1 \otimes C_2$ is locally equivalent to $X_1$, which is isomorphic to $(C_3, \iota_3)$ as described in the statement of the lemma.
\end{proof}

\begin{lem}\label{lem:C1xC4}
Let $C_4$ be the complex in Figure \ref{fig:C4} together with
\[ \iota_4(z_i) = z_{20-i}, \quad 0 \leq i \leq 20. \]
The tensor product of $(C_1, \iota_1)$ and $(C_4, \iota_4)$ is locally equivalent to $C_5$ in Figure \ref{fig:C5}, where $\iota_5$ is given by
\begin{align*}
	\iota_5(w_i) &= w_{20-i}, \quad i \neq 10 \\
	\iota_5(w_{10}) &= w_{10} + U^{-1}m \\
	\iota_5(j) &= j + w_{10} \\
	\iota_5(k) &= l + w_9 \\
	\iota_5(l) &= k + w_{11} \\
	\iota_5(m) &=  m. \\
\end{align*}
\end{lem}

\begin{proof}
The proof is analogous to the proof of Lemma \ref{lem:C1xC2}.
\end{proof}

\begin{figure}
\[
\begin{tikzpicture}[scale=0.5]
	\draw[step=1, black!30!white, very thin] (-0.5, -0.5) grid (31.5, 31.5);

	\filldraw (30.5, 0.5) circle (2pt) node[label=right:{\lab{z_0}}] (z_0) {};
	\filldraw (29.5, 0.5) circle (2pt) node[label=above right:{\lab{z_1}}] (z_1) {};
	\filldraw (29.5, 5.5) circle (2pt) node[label=right:{\lab{z_2}}] (z_2) {};
	\filldraw (28.5, 5.5) circle (2pt) node[label=above right:{\lab{z_3}}] (z_3) {};
	\filldraw (28.5, 10.5) circle (2pt) node[label=right:{\lab{z_4}}] (z_4) {};
	\filldraw (26.5, 10.5) circle (2pt) node[label=above right:{\lab{z_5}}] (z_5) {};
	\filldraw (26.5, 14.5) circle (2pt) node[label=right:{\lab{z_6}}] (z_6) {};
	\filldraw (24.5, 14.5) circle (2pt) node[label=above right:{\lab{z_7}}] (z_7) {};
	\filldraw (24.5, 18.5) circle (2pt) node[label=right:{\lab{z_8}}] (z_8) {};
	\filldraw (21.5, 18.5) circle (2pt) node[label=above right:{\lab{z_9}}] (z_9) {};
	\filldraw (21.5, 21.5) circle (2pt) node[label=right:{\lab{z_{10}}}] (z_10) {};
	\filldraw (18.5, 21.5) circle (2pt) node[label=above right:{\lab{z_{11}}}] (z_11) {};
	\filldraw (18.5, 24.5) circle (2pt) node[label=right:{\lab{z_{12}}}] (z_12) {};
	\filldraw (14.5, 24.5) circle (2pt) node[label=above right:{\lab{z_{13}}}] (z_13) {};
	\filldraw (14.5, 26.5) circle (2pt) node[label=right:{\lab{z_{14}}}] (z_14) {};
	\filldraw (10.5, 26.5) circle (2pt) node[label=above right:{\lab{z_{15}}}] (z_15) {};
	\filldraw (10.5, 28.5) circle (2pt) node[label=right:{\lab{z_{16}}}] (z_16) {};
	\filldraw (5.5, 28.5) circle (2pt) node[label=left:{\lab{z_{17}}}] (z_17) {};
	\filldraw (5.5, 29.5) circle (2pt) node[label=above:{\lab{z_{18}}}] (z_18) {};
	\filldraw (0.5, 29.5) circle (2pt) node[label=left:{\lab{z_{19}}}] (z_19) {};
	\filldraw (0.5, 30.5) circle (2pt) node[label=above :{\lab{z_{20}}}] (z_20) {};
	\draw[->](z_0) to (z_1);	
	\draw[->] (z_2) to (z_1);	
	\draw[->] (z_2) to (z_3);	
	\draw[->] (z_4) to (z_3);	
	\draw[->](z_4) to (z_5);	
	\draw[->] (z_6) to (z_5);	
	\draw[->] (z_6) to (z_7);	
	\draw[->] (z_8) to (z_7);	
	\draw[->] (z_8) to (z_9);	
	\draw[->] (z_10) to (z_9);	
	\draw[->](z_10) to (z_11);	
	\draw[->] (z_12) to (z_11);	
	\draw[->] (z_12) to (z_13);	
	\draw[->] (z_14) to (z_13);	
	\draw[->](z_14) to (z_15);	
	\draw[->] (z_16) to (z_15);	
	\draw[->] (z_16) to (z_17);	
	\draw[->] (z_18) to (z_17);	
	\draw[->] (z_18) to (z_19);	
	\draw[->] (z_20) to (z_19);	
	
\end{tikzpicture}
\]
\caption{The complex $C_4 = \CFKi(-T_{6,13})$.}
\label{fig:C4}
\end{figure}

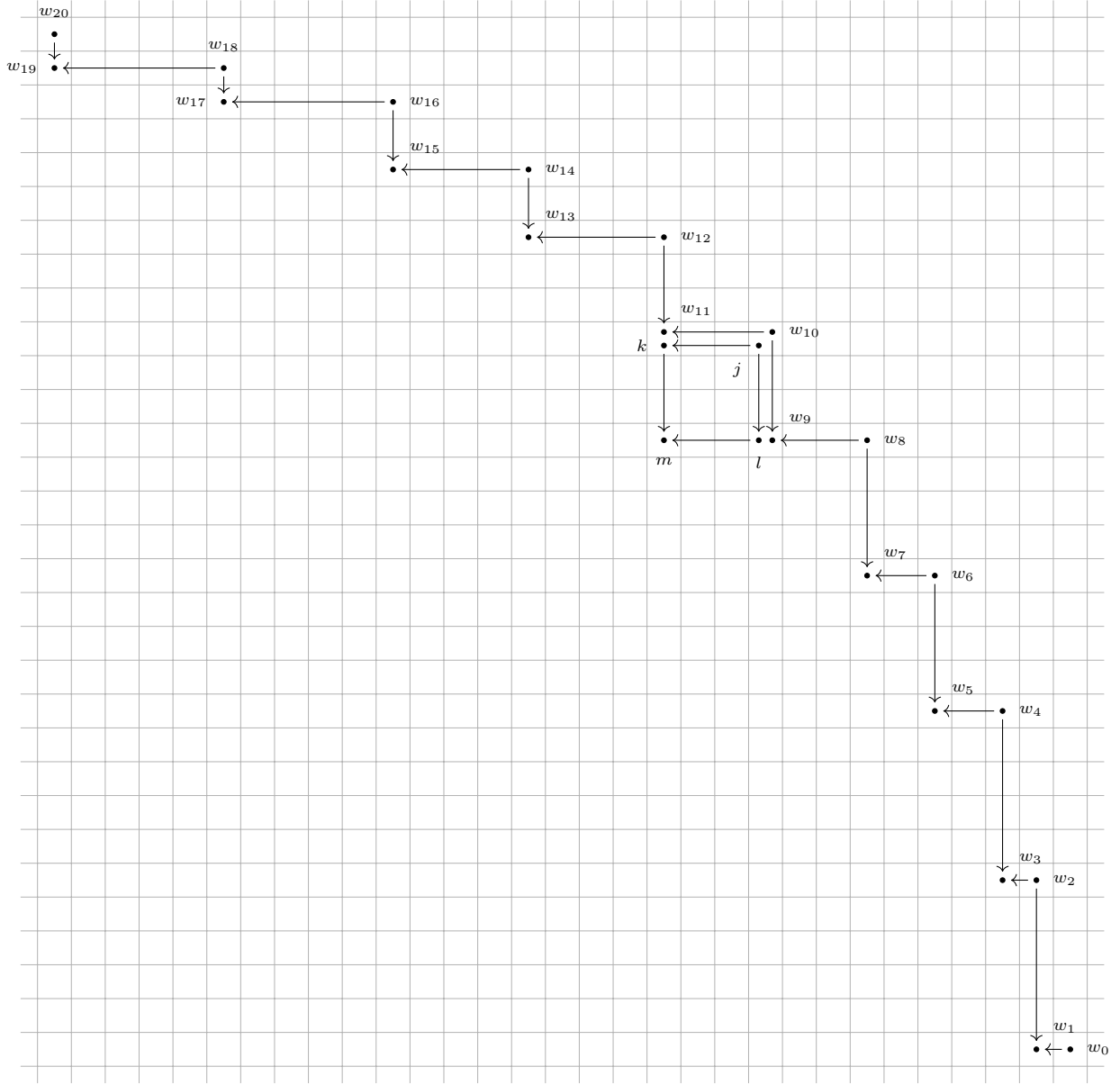
\begin{figure}
\[
\begin{tikzpicture}[scale=0.5]
	\draw[step=1, black!30!white, very thin] (-0.5, -0.5) grid (31.5, 31.5);

	\filldraw (30.5, 0.5) circle (2pt) node[label=right:{\lab{w_0}}] (w_0) {};
	\filldraw (29.5, 0.5) circle (2pt) node[label=above right:{\lab{w_1}}] (w_1) {};
	\filldraw (29.5, 5.5) circle (2pt) node[label=right:{\lab{w_2}}] (w_2) {};
	\filldraw (28.5, 5.5) circle (2pt) node[label=above right:{\lab{w_3}}] (w_3) {};
	\filldraw (28.5, 10.5) circle (2pt) node[label=right:{\lab{w_4}}] (w_4) {};
	\filldraw (26.5, 10.5) circle (2pt) node[label=above right:{\lab{w_5}}] (w_5) {};
	\filldraw (26.5, 14.5) circle (2pt) node[label=right:{\lab{w_6}}] (w_6) {};
	\filldraw (24.5, 14.5) circle (2pt) node[label=above right:{\lab{w_7}}] (w_7) {};
	\filldraw (24.5, 18.5) circle (2pt) node[label=right:{\lab{w_8}}] (w_8) {};
	\filldraw (21.7, 18.5) circle (2pt) node[label=above right:{\lab{w_9}}] (w_9) {};
	\filldraw (21.7, 21.7) circle (2pt) node[label=right:{\lab{w_{10}}}] (w_10) {};
	\filldraw (18.5, 21.7) circle (2pt) node[label=above right:{\lab{w_{11}}}] (w_11) {};
	\filldraw (18.5, 24.5) circle (2pt) node[label=right:{\lab{w_{12}}}] (w_12) {};
	\filldraw (14.5, 24.5) circle (2pt) node[label=above right:{\lab{w_{13}}}] (w_13) {};
	\filldraw (14.5, 26.5) circle (2pt) node[label=right:{\lab{w_{14}}}] (w_14) {};
	\filldraw (10.5, 26.5) circle (2pt) node[label=above right:{\lab{w_{15}}}] (w_15) {};
	\filldraw (10.5, 28.5) circle (2pt) node[label=right:{\lab{w_{16}}}] (w_16) {};
	\filldraw (5.5, 28.5) circle (2pt) node[label=left:{\lab{w_{17}}}] (w_17) {};
	\filldraw (5.5, 29.5) circle (2pt) node[label=above:{\lab{w_{18}}}] (w_18) {};
	\filldraw (0.5, 29.5) circle (2pt) node[label=left:{\lab{w_{19}}}] (w_19) {};
	\filldraw (0.5, 30.5) circle (2pt) node[label=above :{\lab{w_{20}}}] (w_20) {};
	\draw[->](w_0) to (w_1);	
	\draw[->] (w_2) to (w_1);	
	\draw[->] (w_2) to (w_3);	
	\draw[->] (w_4) to (w_3);	
	\draw[->](w_4) to (w_5);	
	\draw[->] (w_6) to (w_5);	
	\draw[->] (w_6) to (w_7);	
	\draw[->] (w_8) to (w_7);	
	\draw[->] (w_8) to (w_9);	
	\draw[->] (w_10) to (w_9);	
	\draw[->](w_10) to (w_11);	
	\draw[->] (w_12) to (w_11);	
	\draw[->] (w_12) to (w_13);	
	\draw[->] (w_14) to (w_13);	
	\draw[->](w_14) to (w_15);	
	\draw[->] (w_16) to (w_15);	
	\draw[->] (w_16) to (w_17);	
	\draw[->] (w_18) to (w_17);	
	\draw[->] (w_18) to (w_19);	
	\draw[->] (w_20) to (w_19);	

	\filldraw (21.3, 21.3) circle (2pt) node[label=below left:{\lab{j}}] (j) {};
	\filldraw (18.5, 21.3) circle (2pt) node[label= left : {\lab{k}}] (k) {};
	\filldraw (21.3, 18.5) circle (2pt) node[label= below : {\lab{l}}] (l) {};
	\filldraw (18.5, 18.5) circle (2pt) node[label= below : {\lab{m}}] (m) {};
	\draw[->] (j) to (k);	
	\draw[->] (j) to (l);	
	\draw[->] (k) to (m);	
	\draw[->] (l) to (m);
	
\end{tikzpicture}
\]
\caption{The complex $C_5$.}
\label{fig:C5}
\end{figure}

\begin{lem}\label{lem:C6xC6}
Let $C_6$ be the complex in Figure \ref{fig:C6} with
\[ \iota_6(v_i) = v_{10-i}, \quad 0 \leq i \leq 10. \]
Then $C_6 \otimes C_6$ is locally equivalent to $(C_5, \iota_5)$, where $\iota_5$ is as described in Lemma \ref{lem:C1xC4}.
\end{lem}

\begin{proof}
After a change of basis, we have that $C_6 \otimes C_6$ splits as in Figure \ref{fig:C6xC6}. Moreover, this splitting may be chosen to be equivariant with respect to $\iota_{C_6 \otimes C_6} = \iota_6 \otimes \iota_6 + (\Phi_6 \otimes \Psi_6) \circ (\iota_6 \otimes \iota_6)$. The ``staircase'' together with the three-by-three box along the diagonal form a direct summand that is isomorphic to $C_5$, under the identification 
\begin{align*}
	w_0 &= v_0v_0	 	\\
	w_1 &= v_0v_1+v_1v_0 	\\							
	w_2 &= v_0v_2+v_2v_0 + U^{-1}v_1v_1					 \\
	w_3 &= v_0v_3+v_1v_2+v_2v_1+v_3v_0 					 \\
	w_4 &= v_0v_4+v_2v_2+v_4v_0 						 \\
	w_5 &= v_0v_5+v_1v_4+v_2v_3+v_3v_2+v_4v_1+v_5v_0 	 \\
	w_6 &= v_0v_6+v_2v_4+v_4v_2+v_6v_0 	+ U^{-1}(v_1v_5+v_5v_1)			 \\
	w_7 &= v_0v_7+v_1v_6+v_2v_5+v_3v_4+v_4v_3+v_5v_2+v_6v_1+v_7v_0 \\
	w_8 &= v_0v_8+v_2v_6+v_4v_4+v_6v_2+v_8v_0 			 \\
	w_9 &= v_0v_9+v_1v_8+v_2v_7+v_3v_6+v_4v_5+v_5v_4+v_6v_3+v_7v_2+v_8v_1+v_9v_0 		 \\
	w_{10} &=  v_0v_{10}+v_2v_8+v_4v_6+v_6v_4+v_8v_2+v_{10}v_0	+ U^{-1}(v_1v_9+v_5v_5)				\\
	w_{11}&= v_1v_{10}+v_2v_9+v_3v_8+v_4v_7+v_5v_6+v_6v_5+v_7v_4+v_8v_3+v_9v_2+v_{10}v_1 \\
	w_{12} &= v_2v_{10}+v_4v_8+v_6v_6+v_8v_4+v_{10}v_2 \\
	w_{13} &= v_3v_{10}+v_4v_9+v_5v_8+v_6v_7+v_7v_6+v_8v_5+v_9v_4+v_{10}v_3 \\
	w_{14} &= v_4v_{10}+v_6v_8+v_8v_6+v_{10}v_4 \\
	w_{15} &= v_5v_{10}+v_6v_9+v_7v_8+v_8v_7+v_9v_6+v_{10}v_5 \\
	w_{16} &= v_6v_{10}+v_8v_8+v_{10}v_6 \\
	w_{17} &= v_{7}v_{10}+v_8v_9+v_9v_8+v_{10}v_7 \\
	w_{18} &= v_{8}v_{10}+v_{10}v_8 \\
	w_{19} &= v_9v_{10}+v_{10}v_9 \\
	w_{20} &= v_{10}v_{10} \\
	j &= v_6v_4+v_8v_2+v_{10}v_0\\
	k &= v_6v_5+v_7v_4+v_8v_3+v_9v_2+v_{10}v_1 \\
	l &= v_5v_4+v_6v_3+v_7v_2+v_8v_1+v_9v_0 \\
	m &= v_5v_5,
\end{align*}
as shown in Figure \ref{fig:C6xC6=C5}. Indeed, it is straightforward to verify that $\iota_{C_6 \otimes C_6}$ maps the span of $\{w_1, \dots, w_{20}, j, k, l, m\}$ into itself. It is then sufficient to show that the basis in Figure \ref{fig:C6xC6} can be chosen so that none of $w_1, \dots, w_{20}, j, k, l, m$ appear in $\iota_{C_6 \otimes C_6}$ of any of the other basis elements. Since $\iota_{C_6 \otimes C_6}$ is grading-preserving and skew-filtered, it is sufficient to verify this for the vertices along the top and right edges of the rectangles shown in Figure \ref{fig:C6xC6-4boxes}, as well as for $v_2v_8$ and $v_8v_2$, the top right corners of the two squares in Figure \ref{fig:C6xC6-4boxes}. By a direct computation, one sees that the span of
\[ \{ v_1v_0, v_1v_1,  v_2v_0, v_2v_1+v_3v_0, v_7v_{10}+v_8v_9, v_8v_{10}, v_9v_9, v_9v_{10} \}, \]
that is, the span of the vertices of the upper left and bottom right rectangles in Figure \ref{fig:C6xC6-4boxes}, forms an $\iota_{C_6 \otimes C_6}$ equivariant subcomplex. Similarly, one sees that the span of
\begin{align*}
 \{ v_3v_2+v_4v_1+v_5v_0, v_3v_3, v_4v_2+v_6v_0+U^{-1}v_5v_1, v_4v_3+v_5v_2v_6v_1+v_7v_0, v_3v_3, \\ v_3v_{10}+v_4v_9+v_5v_8+v_6v_7,  v_4v_{10}+v_6v_8, v_5v_{10}v_6v_9+v_7v_8, v_7v_7, \}
\end{align*}
that is, the span of the vertices of the middle two rectangles in Figure \ref{fig:C6xC6-4boxes}, also forms an $\iota_{C_6 \otimes C_6}$-equivariant subcomplex. Lastly, one sees that $\iota_{C_6 \otimes C_6}(v_2v_8) = v_8v_2+U^{-1}v_9v_1$ and $\iota_{C_6 \otimes C_6}(v_8v_2) = v_2v_8$. (Alternatively, there is a change of basis such that all of the rectangles in Figure \ref{fig:C6xC6-4boxes} except the central $3 \times 3$ square are ``paired'' with their symmetric mate by $\iota_{C_6 \otimes C_6}$.)

Thus, the span of $\{w_1, \dots, w_{20}, j, k, l, m\}$, which is isomorphic to $C_5$ in Figure \ref{fig:C5}, forms an $\iota_{C_6 \otimes C_6}$-equivariant summand of $C_6 \otimes C_6$. Since this summand supports $\HF^\infty$, it follows that $C_6 \otimes C_6$ is locally equivalent to $C_5$, as desired.
\end{proof}

\begin{figure}
\[
\begin{tikzpicture}[scale=0.5]
	\draw[step=1, black!30!white, very thin] (-0.5, -0.5) grid (16.5, 16.5);

	\filldraw (15.5, 0.5) circle (2pt) node[label=right:{\lab{v_0}}] (v_0) {};
	\filldraw (14.5, 0.5) circle (2pt) node[label=above right:{\lab{v_1}}] (v_1) {};
	\filldraw (14.5, 5.5) circle (2pt) node[label=right:{\lab{v_2}}] (v_2) {};
	\filldraw (12.5, 5.5) circle (2pt) node[label=above right:{\lab{v_3}}] (v_3) {};
	\filldraw (12.5, 9.5) circle (2pt) node[label=right:{\lab{v_4}}] (v_4) {};
	\filldraw (9.5, 9.5) circle (2pt) node[label=above right:{\lab{v_5}}] (v_5) {};
	\filldraw (9.5, 12.5) circle (2pt) node[label=right:{\lab{v_6}}] (v_6) {};
	\filldraw (5.5, 12.5) circle (2pt) node[label=above right:{\lab{v_7}}] (v_7) {};
	\filldraw (5.5, 14.5) circle (2pt) node[label=right:{\lab{v_8}}] (v_8) {};
	\filldraw (0.5, 14.5) circle (2pt) node[label=left:{\lab{v_9}}] (v_9) {};
	\filldraw (0.5, 15.5) circle (2pt) node[label=above:{\lab{v_{10}}}] (v_10) {};
	\draw[->](v_0) to (v_1);	
	\draw[->] (v_2) to (v_1);	
	\draw[->] (v_2) to (v_3);	
	\draw[->] (v_4) to (v_3);	
	\draw[->](v_4) to (v_5);	
	\draw[->] (v_6) to (v_5);	
	\draw[->] (v_6) to (v_7);	
	\draw[->] (v_8) to (v_7);	
	\draw[->] (v_8) to (v_9);	
	\draw[->] (v_10) to (v_9);	

\end{tikzpicture}
\]
\caption{The complex $C_6 = \CFKi(-T_{6,7})$.}
\label{fig:C6}
\end{figure}

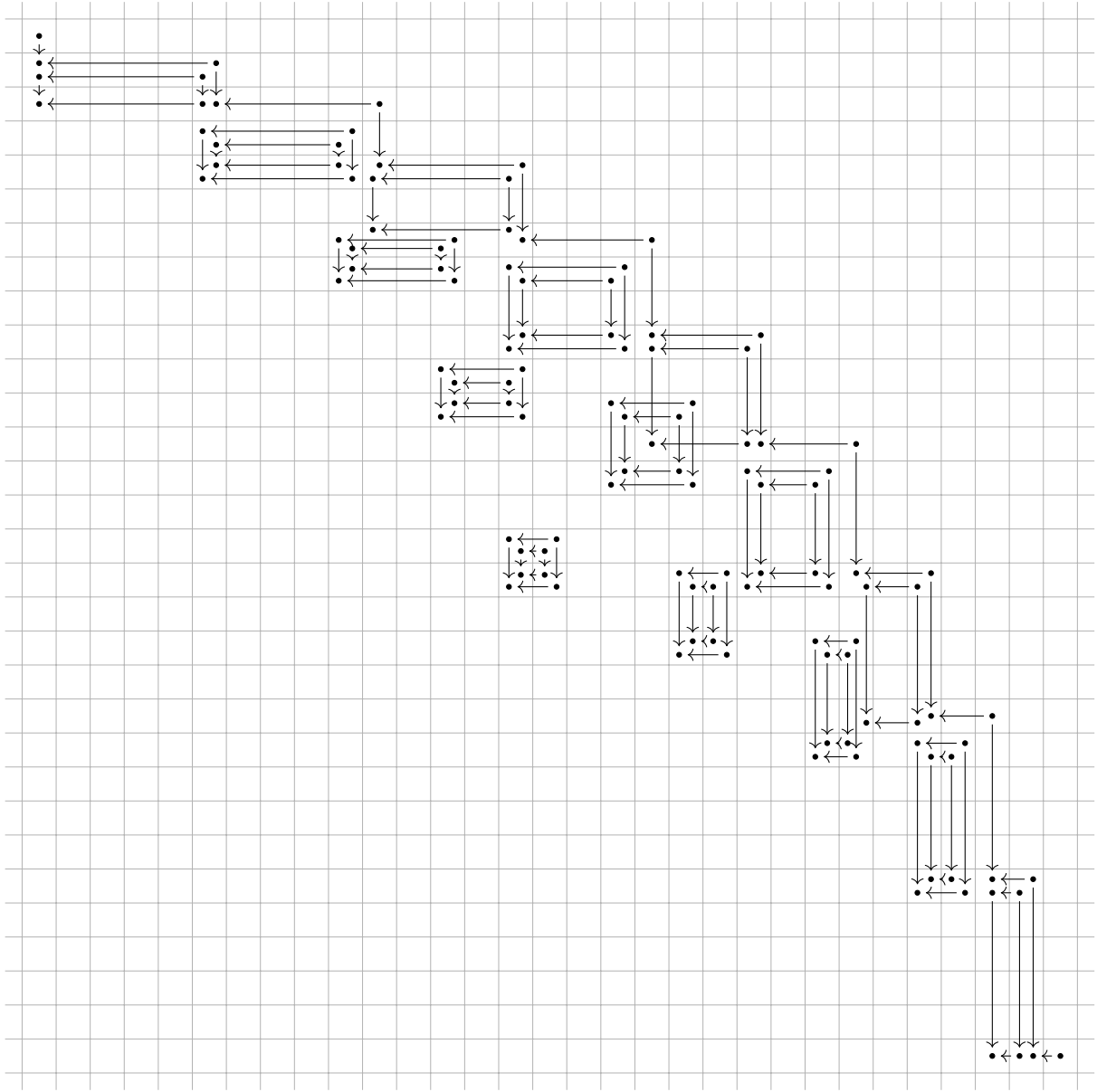
\begin{figure}
\[
\begin{tikzpicture}[scale=0.5]
	\draw[step=1, black!30!white, very thin] (-0.5, -0.5) grid (31.5, 31.5);

	\filldraw (30.5, 0.5) circle (2pt) node[label=right:{\lab{}}] (w_0) {};
	\filldraw (29.7, 0.5) circle (2pt) node[label=above right:{\lab{}}] (w_1) {};
	\filldraw (29.7, 5.7) circle (2pt) node[label=right:{\lab{}}] (w_2) {};
	\filldraw (28.5, 5.7) circle (2pt) node[label=left:{\lab{}}] (w_3) {};
	\filldraw (28.5, 10.5) circle (2pt) node[label=right:{\lab{}}] (w_4) {};
	\filldraw (26.7, 10.5) circle (2pt) node[label=left:{\lab{}}] (w_5) {};
	\filldraw (26.7, 14.7) circle (2pt) node[label=right : {\lab{}}] (w_6) {};
	\filldraw (24.5, 14.7) circle (2pt) node[label=left:{\lab{}}] (w_7) {};
	\filldraw (24.5, 18.5) circle (2pt) node[label=right:{\lab{}}] (w_8) {};
	\filldraw (21.7, 18.5) circle (2pt) node[label= {[xshift=60pt, yshift=0pt] \lab{}}] (w_9) {};
	\filldraw (21.7, 21.7) circle (2pt) node[label=right:{\lab{}}] (w_10) {};
	\filldraw (18.5, 21.7) circle (2pt) node[label= above left:{\lab{}}] (w_11) {};
	\filldraw (18.5, 24.5) circle (2pt) node[label=right:{\lab{}}] (w_12) {};
	\filldraw (14.7, 24.5) circle (2pt) node[label=above right:{\lab{}}] (w_13) {};
	\filldraw (14.7, 26.7) circle (2pt) node[label=right:{\lab{}}] (w_14) {};
	\filldraw (10.5, 26.7) circle (2pt) node[label=above right:{\lab{}}] (w_15) {};
	\filldraw (10.5, 28.5) circle (2pt) node[label=right:{\lab{}}] (w_16) {};
	\filldraw (5.7, 28.5) circle (2pt) node[label=left:{\lab{}}] (w_17) {};
	\filldraw (5.7, 29.7) circle (2pt) node[label=above:{\lab{}}] (w_18) {};
	\filldraw (0.5, 29.7) circle (2pt) node[label=left:{\lab{}}] (w_19) {};
	\filldraw (0.5, 30.5) circle (2pt) node[label=above :{\lab{}}] (w_20) {};
	\draw[->](w_0) to (w_1);	
	\draw[->] (w_2) to (w_1);	
	\draw[->] (w_2) to (w_3);	
	\draw[->] (w_4) to (w_3);	
	\draw[->](w_4) to (w_5);	
	\draw[->] (w_6) to (w_5);	
	\draw[->] (w_6) to (w_7);	
	\draw[->] (w_8) to (w_7);	
	\draw[->] (w_8) to (w_9);	
	\draw[->] (w_10) to (w_9);	
	\draw[->](w_10) to (w_11);	
	\draw[->] (w_12) to (w_11);	
	\draw[->] (w_12) to (w_13);	
	\draw[->] (w_14) to (w_13);	
	\draw[->](w_14) to (w_15);	
	\draw[->] (w_16) to (w_15);	
	\draw[->] (w_16) to (w_17);	
	\draw[->] (w_18) to (w_17);	
	\draw[->] (w_18) to (w_19);	
	\draw[->] (w_20) to (w_19);	

	\filldraw (21.3, 21.3) circle (2pt) node[label=below :{\lab{}}] (j) {};
	\filldraw (18.5, 21.3) circle (2pt) node[label= left : {\lab{}}] (k) {};
	\filldraw (21.3, 18.5) circle (2pt) node[label= below : {\lab{}}] (l) {};
	\filldraw (18.5, 18.5) circle (2pt) node[label= left : {\lab{}}] (m) {};
	\draw[->] (j) to (k);	
	\draw[->] (j) to (l);	
	\draw[->] (k) to (m);	
	\draw[->] (l) to (m);
	
	\filldraw (29.3, 5.3) circle (2pt) node[label=below :{\lab{}}] (n) {};
	\filldraw (29.3, 0.5) circle (2pt) node[label= left : {\lab{}}] (o) {};
	\filldraw (28.5, 5.3) circle (2pt) node[label= below : {\lab{}}] (p) {};
	\filldraw (28.5, 0.5) circle (2pt) node[label= left : {\lab{}}] (q) {};
	\draw[->] (n) to (o);	
	\draw[->] (n) to (p);	
	\draw[->] (o) to (q);	
	\draw[->] (p) to (q);	

	\filldraw (27.7, 9.7) circle (2pt) node[label=below :{\lab{}}] (r) {};
	\filldraw (27.7, 5.3) circle (2pt) node[label= left : {\lab{}}] (s) {};
	\filldraw (26.3, 9.7) circle (2pt) node[label= below : {\lab{}}] (t) {};
	\filldraw (26.3, 5.3) circle (2pt) node[label= left : {\lab{}}] (u) {};
	\draw[->] (r) to (s);	
	\draw[->] (r) to (t);	
	\draw[->] (s) to (u);	
	\draw[->] (t) to (u);

	\filldraw (27.3, 9.3) circle (2pt) node[label=below :{\lab{}}] (r) {};
	\filldraw (27.3, 5.7) circle (2pt) node[label= left : {\lab{}}] (s) {};
	\filldraw (26.7, 9.3) circle (2pt) node[label= below : {\lab{}}] (t) {};
	\filldraw (26.7, 5.7) circle (2pt) node[label= left : {\lab{}}] (u) {};
	\draw[->] (r) to (s);	
	\draw[->] (r) to (t);	
	\draw[->] (s) to (u);	
	\draw[->] (t) to (u);	

	\filldraw (26.3, 14.3) circle (2pt) node[label=below :{\lab{}}] (r) {};
	\filldraw (26.3, 10.3) circle (2pt) node[label= left : {\lab{}}] (s) {};
	\filldraw (24.8, 14.3) circle (2pt) node[label= below : {\lab{}}] (t) {};
	\filldraw (24.8, 10.3) circle (2pt) node[label= left : {\lab{}}] (u) {};
	\draw[->] (r) to (s);	
	\draw[->] (r) to (t);	
	\draw[->] (s) to (u);	
	\draw[->] (t) to (u);
	
	\filldraw (24.5, 12.7) circle (2pt) node[label=below :{\lab{}}] (r) {};
	\filldraw (24.5, 9.3) circle (2pt) node[label= left : {\lab{}}] (s) {};
	\filldraw (23.3, 12.7) circle (2pt) node[label= below : {\lab{}}] (t) {};
	\filldraw (23.3, 9.3) circle (2pt) node[label= left : {\lab{}}] (u) {};
	\draw[->] (r) to (s);	
	\draw[->] (r) to (t);	
	\draw[->] (s) to (u);	
	\draw[->] (t) to (u);
	
	\filldraw (24.25, 12.3) circle (2pt) node[label=below :{\lab{}}] (r) {};
	\filldraw (24.25, 9.7) circle (2pt) node[label= left : {\lab{}}] (s) {};
	\filldraw (23.65, 12.3) circle (2pt) node[label= below : {\lab{}}] (t) {};
	\filldraw (23.65, 9.7) circle (2pt) node[label= left : {\lab{}}] (u) {};
	\draw[->] (r) to (s);	
	\draw[->] (r) to (t);	
	\draw[->] (s) to (u);	
	\draw[->] (t) to (u);	
	
	\filldraw (23.7, 17.7) circle (2pt) node[label=below :{\lab{}}] (r) {};
	\filldraw (23.7, 14.3) circle (2pt) node[label= left : {\lab{}}] (s) {};
	\filldraw (21.3, 17.7) circle (2pt) node[label= below : {\lab{}}] (t) {};
	\filldraw (21.3, 14.3) circle (2pt) node[label= left : {\lab{}}] (u) {};
	\draw[->] (r) to (s);	
	\draw[->] (r) to (t);	
	\draw[->] (s) to (u);	
	\draw[->] (t) to (u);	

	\filldraw (23.3, 17.3) circle (2pt) node[label=below :{\lab{}}] (r) {};
	\filldraw (23.3, 14.7) circle (2pt) node[label= left : {\lab{}}] (s) {};
	\filldraw (21.7, 17.3) circle (2pt) node[label= below : {\lab{}}] (t) {};
	\filldraw (21.7, 14.7) circle (2pt) node[label= left : {\lab{}}] (u) {};
	\draw[->] (r) to (s);	
	\draw[->] (r) to (t);	
	\draw[->] (s) to (u);	
	\draw[->] (t) to (u);	

	\filldraw (20.7, 14.7) circle (2pt) node[label=below :{\lab{}}] (r) {};
	\filldraw (20.7, 12.3) circle (2pt) node[label= left : {\lab{}}] (s) {};
	\filldraw (19.3, 14.7) circle (2pt) node[label= below : {\lab{}}] (t) {};
	\filldraw (19.3, 12.3) circle (2pt) node[label= left : {\lab{}}] (u) {};
	\draw[->] (r) to (s);	
	\draw[->] (r) to (t);	
	\draw[->] (s) to (u);	
	\draw[->] (t) to (u);
	
	\filldraw (20.3, 14.3) circle (2pt) node[label=below :{\lab{}}] (r) {};
	\filldraw (20.3, 12.7) circle (2pt) node[label= left : {\lab{}}] (s) {};
	\filldraw (19.7, 14.3) circle (2pt) node[label= below : {\lab{}}] (t) {};
	\filldraw (19.7, 12.7) circle (2pt) node[label= left : {\lab{}}] (u) {};
	\draw[->] (r) to (s);	
	\draw[->] (r) to (t);	
	\draw[->] (s) to (u);	
	\draw[->] (t) to (u);
	
	\filldraw (19.7, 19.7) circle (2pt) node[label=below :{\lab{}}] (r) {};
	\filldraw (19.7, 17.3) circle (2pt) node[label= left : {\lab{}}] (s) {};
	\filldraw (17.3, 19.7) circle (2pt) node[label= below : {\lab{}}] (t) {};
	\filldraw (17.3, 17.3) circle (2pt) node[label= left : {\lab{}}] (u) {};
	\draw[->] (r) to (s);	
	\draw[->] (r) to (t);	
	\draw[->] (s) to (u);	
	\draw[->] (t) to (u);

	\filldraw (19.3, 19.3) circle (2pt) node[label=below :{\lab{}}] (r) {};
	\filldraw (19.3, 17.7) circle (2pt) node[label= left : {\lab{}}] (s) {};
	\filldraw (17.7, 19.3) circle (2pt) node[label= below : {\lab{}}] (t) {};
	\filldraw (17.7, 17.7) circle (2pt) node[label= left : {\lab{}}] (u) {};
	\draw[->] (r) to (s);	
	\draw[->] (r) to (t);	
	\draw[->] (s) to (u);	
	\draw[->] (t) to (u);

	\filldraw (15.7, 15.7) circle (2pt) node[label=below :{\lab{}}] (r) {};
	\filldraw (15.7, 14.3) circle (2pt) node[label= left : {\lab{}}] (s) {};
	\filldraw (14.3, 15.7) circle (2pt) node[label= below : {\lab{}}] (t) {};
	\filldraw (14.3, 14.3) circle (2pt) node[label= left : {\lab{}}] (u) {};
	\draw[->] (r) to (s);	
	\draw[->] (r) to (t);	
	\draw[->] (s) to (u);	
	\draw[->] (t) to (u);
	
	\filldraw (15.35, 15.35) circle (2pt) node[label=below :{\lab{}}] (r) {};
	\filldraw (15.35, 14.65) circle (2pt) node[label= left : {\lab{}}] (s) {};
	\filldraw (14.65, 15.35) circle (2pt) node[label= below : {\lab{}}] (t) {};
	\filldraw (14.65, 14.65) circle (2pt) node[label= left : {\lab{}}] (u) {};
	\draw[->] (r) to (s);	
	\draw[->] (r) to (t);	
	\draw[->] (s) to (u);	
	\draw[->] (t) to (u);

	\filldraw (5.3, 29.3) circle (2pt) node[label=below :{\lab{}}] (n) {};
	\filldraw (0.5, 29.3) circle (2pt) node[label= left : {\lab{}}] (o) {};
	\filldraw (5.3, 28.5) circle (2pt) node[label= below : {\lab{}}] (p) {};
	\filldraw (0.5, 28.5) circle (2pt) node[label= left : {\lab{}}] (q) {};
	\draw[->] (n) to (o);	
	\draw[->] (n) to (p);	
	\draw[->] (o) to (q);	
	\draw[->] (p) to (q);	

	\filldraw (9.7, 27.7) circle (2pt) node[label=below :{\lab{}}] (r) {};
	\filldraw (5.3, 27.7) circle (2pt) node[label= left : {\lab{}}] (s) {};
	\filldraw (9.7, 26.3) circle (2pt) node[label= below : {\lab{}}] (t) {};
	\filldraw (5.3, 26.3) circle (2pt) node[label= left : {\lab{}}] (u) {};
	\draw[->] (r) to (s);	
	\draw[->] (r) to (t);	
	\draw[->] (s) to (u);	
	\draw[->] (t) to (u);

	\filldraw (9.3, 27.3) circle (2pt) node[label=below :{\lab{}}] (r) {};
	\filldraw (5.7, 27.3) circle (2pt) node[label= left : {\lab{}}] (s) {};
	\filldraw (9.3, 26.7) circle (2pt) node[label= below : {\lab{}}] (t) {};
	\filldraw (5.7, 26.7) circle (2pt) node[label= left : {\lab{}}] (u) {};
	\draw[->] (r) to (s);	
	\draw[->] (r) to (t);	
	\draw[->] (s) to (u);	
	\draw[->] (t) to (u);	

	\filldraw (14.3, 26.3) circle (2pt) node[label=below :{\lab{}}] (r) {};
	\filldraw (10.3, 26.3) circle (2pt) node[label= left : {\lab{}}] (s) {};
	\filldraw (14.3, 24.8) circle (2pt) node[label= below : {\lab{}}] (t) {};
	\filldraw (10.3, 24.8) circle (2pt) node[label= left : {\lab{}}] (u) {};
	\draw[->] (r) to (s);	
	\draw[->] (r) to (t);	
	\draw[->] (s) to (u);	
	\draw[->] (t) to (u);
	
	\filldraw (12.7, 24.5) circle (2pt) node[label=below :{\lab{}}] (r) {};
	\filldraw (9.3, 24.5) circle (2pt) node[label= left : {\lab{}}] (s) {};
	\filldraw (12.7, 23.3) circle (2pt) node[label= below : {\lab{}}] (t) {};
	\filldraw (9.3, 23.3) circle (2pt) node[label= left : {\lab{}}] (u) {};
	\draw[->] (r) to (s);	
	\draw[->] (r) to (t);	
	\draw[->] (s) to (u);	
	\draw[->] (t) to (u);
	
	\filldraw (12.3, 24.25) circle (2pt) node[label=below :{\lab{}}] (r) {};
	\filldraw (9.7, 24.25) circle (2pt) node[label= left : {\lab{}}] (s) {};
	\filldraw (12.3, 23.65) circle (2pt) node[label= below : {\lab{}}] (t) {};
	\filldraw (9.7, 23.65) circle (2pt) node[label= left : {\lab{}}] (u) {};
	\draw[->] (r) to (s);	
	\draw[->] (r) to (t);	
	\draw[->] (s) to (u);	
	\draw[->] (t) to (u);	
	
	\filldraw (17.7, 23.7) circle (2pt) node[label=below :{\lab{}}] (r) {};
	\filldraw (14.3, 23.7) circle (2pt) node[label= left : {\lab{}}] (s) {};
	\filldraw (17.7, 21.3) circle (2pt) node[label= below : {\lab{}}] (t) {};
	\filldraw (14.3, 21.3) circle (2pt) node[label= left : {\lab{}}] (u) {};
	\draw[->] (r) to (s);	
	\draw[->] (r) to (t);	
	\draw[->] (s) to (u);	
	\draw[->] (t) to (u);	

	\filldraw (17.3, 23.3) circle (2pt) node[label=below :{\lab{}}] (r) {};
	\filldraw (14.7, 23.3) circle (2pt) node[label= left : {\lab{}}] (s) {};
	\filldraw (17.3, 21.7) circle (2pt) node[label= below : {\lab{}}] (t) {};
	\filldraw (14.7, 21.7) circle (2pt) node[label= left : {\lab{}}] (u) {};
	\draw[->] (r) to (s);	
	\draw[->] (r) to (t);	
	\draw[->] (s) to (u);	
	\draw[->] (t) to (u);	

	\filldraw (14.7, 20.7) circle (2pt) node[label=below :{\lab{}}] (r) {};
	\filldraw (12.3, 20.7) circle (2pt) node[label= left : {\lab{}}] (s) {};
	\filldraw (14.7, 19.3) circle (2pt) node[label= below : {\lab{}}] (t) {};
	\filldraw (12.3, 19.3) circle (2pt) node[label= left : {\lab{}}] (u) {};
	\draw[->] (r) to (s);	
	\draw[->] (r) to (t);	
	\draw[->] (s) to (u);	
	\draw[->] (t) to (u);
	
	\filldraw (14.3, 20.3) circle (2pt) node[label=below :{\lab{}}] (r) {};
	\filldraw (12.7, 20.3) circle (2pt) node[label= left : {\lab{}}] (s) {};
	\filldraw (14.3, 19.7) circle (2pt) node[label= below : {\lab{}}] (t) {};
	\filldraw (12.7, 19.7) circle (2pt) node[label= left : {\lab{}}] (u) {};
	\draw[->] (r) to (s);	
	\draw[->] (r) to (t);	
	\draw[->] (s) to (u);	
	\draw[->] (t) to (u);
	
\end{tikzpicture}
\]
\caption{The complex $C_6 \otimes C_6$, after a change of basis.}
\label{fig:C6xC6}
\end{figure}


\begin{figure}
\[
\begin{tikzpicture}[scale=0.5]
	\draw[step=1, black!30!white, very thin] (-0.5, -0.5) grid (31.5, 31.5);

	\filldraw (30.5, 0.5) circle (2pt) node[label=right:{\lab{u}}] (w_0) {};
	\filldraw (29.5, 0.5) circle (2pt) node[label=left:{\lab{t}}] (w_1) {};
	\filldraw (29.5, 5.5) circle (2pt) node[label=right:{\lab{s}}] (w_2) {};
	\filldraw (28.5, 5.5) circle (2pt) node[label=left:{\lab{r}}] (w_3) {};
	\filldraw (28.5, 10.5) circle (2pt) node[label=right:{\lab{q}}] (w_4) {};
	\filldraw (26.5, 10.5) circle (2pt) node[label=left:{\lab{p}}] (w_5) {};
	\filldraw (26.5, 14.5) circle (2pt) node[label= right:{\lab{o}}] (w_6) {};
	\filldraw (24.5, 14.5) circle (2pt) node[label=left:{\lab{n}}] (w_7) {};
	\filldraw (24.5, 18.5) circle (2pt) node[label=right:{\lab{m}}] (w_8) {};
	\filldraw (21.7, 18.5) circle (2pt) node[label= below:{\lab{\ell}}] (w_9) {};
	\filldraw (21.7, 21.7) circle (2pt) node[label= above:{\lab{k}}] (w_10) {};
	\filldraw (18.5, 21.7) circle (2pt) node[label= left:{\lab{j}}] (w_11) {};
	\filldraw (18.5, 24.5) circle (2pt) node[label=above:{\lab{i}}] (w_12) {};
	\filldraw (14.5, 24.5) circle (2pt) node[label=below:{\lab{h}}] (w_13) {};
	\filldraw (14.5, 26.5) circle (2pt) node[label=above:{\lab{g}}] (w_14) {};
	\filldraw (10.5, 26.5) circle (2pt) node[label=below:{\lab{f}}] (w_15) {};
	\filldraw (10.5, 28.5) circle (2pt) node[label=above:{\lab{e}}] (w_16) {};
	\filldraw (5.5, 28.5) circle (2pt) node[label=below:{\lab{d}}] (w_17) {};
	\filldraw (5.5, 29.5) circle (2pt) node[label=above:{\lab{c}}] (w_18) {};
	\filldraw (0.5, 29.5) circle (2pt) node[label=below:{\lab{b}}] (w_19) {};
	\filldraw (0.5, 30.5) circle (2pt) node[label=above :{\lab{a}}] (w_20) {};
	\draw[->](w_0) to (w_1);	
	\draw[->] (w_2) to (w_1);	
	\draw[->] (w_2) to (w_3);	
	\draw[->] (w_4) to (w_3);	
	\draw[->](w_4) to (w_5);	
	\draw[->] (w_6) to (w_5);	
	\draw[->] (w_6) to (w_7);	
	\draw[->] (w_8) to (w_7);	
	\draw[->] (w_8) to (w_9);	
	\draw[->] (w_10) to (w_9);	
	\draw[->](w_10) to (w_11);	
	\draw[->] (w_12) to (w_11);	
	\draw[->] (w_12) to (w_13);	
	\draw[->] (w_14) to (w_13);	
	\draw[->](w_14) to (w_15);	
	\draw[->] (w_16) to (w_15);	
	\draw[->] (w_16) to (w_17);	
	\draw[->] (w_18) to (w_17);	
	\draw[->] (w_18) to (w_19);	
	\draw[->] (w_20) to (w_19);	

	\filldraw (21.3, 21.3) circle (2pt) node[label=below left:{\lab{v}}] (j) {};
	\filldraw (18.5, 21.3) circle (2pt) node[label= left : {\lab{w}}] (k) {};
	\filldraw (21.3, 18.5) circle (2pt) node[label= below : {\lab{x}}] (l) {};
	\filldraw (18.5, 18.5) circle (2pt) node[label= left : {\lab{y}}] (m) {};
	\draw[->] (j) to (k);	
	\draw[->] (j) to (l);	
	\draw[->] (k) to (m);	
	\draw[->] (l) to (m);
	
	\node[align=left] at (10, 10) {\lab{a=v_{10}v_{10}} \\ 
						\lab{b=v_9v_{10}+v_{10}v_9} \\
						\lab{c=v_{8}v_{10}+v_{10}v_8} \\ 
						\lab{d=v_{7}v_{10}+v_8v_9+v_9v_8+v_{10}v_7} \\
						\lab{e=v_6v_{10}+v_8v_8+v_{10}v_6} \\ 
						\lab{f=v_5v_{10}+v_6v_9+v_7v_8+v_8v_7+v_9v_6+v_{10}v_5} \\
						\lab{g=v_4v_{10}+v_6v_8+v_8v_6+v_{10}v_4} \\ 
						\lab{h=v_3v_{10}+v_4v_9+v_5v_8+v_6v_7+v_7v_6+v_8v_5+v_9v_4+v_{10}v_3} \\
						\lab{i=v_2v_{10}+v_4v_8+v_6v_6+v_8v_4+v_{10}v_2} \\ 
						\lab{j=v_1v_{10}+v_2v_9+v_3v_8+v_4v_7+v_5v_6+v_6v_5+v_7v_4+v_8v_3+v_9v_2+v_{10}v_1} \\
						\lab{k=v_0v_{10}+v_2v_8+v_4v_6+v_6v_4+v_8v_2+v_{10}v_0+U^{-1}(v_1v_9+v_5v_5)} \\ 
						\lab{\ell=v_0v_9+v_1v_8+v_2v_7+v_3v_6+v_4v_5+v_5v_4+v_6v_3+v_7v_2+v_8v_1+v_9v_0} \\
						\lab{m=v_0v_8+v_2v_6+v_4v_4+v_6v_2+v_8v_0} \\ 
						\lab{n=v_0v_7+v_1v_6+v_2v_5+v_3v_4+v_4v_3+v_5v_2+v_6v_1+v_7v_0} \\
						\lab{o=v_0v_6+v_2v_4+v_4v_2+v_6v_0+U^{-1}(v_1v_5+v_5v_1)} \\ 
						\lab{p=v_0v_5+v_1v_4+v_2v_3+v_3v_2+v_4v_1+v_5v_0} \\
						\lab{q=v_0v_4+v_2v_2+v_4v_0} \\ 
						\lab{r=v_0v_3+v_1v_2+v_2v_1+v_3v_0} \\
						\lab{s=v_0v_2+v_2v_0+U^{-1}v_{11}} \\ 
						\lab{t=v_0v_1+v_1v_0} \\
						\lab{u=v_0v_0} \\ 
						\lab{v=v_6v_4+v_8v_2+v_{10}v_0} \\
						\lab{w=v_6v_5+v_7v_4+v_8v_3+v_9v_2+v_{10}v_1} \\ 
						\lab{x=v_5v_4+v_6v_3+v_7v_2+v_8v_1+v_9v_0} \\
						\lab{y=v_5v_5}
						};
	
\end{tikzpicture}
\]
\caption{The summand of $C_6 \otimes C_6$ that is isomorphic to $C_5$.}
\label{fig:C6xC6=C5}
\end{figure}
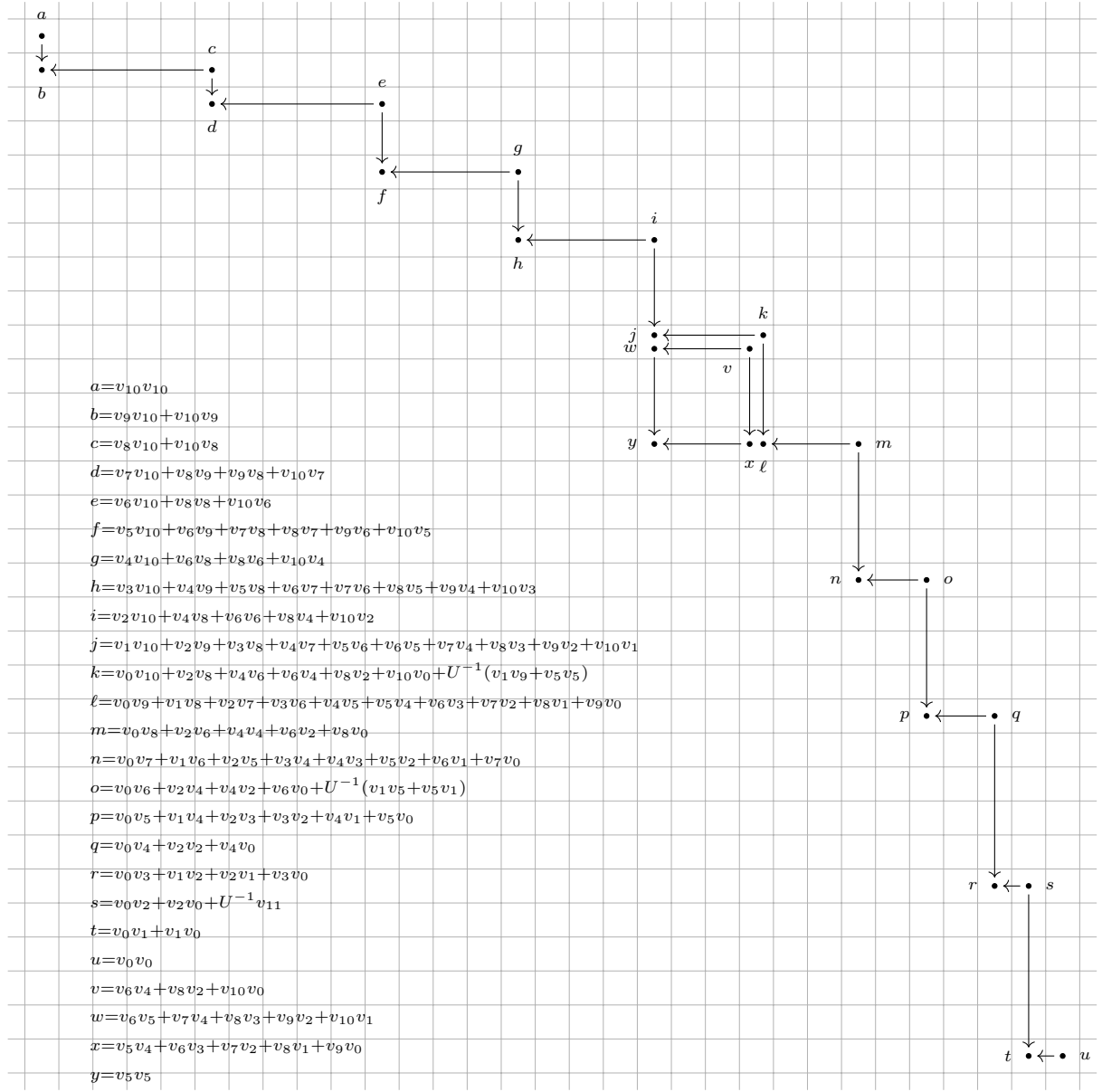


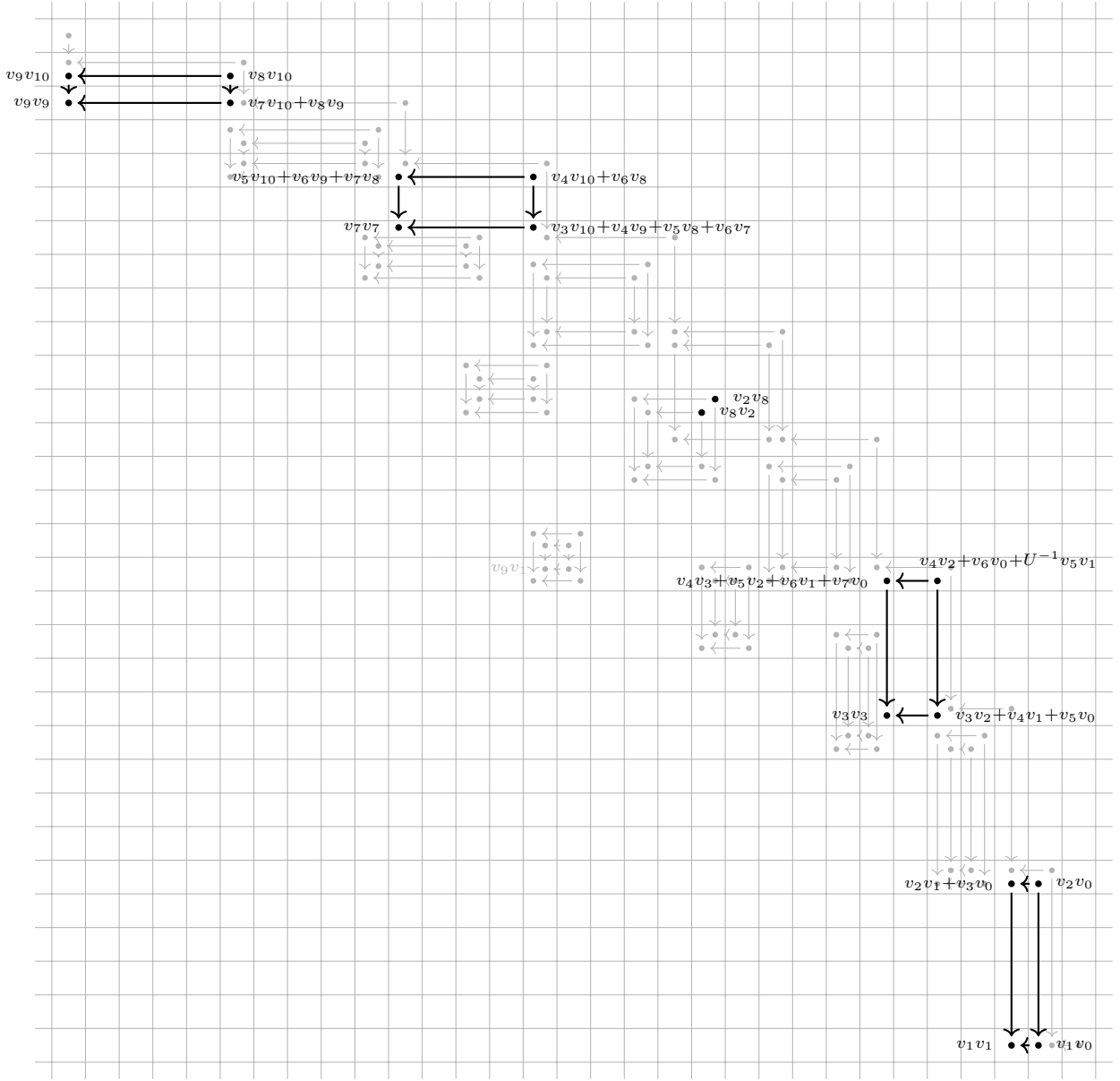
\begin{figure}
\[
\begin{tikzpicture}[scale=0.5]
	\draw[step=1, black!30!white, very thin] (-0.5, -0.5) grid (31.5, 31.5);

\begin{scope}[black!30!white]

	\filldraw (30.5, 0.5) circle (2pt) node[label=right:{\lab{}}] (w_0) {};
	\filldraw (29.7, 0.5) circle (2pt) node[label=above right:{\lab{}}] (w_1) {};
	\filldraw (29.7, 5.7) circle (2pt) node[label=right:{\lab{}}] (w_2) {};
	\filldraw (28.5, 5.7) circle (2pt) node[label=left:{\lab{}}] (w_3) {};
	\filldraw (28.5, 10.5) circle (2pt) node[label=right:{\lab{}}] (w_4) {};
	\filldraw (26.7, 10.5) circle (2pt) node[label=left:{\lab{}}] (w_5) {};
	\filldraw (26.7, 14.7) circle (2pt) node[label=right : {\lab{}}] (w_6) {};
	\filldraw (24.5, 14.7) circle (2pt) node[label=left:{\lab{}}] (w_7) {};
	\filldraw (24.5, 18.5) circle (2pt) node[label=right:{\lab{}}] (w_8) {};
	\filldraw (21.7, 18.5) circle (2pt) node[label= {[xshift=60pt, yshift=0pt] \lab{}}] (w_9) {};
	\filldraw (21.7, 21.7) circle (2pt) node[label=right:{\lab{}}] (w_10) {};
	\filldraw (18.5, 21.7) circle (2pt) node[label= above left:{\lab{}}] (w_11) {};
	\filldraw (18.5, 24.5) circle (2pt) node[label=right:{\lab{}}] (w_12) {};
	\filldraw (14.7, 24.5) circle (2pt) node[label=above right:{\lab{}}] (w_13) {};
	\filldraw (14.7, 26.7) circle (2pt) node[label=right:{\lab{}}] (w_14) {};
	\filldraw (10.5, 26.7) circle (2pt) node[label=above right:{\lab{}}] (w_15) {};
	\filldraw (10.5, 28.5) circle (2pt) node[label=right:{\lab{}}] (w_16) {};
	\filldraw (5.7, 28.5) circle (2pt) node[label=left:{\lab{}}] (w_17) {};
	\filldraw (5.7, 29.7) circle (2pt) node[label=above:{\lab{}}] (w_18) {};
	\filldraw (0.5, 29.7) circle (2pt) node[label=left:{\lab{}}] (w_19) {};
	\filldraw (0.5, 30.5) circle (2pt) node[label=above :{\lab{}}] (w_20) {};
	\draw[->](w_0) to (w_1);	
	\draw[->] (w_2) to (w_1);	
	\draw[->] (w_2) to (w_3);	
	\draw[->] (w_4) to (w_3);	
	\draw[->](w_4) to (w_5);	
	\draw[->] (w_6) to (w_5);	
	\draw[->] (w_6) to (w_7);	
	\draw[->] (w_8) to (w_7);	
	\draw[->] (w_8) to (w_9);	
	\draw[->] (w_10) to (w_9);	
	\draw[->](w_10) to (w_11);	
	\draw[->] (w_12) to (w_11);	
	\draw[->] (w_12) to (w_13);	
	\draw[->] (w_14) to (w_13);	
	\draw[->](w_14) to (w_15);	
	\draw[->] (w_16) to (w_15);	
	\draw[->] (w_16) to (w_17);	
	\draw[->] (w_18) to (w_17);	
	\draw[->] (w_18) to (w_19);	
	\draw[->] (w_20) to (w_19);	

	\filldraw (21.3, 21.3) circle (2pt) node[label=below :{\lab{}}] (j) {};
	\filldraw (18.5, 21.3) circle (2pt) node[label= left : {\lab{}}] (k) {};
	\filldraw (21.3, 18.5) circle (2pt) node[label= below : {\lab{}}] (l) {};
	\filldraw (18.5, 18.5) circle (2pt) node[label= left : {\lab{}}] (m) {};
	\draw[->] (j) to (k);	
	\draw[->] (j) to (l);	
	\draw[->] (k) to (m);	
	\draw[->] (l) to (m);
	
	\filldraw (27.7, 9.7) circle (2pt) node[label=below :{\lab{}}] (r) {};
	\filldraw (27.7, 5.3) circle (2pt) node[label= left : {\lab{}}] (s) {};
	\filldraw (26.3, 9.7) circle (2pt) node[label= below : {\lab{}}] (t) {};
	\filldraw (26.3, 5.3) circle (2pt) node[label= left : {\lab{}}] (u) {};
	\draw[->] (r) to (s);	
	\draw[->] (r) to (t);	
	\draw[->] (s) to (u);	
	\draw[->] (t) to (u);

	\filldraw (27.3, 9.3) circle (2pt) node[label=below :{\lab{}}] (r) {};
	\filldraw (27.3, 5.7) circle (2pt) node[label= left : {\lab{}}] (s) {};
	\filldraw (26.7, 9.3) circle (2pt) node[label= below : {\lab{}}] (t) {};
	\filldraw (26.7, 5.7) circle (2pt) node[label= left : {\lab{}}] (u) {};
	\draw[->] (r) to (s);	
	\draw[->] (r) to (t);	
	\draw[->] (s) to (u);	
	\draw[->] (t) to (u);	
	
	\filldraw (24.5, 12.7) circle (2pt) node[label=below :{\lab{}}] (r) {};
	\filldraw (24.5, 9.3) circle (2pt) node[label= left : {\lab{}}] (s) {};
	\filldraw (23.3, 12.7) circle (2pt) node[label= below : {\lab{}}] (t) {};
	\filldraw (23.3, 9.3) circle (2pt) node[label= left : {\lab{}}] (u) {};
	\draw[->] (r) to (s);	
	\draw[->] (r) to (t);	
	\draw[->] (s) to (u);	
	\draw[->] (t) to (u);
	
	\filldraw (24.25, 12.3) circle (2pt) node[label=below :{\lab{}}] (r) {};
	\filldraw (24.25, 9.7) circle (2pt) node[label= left : {\lab{}}] (s) {};
	\filldraw (23.65, 12.3) circle (2pt) node[label= below : {\lab{}}] (t) {};
	\filldraw (23.65, 9.7) circle (2pt) node[label= left : {\lab{}}] (u) {};
	\draw[->] (r) to (s);	
	\draw[->] (r) to (t);	
	\draw[->] (s) to (u);	
	\draw[->] (t) to (u);	
	
	\filldraw (23.7, 17.7) circle (2pt) node[label=below :{\lab{}}] (r) {};
	\filldraw (23.7, 14.3) circle (2pt) node[label= left : {\lab{}}] (s) {};
	\filldraw (21.3, 17.7) circle (2pt) node[label= below : {\lab{}}] (t) {};
	\filldraw (21.3, 14.3) circle (2pt) node[label= left : {\lab{}}] (u) {};
	\draw[->] (r) to (s);	
	\draw[->] (r) to (t);	
	\draw[->] (s) to (u);	
	\draw[->] (t) to (u);	

	\filldraw (23.3, 17.3) circle (2pt) node[label=below :{\lab{}}] (r) {};
	\filldraw (23.3, 14.7) circle (2pt) node[label= left : {\lab{}}] (s) {};
	\filldraw (21.7, 17.3) circle (2pt) node[label= below : {\lab{}}] (t) {};
	\filldraw (21.7, 14.7) circle (2pt) node[label= left : {\lab{}}] (u) {};
	\draw[->] (r) to (s);	
	\draw[->] (r) to (t);	
	\draw[->] (s) to (u);	
	\draw[->] (t) to (u);	

	\filldraw (20.7, 14.7) circle (2pt) node[label=below :{\lab{}}] (r) {};
	\filldraw (20.7, 12.3) circle (2pt) node[label= left : {\lab{}}] (s) {};
	\filldraw (19.3, 14.7) circle (2pt) node[label= below : {\lab{}}] (t) {};
	\filldraw (19.3, 12.3) circle (2pt) node[label= left : {\lab{}}] (u) {};
	\draw[->] (r) to (s);	
	\draw[->] (r) to (t);	
	\draw[->] (s) to (u);	
	\draw[->] (t) to (u);
	
	\filldraw (20.3, 14.3) circle (2pt) node[label=below :{\lab{}}] (r) {};
	\filldraw (20.3, 12.7) circle (2pt) node[label= left : {\lab{}}] (s) {};
	\filldraw (19.7, 14.3) circle (2pt) node[label= below : {\lab{}}] (t) {};
	\filldraw (19.7, 12.7) circle (2pt) node[label= left : {\lab{}}] (u) {};
	\draw[->] (r) to (s);	
	\draw[->] (r) to (t);	
	\draw[->] (s) to (u);	
	\draw[->] (t) to (u);

	\filldraw (15.7, 15.7) circle (2pt) node[label=below :{\lab{}}] (r) {};
	\filldraw (15.7, 14.3) circle (2pt) node[label= left : {\lab{}}] (s) {};
	\filldraw (14.3, 15.7) circle (2pt) node[label= below : {\lab{}}] (t) {};
	\filldraw (14.3, 14.3) circle (2pt) node[label= left : {\lab{}}] (u) {}; 
	\draw[->] (r) to (s);	
	\draw[->] (r) to (t);	
	\draw[->] (s) to (u);	
	\draw[->] (t) to (u);
	
	\filldraw (15.35, 15.35) circle (2pt) node[label=below :{\lab{}}] (r) {};
	\filldraw (15.35, 14.65) circle (2pt) node[label= left : {\lab{}}] (s) {};
	\filldraw (14.65, 15.35) circle (2pt) node[label= below : {\lab{}}] (t) {};
	\filldraw (14.65, 14.65) circle (2pt) node[label= left : {\lab{v_9v_1}}] (u) {};
	\draw[->] (r) to (s);	
	\draw[->] (r) to (t);	
	\draw[->] (s) to (u);	
	\draw[->] (t) to (u);

	\filldraw (9.7, 27.7) circle (2pt) node[label=below :{\lab{}}] (r) {};
	\filldraw (5.3, 27.7) circle (2pt) node[label= left : {\lab{}}] (s) {};
	\filldraw (9.7, 26.3) circle (2pt) node[label= below : {\lab{}}] (t) {};
	\filldraw (5.3, 26.3) circle (2pt) node[label= left : {\lab{}}] (u) {};
	\draw[->] (r) to (s);	
	\draw[->] (r) to (t);	
	\draw[->] (s) to (u);	
	\draw[->] (t) to (u);

	\filldraw (9.3, 27.3) circle (2pt) node[label=below :{\lab{}}] (r) {};
	\filldraw (5.7, 27.3) circle (2pt) node[label= left : {\lab{}}] (s) {};
	\filldraw (9.3, 26.7) circle (2pt) node[label= below : {\lab{}}] (t) {};
	\filldraw (5.7, 26.7) circle (2pt) node[label= left : {\lab{}}] (u) {};
	\draw[->] (r) to (s);	
	\draw[->] (r) to (t);	
	\draw[->] (s) to (u);	
	\draw[->] (t) to (u);	

	\filldraw (12.7, 24.5) circle (2pt) node[label=below :{\lab{}}] (r) {};
	\filldraw (9.3, 24.5) circle (2pt) node[label= left : {\lab{}}] (s) {};
	\filldraw (12.7, 23.3) circle (2pt) node[label= below : {\lab{}}] (t) {};
	\filldraw (9.3, 23.3) circle (2pt) node[label= left : {\lab{}}] (u) {};
	\draw[->] (r) to (s);	
	\draw[->] (r) to (t);	
	\draw[->] (s) to (u);	
	\draw[->] (t) to (u);
	
	\filldraw (12.3, 24.25) circle (2pt) node[label=below :{\lab{}}] (r) {};
	\filldraw (9.7, 24.25) circle (2pt) node[label= left : {\lab{}}] (s) {};
	\filldraw (12.3, 23.65) circle (2pt) node[label= below : {\lab{}}] (t) {};
	\filldraw (9.7, 23.65) circle (2pt) node[label= left : {\lab{}}] (u) {};
	\draw[->] (r) to (s);	
	\draw[->] (r) to (t);	
	\draw[->] (s) to (u);	
	\draw[->] (t) to (u);	
	
	\filldraw (17.7, 23.7) circle (2pt) node[label=below :{\lab{}}] (r) {};
	\filldraw (14.3, 23.7) circle (2pt) node[label= left : {\lab{}}] (s) {};
	\filldraw (17.7, 21.3) circle (2pt) node[label= below : {\lab{}}] (t) {};
	\filldraw (14.3, 21.3) circle (2pt) node[label= left : {\lab{}}] (u) {};
	\draw[->] (r) to (s);	
	\draw[->] (r) to (t);	
	\draw[->] (s) to (u);	
	\draw[->] (t) to (u);	

	\filldraw (17.3, 23.3) circle (2pt) node[label=below :{\lab{}}] (r) {};
	\filldraw (14.7, 23.3) circle (2pt) node[label= left : {\lab{}}] (s) {};
	\filldraw (17.3, 21.7) circle (2pt) node[label= below : {\lab{}}] (t) {};
	\filldraw (14.7, 21.7) circle (2pt) node[label= left : {\lab{}}] (u) {};
	\draw[->] (r) to (s);	
	\draw[->] (r) to (t);	
	\draw[->] (s) to (u);	
	\draw[->] (t) to (u);	

	\filldraw (14.7, 20.7) circle (2pt) node[label=below :{\lab{}}] (r) {};
	\filldraw (12.3, 20.7) circle (2pt) node[label= left : {\lab{}}] (s) {};
	\filldraw (14.7, 19.3) circle (2pt) node[label= below : {\lab{}}] (t) {};
	\filldraw (12.3, 19.3) circle (2pt) node[label= left : {\lab{}}] (u) {};
	\draw[->] (r) to (s);	
	\draw[->] (r) to (t);	
	\draw[->] (s) to (u);	
	\draw[->] (t) to (u);
	
	\filldraw (14.3, 20.3) circle (2pt) node[label=below :{\lab{}}] (r) {};
	\filldraw (12.7, 20.3) circle (2pt) node[label= left : {\lab{}}] (s) {};
	\filldraw (14.3, 19.7) circle (2pt) node[label= below : {\lab{}}] (t) {};
	\filldraw (12.7, 19.7) circle (2pt) node[label= left : {\lab{}}] (u) {};
	\draw[->] (r) to (s);	
	\draw[->] (r) to (t);	
	\draw[->] (s) to (u);	
	\draw[->] (t) to (u);

	\filldraw (19.7, 19.7) circle (2pt) node[label= right :{\lab{}}] (r) {};
	\filldraw (19.7, 17.3) circle (2pt) node[label= right : {\lab{}}] (s) {}; 
	\filldraw (17.3, 19.7) circle (2pt) node[label= left : {\lab{}}] (t) {}; 
	\filldraw (17.3, 17.3) circle (2pt) node[label= left : {\lab{}}] (u) {}; 
	\draw[->] (r) to (s);	
	\draw[->] (r) to (t);	
	\draw[->] (s) to (u);	
	\draw[->] (t) to (u);

	\filldraw (19.3, 19.3) circle (2pt) node[label=right :{\lab{}}] (r) {};
	\filldraw (19.3, 17.7) circle (2pt) node[label= right : {\lab{}}] (s) {}; 
	\filldraw (17.7, 19.3) circle (2pt) node[label= left : {\lab{}}] (t) {}; 
	\filldraw (17.7, 17.7) circle (2pt) node[label= left : {\lab{}}] (u) {}; 
	\draw[->] (r) to (s);	
	\draw[->] (r) to (t);	
	\draw[->] (s) to (u);	
	\draw[->] (t) to (u);
	
\end{scope}

\begin{scope}[thick]
	\filldraw (19.7, 19.7) circle (2pt) node[label= right :{\lab{v_2v_8}}] (r) {};
	\filldraw (19.3, 19.3) circle (2pt) node[label=right :{\lab{v_8v_2}}] (r) {};

	\filldraw (5.3, 29.3) circle (2pt) node[label=right :{\lab{v_8v_{10}}}] (n) {};
	\filldraw (0.5, 29.3) circle (2pt) node[label= left : {\lab{v_9v_{10}}}] (o) {};
	\filldraw (5.3, 28.5) circle (2pt) node[label= right : {\lab{v_7v_{10}+v_8v_9}}] (p) {};
	\filldraw (0.5, 28.5) circle (2pt) node[label= left : {\lab{v_9v_9}}] (q) {};
	\draw[->] (n) to (o);	
	\draw[->] (n) to (p);	
	\draw[->] (o) to (q);	
	\draw[->] (p) to (q);	
	
	\filldraw (29.3, 5.3) circle (2pt) node[label=right :{\lab{v_2v_{0}}}] (n) {};
	\filldraw (29.3, 0.5) circle (2pt) node[label= right : {\lab{v_1v_0}}] (o) {};
	\filldraw (28.5, 5.3) circle (2pt) node[label= left : {\lab{v_2v_1+v_3v_0}}] (p) {};
	\filldraw (28.5, 0.5) circle (2pt) node[label= left : {\lab{v_1v_1}}] (q) {};
	\draw[->] (n) to (o);	
	\draw[->] (n) to (p);	
	\draw[->] (o) to (q);	
	\draw[->] (p) to (q);	

	\filldraw (14.3, 26.3) circle (2pt) node[label= right : {\lab{v_4v_{10}+v_6v_8}}] (r) {};
	\filldraw (10.3, 26.3) circle (2pt) node[label= left : {\lab{v_5v_{10}+v_6v_9+v_7v_8}}] (s) {};
	\filldraw (14.3, 24.8) circle (2pt) node[label= right : {\lab{v_3v_{10}+v_4v_9+v_5v_8+v_6v_7}}] (t) {};
	\filldraw (10.3, 24.8) circle (2pt) node[label= left : {\lab{v_7v_7}}] (u) {};
	\draw[->] (r) to (s);	
	\draw[->] (r) to (t);	
	\draw[->] (s) to (u);	
	\draw[->] (t) to (u);

	\filldraw (26.3, 14.3) circle (2pt) node[label={[xshift=30pt, yshift=-2pt] \lab{v_4v_2+v_6v_0+U^{-1}v_5v_1}}] (r) {};
	\filldraw (26.3, 10.3) circle (2pt) node[label= right : {\lab{v_3v_2+v_4v_1+v_5v_0}}] (s) {};
	\filldraw (24.8, 14.3) circle (2pt) node[label= left : {\lab{v_4v_3+v_5v_2+v_6v_1+v_7v_0}}] (t) {};
	\filldraw (24.8, 10.3) circle (2pt) node[label= left : {\lab{v_3v_3}}] (u) {};
	\draw[->] (r) to (s);	
	\draw[->] (r) to (t);	
	\draw[->] (s) to (u);	
	\draw[->] (t) to (u);

\end{scope}
	
\end{tikzpicture}
\]
\caption{The complex $C_6 \otimes C_6$, after a change of basis, with some basis elements labelled.}
\label{fig:C6xC6-4boxes}
\end{figure}

\subsection{The \texorpdfstring{$\iota_K$}{iota-K}-complex of \texorpdfstring{$-2T_{6,7} \# T_{6,13} \# -T_{2,3;2,5}$}{-2T(6,7) \# T(6,13) \# -T(2,3;2,5) }}\label{sec:iotaKexample}
Recall that $K \subset S^3$ is an \emph{L-space knot} if $S^3_r(K)$ is an L-space for some $r>0$. We say that $K$ is a \emph{negative L-space knot} if $S^3_{-r}(K)$ is an L-space for some $r>0$. The following lemma, combining results of \cite{OSlens} and \cite{HMInvolutive}, describes the $\iota_K$-complex associated to a negative L-space knot $K$ in terms of the Alexander polynomial $\Delta_K(t)$.

\begin{lem}\label{lem:Lspaceknot}
Let $K$ be a negative L-space knot with Alexander polynomial 
\[ \Delta_K(t) = \sum_{i=0}^n (-1)^i t^{a_i}, \]
where $(a_i)_{i=0}^n$ is a decreasing sequence of integers and $n$ is even. Then $\CFKi(K)$ is generated over $\bF[U, U^{-1}]$ by $x_i, 0 \leq i \leq n$, with Maslov gradings
\begin{align*}
	M(x_i) &= 0, \quad i \textup{ even} \\
	M(x_i) &= -1, \quad i \textup{ odd} 	
\end{align*}
where the $\Z \oplus \Z$-filtration level of $x_0$ is $\cF(x_0) = (0, -g(K))$ and relative filtration of $x_i, x_{i+1}$ is
\begin{align*}
	\cF(x_{i}) - \cF(x_{i+1}) &= (a_i-a_{i+1}, 0), \quad i \textup{ even} \\
	\cF(x_{i}) - \cF(x_{i+1}) &= (0, a_{i+1} - a_{i}), \quad i \textup{ odd}.
\end{align*}
The differential is given by
\begin{align*}
	\d x_0 &= x_1 \\
	\d x_i &= x_{i-1} + x_{i+1}, \quad 2 \leq i \leq n-2 \textup{ even} \\
	\d x_n &= x_{n-1} \\
	\d x_i &=0, \quad i \textup{ odd}.
\end{align*}
The endomorphism $\iota_K$ is given by
\[ \iota_K(x_i) = x_{n-i}, \quad 0 \leq i \leq n. \] 
\end{lem}

\begin{proof}
The description of $\CFKi(K)$ follows from \cite[Theorem 1.2]{OSlens} (cf. \cite[Theorem 2.10]{OSSUpsilon}). The description of $\iota_K$ follows from \cite[Section 7]{HMInvolutive}.
\end{proof}

\begin{lem}\label{lem:iotaKcomplexes}
We have the following identifications:
\begin{enumerate}
	\item The complex $(C_2, \iota_2)$ in Figure \ref{fig:C1C2} where $\iota_2(x_i) = x_{4-i}$ is the $\iota_K$-complex associated to $-T_{2,3;2,5}$.
	\item The complex $(C_4, \iota_4)$ in Figure \ref{fig:C4} where $\iota_4(z_i) = z_{20-i}$ is the $\iota_K$-complex associated to $-T_{6,13}$.
	\item The complex $(C_6, \iota_6)$ in Figure \ref{fig:C6} where $\iota_5(v_i) = v_{10-i}$ is the $\iota_K$-complex associated to $-T_{6,7}$.
\end{enumerate}
\end{lem}

\begin{proof}
The Alexander polynomial of $-T_{2,3;2,5}$ is 
\[ \Delta_{-T_{2,3;2,5}}(t) = t^8-t^7+t^4-t+1, \] 
the Alexander polynomial $\Delta_{-T_{6,13}}(t)$ is
\[ t^{60}-t^{59}+t^{54}-t^{53}+t^{48}-t^{46}+t^{42}-t^{40}+t^{36}-t^{33}+t^{30}-t^{27}+t^{24}-t^{20}+t^{18}-t^{14}+t^{12}-t^{7}+t^6-t+1, \]
and the Alexander polynomial of $-T_{6,7}$ is  
\[ \Delta_{-T_{6,7}}(t) = t^{30}-t^{29}+t^{24}-t^{22}+t^{18}-t^{15}+t^{12}-t^{8}+t^5-t+1. \]
Recall that negative torus knots admit negative L-space surgeries, as does $-T_{2,3;2,5}$ \cite[Theorem 1.10]{HeddencablingII} (cf. \cite{HomcablesLspace}). The result now follows from Lemma \ref{lem:Lspaceknot}.
\end{proof}

We now combine Lemma \ref{lem:iotaKcomplexes} with the computations of Section \ref{sec:tensorprods} to compute the $\iota_K$-complex of $K = -2T_{6,7} \# T_{6,13} \# -T_{2,3;2,5}$.

\begin{prop}\label{prop:iotaKcomplex}
The $\iota_K$-complex of $K = -2T_{6,7} \# T_{6,13} \# -T_{2,3;2,5}$ is locally equivalent to the complex $(C_3, \iota_3)$ in Figure \ref{fig:C3}, where $\iota_3$ is given by
\begin{align*}
	\iota_3(y_i) &= y_{4-i} , \quad i \neq 2 \\
	\iota_3(y_2) &= y_2 + U^{-1}i \\
	\iota_3(f) &= f+y_2 \\ 
	\iota_3(g) &= h+y_1 \\ 
	\iota_3(h) &= g+y_3 \\ 
	\iota_3(i) &= i.
\end{align*}
\end{prop}

\begin{proof}
By Lemma \ref{lem:iotaKcomplexes} and the K\"unneth formula for $\iota_K$-complexes, we have that the $\iota_K$-complex of $K$ is $C_6 \otimes C_6 \otimes C_4^* \otimes C_2$.

We have the following sequence of equivalences, where $\simeq$ denotes local equivalence:
\begin{align*}
	C_6 \otimes C_6 \otimes C_4^* \otimes C_2 &\simeq C_6 \otimes C_6 \otimes C_4^* \otimes C_1^* \otimes C_1 \otimes C_2 \\
		&\simeq C_5 \otimes C_5^* \otimes C_1 \otimes C_2 \\
		&\simeq C_1 \otimes C_2 \\
		&\simeq C_3,
\end{align*}
where the first and third equivalences follow from the fact that $C_i$ and $C_i^*$ are inverses in the local equivalence group of $\iota_K$-complexes, the second equivalence follows from  Lemmas \ref{lem:C6xC6} and \ref{lem:C1xC4} (and duality, in the latter case), and the fourth equivalence follows from Lemma \ref{lem:C1xC2}.
\end{proof}

\subsection{Surgery along \texorpdfstring{$-2T_{6,7} \# T_{6,13} \# -T_{2,3;2,5}$}{-2T(6,7)\# T(6,13) \# -T(2,3;2,5)}}\label{sec:surgeryalongK}

In this section, we employ Proposition \ref{prop:local-equivalence-class} to compute the local equivalence class of $(\CF^-(S^3_{+1}(K), \iota)$, where $K = -2T_{6,7} \# T_{6,13} \# -T_{2,3;2,5}$. 

\begin{prop}\label{prop:knotexample}
Let $K=-2T_{6,7} \# T_{6,13} \# -T_{2,3;2,5}$. Then $(\CF^-(S^3_{+1}(K), \iota)$ is locally equivalent to $(C, \iota_C)$ where $C$ is generated by $x_0, x_1, x_2, x_3, x_4$ where
\begin{align*}
	\d x_0 &= 0 		& \omega x_0 &= 0 \\
	\d x_1 &= Ux_2		&\omega  x_1 &= x_0 \\
	\d x_2 &=0 		& \omega  x_2 &=0  \\
	\d x_3 &= U^3x_4	& \omega  x_3 &= x_2 \\
	\d x_4 &= 0		& \omega  x_4 &= 0,
\end{align*}
where $\omega = 1 + \iota_C$.

In particular, in the notation of Section \ref{sec:alec}, the almost local equivalence class of $S^3_{+1}(K)$ is $(+, -1, +, -3)$.
\end{prop}

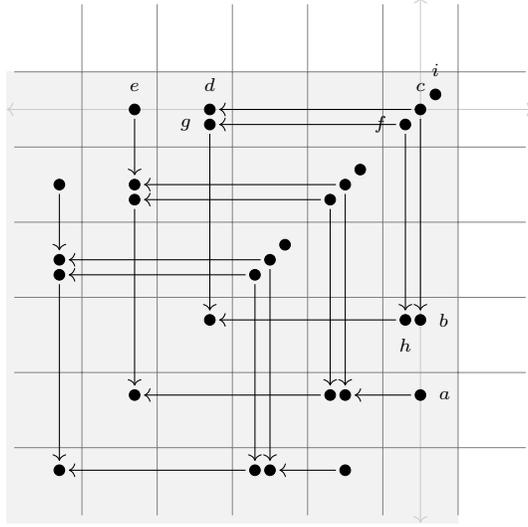
\begin{figure}
\begin{tikzpicture}[scale=1]

\filldraw[black!5!white] (-5, -5) rectangle (1, 1);

	\begin{scope}[thin, black!20!white]
		\draw [<->] (-5, 0.5) -- (2, 0.5);
		\draw [<->] (0.5, -5) -- (0.5, 2);
	\end{scope}
	\draw[step=1, black!50!white, very thin] (-4.9, -4.9) grid (1.9, 1.9);

	\filldraw (-0+0.7, -0+0.7) circle (2pt) node[label= above : {\lab{i}}] (i) {};
	\filldraw (-1+0.7, -1+0.7) circle (2pt) node[label= below : {\lab{}}] (i) {};
	\filldraw (-2+0.7, -2+0.7) circle (2pt) node[label= below : {\lab{}}] (i) {};
	\filldraw (-3+3.5, -3+0.7) circle (2pt) node[label= right:{\lab{b}}] (y_1) {};
	\filldraw (-3+3.5, -3+3.5) circle (2pt) node[label=above :{\lab{c}}] (y_2) {};
	\filldraw (-3+0.7, -3+3.5) circle (2pt) node[label=above:{\lab{d}}] (y_3) {};	
	\draw[->] (y_2) to (y_1);	
	\draw[->] (y_2) to (y_3);		
	\filldraw (-3+3.3, -3+3.3) circle (2pt) node[label= left:{\lab{f}}] (f) {};
	\filldraw (-3+0.7, -3+3.3) circle (2pt) node[label= left : {\lab{g}}] (g) {};
	\filldraw (-3+3.3, -3+0.7) circle (2pt) node[label= below : {\lab{h}}] (h) {};
	\filldraw (-3+0.7, -3+0.7) circle (2pt) node[label= below : {\lab{}}] (i) {};
	\draw[->] (f) to (h);	
	\draw[->] (f) to (g);	
	\draw[->] (g) to (i);	
	\draw[->] (h) to (i);

	\filldraw (-4+4.5, -4+0.7) circle (2pt) node[label=right:{\lab{a}}] (y_0) {};
	\filldraw (-4+3.5, -4+0.7) circle (2pt) node[label=above right:{\lab{}}] (y_1) {};
	\filldraw (-4+3.5, -4+3.5) circle (2pt) node[label=right:{\lab{}}] (y_2) {};
	\filldraw (-4+0.7, -4+3.5) circle (2pt) node[label=above right:{\lab{}}] (y_3) {};
	\filldraw (-4+0.7, -4+4.5) circle (2pt) node[label=above:{\lab{e}}] (y_4) {};
	\draw[->](y_0) to (y_1);	
	\draw[->] (y_2) to (y_1);	
	\draw[->] (y_2) to (y_3);	
	\draw[->] (y_4) to (y_3);	

	\filldraw (-4+3.3, -4+3.3) circle (2pt) node[label=below left:{\lab{}}] (f) {};
	\filldraw (-4+0.7, -4+3.3) circle (2pt) node[label= left : {\lab{}}] (g) {};
	\filldraw (-4+3.3, -4+0.7) circle (2pt) node[label= below : {\lab{}}] (h) {};
	\filldraw (-4+0.7, -4+0.7) circle (2pt) node[label= below : {\lab{}}] (i) {};
	\draw[->] (f) to (h);	
	\draw[->] (f) to (g);	
	\draw[->] (g) to (i);	
	\draw[->] (h) to (i);

	\filldraw (-5+4.5, -5+0.7) circle (2pt) node[label=right:{\lab{}}] (y_0) {};
	\filldraw (-5+3.5, -5+0.7) circle (2pt) node[label=above right:{\lab{}}] (y_1) {};
	\filldraw (-5+3.5, -5+3.5) circle (2pt) node[label=right:{\lab{}}] (y_2) {};
	\filldraw (-5+0.7, -5+3.5) circle (2pt) node[label=above right:{\lab{}}] (y_3) {};
	\filldraw (-5+0.7, -5+4.5) circle (2pt) node[label=above:{\lab{}}] (y_4) {};
	\draw[->](y_0) to (y_1);	
	\draw[->] (y_2) to (y_1);	
	\draw[->] (y_2) to (y_3);	
	\draw[->] (y_4) to (y_3);	

	\filldraw (-5+3.3, -5+3.3) circle (2pt) node[label=below left:{\lab{}}] (f) {};
	\filldraw (-5+0.7, -5+3.3) circle (2pt) node[label= left : {\lab{}}] (g) {};
	\filldraw (-5+3.3, -5+0.7) circle (2pt) node[label= below : {\lab{}}] (h) {};
	\filldraw (-5+0.7, -5+0.7) circle (2pt) node[label= below : {\lab{}}] (i) {};
	\draw[->] (f) to (h);	
	\draw[->] (f) to (g);	
	\draw[->] (g) to (i);	
	\draw[->] (h) to (i);
\end{tikzpicture}
\caption{The complex $A^-_0$ associated to $C_3$ with (appropriate $U$-powers of) $y_1, y_2, y_3, y_4, y_5$ renamed $a, b, c, d, e$ respectively.}
\label{fig:A-0}		
\end{figure}

\begin{proof}
By Proposition \ref{prop:iotaKcomplex}, we have that $(\CFKi(K), \iota_K)$ is locally equivalent to $(C_3, \iota_3)$. Let $(A^-_0, \iota_K)$ denote $C_3\{ \max(i,j) \leq 0 \}$ with the homotopy involution induced by $\iota_3$.  Proposition~\ref{prop:local-equivalence-class} implies that $(A^-_0, \iota_K)$ is locally equivalent to $(\CF^-(S^3_{+1}(K), \iota)$. Thus it suffices to show that $(A^-_0, \iota_K)$ is locally equivalent to the $\iota$-complex $(C, \iota_C)$ in the statement of Proposition \ref{prop:knotexample}.

We have that $(A^-_0, \iota_K)$ is generated over $\bF[U]$ by $a, b, c, d, e, f, g, h, i$ with
\begin{align*}
\d a &= Ub  	&\iota_K(a) &= e 		&\omega(a) &= a+e \\
\d b &= 0  		&\iota_K(b) &= d 		&\omega(b) &= b+d \\
\d c &= b+d 	&\iota_K(c) &= c + U^{2}i 	&\omega(c) &= U^2i \\
\d d &= 0 		& \iota_K(d) &= b 		&\omega(d) &= b+d \\
\d e &= Ud	&\iota_K(e) &= a 		&\omega(e) &= a+e \\
\d f &= h+g 	&\iota_K(f) &= f+c		&\omega(f) &= c  \\
\d g &= U^3i 	&\iota_K(g) &= b+h		&\omega(g) &= b+g+h  \\
\d h &= U^3 i 	&\iota_K(h) &= d+g		&\omega(h) &= d+g+h \\
\d i &= 0		&\iota_K (i) &= i			&\omega(i) &= 0
\end{align*}
as shown in Figure \ref{fig:A-0}, where $y_0, U^{-1}y_1, U^{-1}y_2, U^{-1}y_3, y_4$ have been renamed $a, b, c, d, e$ respectively, $f, g, h$ have been shifted by $U^{-1}$, and $i$ has been shifted by $U^{-4}$. (We include $\omega = \id + \iota_K$ above in order to aid with verifying the calculations below.) 

Consider $\psi \co A^-_0 \to C$ defined to be
\begin{align*}
	\psi(a) &= x_1 \\
	\psi(b) &= x_2 \\
	\psi(c) &= 0 \\
	\psi(d) &= x_2 \\
	\psi(e) &= x_0+x_1 \\
	\psi(f) &= 0 \\
	\psi(g) &= x_3 \\
	\psi(h) &= x_3 \\
	\psi(i) &= x_4.
\end{align*}
It is straightforward to verify that $\psi$ is a chain map and that $\omega \psi = \psi \omega$. Moreover, since $[a+Uc+e]$ is a $U$-nontorsion cycle in $H_*(A^-_0)$ and $[x_0]$ is a $U$-nontorsion cycle in $H_*(C)$, it follows that $\psi$ induces an isomorphism on $U^{-1}H_*$.

Now consider $\phi \co C \to A^-_0$ defined to be
\begin{align*}
	\phi(x_0) &= a+Uc+e \\
	\phi(x_1) &= a \\
	\phi(x_2) &= b \\
	\phi(x_3) &= g \\
	\phi(x_4) &= i.
\end{align*}
It is straightforward to verify that $\phi$ is a chain map. We have that $\omega \phi \simeq \phi \omega$ via $H \co C \to A^-_0$ defined by
\begin{align*}
	H(x_0) &= g \\
	H(x_1) &= 0 \\
	H(x_2) &= c \\
	H(x_3) &= f \\
	H(x_4) &= 0.
\end{align*}
Moreover, since $[x_0]$ is $\bF[U]$-non-torsion  in $H_*(C)$ and $[a+Uc+e]$ is $\bF[U]$-non-torsion in $H_*(A^-_0)$, it follows that $\phi$ induces an isomorphism on $U^{-1}H_*$.

Hence $\psi$ and $\phi$ provide the desired local equivalence between $(A^-_0, \iota_K)$ and $(C, \iota_C)$.
(In fact, $\psi \phi = \id_{C_3}$ and $\phi \psi \simeq \id_{A^-_0}$ via the homotopy $J$ defined by $J(d) = c, J(h) = f$, and $J(y) = 0$ if $y \neq d, h$.)
\end{proof}

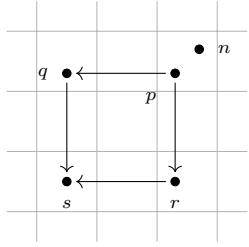
\begin{figure}
\[
\begin{tikzpicture}[scale=0.8]
	\draw[step=1, black!30!white, very thin] (-0.5, -0.5) grid (3.5, 3.5);

	\filldraw (2.7, 2.7) circle (2pt) node[label=right:{\lab{n}}] (a) {};
	\filldraw (2.3, 2.3) circle (2pt) node[label=below left:{\lab{p}}] (b) {};
	\filldraw (0.5, 2.3) circle (2pt) node[label=left:{\lab{q}}] (c) {};
	\filldraw (2.3, 0.5) circle (2pt) node[label=below:{\lab{r}}] (d) {};
	\filldraw (0.5, 0.5) circle (2pt) node[label=below:{\lab{s}}] (e) {};
	
	\draw[->, bend left=0] (b) to node[above,sloped]{\lab{}} (d);	
	\draw[->] (b) to (c);	
	\draw[->] (c) to (e);	
	\draw[->] (d) to (e);		
	
\end{tikzpicture}
\]
\caption{The complex $C_7$.}
\label{fig:C7}
\end{figure}

\begin{rem}
Recall from \cite[Remark 9.7]{HLInvolutiveBordered} that a homology sphere $Y$ is called \emph{$\HFIhat$-trivial} if $\HFIhat (Y) \cong \HFhat(Y) \oplus \bF$ where $Q \cdot \HFhat(Y) =0$ and $Q$ is non-vanishing on the remaining generator.
Let $J = -2T_{4,5} \# T_{4,9}$. Similar computations to those in Sections \ref{sec:tensorprods} and \ref{sec:iotaKexample} can show that $Y=S^3_{+1}(J)$ is not $\HFIhat$-trivial. Indeed, we have that $(\CFKi(J), \iota_J)$ decomposes $\iota_J$-equivariantly as $C_7 \oplus C_8$ where $C_7$ the complex in Figure \ref{fig:C7} where the Maslov gradings of the generators are
\[ M(n) = M(p) = 0, \quad M(q) = M(r) = -1, \quad M(s) = -2, \]
and $\iota_{7}$ is given by
\begin{align*}
	\iota_{7} (n) &= n \\
	\iota_{7} (p) &= p+n \\
	\iota_{7} (q) &= r \\
	\iota_{7} (r) &= q \\
	\iota_{7} (s) &= s \\
\end{align*}
and $C_8$ consists of an even number of rectangles which are paired by $\iota_{8}$ much like the (non-square) rectangles in Figure \ref{fig:C1xC2} are. (Compare this to the computation of $-2T_{6,7} \# T_{6, 13}$, which is locally equivalent to $C_1$ in Figure \ref{fig:C1}.) An application of the involutive mapping cone then shows that $\HFIhat(S^3_{+1}(J)) \cong \HFhat(Y) \oplus \bF^3$ and where $Q \cdot \HFhat(Y) =0$ and $Q$ is non-vanishing on the remaining generators.\end{rem}

\subsection{Almost local equivalence classes and \texorpdfstring{$S^3_{+1}(-2T_{6,7} \# T_{6,13} \# -T_{2,3;2,5})$}{+1 surgery on -2T(6,7) \# T(6,13) \# -T(2,3;2,5)}}\label{sec:alec}
The goal of this section is to prove Theorem \ref{thm:nYThetaSF}. We begin by recalling some some terminology and notation.

Recall from \cite[Section 4.1]{DHSThomcob} that to a sequence $(a_1, b_2, a_3, b_4, \dots, a_{2m-1}, b_{2m})$, where $a_i \in \{ \pm\}$ and $b_i \in \Z \setminus \{0\}$, we may associate an almost $\iota$-complex, called the \emph{standard complex of type $(a_1, b_2, a_3, b_4, \dots, a_{2m-1}, b_{2m})$}, as follows. (We will abuse notation and let $(a_i, b_i)$ denote both a sequence and the associated standard complex.) The standard complex is generated over $\bF[U]$ by $t_0, t_1, \dots, t_{2m}$. For each symbol $a_i$, we introduce an $\omega$-relation between $t_{i-1}$ and $t_i$ as follows:
\begin{itemize}
	\item If $a_i = +$, then $\omega t_{i} = t_{i-1}$.
	\item If $a_i = -$, then $\omega t_{i-1} = t_i$. 
\end{itemize}
For each symbol $b_i$, 
we introduce an $\d$-relation between $t_{i-1}$ and $t_i$ as follows:
\begin{itemize}
	\item If $b_i>0$, then $\d t_{i} = U^{|b_i|} t_{i-1}$.
	\item If $b_i<0$, then $\d t_{i-1} = U^{|b_i|} t_{i}$.
\end{itemize}
We begin with a few lemmas about tensor products of standard complexes.

\begin{lem}\label{lem:Sn}
Let $k \geq 2$ and let 
\[ S_n = (+, -1, +, -k, +, -1, +, -k, \dots, +, -1, +, -k) \]
where the right-hand side consists of $+, -1, +, -k$ repeated $n$ times. There is an almost local map from $S_n$ to $S_1 \otimes S_{n-1}$.
\end{lem}

\begin{proof}
Let the complex $S_n$ be generated by $x_0, x_1, \dots, x_{4n}$, the complex $S_1$ by $y_0, y_1, \dots, y_4$ and the complex $S_{n-1}$ by $z_0, z_1, \dots, z_{4n-4}$. Define $ f \co S_n \to S_1 \otimes S_{n-1}$ to be
\begin{align*}
	f(x_i) = \begin{cases}
		y_i \otimes z_0 \qquad \textup{ if } 0 \leq i \leq 4 \\
		y_4 \otimes z_{i-4} \quad \textup{ if } 5 \leq i \leq 4n.
		\end{cases}
\end{align*}
It is straightforward to verify that $f$ is an almost local map.
\end{proof}

\begin{lem}\label{lem:-Sn}
Let $k \geq 2$ and let 
\[ -S_n = (-, 1, -, k, -, 1, -, k, \dots, -, 1, -, k) \]
where the right-hand side consists of $-, 1, -, k$ repeated $n$ times. There is an almost local map from $-S_n$ to $-S_1 \otimes -S_{n-1}$.
\end{lem}

\begin{proof}
Let the complex $-S_{n}$ be generated by $x_0, x_1, \dots, x_{4n}$, the complex $-S_{1}$ by $y_0, y_1, \dots, y_4$ and the complex $-S_{n-1}$ by $z_0, z_1, \dots, z_{4n-4}$. Define $ f \co -S_{n} \to -S_{1} \otimes -S_{n-1}$ to be
\begin{align*}
	f(x_0) &= y_0 \otimes z_0 \\
	f(x_1) &= y_0 \otimes z_1 + y_1 \otimes z_0 + y_1\otimes z_1 \\	
	f(x_2) &= y_0 \otimes z_2 + y_2 \otimes z_0 + y_1 \otimes z_2 \\
	f(x_3) &= y_0 \otimes z_3 + y_1 \otimes z_2 + y_2 \otimes z_1 + y_3 \otimes z_0 + y_3 \otimes z_1 \\	
	f(x_4) &= y_0 \otimes z_4 + U^{k-1}y_2 \otimes z_2 + y_4 \otimes z_0 + y_4 \otimes z_1 \\	
	&\,\,\,\vdots\\
	f(x_{4i+1}) &= y_0 \otimes z_{4i+1} + y_1 \otimes z_{4i} + U^{k-1}y_2 \otimes z_{4i-1} + U^{k-1}y_3 \otimes z_{4i-2}\\
	&\quad  + y_4 \otimes z_{4i-3} + U^{k-1}y_3 \otimes z_{4i-1} + y_1 \otimes z_{4i+1} \\
	f(x_{4i+2}) &= y_0 \otimes z_{4i+2} + y_2 \otimes z_{4i} + y_4 \otimes z_{4i-2} + y_3 \otimes z_{4i} + y_1 \otimes z_{4i+2} \\	
	f(x_{4i+3}) &= y_0 \otimes z_{4i+3} + y_1 \otimes z_{4i+2} + y_2 \otimes z_{4i+1} + y_3 \otimes z_{4i} + y_4 \otimes z_{4i-1} \\
	f(x_{4i+4}) &= y_0 \otimes z_{4i+4} + U^{k-1}y_2 \otimes z_{4i+2} + y_4 \otimes z_{4i} \\	
	&\,\,\,\vdots\\
	f(x_{4n-3}) &= y_1 \otimes z_{4n-4} + U^{k-1}y_2 \otimes z_{4n-5} + U^{k-1}y_3 \otimes z_{4n-6} + y_4 \otimes z_{4n-7} + U^{k-1}y_3 \otimes z_{4n-5}  \\
	f(x_{4n-2}) &= y_2 \otimes z_{4n-4} + y_4 \otimes z_{4n-6}  + y_3 \otimes z_{4n-4} \\	
	f(x_{4n-1}) &= y_3 \otimes z_{4n-4} + y_4 \otimes z_{4n-5} \\
	f(x_{4n}) &= y_4 \otimes z_{4n-4}		
\end{align*}
where $i=1, \dots, n-2$. It is straightforward to verify that $f$ is an almost local map.
\end{proof}

\begin{lem}\label{lem:aiccomp}
Let $k \geq 2$. Then, in the notation of Lemmas \ref{lem:Sn} and \ref{lem:-Sn}, we have that $S_1^{\otimes n}$ is locally equivalent to $S_n$.
\end{lem}

\begin{proof}
By Lemma \ref{lem:Sn}, we have that $S_n \leq S_1 \otimes S_{n-1}$ and by Lemma \ref{lem:-Sn}, we have that $-S_n \leq -S_1 \otimes -S_{n-1}$, or equivalently, $S_n \geq S_1 \otimes S_{n-1}$. Since $\leq$ is a partial order on the set of almost local equivalence classes (see \cite[Definition 3.19]{DHSThomcob}), it follows that $S_n = S_1 \otimes S_{n-1}$, where equality denotes almost local equivalence. Hence by induction, $S_1^{\otimes n} = S_n$.
\end{proof}

\begin{rem}
In Lemmas \ref{lem:Sn} and \ref{lem:-Sn}, the complexes $S_n$ and $-S_n$ are in fact iota complexes (not just almost iota complexes), and the maps are in fact local maps (not just almost local maps). Hence $S_1^{\otimes n}$ and $S_n$ are in fact locally equivalent (not just almost locally equivalent).
\end{rem}

We are now ready to prove Theorem \ref{thm:nYThetaSF}.

\begin{proof}[Proof of Theorem \ref{thm:nYThetaSF}]
Let $Y = S^3_{+1}(-2T_{6,7} \# T_{6,13} \# -T_{2,3;2,5})$. By Proposition \ref{prop:knotexample}, we have that the almost local equivalence class of $(\CF^-(Y), \iota)$ is $(+, -1, +, -3)$. Then by Lemma \ref{lem:aiccomp}, the almost local equivalence class of $(\CF^-(nY), \iota)$ is $(+, -1, +, -3, +, -1, +, -3, \dots, +, -1, +, -3)$, that is, the sequence $(+, -1, +, -3)$ repeated $n$ times.

Suppose $(a_i, b_i)$ represents the almost local equivalence class of a linear combination of Seifert fibered spaces.
By \cite[Theorem 8.1]{DHSThomcob} (cf. the proof of \cite[Theorem 1.1]{DaiConnected}), if $a_i = +$ for all $i$, then $b_i \geq b_{i+1}$ for all $i$. Hence $nY$, for $n \neq 0$, is not homology cobordant to any linear combination of Seifert fibered spaces.
\end{proof}

\section{Hypercubes and hyperboxes}
\label{sec:hypercubes}

In this section, we present some background on homological algebra.

\subsection{Hyperboxes of chain complexes}
\label{sec:hypercubes-chain-complexes}

We recall the notion of a hyperbox of chain complexes, due to Manolescu and Ozsv\'{a}th \cite{MOIntegerSurgery}. We let $\bE_n$ denote the unit cube in $\R^n$:
\[
\bE_n:=\{0,1\}^n.
\]
More generally, if $\ve{d}=(d_1,\dots, d_n)\in \Z_{\ge 0}^n$, we define the cube of size $\ve{d}$ to be
\[
\bE(\ve{d})=\prod_{i=1}^n \{0,\dots, d_i\}.
\]
If $\varepsilon=(\varepsilon_1,\dots, \varepsilon_n)\in \bE(\ve{d})$, we write
\[
\|\varepsilon\|=\varepsilon_1+\cdots +\varepsilon_n.
\]

\begin{define}\label{def:hyperbox-chain-complexes} An $n$-dimensional \emph{hyperbox of chain complexes} of size $\ve{d}=(d_1,\dots, d_n)\in \Z_{\ge 0}^{n}$ consists of a collection of $\Z$-graded $\bF$-vector spaces
\[
(C^{\varepsilon})_{\varepsilon \in \bE(\ve{d})},
\]
together with a collection of maps 
\[
D^\varepsilon_{\varepsilon_0}\colon C^{\varepsilon_0}_*\to C^{\varepsilon_0+\varepsilon}_{*-1+\|\varepsilon\|},
\]
indexed over pairs $(\varepsilon_0,\varepsilon)\in \bE(\ve{d})\times \bE_n$ such that  $\varepsilon_0+\varepsilon\in \bE(\ve{d})$. Furthermore, if $\varepsilon\in \bE_n$ and $\varepsilon_0\in \bE(\ve{d})$, then
\begin{equation}
\sum_{\substack{\varepsilon'\in \bE_n \\ \varepsilon'\le \varepsilon} } D_{\varepsilon_0+\varepsilon'}^{\varepsilon-\varepsilon'} \circ D_{\varepsilon_0}^{\varepsilon'}=0, \label{eq:D^2=0}
\end{equation}
whenever $\varepsilon_0$ and $\varepsilon_0+\varepsilon$ are both in $\bE(\ve{d})$.
\end{define}

If $H=(C^{\varepsilon}, D_{\veps_0}^{\veps})$ is a hyperbox of chain complexes of size $\ve{d}$, we  write $C$ for the total space
\[
C:=\bigoplus_{\varepsilon\in \bE(\ve{d})} C^\varepsilon.
\]
Similarly, the maps $D_{\veps_0}^{\veps}$ may be summed to give an endomorphism
\[
D\colon C\to C.
\]
Equation~\eqref{eq:D^2=0} implies $D^2=0$. Similarly, if $(C,D)$ is a hyperbox, and $\veps\in \bE_n$, we write $D^{\veps}$ for the direct sum of $D_{\veps_0}^{\veps}$ over all $\veps_0$ where the map is defined.

A hyperbox of chain complexes of size $\ve{d}=(1,\dots, 1)\in \Z_{\ge 0}^n$ is called an $n$-dimensional  \emph{hypercube} of chain complexes.

For this paper, we only consider hyperboxes of dimension $n\le 3$.

\subsection{Stacking and compressing hyperboxes}

In this section, we describe two algebraic operations, called \emph{stacking} and \emph{compressing}.

We begin with the operation of stacking. A hyperbox of size $(1,1,2)$ is the same as two hyperboxes of size $(1,1,1)$, which have chain complexes and maps which coincide along one face. More generally, suppose we have a hyperbox $ H_1$ of size $(x,y,d_1)$ and a hyperbox $ H_2$ of size $(x,y,d_2)$. Suppose further that the restriction of $H_1$ to $\bE(x,y)\times \{d_1\}$ coincides with the restriction of $ H_2$ to $\bE(x,y)\times \{0\}$. We form a new hyperbox $\St( H_1, H_2)$ , which is a hyperbox of size $(x,y,d_1+d_2)$. A similar description holds for stacking two hypercubes of the same dimension, which coincide along a hyperface.

We now describe \emph{compression}, which is described by Manolescu and Ozsv\'{a}th \cite{MOIntegerSurgery}*{Section~5}. The operation takes an $n$-dimensional hyperbox $H=(C,D)$ of size $\ve{d}$, and yields an $n$-dimensional hypercube $\widehat{H}=(\widehat{C},\widehat{D})$. The compression operation is fundamental to our construction of the main hypercube from Theorem~\ref{thm:main-hypercube}. To simplify the exposition, we focus on the compression of hyperboxes size $(d)$, $(1,d)$ or $(1,1,d)$, which will be sufficient for our purposes.

We begin with 1-dimensional hyperboxes of size $(d)$.  A hyperbox $H=(C,D)$ of size $(d)$ is a collection of $d+1$ chain complexes $C^0,\dots, C^d$, together with chain maps
\[
D_i^{1}\colon C^i\to C^{i+1}.
\]
We depict such a hyperbox as
\[
\begin{tikzcd}
C^0 \arrow[r,"D^1"] & C^1\arrow[r, "D^1"]&\cdots\arrow[r,"D^1"]&C^d. 
\end{tikzcd}
\]
The compression $\widehat{H}=(\widehat{C}, \widehat{D})$ has underlying chain complex equal to $C^0\oplus C^d$.  The length 0 map $\widehat{D}^0$ is the restriction of the original length 0 map $D^0$. The length 1 map $\widehat{D}^1$ is the $d$-fold composition of the original length 1 map, $D^1$. We depict $\widehat{H}$ as
\[
\begin{tikzcd}
C^0 \arrow[r,"(D^1)^{\circ d}"] &C^d.
\end{tikzcd}
\]

Next, we move on to 2 dimensional hyperboxes of size $(1,d)$. A 2-dimensional hyperbox of size $(1,d)$ consists of $2(d+1)$ chain complexes $(C^{\varepsilon})_{\varepsilon \in \{0,1\}\times \{0,\dots, d\}}$, as well as chain maps $D^{01}$ and $D^{10}$, and a chain homotopy $D^{11}$ between $D^{10}\circ D^{01}$ and $D^{01}\circ D^{10}$. Schematically, we depict $(C,D)$ as the diagram
\[
\begin{tikzcd}[row sep=huge, column sep=huge]
 C^{00}\arrow[dr,dashed,"D^{11}"]\arrow[d,"D^{01}"]\arrow[r, "D^{10}"]& C^{10}\arrow[d,"D^{01}"]\\
C^{01}\arrow[r, "D^{10}"]\arrow[d,"D^{01}"]\arrow[dr,dashed,"D^{11}"]& C^{11}\arrow[d,"D^{01}"]\\
\vdots\arrow[d,"D^{01}"]\arrow[dr,dashed,"D^{11}"]& \vdots\arrow[d,"D^{01}"]\\
C^{0d}\arrow[r, "D^{10}"]& C^{1d}.
\end{tikzcd}
\]
The compression $\widehat{H}=(\widehat{C},\widehat{D})$ is the 2-dimensional hypercube depicted below:
\[
\begin{tikzcd}[row sep=huge, column sep=huge]
C^{00}\arrow[r, "\widehat{D}^{10}"]\arrow[d, " \widehat{D}^{01}"]\arrow[dr,dashed, "\widehat{D}^{11}"]& C^{10}\arrow[d,"\widehat{D}^{01}"]\\
C^{0d}\arrow[r,"\widehat{D}^{10}"]& C^{1d}.
\end{tikzcd}
\]
Here
\begin{equation}
\begin{split} \widehat{D}^{10}&=D^{10}\\
\widehat{D}^{01}&=(D^{01})^{\circ d}\\
\widehat{D}^{11}&=\sum_{\substack{i+j=d-1\\ i,j\ge 0}} (D^{01})^{\circ j} \circ D^{11} \circ (D^{01})^{\circ i}.
\end{split}\label{eq:compressed-2-d}
\end{equation}
The length 2 compressed map $\widehat{D}^{11}$ is shown schematically in Figure~\ref{fig:80}.

\begin{figure}[ht!]
\[
\begin{tikzpicture}[scale=1.2]
	\draw[step=1, black!30!white] (0, 0) grid (1,5);
	\draw (-1.2,2.5) node[auto] {$\widehat{D}^{11}=\displaystyle{\sum}$};
	\draw[->] (0,4.9) -- (0,4.1) node[midway,fill=white, fill opacity=.8, text opacity=1] {$D^{01}$};
	\draw[->,dashed]  (.1,3.9) -- (.9,3.1) node[midway,fill=white, fill opacity=.8, text opacity=1] {$D^{11}$};
	\draw[->] (1,2.9) -- (1,2.1) node[midway,fill=white, fill opacity=.8, text opacity=1] {$D^{01}$};
	\draw[->] (1,1.9) --(1, 1.1) node[midway,fill=white, fill opacity=.8, text opacity=1] {$D^{01}$};
	\draw[->] (1,.9) -- (1, .1) node[midway,fill=white, fill opacity=.8, text opacity=1] {$D^{01}$};
\end{tikzpicture}
\]
	\caption{The length 2 map $\widehat{D}^{11}$ of the compression of a 2-dimensional hyperbox of size $(1,d)$.}\label{fig:80}
\end{figure}

We now discuss 3-dimensional hyperboxes and their compressions. Suppose that $H=(C,D)$ is a 3-dimensional hyperbox of size $(1,1,d)$. The underlying chain complexes of the compressed $\widehat{C}$ are 
\[
\widehat{C}^{(\veps_1,\veps_2,\veps_3)}:=C^{(\veps_1,\veps_2,\veps_3d)}.
\]
The length 1 maps, $\widehat{D}^{100},$ $\widehat{D}^{010}$ and $\widehat{D}^{001}$, as well as the length 2 maps, $\widehat{D}^{110}$, $\widehat{D}^{101}$ and $\widehat{D}^{011}$, are essentially the same as in the 2-dimensional case. See ~\eqref{eq:compressed-2-d}. The length 3 map, $\widehat{D}^{111}$, is given by the formula
\begin{equation}
\begin{split}
\widehat{D}^{111}&:=\sum_{\substack{i+j=d-1\\ i,j\ge 0}} (D^{001})^{\circ j} \circ D^{111} \circ (D^{001})^{\circ i}\\ 
&+\sum_{\substack{i+j+k=d-2\\ i,j,k\ge 0}} (D^{001})^{\circ i}\circ D^{101}\circ (D^{001})^{\circ j} \circ D^{011}\circ (D^{001})^{\circ k}\\
&+\sum_{\substack{i+j+k=d-2\\ i,j,k\ge 0}} (D^{001})^{\circ i}\circ D^{011}\circ (D^{001})^{\circ j} \circ D^{101}\circ (D^{001})^{\circ k}.
\end{split}
\label{eq:cc111}
\end{equation}

Equation~\eqref{eq:cc111} is schematically illustrated in Figure~\ref{fig:79}.

\begin{figure}[ht!]
\[
\begin{tikzpicture}[scale=1.2]
	\draw[-,black!30!white] (-3,2) -- (-2,2);
	\draw[-,black!30!white] (-3,3) -- (-2,3);
	\draw[-,black!30!white] (-2.5,1.7) -- (-1.5,1.7);
	\draw[-,black!30!white] (-2.5,2.7) -- (-1.5,2.7);
	\draw[-,black!30!white] (-3,2) -- (-2.5,1.7);
	\draw[-,black!30!white] (-2,2) -- (-1.5,1.7);
	\draw[-,black!30!white] (-3,3) -- (-2.5,2.7);
	\draw[-,black!30!white] (-2,3) -- (-1.5,2.7);
	\draw[-,black!30!white] (-3,3) -- (-3,2);
	\draw[-,black!30!white] (-2,3) -- (-2,2);
	\draw[-,black!30!white] (-2.5,2.7) -- (-2.5,1.7);
	\draw[-,black!30!white] (-1.5,2.7) -- (-1.5,1.7);
	\draw[->, dotted,thick] (-2.9,2.9) -- (-1.6,1.8) 
	node[midway,fill=white, fill opacity=.8, text opacity=1] {$\widehat{D}^{111}$};
	\draw (-.8,2.5) node[auto] {$=\displaystyle{\sum}$};
	\draw[step=1, black!30!white] (0, 0) grid (1,5);
	\draw[step=1, black!30!white,shift={(.5,-.3)}]  (0, 0) grid (1,5);
	\draw[-,black!30!white] (0,0) -- (.5,-.3);
	\draw[-,black!30!white] (0,1) -- (.5,.7);
	\draw[-,black!30!white] (0,2) -- (.5,1.7);
	\draw[-,black!30!white] (0,3) -- (.5,2.7);
	\draw[-,black!30!white] (0,4) -- (.5,3.7);
	\draw[-,black!30!white] (0,5) -- (.5,4.7);
	\draw[-,black!30!white] (1,0) -- (1.5,-.3);
	\draw[-,black!30!white] (1,1) -- (1.5,.7);
	\draw[-,black!30!white] (1,2) -- (1.5,1.7);
	\draw[-,black!30!white] (1,3) -- (1.5,2.7);
	\draw[-,black!30!white] (1,4) -- (1.5,3.7);
	\draw[-,black!30!white] (1,5) -- (1.5,4.7);
	\draw[step=1, black!30!white,shift={(3,0)}] (0, 0) grid (1,5);
	\draw[step=1, black!30!white,shift={(3.5,-.3)}]  (0, 0) grid (1,5);
	\draw[-,black!30!white,shift={(3,0)}] (0,0) -- (.5,-.3);
	\draw[-,black!30!white,shift={(3,0)}] (0,1) -- (.5,.7);
	\draw[-,black!30!white,shift={(3,0)}] (0,2) -- (.5,1.7);
	\draw[-,black!30!white,shift={(3,0)}] (0,3) -- (.5,2.7);
	\draw[-,black!30!white,shift={(3,0)}] (0,4) -- (.5,3.7);
	\draw[-,black!30!white,shift={(3,0)}] (0,5) -- (.5,4.7);
	\draw[-,black!30!white,shift={(3,0)}] (1,0) -- (1.5,-.3);
	\draw[-,black!30!white,shift={(3,0)}] (1,1) -- (1.5,.7);
	\draw[-,black!30!white,shift={(3,0)}] (1,2) -- (1.5,1.7);
	\draw[-,black!30!white,shift={(3,0)}] (1,3) -- (1.5,2.7);
	\draw[-,black!30!white,shift={(3,0)}] (1,4) -- (1.5,3.7);
	\draw[-,black!30!white,shift={(3,0)}] (1,5) -- (1.5,4.7);
		\draw[step=1, black!30!white,shift={(6,0)}] (0, 0) grid (1,5);
	\draw[step=1, black!30!white,shift={(6.5,-.3)}]  (0, 0) grid (1,5);
	\draw[-,black!30!white,shift={(6,0)}] (0,0) -- (.5,-.3);
	\draw[-,black!30!white,shift={(6,0)}] (0,1) -- (.5,.7);
	\draw[-,black!30!white,shift={(6,0)}] (0,2) -- (.5,1.7);
	\draw[-,black!30!white,shift={(6,0)}] (0,3) -- (.5,2.7);
	\draw[-,black!30!white,shift={(6,0)}] (0,4) -- (.5,3.7);
	\draw[-,black!30!white,shift={(6,0)}] (0,5) -- (.5,4.7);
	\draw[-,black!30!white,shift={(6,0)}] (1,0) -- (1.5,-.3);
	\draw[-,black!30!white,shift={(6,0)}] (1,1) -- (1.5,.7);
	\draw[-,black!30!white,shift={(6,0)}] (1,2) -- (1.5,1.7);
	\draw[-,black!30!white,shift={(6,0)}] (1,3) -- (1.5,2.7);
	\draw[-,black!30!white,shift={(6,0)}] (1,4) -- (1.5,3.7);
	\draw[-,black!30!white,shift={(6,0)}] (1,5) -- (1.5,4.7);
	\draw[->] (0,4.9) -- (0,4.1)
	node[midway,fill=white, fill opacity=.8, text opacity=1] {$D^{001}$};
	\draw[->,dotted, thick] (.1,3.9) -- (1.4,2.8)
	node[midway,fill=white, fill opacity=.8, text opacity=1] {$D^{111}$};
	\draw[->] (1.5, 2.6) -- (1.5,1.8)
	node[midway,fill=white, fill opacity=.8, text opacity=1] {$D^{001}$};
	\draw[->] (1.5, 1.6) -- (1.5,.8)
	node[midway,fill=white, fill opacity=.8, text opacity=1] {$D^{001}$};
	\draw[->] (1.5, .6) -- (1.5, -.2)
	node[midway,fill=white, fill opacity=.8, text opacity=1] {$D^{001}$};
	\draw[->] (3,4.9) -- (3,4.1)
	node[midway,fill=white, fill opacity=.8, text opacity=1] {$D^{001}$};
	\draw[->,dashed] (3.1, 3.9)--(3.4, 2.8)
	node[midway,fill=white, fill opacity=.8, text opacity=1] {$D^{011}$};
	\draw[->] (3.5, 2.6)--(3.5, 1.8)
	node[midway,fill=white, fill opacity=.8, text opacity=1] {$D^{001}$};
	\draw[->,dashed] (3.6, 1.6)--(4.4, .8)
	node[midway,fill=white, fill opacity=.8, text opacity=1] {$D^{101}$};
	\draw[->] (4.5, .6) -- (4.5, -.2)
	node[midway,fill=white, fill opacity=.8, text opacity=1] {$D^{001}$};
	\draw[->] (6,4.9) -- (6,4.1)
	node[midway,fill=white, fill opacity=.8, text opacity=1] {$D^{001}$};
	\draw[->, dashed] (6.1, 3.9)--(6.9, 3.1)
	node[midway,fill=white, fill opacity=.8, text opacity=1] {$D^{101}$};
	\draw[->, dashed] (7.1, 2.9)--(7.4, 1.8)
	node[midway,fill=white, fill opacity=.8, text opacity=1] {$D^{011}$};
	\draw[->] (7.5, 1.6) -- (7.5, .8)
	node[midway,fill=white, fill opacity=.8, text opacity=1] {$D^{001}$};
	\draw[->] (7.5, .6) -- (7.5, -.2)
	node[midway,fill=white, fill opacity=.8, text opacity=1] {$D^{001}$};
	\draw (2.3,2.3) node[auto] {$+$};
	\draw (5.3,2.3) node[auto] {$+$};
\end{tikzpicture}
\]
	\caption{The length 3 map $\widehat{D}^{111}$ of the compression of a 3-dimensional hyperbox of size $(1,1,d)$.}\label{fig:79}
\end{figure}

\subsection{Mapping cones and hyperboxes}

In this section, we describe the relationship between the hyperbox construction of Manolescu and Ozsv\'{a}th and the more common mapping cone construction. The results of this section will be helpful for deriving the involutive mapping cone formula from the hypercube in Theorem~\ref{thm:main-hypercube}.

\begin{example}Suppose $(C_0,D_0)$ and $(C_1,D_1)$ are chain complexes, and $F_{0}^1\colon C_0\to C_1$ is a chain map. We may form the mapping cone complex
\[
\Cone(F_0^1):=C_0[-1]\oplus C_1,
\]
which has differential $D_0+ D_1+F$. The mapping cone complex is naturally filtered over the set $\{0,1\}$. The mapping cone complex is a 1-dimensional hypercube of chain complexes, of size $(1)$. If $F_0^1\colon C_0\to C_1$ is a chain map, then there is an exact triangle
\[
\cdots \to  H_*(C_0)\to H_*(C_1)\to H_*(\Cone(F_0^1)) \to H_{*-1}(C_0)\to \cdots.
\]
\end{example}

\begin{example}Consider a 2-dimensional hypercube $H=(C,D)$ satisfying $C^{10}=0$. Such a hypercube has the following presentation:
\[
\begin{tikzcd}[row sep=1cm, column sep=1cm]
C^{00}\arrow[d,swap, "D^{01}"]\arrow[dr, dashed, "D^{11}"]& \\
 C^{01} \arrow[r, swap, "D^{10}"]&C^{11}.
\end{tikzcd}
\]
The data of $H$ is equivalent to three chain complexes, two chain maps $D^{01}$ and $D^{10}$, and a null-homotopy $D^{11}$ of the composition $D^{10}\circ D^{01}$. Equivalently, the hypercube is the same as a chain map from $C^{00}$ to $\Cone(D^{10}\colon C^{01}\to C^{11})$.

More generally, suppose $H=(C,D)$ is an arbitrary 2-dimensional hypercube of chain complexes:
\[
\begin{tikzcd}[row sep=huge, column sep=huge]
C^{00}\arrow[d, "D^{01}",swap]\arrow[r, "D^{10}"]\arrow[dr,dashed, "D^{11}"]&C^{10}\arrow[d, "D^{01}"]\\
C^{01}\arrow[r, "D^{10}"]& C^{11}
\end{tikzcd}.
\]
The data of the hypercube $H$ is equivalent to a chain map from $\Cone(D_{00}^{01}\colon C^{00}\to C^{01})$ to $\Cone (D^{01}_{10}\colon C^{10}\to C^{11})$ which is filtered with respect to the mapping cone filtrations. 
\end{example}

\begin{example}
Suppose that $H=(C,D)$ is a 3-dimensional hypercube, which has the property that $C^{010}$ and $C^{011}$ are both zero. See ~\eqref{eq:triple-mapping-cone}.
\begin{equation}
\begin{tikzcd}[row sep=.7cm, column sep=2cm]
C^{000}
	\arrow[ddddrr,dotted]
	\arrow[dr]
	\arrow[ddd]
	\arrow[drr,dashed]&[-1cm] \,&\,
\\[-.3cm]
& C^{100}
	\arrow[r]
	\arrow[dddr,dashed]
&C^{110}\arrow[ddd]\\
\\
C^{001}
	\arrow[dr]
	\arrow[rdr,dashed]&&\,
\\[-.3cm]
& C^{101}
	\arrow[from=uuuul, dashed, crossing over]
	\arrow[from=uuu,crossing over]
	\arrow[r]&
C^{111}
\end{tikzcd}
\label{eq:triple-mapping-cone}
\end{equation}

We  can identify the subcomplex generated by $C^{100}$, $C^{110}$, $C^{101}$ and $C^{111}$ as the cone of a filtered chain map 
\[
G\colon \Cone \left(\begin{tikzcd} C^{100} \arrow[d]\\
C^{101} \end{tikzcd}\right)\to \Cone\left(\begin{tikzcd} C^{110}\arrow[d] \\ C^{111}\end{tikzcd} \right).
\]

The 3-dimensional hypercube in ~\eqref{eq:triple-mapping-cone} may be interpreted as encoding the data of a  chain map
\[
\Phi\colon \Cone\left(\begin{tikzcd}C^{000}\arrow[d]\\
C^{001}\end{tikzcd}\right)\to \Cone(G),
\]
which is filtered with respect the last component of the indexing (the vertical direction in our diagrams).
\end{example}

\subsection{Hypercubes and involutive complexes}

In this section, we describe some helpful algebra for relating the hypercube construction to the algebraic setting of involutive Heegaard Floer homology. Suppose $\cR$ is a ring such that $2=0$, and that we have a hypercube of chain complexes over $\cR$ of following form:
\begin{equation}
\begin{tikzcd}[row sep=huge, column sep=huge]
X\arrow[d, "g",swap]\arrow[r, "f"]\arrow[dr,dashed, "h"]&Y\arrow[d, "k"]\\
X\arrow[r, "f"]& Y
\end{tikzcd}.
\label{eq:involutive-cube}
\end{equation}
In particular, the complexes and morphisms along the top and bottom coincide. We can view $\Cone(g)$ and $\Cone(k)$ as complexes over $\cR[Q]/Q^2$, where $Q$ takes the first copy of $X$ to the second, and similarly for $Y$.

 The diagram in ~\eqref{eq:involutive-cube}
determines a map $f+Q\cdot h$ from $\Cone(g)$ to $\Cone(k)$ which is $Q$-equivariant. We collect the following basic algebraic facts:
\begin{lem}
\label{lem:homotopy-equivalence-iota-complexes} Suppose that $X$, $Y$, $f$, $g$, $k$ and $h$ form a hypercube, as in ~\eqref{eq:involutive-cube}. Suppose that $f\colon X\to Y$ is a homotopy equivalence.
\begin{enumerate}
\item The induced map $f+Q\cdot h$ from $\Cone(g)$ to $\Cone(k)$ is also a homotopy equivalence over $\cR[Q]/Q^2$.
\item If $g^2\simeq \id$, then $k^2\simeq \id$. Similarly, if $g^2\simeq 0$, then $k^2\simeq 0$.
\end{enumerate}
\end{lem}
\begin{proof} Let $j\colon  Y\to X$ be a homotopy inverse to $f$. We consider the first claim. First, we claim that $jk\simeq gj$. To see this, we begin with the relation $kf\simeq fg$, and pre- and post-compose with $j$ to obtain $jkfj\simeq jfgj$. Using the fact that $j$ and $f$ are homotopy inverses, we obtain that $jk\simeq gj$. If $s$ is a map so that $jk+gj=[\d,s]$, then $j+Q\cdot s$ is a chain map from $\Cone(k)$ to $\Cone(g)$. Since $jf\simeq  \id$, the composition $(j+Q\cdot s)\circ (f+Q\cdot h)$ is chain homotopic to a map of the form $\Psi=(\id+t\cdot Q)$, for some $t\colon X\to X$. The map $\Psi$ is a chain isomorphism, with inverse equal to itself. Hence $\Psi\circ (j+Q\cdot s)$ is a left homotopy inverse of $f+Q\cdot h$. A similar argument shows that $(f+Q\cdot h)$ has a right homotopy inverse. An easy argument shows that if $(f+Q\cdot h)$ has both a left and a right homotopy inverse, then they must be homotopic.

To prove the second claim, suppose first that $g^2\simeq \id$. We perform the following computation:
\[
\begin{split}
\id&\simeq f j\\
&\simeq (fg)(gj)\\
&\simeq  (kf)(jk)\\
&\simeq k^2.
\end{split}
\]
A similar argument shows that if $g^2\simeq 0$, then $k^2\simeq 0$.
\end{proof}

\subsection{Hyperboxes of attaching curves}

Manolescu and Ozsv\'{a}th also describe a natural analog of a hyperbox for diagrams in the Fukaya category. The main structural difference is that in the Fukaya category, we must include higher compositions into the hypercube relations. We make the following defintion, focusing on the case relevant to Heegaard Floer homology:

\begin{define}\label{def:hyperbox-attaching-curves}
 An $n$-dimensional \emph{hypercube of beta attaching curves} on $(\Sigma,w)$ is a collection of attaching curves $\bs^{\varepsilon}$, ranging over $\varepsilon\in \bE_n$, together with a distinguished element
 \[
 \Theta_{\b^\veps,\b^{\veps'}}\in \bCF^-(\bs^{\veps}, \bs^{\veps'}),
 \]
 whenever $\veps<\veps'$.
Furthermore, whenever $\veps< \veps'$, the following compatibility condition is satisfied: 
 \begin{equation}
0=\sum_{\veps< \veps_1< \cdots<\veps_k < \veps'}f_{\b^\veps, \b^{\veps_1},\dots,\b^{\veps_k},\b^{\veps'}}(\Theta_{\b^{\veps},\b^{\veps_1}},\dots, \Theta_{\b^{\veps_k},\b^{\veps'}}), \label{eq:compatibility-hypercube-attaching-curves}
 \end{equation}
 where $f_{\b^\veps, \b^{\veps_1},\dots,\b^{\veps_k},\b^{\veps'}}$ is the map which counts $2+k$-gons. An $n$-dimension \emph{hyperbox} of beta attaching curves is defined by the obvious adaptation.
 
 An $n$-dimensional \emph{hypercube of alpha attaching curves} is a collection of attaching curves $(\as^{\nu})_{\nu\in \bE_n}$, together with a distinguished chain $\Theta_{\a^{\nu'}, \a^\nu}\in \bCF^-(\as^{\nu'}, \as^{\nu})$ whenever $\nu<\nu'$. Furthermore, the classes $\Theta_{\a^{\nu'}, \a^{\nu}}$ satisfy the natural adaptation of ~\eqref{eq:compatibility-hypercube-attaching-curves}.
\end{define}

\begin{rem} Manolescu and Ozsv\'{a}th's definition of a hypercube of beta attaching curves \cite{MOIntegerSurgery}*{Section~8.1} assumes that for all $\veps,\veps'\in \bE_n$, the curves $\bs^{\veps}$ and $\bs^{\veps'}$ are  related by a sequence of handleslides and isotopies, whereas ours does not. The additional generality will be useful especially in Section~\ref{sec:central-hypercube} when we construct the central hypercube.
\end{rem}

\begin{example}
We consider a 2-dimensional hypercube of beta attaching curves, which has the following form:
\[
\begin{tikzcd}[row sep=1.5 cm, column sep=2.5 cm,labels=description]
\bs^{00}
	\arrow[r, "\Theta_{\b^{00},\b^{01}}"]
	\arrow[dr,dashed,sloped, "\Theta_{\b^{00},\b^{11}}"]
	\arrow[d,swap, "\Theta_{\b^{00},\b^{10}}"]
& \bs^{01}
	\arrow[d,"\Theta_{\b^{01},\b^{11}}"]\\
\bs^{10}
	\arrow[r, "\Theta_{\b^{10},\b^{11}}"]
& \bs^{11}
\end{tikzcd}.
\]
The compatibility relations are equivalent to the length 1 chains, $\Theta_{\b^{00},\b^{01}},$ $\Theta_{\b^{00},\b^{10}}$, $\Theta_{\b^{10},\b^{11}}$ and $\Theta_{\b^{01},\b^{11}}$, being cycles, and to the relation
\[
\d \Theta_{\b^{00},\b^{11}}+f_{\b^{00}, \b^{01}, \b^{11}}(\Theta_{\b^{00},\b^{01}}, \Theta_{\b^{01},\b^{11}})+f_{\b^{00}, \b^{10}, \b^{11}}(\Theta_{\b^{00},\b^{10}}, \Theta_{\b^{10},\b^{11}})=0.
\]
\end{example}

\begin{rem}
Given a collection of attaching curves $(\bs^\veps)_{\veps\in \bE_n}$, and a set of length 1 chains $\Theta_{\b^{\veps},\b^{\veps'}}$, ranging over $\veps,\veps'\in \bE_n$ such that $|\veps-\veps'|=1$, it is natural to ask when one may pick higher length chains  which satisfy the hypercube relations.  If, for all $\veps,\veps'\in \bE_n$, the attaching curves $\bs^{\veps}$ and $\bs^{\veps'}$ are related by a sequence of handleslides and isotopies, and length 1 cycles are chosen to represent the top degree generators of $\bHF^-(\Sigma,\bs^{\veps}, \bs^{\veps'},w)$, then a procedure of Manolescu and Ozsv\'{a}th \cite{MOIntegerSurgery}*{Lemma~8.6} proves the existence of such a collection of higher length chains which satisfy the hypercube relations. 
\end{rem}

\subsection{Pairing hypercubes of attaching curves}
\label{eq:pairing-hypercube}

Suppose $\cL_{\a}=(\as^{\nu})_{\nu \in \bE_m}$ is a hypercube of alpha attaching curves, and $\cL_{\b}=(\bs^{\veps})_{\veps\in \bE_n}$ is a hypercube of beta attaching curves. Furthermore, suppose that the Heegaard multi-diagram with all $2^{m}+2^n$ attaching curves has only transverse double point intersections, and is weakly admissible. In this situation, there is a naturally associated $n+m$-dimensional hypercube of chain complexes $\bCF^-(\cL_{\a},\cL_{\b})$, as follows.

The chain complex associated to the vertex $(\nu,\veps)\in \bE_{m+n}$ is $\bCF^-(\as^{\nu}, \bs^{\veps})$. To describe the maps, it is helpful to introduce some notation. If $\nu_1<\dots< \nu_i$ and $\veps_1<\cdots <\veps_j$, we define the map
\[
f_{\a^{\nu_1}\to \dots\to  \a^{\nu_i}}^{\b^{\veps_1}\to \cdots \to \b^{\veps_j}}\colon \bCF^-(\as^{\nu_1}, \bs^{\veps_1})\to \bCF^-(\as^{\nu_i}, \bs^{\veps_j})
\]
via the formula
\[
f_{\a^{\nu_1}\to \dots\to  \a^{\nu_i}}^{\b^{\veps_1}\to \cdots \to \b^{\veps_j}}(\xs)=f_{\a^{\nu_i},\dots, \a^{\nu_1},\b^{\veps_1},\dots, \b^{\veps_j}}\left(\Theta_{\a^{\nu_i},\a^{\nu_{i-1}}},\dots, \Theta_{\a^{\nu_2},\a^{\nu_1}}, \ve{x}, \Theta_{\b^{\veps_1},\b^{\veps_2}},\dots, \Theta_{\b^{\veps_{j-1}}, \b^{\veps_{j}}}\right),
\]
extended equivariantly over the action of $U$.
The hypercube map in $\bCF^-(\cL_{\a}, \cL_{\b})$ from $(\nu,\veps)$ to $(\nu',\veps')$ is given by the formula
\begin{equation}
F_{\a^{\nu}\to \a^{\nu'}}^{\b^{\veps}\to \b^{\veps'}} :=\sum_{\substack{\nu<\nu_1<\dots<\nu_i<\nu' \\ \veps<\veps_1<\dots< \veps_j<\veps'}}f_{\a^\nu\to\a^{\nu_1}\to \dots\to  \a^{\nu_i}\to \a^{\nu'}}^{\b^{\veps}\to \b^{\veps_1}\to \cdots \to \b^{\veps_j}\to \b^{\veps'}}. \label{eq:pairing-attaching-curves}
\end{equation}
In the notation of Section~\ref{sec:hypercubes-chain-complexes}, the map $F_{\a^{\nu}\to \a^{\nu'}}^{\b^{\veps}\to \b^{\veps'}}$ would be written $D_{(\nu,\veps)}^{(\nu'-\nu,\veps'-\veps)}$. If $\nu=\nu'$, then we will usually write $F_{\a^{\nu}}^{\b^{\veps}\to \b^{\veps'}}$, and similarly if $\veps=\veps'$.

\begin{lem}\label{lem:hypercube-diagram-to-complex}
 If $\cL_{\a}$ and $\cL_{\b}$ are $m$ and $n$-dimensional hypercubes of attaching curves on $(\Sigma,w)$, such that the Heegaard multi-diagram consisting of all $2^m+2^n$ attaching curves is weakly admissible, then the diagram $\bCF^-(\cL_{\a},\cL_{\b})$ is an $(n+m)$-dimensional hypercube of chain complexes.
\end{lem}

Lemma~\ref{lem:hypercube-diagram-to-complex} follows immediately from the associativity relations for holomorphic polygons, together with the compatibility conditions in~\eqref{eq:compatibility-hypercube-attaching-curves}.

\subsection{Hypercubes with twisted coefficients}
\label{sec:twisted-hypercubes-1}
The above notions naturally adapt to give a notion of a hypercube of attaching curves on a doubly pointed surface  $(\Sigma,w,z)$ with twisted coefficients. Assume first that there is a chosen path $\g$ from $w$ to $z$. A twisted hypercube of alpha attaching curves $\cL_{\a}$ is defined similarly to the untwisted case, except all of the Floer complexes have coefficients in $\bF\llsquare U\rrsquare \otimes \bF[\Z/m]$. A holomorphic curve representing a class $\phi$ is counted with a multiplicity of $U^{n_w(\phi)} T^{n_z(\phi)-n_w(\phi)}$. We may pair two twisted hypercubes $\cL_{a}$ and $\cL_{\b}$ to form a hypercube of chain complexes over $\bF\llsquare U\rrsquare \otimes\bF[\Z/m]$, denoted
\[
\ul{\buCF}^-(\cL_{\a},\cL_{\b}).
\]
We are usually interested in a subcomplex of $\ul{\buCF}^-(\cL_{\a},\cL_{\b})$, which we denote by $\buCF^-(\cL_{\a},\cL_{\b})$. Suppose that $\cL_{\a}$ and $\cL_{\b}$ are $n$ and $m$-dimensional, respectively. If $(\veps,\nu)\in \bE_n\times \bE_m$, we define $\buCF^-(\cL_{\a},\cL_{\b})^{(\veps,\nu)}$ to be $\bCF^-(\as^{\veps}, \bs^{\nu})\otimes T^0\iso \bCF^-(\as^{\veps}, \bs^{\nu})$ if
\[
\as^{\veps'}\cap \g=\bs^{\nu'}\cap \g=\emptyset \quad \text{whenever} \quad \veps'\le\veps \quad  \text{and} \quad \nu'\le \nu,
\]
and we define $\buCF^-(\cL_{\a},\cL_{\b})^{(\veps,\nu)}=\buCF^-(\as^{\veps}, \bs^{\veps})$ otherwise. It is straightforward to see that $\buCF^-(\cL_{\a}, \cL_{\b})$ is a subcomplex of $\ul{\buCF}^-(\cL_{\a},\cL_{\b})$.

\subsection{Low dimensional notation}
\label{sec:low-dim-notation}

Since we will only be dealing with hypercubes of dimension 3 or less, it is convenient to make some special notation.

\begin{enumerate}
\item We use the characters $\Theta,$ $\lambda$ and $\omega$ for length 1, 2 and 3 chains, respectively, in a hypercube of attaching curves. 
\item We use the characters $F$, $H$, and $P$, respectively, for the length 1, 2 and 3 maps in a hypercube of chain complexes.
\item We use the characters $f$, $h$ and $p$ to indicate triangle, quadrilateral, and pentagon counting maps, respectively. 
\end{enumerate}
As particular examples, if $\as'$ and $\as$ differ in a single hypercube coordinate, we write $f_{\a\to \a'}^{\b}$ for the map $f_{\a',\a,\b}(\Theta_{\a',\a},-)$. Similarly, if $\as'$ and $\as$ differ in two coordinates, we write $f_{\a\to \a'}^{\b}$ for the map $f_{\a',\a,\b}(\lambda_{\a',\a},-)$.

\section{Doubled Heegaard diagrams and the involution}
\label{sec:doubled-diagrams}

In this section, we begin by describing the operation of \emph{doubling} a Heegaard diagram of a 3-manifold in Section \ref{sec:doubled-3-manifolds}. There is an analog for diagrams of knots, which is described in Section~\ref{sec:doubled-knots}. In Section~\ref{sec:involution-doubling}, we will see how doubled Heegaard diagrams give a conceptually simple formula for the involution. In Section~\ref{sec:outline}, we sketch our construction of the main hypercube, in terms of doubled diagrams. We prove a technical result regarding admissibility of doubled diagrams in Section \ref{sec:admissibility}.

In this section, we focus on the case that $Y$ is an integer homology 3-sphere, which simplifies the discussion of the twisted complexes $\buCF^-(Y)$. This level of generality will be sufficient for our purposes until Section~\ref{sec:rational}, when we talk about rational surgeries.

\subsection{Doubled Heegaard diagrams}
\label{sec:doubled-3-manifolds}
Suppose that $(\Sigma,\as,\bs,w)$ is a Heegaard diagram for $(Y,w)$.
 Let $D\subset \Sigma$ be an open disk centered at $w$. Note that both $N(\Sigma\setminus D)$ and its complement are handlebodies. We construct a new Heegaard surface $D(\Sigma)$ for $Y$:
\[
D(\Sigma):=\d\left( (\Sigma\setminus D)\times [0,1]\right)\iso \Sigma\# \bar{\Sigma}.
\]

We now describe a method for constructing attaching curves on $D(\Sigma)$. We may take one set of attaching curves to be $\as\cup \bar{\bs}$, where $\as\subset \Sigma$ and $\bar{\bs}\subset \bar{\Sigma}$. (The barring of $\bs$ indicates only that it lives on $\bar{\Sigma}$). A second set of attaching curves may be constructed as follows. Let $I\subset \d D$ be a subinterval. Let $d_1,\dots, d_{2g}$ be a set of pairwise disjoint,  properly embedded arcs in $\Sigma\setminus D$, such that $\d d_i\subset I$, and which induce a basis of $H_1(\Sigma\setminus D, I;\Z)$. If we view $D(\Sigma)$ as $\Sigma\# \bar{\Sigma}$, we may double the arcs $d_1,\dots, d_{2g}$ to obtain a collection of pairwise disjoint simple closed curves $\Ds=(\Delta_1,\dots, \Delta_{2g})$ on $D(\Sigma)$, which form a set of compressing curves for $D(\Sigma)$. On $\Sigma\# \bar{\Sigma}$, we view the basepoint $w$ as living in a subset of the connected sum region which is complementary to $I\times [0,1]$.

The above construction naturally gives two Heegaard diagrams,
\[
 D(\cH):=\left(\Sigma\# \bar{\Sigma}, \as\cup \bar{\bs}, \Ds,w\right) \quad \text{and} \quad  \bar{D}(\cH):=\left(\bar{\Sigma}\# \Sigma, \Ds,\bar{\as}\cup \bs,w\right),
\]
 which are conjugates of each other.

The transition maps $\Psi_{\cH\to D(\cH)}$ and $\Psi_{D(\cH)\to \cH}$ have a simple description, as we now describe. Let $\bs'$ be small Hamiltonian translate of $\bs$, such
that $|\b_i\cap \b_j'|=2\delta_{ij}$. There is a 1-handle map
\[
F_1^{\bar{\b},\bar{\b}'}\colon \CF^-(\Sigma,\as,\bs)\to \CF^-(\Sigma\# \bar{\Sigma}, \as\cup \bar{\bs}, \bs\cup \bar{\bs}'),
\]
defined by the formula
\begin{equation}
F_1^{\bar{\b},\bar{\b}'}(\xs)=\xs\times \Theta_{\bar{\b},\bar{\b}'}^+, \label{eq:compound-1-handle-def}
\end{equation}
extended $\bF[U]$-equivariantly, where $\Theta_{\bar{\b},\bar{\b}'}^+$ is the top degree generator of $\CF^-(\bar{\Sigma}, \bar{\bs}, \bar{\bs}')$.  Next, we note that $(\Sigma\# \bar{\Sigma}, \bs\cup \bar{\bs}', \Ds,w)$ is the double of $(\Sigma,\bs,\bs',w)$ which is a diagram $(S^1\times S^2)^{\# g}$. Hence, assuming the diagram $(\Sigma\# \bar{\Sigma}, \bs\cup \bar{\bs}', \Ds,w)$ is weakly admissible, there is a cycle
\[
\Theta^+_{\b\cup \bar{\b}', \Delta}\in \CF^-(\Sigma\# \bar{\Sigma},\bs\cup \bar{\bs}', \Ds,w)\simeq  \CF^-\left((S^1\times S^2)^{\# g}\right),
\]
representing the top degree generator of homology. See Section~\ref{sec:admissibility} for more on admissibility and doubled diagrams.
 We define
\[
f_{\a\cup \bar{\b}}^{\b\cup \bar{\b}'\to \Delta}(\xs):=f_{\a\cup \bar{\b}, \b\cup \bar{\b}',\Delta}\left(\xs, \Theta^+_{\b\cup \bar{\b}',\Delta}\right).
\]

\begin{prop}\label{prop:doubling-transition-formula} If $\cH=(\Sigma,\as,\bs,w)$ is a Heegaard diagram and $D(\cH)=(\Sigma\# \bar{\Sigma}, \as\cup \bar{\bs}, \Ds,w)$ is a double of $\cH$, then
\[
\Psi_{\cH\to D(\cH)}\simeq f_{\a\cup \bar{\b}}^{\b\cup \bar{\b}'\to \Delta}\circ F_1^{\bar{\b},\bar{\b}'}.
\]
Dually,
\[
\Psi_{D(\cH)\to \cH}\simeq F_3^{\bar{\b},\bar{\b}'}\circ f_{\a\cup \bar{\b}}^{\Delta\to \b\cup \bar{\b}'}.
\]
\end{prop}

In Proposition~\ref{prop:doubling-transition-formula}, the map $F_3^{\bar{\b},\bar{\b}'}$ is the dual of the map $F_1^{\bar{\b},\bar{\b}'}$ in~\eqref{eq:compound-1-handle-def}, i.e., it is $\bF[U]$-equivariant and satisfies
\[
F_3^{\bar{\b},\bar{\b}'}(\xs\times \Theta)=\begin{cases} \xs& \text{if } \Theta=\Theta_{\bar{\b},\bar{\b}'}^-\\
0& \text{otherwise} .
\end{cases}
\]

Proposition~\ref{prop:doubling-transition-formula} is proven by interpreting the first composition as the cobordism map for a canceling sequence of 1-handles and 2-handles. The second map is similarly interpreted as the cobordism map for a canceling sequence of 2-handles and 3-handles. See \cite{ZemDuality}*{Proposition~7.2} for further details.

\subsection{Doubled diagrams for knots}
\label{sec:doubled-knots}

In this section, we describe an analog of the doubling construction from Section~\ref{sec:doubled-3-manifolds}, for doubly pointed diagrams of knots. We also describe an analog of Proposition~\ref{prop:doubling-transition-formula} for the twisted complexes $\buCF^-(\Sigma,\as,\bs,\hat{w},z)$. There is a new source of ambiguity in the construction, since we may add the new 1-handles near either $w$, or near $z$.

Suppose $\cH=(\Sigma,\as,\bs,w,z)$ is a Heegaard diagram for a doubly pointed knot $(K,w,z)$ in $Y$. We form a diagram $(\Sigma\# \bar{\Sigma}, \as\cup \bar{\bs}, \Ds, z, \bar{z})$ by doubling the diagram $\cH$, at the basepoint $w$. In this diagram, we delete $w$, and add a new basepoint $\bar{z}\in \bar{\Sigma}$, which is the image of $z$ on $\bar{\Sigma}$. 

\begin{figure}[ht!]
	\centering
	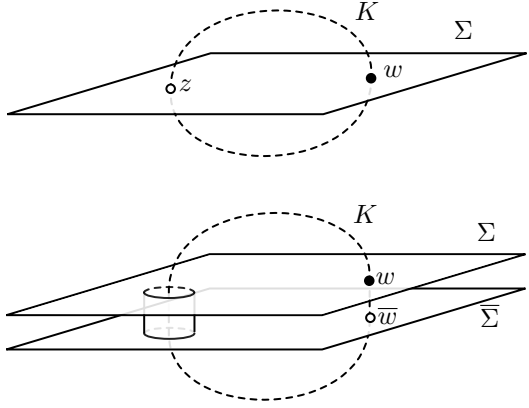
	\caption{A Heegaard diagram for a knot (top), and its double (bottom).}\label{fig:113}
\end{figure}

The diagram $(\Sigma\# \bar{\Sigma}, \as\cup \bar{\bs}, \Ds,z,\bar{z})$ may naturally be viewed as a diagram for $(Y,K)$, as follows. Let $D$ be an open disk containing the basepoint $w$, and we embed $\Sigma\# \bar{\Sigma}$ into $Y$ as $\d (N(\Sigma\setminus D))$. The knot $K$ no longer intersects the Heegaard surface at $w$, and instead intersects it at $z\in \Sigma$ and $\bar{z}\in \bar{\Sigma}$. See Figure~\ref{fig:113}.

 There is a 1-handle map
\[
F_{1}^{\bar{\b},\bar{\b}'} \colon \buCF^-\left(\Sigma,\as,\bs,\hat{w},z\right)\to \buCF^-\left(\Sigma\# \bar{\Sigma}, \as\cup \bar{\bs}, \bs\cup \bar{\bs}',z,\hat{\bar{z}}\right).
\]
 We think of $F_1^{\bar{\b},\bar{\b}'}$ as corresponding to attaching 1-handles near $w$, then moving $w$ to the position of $\bar{z}$.

We have the following analog to Proposition~\ref{prop:doubling-transition-formula}, for twisted complexes:
\begin{prop}\label{prop:doubling-twisted} The composition
\[
f_{\a\cup \bar{\b}}^{\b\cup \bar{\b}'\to \Delta}\circ F_1^{\bar{\b},\bar{\b}'}\colon \buCF^-\left(\Sigma,\as,\bs,\hat{w},z\right)\to \buCF^-\left(\Sigma\#\bar{\Sigma}, \as\cup \bar{\bs},\Ds,z,\hat{\bar{z}}\right)
\]
is chain homotopic to the map for the diffeomorphism which moves $w$ to $\bar{z}$ along a subarc of $K$, while keeping $z$ fixed, followed by the map from naturality. Dually,
\[
F_3^{\bar{\b},\bar{\b}'}\circ f_{\a\cup \bar{\b}}^{\Delta\to \b\cup \bar{\b}'}\colon \buCF^-\left(\Sigma\# \bar{\Sigma}, \as\cup \bar{\bs}, \Ds, z, \hat{\bar{z}}\right)\to \buCF^-(\Sigma,\as,\bs,\hat{w},z),
\]
is chain homotopic to the diffeomorphism map moving $\bar{z}$ to $w$ while fixing $z$, followed by the map from naturality.
\end{prop}

\subsection{The involution via doubling}
\label{sec:involution-doubling}
A doubled diagram $D(\cH)$ may be viewed as the double of both $\cH$ and $\bar{\cH}$. Hence, the involution may be described as a composition of a doubling map, followed by an un-doubling map, using the formulas from Proposition~\ref{prop:doubling-transition-formula} and \ref{prop:doubling-twisted}. More concretely, if $(\Sigma,\as,\bs,w)$ is a Heegaard diagram, then
\begin{equation}
\iota\simeq \eta\circ F_{3}^{\a,\a'}\circ f_{\a\cup \bar{\b}}^{\Delta\to \a'\cup \bar{\a}} \circ f_{\a\cup \bar{\b}}^{\b\cup \bar{\b}'\to \Delta} \circ F_1^{\bar{\b},\bar{\b}'}. \label{eq:decompose-involution-doubling}
\end{equation}
(Recall that $\eta$ is the tautological chain isomorphism from Section \ref{sec:3-manifolds}.)

Similarly, if $(\Sigma,\as,\bs,\hat{w},z)$ is a diagram for an oriented null-homologous knot $K$ in $Y$, with $w$ marked as the special basepoint, then the twisted involution defined in ~\eqref{eq:def-iota-3} satisfies
\[
\underline{\iota}\simeq \frF_{z\to w}\circ \underline{\eta}\circ  F_{3}^{\a,\a'} \circ f_{\a\cup \bar{\b}}^{\Delta\to \a'\cup \bar{\a}}\circ f_{\a\cup \bar{\b}}^{\b\cup \bar{\b}'\to \Delta}\circ F_1^{\bar{\b},\bar{\b}'}.
\]
(Recall that $\underline{\eta}$ is the tautological chain isomorphism from Section \ref{sec:twisted-complexes-def}.)

\begin{figure}[h]
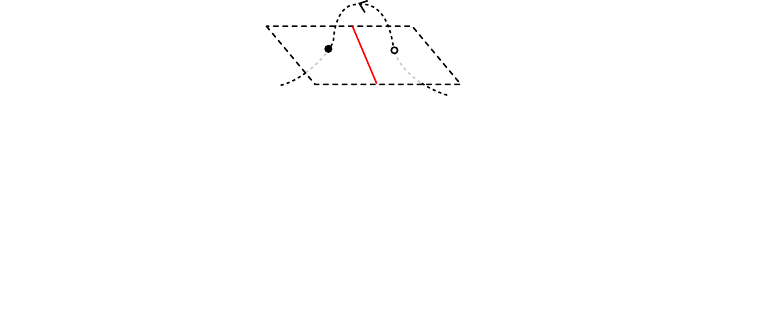
\caption{A schematic of the map $\iota_{\bB}$ defined using the wrong choice of connected sum point. The top image indicates the orientation of $K$, and the bottom row indicates the diagrams which appear in the definition of $\iota_{\bB}$. Arrows indicate movement of the basepoints throughout. The diagram on the left is $(\Sigma,\as,\bs,w,z)$. The middle-left is the doubled diagram. The middle-right is the diagram $(\Sigma,\bar \bs, \bar \as, \bar z, \bar w)$. The final right-hand diagram indicates the path for the flip map. Note that $w$ is sent through a full loop around $K$ as we go from left to right.}
\label{fig:4}
\end{figure}

\begin{rem}
\label{rem:careful-twisting} At this point we pause for one important note concerning the definition of the knot involution. The choice of which basepoint at which we form the connect sum in Figure~\ref{fig:113} is not arbitrary; it is instead determined by our choice of flip map $\frF_{z\to w}$. In particular, the choice of basepoint at which we form the connect sum determines the direction in which the basepoint-moving diffeomorphism in the definition of $\iota_K$ twists along the knot $K$. If we take the connected sum near $w$, the basepoint-moving diffeomorphism from $\iota_K$ twists oppositely to the orientation of $K$, whereas the flip map twists coherently to $K$. The net effect is that the map $\iota_{\bB}:=\frF_{z\to w}\circ \iota_K$ carries the basepoint $w$ through a nullhomotopic loop. If we instead take the connected sum near $z$ then both $\iota_K$ and $\frF_{z\to w}$ move basepoints positively along $K$. Hence $\iota_{\bB}$ sends the basepoints through a full loop around $K$, as in Figure~\ref{fig:4}. Our proof of the surgery exact triangle does not carry through with this choice. Concretely, one of the model counts in Lemma~\ref{lem:ends-r-to-infty} will fail.
\end{rem}

\subsection{Doubling and the involutive hypercube of a triad}
\label{sec:outline}
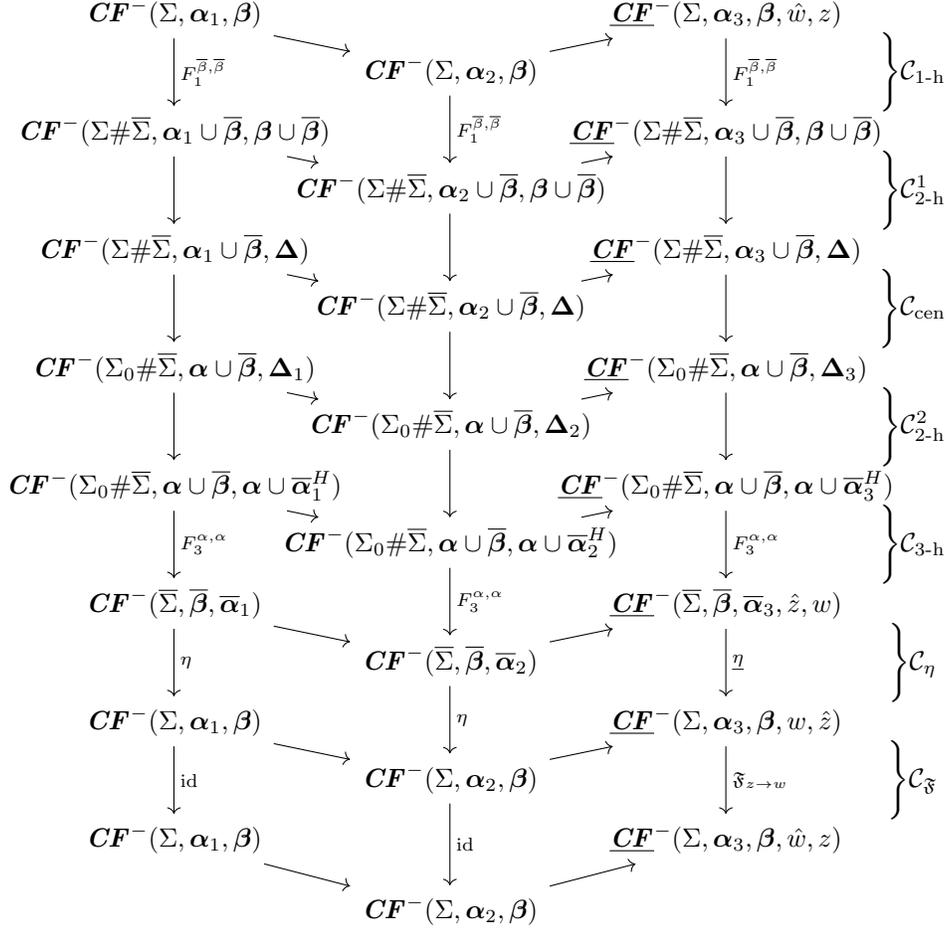
\begin{figure}[h]
	\centering
\[
\begin{tikzcd}[row sep=.3cm, column sep=-1.1 cm]
\bCF^-(\Sigma,\as_1,\bs)\arrow[rd]\arrow[dd,"F_1^{\bar{\b},\bar{\b}}"]
& \,
&\buCF^- (\Sigma,\as_3,\bs,\hat{w},z)\arrow[dd,"F_1^{\bar{\b},\bar{\b}}"]
\\[-.5cm]
&\bCF^-(\Sigma,\as_2,\bs) \arrow[ur] \arrow[dd,"F_1^{\bar{\b},\bar{\b}}"]
&\,&[.6cm] \Bigg\} \cC_{\oneh}
\\[-.5cm]
\bCF^-(\Sigma\# \bar{\Sigma},\as_1\cup \bar{\bs},\bs\cup \bar{\bs})\arrow[dr]\arrow[dd]
& \, &\buCF^-(\Sigma\# \bar{\Sigma},\as_3\cup \bar{\bs},\bs\cup \bar{\bs})\arrow[dd]&
\\[-.5cm]
&
\bCF^-(\Sigma\# \bar{\Sigma},\as_2\cup \bar{\bs},\bs\cup \bar{\bs})\arrow[ur]\arrow[dd]
&\, &\Bigg\}\cC_{\twoh}^1
\\[-.5cm]
\bCF^-(\Sigma\# \bar{\Sigma},\as_1\cup \bar{\bs},\Ds)\arrow[dr]\arrow[dd]
&\,&\buCF^-(\Sigma\# \bar{\Sigma},\as_3\cup \bar{\bs},\Ds)\arrow[dd]&
\\[-.5cm]
&
\bCF^-(\Sigma\# \bar{\Sigma},\as_2\cup \bar{\bs},\Ds) \arrow[ur] \arrow[dd]
&\,&\Bigg\}\cC_{\cen}
\\[-.5cm]
\bCF^-(\Sigma_0\# \bar{\Sigma},\as\cup \bar{\bs},\Ds_1)\arrow[dr]\arrow[dd]
&\, &\buCF^-(\Sigma_0\# \bar{\Sigma},\as\cup \bar{\bs},\Ds_3)\arrow[dd]&
\\[-.5cm]
&\bCF^-(\Sigma_0\# \bar{\Sigma},\as\cup \bar{\bs},\Ds_2)\arrow[ur]\arrow[dd]
&\,&\Bigg\}\cC_{\twoh}^2
\\[-.5cm]
\bCF^-(\Sigma_0\# \bar{\Sigma},\as\cup \bar{\bs},\as\cup \bar{\as}^H_1)\arrow[dr]\arrow[dd,"F_3^{\a,\a}"]
&\,&\buCF^-(\Sigma_0\# \bar{\Sigma},\as\cup \bar{\bs},\as\cup \bar{\as}^H_3)\arrow[dd,"F_3^{\a,\a}"]&
\\[-.5cm]
&\bCF^-(\Sigma_0\# \bar{\Sigma},\as\cup \bar{\bs},\as\cup \bar{\as}^H_2)\arrow[ur] \arrow[dd,"F_3^{\a,\a}"]
&\,&\Bigg\}\cC_{\thrh}
\\[-.5cm]
\bCF^-(\bar{\Sigma},\bar{\bs},\bar{\as}_1)\arrow[dr]\arrow[dd, "\eta"]
& \,&\buCF^-(\bar{\Sigma},\bar{\bs},\bar{\as}_3,\hat{z},w) \arrow[dd, "\underline{\eta}"] &
\\[-.5cm]
&\bCF^-(\bar{\Sigma},\bar{\bs}, \bar{\as}_2)\arrow[ur]\arrow[dd,"\eta"]
&\,& \Bigg\}\cC_{\eta}
\\[-.5cm]
\bCF^-(\Sigma,\as_1,\bs)\arrow[dr]\arrow[dd,"\id"]&\,& \buCF^-(\Sigma,\as_3,\bs,w,\hat{z})\arrow[dd, "\frF_{z\to w}"]&
\\[-.5cm]
& \bCF^-(\Sigma,\as_2,\bs)\arrow[ur]\arrow[dd, "\id"]
&\,&\Bigg\}\cC_{\frF}
\\[-.5cm]
\bCF^-(\Sigma,\as_1,\bs)\arrow[dr]&\,&\buCF^-(\Sigma,\as_3,\bs,\hat{w},z)&\\
&\bCF^-(\Sigma,\as_2,\bs)\arrow[ur]& 
\end{tikzcd}
\]
	\caption{The hyperbox whose compression is the hypercube in Theorem~\ref{thm:main-hypercube}. Length 2 and 3 arrows are not shown.}\label{fig:large-hyperbox}
\end{figure}

We now outline our construction of the main hypercube, from Theorem~\ref{thm:main-hypercube}. We will construct a   hyperbox of chain complexes, of size $(1,1,7)$, which we show in Figure~\ref{fig:large-hyperbox}. In the subsequent sections, we will construct the hypercubes labeled $\cC_{\oneh}$, $\cC_{\twoh}^1$, $\cC_{\cen}$, $\cC_{\twoh}^2$, $\cC_{\thrh}$, $\cC_{\eta}$ and $\cC_{\frF}$.

We have not drawn any of the length 2 or 3 arrows in any of the hypercubes. The compression of the above hyperbox of chain complexes will be our main hypercube. Our construction is organized as follows:
\begin{itemize}
\item  $\cC_{\oneh}$ and $\cC_{\thrh}$ are the \emph{1-handle} and \emph{3-handle hypercubes}. They are special cases of the \emph{hypercubes for stabilization}, constructed in Section~\ref{sec:hypercube-for-stabilization}.
\item  $\cC_{\twoh}^1$ and $\cC_{\twoh}^2$ are the \emph{2-handle hypercubes}. They are constructed in Sections~\ref{sec:peripheral-hypercubes} and~\ref{sec:hypercube-C-3}.
\item  $\cC_{\cen}$ is the \emph{central hypercube}. It is constructed in Section~\ref{sec:central-hypercube}.
\item $\cC_{\eta}$ is the \emph{tautological hypercube}. It is described in Section~\ref{sec:flip-hypercube}. 
\item $\cC_{\frF}$ is the \emph{flip-map hypercube}. It is constructed in Section~\ref{sec:flip-hypercube}.
\end{itemize}

The surface $\Sigma_0\# \bar{\Sigma}$, and the attaching curves $\Ds_i$, $\as$ and $ \bar{\as}_i^H$ will be described in Section~\ref{sec:peripheral-hypercubes}. Note that the columns are not identical to the doubling model described in Section~\ref{sec:involution-doubling} (for example, the model described therein did not feature $\Sigma_0\# \bar{\Sigma}$). We prove in Lemma~\ref{lem:columns=change-of-diagrams} that the vertical maps in $\cC_{\cen}$ coincide with change of diagram maps. Using this fact, the same logic as in Proposition~\ref{prop:doubling-transition-formula} also shows that the composition of each vertical arrow in the columns of Figure~\ref{fig:large-hyperbox} gives the involution.

\subsection{Admissibility and doubled diagrams}
\label{sec:admissibility}

In this section, we prove a technical result concerning admissibility of doubled diagrams, which will be used as the model for more complicated admissibility arguments later.

\begin{define} Suppose $\cD=(\Sigma,\gs_1, \dots, \gs_n,w)$ is a Heegaard multi-diagram.
\begin{enumerate}
\item A \emph{(rational) periodic domain} is an integral (resp. rational) 2-chain $D$ on $\Sigma$ such that $\d D$ is a (rational) linear combination of the curves in $\gs_1\cup \cdots \cup \gs_n$, and furthermore, $n_{w}(D)=0$.
\item $\cD$ is \emph{weakly admissible} if every non-zero periodic domain has both positive and negative multiplicities.
\end{enumerate}
\end{define}

\begin{figure}[ht!]
	\centering
	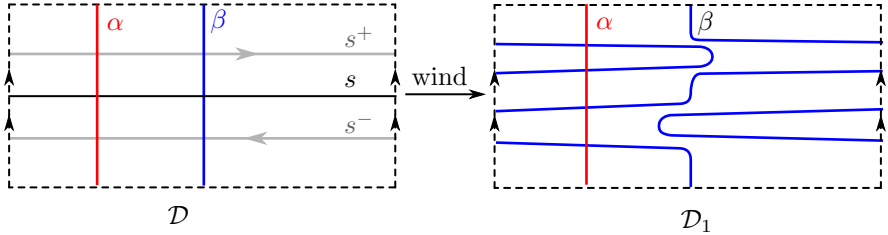
	\caption{Winding along $s_i$ once to obtain $\cD_1$.}\label{fig:admissible-1}
\end{figure}

Since we are working over the power series ring $\bF\llsquare U \rrsquare$, it is sufficient to work with diagrams which are weakly admissible. (See \cite{OSDisks}*{Lemma~4.13}.) Admissibility is usually achieved by the process of \emph{winding} \cite{OSDisks}*{Section~5}, which we now sketch. If $(\Sigma,\as,\bs,w)$ is a Heegaard diagram, we pick pairwise disjoint, embedded curves $s_1,\dots, s_g$ which are dual to the $\bs$ curves, in the sense that
\[
|s_i\cap \beta_j|=\delta_{ij}.
\]
We then consider two pushoffs $s_i^+$ and $s_i^-$ of $s_i$. We form a diagram $\cD_N$ by winding $\b_i$ $N$-times, positively, along $s_i^+$, and $N$-times, negatively, along $s_i^-$. See Figure~\ref{fig:admissible-1}.

Ozsv\'{a}th and Szab\'{o} show that if one winds a Heegaard multi-diagram sufficiently, then weak admissibility is obtained \cite{OSDisks}*{Lemma~5.4}. In particular, if $(\Sigma\# \bar{\Sigma}, \as\cup \bar{\bs}, \Ds,w)$ is a doubled diagram, we may wind $\as\cup \bar{\bs}$ to achieve admissibility. This is not quite suitable for our purposes. Instead, we now show that for a doubled diagram, we need only wind $\bar{\bs}$ to achieve admissibility.

 We now prove a helpful admissibility result for doubled diagrams where $\Ds$ are adapted to $\as$. The following will be a model for proving admissibility of Heegaard multi-diagrams appearing in the construction of the central hypercube.

\begin{lem}\label{lem:admissibility-2-tuple}
 Suppose that $\cH=(\Sigma,\as,\bs,w)$ is a Heegaard diagram, and $D(\cH)=(\Sigma\# \bar{\Sigma}, \as\cup\bar{\bs}, \Ds,w)$ is a double of $\cH$. Then, after winding the $\bar{\bs}$ curves sufficiently on $\bar{\Sigma}$, the diagram $D(\cH)$ becomes weakly admissible.
\end{lem}
\begin{proof} Let $s_1,\dots, s_g\subset \bar{\Sigma}$ denote dual curves to the $\bar{\bs}$ curves, along which to wind. Let $D(\cH)_N$ denote the diagram obtained by winding $\bar{\bs}$ $N$-times, along positive and negative push-offs of the curves $s_1,\dots, s_g$, as in Figure~\ref{fig:admissible-1}.

 Let $\cP_N^{\Q}$ denote the set of rational periodic domains on $D(\cH)_N$. It is clearly sufficient to rule out the existence of a nonzero rational periodic domain with nonnegative coefficients.
 
 Note that there is a canonical isomorphism
 \[
\phi_N\colon \cP^{\Q}\to \cP^{\Q}_N,
 \]
 which preserves the multiplicities which are not in the winding region.

We suppose for contradiction that $D(\cH)_N$ is not weakly admissible, for arbitrarily large $N$. This means there is a sequence of periodic domains, $\{D_i\}_{i\in \N}$ such that $D_i$ is a nonzero periodic domain on $D(\cH)_{N_i}$ with only nonnegative multiplicities, where $N_i\to \infty$. By rescaling, we assume that
\[
\| D_i\|_{L^{\infty}}=1,
\]
for all $i$. We consider the domains $D_i':=\phi_{N_i}^{-1}(D_i)$. Noting that  the set of multiplicities of $D_i'$ is a subset of the multiplicities of $D_i$, it follows that $\|D_i'\|_{L^\infty}\le 1$, for all $i$. Let $D_{\infty}'$ be the limit of a convergent subsequence of the $D_i'$. We make the following claims:
\begin{enumerate}[label=($d$-\arabic*), ref=$d$-\arabic*, leftmargin=*,widest=IIII]
\item\label{d-infty-1} $D_\infty'$ is nonnegative.
\item\label{d-infty-2} $D_\infty'$ is a real periodic domain such that  $\d D_{\infty}'$ is a linear combination of the $\as$ and $\Ds$ curves (in particular, the boundary contains no multiples of $\bar{\bs}$).
\item\label{d-infty-3} $\| D_\infty'\|\ge 1/2$.
\end{enumerate}

Claim~\eqref{d-infty-1} follows since $D_i$ is nonnegative, and the set of multiplicities of $D_i'$ is a subset of the set of multiplicities of $D_i$.

 To prove claim~\eqref{d-infty-2}, it is sufficient to show that $D_\infty'$ has equal multiplicity on either side of a curve $\b$ in $\bar{\bs}$. Assume $D_i'\to D_\infty'$. As the difference of $D_\infty'$ across $\b$ at a point $p\in \b$ is independent of $p$, it is sufficient to consider $p$ adjacent to the winding region, but not contained in it. Let $x_i^l$ and $x_i^r$ denote the multiplicities of $D_i'$ on either side of $\b$ at $p$ (the order is unimportant). Let $x_\infty^l$ and $x_\infty^r$ denote the limits (which are the values of $D_\infty'$ near $p$). By construction, $D_i$ has multiplicities of both $x_i^l+N_i(x_i^l-x_i^r)$ and $x_i^r+N_i(x_i^r-x_i^l)$.  Since for all $i$,
 \[
0\le x_i^l+N_i(x_i^l-x_i^r)\le 1 \quad \text{and} \quad 0 \le x_i^l\le 1,
 \]
 we obtain that $x_\infty^l-x_\infty^r= 0$,  by taking the limit. In particular, $D_\infty'$ has equal multiplicities on both sides of $\b$. Claim~\eqref{d-infty-2} follows.

We now consider claim \eqref{d-infty-3}. We claim that each $D_i'$ has a multiplicity which is at least $1/2$. By construction, each $D_i$ has a multiplicity which is 1. If this multiplicity occurs outside of the winding region, then we are done. Suppose that it occurs inside of the winding region for a curve $s\in \{s_1,\dots, s_g\}$ which is dual to $\b\in \bs$. Let $\delta$ denote the coefficient of $\b$ in $\d D_i'$, which is equivalently the difference in multiplicities on the two sides of $\b$. We may assume $\delta\ge 0$ by changing the orientation of $\beta$, if necessary. The multiplicities of $D_i$ are exactly $m+k\cdot \delta$, where $m$ is a multiplicity of $D_i'$ and $k\in \{-N_i,\dots, N_i\}$. For $m+k\cdot \delta$ to be the maximal multiplicity of $D_i'$, we must have $k=N_i$. However, 
\[
0\le m-N_i\cdot \delta.
\] 
Adding $1=m+N_i\cdot \delta$ and dividing by 2 we obtain $m\ge 1/2$, establishing~\eqref{d-infty-3}.

However,  there are no non-trivial, real 2-chains on $\Sigma\# \bar{\Sigma}$ which have multiplicity 0 at $w$ and satisfy~\eqref{d-infty-1}, \eqref{d-infty-2} and~\eqref{d-infty-3}. Indeed, such a 2-chain would give a non-trivial linear relation amongst the classes $[\a_1],\dots [\a_g], [\Delta_1],\dots, [\Delta_{2g}]$ in $H_1(\Sigma\# \bar{\Sigma};\R)$. However, the subspaces 
\[
\Span([\a_1],\dots, [\a_g]), \, \, \Span([\Delta_1],\dots, [\Delta_{2g}])\subset H_1(\Sigma\# \bar{\Sigma};\R)
\]
 have trivial intersection, and are themselves $g$ and $2g$-dimensional, respectively.
\end{proof}

To ensure certain stabilization results for holomorphic curves, it will be helpful to introduce the following restricted version of a doubled diagram:

\begin{define}\label{def:adapted-to-alpha}
Suppose that $\cH=(\Sigma,\as,\bs,w)$ is a Heegaard diagram, and $(\Sigma\# \bar{\Sigma}, \as\cup \bar{\bs}, \Ds,w)$ is a double of $\cH$. We say that the curves $\Ds$ are \emph{adapted to $\as$} if they are constructed by doubling  a basis of arcs $d_1,\dots d_{2g}$ for $H_1(\Sigma\setminus N(w), I;\Z)$, where $I$ is a subarc of $\d N(w)$, as follows. The arcs $d_{2},d_4,\dots, d_{2g}$ are chosen (arbitrarily) to satisfy
\[
|d_{2i}\cap \a_{j}|=\delta_{ij}.
\]
The arc $d_{2i}$ induces two arcs from $\a_i$ to $\bar{\a}_i$. Let $\lambda_i$ denote one of these arcs (chosen arbitrarily). We let $\Delta_{2i-1}$ be obtained by handlesliding $\a_i$ across $\bar{\a}_i$, using the arc $\lambda_i$ as a guide. We assume $\Delta_{2i-1}$ is chosen to be disjoint from $\a_i$. 
\end{define}

Note that Lemma~\ref{lem:admissibility-2-tuple} implies that $\as$-adapted diagrams may be made weakly admissible by winding $\bar{\bs}$ sufficiently.

 \begin{figure}[ht!]
	\centering
\begingroup%
  \makeatletter%
  \providecommand\color[2][]{%
    \errmessage{(Inkscape) Color is used for the text in Inkscape, but the package 'color.sty' is not loaded}%
    \renewcommand\color[2][]{}%
  }%
  \providecommand\transparent[1]{%
    \errmessage{(Inkscape) Transparency is used (non-zero) for the text in Inkscape, but the package 'transparent.sty' is not loaded}%
    \renewcommand\transparent[1]{}%
  }%
  \providecommand\rotatebox[2]{#2}%
  \newcommand*\fsize{\dimexpr\f@size pt\relax}%
  \newcommand*\lineheight[1]{\fontsize{\fsize}{#1\fsize}\selectfont}%
  \ifx\svgwidth\undefined%
    \setlength{\unitlength}{139.25724382bp}%
    \ifx\svgscale\undefined%
      \relax%
    \else%
      \setlength{\unitlength}{\unitlength * \real{\svgscale}}%
    \fi%
  \else%
    \setlength{\unitlength}{\svgwidth}%
  \fi%
  \global\let\svgwidth\undefined%
  \global\let\svgscale\undefined%
  \makeatother%
  \begin{picture}(1,0.82821884)%
    \lineheight{1}%
    \setlength\tabcolsep{0pt}%
    \put(0,0){\includegraphics[width=\unitlength,page=1]{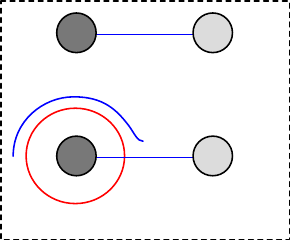}}%
    \put(0.49686761,0.73808201){\color[rgb]{0,0,0.81960784}\makebox(0,0)[t]{\lineheight{1.25}\smash{\begin{tabular}[t]{c}$\Delta_{2i}$\end{tabular}}}}%
    \put(0.72106261,0.52668442){\color[rgb]{0,0,0.81960784}\makebox(0,0)[t]{\lineheight{1.25}\smash{\begin{tabular}[t]{c}$\Delta_{2i-1}$\end{tabular}}}}%
    \put(0.40387334,0.31405485){\color[rgb]{1,0,0}\makebox(0,0)[rt]{\lineheight{1.25}\smash{\begin{tabular}[t]{r}$\a_i$\end{tabular}}}}%
    \put(0,0){\includegraphics[width=\unitlength,page=2]{fig184.pdf}}%
  \end{picture}%
\endgroup%

	\caption{A schematic of a choice of $\Ds\subset\Sigma\# \bar{\Sigma}$, which are adapted to $\as\subset \Sigma$. Here, $\Sigma$ is the genus 1 region on the left, and $\bar{\Sigma}$ is the genus 1 region on the right.}\label{fig:184}
\end{figure}

\section{Moduli spaces and almost complex structures}
\label{sec:moduli-spaces}

In this section, we describe some analytic input we will need. We begin by defining our spaces of almost complex structures and moduli spaces for holomorphic disks in Section~\ref{ac:disks}.  In Section~\ref{sec:matched-moduli-spaces} we describe some results on transversality of holomorphic curves with pointwise constraints, including pointwise tangency constraints, which we use frequently. In Sections~\ref{sec:cylindrical-boundary}--\ref{sec:counting-cylindrical-boundary} we describe the necessary input on boundary degenerations, including Proposition~\ref{prop:relaxed-count}, which has not appeared previously in the literature.

\subsection{Almost complex structures for holomorphic disks}

\label{ac:disks}

Let $\omega=dA+ ds\wedge dt$ be a fixed, split symplectic form on $\Sigma\times [0,1]\times \R$. Suppose $p_1,\dots, p_n$ is a fixed collection of points in $\Sigma\setminus (\as\cup \bs)$, and $D_1,\dots, D_n$ are a collection of small regular neighborhoods of $p_1,\dots, p_n$ in $\Sigma$. 

\begin{define} We say an almost complex structure $J$ on $\Sigma\times [0,1]\times \R$ is \emph{admissible} if it satisfies the following:
\begin{enumerate}[ref=$J$-\arabic*, label=($J$-\arabic*),leftmargin=*, widest=IIII]
\item\label{J1} $J$ is tamed by $\omega$.
\item\label{J2} $J$ is split on  $(D_1\cup \cdots \cup D_n)\times [0,1]\times \R$.
\item\label{J3} $J$ is translation invariant under the $\R$-action.
\item\label{J4} $J \d/\d s=\d/\d t$.
\item\label{J5} $T_p \Sigma\oplus \{0\}$ is a complex line in a neighborhood of $(\as\cup \bs)\times [0,1]\times \R$, and on $\Sigma\times \{0,1\}\times \R$.
\end{enumerate}
\end{define}
Note that our axioms \eqref{J1}--\eqref{J5} are equivalent to the ones that Lipshitz labels ($J$-1)--($J$-4) and ($J$-5') (in particular, our numbering does not coincide with his). A similar family of almost complex structures are called \emph{relaxed} in \cite{LOTBordered}*{p. 229}.

\subsection{Moduli spaces of disks}
\label{sec:moduli-spaces-disks}

Suppose $(\Sigma,\as,\bs,w)$ is a Heegaard diagram. Suppose that $S$ is a smooth (not nodal) Riemann surface with boundary and a finite collection of boundary punctures, designated as $+$ or $-$. Suppose $J$ is an almost complex structure on $\Sigma\times [0,1]\times \R$ which satisfies \eqref{J1}--\eqref{J5}. We are interested in pairs $(u,j)$, where $j$ is an almost complex structure on $S$, and $u$ is a smooth map
\[
u\colon (S,\d S)\to (\Sigma\times [0,1]\times \R, (\as\times \{1\}\times \R) \cup (\bs\times \{0\}\times \R)),
\]
 representing the homology class $\phi$,  and satisfying the following:
\begin{enumerate}[ref=($M$-\arabic*), label=($M$-\arabic*),leftmargin=*, widest=IIII]
\item \label{M-1} $u$ is $(j,J)$-holomorphic.
\item\label{M-2} $u$ is proper.
\item\label{M-3} For each $i$, $u^{-1}(\a_i\times \{1\}\times \R)$ and $u^{-1}(\b_i\times \{0\}\times \R)$ each consist of exactly one component of $\d S$.
\item\label{M-4} If $q$ is a $+$ puncture of $S$, then $\lim_{z\to q} (\pi_{\R}\circ u)(z)=+\infty$.
\item\label{M-5} If $q$ is a $-$ puncture of $S$, then $\lim_{z\to q} (\pi_{\R}\circ u)(z)=-\infty$.
\item\label{M-6} $u$ has finite energy.
\item\label{M-7} $\pi_{\bD}\circ u$ is non-constant on each component of $S$, where $\pi_{\bD}$ denotes the projection onto $[0,1]\times \R$. 
\end{enumerate}
Compare \cite{LipshitzCylindrical}*{p. 960}, \cite{LOTBordered}*{p. 63}. We write $\cM_J(S,\phi)$ for the set of such pairs, modulo the relation that $(u,j)\sim (u',j')$ if there is a $(j,j')$-holomorphic diffeomorphism $\phi\colon S\to S$ such that $u'\circ \phi=u$.  

The next object of interest is the moduli space of holomorphic curves with marked points. Suppose $S$ is a Riemann surface, as above, and $\ve{q}=(q_1,\dots, q_n)$ is an ordered collection of distinct points along $S$. We write $S^{\qs}$ for a surface, equipped with a choice of $\ve{q}$. We write $\cM_J(S^{\qs},\phi)$ for the set of pairs $(u,j)$ satisfying~\ref{M-1}--\ref{M-7}, except we quotient only by diffeomorphisms $\phi\colon S\to S$ which fix $\ve{q}$ pointwise.

We say two marked surfaces $S^{\qs}$ and $T^{\ps}$ are \emph{equivalent} if there is a diffeomorphism $\phi\colon S\to T$ which takes the positive (resp. negative) punctures to positive (resp. negative) punctures, and restricts to an order preserving bijection between $\qs$ and $\ps$.

\subsection{Moduli spaces with matching conditions}
\label{sec:matched-moduli-spaces}

 A common construction is to restrict the values obtained by holomorphic curves on their marked points. A refinement of this construction is to consider holomorphic curves with additional tangency constraints at their marked points. In this section, we provide background on the construction.

Suppose that $(W,J)$ is an almost complex manifold and $p\in W$. A \emph{complex $r$-jet} at $p$ is an equivalence class of germs of $J$-holomorphic maps
\[
g\colon (\C,0)\to (W,p),
\]
where two germs, $g_1$ and $g_2$, are \emph{$r$-equivalent} if  $\left(\tfrac{\d^i}{\d^i x}g_1\right)(0)=\left(\tfrac{\d^i}{\d^i x}g_2\right)(0)$ 
for all $0\le i\le r$, where $(x,y)$ are the standard coordinates on $\C$. We write $\Jet_p^r(W,J)$ for the set of complex $r$-jets at $p$. 

Using the local existence of holomorphic curves, one can identify $\Jet_{p}^r(W^{2n},J)$ with $\C^{rn}$ \cite{ZehmischJets}*{Proposition~2.1}. By convention, we interpret the 0-jet space as a single point. When $W=\C^{n}$, coordinates on the jet spaces may be given by coefficients of complex Taylor series.

It is straightforward to identify $\Jet^1_p(W,J)$ with $T_p W$, and hence $\Jet^1_p(W,J)$ is naturally a complex vector space. For $r>1$, $\Jet^r_p(W,J)$ is not naturally a vector space, however it does have a distinguished point $0$, given by the constant map at $p$. We define
\[
\Jet^r(W,J):=\bigcup_{p\in W} \Jet_p^r(W,J),
\]
which forms a fiber bundle over $W$.

If $(W_1,J_1)$ and $(W_2,J_2)$ are two almost complex manifolds and $u\colon W_1\to W_2$ is $(J_1,J_2)$-holomorphic, there is an $r$-jet evaluation map
\[
J_p^r u\colon \Jet_p^r(W_1,J_1)\to \Jet_{u(p)}^r(W_2,J_2),
\]
given by $(J_p^r u)(g)=u\circ g$.

\begin{define} We say $M=(M_0,\ve{t})$ is a \emph{matching condition} if $M_0\subset ((0,1)\times \R)^k$ is a smoothly embedded submanifold, and $\ve{t}=(n_1,\dots, n_k, m_1,\dots, m_k)$ is a tuple of nonnegative integers.
\end{define}

If $\ve{t}=(n_1,\dots, n_k,m_1,\dots, m_k)$, we write
\[
|\ve{t}|:=n_1+\cdots + n_k+m_1+\cdots +m_k.
\]

Suppose $M$ is a matching condition, $S^{\qs}$ is a marked Riemann surface for disks, $p\in \Sigma\setminus (\as\cup \bs)$ is a chosen point, and $J$ is an admissible almost complex structure on $\Sigma\times [0,1]\times \R$ which is split in a cylindrical neighborhood of $\{p\}\times [0,1]\times \R$. We define the \emph{moduli space of $M$-matched disks} to be
\begin{equation}
\cM_J(S^{\qs},\phi,M):=\left\{ (u,j)\in \cM_J(S,\phi)\middle \vert \begin{array}{c}
(\pi_{\Sigma}\circ u)(q_i)=p\\
\left((\pi_{\bD}\circ u)(q_1),\dots,(\pi_{\bD}\circ u)(q_k)\right)\in M_0,\\ 
J^{n_i}_{q_i}(\pi_\Sigma\circ u)=0,\\
J^{m_i}_{q_i}(\pi_\bD\circ u)=0,
\\
\text{ for } i=1,\dots, k
 \end{array}  \right\}.
 \label{eq:def-matched-moduli-space-most-general}
\end{equation}
In ~\eqref{eq:def-matched-moduli-space-most-general}, if $n_i=0$, we interpret $J_{q_i}^{n_i}(\pi_\Sigma \circ u)=0$ as holding vacuously.

If $p_1,\dots, p_\ell$ is a collection of special points on $\Sigma\setminus (\as\cup \bs)$, and we are equipped with a function $f\colon \qs\to \{p_1,\dots, p_\ell\}$, then we may also define the moduli space of $M$-matched disks for $S^{\qs,f}$ by a simple adaptation of~\eqref{eq:def-matched-moduli-space-most-general}.

\subsection{Expected dimensions}

Given a matching condition $M=(M_0,\ve{t})$, we will write
\[
|M|:=\codim(M_0) +2|\ve{t}|.
\]
We define the \emph{embedded, matched index} of a pair $(\phi,M)$ to be
\[
\ind_{\emb}(\phi,M):=\mu(\phi)-|M|,
\]
where $\mu(\phi)$ denotes the ordinary Maslov index of the class $\phi$.

\begin{prop}\label{prop:dimension-counts} Suppose that $S^{\qs}$ is a Riemann surface with $k$ marked points $\qs=(q_1,\dots, q_k)$, and $M=(M_0,\ve{t})$ is a matching condition at a point $p\in \Sigma$.
\begin{enumerate}
\item\label{expected-dimensions-1}  For generic $J$ satisfying \eqref{J1}--\eqref{J5}, the $M$-matched moduli space $\cM_J(S^{\qs},\phi,M)$ is a smooth manifold of dimension
\[
\ind(u)-|M|,
\]
where $\ind(u)$ denotes the Fredholm index, viewing $u$ as a holomorphic curve with no marked points.
\item\label{expected-dimensions-2} Furthermore, the dimension of $\cM_J(S^{\qs},\phi,M)$ also coincides with the quantity
\[
\ind_{\emb}(\phi,M)-2\Sing(u),
\]
where $\Sing(u)$ denotes the singularity number of $u$.
\end{enumerate}
\end{prop}

In part~\eqref{expected-dimensions-2} of Proposition~\ref{prop:dimension-counts}, the quantity $\Sing(u)$ is the signed count of double points, in an equivalent singularity. An interior (positive) double point contributes $1$, and a boundary double point contributes $\tfrac{1}{2}$.

Proposition~\ref{prop:dimension-counts} follows from a standard line of reasoning, as we now sketch. We begin with part~\eqref{expected-dimensions-1}. In the restricted case that $\ve{t}=(0,\dots, 0)$, the result follows from a straightforward combination of the proofs of \cite{LipshitzCylindrical}*{Proposition~3.7} and \cite{McDuffSalamonSymplectic}*{Proposition~3.4.2}.
For a generic $J$, $\cM_J(S^{\qs},\phi)$ is a smooth manifold, and the map $\ev_{\qs}$ is transverse to $M_0$, giving the claimed dimension count immediately. 

The case when $|\ve{t}|>0$ is slightly more involved, as we now describe. We note that $J_{q_i}^{n_i}(\pi_\Sigma\circ u)$ lives in $\Map (\Jet^{n_i}_{q_i}(S,j), \Jet^{n_i}_{(\pi_\Sigma\circ u)(q_i)}(\Sigma,j_\Sigma))$, which depends on $j$. This is slightly inconvenient for computing expected dimensions, and one solution is to pick complex coordinates carefully on $S$, near $\ve{q}$, as we now describe.

Write $\cJ(S)$ for the set of almost complex structures on $S$. The \emph{Teichm\"{u}ller space} of $S^{\qs}$, which we denote by $\cT(S^{\qs})$, is the the quotient of $\cJ(S)$ by $\Diff_0(S,\qs)$ (automorphisms of $S$, which are isotopic to $\id_S$ via an isotopy which fixes $\qs$). Let us write $\pi\colon \cJ(S)\to \cT(S^{\qs})$ for the natural map.

\begin{define} If $j_0\in \cJ(S)$, a \emph{Teichm\"{u}ller slice} through $j_0$ is a subset $\cT\subset \cJ(S)$ parametrized by a map
\[
\phi\colon (U,0)\to (\cJ(S),j_0)
\]
where $U$ is an open subset of some Euclidean space, containing $0$,  such that the composition $\pi\circ \phi\colon U\to \cT(S^{\qs})$ is a diffeomorphism onto its image. Furthermore, all $j\in \cT$ coincide on a neighborhood of $\ve{q}$.
\end{define}

See \cite{WendlAutomatic}*{Section~3.1} for information on constructing Teichm\"{u}ller slices.

Let $\cB$  denote the set of $W^{p,1}_{\delta}$ maps from $S$ to $\Sigma\times [0,1]\times \R$. Fix a Teichm\"{u}ller slice $\cT$ through $j_0\in \cJ(S)$. There is a Banach space bundle $\cE$ over $\cT\times \cB$, whose fiber over $(j,u)$ consists of $L_{\delta}^p$ sections of the vector bundle  $\Omega^{0,1}(S,j,u)\to S$, whose fiber at $p$ is the vector space
\[
\Omega_p^{0,1}(S,j,u)=\bar{\Hom}_{\C}\left((T_p S, j), ( T_{u(p)}\Sigma \times [0,1]\times \R, J)\right).
\]
We view the operator $\bar{\d}_{J}\colon \cT\times \cB\to \cE,$ given by
\[
\bar{\d}_J(u)=d u+J\circ d u\circ j,
\]
as giving a section of the bundle $\cE$.

Since the almost complex structures $j\in \cT$ on $S$ from a Teichm\"{u}ller slice $\cT$ are constant near $\qs$, we may pick complex coordinates in a neighborhood of each $q_i$. For a tuple $(r_1,\dots, r_k)\in \Z_{\ge 0}^k$, we obtain a jet space evaluation map
\begin{equation}
\ev^{\Jet}\colon \bar{\d}_{J}^{-1}(0) \cap (\cT\times \cB)\to \Jet^{r_1}(\Sigma\times [0,1]\times \R,J)\times \cdots \times \Jet^{r_k}(\Sigma\times [0,1]\times \R,J). \label{eq:jet-evaluation}
\end{equation}

Write $\bar{\d}^{-1}_{J}(0)|_D$ for the set of open subset consisting of marked holomorphic curves whose marked points have evaluation in $D\times [0,1]\times \R$.  Furthermore, we may restrict attention to the set of almost complex structures satisfying \eqref{J1}-\eqref{J5} which additionally coincide with $j_D\times j_{\bD}$ on $D\times [0,1]\times \R$, for some fixed almost complex structure $j_D$ on $D$. In this case, the jet evaluation map from ~\eqref{eq:jet-evaluation} restricts to a map
\begin{equation}
\ev^{\Jet}_D \colon \bar{\d}^{-1}_{J}(0)|_D\to \Jet^{r_1}(D\times [0,1]\times \R)\times \cdots \times \Jet^{r_k}(D\times [0,1]\times \R),
\label{eq:restricted-jet-evaluation-map}
\end{equation}
which has a codomain independent of $J$.

The matched moduli space in~\eqref{eq:def-matched-moduli-space-most-general} is then locally given as the preimage under $\ev^{\Jet}_D$ of a submanifold 
\[
Z(M)\subset \Jet^{r_1}(D\times [0,1]\times \R)\times \cdots \times \Jet^{r_k}(D\times [0,1]\times \R),
\]
where $r_i=\max\{n_i,m_i\}$. We leave it to the reader to write down an explicit formula for $Z(M)$ in terms of $M$, however we note that it has codimension
\[
2k+\codim(M_0)+2|\ve{t}|.
\]

\begin{define}
We say that $J$ is \emph{regular} for $M$ if at every $(u,j)\in \cM_J(S^{\qs}, \phi,M)$  there is a Teichm\"{u}ller slice $\cT$ through $j$, such that the map $\bar{\d}_{J}$ is transverse to the zero section of the bundle $\cE\to \cT\times \cB$, and the map $\ev_D^{\Jet}$ is transverse to $Z(M)$. 
\end{define}

The above regularity condition is independent of the Teichm\"{u}ller slice $\cT$; see \cite{WendlSurvey}*{Lemma~4.3.2}. Finally, an essentially standard argument shows that there is a dense (Baire) subset of almost complex structures satisfying \eqref{J1}--\eqref{J5} which are regular. The argument can be obtained by combining \cite{LipshitzCylindrical}*{Proposition~4.3} together with one of \cite{WendlSuperRigidity}*{Lemma~A.3}, \cite{CieliebakMohnkeGW}*{Lemma~6.6} or \cite{ZehmischJets}. 

Taking the preimage of $Z(M)$ under $\ev^{\Jet}_D$ gives part~\eqref{expected-dimensions-1} of Proposition~\ref{prop:dimension-counts}. Part~\eqref{expected-dimensions-2} of Proposition~\ref{prop:dimension-counts} follows from part~\eqref{expected-dimensions-1}, together with the formula
\begin{equation}
\ind(u)=\mu(\phi)-2\Sing(u). \label{eq:singularity-number-correction}
\end{equation}
Equation~\eqref{eq:singularity-number-correction} follows from \cite{LipshitzErrata}*{Proposition~4.2'}. See also \cite{LOTBordered}*{Proposition~5.69}.

\subsection{Cylindrical boundary degenerations}
\label{sec:cylindrical-boundary}
In this section, we recall some notions about cylindrical boundary degenerations. The material is mostly based on \cite{OSLinks}*{Section~5} and \cite{OSBorderedHFK}*{Section~5.6}, though Proposition~\ref{prop:relaxed-count} is new.

\begin{define}
A \emph{cylindrical beta boundary degeneration} is a smooth map
\[
u\colon (S,\d S)\to (\Sigma\times [0,\infty)\times \R, \bs\times \{0\}\times \R),
\]
satisfying the following:
\begin{enumerate}[ref=$N$-\arabic*, label=($N$-\arabic*),leftmargin=*, widest=IIII]
\item \label{N-1} $u$ is $(j,J)$-holomorphic.
\item \label{N-2} $u$ is proper.
\item \label{N-3} For each $t\in \R$ and $i\in \{1,\dots, n\}$, $u^{-1}(\beta_i\times \{t\})$ consists of a single point.
\item \label{N-4} $u$ has finite energy.
\item \label{N-5} $\pi_{\bH}\circ u$ is non-constant on each component of $S$, where $\pi_{\bH}$ denotes the projection onto $[0,\infty)\times \R$.
\end{enumerate}
We define cylindrical alpha boundary degenerations analogously.
\end{define}

If $B$ is a relative 2-cycle on $(\Sigma,\bs)$, we write $\cN_J(B)$ for the moduli space of cylindrical beta-boundary degenerations representing the homology class $B$. There is an evaluation map
\[
\ev^{\infty}\colon \cN_J(B)\to \bT_{\b},
\]
obtained by evaluating a boundary degeneration at its punctures. If $\ve{x}\in \bT_{\b}$, we write $\cN_J(B,\ve{x})$ for $(\ev^{\infty})^{-1}(\xs)$.

The space $\cN_J(B,\ve{x})$ has an action of $\Aut(\bH)$, the group of conformal autormorphisms of the half plane, which acts pointwise on the $[0,\infty)\times \R$ factor of a curve. We write $\tilde{\cN}_J(B,\ve{x})$ for the quotient
\[
\tilde{\cN}_J(B,\ve{x}):=\cN_J(B,\ve{x})/\Aut(\bH).
\]

\subsection{Boundary degenerations and split almost complex structures}

  Ozsv\'{a}th and Szab\'{o} prove the following:
\begin{thm}[\cite{OSLinks}*{Theorem~5.5}]\label{thm:os-boundary-deg-count} Suppose that $(\Sigma,\as,\bs,\ws)$ is a Heegaard diagram, and $J$ is a generic, split almost complex structure on $\Sigma\times [0,\infty)\times \R$. If $B$ is a Maslov index 2 class of beta boundary degenerations and $\xs\in  \bT_{\b}$, then
\[
\# \tilde{\cN}_{J}(B,\xs)\equiv \begin{cases} 1& \text{if} \quad|\ve{w}|>1\\
0& \text{if}\quad  |\ve{w}|=1.
\end{cases}
\]
\end{thm}

Ozsv\'{a}th and Szab\'{o}'s proof when $|\ws|=1$ is  based on the corresponding result for the symmetric product \cite{OSDisks}*{Theorem~3.15}, via the \emph{tautological correspondence}. This turns out to not be the most natural setting for our purposes, since it does not rule out curves appearing in codimension 1 degenerations which satisfy only \eqref{N-1}--\eqref{N-4}, but not \eqref{N-5}. For example, we may always form a holomorphic curve as the union of $\Sigma\times \{pt\}$ together with a constant curve at $\xs\in \bT_{\b}$. The $\bar{\d}$ operator does not achieve transversality for such curves, as can easily be seen since they violate expected dimension counts; compare Proposition~\ref{prop:dimension-counts}.

In Sections~\ref{sec:generic-boundary-degenerations} and ~\ref{sec:counting-cylindrical-boundary}, we adapt Lipshitz's proof of stabilization invariance \cite{LipshitzCylindrical}*{Section~12} and consider more general families of almost complex structure where, generically, transversality is achieved for all curves appearing in codimension 1 degenerations. Note that when we do this, we get a different count for the case $|\ws|=1$; see Proposition~\ref{prop:relaxed-count}.

\subsection{Boundary degenerations and generic almost complex structures}
\label{sec:generic-boundary-degenerations}
In this section, we define the almost complex structures we consider for boundary degenerations. Suppose that $p_1,\dots, p_n$ is a (possibly empty) collection of distinct points in $\Sigma\setminus \bs$, and $D_1,\dots, D_n$ is a collection of small disks in $\Sigma\setminus \bs$, containing $p_1,\dots, p_n$, respectively.

\begin{define} We say an almost complex structure $J$ on $\Sigma\times [0,\infty)\times \R$ is \emph{admissible} if it satisfies the following:
\begin{enumerate}[ref=$J'$-\arabic*, label=($J'$-\arabic*),leftmargin=*, widest=IIII]
\item\label{J'1} $J$ is tamed by $\omega$.
\item\label{J'2} $J$ is split on  $(D_1\cup \cdots \cup D_n)\times [0,\infty)\times \R$.
\item\label{J'3} $J$ is invariant under the action of $\Aut(\bH)$.
\item\label{J'4} $J \d/\d s=\d/\d t$.
\item\label{J'5} $T_p \Sigma \times \{0\}$ is a complex line in a neighborhood of $\bs\times [0,\infty)\times \R$, and along $\Sigma\times \{0\}\times \R$.
\end{enumerate}
\end{define}

If $J$ is an almost complex structure on $\Sigma\times [0,1]\times \R$, satisfying \eqref{J1}--\eqref{J5}, then there is a naturally associated almost complex structure $J'$ on $\Sigma\times [0,\infty)\times \R$, as we now describe. Fix $y_0\in \{0\}\times \R$. If $r>0$, write
\[
R_{r,y_0}\colon \Sigma\times [0,1]\times \R\hookrightarrow \Sigma\times [0,\infty)\times \R
\]
for the embedding $R_{r,y_0}(x,z)=(x,r\cdot (z-y_0))$. We define
\begin{equation}
J'=\lim_{r\to \infty} d R_{r,y_0}\circ J\circ d (R_{r,y_0})^{-1}. \label{eq:infinitesimal-ac}
\end{equation}
 Note that $J'$ is independent of $y_0$.

There is a more concrete description of $J'$ (which in particular makes it clear that the limit exists). With respect to the isomorphism
\[
T_{(p,y)}(\Sigma\times [0,1]\times \R)\iso T_p\Sigma \oplus T_y ([0,1]\times \R),
\]
we may decompose
\[
J_{(p,y)}=\begin{pmatrix} (j_{\Sigma})_{(p,y)}& 0\\
Y_{(p,y)}& j_{\bD}
\end{pmatrix},
\]
where the following hold:
\begin{itemize}
\item $(j_{\Sigma})_{(p,y)}$ is an almost complex structure on $T_{p} \Sigma$,
\item $j_{\bD}$ is the standard almost complex structure on $[0,1]\times \R$,
\item $Y_{(p,y)}\colon T_{p} \Sigma\to T_y([0,1]\times \R)$ is complex anti-linear.
\end{itemize}
Equation~\eqref{eq:infinitesimal-ac} has the more concrete description
\begin{equation}
J'_{(p,y)}=\begin{pmatrix} (j_{\Sigma})_{(p,y_0)} &0\\
(\d_{y} Y)_{(p,y_0)}& j_{\bH}
\end{pmatrix},
\label{eq:matrix-derivative}
\end{equation}
where $\d_y Y$ denotes the directional derivative, in the direction $y$.

Finally, we note that if a boundary degeneration $u\colon S\to \Sigma \times [0,\infty)\times \R$ appears in the limit of $J$-holomorphic curves on $\Sigma\times [0,1]\times \R$, then $u$ will be $J'$-holomorphic. This holds because a boundary degeneration which appears in the limit of a sequence of curves in $\Sigma\times [0,1]\times \R$ is obtained by rescaling the $[0,1]\times \R$ component. Hence, the description in~\eqref{eq:infinitesimal-ac} makes it clear that the limiting curve will be $J'$-holomorphic.

\subsection{Boundary degenerations for generic almost complex structures}
\label{sec:counting-cylindrical-boundary}

In this section, we count index 2 boundary degenerations, for almost complex structures on $\Sigma\times [0,\infty)\times \R$ which satisfy \eqref{J'1}--\eqref{J'5}; compare Theorem~\ref{thm:os-boundary-deg-count}.

A key observation is that for a generic almost complex structure on $\Sigma\times [0,\infty)\times \R$, satisfying \eqref{J'1}--\eqref{J'5}, the $\bar{\d}$ operator achieves transversality at holomorphic maps satisfying~\eqref{N-1}--\eqref{N-4}, which have no multiply covered, closed components. The proof may be easily adapted from Lipshitz's proof of the analogous statement for holomorphic disks \cite{LipshitzCylindrical}*{Proposition~3.7}.

\begin{figure}[H]
	\centering
	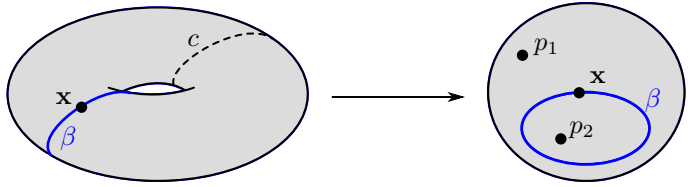
	\caption{Stretching the neck along $c$ identifies $\cN(B,\ve{x})$ with the space of self-matched, Maslov index 4 boundary degenerations on $(S^2,\b)$.}\label{fig:131}
\end{figure}

We now state our main count:

\begin{prop}\label{prop:relaxed-count}Suppose $\bs\subset \Sigma$ is a collection of attaching curves on $\Sigma$, $B$ is a Maslov index 2 class of beta boundary degenerations, and $\ve{x}\in \bT_{\b}$. Then for generic $J$ satisfying ~\eqref{J'1}--\eqref{J'5}, 
\begin{equation}
\# \tilde{\cN}_J(B,\ve{x})\equiv 1\pmod 2, \label{eq:count-boundary-degenerations}
\end{equation}
(regardless of whether $|\ve{w}|=1$ or $|\ve{w}|>1$).
\end{prop}
\begin{proof} First, we argue that the mod 2 count of $\tilde{\cN}_J(B,\ve{x})$ is independent of  the choice of generic $J$. Indeed if $J_0$ and $J_1$ are two generic almost complex structures, then one picks a generic family, $(J_t)_{t\in [0,1]}$, connecting $J_0$ and $J_1$. We consider the parametrized moduli space
\[
\tilde{\cN}_{J_t, t\in [0,1]}(B,\xs):=\bigcup_{t\in [0,1]} \tilde{\cN}_{J_t}(B,\ve{x})\times \{t\}.
\]
We claim the ends of $\tilde{N}_{J_t, t\in [0,1]}(B,\xs)$ are exactly the curves at $t=0$ and $t=1$. To see this, we consider the possible degenerations away from $t=0$ or $t=1$. The source may degenerate along  a collection of closed circles and boundary-to-boundary arcs. If the source degenerates along a boundary-to-boundary arc, then \eqref{N-3} implies that the limit contains either an additional boundary degeneration, or a \emph{ghost curve} component (i.e. a component which has image equal to a single point). The formation of an additional boundary degeneration is prohibited since every non-trivial boundary degeneration has Maslov index at least 2, so the formation of an additional boundary degeneration would force $B$ to have Maslov index at least 4. 

For the moment, consider the components remaining after we trim any ghost components. If there are additional degenerations of the source, then the remaining components of the source have nodal points in their interior. The existence of such curves are prohibited by the analog of Proposition~\ref{prop:dimension-counts} for boundary degenerations. In particular, the expected dimension over all $t$ of the space of non-embedded curves representing $B$ is at most 1. After quotienting by $\Aut(\bH)$, we obtain an expected dimension of $-1$, so a generic family $J_t$ misses the almost complex structures which support such curves. 

Finally, it remains to rule out ghost curves. These are ruled out by an easy adaptation of a standard Fredholm index argument, for which we refer the reader to \cite{LOTBordered}*{Lemma~5.57} for more details.

Next, it remains to establish the stated count of boundary degenerations. We begin with the $g=1$ case, when $(\Sigma,\bs)=(\bT^2,\b)$. In this case, we take an embedded, homologically essential closed curve $c\subset \bT^2\setminus \b$. We stretch the neck along $c$. If we cut $\bT^2$ along $c$, and collapse the two boundary components, we are left with a 2-sphere, with two special points, $p_1$ and $p_2$. For large neck-length parameter, by gluing, we may identify $\tilde{\cN}_{J(T)}(B,\ve{x})$ with the space of \emph{self-matched boundary degenerations} on $(S^2,\b)$, of Maslov index 4. A self-matched boundary degeneration is a boundary degeneration $u\colon S\to S^2\times [0,\infty)\times \R$, with two marked points, $q_1,q_2\in S$, satisfying
 \[
(\pi_\Sigma \circ u)(q_1)=p_1\quad \text{and}\quad (\pi_\Sigma \circ u)(q_2)=p_2\quad  \text{and}\quad (\pi_\bH\circ u)(q_1) = (\pi_{\bH}\circ u)(q_2).
\]
See Figure~\ref{fig:131}.

We deform the space of self-matched boundary degenerations and consider, for $t\in (0,\infty)$, the space of \emph{$t$-self-matched} cylindrical boundary degenerations, which are Maslov index 4 boundary degenerations with two marked points, $q_1$ and $q_2$ on the source, which satisfy
 \[
(\pi_\Sigma \circ u)(q_1)=p_1\quad \text{and}\quad (\pi_\Sigma\circ u)(q_2)=p_2 \quad  \text{and}\quad   (\pi_\bH \circ u)(q_1)=t\cdot (\pi_{\bH}\circ u)(q_2).
\]
Here, if $z\in [0,\infty)\times \R$, then $t\cdot z\in [0,\infty)\times \R$ denotes the ordinary product, viewing $[0,\infty)\times \R$ as the right-half plane in $\C$. The ordinary self-matched moduli space coincides with the 1-matched moduli space. Write $\cN^{t}(B,\ve{x})$ for the moduli space of $t$-self-matched boundary degenerations, which we view as a subspace of $\cN(B_1\cup B_2,\ve{x})$, where $B_1$ and $B_2$ denote the two components of $S^2\setminus \b$. Note that the action of $\Aut(\bH)$ preserves $\cN^t(B,\ve{x})$ for each $t$.

Gluing identifies the moduli space $\tilde{\cN}^t(B,\ve{x})$ with the cartesian product
\[
\tilde{\cN}(B_1,\ve{x})\times \tilde{\cN}(B_2,\ve{x}),
\]
when $t$ is large, where $B_1$ and $B_2$ are the index 2 classes which cover $p_1$ and $p_2$, respectively. Hence, it remains to count $\tilde{\cN}(B_i,\ve{x})$ individually. These may be counted by arguing the same as \cite{OSLinks}*{Theorem~5.5}. Namely, we deform the almost complex structure so that it becomes a product, and then we are counting holomorphic maps from $\bH$ to $\bH\times [0,\infty)\times \R$, which are the identity on $\bH$, modulo the action of $\Aut(\bH)$. The count of such curves is  clearly 1.

The case when $|\ws|=1$, but $g(\Sigma)>1$ is proven repeatedly performing a neck-stretching argument, similar to the above one, to reduce to the $g=1$ case (compare the proof of \cite{OSLinks}*{Theorem~5.5}).
\end{proof}

\begin{rem} 
Proposition~\ref{prop:relaxed-count} gives an alternative to the last step in Lipshitz's proof of stabilization invariance. In his proof, Lipshitz proves several gluing results \cite{LipshitzCylindrical}*{Appendix~A} which reduce the question of stabilization invariance to the claim that
\[
\# \cM(\phi, d)\equiv 1 \pmod{2},
\] 
where $\phi$ is the Maslov index 2 class on a standard genus 1 diagram for $S^3$,  $d\in (0,1)\times \R$. This count may instead be performed by letting $d$ approach a point in $d_0\in \{1\}\times \R$, along some path $\g$. By considering the moduli space of curves which match any point along $\g$, one identifies $\cM(\phi,d)$ with the count of beta-boundary degenerations, which is 1, according to our Proposition~\ref{prop:relaxed-count}.
\end{rem}

Another situation that is helpful to understand is that of \emph{matched boundary degenerations} on pairs of diagrams. Suppose that $\cH_1=(\Sigma_1,\bs_1)$ and $\cH_2=(\Sigma_2,\bs_2)$ are two partial Heegaard diagrams, each with a special point $p_i\in \Sigma_i\setminus \bs_i$.  Suppose $\xs_1\in  \bT_{\b_1}$ and $\xs_2\in  \bT_{\b_2}$ are two points. There are unique, non-negative, index 2 classes, $B_1$ and $B_2$, which cover $p_1$ and $p_2$ once, respectively.  For $t\in \R$, we consider the moduli space $\cN\cN^t(B_1,B_2,\xs_1,\xs_2)$ which consists of pairs $(u_1,u_2)$ of boundary degenerations, each with a marked point, $q_1$ and $q_2$, such that
\begin{equation}
\begin{split}
(\pi_{\Sigma_1}\circ u_1)(q_1)=p_1,& \quad (\pi_{\Sigma_2}\circ u_2)(q_2)=p_2, \quad
(\pi_{\bH}\circ u_1)(q_1)= (\pi_{\bH}\circ u_2)(q_2),\\
  &\ev^{\infty}(u_1)=\xs_1\quad \text{and} \quad \ev^{\infty}(u_2)=\xs_2.
\end{split}
\label{eq:matched-boundary-degenens}
 \end{equation}
We write $\tilde{\cN\cN}^t(B_1, B_2,\xs_1,\xs_2)$ for the quotient by the action of $\Aut(\bH)$.

\begin{lem}\label{lem:t-matched-boundary-degenerations}
Suppose $(\Sigma_1,\bs_1)$ and $(\Sigma_2,\bs_2)$ are two partial Heegaard diagrams, each with a special point $p_i\in \Sigma_i\setminus \bs_i$, as above. If $B_1$ and $B_2$ are classes described above, and $J$ is a generic almost complex structure on $(\Sigma_1\sqcup \Sigma_2)\times [0,\infty)\times \R$ which satisfies \eqref{J'1}--\eqref{J'5}, then
\[
\# \tilde{\cN\cN}_J(B_1,B_2,\xs_1,\xs_2)\equiv 1.
\]
The same also holds for alpha boundary degenerations.
\end{lem}
\begin{proof} The argument follows from the same line of reasoning as Proposition~\ref{prop:relaxed-count}. For $t\in [1,\infty)$, we consider the moduli space of $t$-matched boundary degenerations $\cN\cN^t(B_1,B_2,\xs_1,\xs_2)$, which are defined as in~\eqref{eq:matched-boundary-degenens}, except we require
\[
(\pi_{\bH}\circ u_1)(q_1)=t\cdot (\pi_{\bH}\circ u_2)(q_2),
\]
where $t\cdot$ denotes the standard action of $\R^{>0}$ on $\bH$. We may view $\cN\cN^t(B_1,B_2,\xs_1,\xs_2)$ as being a subspace of the four dimensional moduli space $\cN(B_1,B_2,\xs_1,\xs_2)$. We note that the action of $\Aut(\bH)$ preserves the subspace $\cN\cN^t(B_1,B_2,\xs_1,\xs_2)$, for all $t$, so we may consider the quotient $\tilde{\cN\cN}^t(B_1,B_2,\xs_1,\xs_2)$.
Next, we consider the parametrized moduli space 
\[
\tilde{\cN}=\bigcup_{t\in [1,\infty)} \tilde{\cN\cN}^t(B_1,B_2,\xs_1,\xs_2)\times \{t\}.
\]
One end of $\tilde{\cN}$ may be identified with $\tilde{\cN\cN}(B_1,B_2,\xs_1,\xs_2)$, while the other may be identified with $\tilde{\cN}(B_1,\xs_1)\times \tilde{\cN}(B_2,\xs_2)$. Additional ends are ruled out as in Proposition~\ref{prop:relaxed-count}, completing the proof.
\end{proof}

\subsection{Holomorphic polygons and related constructions}

In this section, we describe the holomorphic polygon counting maps, and their friends, in the cylindrical setting. We write $D_n$ for a disk with $n$ boundary punctures, which we view as an $n$-gon with cylindrical ends.

\begin{define}\label{def:stratified-decomposing-map}
Suppose that $\cA=\{a_1,\dots, a_k\}$ is a collection of compact, properly embedded embedded arcs in $D_n$. We say that $\cA$ is a collection of \emph{decomposing arcs} if no $a_i$ has its endpoints on the same component of $\d D_n$. Let $2^{\cA}$ denote the power set of $\cA$. Suppose $X$ is a smooth, compact manifold with corners. We say $f\colon X\to 2^{\cA}$ is a \emph{stratified decomposing map} if it satisfies the following:
\begin{enumerate}[ref=$s$-\arabic*, label=($s$-\arabic*),leftmargin=*, widest=IIII]
\item $f$ is constant on the interior of each component of the codimension $l$ strata of $X$, for all $l$.
\item $f(x)\equiv \emptyset$ on the codimension 0 strata.
\item If $h$ is a component of the codimension $l$ strata of $X$, $x\in h$ and $y\in \d h$, then $f(x)\subset f(y)$.
\item If $f(x)=\{a_{i_1},\dots, a_{i_l}\}$, then $a_{i_1},\dots, a_{i_l}$ are pairwise disjoint.
\end{enumerate}
\end{define}

If $\cA=\{a_1,\dots, a_k\}$ is a collection of decomposing arcs in $D_n$, we pick regular neighborhoods $N(a_1),\dots, N(a_k)$, equipped with identifications $N(a_i)\iso a_i\times [-1,1]$. Further, we assume that if $a_i$ and $a_j$ are disjoint, then $N(a_i)$ and $N(a_j)$ are as well. Additionally, pick a smooth function $\rho\colon [-1,1]\to [0,1]$ such that 
\begin{enumerate}
\item $\rho(x)=0$ for $x$ in a neighborhood of $-1$ and $1$.
\item $\rho(0)=1$, and $\rho(x)<1$ for $x\neq 0$.
\end{enumerate}

We encourage the reader to compare the following to \cite{LOTDoubleBranchedII}*{Definition~3.5}:

\begin{define}\label{def:stratified-family}
 Suppose $\cA$ is a set of decomposing arcs for $D_n$, as above, $X$ is a smooth manifold with corners, and $f\colon X\to 2^{\cA}$ is a stratified decomposing map. We say  $(J_x)_{x\in X}$ is a \emph{stratified family of almost complex structures} on $\Sigma\times D_n$ if the following hold:
 \begin{enumerate}[ref=$J''$-\arabic*, label=($J''$-\arabic*),leftmargin=*, widest=IIIII]
 \item \label{J1''} For each $x\in X$, $J_x$ is an almost complex structure on $\Sigma\times (D_n\setminus f(x))$, which is tamed by the split symplectic form on $\Sigma\times (D_n\setminus f(x))$.
 \item \label{J2''} If $h$ is a component of the codimension $l$ strata of $X$, and $D_m$ is a component of $D_n\setminus f(x)$,  then $J_x$ is constant on the cylindrical ends of $D_m$, for all $x\in  h$. Furthermore, the induced cylindrical almost complex structure on $\Sigma\times [0,1]\times \R$ satisfies \eqref{J1}--\eqref{J5}.
\item\label{J4''} $J_x$ preserves $T D_n$ for all $x\in X$. In particular, there is an induced family $(j_{D_n,(x,p)})_{(x,p)\in X\times \Sigma}$ of almost complex structures on $D_n$. We assume the family $j_{D_n,(x,p)}$ is constant in $p$.
\item\label{J5''} $T_p\Sigma \oplus \{0\}$ are  complex lines on a neighborhood of $(\as\cup \bs)\times D_n$, and along $\Sigma\times \d D_n$.
 \item\label{J6''} Roughly, if $x\in X$ and $f(x)=\{a_{i_1},\dots, a_{i_k}\}$, then for $y$ near $x$, $J_y$ is obtained by inserting long necks along $a_{i_1},\dots, a_{i_k}$, consistently with the stratified map $f$. Precisely, there is a neighborhood $U\subset X$ of $x$  and a collection of smooth functions $\tau_1,\dots, \tau_k\colon U\to [0,1]$, satisfying the following:
 \begin{enumerate} 
 \item For all $y\in U$, $\tau_j(y)=1$ if and only if $a_{i_j}\in f(y)$. 
 \item On each $N(a_{i_j})$,
 \[
	J_y \frac{\d}{\d s}=(1-\tau_j(y)\cdot \rho (t))\frac{\d}{\d t},
 \] 
 where $(s,t)$ denote the coordinates on $N(a_{i_j})\iso a_{i_j}\times [-1,1]$.
 \item  For each $y\in U$, $J_y$ is invariant on $\Sigma \times N(a_{i_j})$ under the flow of 
 \[
 V_y=(1-\tau_j(y)\cdot \rho(t))\frac{\d}{\d t},
 \]
 and $(J_y)|_{\Sigma\times a_i\times \{-1,1\}}$ is constant for $y\in U$.
 \item If $x_n\to x$, then $J_{x_n}$ approaches $J_x$ on $\Sigma \times D_n\setminus (N(a_{i_1})\cup \cdots \cup  N(a_{i_k}) )$ in the $C^\infty$ topology. 
 \end{enumerate}
 \end{enumerate}
\end{define}

\begin{figure}[ht!]
	\centering
\begingroup%
  \makeatletter%
  \providecommand\color[2][]{%
    \errmessage{(Inkscape) Color is used for the text in Inkscape, but the package 'color.sty' is not loaded}%
    \renewcommand\color[2][]{}%
  }%
  \providecommand\transparent[1]{%
    \errmessage{(Inkscape) Transparency is used (non-zero) for the text in Inkscape, but the package 'transparent.sty' is not loaded}%
    \renewcommand\transparent[1]{}%
  }%
  \providecommand\rotatebox[2]{#2}%
  \newcommand*\fsize{\dimexpr\f@size pt\relax}%
  \newcommand*\lineheight[1]{\fontsize{\fsize}{#1\fsize}\selectfont}%
  \ifx\svgwidth\undefined%
    \setlength{\unitlength}{198.25301962bp}%
    \ifx\svgscale\undefined%
      \relax%
    \else%
      \setlength{\unitlength}{\unitlength * \real{\svgscale}}%
    \fi%
  \else%
    \setlength{\unitlength}{\svgwidth}%
  \fi%
  \global\let\svgwidth\undefined%
  \global\let\svgscale\undefined%
  \makeatother%
  \begin{picture}(1,1.04671342)%
    \lineheight{1}%
    \setlength\tabcolsep{0pt}%
    \put(0,0){\includegraphics[width=\unitlength,page=1]{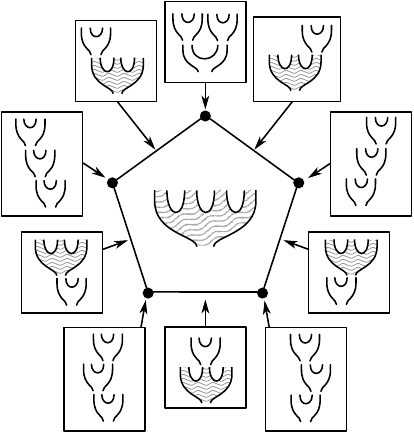}}%
    \put(0.38258852,0.38916607){\color[rgb]{0,0,0}\makebox(0,0)[lt]{\lineheight{1.25}\smash{\begin{tabular}[t]{l}$K_4$\end{tabular}}}}%
  \end{picture}%
\endgroup%

	\caption{A stratified family of almost complex structures indexed by the assocediahedron $K_4$ (a pentagon). We have schematically illustrated how the almost complex structure degenerates on the codimension  1 and 2 strata. The wavy lines indicate where the almost complex structure is allowed to vary along the codimension 0 and 1 strata.}\label{fig:158}
\end{figure}

The main example is when $X$ is an associahedron. The example that $X=K_4$ is shown in Figure~\ref{fig:158}. Note that the associahedron is not the only example we need; see for example Figure~\ref{fig:111}.

\begin{define}\label{def:split}
 We say that a stratified family of almost complex structures $(J_x)_{x\in X}$ on $\Sigma\times D_n$ is \emph{split} if $J_x=j_\Sigma\times j_{D_n,x}$, for a fixed almost complex structure $j_\Sigma$ on $\Sigma$, and a family $(j_{D_n,x})_{x\in X}$ on $D_n$.
\end{define}

If $J=(J_x)_{x\in X}$ is a stratified family of almost complex structures on $\Sigma\times D_n$, and $S$ is a source curve for $n$-gons, we define the moduli space $\cM_{J}(S,\phi)$ to be the set of tuples $(x,u,j)$, where $x$ lies in the codimension 0 strata, $j\in \cJ(S)$ and $u\colon S\to \Sigma\times D_n$ is $(j,J_x)$-holomorphic and satisfies the analogs of \ref{M-2}--\ref{M-7}.

If $\cD=(\Sigma,\gs_1,\dots,\gs_n,w)$ is a Heegaard $n$-tuple, then (assuming $\cD$ is weakly admissible), there is a cylindrical holomorphic polygon counting map, obtained by picking a stratified family $(J_x)_{x\in X}$, where $X$ is the Deligne-Mumford compactification of the moduli space of $n$ marked points on a disk (or equivalently, $X$ is the associahedron $K_{n-1}$) and defining
\[
F_{\g_1,\dots, \g_n;J}\colon \bCF^-(\gs_1,\gs_2)\otimes \cdots \otimes \bCF^-(\gs_{n-1}, \gs_n)\to \bCF^-(\gs_1,\gs_n),
\]
to count curves of Maslov index $3-n$ mapping into $\Sigma\times D_n$ which are $J_x$-holomorphic for some $x\in \Int(X)$. The motivating example is when $X$ is an associahedron, however we also need some variations of this. See Figure~\ref{fig:111}.

\subsection{Compactifications via holomorphic combs}

We now describe a compactification of the moduli spaces in Section~\ref{sec:moduli-spaces} which provides a helpful framework for some of the arguments of this paper (and is used implicitly in Section~\ref{sec:counting-cylindrical-boundary}). It is based on the compactness results of symplectic field theory \cite{EGH-SFT} \cite{BEHWZCompactness}, as well as several related notions in the setting of bordered Floer homology (see especially \cite{LOTBordered}*{Section~5.4} \cite{OSBorderedHFK}*{Section~7.7}).

\begin{define}\label{def:holomorphic-comb}
Suppose that $(\Sigma,\as,\bs,w)$ is a Heegaard diagram of genus $g$.
\begin{enumerate}
\item Let $u$ be a holomorphic curve in $\Sigma\times [0,1]\times \R$ and $\ve{P}$ be a set of boundary punctures on the source of $u$. We say that $\ve{P}$ is an \emph{extended fiber of $u$ over $\d \bD$} if $\ve{P}$ consists of $g$ punctures $p_1,\dots, p_g$, such that $\lim_{z\to p_i}u(z)$ exists and is in $\Sigma\times \{0,1\}\times \R$ for all $i$, and the projection to $[0,1]\times \R$ is independent of $i$. Furthermore, we assume that no two distinct $p_i$ have limit in the same $\as$ or $\bs$ curve. We also allow the case that $\lim_{z\to p_i} (\pi_{\bD}\circ u)(z)=\pm \infty$.

 For a holomorphic curve mapping into $\Sigma\times \bH$, we say a set of $g$ boundary punctures is an \emph{extended fiber over $\d \bH$} if the analogous conditions hold, though we also allow the case that $\lim_{z\to p_i} (\pi_{\bH}\circ u(z))=\infty$.
\item A \emph{holomorphic story} is a collection of holomorphic curves 
\[
\{b_n,\dots, b_1, u,a_1,\dots, a_m\},
\]
 as follows:
\begin{enumerate}
\item For each $i\in \{1,\dots, m\}$, $a_i$ is a finite collection of alpha-boundary degenerations on $(\Sigma,\as)$, equipped with two sets of boundary punctures, $\ve{A}^+(a_i)$ and  $\ve{A}^-(a_i)$.
\item For each $i\in \{1,\dots, n\}$, $b_i$ is a finite collection of beta-boundary degenerations on $(\Sigma,\bs)$, equipped with two sets of boundary punctures $\ve{B}^+(b_i)$ and $\ve{B}^-(b_i)$.
\item $u$ is a holomorphic curve in  $\Sigma\times [0,1]\times \R$. The source of $u$ has two $g$-tuples of  boundary punctures labeled $+\infty$ and $-\infty$, which are asymptotic to intersection points $\xs\times \{-\infty\}$ and $\ys\times \{+\infty\}$ in $(\bT_{\a}\cap \bT_{\b})\times \bar{\bD}$. The source of $u$ also has two sets of boundary punctures, $\ve{A}^-(u)$ and $\ve{B}^-(u)$. After  completing over $\ve{A}^-(u)$ and $\ve{B}^-(u)$, the curve $u$ satisfies \ref{M-1}--\ref{M-6}.
\end{enumerate}
Additionally:
\begin{enumerate}
\item The set $\ve{A}^-(u)$ is  partitioned into extended $\d \bD$ fibers, which all project to $\{1\}\times \R$, and have distinct projection thereto. The set $\ve{B}^-(u)$ are similarly partitioned into extended $\d \bD$ fibers which map to $\{0\}\times \R$, and have distinct projection.
\item For each $a_i$, the punctures $\ve{A}^+(a_i)$ are partitioned into extended $\d \bH$ fibers (one for each boundary degeneration of $a_i$), and the projection $\pi_{\bH}\circ a_i$ is asymptotic to $\infty$ at each puncture. The same holds for the punctures $\ve{B}^+(b_i)$ on $b_i$.
\item  The punctures $\ve{A}^-(a_i)$ are partitioned into extended $\d \bH$ fibers. Furthermore, any two extended $\d \bH$ fibers on the same boundary degeneration have distinct projection to $\d \bH$. The analogous statement holds for each $\ve{B}^-(b_i)$.
\item There is a chosen bijection between  $\ve{A}^-(u)$ and  $\ve{A}^+(a_1)$, which is compatible with the two partitions into extended fibers. Furthermore, each pair of extended fibers which are matched, have equal evaluation in $\bT_{\a}$. A similar identification and matching condition holds between $\ve{A}^-(a_{i+1})$ and $\ve{A}^+(a_i)$. There are analogous identifications and matching conditions between $\ve{B}^-(u)$ and $\ve{B}^+(b_{1})$, and also between $\ve{B}^-(b_i)$ and $\ve{B}^+(b_{i+1})$.
\end{enumerate}
\item A \emph{holomorphic comb} is a sequence of holomorphic stories $\cU_1,\dots, \cU_n$, such that the intersection point at $+\infty$ of $\cU_i$ coincides with the intersection point at $-\infty$ of $\cU_{i+1}$.
\end{enumerate}
\end{define}

\begin{figure}[H]
	\centering
\begingroup%
  \makeatletter%
  \providecommand\color[2][]{%
    \errmessage{(Inkscape) Color is used for the text in Inkscape, but the package 'color.sty' is not loaded}%
    \renewcommand\color[2][]{}%
  }%
  \providecommand\transparent[1]{%
    \errmessage{(Inkscape) Transparency is used (non-zero) for the text in Inkscape, but the package 'transparent.sty' is not loaded}%
    \renewcommand\transparent[1]{}%
  }%
  \providecommand\rotatebox[2]{#2}%
  \newcommand*\fsize{\dimexpr\f@size pt\relax}%
  \newcommand*\lineheight[1]{\fontsize{\fsize}{#1\fsize}\selectfont}%
  \ifx\svgwidth\undefined%
    \setlength{\unitlength}{117.94999129bp}%
    \ifx\svgscale\undefined%
      \relax%
    \else%
      \setlength{\unitlength}{\unitlength * \real{\svgscale}}%
    \fi%
  \else%
    \setlength{\unitlength}{\svgwidth}%
  \fi%
  \global\let\svgwidth\undefined%
  \global\let\svgscale\undefined%
  \makeatother%
  \begin{picture}(1,0.50023309)%
    \lineheight{1}%
    \setlength\tabcolsep{0pt}%
    \put(0,0){\includegraphics[width=\unitlength,page=1]{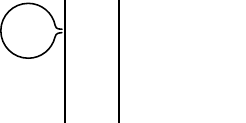}}%
    \put(0.36294313,0.17177405){\color[rgb]{0,0,0}\makebox(0,0)[t]{\lineheight{1.25}\smash{\begin{tabular}[t]{c}$u$\end{tabular}}}}%
    \put(0.11145942,0.36176845){\color[rgb]{0,0,0}\makebox(0,0)[t]{\lineheight{1.25}\smash{\begin{tabular}[t]{c}$b_1$\end{tabular}}}}%
    \put(0,0){\includegraphics[width=\unitlength,page=2]{fig200.pdf}}%
    \put(0.62812286,0.21236923){\color[rgb]{0,0,0}\makebox(0,0)[t]{\lineheight{1.25}\smash{\begin{tabular}[t]{c}$a_1$\end{tabular}}}}%
    \put(0,0){\includegraphics[width=\unitlength,page=3]{fig200.pdf}}%
    \put(0.87222602,0.33000765){\color[rgb]{0,0,0}\makebox(0,0)[t]{\lineheight{1.25}\smash{\begin{tabular}[t]{c}$a_2$\end{tabular}}}}%
    \put(0.88161889,0.07052181){\color[rgb]{0,0,0}\makebox(0,0)[t]{\lineheight{1.25}\smash{\begin{tabular}[t]{c}$a_2$\end{tabular}}}}%
    \put(0,0){\includegraphics[width=\unitlength,page=4]{fig200.pdf}}%
    \put(0.1112845,0.10606833){\color[rgb]{0,0,0}\makebox(0,0)[t]{\lineheight{1.25}\smash{\begin{tabular}[t]{c}$b_1$\end{tabular}}}}%
    \put(0,0){\includegraphics[width=\unitlength,page=5]{fig200.pdf}}%
  \end{picture}%
\endgroup%

	\caption{A holomorphic comb with one story.}\label{fig:200}
\end{figure}

\begin{rem}
\begin{enumerate}
\item Definition~\ref{def:holomorphic-comb} admits a natural extension for holomorphic polygons mapping into $\Sigma\times D_n$. The only complication is that one also has to keep track of degenerations of $D_n$, via a stratified family of almost complex structures on $\Sigma\times D_n$.
\item Definition~\ref{def:holomorphic-comb} admits a natural extension to holomorphic combs with marked points.
\item The set of holomorphic combs (with or without marked points) possesses a natural topology, similar to the one described in \cite{BEHWZCompactness}.  See also \cite{LOTBordered}*{Definition~5.21}.
\end{enumerate}
\end{rem}

\section{Degenerating connected sum tubes}
\label{sec:degeneratingconnectedsumtubes}

A common tool in gauge theory is the neck-stretching argument. In our paper, it becomes helpful to allow the neck-length parameter to go to $\infty$, and work with moduli spaces which are defined as a fibered product over an evaluation map. In this section, we describe some of these techniques.

\subsection{Connected sums and matched moduli spaces}
\label{sec:matched-disks}
Suppose $\cH_1=(\Sigma_1,\as_1,\bs_1,w_1)$ and $\cH_2=(\Sigma_2,\as_2,\bs_2,w_2)$ are two Heegaard diagrams, and suppose that $J_1$ and $J_2$ are almost complex structures on $\Sigma_1\times [0,1]\times \R$ and $\Sigma_2\times [0,1]\times \R$, which are split in a neighborhood of two chosen points $p_1\in \Sigma_1$ and $p_2\in \Sigma_2$. In this section, we focus on the case that $p_1$ and $p_2$ are immediately adjacent to $w_1$ and $w_2$.

Suppose $\phi_1\in \pi_2(\xs_1,\ys_1)$ and $\phi_2\in \pi_2(\xs_2,\ys_2)$ are homotopy classes of disks on $\cH_1$ and $\cH_2$, satisfying
\[
n_{p_1}(\phi_1)=n_{p_2}(\phi_2).
\]
Let us write $k$ for the common value. Suppose $S_1^{\qs_1}$ and $S_2^{\qs_2}$ are Riemann surfaces, each decorated with an ordered $k$-tuple of marked points. We  define the \emph{matched moduli space of disks}, representing $\phi_1$ and $\phi_2$, to be
\begin{equation}
\begin{split}
&\cM\cM_{J_1\wedge J_2}(S_1^{\qs_1}, S_2^{\qs_2}, \phi_1,\phi_2)\\
:=&\left\{(u_1,u_2)\in \cM_{J_1}(S_1^{\qs_1},\phi_1)\times \cM_{J_2}(S_2^{\qs_2}, \phi_2)\middle \vert \begin{array}{c} (\pi_\Sigma\circ u_1)(q_{1,i})=p_1, \\ (\pi_\Sigma\circ u_2)(q_{2,i})=p_2 , \\
(\pi_{\bD}\circ u_1)(q_{1,i}) =(\pi_{\bD}\circ u_2)(q_{2,i}),\\
 \text{ for } i=1,\dots, k
 \end{array}\right\}.
 \end{split}
 \label{eq:matched-moduli-space-disks}
\end{equation}
The construction in ~\eqref{eq:matched-moduli-space-disks} is a special case of the construction from Section~\ref{sec:matched-moduli-spaces}, for a matching condition $M_{\#}$ at the points $p_1$ and $p_2$. In this case, the associated submanifold $M_0$ is the diagonal 
\[
\Delta_{([0,1]\times \R)^k} \subset ([0,1]\times \R)^{2k}.
\]
 The tuple of tangencies $\ve{t}$ is $\{0\}\times \cdots \times  \{0\}$. We note that $M_\#$ has  a dependence on $k$, though it is convenient to suppress this from the notation.

By Proposition~\ref{prop:dimension-counts}, for generic $J_1$ and $J_2$,  the space $\cM\cM_{J_1\wedge J_2}(S_1^{\qs_1}, S_2^{\qs_2},\phi_1,\phi_2)$ is a smooth manifold of dimension
\[
\ind_{\emb}(\phi_1,\phi_2,M_{\#})=\mu(\phi_1)+\mu(\phi_2)-2n_{p_i}(\phi_i),
\]
near a pair $(u_1,u_2)$ of embedded curves.

When no confusion may arise, we write $\cM\cM_{J_1\wedge J_2}(\phi_1,\phi_2)$ for the disjoint union of the matched moduli spaces over all equivalence classes of pairs of marked sources $S_1^{\qs_1}$ and $S_2^{\qs_2}$, where $|\qs_1|=|\qs_2|=n_{p_i}(\phi_i)$.

\subsection{Matched chain complexes for connected sums}
\label{sec:totally-degenerated-ac-complexes}

Suppose that $\cH_1$ and $\cH_2$ are two singly pointed Heegaard diagrams. The matched moduli spaces $\cM\cM_{J_1\wedge J_2}(\phi_1,\phi_2)$ from Section~\ref{sec:matched-disks} may be used to construct a chain complex $\CF^-_{J_1\wedge J_2}(\cH_1, \cH_2)$, as we now describe.

 As a module, $\CF^-_{J_1\wedge J_2}(\cH_1, \cH_2)$ is generated over $\bF[U]$ by pairs $\xs_1\times \xs_2$, where $\xs_i$ is a generator on $\cH_i$. We define  
\[
\d_{J_1\wedge J_2}(\xs_1\times \xs_2):=\sum_{\substack{
\phi_1\in \pi_2(\xs_1,\ys_1)\\
\phi_2\in \pi_2(\xs_2,\ys_2)\\
n_{p_1}(\phi_1)=n_{p_2}(\phi_2)\\
\ind_{\emb}(\phi_1,\phi_2,M_{\#})=1}} \# (\cM\cM_{J_1\wedge J_2}(\phi_1, \phi_2)/\R) \cdot  U^{n_w(\phi_1, \phi_2)}\cdot \ys_1\times \ys_2,
\]
where $n_{w}(\phi_1,\phi_2)$ denotes the multiplicity in the connected sum region.

\begin{lem}\label{lem:matched-d^2=0}
The endomorphism $\d_{J_1\wedge J_2}$ of $\CF^-_{J_1\wedge J_2}(\cH_1, \cH_2)$ squares to zero.
\end{lem}
\begin{proof} The argument follows from the usual strategy of proving that $\d^2=0$ in Floer homology. One considers the ends of the 1-dimensional moduli spaces $\cM\cM_{J_1\wedge J_2}(\phi_1,\phi_2)/\R$, ranging 
over pairs $(\phi_1,\phi_2)$ where $\ind_{\emb}(\phi_1,\phi_2,M_{\#})=2$.

Using Proposition~\ref{prop:dimension-counts}, it is straightforward to adapt \cite{LipshitzCylindrical}*{Lemma~8.2} to show that the ends of the spaces of matched, index 2 flowlines consist exactly of the following broken curves:
\begin{enumerate}
\item 2-level, broken holomorphic disks, where each level satisfies \ref{M-1}--\ref{M-7} and has matched index 1.
\item Maslov index 4, matched alpha or beta boundary degenerations, together with a constant flowline.
\end{enumerate}
The count of the first type of broken curve gives $\d_{J_1\wedge J_2}^2(\xs_1\times \xs_2)$. Hence,
\[
\d_{J_1\wedge J_2}^2(\xs_1\times \xs_2)=U\cdot \left(\# \tilde{\cN}^{\alpha}(\xs_1\times \xs_2)+\# \tilde{\cN}^{\beta}(\xs_1\times \xs_2)\right) \cdot \xs_1\times \xs_2,
\]
where $\tilde{\cN}^{\a}(\xs_1\times \xs_2)$ and $\tilde{\cN}^{\b}(\xs_1\times \xs_2)$ denote the moduli spaces of Maslov index 4, matched alpha and beta degenerations at $\xs_1\times \xs_2$, modulo the action of $\Aut(\bH)$. By Lemma~\ref{lem:t-matched-boundary-degenerations},
\[
\# \tilde{\cN}^{\a}(\xs_1\times \xs_2)\equiv \# \tilde{\cN}^{\b}(\xs_1\times \xs_2)\equiv 1.
\]
Hence $\d^2_{J_1\wedge J_2}=0$, concluding the proof.
\end{proof}

If $J$ is a (non-singular) almost complex structure on $\Sigma_1\# \Sigma_2\times [0,1]\times \R$, satisfying \eqref{J1}--\eqref{J5}, we may define a transition map
\[
\Psi_{J\to J_1\wedge J_2}\colon \CF^-_{J}(\cH_1\# \cH_2)\to \CF^-_{J_1\wedge J_2}(\cH_1, \cH_2),
\]
by picking a non-cylindrical almost complex structure  $\tilde{J}$ on $\Sigma_1\# \Sigma_2\times [0,1]\times \R$, which for large $T\in \R$, coincides with an almost complex structure $J(T)$ along $\Sigma_1\# \Sigma_2\times [0,1]\times \{T\}$, which has a neck of length $T$.  We define $\Psi_{J\to J_1\wedge J_2}$ to count index $0$ $\tilde{J}$-holomorphic curves in $\Sigma\times [0,1]\times \R$.  A  map 
\[
\Psi_{J_1\wedge J_2\to J}\colon \CF^-_{J_1\wedge J_2}(\cH_1, \cH_2)\to \CF^-_J(\cH_1\# \cH_2),
\]
 is defined similarly.

\begin{lem}\label{lem:desingularization-maps} The maps $\Psi_{J\to J_1\wedge J_2}$ and  $\Psi_{J_1\wedge J_2 \to J}$
are chain maps, and furthermore, are homotopy inverses.
\end{lem}
\begin{proof}  Suppose $\tilde{J}$ is an almost complex structure, as above. To see that $\Psi_{J\to J_1\wedge J_2}$ is a chain map, we consider the ends of index 1 families of $\tilde{J}$-holomorphic matched curves. Using Proposition~\ref{prop:dimension-counts}, it is straightforward to see that the ends of $\tilde{J}$-holomorphic matched moduli spaces consist of the following configurations:
\begin{enumerate}
\item An index $1$ $J_1\wedge J_2$-holomorphic matched curve $(u_1,u_2)$, and an index 0 $\tilde{J}$-holomorphic curve.
\item An index $0$ $\tilde{J}$-holomorphic curve, and an index 1 $J$-holomorphic curve.
\end{enumerate}
This implies that $\Psi_{J\to J_1\wedge J_2}$ is a chain map.

Let $\tilde{I}$ be an almost complex structure which can similarly be used to compute $\Psi_{J_1\wedge J_2\to J}$. A homotopy between $\Psi_{J_1\wedge J_2\to J}\circ \Psi_{J\to J_1\wedge J_2}$ and $\id_{\CF^-_{J}(\cH_1\# \cH_2)}$ is constructed by picking a 1-parameter family of almost complex structures $(J_t)_{t\in (0,1)}$ on $\Sigma\times [0,1]\times \R$, such that as $t\to 0$, $J_t$  approaches $J$, and as $t\to 1$, the family $J_t$ approaches the two-level almost complex structure with $\tilde{J}$ on one level, and $\tilde{I}$ on the other. Counting the ends of the moduli spaces of index 0, $J_t$-holomorphic curves, ranging over $t\in (0,1)$, we obtain
\[
\Psi_{J_1\wedge J_2\to J} \circ \Psi_{J\to J_1\wedge J_2}+\id=\d \circ H+H\circ \d,
\]
where $H$ counts index $-1$ curves which are $J_t$-holomorphic for some $t\in (0,1)$. A homotopy between $\id_{\CF^-_{J_1\wedge J_2}(\cH_1, \cH_2)}$ and $\Psi_{J\to J_1\wedge J_2}\circ \Psi_{J_1\wedge J_2\to J}$ is constructed in an analogous fashion.
\end{proof}

 \begin{rem}
 When $\mu(\phi_1\# \phi_2)=1$, one can identify $\cM\cM_{J_1\wedge J_2}(\phi_1,\phi_2)$ with the ordinary moduli space $\cM_{J(T)}(\phi_1\# \phi_2)$ for large $T$, where $J(T)$ is obtained by joining $J_1$ and $J_2$ with neck-length $T$. Hence we could avoid Lemma~\ref{lem:matched-d^2=0} by just using the normal proof that $\CF^-_{J(T)}(\cH_1\# \cH_2)$ is a chain complex, where $J(T)$ is an almost complex structure with sufficiently long neck. This becomes impractical when working with higher polygon counts, since the polygon counting maps count infinitely many classes, and it is not clear we can pick a single value of $T$ which works for all classes. 
\end{rem}

\subsection{Maximally pinched families of almost complex structures}

\label{sec:omega-a-compatible}
\label{sec:maximally-pinched}

In Section~\ref{sec:matched-disks}, we described moduli spaces of disks for almost complex structures which were stretched along a simple closed curve in the Heegaard surface. The construction therein adapts easily to connected sums of Heegaard $n$-tuples, however the most straightforward generalization is slightly too restrictive for our purposes. In this section, we describe a more subtle generalization of the connected sum operation which we refer to as a stratified family of almost complex structures on $\Sigma\times D_n$ which is \emph{maximally pinched} along a simple closed curve $S\subset \Sigma$. Roughly speaking, this means that if $f\colon X\to 2^{\cA}$ is the stratifying map, then $J_x$ has infinite neck length along $S$ on any component of $D_n\setminus f(x)$ whose associated subdiagram of $\cD$ has no attaching curves which intersect $S$. For Heegaard diagrams (with two attaching curves) this reduces to the infinite neck length construction found in Section~\ref{sec:matched-disks}, but for Heegaard $n$-tuples with $n\ge 3$, it is a more flexible generalization.

 We begin with some notation. First, let $j_0$ denote the almost complex structure on $S^1\times [-1,1]$ given by $j_0 (\d/\d \theta)=\d/\d t$. If $\delta\in [0,1)$, we write $j_\delta$ for an almost complex structure which is conformally equivalent to $(S^1\times [-1/(1-\delta),1/(1-\delta)],j_0)$. Furthermore, we assume that $j_{\delta}$ is chosen to coincide with $j_0$ in a neighborhood of $S^1\times \{-1,1\}$. We do not need a concrete formula for $j_\delta$, however we assume that sending $\delta\to 1$ corresponds to the almost complex structure collapsing along the circle $S^1\times \{0\}$.

 Next, suppose that $S\subset \Sigma$ is a simple closed curve, and $(J_x)_{x\in X}$ is a stratified family of almost complex structures on $\Sigma\times D_n$, which is split on $N(S)\times U$, for some open set $U\subset D_n\times X$, where $N(S)$ denotes a regular neighborhood of $S$. In particular, over $N(S)\times U$, there is an induced almost complex structure $j_\Sigma$ on $N(S)$, and we assume that $(N(S),j_\Sigma)$ is conformally identified with $(S^1\times [-1,1],j_0)$. If $\delta\colon D_n\times X\to [0,1)$ is a smooth function with support in $U$, we define a family $(J_x(\delta))_{x\in X}$ by modifying $J_x$ so that it coincides with $j_{\delta(x,y)}\times j_{D_n,x}$ on $N(S)\times U$, and coincides with $J_x$, outside of $N(S)\times U$.

\begin{define}\label{def:maximally-pinched}
Suppose that $\cD=(\Sigma,\gs_1,\dots, \gs_n)$ is a Heegaard $n$-tuple, and $S$ is a closed curve on $\Sigma$ (possibly intersecting some of the $\gs_i$). Suppose that $\cA$ is a set of decomposing arcs of $D_n$, $X$ is a compact, smooth manifold with corners, and $f\colon X\to 2^{\cA}$ is a stratified decomposing map, as in Definition~\ref{def:stratified-decomposing-map}. We say a family of (singular) almost complex structures $(J_x)_{x\in X}$ on $\Sigma\times D_n$ is a \emph{stratified family which is maximally pinched along} $S$ if it is constructed via the following procedure:
\begin{enumerate}[ref=Step-\arabic*, label=(Step-\arabic*),leftmargin=*,widest=IIII]
\item Pick a stratified family of almost complex structures $(J_x^0)_{x\in X}$ on $\Sigma\times D_n$, such that there is an open set $U\subset D_n\times X$,  satisfying the following:
\begin{enumerate}
\item $J_x$ coincides with $j_0\times j_{D_n,x}$ on $N(S)\times U$.
\item Suppose $x\in X$ and $f(x)=\{a_{i_1},\dots, a_{i_l}\}$. If $D_l$ is a component of $D_n \setminus f(x)$ corresponding to a subdiagram of $\cD$, all of whose attaching curves are disjoint from $S$, and whose boundary contains the arcs arcs $a_{j_1},\dots, a_{j_s}$, then 
\[
 (D_l\cup N(a_{j_1})\cup \cdots \cup N(a_{j_s}) )\times \{x\}\subset U.
\]
\item If $\gs_i$ and $\gs_{i+1}$ are cyclically consecutive attaching curves on $\cD$ which are disjoint from $S$, then $U$ contains $ E_{\g_i,\g_{i+1}}\times X$, where $E_{\g_i,\g_{i+1}}$ denotes the cylindrical end of $D_n$ corresponding to $\gs_i$ and $\gs_{i+1}$.
\end{enumerate}
\item Pick a smooth $\delta\colon D_n\times X\to [0,1]$ (the neck-length parameter) which is supported in $U$. We assume that $\delta(y,x)=1$ if and only if $f(x)=\{a_{i_1},\dots, a_{i_m}\}$ and $y$ is in the closure of a component of $D_n\setminus f(x)$ which corresponds to a subdiagram of $\cD$, all of whose attaching curves are disjoint from $S$. We further assume the following:
\begin{enumerate}
\item $\delta\to 1$ in the cylindrical ends $E_{\g_i,\g_{i+1}}$ of $D_n$, whenever $\gs_i$ and $\gs_{i+1}$ are disjoint from $S$.
\item $\delta$ is constant along the $[0,1]$ component of $E_{\g_i,\g_{i+1}}\iso [0,1]\times [0,\infty)$.
\item For each $x$, the function $\delta(-,x):D_n\to [0,1]$ is constant on $a_i\times \{t\}\subset N(a_i)$.
\end{enumerate}
\item Set $(J_x)_{x\in X}$ to be $(J_x^0(\delta))_{x\in X}$.
\end{enumerate}
\end{define}

If $S_1,\dots, S_k$ is a collection of pairwise disjoint, simple closed  curves on $\Sigma$, we may similarly define a family $(J_x)_{x\in X}$ to be maximally pinched along $(S_1,\dots, S_k)$ if it constructed by performing the above procedure along each $S_i$.

\section{Stabilizing  Heegaard multi-diagrams}
\label{sec:stabilization}
In this section, we describe several  stabilization operations for Heegaard diagrams, triples and quadruples. In Section~\ref{sec:stabilizations-holomorphic}, we describe their interaction with the counts of holomorphic curves.

\subsection{Stabilizing Heegaard diagrams}

\begin{figure}[ht!]
	\centering
\begingroup%
  \makeatletter%
  \providecommand\color[2][]{%
    \errmessage{(Inkscape) Color is used for the text in Inkscape, but the package 'color.sty' is not loaded}%
    \renewcommand\color[2][]{}%
  }%
  \providecommand\transparent[1]{%
    \errmessage{(Inkscape) Transparency is used (non-zero) for the text in Inkscape, but the package 'transparent.sty' is not loaded}%
    \renewcommand\transparent[1]{}%
  }%
  \providecommand\rotatebox[2]{#2}%
  \newcommand*\fsize{\dimexpr\f@size pt\relax}%
  \newcommand*\lineheight[1]{\fontsize{\fsize}{#1\fsize}\selectfont}%
  \ifx\svgwidth\undefined%
    \setlength{\unitlength}{182.14771615bp}%
    \ifx\svgscale\undefined%
      \relax%
    \else%
      \setlength{\unitlength}{\unitlength * \real{\svgscale}}%
    \fi%
  \else%
    \setlength{\unitlength}{\svgwidth}%
  \fi%
  \global\let\svgwidth\undefined%
  \global\let\svgscale\undefined%
  \makeatother%
  \begin{picture}(1,0.72092426)%
    \lineheight{1}%
    \setlength\tabcolsep{0pt}%
    \put(0,0){\includegraphics[width=\unitlength,page=1]{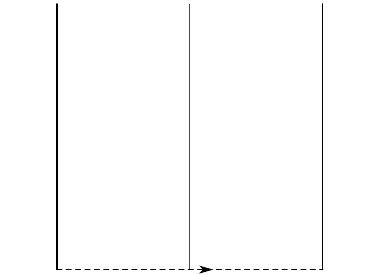}}%
    \put(0.47238638,0.32550139){\color[rgb]{1,0,0}\makebox(0,0)[rt]{\lineheight{1.25}\smash{\begin{tabular}[t]{r}$\xi$\end{tabular}}}}%
    \put(0.59395486,0.32550139){\color[rgb]{0,0,1}\makebox(0,0)[lt]{\lineheight{1.25}\smash{\begin{tabular}[t]{l}$\zeta$\end{tabular}}}}%
    \put(0.13590017,0.25125694){\color[rgb]{0,0,0}\makebox(0,0)[rt]{\lineheight{1.25}\smash{\begin{tabular}[t]{r}$p_0$\end{tabular}}}}%
    \put(0.86386393,0.25125694){\color[rgb]{0,0,0}\makebox(0,0)[lt]{\lineheight{1.25}\smash{\begin{tabular}[t]{l}$p_1$\end{tabular}}}}%
    \put(0,0){\includegraphics[width=\unitlength,page=2]{fig127.pdf}}%
  \end{picture}%
\endgroup%

	\caption{The standard diagram $(S^2,\xi,\zeta,p_1,p_2)$ for $(S^3,p_1,p_2)$ used for a 1-handle stabilization. The Heegaard surface is obtained by identifying the top and bottom sides, and collapsing the left and right sides each to a point.}
	\label{fig:127}
\end{figure}

The most basic way to stabilize a Heegaard diagram is to take the connected sum with a standard genus 1 diagram for $S^3$. More generally, we may take the connected sum of a Heegaard diagram with a standard diagram for the lens space $L(m,1)$, shown in Figure~\ref{fig:179}.  A third operation is related to 1-handles, and involves taking the connected sum at two points with a standard, genus 0, doubly pointed diagram of $(S^3,p_1,p_2)$. See Figure~\ref{fig:127}. We call these \emph{elementary stabilizations}, \emph{elementary lens space stabilizations}, and \emph{elementary 1-handle stabilizations}, respectively.

It is convenient to perform multiple stabilizations at once. We make the following definition:

\begin{define}\label{def:multi-stabilization-diagram}
 We say a Heegaard diagram $\cH_0=(\Sigma_0,\xis,\zetas,p_0)$ is a \emph{multi-stabilizing diagram} if it is constructed by starting with the diagram $(S^2,p_0)$ (with no attaching curves) and inductively attaching elementary stabilizations and 1-handle stabilizations.
\end{define}

One could also consider multi-stabilizing diagrams which have standard lens space summands, though we will not have a need for such diagrams.

\subsection{Stabilizing Heegaard triples}

We now consider stabilizations of Heegaard triples.

\begin{figure}[ht!]
	\centering
	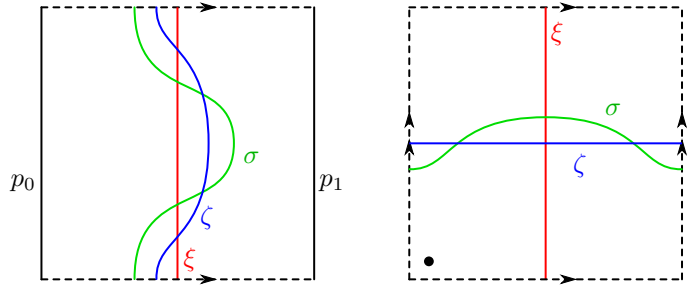
	\caption{Left: An elementary 1-handle triple. Right: An elementary stabilization triple. We collapse the solid lines to basepoints, and identify the dashed sides.}
	\label{fig:163}
\end{figure}

\begin{define}\label{def:elem-stabilization-triple} 
\begin{enumerate}
\item  We say a genus 0 triple $(S^2,\xi,\zeta,\sigma,p_0,p_1)$ is an \emph{elementary 1-handle triple}, if $\xi$, $\zeta$ and $\sigma$ are transverse, pairwise isotopic, and each pair of distinct curves in $\{\xi,\zeta,\sigma\}$ intersects in exactly 2 points.
\item We say that a genus 1 triple $(\bT^2,\xi,\zeta,\sigma,p_0)$ is an \emph{elementary stabilization triple} if (up to cyclic permutation) 
$\zeta$ and $\sigma$ are pairwise Hamiltonian isotopic, and
$|\zeta\cap \sigma|=2,$ while $|\xi\cap \zeta|=|\xi\cap \sigma|=1.$ See Figure~\ref{fig:163}.
\end{enumerate}
\end{define}

\begin{define}\label{def:stabilization-triple}
 We say a Heegaard triple $\cT_0=(\Sigma_0,\xis,\zetas,\sigmas,p_0)$ is a \emph{multi-stabilizing triple} if it is constructed by starting with the diagram $(S^2,p_0)$, and inductively attaching elementary stabilization triples and elementary 1-handle triples.
\end{define}

\subsection{Stabilizing Heegaard quadruples}

Finally, we consider stabilizations of quadruples.

\begin{define}\label{def:elem-stabilization-quad} 
\begin{enumerate}
\item  We say a genus 0 quadruple $(S^2,\xi,\zeta,\sigma,\tau,p_0,p_1)$ is an \emph{elementary 1-handle quadruple} if $\xi$, $\zeta$, $\sigma$ and $\tau$ are transverse, pairwise isotopic, and each pair of distinct curves in $\{\xi,\zeta,\sigma,\tau\}$ intersects in exactly 2 points. 
\item We say a genus 1 quadruple $(\bT^2,\xi,\zeta,\sigma,\tau,p_0)$ is an \emph{elementary stabilization quadruple} if (up to cyclic permutation) one of the following holds:
\begin{enumerate}
\item $\xi$, $\zeta$ and $\sigma$ are pairwise Hamiltonian isotopic, and
\[
|\xi\cap \zeta|=|\xi\cap \sigma|=|\zeta\cap \sigma|=2, \quad \text{while} \quad |\xi\cap \tau|=|\zeta\cap \tau|=|\sigma\cap \tau|=1.
\]
\item $\xi$ and $\zeta$ are Hamiltonian isotopic, and $\sigma$ and $\tau$ are Hamiltonian isotopic. Furthermore,
\[
|\xi\cap \zeta|=|\sigma\cap \tau|=2\quad \text{and} \quad  |\xi\cap \sigma|=|\zeta\cap \sigma|=|\xi\cap \tau|=|\zeta\cap \tau|=1.
\]
\end{enumerate}
\end{enumerate}
\end{define}

\begin{figure}[ht!]
	\centering
	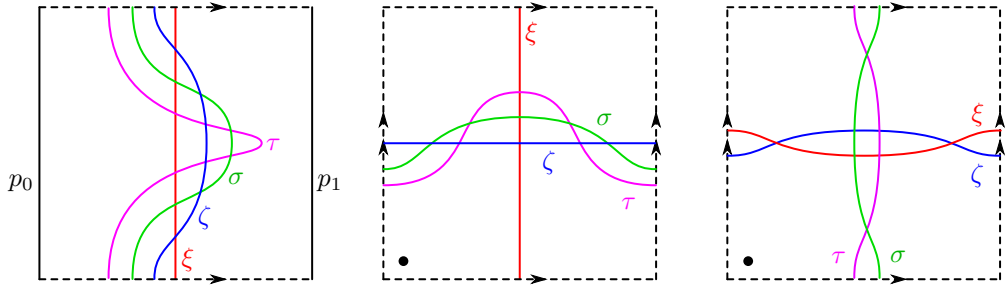
	\caption{Left: An elementary 1-handle quadruple. Center and right: Elementary stabilization quadruples.}
	\label{fig:162}
\end{figure}

See Figure~\ref{fig:162} for examples of elementary 1-handle and stabilization quadruples.

\begin{define}
We say that $\cQ_0=(\Sigma_0, \xis, \zetas,\sigmas,\taus,p_0)$ is a \emph{multi-stabilizing Heegaard quadruple} if it is constructed by starting with the diagram $(S^2,p_0)$ and inductively attaching elementary 1-handle quadruples and simple 1-handle stabilization quadruples.
\end{define}

\subsection{Index tuples}
\label{sec:index-tuple}
In this section, we define the \emph{index tuple} of a multi-stabilizing triple or quadruple. The index tuple will allow us to efficiently state the holomorphic stabilization formulas in Section~\ref{sec:stabilizations-holomorphic}.

For a multi-stabilizing triple $\cT_0$, the index tuple is a triple
\[
\ve{m}(\cT_0)=(m_1(\cT_0),m_2(\cT_0), m_3(\cT_0)).
\]
For a multi-stabilizing quadruple $\cQ_0$, the index tuple is a 4-tuple
\[
\ve{m}(\cQ_0)=(m_1(\cQ_0),m_2(\cQ_0),m_3(\cQ_0),m_4(\cQ_0)).
\]

If $\cT_0=(\Sigma_0,\xis,\zetas,\sigmas,p_0)$ is a multi-stabilizing triple, we set $m_1(\cT_0)$ to be $\mu(\psi)$ for any triangle class $\psi\in \pi_2(\Theta_{\xi,\zeta}^-, \Theta_{\zeta,\sigma}^+, \Theta_{\sigma,\xi}^+)$ with $n_{p_0}(\psi)=0$. (Note that we are using an unusual ordering of the subscripts, and $\Theta_{\sigma,\xi}^{\pm}=\Theta_{\xi,\sigma}^{\mp}$).
Similarly, we define $m_2(\cT_0)$ to be the Maslov index of any triangle class $\psi\in \pi_2(\Theta_{\xi,\zeta}^+, \Theta_{\zeta,\sigma}^-,\Theta_{\sigma,\xi}^+)$ with $n_{p_0}(\psi)=0$, and $m_3(\cT_0)$ is the Maslov index of any triangle class $\psi\in \pi_2(\Theta_{\xi,\zeta}^+,\Theta_{\zeta,\sigma}^+,\Theta_{\sigma,\xi}^-)$.

We define the index 4-tuple of a multi-stabilizing quadruple $\cQ_0=(\Sigma_0,\xis,\zetas,\sigmas,\taus,p_0)$ as follows. We define $m_1(\cQ_0)$ to be the Maslov index of any rectangle class $\psi\in \pi_2(\Theta_{\xi,\zeta}^-,\Theta_{\zeta,\sigma}^+,\Theta_{\sigma,\tau}^+, \Theta_{\tau,\xi}^+)$ with $n_{p_0}(\psi)=0$.  Similarly $m_2(\cQ_0)$ is the Maslov index of any rectangle class $\psi\in \pi_2(\Theta_{\xi,\zeta}^+, \Theta_{\zeta,\sigma}^-,\Theta_{\sigma,\tau}^+, \Theta_{\tau,\xi}^+)$ with $n_{p_0}(\psi)=0$. The integers $m_3(\cQ_0)$ and $m_4(\cQ_0)$ are defined similarly.

In the subsequent Lemma~\ref{lem:Maslov-triangle-multi-stabilization}, we prove that the definition of index tuples given above is well defined. Before stating the lemma, we note that we use the grading convention that the top degree generator of $\CF^-((S^1\times S^2)^{\# k},p_0)$ is supported in grading $k/2$.

\begin{lem}\label{lem:Maslov-triangle-multi-stabilization}
 \begin{enumerate}
\item  If $\cT_0=(\Sigma_0,\xis,\zetas,\sigmas,p_0)$ is a multi-stabilizing Heegaard triple, and $\psi\in \pi_2(\xs,\ys,\zs)$ is a homology class of triangles, then
\[
\mu(\psi)=2n_{p_0}(\psi)+\gr(\xs)+\gr(\ys)-\gr(\zs)-n/2,
\]
 where $\gr(\zs)$ is the absolute grading of $\zs$, as an element of $\bCF^-(\Sigma_0,\xis,\sigmas)$, and $n=|\xis|=|\zetas|=|\sigmas|$.
 \item If $\cQ_0=(\Sigma_0,\xis,\zetas,\sigmas,\taus,p_0)$ is a multi-stabilizing Heegaard quadruple, and $\psi\in \pi_2(\xs,\ys,\zs,\ws)$ is a homology class of quadrilaterals, then
\[
\mu(\psi)=2n_{p_0}(\psi)+\gr(\xs)+\gr(\ys)+\gr(\zs)-\gr(\ws)-n,
\]
where $\gr(\ws)$ denotes the Maslov grading of $\ws$ as an element of $\bCF^-(\Sigma_0,\xis,\taus)$, and $n=|\xis|=|\zetas|=|\sigmas|=|\taus|$.
\end{enumerate}
\end{lem}
\begin{proof} Both claims are proven using the same strategy. First, one checks directly that they hold for elementary 1-handle and stabilization triples and quadruples. Next, one argues by induction that the formulas hold when we add an extra elementary 1-handle or stabilization triple or quadruple. We leave the details to the reader.
\end{proof}

\begin{lem}\label{lem:monotonicity-index-tuple}
 Suppose $\cT_0=(\Sigma_0,\xis,\zetas,\sigmas,p_0)$ is a multi-stabilizing triple which is built using elementary triples $\cT_1,\dots, \cT_n$.
 \begin{enumerate}
 \item  $m_i(\cT_0)$ is non-positive, for $i\in \{1,2,3\}$.
 \item $\ve{m}(\cT_0)=\sum_{i=1}^n \ve{m}(\cT_i)$.
 \item $m_i(\cT_0)=0$ if and only if $m_i(\cT_j)=0$ for all $j\in \{1,\dots, n\}$.
\end{enumerate} 
  The analogous statements hold for multi-stabilizing quadruples.
\end{lem}
\begin{proof}
The same argument works for both triples and quadruples, so we focus on  triples. Firstly, additivity of $\ve{m}(\cT_0)$ is immediate from the construction, as the grading is additive under tensor products. For elementary tuples, the relation $m_i(\cT_0)\le 0$ is easily checked. It follows that $m_i(\cT_0)\le 0$ and $m_i(\cT_0)=0$ if and only if each $m_i(\cT_j)=0$ for all $j$.
\end{proof}

\begin{rem}\label{rem:cyclically-permute-triple}
Cyclic permutation of the attaching curves in a multi-stabilizing triple $\cT_0$, has the effect on $\ve{m}(\cT_0)$ of cyclic permuation. The same holds for multi-stabilizing quadruples.
\end{rem}

\begin{figure}[ht!]
	\centering
	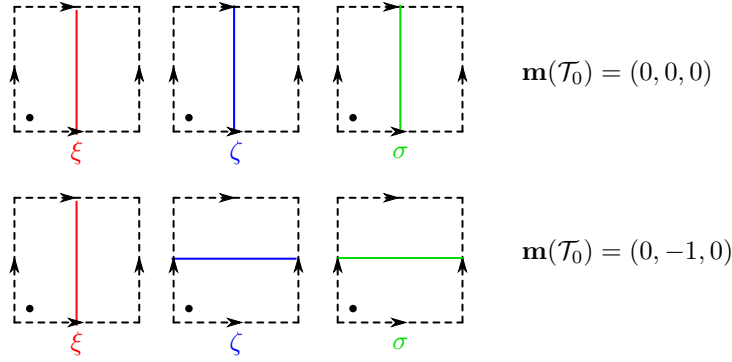
	\caption{Top: An elementary stabilization triple. Bottom: An elementary 1-handle triple. The index tuples are labeled.}\label{fig:117}
\end{figure}

\begin{figure}[ht!]
	\centering
	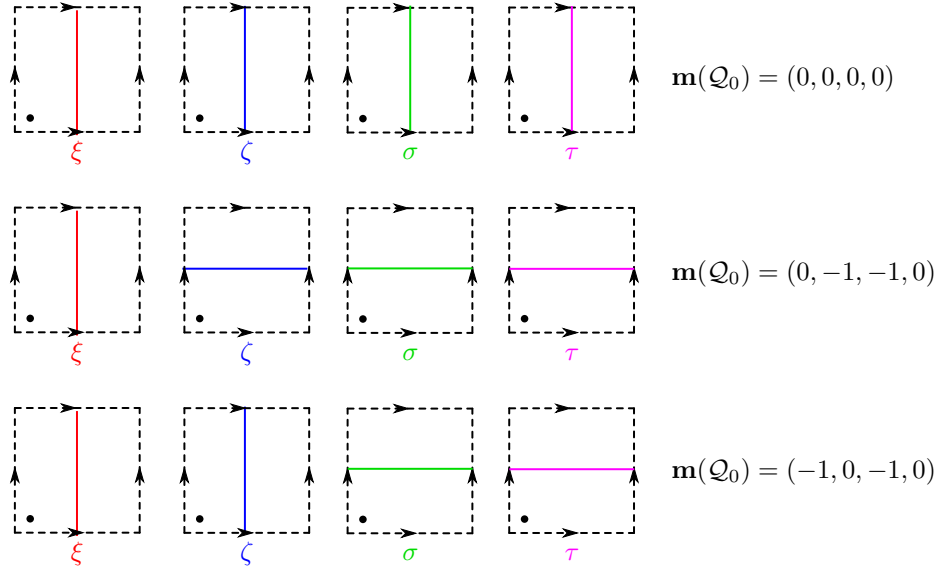
	\caption{Elementary 1-handle and stabilization quadruples, and their index tuples.}\label{fig:118}
\end{figure}

The index tuples of elementary stabilization and 1-handle triples and quadruples are shown in Figures~\ref{fig:117} and~\ref{fig:118}.

We now state an associativity relation for the index tuples:

\begin{lem}\label{lem:index-associativity} Let $\cQ_0=(\Sigma_0,\xis,\zetas,\sigmas,\taus,p_0)$ be a multi-stabilizing quadruple, and let $\cT_{\xi,\zeta,\sigma}$ denote the triple $(\Sigma_0,\xis,\zetas,\sigmas,p_0)$ (and so forth for the other sub-triples). The following hold:
\begin{enumerate}
\item $m_1(\cQ_0)=m_1(\cT_{\xi,\zeta,\sigma})+m_1(\cT_{\xi,\sigma,\tau})=m_3(\cT_{\zeta,\sigma,\tau})+m_1(\cT_{\xi,\zeta,\tau}).$
\item $m_2(\cQ_0)=m_2(\cT_{\xi,\zeta,\sigma})+m_1(\cT_{\xi,\sigma,\tau})=m_1(\cT_{\zeta,\sigma,\tau})+m_2(\cT_{\xi,\zeta,\tau})$.
\item $m_3(\cQ_0)=m_3(\cT_{\xi,\zeta,\sigma})+m_2(\cT_{\xi,\sigma,\tau})=m_2(\cT_{\zeta,\sigma,\tau})+m_2(\cT_{\xi,\zeta,\tau})$.
\item $m_4(\cQ_0)=m_3(\cT_{\xi,\zeta,\sigma})+m_3(\cT_{\xi,\sigma,\tau})=m_3(\cT_{\zeta,\sigma,\tau})+m_3(\cT_{\xi,\zeta,\tau}).$
\end{enumerate}
\end{lem}
\begin{proof} We focus on the first equation of the first claim. We pick triangle classes 
\[
\psi_{\xi,\zeta,\sigma}\in \pi_2(\Theta_{\xi,\zeta}^-,\Theta_{\zeta,\sigma}^+,\Theta_{\sigma,\xi}^+)\quad \text{and}\quad \psi_{\xi,\sigma,\tau}\in \pi_2(\Theta_{\xi,\sigma}^-, \Theta_{\sigma,\tau}^+, \Theta_{\tau,\xi}^+),
\]
which both have zero multiplicity on $p_0$. The classes $\psi_{\xi,\zeta,\sigma}$ and $\psi_{\xi,\sigma,\tau}$ may be spliced together since $\Theta_{\sigma,\xi}^+=\Theta_{\xi,\sigma}^-$. Using additivity of the index under juxtoposition, as well as the definition of $m_i$, we obtain the stated equation. The same argument applies for all of the other equations.
\end{proof}

\section{Stabilizations and holomorphic curves}
\label{sec:stabilizations-holomorphic}
In this section, we describe how the counts of holomorphic curves interact with the stabilization operations from Section~\ref{sec:stabilization}.

\subsection{Stabilizations and holomorphic disks}

We now consider multi-stabilizations of Heegaard diagrams. Suppose that $\cH_0=(\Sigma_0,\xis,\zetas,p_0)$ is a multi-stabilizing Heegaard diagram. By construction, there is a top graded intersection point $\Theta_{\xi,\zeta}^+\in \bT_{\xi}\cap \bT_{\zeta}$, as well as a bottom graded intersection point $\Theta_{\xi,\zeta}^-$. We define
\[
F_1^{\xi,\zeta}\colon \bCF^-_J(\cH)\to \bCF^-_{J\wedge J_0}(\cH,\cH_0)\quad \text{and} \quad F_3^{\xi,\zeta}\colon \bCF^-_{J\wedge J_0}(\cH, \cH_0)\to \bCF^-_J(\cH)
\]
via the formulas
\[
F_1^{\xi,\zeta}(\xs)=\xs\times \Theta_{\xi,\zeta}^+\quad \text{and} \quad F_3^{\xi,\zeta}(\xs\times \Theta)=\begin{cases} \xs& \text{if } \Theta=\Theta_{\xi,\zeta}^-\\
0& \text{otherwise}
\end{cases},
\]
extended equivariantly over $U$.  Similar to Ozsv\'{a}th and Szab\'{o}'s original 1- and 3-handle maps \cite{OSTriangles}*{Section~4.3}, we have the following:

\begin{lem}\label{lem:singular-1-handle-chain-maps}
The maps $F_1^{\xi,\zeta}$ and $F_3^{\xi,\zeta}$ are chain maps.
\end{lem}
\begin{proof} 
We focus on the claim for $F_1^{\xi,\zeta}$, as $F_3^{\xi,\zeta}$ is algebraically dual. Suppose  $\phi_0\in \pi_2(\Theta_1,\Theta_2)$ is a class on $\cH_0$. By the grading formula, we have
\[
\mu(\phi_0)=\gr(\Theta_1,\Theta_2)+2n_{p_0}(\phi_0).
\]
Hence, if $\phi\in \pi_2(\xs,\ys)$ is a class on $\cH$, then
\begin{equation}
\ind_{\emb}(\phi,\phi_0,M_{\#})=\mu(\phi)+\gr(\Theta_1,\Theta_2).
\end{equation}
(We remind the reader that $\ind_{\emb}(\phi,\phi_0,M_{\#})=\mu(\phi\# \phi_0)$). The differential on $\bCF^-_{J\wedge J_0}(\cH,\cH_0)$ counts representatives of pairs  of classes $(\phi,\phi_0)$ with $\ind_{\emb}(\phi,\phi_0,M_{\#})=1$. We note that $\gr(\Theta_{\xi,\zeta}^+,\Theta_2)\ge 0$. Furthermore,  transversality for holomorphic curves representing $\phi$ implies that $\mu(\phi)\ge 0$, with equality if and only if $\phi$ is a constant class. Hence, any index 1 pair $(\phi,\phi_0)\in \pi_2(\xs,\ys)\times \pi_2(\Theta_{\xi,\zeta}^+,\Theta_2)$ which admits holomorphic representatives must have one of the following configurations: either $\mu(\phi)=1$ and $\Theta_{\xi,\zeta}^+=\Theta_2$, or $\mu(\phi)=0$ and $\gr(\Theta_{\xi,\zeta}^+,\Theta_2)=1$. The contribution to $\d(\xs\times \Theta^+_{\xi,\zeta})$ of curves satisfying the latter condition coincides with $\xs\otimes \hat{\d}_0(\Theta^+_{\xi,\zeta})$, where $\hat{\d}_0$ denotes the ordinary differential on the complex $\widehat{\CF}(\cH_0)$. Since $\Theta^+_{\xi,\zeta}$ is a cycle in $\widehat{\CF}(\cH_0)$, this contribution is zero. It remains to consider the case where $\mu(\phi)=1$ and $\Theta_1=\Theta_2=\Theta_{\xi,\zeta}^+$. In this case, the main claim amounts to showing that if $\ve{d}\in ((0,1)\times \R)^{k}$ is a generic $k$-tuple of points, then
\begin{equation}
\sum_{\substack{\phi_0\in \pi_2(\Theta_{\xi,\zeta}^+, \Theta_{\xi,\zeta}^+)\\ 
n_{p_0}(\phi_0)=k}} \# \cM(\phi_0,\ve{d}) \equiv 1 \pmod{2}, \label{eq:total-count-matched-disks}
\end{equation}
which is verified in \cite{ZemDuality}*{Proposition~6.1}, by degenerating the set $\ve{d}$, adapting an argument of Ozsv\'{a}th and Szab\'{o} \cite{OSLinks}*{Lemma~6.4}.
\end{proof}

\subsection{Stabilizations and holomorphic triangles}
\label{sec:stabilizations-hol-triangles}

We now describe the interaction of the holomorphic triangle maps with multi-stabilizations.

\begin{prop}\label{prop:stabilize-triangles-general}
 Suppose that $\cT=(\Sigma,\as,\bs,\gs,w)$ is a Heegaard triple and $\cT_0=(\Sigma_0,\xis,\zetas,\sigmas,p_0)$ is a multi-stabilizing triple. We form the wedge product triple $\cT\wedge \cT_0$ (thought of as a connected sum with infinite neck length), at $p_0$ and some point in $\Sigma\setminus (\as\cup \bs\cup \gs)$. We delete the basepoint $p_0$, and leave $w$ as a basepoint. With respect to the almost complex structures $J$ and $J\wedge J_0$, we have the following:
 \begin{enumerate}
 \item If $m_1(\cT_0)=0$, then
 \[
 \begin{split}
 &f_{\cT\wedge \cT_0;\frs\# \frs_0}(\xs\times \Theta_{\xi,\zeta}^-,\ys\times \Theta_{\zeta,\sigma}^+)\\
 =&f_{\cT;\frs}(\xs,\ys)\otimes \Theta_{\xi,\sigma}^-+\sum_{\substack{\zs\in \bT_{\a}\cap \bT_{\g}\\ \Theta\in \bT_{\xi}\cap \bT_{\sigma}\setminus \{\Theta_{\xi,\sigma}^-\} }} C_{\xs,\ys,\zs,\Theta}\cdot \zs\times \Theta,
 \end{split}
 \]
 for some $C_{\xs,\ys,\zs,\Theta}\in \bF[U]$.
 \item If $m_2(\cT_0)=0$, then
 \[
 \begin{split}
 &f_{\cT\wedge \cT_0;\frs\#\frs_0}(\xs\times \Theta_{\xi,\zeta}^+,\ys\times \Theta_{\zeta,\sigma}^-)\\
 =&f_{\cT;\frs}(\xs,\ys)\otimes \Theta_{\xi,\sigma}^-+\sum_{\substack{\zs\in \bT_{\a}\cap \bT_{\g}\\ \Theta\in \bT_{\xi}\cap \bT_{\sigma}\setminus \{\Theta_{\xi,\sigma}^-\} }} C_{\xs,\ys,\zs,\Theta}\cdot \zs\times \Theta,
 \end{split}
 \]
 for some $C_{\xs,\ys,\zs,\Theta}\in \bF[U]$.
 \item If $m_3(\cT_0)=0$, then
  \[
  \begin{split}
& f_{\cT\wedge \cT_0;\frs\#\frs_0}(\xs\times \Theta_{\xi,\zeta}^+,\ys\times \Theta_{\zeta,\sigma}^+)\\
 =&f_{\cT;\frs}(\xs,\ys)\otimes \Theta_{\xi,\sigma}^+.
 \end{split}
 \]
 \end{enumerate}
 Here, $\frs_0\in \Spin^c(X_{\xi,\zeta,\tau})$ is the one which restricts to the torsion $\Spin^c$ structure on each end.
\end{prop}

\begin{proof} We focus on the third claim, since the other two are obtained by cyclically permuting the attaching curves. 

Suppose that $\psi\# \psi_0\in \pi_2(\xs\times \Theta^+_{\xi,\zeta}, \ys\times \Theta^+_{\zeta,\sigma},\zs\times \Theta)$, for some $\Theta\in \bT_{\xi}\cap \bT_{\sigma}$. By definition of $m_3(\cT_0)$,
\[
\mu(\psi_0)=m_3(\cT_0)+2n_{p_0}(\psi_0)+\gr(\Theta_{\xi,\sigma}^+,\Theta).
\]
Since $m_3(\cT_0)=0$ by assumption, we conclude that
\begin{equation}
\mu(\psi_0)=2n_{p_0}(\psi_0)+\gr(\Theta_{\xi,\sigma}^+,\Theta).
\end{equation}
\label{eq:Maslov-index-psi-0-triangles}
Using ~\eqref{eq:Maslov-index-psi-0-triangles}, we have
\begin{equation}
\ind_{\emb}(\psi,\psi_0,M_{\#})=\mu(\psi)+\gr(\Theta_{\xi,\sigma}^+,\Theta).\label{eq:maslov-connected-sum-triangles}
\end{equation}

If $\cM\cM_{J\wedge J_0}(\psi,\psi_0)$ is non-empty, then $\cM_J(\psi)$ is non-empty by definition, so $\mu(\psi)\ge 0$ by transversality. Equation~\eqref{eq:maslov-connected-sum-triangles} now gives
\[
\mu(\psi)=\gr(\Theta_{\xi,\sigma}^+,\Theta)=0.
\]
Hence, it remains to show that
\begin{equation}
\sum_{\substack{\psi_0\in \pi_2(\Theta^+_{\xi,\zeta},\Theta_{\zeta,\sigma}^+,\Theta_{\xi,\sigma}^+)\\ n_{p_0}(\psi_0)=k}} \# \cM(\psi,\ve{d})\equiv 1, \label{eq:main-triangle-count-triangles}
\end{equation}
for a generic $\ve{d}\in \Delta^{k},$ which is disjoint from the fat diagonal. Equation~\eqref{eq:main-triangle-count-triangles} is proven by considering a path $(\ve{d}_t)_{t\in [0,\infty)}$, satisfying the following:
\begin{enumerate}
\item $\ve{d}_0=\ve{d}$.
\item The $k$ components of $\ve{d}_t$ approach $\infty$ in the $\a$-$\b$ cylindrical end of $\Delta$. Furthermore, with respect to an identification of the cylindrical end with $[0,1]\times [0,\infty)$, the points of $\ve{d}_t$ approach some fixed 
\[
\ve{d}'\in ((0,1)\times \R)^k,
\]
modulo overall translation by $\R$. Furthermore, $\ve{d}'$ lies in the complement of the fat diagonal.
\end{enumerate}

If $\zs\in \bT_{\xi}\cap \bT_{\sigma}$, we consider the 1-dimensional moduli space
\[
\cM:= \bigcup_{t\in [0,\infty)} \coprod_{\substack{\psi\in \pi_2(\Theta_{\xi,\zeta}^+,\Theta_{\zeta,\sigma}^+,\Theta_{\xi,\sigma}^+)\\ n_{p_0}(\psi)=k }} \cM(\psi,\ve{d}_t).
\]
The count of the ends of $\cM$ at $t=0$ coincides with the left-hand side of ~\eqref{eq:main-triangle-count-triangles}. The ends of $\cM$ at finite $t$ occur in canceling pairs, since the differentials vanish on $\widehat{\CF}(\Sigma_0, \xis,\zetas)$, $\widehat{\CF}(\Sigma_0,\zetas,\sigmas)$ and $\widehat{\CF}(\Sigma_0,\xis,\zetas)$. The remaining ends of $\cM$ occur as $t\to\infty$, and are identified with the Cartesian product
\begin{equation}
\bigg(\coprod_{\substack{\psi\in \pi_2(\Theta_{\xi,\zeta}^+, \Theta_{\zeta,\sigma}^+,\Theta_{\xi,\sigma}^+) \\ n_{p_0}(\psi)=0\\ \mu(\psi)=0}} \cM(\psi) \bigg)\times \bigg( \coprod_{\substack{\phi\in \pi_2(\Theta_{\xi,\zeta}^+,\Theta_{\xi,\zeta}^+) \\ n_{p_0}(\psi_0)=k}} \cM(\phi,\ve{d}')\bigg).\label{eq:cobordism-triangles-d}
\end{equation}

 The left-hand side of ~\eqref{eq:cobordism-triangles-d}, is the $\Theta_{\xi,\sigma}^+$ coefficient of $\hat{f}_{\xi,\zeta,\sigma}(\Theta_{\xi,\zeta}^+,\Theta_{\zeta,\sigma}^+)$. The right-hand factor of ~\eqref{eq:cobordism-triangles-d} has total count 1 (modulo 2) by ~\eqref{eq:total-count-matched-disks}.
 
 The main claim thus amounts to showing that the left-hand factor of ~\eqref{eq:cobordism-triangles-d} has total count congruent to 1. This is equivalent to showing
\begin{equation}
\hat{f}_{\xi,\zeta,\sigma}(\Theta_{\xi,\zeta}^+,\Theta_{\zeta,\sigma}^+)=\Theta_{\xi,\sigma}^+,\label{eq:model-computation-triangles}
\end{equation}
where $\hat{f}_{\xi,\zeta,\sigma}$ denotes the map which counts triangles on $\cT_0$ which do not cross the basepoints. Equation~\eqref{eq:model-computation-triangles} is verified in the subsequent Lemma~\ref{lem:model-computation-stabilize-triangle}.
\end{proof}

\begin{rem} Note that if we were to verify the first two statements of Proposition \ref{prop:stabilize-triangles-general} directly, the main difference would be that $\gr(\Theta_{\xi,\sigma}^-,\Theta)$ may be negative, allowing for the existence of curves with $\mu(\psi) \geq 1$. \end{rem}

\begin{lem} \label{lem:model-computation-stabilize-triangle}
 Suppose that $\cT_0=(\Sigma_0,\xis,\zetas,\taus,\ve{p}_0)$ is a multi-stabilizing Heegaard triple.
 \begin{enumerate}
 \item If $m_1(\cT_0)=0$, then
 \[
	\hat{f}_{\cT_0}(\Theta_{\xi,\zeta}^-,\Theta_{\zeta,\sigma}^+)=\Theta_{\xi,\sigma}^-+\sum_{\Theta\in \bT_{\xi}\cap \bT_{\sigma}\setminus \{\Theta_{\xi,\sigma}^-\}} C_{\Theta}\cdot \Theta, 
 \]
 for some $C_{\Theta}\in \bF$.
 \item If $m_2(\cT_0)=0$, then
  \[
	\hat{f}_{\cT_0}(\Theta_{\xi,\zeta}^+,\Theta_{\zeta,\sigma}^-)=\Theta_{\xi,\sigma}^-+\sum_{\Theta\in \bT_{\xi}\cap \bT_{\sigma}\setminus \{\Theta_{\xi,\sigma}^-\}} C_{\Theta}\cdot \Theta, 
 \]
 for some $C_{\Theta}\in \bF$.
 \item If $m_3(\cT_0)=0$, then
 \[
\hat{f}_{\cT_0}(\Theta_{\xi,\zeta}^+,\Theta_{\zeta,\sigma}^+)=\Theta_{\xi,\sigma}^+.
 \]
 \end{enumerate}
\end{lem}
\begin{proof} We focus on the third claim, when $m_3(\cT_0)=0$, since it is the simplest notationally, and all claims are equivalent upon dualizing, as in~\eqref{eq:symmetrized-triangle-maps}.

According to Definition~\ref{def:stabilization-triple}, a multi-stabilizing triple is built by inductively attaching elementary stabilizing triples. Hence, we prove the claim by induction on the number $n$ of elementary stabilizing triples used in the construction. If $n=1$, the computation is easy to perform directly. We suppose now that the claim holds for all multi-stabilizing triples built with $n$ or fewer elementary triples. This implies that Proposition~\ref{prop:stabilize-triangles-general} holds whenever the stabilizing triple is built with $n$ or fewer elementary triples. In particular, Proposition~\ref{prop:stabilize-triangles-general} holds when $\cT_0$ is an elementary stabilizing triple.
 
We suppose now that $\cT$ is a general multi-stabilizing triple.  It is sufficient to show the present claim for some almost complex structure, since if $J$ and $J'$ are two almost complex structures on $\Sigma\times \Delta$, then we may pick a stratified family of almost complex structures  $(J_s)_{s\in [0,1]}$ on $\Sigma_0\times \Delta$ such that $J_0=J$, and as $s\to 1$, the family $J_s$ limits to $J'$ on $\Sigma\times \Delta$, as well as three non-cylindrical almost complex structures on the ends of $\Sigma\times \Delta$. Counting the ends of index 0 families of $J_s$-holomorphic triangles, we obtain 
  \begin{equation}
  \begin{split} &\hat{f}_{\cT_0,J}(\Theta_{\xi,\zeta}^+,\Theta_{\zeta,\sigma}^+)+\Psi_{J'\to J} (f_{\cT_0,J'}(\Psi_{J\to J'}(\Theta_{\xi,\zeta}^+),\Psi_{J\to J'}(\Theta_{\zeta,\sigma}^+)))
  \\ &+\hat{\d}_{\xi,\sigma} H(\Theta^+_{\xi,\zeta},\Theta^+_{\zeta,\sigma})+H(\hat{\d}_{\xi,\zeta}(\Theta_{\xi,\zeta}^+),\Theta_{\zeta,\sigma}^+)+H(\Theta_{\xi,\zeta}^+,\hat{\d}_{\zeta,\sigma}(\Theta_{\zeta,\sigma}^+))=0,
\end{split}  
  \label{eq:change-almost-complex-triangles}
  \end{equation}
  where $\Psi_{J\to J'}$ denotes the change of almost complex structure map, and $H$ is the map which counts index $-1$ triangles which are $J_s$-holomorphic for some $s\in (0,1)$. (Note that we are abusing notation slightly by writing $\Psi_{J\to J'}$ and writing $J$ for an almost complex structure on $\Sigma_0\times \Delta$, as well as an almost complex structure on the cylindrical ends).

Since the change of almost complex structure maps are graded, they must fix $\Theta_{\xi,\zeta}^+$, $\Theta_{\zeta,\sigma}^+$ and $\Theta_{\xi,\sigma}^+$. Since also the differentials on $\widehat{\CF}(\Sigma_0,\xis,\zetas)$, $\widehat{\CF}(\Sigma_0,\zetas,\sigmas)$ and $\widehat{\CF}(\Sigma_0,\xis,\sigmas)$ all vanish, we conclude from ~\eqref{eq:change-almost-complex-triangles} that if the main claim holds using $J$, then it also holds using $J'$.
 
Suppose now that $\cT_0$ is a multi-stabilizing triple with $m_3(\cT_0)=0$, obtained by adjoining $n+1$ elementary stabilizations. Let us say that $\cT_0$ is obtained by adjoining the elementary triple $\cE$ to a triple $\cT'$, obtained from $n$ elementary stabilizations. We pick the almost complex structure to be infinitely degenerated between $\cE$ and $\cT'$.  By the inductive hypothesis, the present claim holds for $\cT'$. We now view $\cT_0$ as an elementary stabilization of $\cT'$, using the triple $\cE$, and apply Proposition~\ref{prop:stabilize-triangles-general} to obtain the claim for $\cT_0$. The proof is complete. 
\end{proof}

\begin{rem}
If we symmetrize the triangle map,  by viewing it instead as a map
\begin{equation}
f_{\a,\b,\g;\frs}\colon\bCF^-(\Sigma,\as,\bs,\frs|_{Y_{\a,\b}})\otimes \bCF^-(\Sigma,\bs,\gs,\frs|_{Y_{\b,\g}})\otimes \bCF^-(\Sigma,\gs,\as, \frs|_{Y_{\g,\a}})\to \bF \llsquare U\rrsquare, \label{eq:symmetrized-triangle-maps}
\end{equation}
 then the three statements of Proposition~\ref{prop:stabilize-triangles-general} become equivalent.
\end{rem}

\subsection{Stabilizations and holomorphic quadrilaterals}

We now describe the analogs of the stabilization results from Section~\ref{sec:stabilizations-hol-triangles} for holomorphic quadrilaterals.

\begin{prop}\label{prop:stabilize-quadrilateral}
 Suppose that $\cQ=(\Sigma,\as,\bs,\gs,\ds,w)$ is a Heegaard quadruple, and $\cQ_0=(\Sigma_0,\xis,\zetas,\sigmas,\taus,p_0)$ is a multi-stabilizing Heegaard quadruple. We form the wedge sum $\Sigma\wedge \Sigma_0$ at $p_0$ and a point in $\Sigma\setminus (\as\cup \bs\cup \gs\cup \ds)$. Suppose and $J=(J_s)_{s\in [0,1]}$ and $I=(I_s)_{s\in [0,1]}$ are stratified families of almost complex structures on $\Sigma\times \Box$ and $\Sigma_0\times \Box$, respectively, for counting rectangles. With respect to the families $J$ and $J\wedge I$, we have the following:
 \begin{enumerate}
 \item \label{prop:stabilize:part-2} If $m_1(\cQ_0)=0$, then
 \begin{align*}
&h_{\cQ\wedge \cQ_0;\frs\# \frs_0}(\xs\times \Theta^-_{\xi,\zeta}, \ys\times \Theta^+_{\zeta,\sigma}, \zs\times \Theta^+_{\sigma,\tau})\\
=& h_{\cQ;\frs}(\xs, \ys, \zs)\otimes \Theta_{\xi,\tau}^-+\sum_{\substack{\ws\in \bT_{\a}\cap \bT_{\dt}\\
\Theta \in \bT_{\xi}\cap \bT_{\tau}\setminus \Theta_{\xi,\tau}^-}} C_{\xs,\ys,\zs,\ws,\Theta}\cdot \ws\times \Theta,
\end{align*}
for some $C_{\xs,\ys,\zs,\ws,\Theta}\in \bF[U]$. 
 \item If $m_2(\cQ_0)=0$,  then
 \begin{align*}
&h_{\cQ\wedge \cQ_0;\frs\# \frs_0}(\xs\times \Theta^+_{\xi,\zeta}, \ys\times \Theta^-_{\zeta,\sigma}, \zs\times \Theta^+_{\sigma,\tau})\\
=& h_{\cQ;\frs}(\xs, \ys, \zs)\otimes \Theta_{\xi,\tau}^-+\sum_{\substack{\ws\in \bT_{\a}\cap \bT_{\dt}\\
\Theta \in \bT_{\xi}\cap \bT_{\tau}\setminus \Theta_{\xi,\tau}^-}} C_{\xs,\ys,\zs,\ws,\Theta}\cdot \ws\times \Theta,
\end{align*}
for some $C_{\xs,\ys,\zs,\ws,\Theta}\in \bF[U]$. 
 \item If $m_3(\cQ_0)=0$, then
 \begin{align*}
&h_{\cQ\wedge \cQ_0;\frs\# \frs_0}(\xs\times \Theta^+_{\xi,\zeta}, \ys\times \Theta^+_{\zeta,\sigma}, \zs\times \Theta^-_{\sigma,\tau})\\
=& h_{\cQ;\frs}(\xs, \ys, \zs)\otimes \Theta_{\xi,\tau}^-+\sum_{\substack{\ws\in \bT_{\a}\cap \bT_{\dt}\\
\Theta \in \bT_{\xi}\cap \bT_{\tau}\setminus \Theta_{\xi,\tau}^-}} C_{\xs,\ys,\zs,\ws,\Theta}\cdot \ws\times \Theta,
\end{align*}
for some $C_{\xs,\ys,\zs,\ws,\Theta}\in \bF[U]$. 
 \item If $m_4(\cQ_0)=0$, then
 \[
h_{\cQ\wedge \cQ_0;\frs\# \frs_0}(\xs\times \Theta^+_{\xi,\zeta}, \ys\times \Theta^+_{\zeta,\sigma}, \zs\times \Theta^+_{\sigma,\tau})=h_{\cQ,J,\frs}(\xs,\ys,\zs)\otimes \Theta^+_{\xi,\tau}. 
 \]
 \end{enumerate}
\end{prop}
\begin{proof} We focus on the proof of the fourth claim, since it is notationally the simplest. The proof follows similarly to the proof of Proposition~\ref{prop:stabilize-triangles-general}, until~\eqref{eq:main-triangle-count-triangles} therein. The necessary substitute is 
\begin{equation}
\sum_{\substack{\psi_0\in \pi_2(\Theta_{\xi,\zeta}^+,\Theta_{\zeta,\sigma}^+,\Theta_{\sigma,\tau}^+,\Theta_{\xi,\tau}^+) \\ n_{p_0}(\psi_0)=k}} \# \cM(\psi_0,\tilde{\ve{d}})\equiv 1 \pmod{2}, \label{eq:main-count}
\end{equation}
for generic $\tilde{\ve{d}}\in \Int(\Box)^k\times (0,1)$.

Write $\cM(\tilde{\ve{d}})$ for the union over such $\psi_0$ of the matched moduli spaces $\cM(\psi_0,\tilde{\ve{d}})$ so $\# \cM(\tilde{\ve{d}})$ is equal to the left side of ~\eqref{eq:main-count}. We will use a cobordism argument to show ~\eqref{eq:main-count}.

We consider a path $(\tilde{\ve{d}}_t)_{t\in [0,\infty]}$ in $(\Int(\Box))^k\times (0,1)$ satisfying the following:
\begin{enumerate}
\item $\tilde{\ve{d}}_0=\tilde{\ve{d}}$.
\item The $(0,1)$-component (i.e. the almost complex structure parameter) of $\tilde{\ve{d}}_t$ approaches $1$ as $t\to \infty$.
\item The $k$ $\Box$-components of $\tilde{\ve{d}}_t$ approach $\infty$ in the $\a$-$\b$ cylindrical end of $\Box$. Furthermore, with respect to an identification of the cylindrical end with $[0,1]\times [0,\infty)$, the points of $\tilde{\ve{d}}_t$ approach some fixed, generic
\[
\ve{d}\in ((0,1)\times \R)^k,
\]
modulo overall translation by $\R$.
\end{enumerate}

There are three types of ends of the 1-dimensional moduli space $\bigcup_{t\in [0,\infty)} \cM(\tilde{\ve{d}}_t)$. The first type are the curves in $\cM(\tilde{\ve{d}})$ (i.e. the ones we want to count). The next type of end is curves which break off at finite $t\in (0,\infty)$. The final type is curves which appear as $t\to \infty$.

The degenerations appearing at finite $t$ can only correspond to strips breaking off (degenerations into a pair of triangles are prohibited since the aspect ratio is determined by $\tilde{\ve{d}}_t$ and is finite for $t\in [0,\infty)$). Ends corresponding to strip breaking occur in canceling pairs, since $\Theta_{\xi,\zeta}^+$, $\Theta_{\zeta,\sigma}^+$, $\Theta_{\sigma,\tau}^+$ and $\Theta_{\xi,\tau}^+$ are all cycles in the hat version of the complexes which contain them.

As $t\to \infty$, a family of curves approaches a pair of triangles not crossing the basepoints, as well as a holomorphic disk which matches $\ve{d}$.  As such, these ends correspond to the space
\begin{equation}
\begin{split}
\bigg(\coprod_{\zs\in \bT_{\xi}\cap \bT_{\sigma}} \bigg(\coprod_{\substack{\psi_{\xi,\zeta,\sigma}\in \pi_2(\Theta_{\xi,\zeta}^+,\Theta_{\zeta,\sigma}^+,\zs)\\
\mu(\psi_{\xi,\zeta,\sigma})=0\\
n_{p_0}(\psi_{\xi,\zeta,\sigma})=0}} \cM(\psi_{\xi,\zeta,\sigma}) \bigg)&\times  \bigg( \coprod_{\substack{\psi_{\xi,\sigma,\tau}\in \pi_2(\zs,\Theta_{\sigma,\tau}^+,\Theta_{\xi,\tau}^+)\\
\mu(\psi_{\xi,\sigma,\tau})=0\\
n_{p_0}(\psi_{\xi,\sigma,\tau})=0}} \cM(\psi_{\xi,\sigma,\tau})\bigg)\bigg) 
\\
&\times \bigg(\coprod_{\substack{\phi\in \pi_2(\Theta_{\xi,\zeta}^+,\Theta_{\xi,\zeta}^+)\\ n_{p_0}(\phi)=k }} \cM(\phi,\ve{d}) \bigg).
\end{split}
\label{eq:Md-cobordant-to}
\end{equation}
We wish to show the above space has an odd number of points. An identical argument to \cite{ZemDuality}*{Equation~24} establishes that the last factor of ~\eqref{eq:Md-cobordant-to} has an odd number of elements.

The first pair of factors corresponds exactly to the $\Theta_{\xi,\tau}^+$ component of
\[
\hat{f}_{\xi,\sigma,\tau}(\hat{f}_{\xi,\zeta,\sigma}(\Theta_{\xi,\zeta}^+,\Theta_{\zeta,\sigma}^+),\Theta_{\sigma,\tau}^+),
\]
where $\hat{f}_{\xi,\sigma,\tau}$ denotes the map which counts index 0 holomorphic triangles with zero multiplicity on $p_0$. Consequently, it is sufficient to show that
\begin{equation}
\hat{f}_{\xi,\sigma,\tau}(\hat{f}_{\xi,\zeta,\sigma}(\Theta_{\xi,\zeta}^+,\Theta_{\zeta,\sigma}^+),\Theta_{\sigma,\tau}^+)=\Theta_{\xi,\tau}^+.\label{eq:final-model-comp-triangle}
\end{equation}
To prove ~\eqref{eq:final-model-comp-triangle} it is sufficient to prove 
\begin{equation}
\hat{f}_{\xi,\sigma,\tau}(\Theta_{\xi,\sigma}^+,\Theta_{\sigma,\tau}^+)=\Theta_{\xi,\tau}^+,\quad \text{and}\quad \hat{f}_{\xi,\zeta,\sigma}(\Theta_{\xi,\zeta}^+,\Theta_{\zeta,\sigma}^+)=\Theta_{\zeta,\sigma}^+,
\end{equation}
which is proven in Lemma~\ref{lem:model-computation-stabilize-triangle}.
\end{proof}

\subsection{Lens space stabilizations of Heegaard diagrams}

Let $\cH$ be a Heegaard diagram, and let $\cH_0=(\bT^2,c_1,c_2,p_0)$ be a standard genus 1 diagram for $L(m,1)$. Let $J$ be an almost complex structure for $\cH$, and let $J_0$ be an almost complex structure on $\bT^2\times [0,1]\times \R$.

If $i\in \Z/m$, we define a map 
\[
\sigma_{i}\colon \bCF^-_{J}(\cH,\frs)\to \bCF^-_{J\wedge J_0}(\cH, \cH_0, \frs\# \frs_i),
\]
via the formula
\[
\sigma_i(\xs):=  \xs\times x_i,
\]
extended equivariantly over $U$.

\begin{figure}[ht!]
	\centering
	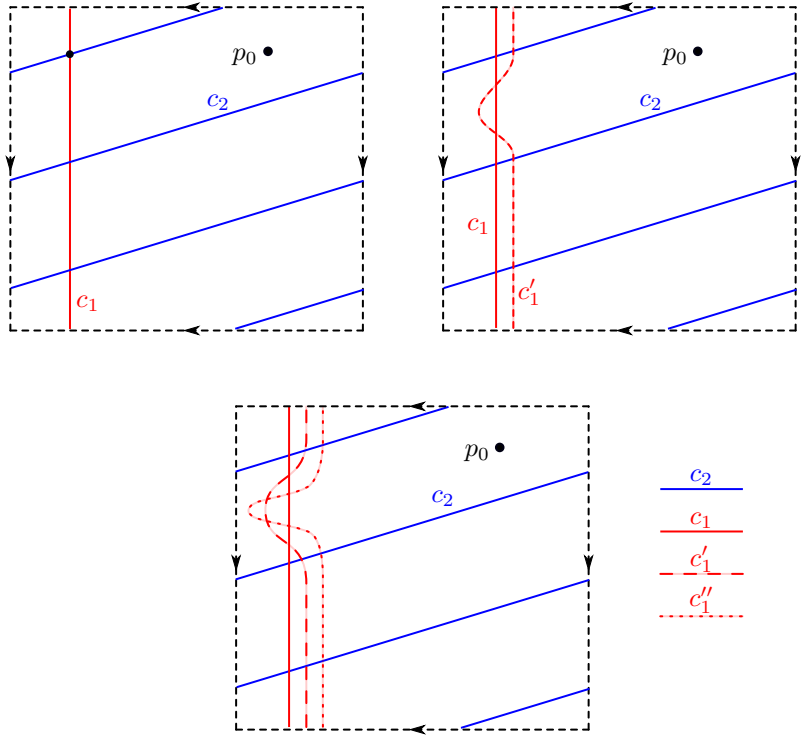
	\caption{Top left: a diagram for $L(m,1)$. Top right: a triple obtained by adding an isotopic copy of the alpha curve. Bottom: a triple obtained by adding two isotopic alpha curves.}\label{fig:179}
\end{figure}

\begin{lem} \label{lem:lens-space-stabilization-disks}
For each $i$, $\sigma_i$ is a chain map.
\end{lem}
\begin{proof}
Note that each intersection point on $\cH_0$ represents a different $\Spin^c$ structure. Hence any class of disks on $(\bT^2,c_1,c_2,p_0)$ must have the same incoming and outgoing intersection point. Furthermore, as $H_2(L(m,1);\Z)=0$, the set of doubly periodic domains vanishes, so $\pi_2(x_i,x_i)\iso \Z$. Furthermore, every class in $\pi_2(x_i,x_i)$ is of the form $k\cdot [\bT_2]+e_{x_i}$, where $e_{x_i}$ denotes the constant class. Write $\phi_{k,x_i}$ for this class. We note
\[
\mu(\phi_{k,x_i})=2k.
\]

Write $p$ for the connected sum point on $\cH$. If $\phi\in \pi_2(\xs,\ys)$ is a class on $\cH$, and $n_p(\phi)=k$,
\[
\mu(\phi\# \phi_{k,x_i})=\mu(\phi),
\]
using Lipshitz's formula for the Maslov index \cite{LipshitzCylindrical}*{Corollary~4.3}.

Hence, the main claim follows if we can show
\[
\# \cM\cM(\phi, \phi_{k,x_i})/\R\equiv \# \cM(\phi)/\R
\]
for any index 1 class $\phi$ with $n_{p}(\phi)=k$.

Hence, it is sufficient to show that if $\ve{d}=(d_1,\dots, d_k)$ is a $k$-tuple of distinct points in $[0,1]\times \R$, then
\begin{equation}
\# \cM(\phi_{k,x_i}, \ve{d})\equiv 1 \pmod 2, \label{eq:lens-space-stabilization-count}
\end{equation}
where $\cM(\phi_{k,x_i}, \ve{d})$ denotes the disjoint union of the moduli spaces
 $\{u\in \cM(S^{\qs},\phi_{k,x_i}): \ev_{q_i}(u)=d_i\}$, over equivalence classes of source curves $S$ with $k$ marked points $\qs$.

Equation~\eqref{eq:lens-space-stabilization-count} may be established using the same argument as in the $k=1$ case (i.e. ordinary stabilizations), as described by Lipshitz \cite{LipshitzCylindrical}*{Appendix A}. Alternatively, one may establish ~\eqref{eq:lens-space-stabilization-count} by degenerating $\ve{d}$ towards $k$ distinct points along $\{0\}\times \R$, and using the count of index 2 boundary degenerations for generic almost complex structures given in Proposition~\ref{prop:relaxed-count}.
\end{proof}

\subsection{Lens space stabilizations and holomorphic triangles}
\label{sec:lens-space-stabilizations}

Next, we consider triangles and lens space stabilizations. We consider taking the connected sum of a Heegaard triple with a triple of the form shown in Figure~\ref{fig:179}.

\begin{lem}
\label{lem:lens-stabilization-triangles} Suppose that $c_1',$ $c_1$ and $c_2$ are the three curves on the torus shown in Figure~\ref{fig:179}, and suppose that $\cT=(\Sigma,\as,\bs,\gs,w)$ is an arbitrary Heegaard triple. Write $\cT_{c_1',c_1,c_2}$ for the triple $(\bT^2,c_1',c_1,c_2,p_0)$, and write $\cT_{c_1,c_2,c_1'}$ and $\cT_{c_2,c_1',c_1}$ for the triples obtained by cyclically permuting the attaching curves. Suppose $J$ and $J_0$ are almost complex structures on $\Sigma\times \Delta$ and $\Sigma_0\times \Delta$, and let $J\wedge J_0$ denote their wedge. Then, counting triangles with respect to $J$ and $J\wedge J_0$, we have the following:
\begin{align*}
f_{\cT\wedge \cT_{c_1',c_1,c_2};\frs\# \frs_i}(\xs\times \Theta_{c_1',c_1}^+, \ys\times x_{i})&=f_{\cT,\frs}(\xs,\ys)\otimes x_i';\\
f_{\cT\wedge \cT_{c_1,c_2,c_1'};\frs\# \frs_i}(\xs\times x_i,\ys\times x_i')&=f_{\cT;\frs}(\xs,\ys)\otimes \Theta_{c_1,c_1'}^-\\
&+ \sum_{\zs\in \bT_{\a}\cap \bT_{\g}} C_{\xs,\ys,\zs}\cdot \zs\times \Theta_{c_1,c_1'}^+, \quad  \text{and}\\
f_{\cT\wedge \cT_{c_2,c_1',c_1};\frs\# \frs_i}(\xs\times x_i',\ys\times \Theta_{c_1',c_1}^+)&=f_{\cT;\frs}(\xs,\ys)\otimes x_i,
\end{align*}
for some $C_{\xs,\ys,\zs}\in \bF[U]$, and where $x_i'$ is the unique intersection point in the corresponding $\Spin^c$ structure.
\end{lem}
\begin{proof} All of the equations are equivalent, by reordering the attaching curves, and dualizing, so we focus on the first one. The $m=1$ case is used in the proof that the cobordism maps are invariants \cite{OSTriangles}*{Theorem~2.14}. Our proof is a modification of the proof of an analogous result in the cylindrical setting \cite{ZemGraphTQFT}*{Theorem~9.7}.

Suppose $\psi_0\in \pi_2(\Theta,x_i,x_i')$ is a class of triangles on $\cT_{c_1',c_1,c_2}$. We have
\[
\mu(\psi_0)=2n_{p_0}(\psi_0)+\gr(\Theta^+_{c_1',c_1},\Theta).
\]
 If $\psi\in \pi_2(\xs,\ys,\zs)$ is a class, it follows that
\begin{equation}
\ind_{\emb}(\psi,\psi_0,M_{\#})=\mu(\psi)+\gr(\Theta^+_{c_1',c_1},\Theta).\label{eq:index-connected-sum-stabilization}
\end{equation}

If the matched moduli space $\cM\cM(\psi,\psi_0)$ is non-empty, then $\cM(\psi)$ is non-empty, so we may assume $\mu(\psi)\ge 0$, by transversality. Equation~\eqref{eq:index-connected-sum-stabilization} implies that if $(\psi, \psi_0)$ has a matched representative, then $\Theta=\Theta^+_{c_1',c_1}$ and $\mu(\psi)=0$. The proof now follows by showing that if $\ve{d}\in (\mathrm{int}\;(\Delta))^k$ is generic, then
\begin{equation}
 \sum_{
\substack{\psi_0\in \pi_2(\Theta^+_{c_1',c_1}, x_i,x_i')\\
n_{p_0}(\psi_0)=k}} \# \cM(\psi_0,\ve{d})\equiv 1 \pmod{2}.\label{eq:matched-disks-stabilization}
\end{equation}
Equation~\eqref{eq:matched-disks-stabilization} is obtained by degenerating $\ve{d}$ into one of the cylindrical ends of $\bT^2\times \Delta$, as in the proof of ~\eqref{eq:main-triangle-count-triangles}. Equation~\eqref{eq:matched-disks-stabilization} then follows from ~\eqref{eq:total-count-matched-disks}, as well as the computation for the hat-triangle map (i.e. the count of holomorphic triangles not covering the basepoint), which is straightforward to perform by hand. The proof is complete.
\end{proof}

\subsection{Lens space stabilizations and holomorphic rectangles}

We will also need a version of the previous lens space stabilization results for rectangles.

\begin{lem}\label{lem:len-space-stabilization-quadrilaterals} Suppose that $\cQ$ is a Heegaard quadruple. Write $\cQ_{c_1'',c_1',c_1,c_2}$ for the genus 1 quadruple shown in Figure~\ref{fig:179}. Suppose that $J=(J_s)_{s\in [0,1]}$ and $I=(I_s)_{s\in [0,1]}$ are stratified families of almost complex structures, for counting rectangles on $\cQ$ and $\cQ_0$, respectively. With respect to the families of almost complex structures $J$ and $J\wedge I$, we have
\[
\begin{split}
&h_{\cQ\wedge \cQ_{c_1',c_1,c_2,c_1''}; \frs\# \frs_{0}}(\xs\times \Theta_{c_1',c_1}^+, \ys\times \Theta^{\can}_{c_1,c_2}, \zs\times \Theta_{c_2,c_1''}^+)\\
=&h_{\cQ;\frs}(\xs,\ys,\zs)\otimes \Theta_{c_1',c_1''}^-+\sum_{\ws\in \bT_{\a}\cap \bT_{\dt}} C_{\xs,\ys,\zs,\ws}\cdot \ws\times \Theta^+_{c_1',c_1''},
\end{split}
\]
for $C_{\xs,\ys,\zs,\ws}\in \bF[U]$.
Similar statements hold for quadruples which are stabilized by the genus 1 diagrams obtained by cyclically permuting $c_1''$, $c_1'$, $c_1$ and $c_2$.
\end{lem}

The proof is not substantially different than the proof of Lemma~\ref{lem:lens-stabilization-triangles}, so we leave it as an exercise to the reader.

\section{Holomorphic triangles and nearest point maps}
\label{sec:holotriangles}

In this section, we consider Heegaard triples where one set of attaching curves is obtained via a small isotopy of another. We show that in favorable circumstances, the counts of holomorphic triangles are appropriately simple. The main result is known to experts. An extremely similar statement and argument was given in a closely related context by Lipshitz, Ozsv\'{a}th and Thurston \cite{LOTDoubleBranchedII}*{Section~3.4}, on which our argument is modeled.

\begin{prop} \label{prop:nearest-point-triangles}
 Suppose that $(\Sigma,\as',\as,\bs)$ is a Heegaard triple and $\as'$ are Hamiltonian isotopies of $\as$, such that 
\[
|\a'\cap \a|=2\delta_{ij}.
\]
If $\as'$ are sufficiently close to $\as$, then the holomorphic triangle map $f_{\a',\a,\b}(\Theta^+_{\a',\a},-)$ counts only small triangles, and is equal to the nearest point map.
\end{prop}

\subsection{Vertex multiplicities}

We now review some basic facts about the domains of triangle classes. Suppose that $(\Sigma,\gs_1,\gs_2,\gs_3)$ is a Heegaard triple. If $\psi\in \pi_2(\xs,\ys,\zs)$ is a homology class of triangles, then we write $D(\psi)$ for the \emph{domain} of $\psi$, which is an integral 2-chain on $\Sigma$, with boundary on $\gs_1\cup \gs_2\cup \gs_3$.

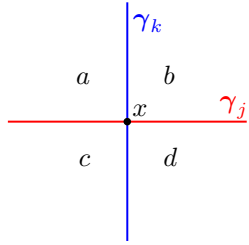
\begin{figure}[ht!]
	\centering
\begingroup%
  \makeatletter%
  \providecommand\color[2][]{%
    \errmessage{(Inkscape) Color is used for the text in Inkscape, but the package 'color.sty' is not loaded}%
    \renewcommand\color[2][]{}%
  }%
  \providecommand\transparent[1]{%
    \errmessage{(Inkscape) Transparency is used (non-zero) for the text in Inkscape, but the package 'transparent.sty' is not loaded}%
    \renewcommand\transparent[1]{}%
  }%
  \providecommand\rotatebox[2]{#2}%
  \newcommand*\fsize{\dimexpr\f@size pt\relax}%
  \newcommand*\lineheight[1]{\fontsize{\fsize}{#1\fsize}\selectfont}%
  \ifx\svgwidth\undefined%
    \setlength{\unitlength}{95.0530061bp}%
    \ifx\svgscale\undefined%
      \relax%
    \else%
      \setlength{\unitlength}{\unitlength * \real{\svgscale}}%
    \fi%
  \else%
    \setlength{\unitlength}{\svgwidth}%
  \fi%
  \global\let\svgwidth\undefined%
  \global\let\svgscale\undefined%
  \makeatother%
  \begin{picture}(1,0.94683991)%
    \lineheight{1}%
    \setlength\tabcolsep{0pt}%
    \put(0,0){\includegraphics[width=\unitlength,page=1]{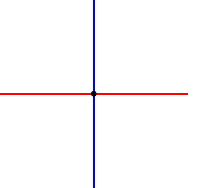}}%
    \put(0.49385036,0.49525954){\color[rgb]{0,0,0}\makebox(0,0)[lt]{\lineheight{1.25}\smash{\begin{tabular}[t]{l}$x$\end{tabular}}}}%
    \put(0.32653015,0.62206836){\color[rgb]{0,0,0}\makebox(0,0)[rt]{\lineheight{1.25}\smash{\begin{tabular}[t]{r}$a$\end{tabular}}}}%
    \put(0.61643227,0.62206836){\color[rgb]{0,0,0}\makebox(0,0)[lt]{\lineheight{1.25}\smash{\begin{tabular}[t]{l}$b$\end{tabular}}}}%
    \put(0.32653015,0.29941015){\color[rgb]{0,0,0}\makebox(0,0)[rt]{\lineheight{1.25}\smash{\begin{tabular}[t]{r}$c$\end{tabular}}}}%
    \put(0.61643227,0.29941015){\color[rgb]{0,0,0}\makebox(0,0)[lt]{\lineheight{1.25}\smash{\begin{tabular}[t]{l}$d$\end{tabular}}}}%
    \put(0.84079606,0.52567875){\color[rgb]{1,0,0}\makebox(0,0)[lt]{\lineheight{1.25}\smash{\begin{tabular}[t]{l}$\gs_j$\end{tabular}}}}%
    \put(0.4970707,0.8537836){\color[rgb]{0,0,1}\makebox(0,0)[lt]{\lineheight{1.25}\smash{\begin{tabular}[t]{l}$\gs_{k}$\end{tabular}}}}%
  \end{picture}%
\endgroup%

	\caption{Multiplicities near an intersection point $x$. We assume $j<k$.}\label{fig:45}
\end{figure}

The domain $D(\psi)$ satisfies \emph{vertex relations} at the intersection points on the Heegaard diagram, which we now describe. Suppose $x\in \Sigma$ is in the intersection of a curve in $\gs_j$ and a curve in $\gs_k$. We can label the multiplicities of the four quadrants near $x$ by $a$, $b$, $c$ and $d$, as shown in Figure~\ref{fig:45}.

The vertex relations are
\[
a+d=c+b+\delta_x(\xs)+\delta_x(\ys)-\delta_x(\zs),
\]
where 
\[
\delta_{x}(\xs)=
\begin{cases}
1& \text{ if } x\in \xs\\
0& \text{ otherwise},
\end{cases}
\]
and similarly for $\delta_x(\ys)$ and $\delta_x(\zs)$.

Any domain which satisfies the vertex relations is a homology class of triangles.

\subsection{Approximating triangle classes}

Suppose $(\as',\as,\bs)$ is a Heegaard triple where $\as'$ consists of Hamiltonian isotopies of $\as$, satisfying
\[
|\a'_i\cap \a_j|=2\delta_{ij}.
\]
We assume all intersections of the attaching curves are transverse. Furthermore, we assume that each curve in $\bs$ intersects the small regions between $\a_k$ and $\a_k'$ in a collection of arcs which each  connect $\a_k$ and $\a_k'$.

If $\xs\in \bT_{\a}\cap \bT_{\b}$, then there is a nearest point $\xs_{np}\in \bT_{\a'}\cap \bT_{\b}$. Given a homology class of triangles $\psi\in \pi_2(\Theta,\xs,\ys_{np})$, there is a unique domain $D_{\app}(\psi)$ on $(\Sigma,\as,\bs)$ which agrees with $D(\psi)$ outside of the regions between $\as$ and $\as'$.

The following is a straightforward computation:
\begin{lem}
If $\xs\in \bT_{\a}\cap \bT_{\b}$ and $\Theta\in \bT_{\a'}\cap \bT_{\a}$, and $\psi\in \pi_2(\Theta,\xs,\ys_{np})$, then the domain $D_{\app}(\psi)$ satisfies the vertex relations for a class in $\pi_2(\xs,\ys)$.
\end{lem}

\begin{figure}[ht!]
	\centering
\begingroup%
  \makeatletter%
  \providecommand\color[2][]{%
    \errmessage{(Inkscape) Color is used for the text in Inkscape, but the package 'color.sty' is not loaded}%
    \renewcommand\color[2][]{}%
  }%
  \providecommand\transparent[1]{%
    \errmessage{(Inkscape) Transparency is used (non-zero) for the text in Inkscape, but the package 'transparent.sty' is not loaded}%
    \renewcommand\transparent[1]{}%
  }%
  \providecommand\rotatebox[2]{#2}%
  \newcommand*\fsize{\dimexpr\f@size pt\relax}%
  \newcommand*\lineheight[1]{\fontsize{\fsize}{#1\fsize}\selectfont}%
  \ifx\svgwidth\undefined%
    \setlength{\unitlength}{77.47314872bp}%
    \ifx\svgscale\undefined%
      \relax%
    \else%
      \setlength{\unitlength}{\unitlength * \real{\svgscale}}%
    \fi%
  \else%
    \setlength{\unitlength}{\svgwidth}%
  \fi%
  \global\let\svgwidth\undefined%
  \global\let\svgscale\undefined%
  \makeatother%
  \begin{picture}(1,2.19209343)%
    \lineheight{1}%
    \setlength\tabcolsep{0pt}%
    \put(0,0){\includegraphics[width=\unitlength,page=1]{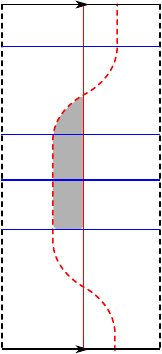}}%
    \put(0.73645712,0.18365263){\color[rgb]{1,0,0}\makebox(0,0)[lt]{\lineheight{1.25}\smash{\begin{tabular}[t]{l}$\a'_j$\end{tabular}}}}%
    \put(0.49776683,0.18243107){\color[rgb]{1,0,0}\makebox(0,0)[rt]{\lineheight{1.25}\smash{\begin{tabular}[t]{r}$\a_j$\end{tabular}}}}%
    \put(0.55870391,1.53341034){\color[rgb]{0,0,0}\makebox(0,0)[lt]{\lineheight{1.25}\smash{\begin{tabular}[t]{l}$\Theta_{\a',\a}^+$\end{tabular}}}}%
    \put(0.29676298,0.8191788){\color[rgb]{0,0,0}\makebox(0,0)[rt]{\lineheight{1.25}\smash{\begin{tabular}[t]{r}$\xs_{np}$\end{tabular}}}}%
    \put(0.54704251,0.8191788){\color[rgb]{0,0,0}\makebox(0,0)[lt]{\lineheight{1.25}\smash{\begin{tabular}[t]{l}$\xs$\end{tabular}}}}%
    \put(0,0){\includegraphics[width=\unitlength,page=2]{fig157.pdf}}%
    \put(0.78912658,0.64243344){\color[rgb]{0,0,1}\makebox(0,0)[lt]{\lineheight{1.25}\smash{\begin{tabular}[t]{l}$\b$\end{tabular}}}}%
  \end{picture}%
\endgroup%

	\caption{A small triangle class.}\label{fig:157}
\end{figure}

Hence, $D_{\app}(\psi)$ determines a class of disks  $\psi_{\app}\in \pi_2(\xs,\ys)$. The following index formula is key to proving Proposition~\ref{prop:nearest-point-triangles}:

\begin{lem}\label{lem:Maslov-index}
Suppose that $(\Sigma,\as',\as,\bs,w)$ is a Heegaard triple, and $\as'$ are small Hamiltonian isotopies of $\as$, as above. Let $\psi\in \pi_2(\Theta,\xs,\ys_{np})$ be a homology class of triangles, and $\psi_{\app}\in \pi_2(\xs,\ys)$ be its approximation. Then
\[
\mu(\psi_{\app})=\mu(\psi)+\gr(\Theta_{\a',\a}^+,\Theta).
\]
\end{lem}
\begin{proof}  Firstly, it is sufficient to show the claim when $\Theta=\Theta_{\a',\a}^+$, since the general claim follows by splicing bigons into $\psi$ to reduce to this case. Hence, we assume $\Theta=\Theta_{\a',\a}^+$.

Noting that $\psi_{\app}\in \pi_2(\xs,\ys)$, we may splice to obtain a triangle class
\begin{equation}
\psi_0:=\psi-\psi_{\app}\in \pi_2(\Theta^+_{\a',\a},\ys,\ys_{np}). \label{eq:maslov-index-approximation}
\end{equation}
The class $\psi_0$ is supported in the regions between $\as$ and $\as'$.  The difference between $\psi_0$ and a canonical small triangle class is a doubly periodic domain which is contained in the region between $\as$ and $\as'$. See Figure~\ref{fig:157}. Such doubly periodic domains have index 0, so
\[
\mu(\psi_0)=0.
\]
Using additivity of the index and~\eqref{eq:maslov-index-approximation}, we obtain $\mu(\psi)=\mu(\psi_0)$, completing the proof.
\end{proof}

\begin{figure}[ht!]
	\centering
	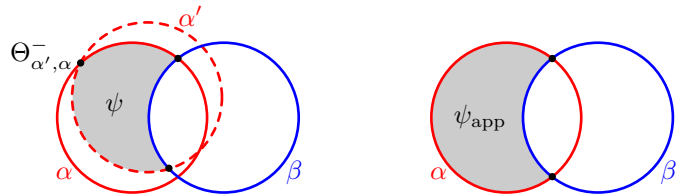
	\caption{A Maslov index 0 triangle class (left), whose approximation has Maslov index 1 (right). Note that the triangle has a vertex at $\Theta_{\a',\a}^-$, and not $\Theta_{\a',\a}^+$.}\label{fig:46}
\end{figure}

In Figure~\ref{fig:46}, we  give an example of a Maslov index 0 class of triangles whose approximation is a Maslov index 1 class of disks. 

\subsection{Proof of Proposition~\ref{prop:nearest-point-triangles}}

\begin{proof}
The proof follows from essentially the same argument as \cite{LOTDoubleBranchedII}*{Lemma~3.38}, as we now sketch. Let $\as'_{i}$ be a sequence of Hamiltonian translates of $\as$, each of which is sufficiently close to $\as$ for the nearest point map to be well defined, and such that $\as'_i\to \as'$ as $i\to \infty$.
 
By picking $\as_i$ suitably, we may assume that the sets of homology classes of triangles on $(\Sigma,\as_i',\as,\bs)$ are canonically identified for different $i$. Let $\Theta^+_i$ denote the top intersection point of $\bT_{\a_i'}\cap \bT_{\a}$. Suppose that $\ys_i'\in \bT_{\a_i'}\cap \bT_{\b}$ are similarly identified, and $\xs\in \bT_{\a}\cap \bT_{\b}$. Suppose that $\psi\in \pi_2(\Theta_i^+,\xs,\ys_i')$ is a Maslov index 0 class of triangles, which has a holomorphic representative for all $i$.

Let us label the edges of the triangle $\Delta$ as $e_{\a'}$, $e_{\a}$ and $e_{\b}$. Analogously to the argument in \cite{LOTDoubleBranchedII}*{Lemma~3.38}, we may extract a Gromov limit of $u_i$, consisting of a collection of curves mapping into $\Sigma\times \Delta$, with boundary on $\as\times e_{\a'}$, $\as\times e_{\a}$ and $\bs\times e_{\b}$, as well as a collection of curves mapping into $\Sigma\times [0,1]\times \R$, corresponding to the three cylindrical ends. We organize all of the components whose sources are not closed, and which project non-constantly to $\Delta$ into a collection $u$. If $u$ is a $J$-holomorphic curve mapping into $\Sigma\times \Delta$ (or $\Sigma\times [0,1]\times \R$) for some $J$ satisfying \eqref{J1}--\eqref{J5}, and $\pi_{\Delta}\circ u$ is locally non-constant, then it follows from \cite{LipshitzCylindrical}*{Lemma~3.1} that $\pi_{\Delta}\circ u$ is an open map (in the notation therein, $\Sigma\times \Delta$ is $E$ and $\Sigma$ is $B$). Hence any closed components  or components with constant projection to $\Sigma$ contribute domain $k[\Sigma]$ for some $k\in \Z$. Positivity of intersections implies that $k\ge 0$.  Similarly, the curves mapping into the cylindrical ends also have to consist of curves satisfying~\ref{M-1}--\ref{M-7}, as well as curves whose domain consists of nonnegative multiples of $[\Sigma]$. Using transversality for holomorphic disks, we may conclude that curves mapping into $\Sigma\times [0,1]\times \R$ have nonnegative Maslov index. 

We may view the curve $u$, mapping into $\Sigma\times \Delta$, as instead mapping into $\Sigma\times [0,1]\times \R\setminus \Sigma\times \{(1,0)\}$, with additional punctures which are asymptopic to points along $\as\times \{(1,0)\}$. Using  removal of singularities \cite{McDuffSalamonSymplectic}*{Theorem~4.1.2}, $u$ may be completed over $\Sigma\times \{(1,0)\}$, to obtain a holomorphic curve mapping into $\Sigma\times [0,1]\times \R$, satisfying~\ref{M-1}--\ref{M-7}. We write $u$ also for this holomorphic curve. We write $\phi_0$ for the class of $u$. The domain of $\phi_0$ differs from $\psi_{\app}$ by $k[\Sigma]$, for some $k\ge 0$, as well as a collection of disk classes corresponding to the curves in the $\alpha$-$\beta$ and $\alpha'$-$\beta$ cylindrical ends of $\Sigma\times \Delta$, which have nonnegative Maslov index. By Lemma~\ref{lem:Maslov-index}, $\mu(\psi_{\app})=0$. We conclude that $\mu(\phi_0)\le 0$, with equality if and only if the remaining classes have trivial domain. Since $\phi_0$ has a holomorphic representative, $\phi(\phi_0)\ge 0$. Transversality considerations imply that $\phi_0$ has trivial domain, and that the remaining curves also have trivial domain. 

It follows that the original triangle class $\psi$ had domain supported in the thin regions between $\as'$ and $\as$, and consequently must be a small triangle class. On the other hand, it is straightforward to enumerate small triangles by hand, and see that the count of small triangles gives exactly the nearest point map. It follows that $f_{\a',\a,\b}(\Theta_{\a',\a}^+,-)$ is the nearest point map.
\end{proof}

\begin{rem}
 The nearest point count of Proposition~\ref{prop:nearest-point-triangles} will be used in the construction of the central hypercube in Sections~\ref{sec:central-hypercube} and~\ref{sec:N-null-homotopy}. We note that a result of Lipshitz \cite{LipshitzCylindrical}*{Proposition~11.4} implies that the nearest point map is chain homotopic to the transition map obtained by counting triangles. It is possible to describe our argument without the nearest point map counts, by adding these homotopies as length 2 maps in the central hypercube. However, a tradeoff to including these homotopies in the central hypercube is a further increase in complexity of the hypercube in Section~\ref{sec:N-null-homotopy}. Hence, we opt to use the nearest point map to simplify the presentation.
\end{rem}

\subsection{Small isotopies and Heegaard quadruples}

Proposition~\ref{prop:nearest-point-triangles} has a natural analog for higher counts of holomorphic curves, and the proof is essentially the same.

\begin{prop}\label{prop:small-rectangles}
 Suppose that $(\Sigma,\as',\as,\bs,\gs,w)$ is a Heegaard quadruple such that $\as'$ are small Hamiltonian translates of $\as$, satisfying
\[
|\a'\cap \a|=2\delta_{ij}.
\]
If $\as'$ are sufficiently small translates of $\as$, then for a suitable family of almost complex structures on $\Sigma\times \Box$, we have $h_{\a',\a,\b,\g;\frs}(\Theta^+_{\a',\a},\xs,\ys)=0$, for all $\xs$ and $\ys$. In fact, the moduli spaces counted by $h_{\a',\a,\b,\g;\frs}(\Theta_{\a',\a}^+,\xs,\ys)$ are empty.
\end{prop}
The proof of Proposition~\ref{prop:small-rectangles} is largely the same as the proof of Proposition~\ref{prop:nearest-point-triangles}, as we sketch presently. We consider a sequence $\as_i'$ which approaches $\as$. The moduli space of rectangles may be identified with the moduli space of triangles with an extra marked point on one of the boundary components. Hence,  from a sequence of Maslov index $-1$ rectangles which have boundary on $\as_{i}',$ $\as,$ $\bs$ and $\gs$, we may extract a Gromov limit, which consists of a broken triangle on $(\Sigma,\as,\bs,\gs,w)$, with Maslov index $-1$. Generically,  no such triangles exist, by transversality, so the moduli spaces counted by $h_{\a_i',\a,\b,\g;\frs}(\Theta_{\a'_i,\a}^+,\xs,\ys)$ must be empty for large $i$.

\section{Input from Ozsv\'{a}th and Szab\'{o}'s surgery exact sequence}
\label{sec:previous-results}
In this section, we review some background on Ozsv\'{a}th and Szab\'{o}'s surgery exact sequence \cite{OSIntegerSurgeries}. In particular,  we describe the maps which appear along the top and bottom faces of our main hypercube.

\subsection{Surgery quadruples}

\begin{define}\label{def:surgery-quadruple} Suppose $K$ is a null-homologous knot in a 3-manifold $Y$. We say a doubly pointed Heegaard quadruple $(\Sigma,\as_3,\as_2,\as_1,\bs,w,z)$ is an $(n,n+m,\infty)$-\emph{surgery quadruple diagram} for $(Y,K)$ if the following are satisfied:
\begin{enumerate}[ref=$D$-\arabic*, label=($D$-\arabic*),leftmargin=*,widest=IIII]
\item We can write \[
\as_1=\as\cup \{c_1\}, \quad \as_2=\as'\cup \{c_2\}, \quad\text{and} \quad \as_3=\as''\cup \{c_3\},
\]
where $\as$, $\as'$ and $\as''$ are pairwise Hamiltonian translates of each other. Furthermore,
\[
|\a_i\cap \a_j'|=|\a'_i\cap \a''_j|=|\a_i\cap \a''_j|=2\delta_{ij}.
\]
\item  There is a once punctured torus $F\subset \Sigma$, which contains $c_1$, $c_2$ and $c_3$ and is disjoint from $\as$, $\as'$ and $\as''$. We call $F$ the \emph{surgery region}.
\item The diagram $(\Sigma,\as_1,\bs,w,z)$ represents a doubly based unknot in $Y_n(K)$. Furthermore, $w$ and $z$ are immediately adjacent. 
\item The diagram $(\Sigma,\as_2,\bs,w,z)$ represents a doubly based unknot in $Y_{n+m}(K)$. Furthermore, $w$ and $z$ are immediately adjacent.
\item The diagram $(\Sigma,\as_3,\bs,w,z)$  represents $K$ in $Y$, decorated with two basepoints.
\item The curves $c_1$, $c_2$ and $c_3$ may be oriented so that
\[
\#(c_1\cap c_2)=m,\quad\text{and} \quad \#(c_2\cap c_3)=\#(c_3\cap c_1)=1.
\]
\end{enumerate}
\end{define}

\begin{figure}[H]
	\centering
	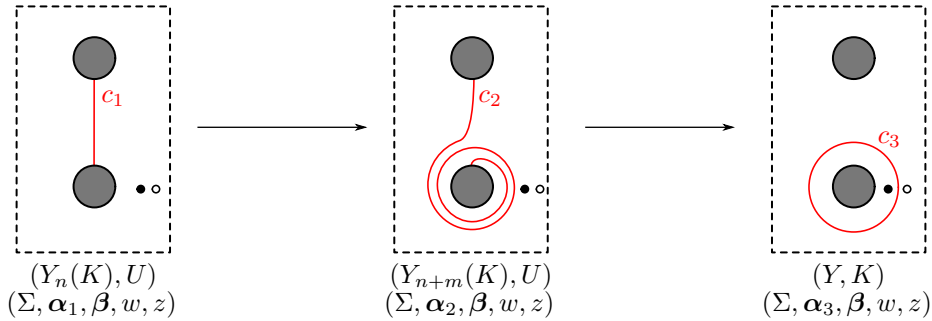
	\caption{The  genus 1 surgery region of an $(n,n+m,\infty)$-surgery quadruple $(\Sigma,\as_3,\as_2,\as_1,\bs,w,z)$. The curves $\bs$ are not shown.}\label{fig:70}
\end{figure}

 See Figure~\ref{fig:70} for a schematic of the surgery region of $\Sigma$.

\subsection{Heegaard Floer homology of lens spaces}

In this section, we recall some facts about the Floer homology of lens spaces, which Ozsv\'{a}th and Szab\'{o} use to build their exact sequence in \cite{OSIntegerSurgeries}.

\begin{define} View $L(m,1)$ as the boundary of a disk bundle $D(m,1)$ over $S^2$ with Euler number $m$. 
 The \emph{canonical $\Spin^c$ structure}, $\frs_0\in \Spin^c(L(m,1))$, is the restriction to $L(m,1)$ of the two $\Spin^c$ structures $\frt\in \Spin^c(D(m,1))$  satisfying
\[
\langle c_1(\frt),S^2 \rangle=\pm m.
\]
\end{define}

 The Floer homology of a lens space has a simple form:
 \[
\HF^-(L(m,1), \frs)\iso \bF[U],
 \]
 for each $\frs\in \Spin^c(L(m,1))$. We write $\Theta^{\can}$ for the generator of $\HF^-(L(m,1),\frs_0)$.

\subsection{Cobordism maps for disk bundles over \texorpdfstring{$S^2$}{S2}}
\label{sec:surgery-m-framed-unknot-1}

Let $m$ be a positive integer. We recall that the induced map on Heegaard Floer homology for $-m$ surgery on an unknot has a simple form. See \cite{OSIntegerSurgeries}*{Section~3}. The blow-up formula corresponds to $m=1$ \cite{OSTriangles}*{Theorem~3.7}.

\begin{figure}[ht!]
	\centering
\begingroup%
  \makeatletter%
  \providecommand\color[2][]{%
    \errmessage{(Inkscape) Color is used for the text in Inkscape, but the package 'color.sty' is not loaded}%
    \renewcommand\color[2][]{}%
  }%
  \providecommand\transparent[1]{%
    \errmessage{(Inkscape) Transparency is used (non-zero) for the text in Inkscape, but the package 'transparent.sty' is not loaded}%
    \renewcommand\transparent[1]{}%
  }%
  \providecommand\rotatebox[2]{#2}%
  \newcommand*\fsize{\dimexpr\f@size pt\relax}%
  \newcommand*\lineheight[1]{\fontsize{\fsize}{#1\fsize}\selectfont}%
  \ifx\svgwidth\undefined%
    \setlength{\unitlength}{156.98577425bp}%
    \ifx\svgscale\undefined%
      \relax%
    \else%
      \setlength{\unitlength}{\unitlength * \real{\svgscale}}%
    \fi%
  \else%
    \setlength{\unitlength}{\svgwidth}%
  \fi%
  \global\let\svgwidth\undefined%
  \global\let\svgscale\undefined%
  \makeatother%
  \begin{picture}(1,0.91997254)%
    \lineheight{1}%
    \setlength\tabcolsep{0pt}%
    \put(0,0){\includegraphics[width=\unitlength,page=1]{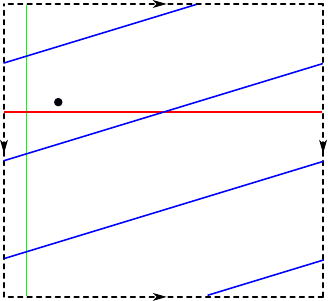}}%
    \put(0.11107228,0.33183338){\color[rgb]{0,0.8627451,0}\makebox(0,0)[lt]{\lineheight{1.25}\smash{\begin{tabular}[t]{l}$c_1$\end{tabular}}}}%
    \put(0.7317077,0.52433252){\color[rgb]{1,0,0}\makebox(0,0)[t]{\lineheight{1.25}\smash{\begin{tabular}[t]{c}$c_3$\end{tabular}}}}%
    \put(0.66810876,0.34735197){\color[rgb]{0,0,1}\makebox(0,0)[rt]{\lineheight{1.25}\smash{\begin{tabular}[t]{r}$c_2$\end{tabular}}}}%
    \put(0.09945948,0.69520275){\color[rgb]{0,0,0}\makebox(0,0)[lt]{\lineheight{1.25}\smash{\begin{tabular}[t]{l}$\Theta_{c_2,c_1}^{\can}$\end{tabular}}}}%
    \put(0,0){\includegraphics[width=\unitlength,page=2]{fig63.pdf}}%
    \put(0.2130278,0.59931374){\color[rgb]{0,0,0}\makebox(0,0)[lt]{\lineheight{1.25}\smash{\begin{tabular}[t]{l}$w$\end{tabular}}}}%
    \put(0,0){\includegraphics[width=\unitlength,page=3]{fig63.pdf}}%
    \put(0.2130278,0.52887315){\color[rgb]{0,0,0}\makebox(0,0)[lt]{\lineheight{1.25}\smash{\begin{tabular}[t]{l}$z$\end{tabular}}}}%
  \end{picture}%
\endgroup%

	\caption{The triple  $(\bT^2,c_3,c_2,c_1,w,z)$. The generator of the canonical $\Spin^c$ structure is labeled.}\label{fig:63}
\end{figure}

\begin{figure}[ht!]
	\centering
	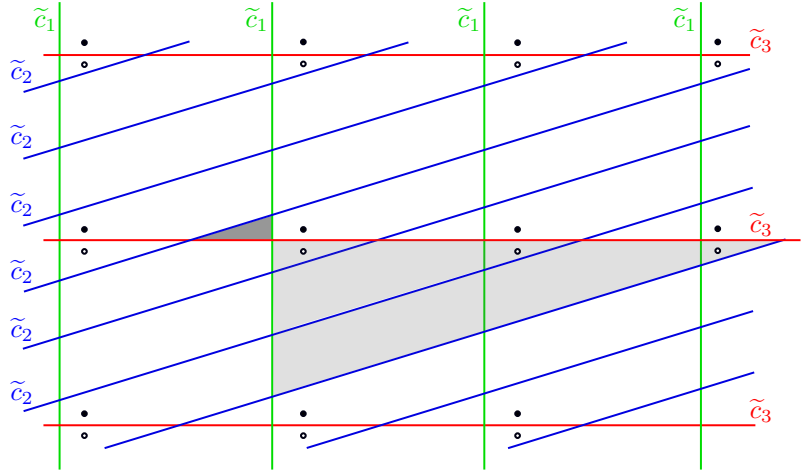
	\caption{The universal cover of Figure~\ref{fig:63}, as well as the lifts of two triangle classes. The curves $\tilde{c}_1$, $\tilde{c}_2$ and $\tilde{c}_3$ are the lifts of $c_1$, $c_2$ and $c_3$.}\label{fig:112}
\end{figure}

\begin{lem}\label{lem:model-computations-m-surgery} Let $(\bT^2,c_3,c_2,c_1,w,z)$ denote the genus 1 Heegaard triple shown in Figure~\ref{fig:63}. The twisted triangle map satisfies
\[
f_{c_3,c_2,c_1;\frs_k^{\pm}}(\Theta_{c_3,c_2},\Theta_{c_2,c_1}^{\can})=U^{mk(k+1)/2}\cdot T^0\cdot  \Theta_{c_3,c_1},
\]
where $\frs_k^{\pm}$  satisfies $\langle c_1(\frs_k^{\pm}),S \rangle=\pm (2k+1)m$, for $k\ge 0$. In particular, the sum over all $\frs_k^{+}$ and $\frs_k^-$ is 0.
\end{lem}
\begin{proof} 
A holomorphic triangle counted by $f_{c_3,c_2,c_1,\frs_k^{\pm}}$ may be lifted to the universal cover of $\bT^2$, where the count is easy to perform. See Figure~\ref{fig:112}.
\end{proof}

It will sometimes be helpful to index $\Spin^c$ structures by the odd integers $s\in 2 \Z+1$. If $u_s$ denotes the triangle representing the $\Spin^c$ structure with $\langle c_1(\frs), S\rangle= sm$, then one may compute directly that
\[
n_{w}(u_s)=\frac{m(s^2-1)}{8}\quad \text{and} \quad n_{z}(u_s)=\frac{m((s+2)^2-1)}{8}.
\]

\subsection{The length 1 maps on the top and bottom  of the main hypercube}

We now describe the length 1 maps along the top and bottom of the involutive hypercube of a triad in Theorem~\ref{thm:main-hypercube}. Suppose that $(\Sigma,\as_3,\as_2,\as_1,\bs,w,z)$ is an $(n,n+m,\infty)$ surgery quadruple for $K\subset Y$. We will write $f_{\a_1\to \a_2}^{\b}$ for the map labeled $f_1$ in Theorem~\ref{thm:main-hypercube}, which we define to be
\[
f_{\a_1\to \a_2}^{\b}(\xs)=\sum_{\substack{\ys\in \bT_{\a_2}\cap \bT_{\b} \\ \psi\in \pi_2(\Theta^{\can}_{\a_2,\a_1},\xs, \ys)
\\ \mu(\psi)=0}} \# \cM(\psi)\cdot U^{n_{w}(\psi)}\cdot \ys,
\]
extended equivariantly over the action of $U$. 

Similarly, we write $f_{\a_2\to \a_3}^{\b}$ for the map labeled $f_2$, which we define as
\[
f_{\a_2\to \a_3}^{\b}(\xs)=\sum_{\substack{\ys\in \bT_{\a_3}\cap \bT_{\b} \\ \psi\in \pi_2(\Theta_{\a_3,\a_2}^+,\xs, \ys)
\\ \mu(\psi)=0}}\# \cM(\psi)\cdot U^{n_{w}(\psi)} T^{n_{z}(\psi)-n_w(\psi)}\cdot \ys.
\]

If we wish to emphasize the almost complex structure, we will write
\[
{}^J\! f_{\a_1\to \a_2}^{\b}\quad \text{and} \quad {}^J\! f_{\a_2\to \a_3}^{\b}.
\]

\subsection{The length 2 maps along the top and bottom of the main hypercube}

The following is a well known fact about Heegaard Floer homology; see, e.g., \cite{OSIntegerSurgeries}:

\begin{lem}\label{lem:omega-a-null-homotopy} Suppose that $(\Sigma,\as_3,\as_2,\as_1,\bs,w,z)$ is an $(n,n+m,\infty)$-surgery quadruple. Let $S$ denote the boundary of the surgery region.  If $J=(J_t)_{t\in [0,1]}$ is a stratified  family of almost complex structures on $\Sigma\times \Box$ which is maximally pinched along $S$, then the diagram in ~\eqref{eq:2-d-hypercube-null-homotopy} is a 2-dimensional hypercube of chain complexes:
\begin{equation}
\begin{tikzcd}[column sep=2cm,row sep=2cm,labels=description]
\bCF^-_{J}(\as_1, \bs)
\arrow[d,swap, "{}^{J}\!f_{\a_1\to \a_2}^{\b} "]
\arrow[dr,dashed, "{}^J \!h_{\a_1\to \a_2\to \a_3}^{\b}",sloped]
&
\,
  \\
\bCF^-_{J}(\as_2,\bs) 
\arrow[r,swap, "{}^J\! f_{\a_2\to \a_3}^{\b}"]
&  \buCF^-_{J}(\as_3,\bs).
\end{tikzcd}
\label{eq:2-d-hypercube-null-homotopy}
\end{equation}
Here ${}^{J}\!h_{\a_1\to \a_2\to \a_3}^{\b}$ denotes the twisted rectangle counting map
\[
{}^J\! h_{\a_1\to \a_2\to \a_3}^{\b}(\xs):=h_{\a_3,\a_2,\a_1,\b}^J(\Theta_{\a_3,\a_2}^+,\Theta_{\a_2,\a_1}^{\can}, \xs).
\]
\end{lem}
\begin{proof} Using the associativity relations for quadrilaterals, it is sufficient to show
\begin{equation}
f_{\a_3,\a_2,\a_1}^J(\Theta_{\a_3,\a_2}^+, \Theta_{\a_2,\a_1}^{\can})=0. \label{eq:Fa3a2a1=0}
\end{equation}
Equation~\eqref{eq:Fa3a2a1=0} is proven as follows. Write $\as_3=\{c_3\}\cup \as''$, $\as_2=\{c_2\}\cup \as'$ and $\as_1=\{c_1\}\cup \as$, where $\as$, $\as'$ and $\as''$ are small Hamiltonian translates of each other. Using Proposition~\ref{prop:stabilize-triangles-general} we obtain
\[
f_{\a_3,\a_2,\a_1}^J(\Theta_{\a_3,\a_2}^+, \Theta_{\a_2,\a_1}^{\can})=f_{c_3,c_2,c_1}^J(\Theta_{c_3,c_2}^+,\Theta_{c_2,c_1}^{\can})\otimes
\Theta_{\a'',\a}^+,
\]
which implies ~\eqref{eq:Fa3a2a1=0} when combined with Lemma~\ref{lem:model-computations-m-surgery}.
\end{proof}

\begin{rem}
 Lemma~\ref{lem:omega-a-null-homotopy} holds without the family $J$ being maximally pinched, though the proof is easier if we assume that $J$ is maximally pinched. Despite being unnecessary for Lemma~\ref{lem:omega-a-null-homotopy}, the maximally pinched condition makes the construction of other portions of the main hypercube substantially easier. As a specific example, see the change of almost complex structure hypercube in  Section~\ref{sec:hypercubes-almost-complex-structure}, where the maximally pinched condition prevents curves from appearing along the left side of Figure~\ref{fig:111}.
\end{rem}

\subsection{Heegaard quadruples related to disk bundles}

Analogous to Lemma~\ref{lem:model-computations-m-surgery}, we will need a model count of holomorphic rectangles appearing on quadruples related to those in Figure~\ref{fig:112}.

\begin{lem}\label{lem:genus-1-quadruple-count}
 Suppose that $\cQ=(\bT^2,c_3,c_2,c_1,c_1',w,z)$ is the quadruple obtained from the triple $(\bT^2,c_3,c_2,c_1,w,z)$ shown in Figure~\ref{fig:63} by adding a curve $c_1'$ which is a small Hamiltonian translate of $c_1$, and intersects $c_1$ in two points.
 Then
 \[
h_{\cQ;\frs}(\Theta_{c_3,c_2}, \Theta_{c_2,c_1}^{\can}, \Theta_{c_1,c_1'}^+)=0,
 \]
 for all $\frs\in \Spin^c(X_{c_3,c_2,c_1,c_1'})$. In fact, there are no classes of index $-1$ in 
 \[
 \pi_2(\Theta_{c_3,c_2},\Theta_{c_2,c_1'}^{\can},\Theta_{c_1,c_1'}^+, \Theta_{c_3,c_1'}).
 \]
  Analogous claims hold for the quadruples $(\bT^2,c_3,c_2',c_2,c_1)$ and $(\bT^2,c_3',c_3,c_2,c_1)$.
\end{lem}
\begin{proof} First, we note that if we fill the $(c_3,c_2)$, $(c_1,c_1')$ and $(c_3,c_1')$ ends of $X_{c_3,c_2,c_1,c_1'}$ with 3- and 4-handles, we obtain a disk bundle over $S^2$ with Euler number $-m$. Hence, there is a natural isomorphism between the $\Spin^c$ structures on $X_{c_3,c_2,c_1,c_1'}$ which have torsion restrictions to the boundary, and $\Spin^c(D(-m,1))$. Hence, writing $\frs_k^{\pm}$ for the $\Spin^c$ structures in Lemma~\ref{lem:omega-a-null-homotopy}, we have that
 $h_{\Q;\frs_k^{\pm}}(\Theta_{c_3,c_2}, \Theta_{c_2,c_1}^{\can}, \Theta_{c_1,c_1'}^+)$ is of homogeneous grading
\[
-mk(k+1)+1,
\]
which is always odd. However $\CF^-(\bT^2, c_3,c_1)\iso \bF[U]$ is supported in even gradings, so we conclude that $h_{\Q;\frs}(\Theta_{c_3,c_2}, \Theta_{c_2,c_1}^{\can}, \Theta_{c_1,c_1'}^+)=0$ for all $\frs$ (and in fact, there are no rectangle classes of index $-1$ spanning the given intersection points).
\end{proof}

\section{Hypercubes for changing the almost complex structure}
\label{sec:hypercubes-almost-complex-structure}

We now construct a 3-dimensional hypercube for changing the almost complex structure. We can think of this 3-dimensional hypercube as relating the 2-dimensional ones considered in Lemma~\ref{lem:omega-a-null-homotopy}, for different 1-dimensional families of almost complex structures.

\begin{figure}[ht!]
	\centering
	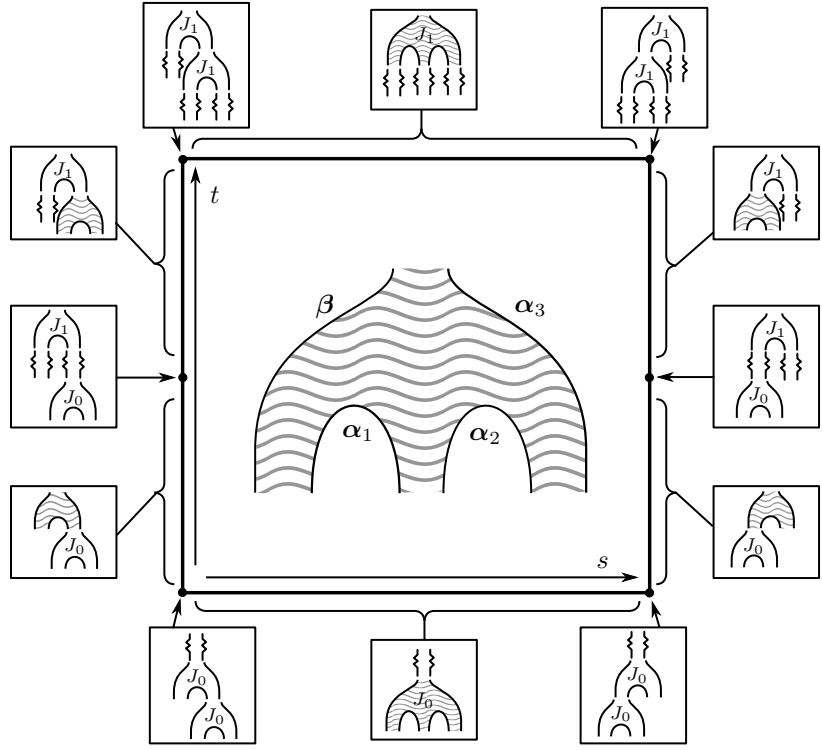
	\caption{A 2-parameter stratified family $(J_{s,t})_{(s,t)\in [0,1]\times [0,1]}$ of almost complex structures on $\Sigma\times \Box$ used to relate two 1-parameter families of almost complex structures. The behavior of $J_{s,t}$ along the codimension 1 and 2 strata is illustrated. Wavy gray lines along the interior of an $n$-gon indicate where the almost complex structure may vary along a given strata.  Wavy boundary lines indicate non-cylindrical almost complex structures on $\Sigma\times [0,1]\times \R$.}\label{fig:111}
\end{figure}

\begin{lem}\label{lem:hypercubes-changing-almost-complex-structure}
Suppose $(\Sigma,\as_3,\as_2,\as_1,\bs,w,z)$ is an $(n,n+m,\infty)$-surgery quadruple for a knot $K\subset Y$, and $J_0=(J_{s,0})_{s\in [0,1]}$ and $J_1=(J_{s,1})_{s\in [0,1]}$ are two stratified families of almost complex structures on $\Sigma\times \Box$, for counting holomorphic rectangles, which are maximally pinched along the boundary of the surgery region. There are maps which make the following diagram a hypercube of chain complexes:
\begin{equation}
\begin{tikzcd}[column sep=2.5cm, labels=description]
\bCF^-_{J_{0}}(\Sigma,\as_1,\bs)
	\arrow[ddddrr,dotted, sloped, pos=.35, swap,
		"{}^{J_0\to J_1}\!h_{\a_1\to \a_2\to \a_3}^{\b}"]
	\arrow[ddd, "\Psi_{J_0\to J_1}",swap]
	\arrow[dr, "{}^{J_0}\!f_{\a_1\to \a_2}^{\b}",sloped,swap]
	\arrow[rrd,sloped, dashed,
		"{}^{J_0}\! h^{\b}_{\a_1\to \a_2\to \a_3}"]&& 
\,
\\
&
\bCF^-_{J_0}(\Sigma,\as_2,\bs)
	\arrow[r, "{}^{J_0}\!f^{\b}_{\a_2\to \a_3}",swap, sloped]
	\arrow[dddr, "{}^{J_0\to J_1}\! f_{\a_2\to \a_3}^{\b}", sloped, dashed] 
&
\buCF^-_{J_0}(\Sigma,\as_3,\bs)
	\arrow[ddd, "\Psi_{J_0\to J_1}"]
 \\
&&
\\
\bCF^-_{J_1}(\Sigma,\as_1,\bs)
	\arrow[drr, dashed,pos=.6, "{}^{J_1}\! h_{\a_1\to \a_2\to \a_3}^{\b}",sloped]
	\arrow[dr, sloped, swap, "{}^{J_1}\!f_{\a_1\to \a_2}^{\b}"]
&& \, \\
&
\bCF^-_{J_1}(\Sigma,\as_2,\bs)
	\arrow[r, sloped,swap, "{}^{J_1}\!f_{\a_2\to \a_3}^{\b}"]
	\arrow[from=uuu, crossing over, "\Psi_{J_0\to J_1}", swap]
	\arrow[from=uuuul,dashed,crossing over, 
	"{}^{J_0\to J_1}\!f_{\a_1\to \a_2}^{\b}", sloped,swap]
& \
\buCF^-_{J_1}(\Sigma,\as_3,\bs)
\end{tikzcd}
\label{eq:hypercube-change-a-c}
\end{equation}
\end{lem}
\begin{proof} We pick a family $\tilde{J}=(J_{s,t})_{(s,t)\in [0,1]\times [0,1]}$ which has the schematic behavior illustrated in Figure~\ref{fig:111}, and which is maximally pinched along $S$.

In equation~\eqref{eq:hypercube-change-a-c}, $\Psi_{J_0\to J_1}$ counts index 0 holomorphic strips for a non-cylindrical almost complex structure (the choice of non-cylindrical almost complex structure is determined by the value of $J_{s,t}$ on the codimension 1 strata of $[0,1]\times [0,1]$).
Similarly, the maps ${}^{J_0\to J_1}\!f_{\a_1\to \a_2}^{\b}$
 and ${}^{J_0\to J_1}\! f_{\a_2\to \a_3}^{\b}$ count index $-1$ triangles which are $J_{s,t}$-holomorphic for $(s,t)$ in a particular  subset of $\d( [0,1]\times [0,1])$, and ${}^{J_0\to J_1}\! h_{\a_1\to \a_2\to \a_3}^{\b}$ counts index $-2$, $J_{s,t}$-holomorphic rectangles, for $(s,t)\in (0,1)\times (0,1)$.

 The length 0, 1 and 2 hypercube relations are easy to verify, so we focus on the length 3 relation. If $\psi\in \pi_2(\Theta_{\a_3,\a_2}^+,\Theta_{\a_2,\a_1}^{\can}, \xs)$ is a Maslov index $-1$ class of rectangles on the quadruple $(\Sigma,\as_3,\as_2,\as_1,\bs)$, then we consider the 1-dimensional moduli space of rectangles
\begin{equation}
\cM_{\tilde{J}}(\psi):=\bigcup_{(s,t)\in (0,1)\times (0,1)} \cM_{J_{s,t}}(\psi)\times \{(s,t)\}. \label{eq:moduli-index-1-triangles}
\end{equation}
 The boundary of the compactification of $\cM_{\tilde{J}}(\psi)$ consists of the following broken curves:
\begin{enumerate}
\item An index $1$ holomorphic strip in one of the cylindrical ends, and an index $-2$, $J_{s,t}$-holomorphic quadrilateral for some $(s,t)\in (0,1)\times (0,1)$.
\item A broken curve consisting of $J_{s,t}$-holomorphic curves for $(s,t)$ in the codimension 1 strata of  $[0,1]\times [0,1]$. Such a curve consists of an index $-1$ curve in the portion with interior wavy lines in Figure~\ref{fig:111}, and index 0 curves in all of the remaining regions.
\end{enumerate}

The count of the ends of the moduli spaces $\cM_{\tilde{J}}(\psi)$ at $(s,t)\in (0,1)\times (0,1)$ give 
\begin{equation}
\left[\d, {}^{J_0\to J_1}\!h_{\a_1\to \a_2\to \a_3}^{\b}\right](\xs). \label{eq:len-3-hypercube-ac-1}
\end{equation}
 Ends of $\cM_{\tilde{J}}(\psi)$ which occur along  $\{1\}\times (0,1)$ contribute the expression
\begin{equation}
\left({}^{J_0\to J_1}\! f_{\a_2\to \a_3}^{\b}\circ {}^{J_0}\! f_{\a_1\to \a_2}^{\b} +{}^{J_1}\! f_{\a_2\to \a_3}^{\b}\circ {}^{J_0\to J_1}\! f_{\a_1\to \a_2}^{\b}\right)(\xs). \label{eq:len-3-hypercube-ac-2}
\end{equation}
Ends which occur along $(0,1)\times \{0,1\}$ contribute
\begin{equation}
\left( \Psi_{J_0\to J_1}\circ {}^{J_0}\!h_{\a_1\to \a_2\to \a_3}^{\b}+{}^{J_1}\! h_{\a_1\to \a_2\to \a_3}^{\b}\circ \Psi_{J_0\to J_1}\right)(\xs). \label{eq:len-3-hypercube-ac-3}
\end{equation}
The sum of equations~\eqref{eq:len-3-hypercube-ac-1}, ~\eqref{eq:len-3-hypercube-ac-2} and~\eqref{eq:len-3-hypercube-ac-3} is exactly the length 3 hypercube relation, so it remains to show that the count of ends appearing along $\{0\}\times (0,1)$ make no total contribution.

The ends appearing along $\{0\}\times (0,1/2)$ contribute
\[
f^{J_0\to J_1}_{\a_3,\a_1,\b}\left(f^{J_0}_{\a_3,\a_2,\a_1}\left(\Theta_{\a_3,\a_2}^+,\Theta_{\a_2,\a_1}^{\can}\right),\xs\right).
\]
Proposition~\ref{prop:stabilize-triangles-general} implies
 \begin{equation}
f_{\a_3,\a_2,\a_1}^{J_0}\left(\Theta_{\a_3,\a_2}^+,\Theta_{\a_2,\a_1}^{\can}\right)=f_{c_3,c_3,c_1}\left(\Theta_{c_3,c_2}^+,\Theta_{c_2,c_1}^{\can}\right)\otimes \Theta_{\a_3',\a_1'}, \label{eq:change-of-ac-hypercube-eq1}
\end{equation}
where $\as_i'=\as_i\setminus \{c_i\}$. By Lemma~\ref{lem:model-computations-m-surgery}, the right side of equation~\eqref{eq:change-of-ac-hypercube-eq1} vanishes. Hence the count of curves along  $\{0\}\times (0,1/2)$ makes no total contribution.

We now consider the ends which occur along $\{0\}\times (1/2,1)$, which contribute the expression
\[
f^{J_1}_{\a_3,\a_1,\b}\left(f^{J_0\to J_1}_{\a_3,\a_2,\a_1}\left(\Theta_{\a_3,\a_2}^+,\Theta_{\a_2,\a_1}^{\can}\right),\Psi_{J_0\to J_1}(\xs)\right).
\]
We claim that 
\begin{equation}
f^{J_0\to J_1}_{\a_3,\a_2,\a_1}\left(\Theta_{\a_3,\a_2}^+,\Theta_{\a_2,\a_1}^{\can}\right)=0.\label{eq:degenerated-ac-a3a2a1}
\end{equation}
By definition, $f^{J_0\to J_1}_{\a_3,\a_2,\a_1}(\Theta_{\a_3,\a_2}^+, \Theta_{\a_2,\a_1}^{\can})$ counts index $-1$, matched triangles over a 1-parameter family of almost complex structures. Since the family $\tilde{J}$ is maximally pinched, the $(\as_3,\as_2,\as_1)$-component of the family $(J_{0,t})_{t\in [0,1]}$ has infinite neck length along $S$, i.e., it is the wedge of a family $(I_t)_{t\in [0,1]}$ on $\Sigma\times \Delta$, and a family $(I'_t)_{t\in [0,1]}$ on $\bT^2\times \Delta$, and by construction, the moduli spaces appearing along $\{0\}\times (0,1)$ are fibered products over the evaluation map at the connected sum point, as well as the parameter $t\in [0,1]$. We claim that the matched moduli spaces counted by $f^{J_0\to J_1}_{\a_3,\a_2,\a_1}(\Theta_{\a_3,\a_2}^+, \Theta_{\a_2,\a_1}^{\can})$ are always empty. To see this, suppose that $(\psi,\psi_0)\in \pi_2(\Theta_{\a'',\a'}^+,\Theta^+_{\a',\a},\zs)\times  \pi_2(\Theta_{c_3,c_2}, \Theta^{\can}_{c_2,c_1},\Theta_{c_3,c_1})$ is a pair of classes, which have the same multiplicity at the connected sum point. The matched index is given by
\begin{equation}
\ind_{\emb}(\psi, \psi_0,M_{\#})=\mu(\psi_0)+\gr(\Theta_{\a'',\a}^+,\zs)\label{eq:index-matched-triangles-a3a2a1}
\end{equation}

 Transversality for the family $I_s$ implies that if $(\psi, \psi_0)$ has a matched representative, then $\mu(\psi_0)\ge -1$. Hence equation~\eqref{eq:index-matched-triangles-a3a2a1} implies that if $(\psi,\psi_0)$ has a representative and $\ind_{\emb}(\psi,\psi_0,M_{\#})=-1$, then
\[
\mu(\psi_0)=-1\quad \text{and} \quad \zs=\Theta_{\a'',\a}^+.
\]
However the diagram $(\bT^2,c_3,c_2,c_1)$ has only triangles of even index, so $\mu(\psi_0)\neq -1$. Hence
 \[
 f_{\a_3,\a_2,\a_1}^{J_0\to J_1}\left(\Theta_{\a_3,\a_2}^+, \Theta_{\a_2,\a_1}^{\can}\right)=0,
 \]
  implying equation~\eqref{eq:degenerated-ac-a3a2a1}, and completing the proof.
\end{proof}

\section{Hypercubes for stabilization}
\label{sec:hypercube-for-stabilization}
We now describe 3-dimensional hypercubes for stabilization and destabilization. Suppose that
\[
\cQ=(\Sigma,\as_3,\as_2,\as_1,\bs,w,z) \quad \text{and} \quad  \cQ_0=(\Sigma_0,\xis,\zetas,\sigmas,\taus,p_0)
\]
are Heegaard quadruples, such that $\cQ$ is an $(n,n+m,\infty)$-surgery quadruple, as in Definition~\ref{def:surgery-quadruple}, and $\cQ_0$ is a multi-stabilizing Heegaard quadruple, as defined in Section~\ref{sec:stabilization}. In this section, we describe a 3-dimensional hypercube of chain complexes relating the 2-dimensional hypercubes for $\cQ$ and $\cQ\# \cQ_{0}$.

 The hypercube for stabilization is constructed as the compression of two hypercubes:
 \begin{enumerate}
 \item A \emph{singular stabilization hypercube}, which involves stabilizing using singularized almost complex structures.
 \item A \emph{desingularization hypercube}, which is a combination of the change of almost complex structure hypercube from Section~\ref{sec:hypercubes-almost-complex-structure}, and the transition maps for desingularizing the almost complex structure from Lemma~\ref{lem:desingularization-maps}.
 \end{enumerate}

\subsection{Singular stabilization hypercubes}
\label{sec:hypercubes-singular-stabilization}

We first construct a stabilization hypercube relating $\cQ$ and $\cQ\# \cQ_0$ which involves singular almost complex structures.

Let $J=(J_s)_{s\in [0,1]}$ be a stratified family of almost complex structures on $\Sigma\times \Box$, for counting rectangles, which is maximally pinched along $\d F$, where $F\subset \Sigma$ denotes the surgery region. Let $I=(I_s)_{s\in [0,1]}$ be a stratified family of maximally pinched almost complex structures on $\Sigma_0\times \Box$, for counting rectangles. We write $J\wedge I=(J_s\wedge I_s)_{s\in [0,1]}$  for the wedge product family. We will always assume that the almost complex structures are split in a neighborhood of the connected sum points.

Similar to section~\ref{sec:totally-degenerated-ac-complexes}, we may define the matched moduli spaces of triangles and rectangles. The matched moduli spaces of triangles may be defined by an easy adaptation of the construction for disks, in ~\eqref{eq:matched-moduli-space-disks}.
The matched moduli space of rectangles is defined similarly, however we also require the rectangles to be matched with respect to the aspect ratio parameter. More concretely, if $S^{\qs}$ and $S_0^{\qs_0}$ are two source curves for rectangles, and $\psi$ and $\psi_0$ are two classes of rectangles on $\cQ$ and $\cQ_0$, then we define the matched moduli space of rectangles to be 
\begin{equation}
\begin{split}
&\cM\cM_{J\wedge I}(S^{\qs}, S_0^{\qs_0}, \psi,\psi_0)\\
:=&\left\{(u,s,u_0,s_0)\in \cM_{J}(S^{\qs},\phi)\times \cM_{I}(S_0^{\qs_0}, \psi_0)\middle \vert \begin{array}{c} (\pi_\Sigma\circ u)(q_i)=p,\\
(\pi_{\Sigma_0}\circ u_0)(q_{0,i})=p_0,\\
(\pi_\Box\circ u)(q_i)=(\pi_\Box\circ u)(q_{i,0}),\\
 \text{ for } i=1,\dots, k, \\
s=s_0
 \end{array}\right\}
 \end{split}
 \label{eq:matched-moduli-space-rectangles}
\end{equation}

We write $\cM\cM_{J\wedge I}(\psi,\psi_0)$ for the disjoint union of the moduli spaces in ~\eqref{eq:matched-moduli-space-rectangles}, over pairs of equivalence classes of marked sources.

Recall that in Section~\ref{sec:index-tuple}, to a multi-stabilizing quadruple $\cQ_0$, we associated a 4-tuple of integers $\ve{m}(\cQ_0)=(m_1(\cQ_0),m_2(\cQ_0),m_3(\cQ_0),m_4(\cQ_0))$, which we called the \emph{index tuple}.

\begin{lem}\label{lem:singular-hypercube-stabilization} Suppose $\cQ=(\Sigma,\as_3,\as_2,\as_1,\bs,w,z)$ is an $(n,n+m,\infty)$-surgery quadruple, with genus 1 surgery region bounded by a curve $S$, and $\cQ_0=(\Sigma_0,\xis,\zetas,\sigmas,\taus,\ps_0)$ is a multi-stabilization quadruple, satisfying
\[
m_4(\cQ_0)=0.
\] Suppose $J=(J_s)_{s\in [0,1]}$ and $I=(I_s)_{s\in [0,1]}$ are stratified families of almost complex structures on $\Sigma\times \Box$ and $\Sigma_0\times \Box$, both for counting  rectangles, and $J$ is maximally pinched along $S$. The following diagram is a hypercube:
\[
\begin{tikzcd}[row sep=.6cm, column sep=1.7cm]
\bCF^-_{J}(\as_1,\bs)
	\arrow[ddd, "F_1^{\sigma,\tau}",description]
	\arrow[rd, "{}^J\! f_{\a_1\to \a_2}^{\b}" description,sloped]
	\arrow[drr, dashed, "{}^{J}\!h_{\a_1\to \a_2\to \a_3}^{\b}" description,sloped]&& 
\,
\\
&
\bCF^-_{J}(\as_2,\bs)
	\arrow[r, "{}^J\!f_{\a_2\to \a_3}^{\b}" description,sloped, swap]
& \buCF^-_{J}(\as_3,\bs)
	\arrow[ddd, "F_1^{\xi,\tau}" description]
 \\
&&&
\\
\bCF^-_{J\wedge I}(\as_1\cup \sigmas,\bs\cup \taus)
	\arrow[drr, dashed,near end, sloped, 
		"{}^{J\wedge I}\! h_{\a_1\cup \sigma\to\a_2\cup \zeta\to \a_3\cup \xi}^{\b\cup \tau}" description]
	\arrow[dr, "{}^{J\wedge I}\!f_{\a_1\cup \sigma\to \a_2\cup \zeta}^{\b\cup \tau}",swap, sloped]
&& \,\\
&
\bCF^-_{J\wedge I}(\as_2\cup \zetas,\bs\cup \taus)
	\arrow[from=uuu, "F_1^{\zeta,\tau}" description, swap, crossing over]
	\arrow[r, swap, 
		"{}^{J\wedge I}\!f_{\a_2\cup \zeta\to \a_3\cup \xi}^{\b\cup \tau}"]
&
\buCF^-_{J\wedge I}(\as_3\cup \xis,\bs\cup \taus)
\end{tikzcd}.
\]
\end{lem}

\begin{proof}
We verify the hypercube relations. The length 0 hypercube relations are immediate. We now consider the length 1 hypercube relations. For the horizontal length 1 maps, the hypercube relations follow by counting the ends of moduli spaces of triangles of index 1. For the vertical length 1 maps, the hypercube relations follow from Lemma~\ref{lem:singular-1-handle-chain-maps}.

We now consider the length 2 relations. Lemma~\ref{lem:omega-a-null-homotopy} imples that ${}^{J}\!h_{\a_1\to \a_2\to \a_3}^{\b}$ is a null-homotopy of ${}^J\!f_{\a_2\to \a_3}^{\b}\circ {}^J \! f_{\a_1\to \a_2}^{\b}$. The same argument, applied to matched triangles and rectangles, implies that ${}^{J\wedge I}\!h_{\a_1\cup \sigma\to \a_2\cup \zeta\to \a_3\cup \xi}^{\b\cup \tau}$ is a null-homotopy of ${}^{J\wedge I}\!f_{\a_2\cup \zeta\to \a_3\cup \xi}^{\b\cup \tau} \circ {}^{J\wedge I}\!f_{\a_1\cup \sigma\to \a_2\cup \zeta}^{\b\cup \tau}$. Additionally, we need to show that
\begin{equation}
\begin{split}
F_1^{\zeta,\tau}\circ {}^J\!f_{\a_1\to \a_2}^{\b}&={}^{J\wedge I}\!f_{\a_1\cup \sigma\to \a_2\cup \zeta}^{\b\cup \tau} \circ F_1^{\sigma,\tau}, \quad \text{and}\\
F_1^{\sigma,\tau}\circ {}^J\!f_{\a_2\to \a_3}^{\b}&={}^{J\wedge I}\!f_{\a_2\cup \zeta\to \a_3\cup \xi}^{\b\cup \tau} \circ F_1^{\zeta, \tau}.
\end{split}
\label{eq:stab-hypercube-triangle-1-handle-commute}
\end{equation}
By hypothesis, $m_4(\cQ_0)=0$. By Lemma~\ref{lem:index-associativity}, we have
\begin{equation}
m_4(\cQ)=m_3(\cT_{\zeta,\sigma,\tau})+m_3(\cT_{\xi,\zeta,\tau}). \label{eq:index-stabilization-associativity}
\end{equation}
Using Lemma~\ref{lem:monotonicity-index-tuple} and~\eqref{eq:index-stabilization-associativity}, we obtain
\[
0=m_3(\cT_{\zeta,\sigma,\tau})=m_3(\cT_{\xi,\zeta,\tau}).
\]
Proposition~\ref{prop:stabilize-triangles-general} now gives~\eqref{eq:stab-hypercube-triangle-1-handle-commute}.

Lastly, we consider the length 3 hypercube relation. The length 3 hypercube relation is equivalent to
\[
F_1^{\xi,\tau}\circ {}^{J}\!h_{\a_1\to \a_2\to \a_3}^{\b}={}^{J\wedge I}\!h_{\a_1\cup \sigma\to\a_2\cup \zeta\to  \a_3\cup \xi}^{\b\cup \tau}\circ F_1^{\xi,\tau},
\]
which follows from Proposition~\ref{prop:stabilize-quadrilateral}.
\end{proof}

Dually, we obtain the following:

\begin{lem}\label{lem:singular-hypercube-destabilization} Suppose that $\cQ=(\Sigma,\as_3,\as_2,\as_1,\bs,w,z)$ is an $(n,n+m,\infty)$-surgery quadruple, with a genus 1 surgery region bounded by the curve $S$. Suppose that $\cQ_0=(\Sigma_0,\xis,\zetas,\sigmas,\taus,\ps_0)$
 is a multi-stabilizing quadruple with 
 \[
 m_3(\cQ_0)=0.
 \]
  Suppose $J=(J_s)_{s\in [0,1]}$ is a family of almost complex structures on $\Sigma\times \Box$ which is maximally pinched along $S$, and $I=(I_s)_{s\in [0,1]}$ is a family of almost complex structures on $\Sigma_0\times \Box$, both for counting holomorphic rectangles.  The following diagram is a hypercube of chain complexes
\[
\begin{tikzcd}[row sep=.6cm, column sep=1.7cm]
\bCF^-_{J\wedge I}(\as_1\cup \sigmas,\bs\cup \taus)
\arrow[drr, dashed,sloped, "{}^{J\wedge I} \! h_{\a_1\cup \sigma\to \a_2\cup \zeta\to \a_3\cup \xi}^{\b\cup \tau}"]
\arrow[dr, "{}^{J\wedge I}\! f_{\a_1\cup \sigma\to \a_2\cup \zeta} ^{\b\cup \tau}",swap, sloped]
\arrow[ddd, "F_3^{\sigma,\tau}"]
&&
\,
\\
& \bCF^-_{J\wedge I}(\as_2\cup \zetas,\bs\cup \taus)
\arrow[r, "{}^{J\wedge I} \!f_{\a_2\cup \zeta\to \a_3\cup \xi}^{\b\cup \tau}", near end, swap]
& \buCF^-_{J\wedge I}(\as_3\cup \xis,\bs\cup \taus) 
\arrow[ddd, "F_3^{\xi,\tau}"]
\\
&&&
\\
\bCF^-_{J}(\as_1,\bs)
\arrow[drr,sloped, dashed,near end, "{}^{J}\!h_{\a_1\to \a_2\to \a_3}^{\b}"]
\arrow[dr, "{}^{J}\!f_{\a_1\to \a_2}^{\b}",swap, sloped]
&& \,\\
& \bCF^-_{J}(\as_2,\bs)
\arrow[r, "{}^{J}\!f_{\a_2\to \a_3}^{\b}", swap]
\arrow[from=uuu, "F_3^{\zeta,\tau}", crossing over, near start]
&\buCF^-_{J}(\as_3,\bs)
\end{tikzcd}
\]
\end{lem}

\subsection{Desingularization hypercubes}

In Section~\ref{sec:hypercubes-singular-stabilization} we described a hypercube for stabilization which involved matched moduli spaces (i.e. fibered products along evaluation maps). We now relate construct a hypercube which relates the singular almost complex structures to ordinary ones:

\begin{lem}\label{lem:hypercube-desingunarlize-stabilization} Suppose that $\cQ=(\Sigma,\as_3,\as_2,\as_1,\bs,\ws)$ is an $(n,n+m,\infty)$ surgery quadruple, with a genus 1 surgery region bounded by a curve $S$, and $\cQ_0=(\Sigma_0,\xis,\zetas,\sigmas,\taus,p)$ is a multi-stabilizing Heegaard quadruple. Suppose that $J=(J_s)_{s\in [0,1]}$ and $I=(I_s)_{s\in [0,1]}$ are stratified families of almost complex structures on $\Sigma\times \Box$ and $\Sigma_0\times \Box$, and $J$ is maximally pinched along $S$. Suppose further that $K=(K_s)_{s \in [0,1]}$ is a family of almost complex structure on $\Sigma\# \Sigma_0\times \Box$ which is maximally pinched along $S$. Then there is a stratified family of almost complex structures $L=(L_{s,t})_{(s,t)\in [0,1]\times [0,1]}$ on $\Sigma\# \Sigma_0\times \Box$, whose associated moduli spaces give the following hypercube of chain complexes:
\[
\begin{tikzcd}[row sep=1 cm, column sep=1.5cm]
\bCF^-_{J\wedge I}(\as_1\cup \sigmas,\bs\cup \taus)
	\arrow[ddd, "\Psi_{J\wedge I\to K}",swap]
	\arrow[dr, "{}^{J\wedge I}\!
	 	f_{\a_1\cup \sigma\to \a_2\cup \zeta}^{\b\cup \xi}",sloped, swap]
	\arrow[drr, dashed, "{}^{J\wedge I}\!
		 h_{\a_1\cup \sigma\to\a_2\cup \zeta\to  \a_3\cup\xi}^{\b\cup \tau}",sloped]
	\arrow[ddddrr,dotted, "{}^{J\wedge I\to K}\!
		h_{\a_1\cup \sigma\to\a_2\cup \zeta\to \a_3\cup \xi}^{\b\cup \tau}\hspace{.7cm}",
		 sloped, near end,swap]
&& 
\,
\\
&\bCF^-_{J\wedge I}(\as_2\cup \zetas,\bs\cup \taus)
 	\arrow[r, "{}^{J\wedge I}\!f_{\a_2\cup \zeta\to 
 		\a_3\cup \xi}^{\b\cup \tau}",near end, swap]
	\arrow[dddr,dashed, sloped, "{}^{J\wedge I\to K}
		\!f_{\a_2\cup \zeta\to \a_3\cup \xi}^{\b\cup \tau}"] 
&\buCF^-_{I\wedge J}(\as_3\cup \xis,\bs\cup \taus)
	\arrow[ddd, "\Psi_{J\wedge I\to K}"]
 \\
&&&
\\
\bCF^-_{K}(\as_1\cup \sigmas,\bs\cup \taus)
	\arrow[drr, dashed,sloped,near end, 
		"{}^K\!h_{\a_1\cup \sigma\to\a_2\cup \zeta\to  \a_3\cup \xi}^{\b\cup \tau}"]
	\arrow[dr, "{}^{K}\! f_{\a_1\cup \sigma\to \a_2\cup \zeta}^{\b\cup \tau}"
		,swap, sloped]
&& \,\\
& 
\bCF^-_{K}(\as_2\cup \zetas,\bs\cup \taus)
	\arrow[from=uuu,crossing over, "\Psi_{J\wedge I\to K}", swap]
	\arrow[r, "{}^{K}\!f_{\a_2\cup \zeta\to \a_3\cup \xi}^{\b\cup \tau}",
		sloped, swap,near end]
	\arrow[from=uuuul, crossing over, dashed, 
		"{}^{J\wedge I\to K}\!f_{\a_1\cup \sigma\to \a_2\cup \zeta}^{\b\cup \tau}",sloped]
&
\buCF^-_{K}(\as_3\cup \xis,\bs\cup \taus)
\end{tikzcd}
\]
\end{lem}
\begin{proof}
The proof follows from the same reasoning as the proof of Lemma~\ref{lem:hypercubes-changing-almost-complex-structure}, which covers more general changes of the almost complex structure.
\end{proof} 

By stacking and compressing the hypercubes in Lemmas~\ref{lem:singular-hypercube-stabilization} and ~\ref{lem:hypercube-desingunarlize-stabilization}, we obtain the general hypercube for stabilization. Dualizing, as in Lemma~\ref{lem:singular-hypercube-destabilization}, yields a hypercube for destabilization.

\section{The first 2-handle hypercube} 

\label{sec:peripheral-hypercubes}

We now construct the first 2-handle hypercube, $\cC_{\twoh}^1$.  For all of the complexes appearing in this hypercube, we use the underlying Heegaard surface  $\Sigma\# \bar{\Sigma}$. We make the following definitions:
\begin{itemize}
\item We write $D(F)$ for the genus 2 subset of $\Sigma\# \bar{\Sigma}$ consisting of the surgery region on $\Sigma$, and its image on $\bar{\Sigma}$.
\item We write $\Ds\subset \Sigma\# \bar{\Sigma}$ for a collection of curves which are obtained by doubling a basis of arcs $d_1,\dots, d_{2g}$ for $H_1(\Sigma\setminus N(w),I;\Z)$. Furthermore, we assume that $\Ds$ is adapted to $\as_3$, in the sense of Definition~\ref{def:adapted-to-alpha}. 
\item We write $\Delta_1$ and $\Delta_2$ for the two curves of $\Ds$ which are in $D(F)$. And we write $\hat{\Ds}$ for $\Ds\setminus \{\Delta_1,\Delta_2\}$.
\end{itemize}

We consider the pentuple
\[
\cP=\left(\Sigma\# \bar{\Sigma}, \as_3\cup \bar{\bs}''', \as_2\cup \bar{\bs}'', \as_1\cup \bar{\bs}', \bs\cup \bar{\bs}, \Ds,w\right),
\]
where $\bar{\bs},$ $\bar{\bs}'$, $\bar{\bs}''$, $\bar{\bs}'''$ are small Hamiltonian translates.

\begin{lem}\label{lem:admissibility-C-2-h-1}
If  $\bs$ are sufficiently wound on $\Sigma$ (with the winding also performed on $\bar{\bs}$ and its small translates), and $\Delta_1$ and $\Delta_2$ are sufficiently wound on $D(F)$, then $\cP$ is weakly admissible.
\end{lem}
\begin{proof}
Let $\Pi_{\cP}^{\Q}$ denote the set of rational periodic domains on $\cP$. Let $\Pi_{\cP,L}^{\Q}$ denote the set of periodic domains with boundary equal to a rational linear combination of the curves in 
\[
\as\cup \{c_1,c_2,c_3\}\cup \bs \cup \bar{\bs}\cup \Ds.
\]
(That is, we remove all but one curve in any collection of pairwise isotopic curves). There is a canonical map
\[
L\colon \Pi^{\Q}_{\cP}\to \Pi_{\cP,L}^{\Q}.
\] 
The map $L$ is clearly surjective. See Figure~\ref{fig:185} for a schematic of $L$.

\begin{figure}[h]
	\centering
\begingroup%
  \makeatletter%
  \providecommand\color[2][]{%
    \errmessage{(Inkscape) Color is used for the text in Inkscape, but the package 'color.sty' is not loaded}%
    \renewcommand\color[2][]{}%
  }%
  \providecommand\transparent[1]{%
    \errmessage{(Inkscape) Transparency is used (non-zero) for the text in Inkscape, but the package 'transparent.sty' is not loaded}%
    \renewcommand\transparent[1]{}%
  }%
  \providecommand\rotatebox[2]{#2}%
  \newcommand*\fsize{\dimexpr\f@size pt\relax}%
  \newcommand*\lineheight[1]{\fontsize{\fsize}{#1\fsize}\selectfont}%
  \ifx\svgwidth\undefined%
    \setlength{\unitlength}{183.24796983bp}%
    \ifx\svgscale\undefined%
      \relax%
    \else%
      \setlength{\unitlength}{\unitlength * \real{\svgscale}}%
    \fi%
  \else%
    \setlength{\unitlength}{\svgwidth}%
  \fi%
  \global\let\svgwidth\undefined%
  \global\let\svgscale\undefined%
  \makeatother%
  \begin{picture}(1,0.38404888)%
    \lineheight{1}%
    \setlength\tabcolsep{0pt}%
    \put(0,0){\includegraphics[width=\unitlength,page=1]{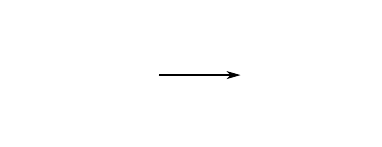}}%
    \put(-0.00165205,0.1480887){\color[rgb]{0,0,0}\makebox(0,0)[lt]{\lineheight{1.25}\smash{\begin{tabular}[t]{l}$n$\end{tabular}}}}%
    \put(0.27868205,0.1480887){\color[rgb]{0,0,0}\makebox(0,0)[lt]{\lineheight{1.25}\smash{\begin{tabular}[t]{l}$m$\end{tabular}}}}%
    \put(0,0){\includegraphics[width=\unitlength,page=2]{fig185.pdf}}%
    \put(0.79389168,0.1480887){\color[rgb]{0,0,0}\makebox(0,0)[lt]{\lineheight{1.25}\smash{\begin{tabular}[t]{l}$n$\end{tabular}}}}%
    \put(0.94201101,0.1480887){\color[rgb]{0,0,0}\makebox(0,0)[lt]{\lineheight{1.25}\smash{\begin{tabular}[t]{l}$m$\end{tabular}}}}%
    \put(0,0){\includegraphics[width=\unitlength,page=3]{fig185.pdf}}%
    \put(0.51471705,0.21217221){\color[rgb]{0,0,0}\makebox(0,0)[t]{\lineheight{1.25}\smash{\begin{tabular}[t]{c}$K$\end{tabular}}}}%
    \put(0,0){\includegraphics[width=\unitlength,page=4]{fig185.pdf}}%
  \end{picture}%
\endgroup%

	\caption{A schematic of the homomorphism $L$. The isotopic curves denote $\bar{\bs}$, $\bar{\bs}'$, $\bar{\bs}''$ and $\bar{\bs}'''$. The numbers $n$ and $m$ denote the multiplicities of a domain $D$. On the right is the domain $L(D)$.}\label{fig:185}
\end{figure}

We will show that every non-zero element of  $\Pi^{\Q}_{\cP,L}$ has both positive and negative multiplicities. This will imply that any nonnegative periodic domains in $\Pi^{\Q}_{\cP}$ must be in the kernel of $L$. By picking the small Hamiltonian translations which appear in the construction of $\cP$ appropriately, one can ensure that any elements in the kernel of $L$ have both positive and negative multiplicities, which will imply the desired result.

By applying the same argument as in Lemma~\ref{lem:admissibility-2-tuple}, for sufficiently wound $\bs$, $\bar{\bs}$, $\Delta_1$ and $\Delta_2$ it is sufficient to show that there are no non-zero periodic domains in $\Pi_{\cP,L}^{\Q}$ with only nonnegative coefficients, whose boundary is a linear combination of only the curves in $\as\cup \{c_1,c_2,c_3\}\cup \hat{\Ds}$. Write $\Pi_{\a\cup \{c_1,c_2,c_3\}\cup \hat{\Dt}}^\Q$ for the set of such periodic domains.

Note that $\{c_1,c_2,c_3\}$ are contained on a genus $1$ summand of $\Sigma\# \bar{\Sigma}$, and $\as\cup \hat{\Ds}$ are contained on a disjoint, genus $2g-2$ summand. Hence the subspace $\Span([c_1],[c_2],[c_3])\in H_1(\Sigma\# \bar{\Sigma};\Q)$ has trivial intersection with each of $\Span([\alpha]:\a\in  \as)$ and $\Span ([\Delta]:\Delta\in \hat{\Ds})$. Furthermore, the subspaces $\Span([\a]:\a\in \as)$ and $\Span([\Delta]:\Delta\in \hat{\Ds})$ have trivial intersection in $H_1(\Sigma\# \bar{\Sigma};\Q)$ by the definition of a doubled diagram. Since $\Span([\a]:\a\in \as)$ and $\Span([\Delta]:\Delta\in \hat{\Ds})$ each have dimension $|\as|$ and $|\hat{\Ds}|$, it follows that an element of $\Pi^{\Q}_{\a\cup \{c_1,c_2,c_3\}\cup \hat{\Dt}}$ must have boundary equal to a linear combination of the curves $c_1$, $c_2$ and $c_3$, so $\Pi^{\Q}_{\a\cup \{c_1,c_2,c_3\}\cup \hat{\Dt}}$ is naturally isomorphic to $\Pi^{\Q}_{\{c_1,c_2,c_3\}}$. The main claim now follows from the fact that the diagram $(\bT^2,c_3,c_2,c_1,w)$ is weakly admissible, as can be seen by explicit examination of the diagram in Figure~\ref{fig:63}.
\end{proof}

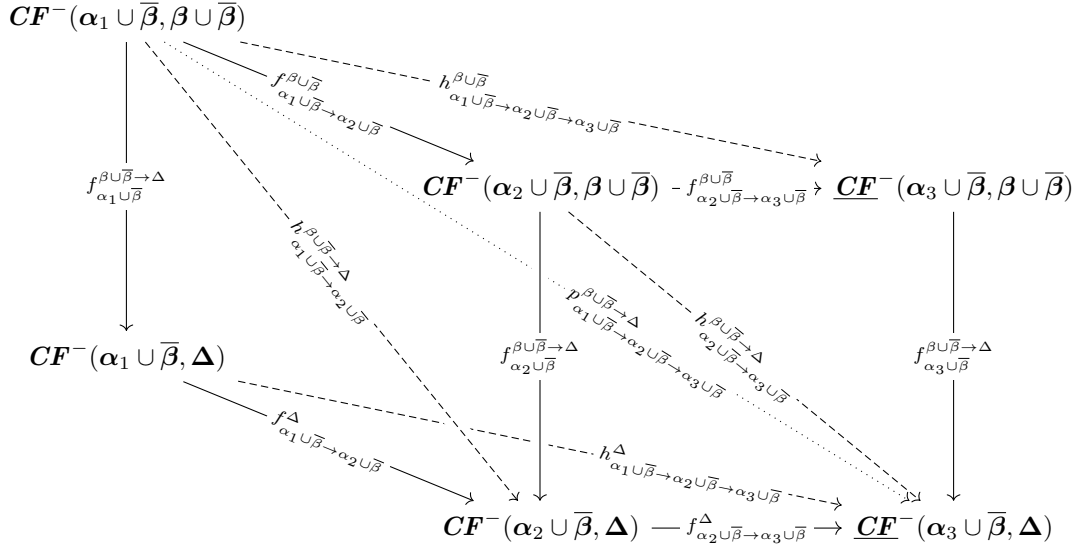
\begin{figure}[h]
\centering
\[
\begin{tikzcd}[column sep=2cm,row sep=1.6cm,labels=description]
\bCF^-( \as_1\cup \bar{\bs}, \bs\cup \bar{\bs})
	\arrow[dr, "f_{\a_1\cup \bar{\b}\to \a_2\cup \bar{\b}}^{ \b\cup \bar{\b}}"
		, sloped]
	\arrow[dd, "f_{\a_1\cup \bar{\b}}^{\b\cup \bar{\b}\to \Delta}",swap ]
	\arrow[drr,dashed, sloped,
		"h_{\a_1\cup \bar{\b} \to \a_2\cup \bar{\b}\to \a_3\cup \bar{\b}}^{\b\cup \bar{\b}}"] 
	\arrow[dddrr,dotted, swap, sloped,
		"p_{\a_1\cup \bar{\b} \to \a_2\cup \bar{\b}\to \a_3\cup \bar{\b}}^{\b\cup \bar{\b}\to \Delta}", pos=.65]
		 &&
 \\
& 
\bCF^-(\as_2\cup \bar{\bs}, \bs\cup \bar{\bs})
	\arrow[r, "f_{\a_2\cup \bar{\b}\to \a_3\cup \bar{\b}}^{\b\cup \bar{\b}}",swap]
	\arrow[ddr,dashed, sloped,
		"h_{ \a_2\cup \bar{\b}\to \a_3\cup \bar{\b}}^{\b\cup \bar{\b}\to \Delta}"]&
\buCF^-( \as_3\cup \bar{\bs}, \bs\cup \bar{\bs}) 
  	\arrow[dd, "f_{\a_3\cup \bar{\b}}^{\b\cup \bar{\b}\to  \Delta}"]
 \\
\bCF^-(\as_1\cup\bar{\bs}, \Ds)
 	\arrow[dr,swap,sloped, "f_{\a_1\cup \bar{\b}\to \a_2\cup \bar{\b}}^{\Delta}"]
	\arrow[drr,dashed,sloped, near end, 
		"h_{ \a_1\cup \bar{\b}\to  \a_2\cup \bar{\b}\to \a_3\cup \bar{\b}}^{ \Delta}"]
 \\
 & 
\bCF^-(\as_2\cup \bar{\bs}, \Ds)
	\arrow[r,swap, "f_{ \a_2\cup \bar{\b}\to \a_3\cup \bar{\b}}^{\Delta}"] 
	\arrow[from=uuul, dashed,crossing over,sloped,
	 "h_{\a_1\cup \bar{\b}\to \a_2\cup \bar{\b}}^{\b\cup \bar{\b}\to \Delta}"]
	\arrow[from=uu, crossing over,
		 "f_{\a_2\cup \bar{\b}}^{\b\cup \bar{\b}\to \Delta}"]
 & 
\buCF^-(\as_3\cup \bar{\bs}, \Ds)
 \end{tikzcd}
\]
\caption{The hypercube $\cC_{\twoh}^1$. The Heegaard surface is $\Sigma\# \bar{\Sigma}$, for all Floer complexes. We have dropped the primes on the various copies of $\bar{\bs}$, to simplify the notation.}
\label{fig:hypercube-C-2}
\end{figure}

Using the pentuple $\cP$, we build the hypercube $\cC_{\twoh}^1$, as shown in Figure~\ref{fig:hypercube-C-2}.

\begin{prop}
 The diagram $\cC_{\twoh}^1$ is a hypercube of chain complexes.
\end{prop}
\begin{proof} The hypercube $\cC_{\twoh}^1$ may equivalently be described as the pairing of the following two hypercubes of attaching curves:
\begin{equation}
\cL_{\a}:=\begin{tikzcd}[row sep=2 cm, column sep=2cm]\as_1\cup \bar{\bs}'\arrow[d, "\Theta_{\a_2\cup \bar{\b}'', \a_1\cup \bar{\b}'}^{\can}"]\\
\as_2\cup \bar{\bs}''\arrow[r, "\Theta_{\a_3\cup \bar{\b}''',\a_2\cup \bar{\b}''}^+"]& \as_3\cup \bar{\bs}'''
\end{tikzcd}\qquad 
\cL_{\b}:=
 \begin{tikzcd}[column sep=2cm]
  \bs\cup \bar{\bs}\arrow[r, "\Theta_{\b\cup \bar{\b},\Dt}^+"]&\Ds.
  \end{tikzcd}
  \label{eq:C-2h-1-as-pairing}
\end{equation}
Using Lemma~\ref{lem:hypercube-diagram-to-complex}, it remains only to show the hypercube relations for $\cL_{\a}$ and $\cL_{\b}$. The only non-trivial relation is the length 2 relation for $\cL_{\a}$, which follows from Lemma~\ref{lem:omega-a-null-homotopy}, completing the proof.
\end{proof}

\section{The central hypercube}
\label{sec:central-hypercube}

In this section, we construct the central hypercube. We construct three hypercubes, $\cC_{\cen}^{(1)}$, $\cC_{\cen}^{(2)}$ and $\cC_{\cen}^{(3)}$, which we stack and compress to obtain a hypercube $\hat{\cC}_{\cen}$. In Section~\ref{sec:constructing-C-cen}, we describe how to modify $\hat{\cC}_{\cen}$ to obtain the hypercube $\cC_{\cen}$.

\subsection{Almost complex structures}

\begin{figure}[ht!]
	\centering
\begingroup%
  \makeatletter%
  \providecommand\color[2][]{%
    \errmessage{(Inkscape) Color is used for the text in Inkscape, but the package 'color.sty' is not loaded}%
    \renewcommand\color[2][]{}%
  }%
  \providecommand\transparent[1]{%
    \errmessage{(Inkscape) Transparency is used (non-zero) for the text in Inkscape, but the package 'transparent.sty' is not loaded}%
    \renewcommand\transparent[1]{}%
  }%
  \providecommand\rotatebox[2]{#2}%
  \newcommand*\fsize{\dimexpr\f@size pt\relax}%
  \newcommand*\lineheight[1]{\fontsize{\fsize}{#1\fsize}\selectfont}%
  \ifx\svgwidth\undefined%
    \setlength{\unitlength}{103.89393292bp}%
    \ifx\svgscale\undefined%
      \relax%
    \else%
      \setlength{\unitlength}{\unitlength * \real{\svgscale}}%
    \fi%
  \else%
    \setlength{\unitlength}{\svgwidth}%
  \fi%
  \global\let\svgwidth\undefined%
  \global\let\svgscale\undefined%
  \makeatother%
  \begin{picture}(1,0.99449228)%
    \lineheight{1}%
    \setlength\tabcolsep{0pt}%
    \put(0,0){\includegraphics[width=\unitlength,page=1]{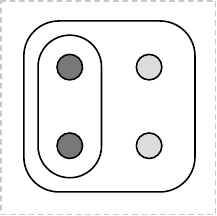}}%
    \put(0.48337353,0.54551636){\color[rgb]{0,0,0}\makebox(0,0)[lt]{\lineheight{1.25}\smash{\begin{tabular}[t]{l}$S_1$\end{tabular}}}}%
    \put(0.59200925,0.13964886){\color[rgb]{0,0,0}\makebox(0,0)[lt]{\lineheight{1.25}\smash{\begin{tabular}[t]{l}$S_2$\end{tabular}}}}%
  \end{picture}%
\endgroup%

	\caption{The closed curves $S_1$ and $S_2$, along which we assume all almost complex structures are maximally pinched (Definition~\ref{def:maximally-pinched}).}\label{fig:129}
\end{figure}

Let $S_1$ and $S_2$ be the two closed curves on $\Sigma\# \bar{\Sigma}$ shown in Figure~\ref{fig:129}. The curve $S_1$ bounds the surgery region of $\Sigma$, and $S_2$ bounds this region as well as its image on $\bar{\Sigma}$. Throughout this section, we assume that all almost complex structures are maximally pinched along $S_1$ and $S_2$ (Definition~\ref{def:maximally-pinched}).

\subsection{Standard curves in the surgery region}

\begin{figure}[ht!]
	\centering
\begingroup%
  \makeatletter%
  \providecommand\color[2][]{%
    \errmessage{(Inkscape) Color is used for the text in Inkscape, but the package 'color.sty' is not loaded}%
    \renewcommand\color[2][]{}%
  }%
  \providecommand\transparent[1]{%
    \errmessage{(Inkscape) Transparency is used (non-zero) for the text in Inkscape, but the package 'transparent.sty' is not loaded}%
    \renewcommand\transparent[1]{}%
  }%
  \providecommand\rotatebox[2]{#2}%
  \newcommand*\fsize{\dimexpr\f@size pt\relax}%
  \newcommand*\lineheight[1]{\fontsize{\fsize}{#1\fsize}\selectfont}%
  \ifx\svgwidth\undefined%
    \setlength{\unitlength}{162.7560454bp}%
    \ifx\svgscale\undefined%
      \relax%
    \else%
      \setlength{\unitlength}{\unitlength * \real{\svgscale}}%
    \fi%
  \else%
    \setlength{\unitlength}{\svgwidth}%
  \fi%
  \global\let\svgwidth\undefined%
  \global\let\svgscale\undefined%
  \makeatother%
  \begin{picture}(1,0.72684256)%
    \lineheight{1}%
    \setlength\tabcolsep{0pt}%
    \put(0,0){\includegraphics[width=\unitlength,page=1]{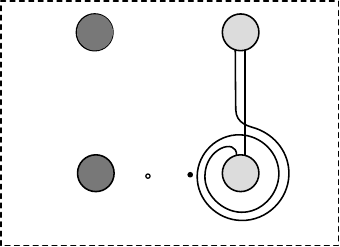}}%
    \put(0.2516975,0.50486162){\color[rgb]{0,0,0}\makebox(0,0)[rt]{\lineheight{1.25}\smash{\begin{tabular}[t]{r}$c_1$\end{tabular}}}}%
    \put(0,0){\includegraphics[width=\unitlength,page=2]{fig186.pdf}}%
    \put(0.3690086,0.37516845){\color[rgb]{0,0,0}\makebox(0,0)[lt]{\lineheight{1.25}\smash{\begin{tabular}[t]{l}$c_3$\end{tabular}}}}%
    \put(0.62788823,0.3755402){\color[rgb]{0,0,0}\makebox(0,0)[rt]{\lineheight{1.25}\smash{\begin{tabular}[t]{r}$\bar{c}_3$\end{tabular}}}}%
    \put(0.31656142,0.50486162){\color[rgb]{0,0,0}\makebox(0,0)[lt]{\lineheight{1.25}\smash{\begin{tabular}[t]{l}$c_2$\end{tabular}}}}%
    \put(0,0){\includegraphics[width=\unitlength,page=3]{fig186.pdf}}%
    \put(0.68185107,0.50486162){\color[rgb]{0,0,0}\makebox(0,0)[rt]{\lineheight{1.25}\smash{\begin{tabular}[t]{r}$\bar{c}_2$\end{tabular}}}}%
    \put(0.73869491,0.50486162){\color[rgb]{0,0,0}\makebox(0,0)[lt]{\lineheight{1.25}\smash{\begin{tabular}[t]{l}$\bar{c}_1$\end{tabular}}}}%
  \end{picture}%
\endgroup%

	\caption{The curves $c_1$, $c_2$, $c_3$, $\bar{c}_1$, $\bar{c}_2$ and $\bar{c}_3$ in the genus 2 region $D(F)$. (Notice the slight asymmetry between the left and right hand sides). The basepoint $w$ is the solid dot on the left, and $z$ is the open dot on the right.}\label{fig:186}
\end{figure}

We will encounter the curves $c_1$, $c_2$ and $c_3$ on $\Sigma$, as well as their images on $\bar{\Sigma}$. We label the curves on $\bar{\Sigma}$ by $\bar{c}_1$, $\bar{c}_2$ and $\bar{c}_3$, as shown in Figure~\ref{fig:186}. Note that to achieve admissibility, we do not pick $\bar{c}_i$ to be an exact mirror of $c_i$, but instead to have the configuration shown in Figure~\ref{fig:186}.

Let $\cP_{D(F)}$ denote the set of rational 2-chains $D$ on $D(F)$ which have boundary equal to a linear combination of the $c_i$ and $\bar{c}_i$ curves, and which have $n_w(D)=0$. The following lemma is helpful to achieve admissibility on the multi-diagrams appearing in the central hypercube:

\begin{lem}\label{lem:admissible-key-computation}
 All non-zero domains $\cP_{D(F)}$ have both positive and negative coefficients.
\end{lem}
\begin{proof} The proof is by direct inspection. The set of triply periodic domains in $\cP_{D(F)}$ is isomorphic to $\Z\oplus \Z$. The first $\Z$ summand is supported on the left side of $D(F)$, while the second $\Z$ summand is supported on the right. Furthermore, both generators have both positive and negative coefficients, so any linear combination does as well.
\end{proof}

\subsection{The hypercube \texorpdfstring{$\cC_{\cen}^{(1)}$}{C-cen-(1)}}
\begin{figure}[ht!]
	\centering
	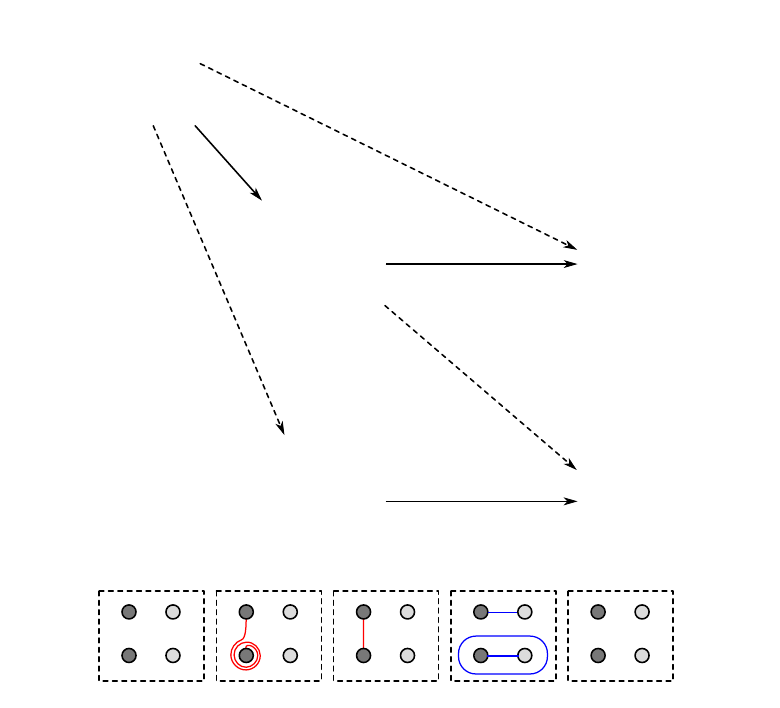
	\caption{The auxiliary hypercube $\cC_{\aux}^{(1)}$. A $\twind$ along a curve indicates that it is wound.}\label{fig:175}
\end{figure}

 We begin by constructing an auxiliary hypercube, $\cC_{\aux}^{(1)}$, shown in Figure~\ref{fig:175}, which is built from a pentuple $(\Sigma\# \bar{\Sigma},\gs'',\gs',\gs,\ds,\ds')$. In the notation of Figure~\ref{fig:large-hyperbox}, the attaching curves of $\cP^{(1)}$ are as follows:
 \begin{itemize}
 \item $\gs=\as_1\cup \bar{\bs}$, $\gs'=\as_2\cup \bar{\bs}$, and $\gs''=\as_3\cup \bar{\bs}$.
 \item $\ds$ coincide with the curves $\Ds$, described in Section~\ref{sec:peripheral-hypercubes}.
 \item  $\ds'$ coincides with small translates of $\hat{\Ds}$ outside of the special genus 2 region $D(F)$, and is equal to small translates of $c_3\cup \bar{c}_3$ in the special genus two region.
 \end{itemize}
 
Recall that we write $\Delta_1$ and $\Delta_2$ for the two curves of $\Ds$ appearing in $D(F)$. We assume that the curves, $c_1$, $c_2$, $c_3$ and $\bar{c}_3$ are in standard position in $D(F)$, i.e. they appear with the configuration of Figure~\ref{fig:70}. The curves $\Delta_1$ and $\Delta_2$ will be wound in $D(F)$, so in general will not be in a simple configuration. Note that we will usually omit the winding of $\Delta_1$ and $\Delta_2$ from the figures.

The maps labeled in Figure~\ref{fig:175} count all triangles, quadrilaterals and pentagons which have special inputs equal to cycles which induce the top graded or canonical elements of the Floer homology of $L(m,1)\# (S^1\times S^2)^{\# k}$ or $(S^1\times S^2)^{\#k}$. For $(\gs,\gs')$ and $(\gs',\gs'')$, the top graded and canonical cycles are each represented by a unique intersection point. For $(\ds,\ds')$, this will usually not be the case, since we wind $\Delta_1,\Delta_2\subset D(F)$ to achieve admissibility.

\begin{lem}\label{lem:P-1-admissible} If $\bar{\bs}$ are sufficiently wound on $\bar{\Sigma}$, and  $\Delta_1$ and $\Delta_2$ are sufficiently wound in $D(F)$, then the tuple $(\Sigma\# \bar{\Sigma},\gs'',\gs',\gs,\ds,\ds',w)$ is weakly admissible.
\end{lem}
\begin{proof} Following the proof of Lemma~\ref{lem:admissibility-C-2-h-1}, if we wind $\bar{\bs}$, $\Delta_1$ and $\Delta_2$ sufficiently, we reduce the main claim to showing that every non-zero rational periodic domain with boundary equal to a linear combination of the curves in $\as\cup \{c_1,c_2, c_3, \bar{c}_3\}\cup \hat{\Ds}$ has both positive and negative multiplicities. Since $\hat{\Ds}$ are adapted to $\as$, it is easy to see that this is equivalent to showing that any non-zero domain with equal to a linear combination of $c_1$, $c_2$, $c_3$, and $\bar{c}_3$ has positive and negative coefficients. This follows immediately from Lemma~\ref{lem:admissible-key-computation}.
\end{proof}

We may equivalently think of $\cC_{\aux}^{(1)}$ as the pairing of alpha and beta hypercubes of attaching curves shown below:
\begin{equation}
\cL_{\a}^{(1)}:=\begin{tikzcd}[row sep=1.5 cm, column sep=1.5cm]\gs\arrow[d, "\Theta_{\g',\g}^{\can}"]\\
\gs'\arrow[r, "\Theta_{\g'',\g'}^+"]& \gs''
\end{tikzcd}\qquad 
\cL_{\b}^{(1)}:=
 \begin{tikzcd}
  \ds\arrow[r, "\Theta_{\dt,\dt'}^+"]&\ds'
  \end{tikzcd}
  \label{eq:def:C-1-aux-as-pairing}
\end{equation}

\begin{lem} 
The diagram $\cC_{\aux}^{(1)}$ is a hypercube of chain complexes.
\end{lem}
\begin{proof}
Given the description of $\cC_{\aux}^{(1)}$ as the pairing of the two diagrams in equation~\eqref{eq:def:C-1-aux-as-pairing}, by Lemma~\ref{lem:hypercube-diagram-to-complex} it remains only to show that the diagrams $\cL_{\a}^{(1)}$ and $\cL_{\b}^{(1)}$ are hypercubes of attaching curves. The only non-trivial relation is
\[
f_{\g'',\g',\g}(\Theta_{\g'',\g'}^+, \Theta_{\g',\g}^{\can})=0,
\]
which follows from Lemma~\ref{lem:model-computations-m-surgery}.
\end{proof}

\begin{figure}[ht!]
	\centering
	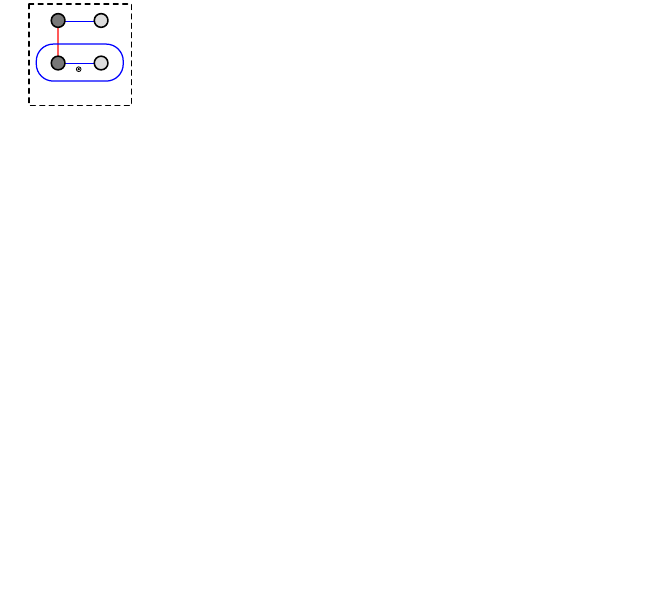
	\caption{The hypercube $\cC_{\cen}^{(1)}$. The character $\twind$ indicates winding.}\label{fig:176}
\end{figure}

Let $\Sigma_0$ denote $\Sigma$ with the surgery region removed, and let $\gs_0$ denote the curves obtained from $\gs'$ by removing the surgery region. Let $\ds_0$ denote the curves on $\Sigma_0$ obtained from $\ds'$ by removing the surgery region.
There is a stabilization map on the twisted complexes
\[
\sigma\colon \buCF^-(\Sigma_0\# \bar{\Sigma}, \gs_0,\ds_0)\to \buCF^-(\Sigma\# \bar{\Sigma}, \gs',\ds'),
\]
defined in the same manner as on the untwisted complexes. The hypercube $\cC_{\cen}^{(1)}$ is shown in Figure~\ref{fig:176}.

\begin{lem} \label{lem:C-c-1-hypercube}
The diagram $\cC_{\cen}^{(1)}$ in Figure~\ref{fig:176} is a hypercube of chain complexes.
\end{lem}
\begin{proof} Before proving the hypercube relations, we make the following claim:
\begin{equation}
\sigma F_3 f_{\g'\to \g''}^{\dt'}=\id_{\buCF^-(\gs',\ds')} \label{eq:C-1-retract}
\end{equation}
Equation~\eqref{eq:C-1-retract} is proven by using Proposition~\ref{prop:stabilize-triangles-general} to identify the composition $\sigma F_3 f_{\g'\to \g''}^{\dt'}$ with the count of triangles after destabilizing the triple for $f_{\g'\to \g''}^{\dt'}$. Proposition~\ref{prop:nearest-point-triangles} identifies this count with a nearest point map.

We now verify the hypercube relations. The length 0 and 1 relations are straightforward.  The length 2 relation along the left and back faces are trivial. The length 2 relation along the top face follows from Lemma~\ref{lem:omega-a-null-homotopy}. The length 2 relation along the bottom face is a consequence of the hypercube relations for $\cC_{\aux}^{(1)}$.

We now consider the length 2 relation for the front face. The desired relation reads
\begin{equation}
\left[\d, \sigma F_3 h_{\g'\to \g''}^{\dt\to \dt'}\right]=\sigma F_3 f_{\g''}^{\dt\to \dt'} f^{\dt}_{\g'\to \g''}+ f_{\g'}^{\dt\to \dt'}.\label{eq:hyperbox1-1n}
\end{equation}
To establish ~\eqref{eq:hyperbox1-1n}, we note that associativity of holomorphic rectangles gives
\begin{equation}
\left[\d, h_{\g'\to\g''}^{\dt\to \dt'}\right] =f_{\g'\to \g''}^{\dt'}f_{\g'}^{\dt\to \dt'}+f_{\g''}^{\dt\to \dt'}f_{\g'\to\g''}^{\dt}.\label{eq:hyperbox1-2n}
\end{equation}
We compose ~\eqref{eq:hyperbox1-2n} with $\sigma F_3$, and then apply~\eqref{eq:C-1-retract}, and immediately obtain ~\eqref{eq:hyperbox1-1n}.

Next, we verify the length 2 hypercube relation along the right face. The desired relation is
\begin{equation}
\left[\d, \sigma F_3 h_1 f_{\g}^{\dt\to \dt'}\right]= f_{\g\to \g'}^{\dt'} f_{\g}^{\dt\to \dt'}.\label{eq:hyperbox1-3n}
\end{equation}
To prove ~\eqref{eq:hyperbox1-3n}, we begin with by the associativity relations, applied to the map $h_1:=h_{\g\to \g'\to \g''}^{\dt'}$, which read
\begin{equation}
\left[\d, h_1\right]=f_{\g'\to \g''}^{\dt'}f_{\g\to \g'}^{\dt'}\label{eq:hyperbox1-4n}.
\end{equation}
We now postcompose ~\eqref{eq:hyperbox1-4n} with $\sigma F_3$, and precompose with $f_{\g}^{\dt\to \dt'}$ to obtain
\begin{equation}
\left[\d, \sigma F_3 h_1 f_{\g}^{\dt\to \dt'}\right]= \sigma F_3 f_{\g'\to \g''}^{\dt'} f_{\g\to \g'}^{\dt'}f_{\g}^{\dt\to \dt'}. \label{eq:hyperbox1-5n}
\end{equation}
Applying~\eqref{eq:C-1-retract} to ~\eqref{eq:hyperbox1-5n} gives the desired~\eqref{eq:hyperbox1-3n}.

Finally, we note by post-composing the length 3 hypercube relation of $\cC_{\aux}^{(1)}$ with $\sigma F_3$, we obtain
\begin{equation}
\begin{split}
&\sigma F_3 f_{\g'\to \g''}^{\dt'}h_{\g\to \g'}^{\dt\to \dt'}\\
+& \sigma F_3 f_{\g''}^{\dt\to\dt'} h_{\g\to \g'\to \g''}^{\dt}\\
+& \sigma F_3 h_{\g'\to \g''}^{\dt\to \dt'} f_{\g\to \g'}^{\dt}\\
+&\sigma F_3 h_{\g\to \g'\to \g''}^{\dt'} f_{\g}^{\dt\to \dt'}\\
=&\left[\d, \sigma F_3 p_{\g\to \g'\to \g''}^{\dt\to \dt'}\right],
\end{split}
\label{eq:hyperbox1-6n}
\end{equation}
which is exactly the desired length 3 relation for $\cC_{\cen}^{(1)}$ (we remind the reader that $h_{\g\to \g'\to \g''}^{\dt'}:=h_1$).
\end{proof}

\subsection{The hypercube \texorpdfstring{$\cC_{\cen}^{(2)}$}{C-cen-(2)}}

\begin{figure}[h]
	\centering
	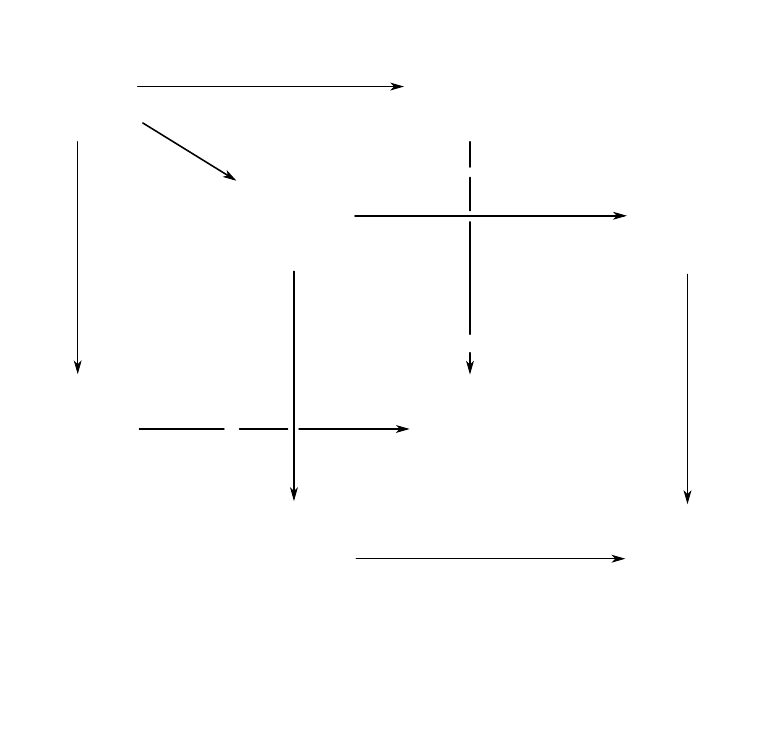
	\caption{The auxiliary hypercube $\cC_{\aux}^{(2)}$. The character $\mathrm{W}$ indicates winding.}\label{fig:177}
\end{figure}

 We consider the attaching curves $\gs'$, $\gs$, $\ds$, $\ds'$, $\tilde{\ds}'$ and $\ds''$, shown in Figure~\ref{fig:177}. We build an auxiliary hypercube $\cC^{(2)}_{\aux}$, shown schematically in Figure~\ref{fig:177}. Using the maps therein, we build the hypercube $\cC_{\cen}^{(2)}$, shown in Figure~\ref{fig:178}.

\begin{lem} If $\bar{\bs}\subset \bar{\Sigma}$ and $\Delta_1,\Delta_2\subset D(F)$ are sufficiently wound, then the tuple
\[
(\Sigma\# \bar{\Sigma}, \gs',\gs,\ds,\ds',\tilde{\ds}',\ds'',w)
\]
is weakly admissible.
\end{lem}
\begin{proof} The proof is not substantially different than the proof of Lemma~\ref{lem:P-1-admissible}.
\end{proof}

We will define $\cC_{\aux}^{(2)}$ as the pairing of two hypercubes of attaching curves, taking the following forms:
\begin{equation}
\begin{tikzcd}[row sep=1.5 cm]
\cL_{\a}^{(2)}:=\gs
\arrow[r, "\Theta_{\g',\g}^{\can}"]& \gs'
\end{tikzcd}
\quad \text{and}\quad 
\cL_{\b}^{(2)}:=
\begin{tikzcd}[row sep=1.5 cm, column sep =1.5 cm] 
\ds
	\arrow[r, "\Theta_{\dt,\dt'}^+"]
	\arrow[d, "\Theta_{\dt,\tilde{\dt}'}^+"]
	\arrow[dr,dashed, "\lambda_{\dt, \dt''}"]	
&\ds'
	\arrow[d, "\Theta_{\dt',\dt''}^+"]
\\
\tilde{\ds}'
	\arrow[r, "\Theta_{\tilde{\dt}',\dt''}^+"]
& \ds''
\end{tikzcd}
\label{eq:C-2-aux-hypercubes-attaching-curves}
\end{equation}
We now explain how to pick the chains appearing in~\eqref{eq:C-2-aux-hypercubes-attaching-curves}.  First, we pick cycles
  \[
  \Theta_{\g',\g}^{\can}, \quad \Theta_{\dt,\dt'}^+,\quad \Theta_{\dt',\dt''}^+, \quad \Theta^+_{\dt,\tilde{\dt}'},\quad \Theta^+_{\tilde{\dt}',\dt''},\quad \text{and}\quad \Theta_{\dt,\dt''}^+,
  \]
	each representing the top degree or canonical elements of their homology groups. Note that since we wound the curves $\Delta_1$ and $\Delta_2$, the chains $\Theta_{\dt,\dt'}^+$, $\Theta_{\dt,\tilde{\dt}'}^+$ and $\Theta_{\dt,\dt''}^+$  may not be represented by a canonical intersection point. For these chains, we pick a cycle which represents the top degree or canonical class arbitrarily.
  
\begin{figure}[h]
	\centering
\begingroup%
  \makeatletter%
  \providecommand\color[2][]{%
    \errmessage{(Inkscape) Color is used for the text in Inkscape, but the package 'color.sty' is not loaded}%
    \renewcommand\color[2][]{}%
  }%
  \providecommand\transparent[1]{%
    \errmessage{(Inkscape) Transparency is used (non-zero) for the text in Inkscape, but the package 'transparent.sty' is not loaded}%
    \renewcommand\transparent[1]{}%
  }%
  \providecommand\rotatebox[2]{#2}%
  \newcommand*\fsize{\dimexpr\f@size pt\relax}%
  \newcommand*\lineheight[1]{\fontsize{\fsize}{#1\fsize}\selectfont}%
  \ifx\svgwidth\undefined%
    \setlength{\unitlength}{251.92445382bp}%
    \ifx\svgscale\undefined%
      \relax%
    \else%
      \setlength{\unitlength}{\unitlength * \real{\svgscale}}%
    \fi%
  \else%
    \setlength{\unitlength}{\svgwidth}%
  \fi%
  \global\let\svgwidth\undefined%
  \global\let\svgscale\undefined%
  \makeatother%
  \begin{picture}(1,0.51337649)%
    \lineheight{1}%
    \setlength\tabcolsep{0pt}%
    \put(0,0){\includegraphics[width=\unitlength,page=1]{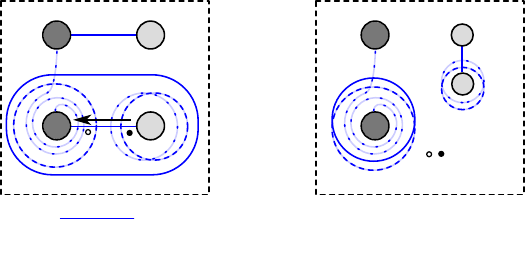}}%
    \put(0.27252864,0.0902096){\color[rgb]{0,0,0}\makebox(0,0)[lt]{\lineheight{1.25}\smash{\begin{tabular}[t]{l}$\ds$\end{tabular}}}}%
    \put(0,0){\includegraphics[width=\unitlength,page=2]{fig189.pdf}}%
    \put(0.27252864,0.04650401){\color[rgb]{0,0,0}\makebox(0,0)[lt]{\lineheight{1.25}\smash{\begin{tabular}[t]{l}$\ds'$\end{tabular}}}}%
    \put(0,0){\includegraphics[width=\unitlength,page=3]{fig189.pdf}}%
    \put(0.27252864,0.00279877){\color[rgb]{0,0,0}\makebox(0,0)[lt]{\lineheight{1.25}\smash{\begin{tabular}[t]{l}$\ds''$\end{tabular}}}}%
    \put(0,0){\includegraphics[width=\unitlength,page=4]{fig189.pdf}}%
  \end{picture}%
\endgroup%

	\caption{Manipulating the triple $(\Sigma\# \bar{\Sigma}, \ds, \ds',\ds'',w,z)$. The arrow indicates that we move the bottom right handle (light gray), across the bottom left handle (dark gray), via a sequence of handleslides. }\label{fig:189}
\end{figure}

\begin{figure}[h]
	\centering
	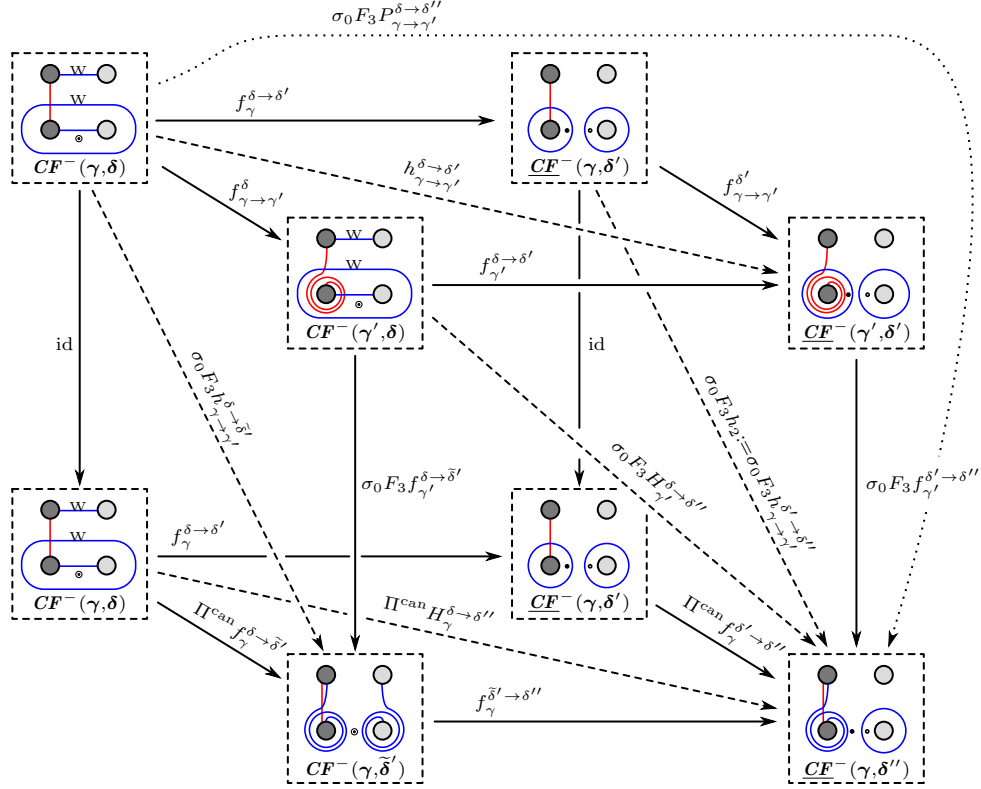
	\caption{The central hypercube $\cC_{\cen}^{(2)}$.}\label{fig:178}
\end{figure}

 \begin{lem}\label{lem:C-2-equal-on-homology}
\begin{equation}
\left[f_{\dt,\dt',\dt''}\left(\Theta_{\dt,\dt'}^+,\Theta_{\dt',\dt''}^+\right)\right]=\left[f_{\dt,\tilde{\dt}',\dt''}\left(\Theta_{\dt,\tilde{\dt}'}^+,\Theta_{\tilde{\dt}',\dt''}^+\right)\right]=\left[\Theta_{\dt,\dt''}^+\right]. \label{eq:model-computation-triangles-C-2}
\end{equation}
In~\eqref{eq:model-computation-triangles-C-2}, the brackets denote the induced elements of the homology group $\buHF^-(\ds,\ds'')$.
\end{lem}
\begin{proof}
The triple $(\Sigma\# \bar{\Sigma}, \ds, \ds',\ds'',w)$ may be interpreted as corresponding to the link cobordism obtained by attaching a 2-handle to $((S^1\times S^2)^{\# 2g-1},U)$, (where $U$ denotes an unknot), which cancels one of the $S^1\times S^2$ summands, and is unlinked from the unknot. This can be seen by handlesliding curves and manipulating the diagram as illustrated in Figure~\ref{fig:189}.  A similar manipulation works for the triple $(\Sigma\# \bar{\Sigma}, \ds, \tilde{\ds}', \ds'')$.
 \end{proof}
 
Lemma~\ref{lem:C-2-equal-on-homology} implies that there is a chain $\lambda_{\dt,\dt''}\in \buCF^-(\ds,\ds'')$, with grading 1 higher than $\Theta^+_{\dt,\dt''}$, such that
\begin{equation}
\d (\lambda_{\dt,\dt''})=f_{\dt,\dt',\dt''}\left(\Theta_{\dt,\dt'}^+,\Theta_{\dt',\dt''}^+\right)+f_{\dt,\tilde{\dt}',\dt''}\left(\Theta_{\dt,\tilde{\dt}'}^+,\Theta_{\tilde{\dt}',\dt''}^+\right). \label{eq:def-eta-d-d''}
\end{equation}
The chain $\lambda_{\dt,\dt''}$ satisfying equation~\eqref{eq:def-eta-d-d''} is the one which features in the hypercube $\cL^{(2)}_{\b}$ in~\eqref{eq:C-2-aux-hypercubes-attaching-curves}. We define the maps in $\cC_{\aux}^{(2)}$ to be the ones obtained by pairing $\cL_{\a}^{(2)}$ and $\cL_{\b}^{(2)}$, via~\eqref{eq:pairing-attaching-curves}.

\begin{lem}\label{lem:C-aux-2-hypercube}
 The diagram $\cC_{\aux}^{(2)}$ in Figure~\ref{fig:177} is a hypercube of chain complexes.
\end{lem}
 Lemma~\ref{lem:C-aux-2-hypercube} follows immediately from the hypercube relations for $\cL_{\a}^{(2)}$ and $\cL_{\b}^{(2)}$, as well as Lemma~\ref{lem:hypercube-diagram-to-complex}.

We now give a more concrete description of the maps appearing in $\cC_{\aux}^{(2)}$. All of the length 1 maps are elementary triangle maps. For example,
\[
f_{\g}^{\dt\to \dt'}=f_{\g,\dt,\dt'}(-,\Theta_{\dt,\dt'}^+).
\]
The length 2 maps $h_{\g\to \g'}^{\dt\to \tilde{\dt}'}$, $h_{\g\to \g'}^{\dt\to \dt'}$, $h_{\g\to \g'}^{\dt\to \dt'}$ and $h_{\g\to \g'}^{\tilde{\dt}'\to \dt''}$ are each given by a single quadrilateral counting map. 

The maps $H_{\g}^{\dt\to \dt''}$ and $H_{\g'}^{\dt\to \dt''}$ are more complicated. They are given by the formulas
\[
H_{\g}^{\dt\to \dt''}=h_{\g}^{\dt\to \dt'\to \dt''}+h_{\g}^{\dt\to \tilde{\dt}'\to \dt''}+f_{\g}^{\dt\to \dt''},
\]
where $f_{\g}^{\dt\to \dt''}:=f_{\g,\dt,\dt''}(-,\lambda_{\dt,\dt''})$. The map $H_{\g'}^{\dt\to \dt''}$ is defined similarly.

Finally, the length 3 map is given by
\[
P_{\g\to \g'}^{\dt\to \dt''}=p_{\g\to \g'}^{\dt\to \dt'\to \dt''}+p_{\g\to \g'}^{\dt\to \tilde{\dt}'\to \dt''}+h_{\g\to \g'}^{\dt\to \dt''},
\]
where $h_{\g\to \g'}^{\dt\to \dt''}=h_{\g',\g,\dt,\dt''}(\Theta_{\g',\g}^{\can},-,\lambda_{\dt,\dt''})$.

In Figure~\ref{fig:178}, we define $\cC_{\cen}^{(2)}$ using the maps from $\cC^{(2)}_{\aux}$. The map $F_3$ is the 3-handle map, and $\sigma_0$ is the lens space stabilization map.
Also, $\Pi^{\can}$ denotes projection onto the canonical $\Spin^c$ structure of the $L(m,1)$ summand. More explicitly,
\begin{equation}
\Pi^{\can}(\xs\times x_i)=\begin{cases} \xs\times x_i& \text{if } i=0\\
0& \text{otherwise},
\end{cases}
\label{eq:Pi-can}
\end{equation}
extended equivariantly over $T$ and $U$.

\begin{lem}
 The diagram $\cC_{\cen}^{(2)}$ in Figure~\ref{fig:177} is a hypercube of chain complexes.
\end{lem}

\begin{proof}
There are two main inputs to the proof. The first is, naturally, the hypercube relations for $\cC_{\aux}^{(2)}$. The second input consists of several results concerning the simplicity of the maps appearing along the bottom face of the hypercube $\cC_{\aux}^{(2)}$. More precisely, we claim
\begin{equation}
\begin{split}
\sigma_0 F_3 f_{\g\to \g'}^{\tilde{\dt}'}&=\Pi^{\can}\colon \bCF^-(\gs,\tilde{\ds}')\to \bCF^-(\gs,\tilde{\ds}'), \\
\sigma_0 F_3 f_{\g\to \g'}^{\dt''}&=\Pi^{\can}\colon \buCF^-(\gs,\ds'')\to \buCF^-(\gs,\ds''),
\label{eq:left-inverse}
\end{split}
\end{equation}
\begin{equation}
f_{\g}^{\tilde{\dt}'\to \dt''} \sigma_0 F_3=  \sigma_0 F_3 f_{\g'}^{\tilde{\dt'}\to \dt''}\colon \bCF^-(\gs',\tilde{\ds}')\to \bCF^-(\gs,\ds''),
\label{eq:inverses-commute}
\end{equation}
and
\begin{equation}
F_3 h_{\g\to \g'}^{\tilde{\dt}'\to \dt''}=0. \label{eq:quadrilateral-count-vanishes}
\end{equation}
 Equation~\eqref{eq:left-inverse} is proven by using Lemma~\ref{lem:lens-stabilization-triangles} to relate the triangle counts for $f_{\g\to \g'}^{\tilde{\dt}'}$ and $f_{\g\to \g'}^{\dt''}$ with triangle counts on the destabilized complex. Proposition~\ref{prop:nearest-point-triangles} identifies the destabilized triangle counts with nearest point maps, giving the stated formula. Equation~\eqref{eq:inverses-commute} is similarly proven by applying Lemma~\ref{lem:lens-stabilization-triangles} and Proposition~\ref{prop:stabilize-triangles-general} to the genus 1 region bounded by the curve $S_1$ to identify both compositions with a triangle count on the same genus $2g(\Sigma)-1$ Heegaard triple. Finally,~\eqref{eq:quadrilateral-count-vanishes} is proven by using Lemma~\ref{lem:len-space-stabilization-quadrilaterals} to identify the map $F_3 h_{\g\to \g'}^{\tilde{\dt}'\to \dt''}$ with the count of quadrilaterals on the quadruple obtained by removing the genus 1 region bounded by $S_1$. Note that outside of the genus 1 region bounded by $S$, the curves $\gs$ and $\gs'$ are pairwise Hamiltonian isotopic (they coincide with small translates of $\as\cup \bar{\bs}$). Hence Proposition~\ref{prop:small-rectangles} implies that the count is zero, for $\gs$ and $\gs'$ chosen to be suitably small translates of each other outside of the surgery region.

We now prove the hypercube relations. The length 0 and 1 hypercube relations are immediate. The length 2 relations along the top and back faces follow immediately from the hypercube relations for $\cC_{\aux}^{(2)}$. The length 2 relation along the front face follows by taking the length 2 relation along the front face of $\cC_{\aux}^{(2)}$ and post-composing with $\sigma_0 F_3$.

The length 2 relation along the bottom face follows immediately from the relation
\begin{equation}
f_{\g}^{\tilde{\dt}'\to \dt''}\Pi^{\can}=\Pi^{\can} f_{\g}^{\tilde{\dt}'\to \dt''},
\label{eq:commute-Pi-can}
\end{equation}
which follows because the map $f_{\g}^{\tilde{\dt}'\to \dt''}$ is the cobordism map for a 2-handle whose attaching cycle is disjoint from the $L(m,1)$ summand of $(\gs,\tilde{\ds}')$, and hence we may use the lens space stabilization result for triangles in Lemma~\ref{lem:lens-stabilization-triangles}.

We now consider the length 2 relation along the left face. The desired relation is
\begin{equation}
\left[\d, \sigma_0 F_3 h_{\g\to \g'}^{\dt\to \tilde{\dt}'}\right]=\sigma_0 F_3 f_{\g'}^{\dt\to \tilde{\dt}'} f_{\g\to \g'}^{\dt}+\Pi^{\can}f_{\g}^{\dt\to \tilde{\dt}'}.
\label{eq:hypercube-2-1n}
\end{equation}
To prove equation~\eqref{eq:hypercube-2-1n}, we begin with the relation
\begin{equation}
\left[\d, h_{\g\to \g'}^{\dt\to \tilde{\dt}'}\right]=f_{\g'}^{\dt\to \tilde{\dt}'}f_{\g\to \g'}^{\dt}+f_{\g\to \g'}^{\tilde{\dt}'} f_{\g}^{\dt\to \tilde{\dt}'}. \label{eq:hypercube-2-2n}
\end{equation} Post-composing with $\sigma_0 F_3$, ~\eqref{eq:hypercube-2-2n} becomes
\[
\left[\d, \sigma_0 F_3 h_{\g\to \g'}^{\dt\to \tilde{\dt}'}\right]=\sigma_0 F_3 f_{\g'}^{\dt\to \tilde{\dt}'}f_{\g\to \g'}^{\dt}+\sigma_0 F_3f_{\g\to \g'}^{\tilde{\dt}'} f_{\g}^{\dt\to \tilde{\dt}'}.
\]
Applying ~\eqref{eq:left-inverse} immediately gives~\eqref{eq:hypercube-2-1n}. The length 2 relation along the right face of $\cC_{\cen}^{(2)}$ is proved entirely analogously.

We now consider the length 3 relation. The desired relation is
\begin{equation}
\begin{split}
\left[\d, \sigma_0 F_3 P_{\g\to \g'}^{\dt\to \dt''}\right]&=f_{\g}^{\tilde{\dt}'\to \dt''} \sigma_0 F_3 h_{\g\to \g'}^{\dt\to \tilde{\dt}'}\\
&+\sigma_0F_3 f_{\g'}^{\dt'\to \dt''} h_{\g\to \g'}^{\dt\to \dt'}\\
&+\Pi^{\can}H_{\g}^{\dt\to \dt''}\\
&+\sigma_0 F_3H_{\g'}^{\dt\to  \dt''}f_{\g\to \g'}^{\dt}\\
&+\sigma_0 F_3 h_{\g\to \g'}^{\dt'\to \dt''} f_{\g}^{\dt\to \dt'}
\end{split}
\label{eq:hypercube-C-2-len-3n}
\end{equation}
To obtain ~\eqref{eq:hypercube-C-2-len-3n}, we  post-compose the length 3 relation of $\cC_{\aux}^{(2)}$ with $\sigma_0 F_3$ to obtain
\begin{equation}
\begin{split}
\left[\d, \sigma_0 F_3 P_{\g\to \g'}^{\dt\to \dt''} \right]&=\sigma_0 F_3 f_{\g'}^{\tilde{\dt}'\to \dt''} h_{\g\to \g'}^{\dt\to \tilde{\dt}'}+\sigma_0 F_3 f_{\g'}^{\dt'\to \dt''} h_{\g\to \g'}^{\dt\to \dt'}\\
&+\sigma_0 F_3 h_{\g\to \g'}^{\dt'\to \dt''} f_{\g}^{\dt\to \dt'}+ \sigma_0 F_3 h_{\g\to \g'}^{\tilde{\dt}'\to \dt''} f_{\g}^{\dt\to \tilde{\dt}'}\\
&+\sigma_0 F_3 f_{\g\to \g'}^{\dt''} H_{\g}^{\dt\to \dt''}\\
&+\sigma_0 F_3H_{\g'}^{\dt\to \dt''}f_{\g\to \g'}^{\dt}.
\end{split}
\label{eq:hypercube-C-2-3n}
\end{equation}
The equivalence of~\eqref{eq:hypercube-C-2-3n} and~\eqref{eq:hypercube-C-2-len-3n} follows quickly from ~\eqref{eq:left-inverse}, \eqref{eq:quadrilateral-count-vanishes} and~\eqref{eq:inverses-commute}.
\end{proof}

\subsection{The hypercube \texorpdfstring{$\cC^{(3)}_{\cen}$}{C-cen-(3)}}
\label{sec:hypercube-Cc3}

\begin{figure}[ht!]
	\centering
	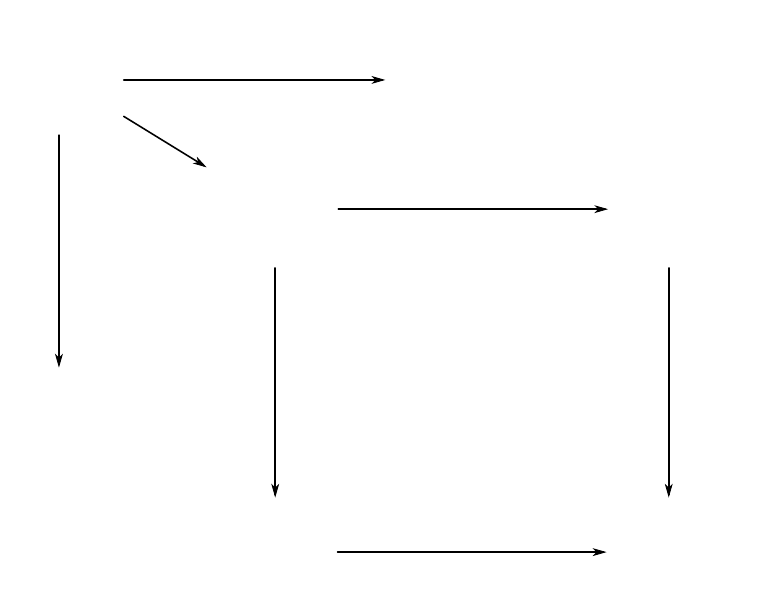
	\caption{The auxiliary hypercube $\cC_{\aux}^{(3)}$.}\label{fig:182}
\end{figure}

In this section, we construct the hypercube $\cC_{\cen}^{(3)}$. We first construct an auxiliary hypercube $\cC_{\aux}^{(3)}$, shown schematically in Figure~\ref{fig:182}.

\begin{lem}\label{lem:C-3-aux-admissible}
 If $\bar{\bs}$ are sufficiently wound on $\bar{\Sigma}$, and $\Delta_1$ and $\Delta_2$ are sufficiently wound on $\bar{\Sigma}$, then the diagram
\[
(\Sigma\# \bar{\Sigma},\gs,\ds, \ds',\ds'',\tilde{\ds}',\ds_1,\ds_2, \ds_3,w)
\]
is weakly admissible.
\end{lem}
The proof of Lemma~\ref{lem:C-3-aux-admissible} is essentially the same as the proof of Lemma~\ref{lem:P-1-admissible}, so we omit it.

The hypercube $\cC_{\aux}^{(3)}$ is the pairing of the 0-dimensional hypercube $\cL_\a^{(3)}$ of alpha attaching curves consisting of $\gs$, as well the following 3-dimensional hypercube of beta attaching curves:
\begin{equation}
\cL_{\b}^{(3)}:=\begin{tikzcd}
[column sep={2.4cm,between origins},
row sep={1cm,between origins},
labels=description]
\ds
	\arrow[dd, swap,"\Theta_{\dt,\dt_1}^+"]
	\arrow[dddrrr,dotted,pos=.55, "\omega_{\dt,\dt_3}"]
	\arrow[dr, swap,"\Theta^+_{\dt,\tilde{\dt}'}"]
	\arrow[rr, "\Theta^+_{\dt,\dt'}"]
	\arrow[drrr,dashed, "\lambda_{\dt,\dt_2}"]
&&[-1.5cm]
\ds'
	\arrow[dr, "\Theta_{\dt',\dt''}^+"]
&
\\
&\tilde{\ds}'
	\arrow[ddrr,dashed, "\lambda_{\tilde{\dt}',\dt_3}"]
	\arrow[rr,crossing over, "\Theta_{\tilde{\dt}',\dt''}^+",swap]
&&
\ds''
	\arrow[dd, "\Theta_{\dt'',\dt_3}^{\can}"]
\\[1.5cm]
\ds_1
	\arrow[dr, "\Theta_{\dt_1,\dt_2}^{\can}"]
	&&\,&
	\\
&
\ds_2
	\arrow[rr, "\Theta_{\dt_2,\dt_3}^+"]
	\arrow[from =uu, crossing over, "\Theta_{\tilde{\dt}',\dt_2}^{\can}"]
	\arrow[from=uuul,dashed,crossing over, "\lambda_{\dt,\dt_2}"]
	&&
\ds_3
\end{tikzcd}
\label{eq:L-b-(3)}
\end{equation}
We now describe the chains appearing in $\cL^{(3)}_{\b}$. First, we pick a collection of cycles
\[
\Theta_{\dt,\tilde{\dt}'}^+,\quad  \Theta_{\dt,\dt'}^+,\quad \Theta_{\dt',\dt''}^+,\quad  \Theta_{\tilde{\dt}',\dt''}^+, \quad \Theta_{\dt,\dt_1}^+,\quad  \Theta_{\dt',\dt_2}^{\can},\quad  \Theta_{\dt'', \dt_3}^{\can},\quad \Theta_{\dt_1,\dt_2}^{\can},\quad \text{and} \quad  \Theta_{\dt_2,\dt_3}^+,
\]
which each represent the top degree generator or canonical generator of their homology groups. We use these cycles for the length 1 cycles in ~\eqref{eq:L-b-(3)}.

\begin{figure}[ht!]
	\centering
	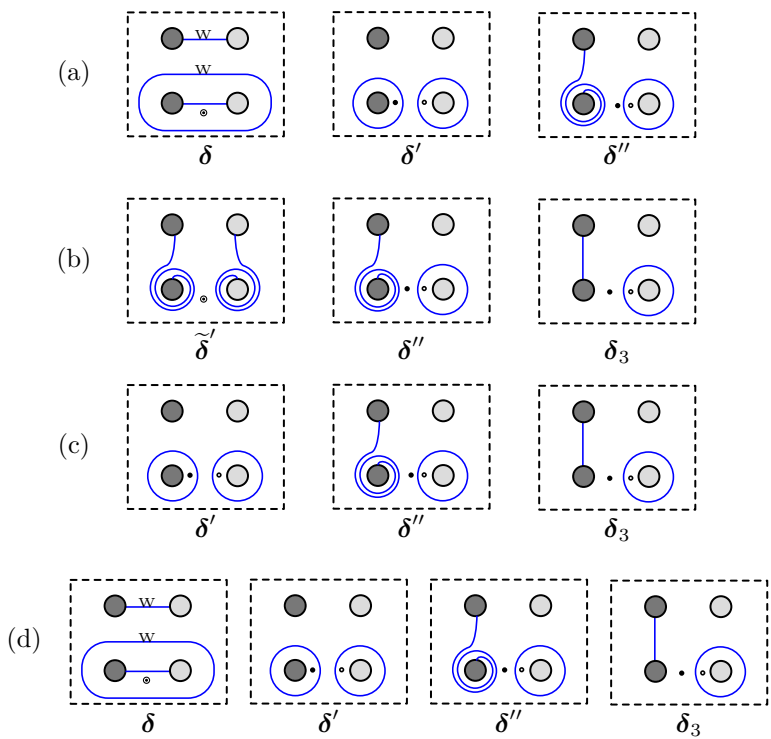
	\caption{Examples of the triples and quadruples from Lemma~\ref{lem:model-computations-for-C3}. (a) The triple $(\ds,\ds',\ds'')$, which represents surgery on  $S^1\times \{pt\}\subset S^1\times S^2$. (b) The triple $(\tilde{\ds}',\ds'',\ds_3)$ represents surgery on $S^1\times \{pt\}\subset S^1\times S^2\# L(m,1)$. (c) The triple $(\ds',\ds'',\ds_3)$ represents surgery on a knot in $L(m,1)$ which gives $S^3$. (d) The quadruple $(\ds,\ds',\ds'',\ds_3)$.}\label{fig:94}
\end{figure}

\begin{lem}\label{lem:model-computations-for-C3} Suppose $\ds$, $\tilde{\ds}'$, $\ds'$, $\ds''$, $\ds_1$, $\ds_2$ and $\ds_3$ are the attaching curves appearing in $\cL_{\b}^{(3)}$. Then \sloppy
\[
\begin{split}
\left[f_{\dt,\dt',\dt''}\left(\Theta_{\dt,\dt'}^+,\Theta_{\dt',\dt''}^+\right)\right]&=\left[\Theta_{\dt,\dt''}^+\right],\\
f_{\dt',\dt'',\dt_3}\left(\Theta_{\dt',\dt''}^+,\Theta_{\dt'',\dt_3}^{\can}\right)&=0,\\
\\
\left[f_{\dt,\tilde{\dt}',\dt''}\left(\Theta_{\dt,\tilde{\dt}'}^+,\Theta_{\tilde{\dt}',\dt''}^+\right)\right]&=\left[\Theta_{\dt,\dt''}^+\right],\\
\left[f_{\tilde{\dt}',\dt'',\dt_3}\left(\Theta_{\tilde{\dt}',\dt''}^+,\Theta_{\dt'',\dt_3}^{\can}\right)\right]&=\left[\Theta_{\tilde{\dt}',\dt_3}^{\can}\right],\\
\\
\left[f_{\dt,\tilde{\dt}',\dt_2}\left(\Theta_{\dt,\tilde{\dt}'}^+,\Theta_{\tilde{\dt}',\dt_2}^{\can}\right)\right]&=\left[\Theta_{\dt,\dt_2}^{\can}\right],\\
\left[f_{\tilde{\dt}',\dt_2,\dt_3}\left(\Theta_{\tilde{\dt}',\dt_2}^{\can},\Theta_{\dt_2,\dt_3}^+\right)\right]&=\left[\Theta_{\tilde{\dt}',\dt_3}^{\can}\right],\\
\\
\left[f_{\dt,\dt_1,\dt_2}\left(\Theta_{\dt,\dt_1}^+,\Theta_{\dt_1,\dt_2}^{\can}\right)\right]&=\left[\Theta_{\dt,\dt_2}^{\can}\right],\\
f_{\dt_1,\dt_2,\dt_3}\left(\Theta_{\dt_1,\dt_2}^{\can},\Theta_{\dt_2,\dt_3}^+\right)&=0,
\end{split}
\]
where the brackets denote the induced elements of homology.
\end{lem}
\begin{proof} First, recall that all almost complex structures in this section are maximally pinched along $S_1$ and $S_2$. Furthermore, outside of the genus two region pictured in Figure~\ref{fig:182}, all of the attaching curves appearing in the statement coincide with small translates of the $\hat{\Ds}$ curves. Hence, we can view all triples appearing in the claim as being stabilizations of the genus two region shown. Furthermore, the stabilizations have index tuple $(0,0,0)$. Hence, using Proposition~\ref{prop:stabilize-triangles-general}, it suffices to prove the claim for the corresponding genus 2 triples.

The computations are handled by one of three cases.
 
The first case occurs when the associated Heegaard triple represents surgery on $S^1\times \{pt\}\subset S^1\times S^2$, yielding $S^3$. The triples $(\ds,\ds',\ds'')$, $(\ds,\tilde{\ds}',\ds_2)$, and $(\ds,\tilde{\ds}',\ds'')$ have this configuration. By invariance of the cobordism maps under index 1 and 2 handle cancelations, the top degree generator is preserved, on homology. The triple $(\ds,\ds',\ds'')$ is shown in Figure~\ref{fig:94}.

The second case  occurs when the triple represents surgery on $S^1\times \{pt\}\subset (S^1\times S^2)\# L(m,1)$. By the same argument as in the first case, the element $\Theta^{\can}$ is preserved. The triples $(\tilde{\ds}',\ds'',\ds_3)$,   $(\ds,\ds_1,\ds_2)$ and $(\tilde{\ds}',\ds_2,\ds_3)$ satisfy this configuration. 

The third case is when the triple represents $-m$-framed surgery on an unknot in $S^3$, giving $L(m,1)$. The triples $(\ds_1,\ds_2,\ds_3)$ and $(\ds',\ds'',\ds_3)$ satisfy this condition. Furthermore, these triples are in standard position (i.e. the triples are not wound), so holomorphic triangles cancel modulo 2 by Lemma~\ref{lem:omega-a-null-homotopy}.
\end{proof}

Lemma~\ref{lem:model-computations-for-C3} implies that there are chains $\lambda_{\dt,\dt''}$, $\lambda_{\tilde{\dt}',\dt_3}$, and $\lambda_{\dt,\dt_2}$  satisfying
\begin{equation}
\begin{split}
\d(\lambda_{\dt,\dt''})&=f_{\dt,\dt',\dt''}\left(\Theta_{\dt,\dt'}^+,\Theta_{\dt',\dt''}^+\right)+f_{\dt,\tilde{\dt}',\dt''}\left(\Theta_{\dt,\tilde{\dt}'}, \Theta^+_{\tilde{\dt}',\dt''}\right),\\
\d(\lambda_{\tilde{\dt}',\dt_3})&=f_{\tilde{\dt}',\dt'',\dt_3}\left(\Theta_{\tilde{\dt}',\dt''}^+, \Theta^{\can}_{\dt'',\dt_3}\right) +f_{\tilde{\dt}',\dt_2,\dt_3}\left(\Theta^{\can}_{\tilde{\dt}',\dt_2}, \Theta_{\dt_2,\dt_3}^+\right),\\
\d(\lambda_{\dt,\dt_2})&=f_{\dt,\tilde{\dt}',\dt_2}\left(\Theta_{\dt,\tilde{\dt}'}^+, \Theta_{\tilde{\dt}',\dt_2}^{\can}\right)+f_{\dt,\dt_1,\dt_2}\left(\Theta_{\dt,\dt_1}^+,\Theta_{\dt_1,\dt_2}^{\can}\right).
\end{split}
\label{def:eta}
\end{equation}
Furthermore, since all of the triangles maps in~\eqref{def:eta} are homogeneously graded with respect to the Maslov grading, we can assume these chains live in Maslov grading one greater than $\Theta_{\dt,\dt_2}^{\can}$, $\Theta_{\dt, \dt''}^+$ and $\Theta^{\can}_{\tilde{\dt}',\dt_3}$. We use the chains from~\eqref{def:eta} for the length 2 chains in $\cL_{\b}^{(3)}$, appearing in ~\eqref{eq:L-b-(3)}.

Next, we describe how to construct the length 3-chain $\omega_{\dt,\dt_3}$ appearing in $\cL_{\b}^{(3)}$. To make our analysis of gradings easier, we will define $\omega_{\dt,\dt_3}$ as an infinite sum
\begin{equation}
\omega_{\dt,\dt_3}=\sum_{k=0}^\infty \omega_{\dt,\dt_3;\frs_k^+}, \label{eq:omega-d-d3-def}
\end{equation}
where each $\omega_{\dt,\dt_3;\frs_k^+}$ is homogeneously graded. Here, $\frs_k^+$ denotes the element of $\Spin^c(D(-m,1))$ defined in Lemma~\ref{lem:model-computations-m-surgery}. We think of $\omega_{\dt,\dt_3;\frs_k^-}$ as also being defined, but being 0. We write $\fru_k$ for $\{\frs_k^+,\frs_k^-\}$. 

Consider the four-ended 4-manifolds constructed from the genus 2 surgery regions of the Heegaard quadruples $(\ds,\ds',\ds'',\ds_3)$, $(\ds,\tilde{\ds}',\ds'',\ds_3)$, $(\ds,
\tilde{\ds}',\ds_2,\ds_3)$ and $(\ds,\ds_1,\ds_2,\ds_3)$. Note that outside of the special genus 2 subregion, the attaching curves are all small translations of $\hat{\Ds}$.  It is a straightforward exercise to see that each of the aforementioned 4-manifolds is obtained by attaching two 2-handles along a 2-component unlink in $S^3$, with framings $0$ and $-m$, and then removing the neighborhood of two 4-dimensional genus 1 handlebodies. Phrased another way, each of the 4-manifolds above have 2 boundary components equal to $S^1\times S^2$, such that when we fill in those components with 3 and 4-handles, we obtain the 2-handle cobordism described above. In particular, $\Spin^c$ structures on the 4-manifolds above, which restrict to the torsion and canonical $\Spin^c$ structures on the boundary, are naturally identified with $\Spin^c$ structures on the disk bundle $D(-m,1)$ which restrict to the canonical one on the boundary. 

\begin{lem}\label{lem:C-d-d3-boundary} For each $k\ge 0$, the chain
\[
\begin{split}
C_{\dt,\dt_3;\fru_k}:=&h_{\dt,\dt',\dt'',\dt_3;\fru_k}\left (\Theta_{\dt,\dt'}^+,\Theta_{\dt',\dt''}^+,\Theta_{\dt'',\dt_3}^{\can}\right)+h_{\dt,\tilde{\dt}',\dt'',\dt_3;\fru_k}\left(\Theta_{\dt,\tilde{\dt}'}^+,\Theta_{\tilde{\dt}',\dt''}^+,\Theta_{\dt'',\dt_3}^{\can}\right)\\
&+h_{\dt,\tilde{\dt}',\dt_2,\dt_3;\fru_k}\left(\Theta_{\dt,\tilde{\dt}'}^+,\Theta_{\tilde{\dt}',\dt_2}^{\can},\Theta_{\dt_2,\dt_3}^+\right)+h_{\dt,\dt_1,\dt_2,\dt_3;\fru_k}\left(\Theta_{\dt,\dt_1}^+,\Theta_{\dt_1,\dt_2}^{\can},\Theta_{\dt_2,\dt_3}^+\right)\\
&+f_{\dt,\dt'',\dt_3;\fru_k}(\lambda_{\dt,\dt''},\Theta_{\dt'',\dt_3}^{\can})+f_{\dt,\tilde{\dt}',\dt_3;\fru_k}\left(\Theta_{\dt,\tilde{\dt}'}^+,\lambda_{\tilde{\dt}',\dt_3}\right)+ f_{\dt,\dt_2,\dt_3;\fru_k}\left(\lambda_{\dt,\dt_2}, \Theta_{\dt_2,\dt_3}^+\right)
\end{split}
\]
is a boundary in $\buCF^-(\ds,\ds_3)$.
\end{lem}
\begin{proof} First, note that using our stabilization formulas for triangles and quadrilaterals in Proposition~\ref{prop:stabilize-triangles-general} and \ref{prop:stabilize-quadrilateral}, it is sufficient to show the claim for the genus 2 regions shown in Figure~\ref{fig:182}. Next, we claim $C_{\dt,\dt_3;\fru_k}$ is a cycle. This is proven by an easy associativity argument, together with the definition of the chains $\lambda_{\dt,\dt''}$, $\lambda_{\tilde{\dt}',\dt_3}$ and $\lambda_{\dt,\dt_2}$ in ~\eqref{def:eta}. The absolute grading formula of \cite{OSIntersectionForms} implies that $C_{\dt,\dt_3;\fru_k}$ has grading
\[
-mk(k+1)+1.
\]
(Compare Lemma~\ref{lem:model-computations-m-surgery}). However $\buHF^-(\ds,\ds_3)\iso \bF\llsquare U\rrsquare\otimes \bF[\Z/m]$ is non-trivial only in even gradings. Hence $C_{\dt,\dt_3;\fru_k}$ must be a boundary.
\end{proof}

Let $\omega_{\dt,\dt_3;\frs_k^+}\in \buCF^-(\ds,\ds_3)$ be any chain satisfying
\begin{equation}
\d(\omega_{\dt,\dt_3;\frs_k^+})=C_{\dt,\dt_3;\fru_k},\label{eq:def-omega-d-d3}
\end{equation}
which has homogeneous grading 1 greater than $C_{\dt,\dt_3;\fru_k}$. As indicated in ~\eqref{eq:omega-d-d3-def}, we define $\omega_{\dt,\dt_3}$ to be the sum over $k\ge 0$ of $\omega_{\dt,\dt_3;\frs_k^+}$. We use $\omega_{\dt,\dt_3}$ as the length 3 chain in the construction of $\cL_{\b}^{(3)}$.

\begin{lem} The diagram $\cL_{\b}^{(3)}$ is a hypercube of beta attaching curves. Furthermore, $\cC_{\aux}^{(3)}$ is a hypercube of chain complexes.
\end{lem}
\begin{proof} The length 1 relations for $\cL_{\b}^{(3)}$ are immediate.
 Lemma~\ref{lem:model-computations-for-C3} and ~\eqref{def:eta} together give the length 2 relations.
Lemma~\ref{lem:C-d-d3-boundary} and~\eqref{eq:def-omega-d-d3} give the length 3 relations of $\cL_{\b}^{(3)}$, concluding the proof.
\end{proof}

\begin{rem}
 Our construction of the chains appearing in $\cL_\b^{(3)}$ is similar to the procedure of \emph{filling} an empty $\beta$-hyperbox, in the terminology of Manolescu and Ozsv\'{a}th \cite{MOIntegerSurgery}*{Section~8}.
\end{rem}

Next, we modify $\cC_{\aux}^{(3)}$ slightly, and construct a hypercube $\cC_{\aux,\Pi^{\can}}^{(3)}$. The hypercube $\cC_{\aux,\Pi^{\can}}^{(3)}$ has the same complexes as $\cC_{\aux}^{(3)}$, except the maps along the top face are composed with $\Pi^{\can}$, as in the bottom face of $\cC_{\cen}^{(2)}$. We claim that $\cC_{\aux, \Pi^{\can}}^{(3)}$ is also a hypercube of chain complexes. To establish this, we note that $\cC_{\aux, \Pi^{\can}}^{(3)}$ coincides with the hypercube obtained by stacking $\cC_{\aux}^{(3)}$ with the hypercube shown in~\eqref{eq:hypercube-Pi-can}, and then compressing. Hence, it is sufficient to show the hypercube relations for the diagram in ~\eqref{eq:hypercube-Pi-can}. The  length 1 relations are clear. The length 2 relations along the back, front and bottom faces are also clear. The length 2 relation along the top face follows from~\eqref{eq:commute-Pi-can}. For the length 2 relations on the left and right faces, we note that the intersection points in the non-canonical $\Spin^c$ structures are in the kernel of $f_{\g}^{\tilde{\dt}'\to \dt_2}$ and $f_{\g}^{\dt''\to \dt_3}$. Indeed these two maps count holomorphic triangles on Heegaard triples which have a genus 1 lens space stabilizations, and by construction they have the intersection point in the canonical $\Spin^c$ structure as a special input, so there are no homology classes which could be counted unless the other input is also in the canonical $\Spin^c$ structure. The length 3 relation follows from the same argument, applied to the rectangle counting map for the front and back face.
\begin{equation}
\begin{tikzcd}[column sep={3.5cm,between origins},
row sep=.3cm,labels=description]
\CF^-(\gs,\tilde{\ds}')
	\ar[dd,"f_\g^{\tilde{\dt}'\to \dt_2}"]
	\ar[dr, "\Pi^{\can}"]
	\ar[rr, "f_{\g}^{\tilde{\dt}'\to \dt''}"]
	\ar[ddrr,dashed, sloped, "H_{\g}^{\tilde{\dt}'\to \dt_3}", pos=.4]
&&[-1.5cm]
\uCF^-(\gs,\ds'')
	\ar[dd, "f_{\g}^{\dt''\to \dt_3}"]
	\ar[dr, "\Pi^{\can}"]
&
\\
&\CF^-(\gs,\tilde{\ds}')
	\ar[rr,crossing over, "f_{\g}^{\tilde{\dt}'\to \dt''}"]
&&
\uCF^-(\gs,\ds'')
	\ar[dd, "f_{\g}^{\dt''\to \dt_3}"]
\\[2cm]
\CF^-(\gs,\ds_2)
	\ar[rr, "f_{\g}^{\dt_2\to \dt_3}", pos=.25]
	\ar[dr, "\id"]
	&
&\uCF^-(\gs,\ds_3) 
	\ar[dr, "\id"]	
&
\\
&
\CF^-(\gs,\ds_2)
	\ar[rr, "f_{\g}^{\dt_2\to \dt_3}"]
	\ar[from =uu, crossing over, "f_{\g}^{\tilde{\dt}'\to \dt_2}"]
	&&
\uCF^-(\gs,\ds_3) 
	\ar[from=uull,dashed, crossing over, "H_{\g}^{\tilde{\dt}'\to \dt_3}",sloped]
\end{tikzcd} 
\label{eq:hypercube-Pi-can}
\end{equation}

\begin{figure}[h]
	\centering
	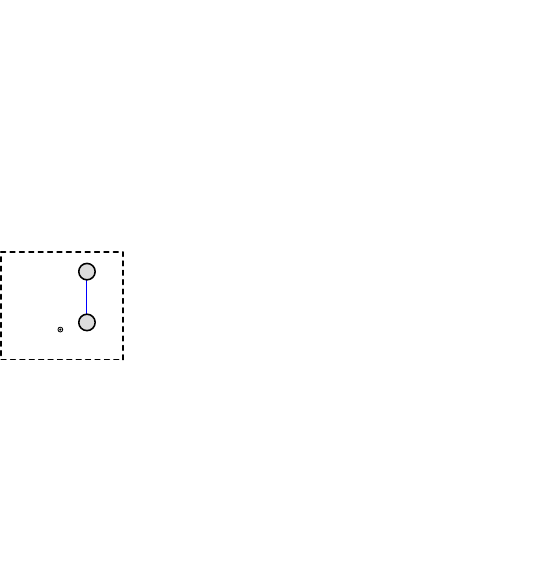
	\caption{The 3-handle hypercube $\cC^{(3)}_{\cen,\thrh}$. }\label{fig:181}
\end{figure}

Next, we construct the additional hypercube $\cC^{(3)}_{\cen,\thrh}$ shown in Figure~\ref{fig:181}. Note that $\cC^{(3)}_{\cen,\thrh}$ is a hypercube for destabilization, as constructed in Section~\ref{sec:hypercube-for-stabilization}.

Finally, we define $\cC_{\cen}^{(3)}$ by stacking and then compressing $\cC_{\aux,\Pi^{\can}}^{(3)}$ and $\cC_{\cen,\thrh}^{(3)}$.

\subsection{Constructing \texorpdfstring{$\cC_{\cen}$}{C-cen}}
\label{sec:constructing-C-cen}
We now stack and compress $\cC_{\cen}^{(1)}$, $\cC_{\cen}^{(2)}$ and $\cC_{\cen}^{(3)}$. Let us write $\hat{\cC}_{\cen}$ for the hypercube obtained by stacking $\cC_{\cen}^{(1)}$, $\cC_{\cen}^{(2)}$ and $\cC_{\cen}^{(3)}$, and then compressing. The hypercube $\hat{\cC}_{\cen}$ has the form shown in Figure~\ref{fig:C-cen-hat}.

\begin{figure}[h]
\[
\begin{tikzcd}[column sep=.5cm, row sep=.4cm]
\bCF^-(\as_1\cup \bar{\bs},\Ds)
	\arrow[dd]
	\arrow[dddrrr,dotted]
	\arrow[dr]
	\arrow[rr,"\id"]
&&[-1.5cm]
\bCF^-(\as_1\cup \bar{\bs},\Ds)
	\arrow[dddr,dashed, "\hat{E}"]
&
\\
&\bCF^-(\as_2\cup \bar{\bs}, \Ds)
	\arrow[ddrr,dashed]
	\arrow[rr,crossing over]
&&
\buCF^-(\as_3\cup \bar{\bs}, \Ds)
	\arrow[dd]
	\arrow[from=ulll,dashed,crossing over]
\\[2.5cm]
\bCF^-(\as\cup \bar{\bs}, \Ds_1)
	\arrow[drrr,dashed]
	\arrow[dr]
	&&\,&
	\\
&
\bCF^-(\as\cup \bar{\bs}, \Ds_2)
	\arrow[rr]
	\arrow[from =uu, crossing over]
	\arrow[from=uuul,dashed,crossing over]
	&&
\buCF^-(\as\cup \bar{\bs}, \Ds_3)
\end{tikzcd}
\]
\caption{The hypercube $\hat{\cC}_{\cen}$.}
\label{fig:C-cen-hat}
\end{figure}

In the cube $\hat{\cC}_{\cen}$, we would like to delete the extra complex labeled $\bCF^-(\as_1\cup \bar{\bs}, \Ds)$ along the top face. If we remove this complex, the length 0, 1 and 2 relations are preserved, though the length 3 relation may not be. However if $\hat{E}$ is itself null-homotopic, i.e. $\hat{E}=\d H+H\d,$
for some $H$, then $H$ may be added to the length 3 arrow of $\hat{\cC}_{\cen}$, and the hypercube relations will be satisfied. By examination of the compression operation, the map $\hat{E}$ is given by
\begin{equation}
\begin{split}
\hat{E}&=F_3 f_{\g}^{\dt''\to \dt_3} \sigma_0 F_3 f_{\g'}^{\dt'\to \dt''} \sigma F_3 h_1 f_{\g}^{\dt\to \dt'}\\
&+F_3 f_{\g}^{\dt''\to \dt_3}\sigma_0F_3 h_2f_{\g}^{\dt\to \dt'}\\
&+F_3 h_3 f_{\g}^{\dt\to \dt'}.
\end{split}
\label{eq:expanded-E-big}
\end{equation}
Equation~\eqref{eq:expanded-E-big} can be simplified as follows. Note first that~\eqref{eq:C-1-retract} is equivalent to 
\begin{equation}
F_3 f_{\g'}^{\dt'\to \dt''}\sigma=\id,
\label{eq:canceling-handles-E-big}
\end{equation} since $\sigma$ is a chain isomorphism.  Similarly~\eqref{eq:left-inverse} implies 
\begin{equation}
F_3 f_{\g}^{\dt''\to \dt_3} \sigma_0=\id.
\label{eq:canceling-handles-E-big-2}
\end{equation}
Hence ~\eqref{eq:expanded-E-big} simplifies to
\begin{equation}
\hat{E}=\left(F_3 h_1 +F_3 h_2+F_3 h_3\right)f_{\g}^{\dt\to \dt'}.
\end{equation}

It turns out that the most natural object of study is the map 
\[
E:=F_3 h_1+F_3 h_2+F_3 h_3.
\]
In Section~\ref{sec:N-null-homotopy}, we will construct a map $H_E$ satisfying
\[
E=\d H_E+H_E \d,
\]
which, of course, also gives a null-homotopy $H=H_E f_{\g}^{\dt\to \dt'}$ of $\hat{E}$. We define $\cC_{\cen}$ to be the hypercube obtained by deleting the extra copy of $\bCF^-(\as_1\cup \bar{\bs}, \Ds)$ from $\hat{\cC}_{\cen}$, and adding $H$ to the length 3 map of $\hat{\cC}_{\cen}$.

\subsection{Compatibility with the involution}

 In this section, we verify the following:

\begin{lem}\label{lem:columns=change-of-diagrams}
 The vertical maps in the central hypercube coincide with the change of diagram maps from naturality.
\end{lem}
\begin{proof}Let us write $\Phi_1$, $\Phi_2$ and $\Phi_3$ for the vertical maps (in the same order as the central hypercube). By construction,
\begin{equation}
\begin{split}
\Phi_1&=F_3f_{\g}^{\dt\to \dt'}\\
\Phi_2&=F_3 f_{\g}^{\tilde{\dt}'\to \dt_2}\sigma_0 F_3 f_{\g'}^{\dt\to \tilde{\dt}'}\\
\Phi_3&=F_3f_{\g}^{\dt''\to \dt_3}\sigma_0 F_3 f_{\g'}^{\dt'\to \dt''}\sigma F_3 f_{\g''}^{\dt\to \dt'}
\end{split}
\label{eq:defs-column-maps}
\end{equation}

We first consider $\Phi_1$. Note that $f_{\g}^{\dt\to \dt'}$ represents the cobordism map for a 0-framed unknot, so the formula for $\Phi_1$ consists of the cobordism map for surgery on a 0-framed unknot, followed by a canceling 3-handle. Hence, $\Phi_1$ represents the map from naturality.

Next, we consider $\Phi_2$. Firstly, we claim that $F_3 f_{\g}^{\tilde{\dt}'\to \dt_2} \sigma_0=\id$, which is proven by using Lemma~\ref{lem:lens-stabilization-triangles} and Proposition~\ref{prop:nearest-point-triangles}, similar to equation~\eqref{eq:left-inverse}. Hence the formula reduces to
\[
\Phi_2=F_3 f_{\g'}^{\dt\to \tilde{\dt}'}.
\]
As with $\Phi_1$, this expression for $\Phi_2$ represents the cobordism map for surgery on a 0-framed unknot, followed by the 3-handle map. Hence $\Phi_2$ coincides with the map from naturality.

Finally, we consider $\Phi_3$. Using ~\eqref{eq:canceling-handles-E-big} and~\eqref{eq:canceling-handles-E-big-2}, the expression for $\Phi_3$ in ~\eqref{eq:defs-column-maps} reduces to
\[
\Phi_3=F_3 f_{\g''}^{\dt\to \dt'},
\]
which similarly represents the cobordism map for a canceling pair of index 2 and 3 handles. The proof is complete.
\end{proof}

\section{The second 2-handle hypercube}
\label{sec:hypercube-C-3}

We now construct the second 2-handle hypercube, $\cC_{\twoh}^2$. We first describe the Heegaard diagrams that appear therein, and then subsequently define the maps, and prove the hypercube relations.

\subsection{Heegaard diagrams for \texorpdfstring{$\cC_{\twoh}^2$}{C-2-h-2}}

All of the complexes used to construct $\cC_{\twoh}^2$ use the Heegaard surface
\[
\Sigma_0\# \bar{\Sigma},
\]
where $\Sigma_0$ denotes the surface obtained by removing the genus 1 surgery region from $\Sigma$. We define the following sets of attaching curves:
\begin{itemize}
\item $\as\subset \Sigma_0$ denotes the original alpha curves which are outside of the genus 1 surgery region of $\Sigma$.
\item $\widehat{\Ds}\subset \Sigma_0\# \bar{\Sigma}_0$ denotes the curves obtained by taking the curves $\Ds$ constructed in Section~\ref{sec:peripheral-hypercubes}, and deleting the two curves which intersect the surgery region of $\Sigma$.
\item We define
\[
\Ds_i:=\widehat{\Ds}\cup \bar{c}_i,
\]
where $\bar{c}_i$ denotes the image of $c_i$ on $\bar{\Sigma}$.
\item $\bar{\as}^H$ denotes the curves obtained by handlesliding each curve in $\bar{\as}$ over the corresponding curve in $\as$. 
Furthermore, we choose each curve in $\bar{\as}^H$ to be a small translate of a curve in $\widehat{\Ds}$.
\item We define 
\[
\bar{\as}_i^H:=\bar{\as}^H\cup \bar{c}_i.
\]
\end{itemize}
See Figure~\ref{fig:116} for a schematic of the curves $\Ds_i$, $\as$ and $\bar{\as}_i^H$ outside of the genus 1 surgery region of $\Sigma_0\# \bar{\Sigma}$. 

\begin{lem}\label{lem:admissibility-C-2h-2}
 If the $\bar{\bs}$ curves are sufficiently wound, then the tuple
\[
(\Sigma_0\# \bar{\Sigma}, \as\cup \bar{\bs}, \Ds_1,\Ds_2,\Ds_3,\as'\cup \bar{\as}_1^H, \as''\cup \bar{\as}_2^H, \as'''\cup \bar{\as}_3^H,w)
\] 
is weakly admissible, where $\as'$, $\as''$ and $\as'''$ are appropriately chosen,  Hamiltonian translates of $\as$.
\end{lem}
 The proof of Lemma~\ref{lem:admissibility-C-2h-2} is an easy adaptation of the proof of Lemma~\ref{lem:admissibility-C-2-h-1}.

\subsection{Constructing \texorpdfstring{$\cC_{\twoh}^2$}{C-2-h-2}}

We first construct a diagram $\cL_{\twoh,\b}^2$ of beta attaching curves:
\[
\cL_{\twoh,\b}^{2}:=\begin{tikzcd}[column sep=2cm]
\Ds_1\arrow[dd, "\Theta_{\Dt_1,\a\cup \bar{\a}_1^H}^+"]
	\arrow[dr, "\Theta_{\Dt_1,\Dt_2}^{\can}",sloped]
	&&\\
&\Ds_2
	\arrow[r, "\Theta_{\Dt_2,\Dt_3}^+"]
	\arrow[dd, "\Theta_{\Dt_2,\a\cup \bar{\a}_2^H}^+"]
&\Ds_3
	\arrow[dd, "\Theta_{\Dt_3,\a\cup \bar{\a}_3^H}^+"]
	\\
\as\cup \bar{\as}_1^H
	\arrow[dr, "\Theta_{\a\cup \bar{\a}_1^H,\a\cup \bar{\a}_2^H}^{\can}",sloped]&&\\
&\as\cup \bar{\as}_2^H
	\arrow[r, "\Theta_{\a\cup \bar{\a}_2^H,\a\cup \bar{\a}_3^H}^+"]&
 \as\cup \bar{\as}_3^H
\end{tikzcd}.
\]
Here, we are omitting writing primes on the various curves labeled $\as$. 

\begin{figure}[H]
	\centering
\begingroup%
  \makeatletter%
  \providecommand\color[2][]{%
    \errmessage{(Inkscape) Color is used for the text in Inkscape, but the package 'color.sty' is not loaded}%
    \renewcommand\color[2][]{}%
  }%
  \providecommand\transparent[1]{%
    \errmessage{(Inkscape) Transparency is used (non-zero) for the text in Inkscape, but the package 'transparent.sty' is not loaded}%
    \renewcommand\transparent[1]{}%
  }%
  \providecommand\rotatebox[2]{#2}%
  \newcommand*\fsize{\dimexpr\f@size pt\relax}%
  \newcommand*\lineheight[1]{\fontsize{\fsize}{#1\fsize}\selectfont}%
  \ifx\svgwidth\undefined%
    \setlength{\unitlength}{349.66588699bp}%
    \ifx\svgscale\undefined%
      \relax%
    \else%
      \setlength{\unitlength}{\unitlength * \real{\svgscale}}%
    \fi%
  \else%
    \setlength{\unitlength}{\svgwidth}%
  \fi%
  \global\let\svgwidth\undefined%
  \global\let\svgscale\undefined%
  \makeatother%
  \begin{picture}(1,0.47046373)%
    \lineheight{1}%
    \setlength\tabcolsep{0pt}%
    \put(0,0){\includegraphics[width=\unitlength,page=1]{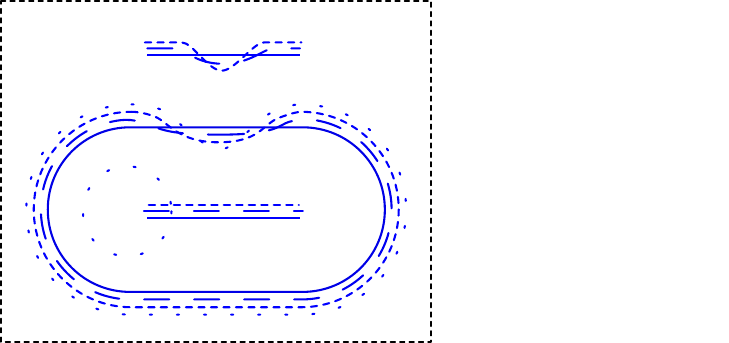}}%
    \put(0.75007522,0.11743955){\color[rgb]{0,0,0}\makebox(0,0)[t]{\lineheight{1.25}\smash{\begin{tabular}[t]{c}$\as\cup \bar{\as}^H_3$\end{tabular}}}}%
    \put(0.75002125,0.31228004){\color[rgb]{0,0,0}\makebox(0,0)[t]{\lineheight{1.25}\smash{\begin{tabular}[t]{c}$\Ds_1$\end{tabular}}}}%
    \put(0.75002125,0.24857853){\color[rgb]{0,0,0}\makebox(0,0)[t]{\lineheight{1.25}\smash{\begin{tabular}[t]{c}$\Ds_2$\end{tabular}}}}%
    \put(0.75002125,0.18217915){\color[rgb]{0,0,0}\makebox(0,0)[t]{\lineheight{1.25}\smash{\begin{tabular}[t]{c}$\Ds_3$\end{tabular}}}}%
    \put(0,0){\includegraphics[width=\unitlength,page=2]{fig172.pdf}}%
  \end{picture}%
\endgroup%

	\caption{A generic genus 2 region of the quadruple $(\Sigma_0\# \bar{\Sigma}, \Ds_1,\Ds_2,\Ds_3,\as\cup \bar{\as}^H_3)$.}\label{fig:116}
\end{figure}

A generic, genus 2 region of the quadruple $(\Sigma_0\# \bar{\Sigma}, \Ds_1,\Ds_2,\Ds_3,\as\cup \bar{\as}_3^H)$ is shown in Figure~\ref{fig:116}.

\begin{lem}
 The diagram $\cL_{\twoh,\b}^{2}$ is a hypercube of attaching curves.
\end{lem}
\begin{proof} The length 1 relations are equivalent to each length 1 class being a cycle, which is straightforward from the construction.

There are four length 2 relations to verify, which follow from the relations
\begin{equation}
\begin{split} 
f_{\Dt_1,\Dt_2,\Dt_3}\left(\Theta_{\Dt_1,\Dt_2}^{\can}, \Theta_{\Dt_2,\Dt_3}^+\right)&=0,\\
f_{\a\cup \bar{\a}_1^H, \a\cup \bar{\a}_2^H,\a\cup \bar{\a}_3^H}\left(\Theta_{\a\cup \bar{\a}_1^H, \a\cup \bar{\a}_2^H}^{\can}, \Theta_{\a\cup \bar{\a}_2^H,\a\cup \bar{\a}_3^H}^+\right)&=0,\\
f_{\Dt_1,\a\cup \bar{\a}_1^H, \a\cup \bar{\a}_2^H}\left(\Theta_{\Dt_1,\a\cup \bar{\a}_1^H}^+, \Theta_{\a\cup \bar{\a}_1^H, \a\cup \bar{\a}_2^H}^{\can}\right)&= \Theta_{\Dt_1,\a\cup \bar{\a}_2^H}^{\can},\\
f_{\Dt_1,\a\cup \bar{\a}_1^H,\a\cup \bar{\a}_2^H}\left(\Theta_{\Dt_1,\a\cup \bar{\a}_1^H}^+, \Theta_{\a\cup \bar{\a}_1^H, \a\cup \bar{\a}_2^H}^{\can}\right)&=\Theta_{\Dt_1,\a\cup \bar{\a}_2^H}^{\can},\\
f_{\Dt_2,\Dt_3,\a\cup \bar{\a}_3^H}\left(\Theta_{\Dt_2,\Dt_3}^+, \Theta_{\Dt_3,\a\cup \bar{\a}_3^H}^+\right)&=\Theta_{\Dt_2,\a\cup \bar{\a}_3^H}^+,\\
f_{\Dt_2,\a\cup \bar{\a}_2^H, \a\cup \bar{\a}_3^H}\left(\Theta_{\Dt_2,\a\cup \bar{\a}_2^H}^+, \Theta_{\a\cup \bar{\a}_2^H, \a\cup \bar{\a}_3^H}^+\right)&=\Theta_{\Dt_2,\a\cup \bar{\a}_3^H}^+.
\end{split}
\label{eq:C-2h-1-length-2-model-computations}
\end{equation}
The first two lines of ~\eqref{eq:C-2h-1-length-2-model-computations} follow from Lemma~\ref{lem:omega-a-null-homotopy}. We now focus on the third equation. The triple $(\Sigma_0\# \bar{\Sigma}, \Ds_1, \as\cup \bar{\as}_1, \as\cup \bar{\as}_2^H)$ is a genus $2g(\Sigma)-2$ stabilization of the genus 1 triple $(\bT^2, \bar{c}'_1, c_1,c_2)$ (where $\bar{c}_1'$ denotes an isotopic copy of $\bar{c}_1$). The index tuple of this stabilization is $(0,-g(\Sigma_0),0)$, by Lemma~\ref{lem:monotonicity-index-tuple} and Figure~\ref{fig:117}. By Proposition~\ref{prop:stabilize-triangles-general}, it suffices to show the claim for the destabilized, genus 1 triple. For the genus 1 triple, the count is immediate. The fourth through sixth model computations are proven similarly.

The length 3 relation fo the diagram $\cL_{\twoh,\b}^{2}$ is reads
\begin{equation}
\begin{split}
0&=h_{\Dt_1,\Dt_2,\Dt_3,\a\cup \bar{\a}_3^H}\left(\Theta^{\can}_{\Dt_1,\Dt_2},\Theta^+_{\Dt_2,\Dt_3},\Theta^+_{\Dt_3,\a\cup \bar{\a}_3^H}\right)\\
&+h_{\Dt_1,\Dt_2,\a\cup \bar{\a}_2^H, \a\cup \bar{\a}_3^H}\left(\Theta^{\can}_{\Dt_1,\Dt_2},\Theta^+_{\Dt_2,\a\cup \bar{\a}^H_2},\Theta^+_{\a\cup \bar{\a}^H_2,\a\cup \bar{\a}^H_3}\right)\\
&+h_{\Dt_1,\a\cup \bar{\a}_1^H, \a\cup \bar{\a}_2^H, \a\cup \bar{\a}_3^H}\left(\Theta^+_{\Dt_1,\a\cup \bar{\a}^H_1},\Theta^{\can}_{\a\cup \bar{\a}^H_1,\a\cup \bar{\a}^H_2},\Theta^+_{\a\cup \bar{\a}^H_2,\a\cup \bar{\a}^H_3}\right).
\end{split}
\label{eq:C-2g-1-length-3-model-computations}
\end{equation}
To prove ~\eqref{eq:C-2g-1-length-3-model-computations}, we will show that each quadrilateral count vanishes. We focus on the count which occurs on $(\Sigma_0\# \bar{\Sigma},\Ds_1,\Ds_2,\Ds_3,\as\cup \bar{\as}_3^H)$. This quadruple is a genus $2g(\Sigma)-2$ stabilization of the genus 1 quadruple $(\bT^2,c_1,c_2,c_3,c_3')$. Furthermore, the index tuple of this stabilization is $(-g(\Sigma_0),-g(\Sigma_0),0,0)$, by Lemma~\ref{lem:monotonicity-index-tuple} and Figure~\ref{fig:118}. By Proposition~\ref{prop:stabilize-quadrilateral} it is sufficient to show that the associated count of quadrilaterals on $(\bT^2,c_1,c_2,c_3,c_3')$ vanishes, which is proven in Lemma~\ref{lem:genus-1-quadruple-count}. The other two quadrilateral counts are proven by an easy modification of this arguement.
\end{proof}

We define $\cC_{\twoh}^{2}$ as the pairing of the $\cL_{\twoh,\beta}^{2}$, and the 0-dimensional hypercube of alpha attaching curves consisting of $\as\cup \bar{\bs}$. The hypercube is schematically shown below:

\begin{figure}[H]
\[\begin{tikzcd}[column sep=1.6cm, row sep=1.3cm]
\bCF^-( \as\cup \bar{\bs},\Ds_1)
\arrow[dr,sloped,swap, "f_{\a\cup \bar{\b}}^{\Delta_1 \to \a\cup \bar{\a}_2^H}"]
\arrow[dd, swap, "f_{\a\cup \bar{\b}}^{\Delta_1\to  \a\cup \bar{\a}^H_1}"]
\arrow[drr,dashed,"h_{\a\cup \bar{\b}}^{ \Delta_1 \to \Dt_2\to  \Dt_3}", sloped]
\arrow[dddrr,dotted,sloped,swap, "\hspace{3.2cm}P_{\a\cup \bar{\b}}^{\Delta_1\to  \a\cup \bar{\a}_3^H}"]&&
\\
& \bCF^-(\as\cup \bar{\bs}, \Ds_2)
\arrow[ddr, dashed, sloped, "H_{\a\cup \bar{\b}}^{\Dt_2\to \a\cup \bar{\a}_3^H}"]
\arrow[r, "f_{\a\cup \bar{\b}}^{\a\cup \bar{\a}^H_2\to  \a\cup \bar{\a}^H_3}",swap]
& \buCF^-(\as\cup \bar{\bs},\Ds_3)
\arrow[dd, "f_{\a\cup \bar{\b}}^{\Dt_3\to \a\cup \bar{\a}_3}"]
\\
\bCF^-(\as\cup \bar{\bs}, \as\cup \bar{\as}^H_1)
\arrow[dr, swap, sloped, "f_{\a\cup \bar{\b}}^{\a\cup\bar{\a}^H_1\to \a\cup \bar{\a}^H_2}"]
\arrow[drr,dashed, sloped, near end, "h_{\a\cup \bar{\b}}^{\a\cup \bar{\a}^H_1\to \a\cup \bar{\a}_2^H\to \a\cup \bar{\a}^H_3}"] &
 \\
& \bCF^-(\as\cup \bar{\bs}, \as\cup \bar{\as}^H_2)
 \arrow[r, "f_{\a\cup \bar{\b}}^{\a\cup \bar{\a}_2^H\to  \a\cup \bar{\a}^H_3}", swap]
 \arrow[from=uu,crossing over, "f_{\a\cup \bar{\b}}^{\Dt_2\to \a\cup \bar{\a}_2^H}",swap]
 \arrow[from=uuul,crossing over,dashed, "H_{\a\cup \bar{\b}}^{\Dt_1\to \a\cup \bar{\a}^H_2}", sloped]
& \buCF^-(\as\cup \bar{\bs}, \as\cup \bar{\as}^H_3)
\end{tikzcd}
\]
\caption{The hypercube $\cC_{\twoh}^{2}$. We omit writing the primes on the various copies of $\as$.}
\label{fig:C2-h-2}
\end{figure}

\section{The tautological and flip-map  hypercubes}
\label{sec:flip-hypercube}

\subsection{The tautological hypercube}
\label{sec:tautological-hypercube}
The bottom face of the hypercube $\cC_{\thrh}$ reads
\begin{equation}
\begin{tikzcd}[row sep=1cm]
\bCF^-\left(\bar{\Sigma}, \bar{\bs}, \bar{\as}_1\right)
	\arrow[d, "f_{\bar{\b}}^{\bar{\a}_1\to \bar{\a}_2}",swap]
	\arrow[dr,dashed, "h_{\bar{\b}}^{\bar{\a}_1\to \bar{\a}_2\to \bar{\a}_3}"]&\\
\bCF^-(\bar{\Sigma}, \bar{\bs}, \bar{\as}_2)
	\arrow[r,"f_{\bar{\b}}^{\bar{\a}_2\to \bar{\a}_3}"]& 
\buCF^-(\bar{\Sigma}, \bar{\bs}, \bar{\as}_3,\hat{z},w).
\end{tikzcd}
\label{eq:almost-mirror-hypercube}
\end{equation}
The diagram in ~\eqref{eq:almost-mirror-hypercube} is very similar to the top of the main hypercube, as described in Section~\ref{sec:previous-results}. There are two subtle differences:
\begin{enumerate}
\item The point $z$ is now the special basepoint of the third complex, whereas at the top of the main hypercube, $w$ is the special basepoint
\item To achieve admissibility for triangles and rectangles, the curves $\bar{c}_1$, $\bar{c}_2$ and $\bar{c}_3$ must have slightly different position than $c_1$, $c_2$ and $c_3$. (Compare the two sides of Figure~\ref{fig:186}). 
\end{enumerate}
We will assume that $\bar{c}_1$ and $\bar{c}_2$ are the mirrors of $c_1$ and $c_2$, but that $\bar{c}_3$ differs slightly from the mirror of $c_3$. We write $\bar{\bar{c}}_3$ for the image of $\bar{c}_3$ on $\Sigma$. We may view $\bar{\bar{c}}_3$ as being obtained by sliding $c_3$ across the intersection of $c_1$ and $c_2$. Compare Figure~\ref{fig:186}. We will write
\[
\bar{\bar{\as}}_3:=\as\cup \bar{\bar{c}}_3.
\]

As defined in Section~\ref{sec:involutive-knots}, there are tautological isomorphism
\[
\begin{split}
\eta_1&\colon \bCF^-(\bar{\Sigma}, \bar{\bs}, \bar{\as}_1)\to \bCF^-(\Sigma,\as_1,\bs)\\
\eta_2&\colon \bCF^-(\bar{\Sigma}, \bar{\bs}, \bar{\as}_2)\to \bCF^-(\Sigma,\as_2,\bs)\\
\underline{\eta}_3& \colon \buCF^-(\bar{\Sigma}, \bar{\bs}, \bar{\as}_3, \hat{z}, w)\to \buCF^-(\Sigma,\bar{\bar{\as}}_3,\bs,w,\hat{z}).
\end{split}
\]
Furthermore, the following diagram is tautologically a hypercube of chain complexes:
\begin{figure}[H]
\[
\begin{tikzcd}[row sep=.8cm, column sep=1.5 cm]
\bCF^-(\bar{\Sigma}, \bar{\bs}, \bar{\as}_1)
	\arrow[dr, "f^{\bar{\a}_1\to \bar{\a}_2}_{\bar{\b}}", sloped, swap]
	\arrow[ddd,"\eta_1"]
	\arrow[rrd,dashed, "h_{\bar{\b}}^{\bar{\a}_1\to \bar{\a}_2\to \bar{\a}_3}",sloped, near start] && \\
& \bCF^-(\bar{\Sigma}, \bar{\bs}, \bar{\as}_2)
	\arrow[r,swap, "f_{\bar{\b}}^{\bar{\a}_2\to \bar{\a}_3}"]
&\buCF^-(\bar{\Sigma}, \bar{\bs}, \bar{\as}_3, \hat{z}, w)
	\arrow[ddd,"\underline{\eta}_3"]\\
\\
\bCF^-(\Sigma,\as_1,\bs)
	\arrow[drr,dashed,sloped,"h^{\b}_{\a_1\to \a_2\to \bar{\bar{\a}}_3}",near start]
	\arrow[dr, "f_{\a_1\to \a_2}^{\b}",sloped, swap]&& \\
& \bCF^-(\Sigma,\as_2,\bs)
	\arrow[from=uuu, crossing over,swap, "\eta_2"]
	\arrow[r,swap,"f_{\a_2\to \bar{\bar{\a}}_3}^{\b}"]
&\buCF^-(\Sigma,\bar{\bar{\as}}_3,\bs,w,\hat{z})
\end{tikzcd}
\]
\caption{The tautological hypercube $\cC_{\eta}$.}
\label{fig:hyperbox-taut}
\end{figure}

\subsection{The flip-map}

In this section, we investigate the \emph{flip-map}
\[
\frF_{z\to w}\colon \buCF^-(\Sigma,\bar{\bar{\as}}_3,\bs,w,\hat{z})\to \buCF^-(\Sigma,\as_3,\bs,\hat{w},z),
\]
in detail. Note that the distinction between $\bar{\bar{\as}}_3$ and $\as_3$ makes no difference to the complexes described above (if $\bar{\bar{c}}_3$ is sufficiently close to $c_3$), however it will make a difference for the triangle and quadrilateral maps appearing in the flip-map hypercube. (Compare~Figure~\ref{fig:186}).

We recall from ~\eqref{eq:def-iota-3} that the flip-map  $\frF_{z\to w}$ is the composition
\[
\frF_{z\to w}:=\theta_w^{-1}\circ (\Psi_{(\Sigma,\bar{\bar{\a}}_3,\b,z)\to (\Sigma,\a_3,\b,w)} \otimes T^n)\circ\theta_z,
\]
where $\Psi_{(\Sigma,\bar{\bar{\a}}_3,\b,z)\to (\Sigma,\a_3,\b,w)}$ is the naturality map induced by moving $z$ to $w$ along a fixed subarc of $K$.

For the purpose of constructing the flip map hypercube, it is helpful to have a concrete model for the naturality map $\Psi_{(\Sigma,\bar{\bar{\a}}_3,\b,z)\to (\Sigma,\a_3,\b,w)}$.  Let us write $(\Sigma,\as_3',\bs,z)$ for the diagram obtained by moving $c_3$ over the basepoint $z$ (this may be achieved by a sequence of handleslides, away from $z$). There is a canonical diffeomorphism
\[
\phi_*\colon (\Sigma,\as_3',\bs,z)\to (\Sigma,\as_3,\bs,w),
\]
which may be identified with the diffeomorphism map which moves $z$ to $w$, along a subarc of $K$.

The map $\Psi_{(\Sigma,\bar{\bar{\a}}_3,\b,z)\to (\Sigma,\a_3,\b,w)}$ is equal to the composition
\[
\Psi_{(\Sigma,\bar{\bar{\a}}_3,\b,z)\to (\Sigma,\a_3,\b,w)}=\phi_*\circ  \Psi_{\bar{\bar{\a}}_3\to \a_3'}^{\b},
\]
where $\Psi_{\bar{\bar{\a}}_3\to \a_3'}^{\b}\colon \bCF^-(\Sigma,\as_3,\bs,z)\to \bCF^-(\Sigma,\as_3',\bs,z)$ is the ordinary map from naturality.

\begin{figure}[H]
	\centering
	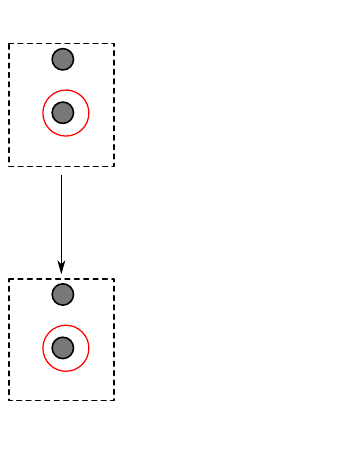
	\caption{Decomposing $\Psi_{\a_3\to \a_3'}^{\b}$ as the composition $f_{\a^\twind_3\to \a_3'}^{\b}\circ f_{\a_3\to \a^\twind_3}^{\b}$. We continue our convention of denoting wound curves with a $\twind$.}\label{fig:73}
\end{figure}

The map $\Psi_{\bar{\bar{\a}}_3\to \a_3'}^{\b}$ cannot in general be computed using a single triangle count, because of admissibility issues. Instead, we pick an intermediate choice of alpha-curves, $\as^{\twind}_3$, by winding $\as_3'$ in the genus 1 surgery region along a curve which is dual to $c_3$. We decompose $\Psi_{\bar{\bar{\a}}_3\to \a_3'}^{\b}$ as 
\[
\Psi_{\bar{\bar{\a}}_3\to \a_3'}^{\b}=f_{\a^{\twind}_3\to \a_3'}^{\b}\circ f_{\bar{\bar{\a}}_3\to \a^{\twind}_3}^{\b},
\]
where $f_{\a^{\twind}_3\to \a_3'}^{\b}$ and $f_{\bar{\bar{\a}}_3\to \a^{\twind}_3}^{\b}$ are holomorphic triangle maps. See Figure~\ref{fig:73}.  Compare \cite{MOIntegerSurgery}*{Section~8.7}. 

\subsection{More hypercubes with twisted coefficients}
\label{sec:twisted-hypercubes-2}

In this section, we describe a construction of hypercubes with twisted coefficients, which is similar to the one from Section~\ref{sec:twisted-hypercubes-1}, but appears in the flip-map hypercube.

Let  $\Sigma$ be a compact, oriented surface, equipped with two points, $p,q\in \Sigma$, as well as a special basepoint $w$ (which may be different from $p$ and $q$). Suppose that $\as_n,\dots, \as_1$ is a collection of attaching curves on $\Sigma$, and suppose that $\g_n,\dots, \g_1$ is a collection of integral 1-chains in $U_{\a_n},\dots, U_{\a_1}$, respectively, such that $\d \g_i=p-q$. The 1-chains $\g_n,\dots, \g_1$ determine a relative 2-cycle
\[
S_{\g_n,\dots, \g_1}\in C_2(X_{\a_n,\dots, \a_1},\d X_{\a_n,\dots, \a_1}).
\]
The 2-chain $S_{\g_n,\dots, \g_1}$ is the sum of $\g_i\times [0,1]\subset U_{\a_i}\times [0,1]$, as well as $\{p,q\}\times D_n\subset \Sigma\times D_n$. If $\psi$ is a homology class of $n$-gons on $(\Sigma,\as_n,\dots, \as_1,w)$ and $\widehat{S}\in H_2(X;\Z)$, we define 
\[
\frm(w,\widehat{S},\psi)=\frac{\left\langle c_1(\frs_{w}(\psi)),[\widehat{S}]\right\rangle -[\widehat{S}]\cdot [\widehat{S}]}{2}.
\]
Note that $\frm(w,\widehat{S},\psi)$ is an integer, since $c_1(\frs_w(\psi))$ is a characteristic vector.

We say $\ul{\cL}_{\a}$ is a hypercube of alpha attaching curves with a \emph{cube of singular twisting chains} on $(\Sigma,p,q,w)$, if $\ul{\cL}_{\a}$ consists of the following:
\begin{enumerate}
\item For each $\veps\in \bE_n$, a set of curves $\as^{\veps}$ on $\Sigma$.
\item For each pair $\veps',\veps\in \bE_n$ such that $\veps'>\veps$, a choice of chain $\Theta_{\a^{\veps'},\a^\veps}\in \bCF^-(\Sigma,\as^{\veps'},\as^{\veps},w)\otimes \bF[\Z/m]$.
\item For each $\veps\in \bE_n$,  a choice of integral 1-chain $\g^\veps\subset U_{\a^{\veps}}$ such that $\d \g^{\veps}=p-q$.
\item For each pair $\veps,\veps'\in \bE_n$ such that $\veps'>\veps$, a choice of integral 2-chain $C_{\veps',\veps}\subset -U_{\a^{\veps'}}\cup U_{\a^{\veps}}$ such that $\d C_{\veps',\veps}=\g_{\veps'}\cup -\g_{\veps}$.
\end{enumerate}
Furthermore, there is a compatibility condition, as we now describe. Write $\ul{\frf}_{\a^{\veps_n},\dots, \a^{\veps_1}}$ for the  map which counts holomorphic $n$-gons of index $3-n$, weighted by $U^{n_w(\psi)} T^{\frm(w,\widehat{S}_{\g_n,\dots, \g_1},\psi)}$. Here, $\widehat{S}_{\g_n,\dots, \g_1}$ is the closed, integral 2-chain
\[
\widehat{S}_{\g_n,\dots, \g_1}=S_{\g_n,\dots,\g_1}-C_{\veps_n,\veps_{n-1}}-\dots-C_{\veps_2,\veps_1}+C_{\veps_n,\veps_1}.
\] 
We assume that the twisted polygon counting maps satisfy ~\eqref{eq:compatibility-hypercube-attaching-curves} (the compatibility condition for a hypercube of attaching curves). A hypercube of \emph{beta} attaching curves with a cube of singular twisting chains is defined by the natural modification.

Suppose $\ul{\cL}_{\a}$ and $\ul{\cL}_{\b}$ are two hypercubes of attaching curves on $(\Sigma,p,q,\hat{w})$ which are equipped with cubes of singular twisting chains. If we additionally pick  a collection of 2-chains $C_{\veps,\nu}$, ranging over $\veps\in \bE_n$ and $\nu\in \bE_{n'}$, we may pair them to form a hypercube $\ul{\buCF}^-(\ul{\cL}_{\a},\ul{\cL}_{\b})$. 
The complex at any $(\veps,\nu)$ is the Floer complex $\bCF^-(\as^{\veps}, \bs^\nu)\otimes \bF[\Z/m]$. The hypercube maps are defined  as in~\eqref{eq:pairing-attaching-curves}, but using the twisted polygon counting maps.

Suppose additionally that whenever $\veps'\in \bE_n$ has the property that the homology class of $\widehat{S}_{\veps_j,\dots, \veps_1}$ is trivial for all sequences $\veps_1< \cdots < \veps_j$ such that $\veps_j\le \veps'$,  the chain $\Theta_{\a^{\veps'},\a^\veps}$ is in $\bCF^-(\as^{\veps'},\as^{\veps})\otimes T^0$, for all $\veps<\veps'$.  We assume a similar condition holds for the beta-chains.
In this situation there is a subcomplex $\buCF^-(\ul{\cL}_{\a},\ul{\cL}_{\b})\subset \ul{\buCF}^-(\ul{\cL}_{\a},\ul{\cL}_{\b})$, as follows.   If $(\veps,\nu)\in \bE_n\times \bE_{n'}$ has the property that
\[
\widehat{S}_{\a^{\veps_i},\dots, \a^{\veps_1}, \b^{\nu_1},\dots, \b^{\nu_j}}=0\in H_2(X_{\a^{\veps_i},\dots, \a^{\veps_1}, \b^{\nu_1},\dots, \b^{\nu_j}}),
\]
whenever $\veps_1< \cdots < \veps_i \le \veps$ and $\nu_1< \cdots < \nu_j\le \nu$, then we set 
\[
\buCF^-(\ul{\cL}_{\a},\ul{\cL}_{\b})^{(\veps,\nu)}=\bCF^-(\as^{\veps}, \bs^{\nu}).
\]
Otherwise, we set $\buCF^-(\ul{\cL}_{\a},\ul{\cL}_{\b})^{(\veps,\nu)}=\bCF^-(\as^{\veps}, \bs^{\nu})\otimes \bF[\Z/m]$. It is straightforward to see that $\buCF^-(\ul{\cL}_{\a}, \ul{\cL}_{\b})^{(\veps,\nu)}$ is a subcomplex of $\ul{\buCF}^-(\ul{\cL}_{\a},\ul{\cL}_{\b})$.

For hypercubes of dimension 3 and lower, we use the following analog of the notation from Section~\ref{sec:low-dim-notation}.  We use lower case letters $\ul{\frf}$, $\ul{\frh}$ and $\ul{\frp}$ to indicate triangle, quadrilateral, and pentagon counting maps, which are weighted using by the quantity $U^{n_w(\psi)}T^{\frm(w,\widehat{S},\psi)}$. We use the upper case letters $\ul{\frH}$ and $\ul{\frP}$ to indicate length 2 and 3 maps in a hypercube, when they involve more than just a single count of quadrilaterals or pentagons.

\subsection{The flip-map hypercube}
\label{sec:flip-map-hypercube-details}
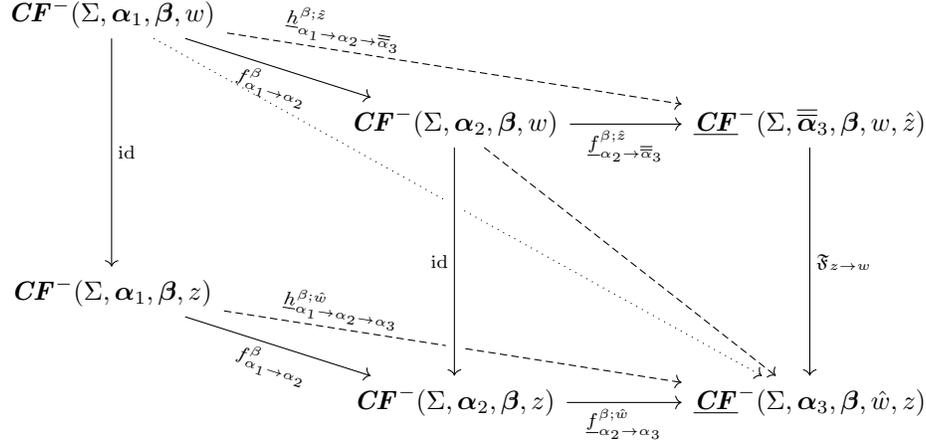
\begin{figure}[H]
\[
\begin{tikzcd}[row sep=.8cm, column sep=1.5 cm]
\bCF^-(\Sigma,\as_1,\bs,w)
	\arrow[dr, "f_{\a_1\to \a_2}^{\b}", sloped, swap]
	\arrow[ddddrr,dotted]
	\arrow[ddd,"\id"]
	\arrow[rrd,dashed, "\underline{h}^{\b;\hat{z}}_{\a_1\to \a_2\to \bar{\bar{\a}}_3}",sloped, near start] && \\
& \bCF^-(\Sigma, \as_2,\bs,w)
	\arrow[r,swap, "\ul{f}_{\a_2\to \bar{\bar{\a}}_3}^{\b;\hat{z}}"]
	\arrow[dddr,dashed]&
\buCF^-(\Sigma,\bar{\bar{\as}}_3,\bs,w,\hat{z})
	\arrow[ddd,"\frF_{z\to w}"]\\
\\
\bCF^-(\Sigma,\as_1,\bs,z)
	\arrow[drr,dashed,sloped,"\underline{h}^{\b;\hat{w}}_{\a_1\to \a_2\to \a_3}",near start]
	\arrow[dr, "f_{\a_1\to \a_2}^{\b}",sloped, swap]&& \\
& \bCF^-(\Sigma,\as_2,\bs,z)
	\arrow[from=uuu, crossing over,swap, "\id"]
	\arrow[r,swap,"\ul{f}_{\a_2\to \a_3}^{\b;\hat{w}}"]
&\buCF^-(\Sigma,\as_3,\bs,\hat{w},z)
\end{tikzcd}
\]
\caption{The flip-map hypercube $\cC_{\frF}$. A superscript containing $\hat{w}$ or $\hat{z}$ indicates which basepoint contributes the power of $U$.}
\label{fig:hyperbox-KFl}
\end{figure}

In this section, we construct the \emph{flip-map hypercube}, shown in Figure~\ref{fig:hyperbox-KFl}.

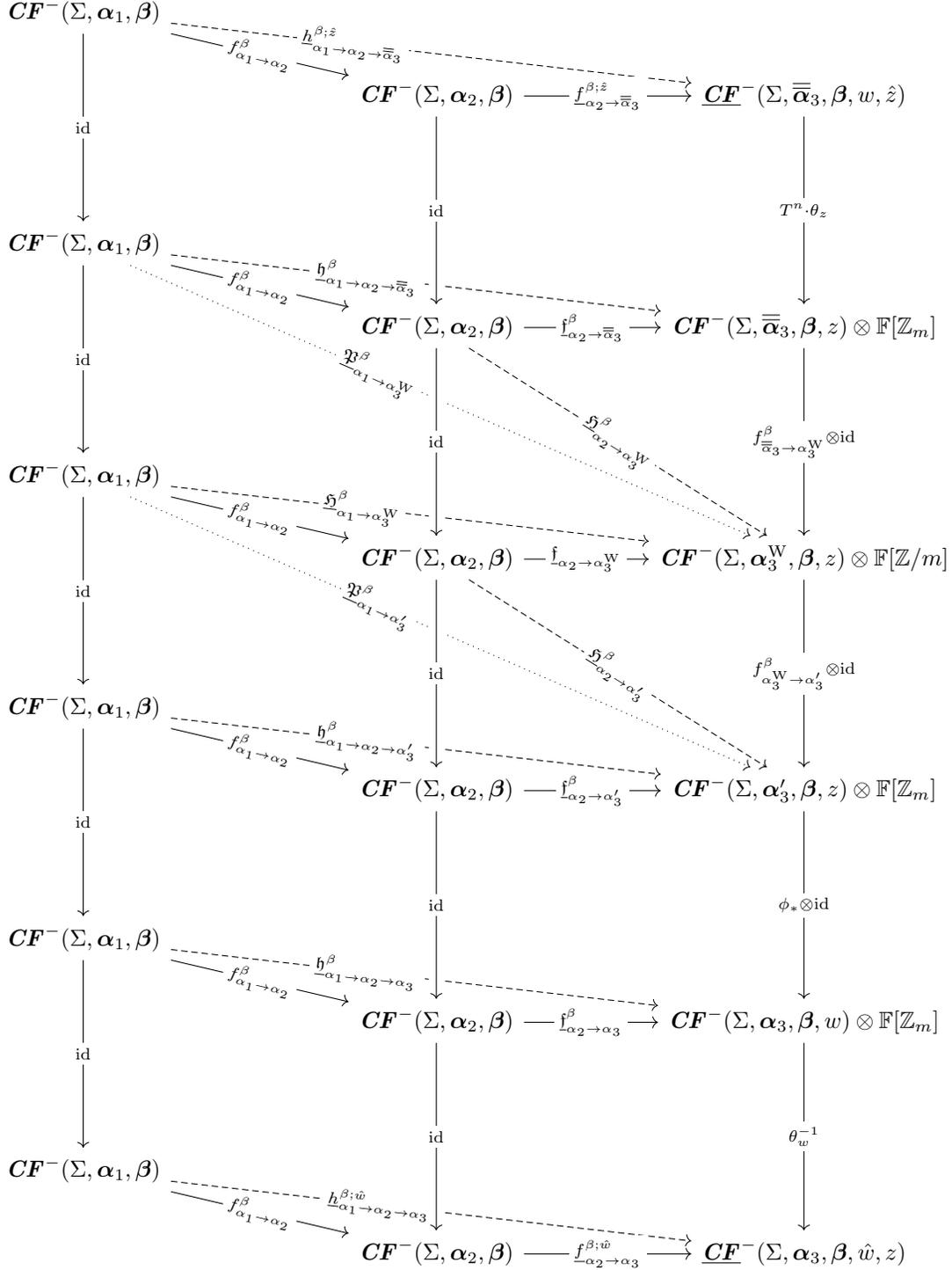
\begin{figure}[p]
\[
\begin{tikzcd}[row sep=.6cm, column sep=1.7 cm,labels=description]
\bCF^-(\Sigma,\as_1,\bs)
	\arrow[dr, "f_{\a_1\to \a_2}^{\b}", sloped, swap]
	\arrow[dd,"\id"]
	\arrow[rrd,dashed, "\underline{h}^{\b;\hat{z}}_{\a_1\to \a_2\to \bar{\bar{\a}}_3}",sloped, pos=.35]
	 &[1cm]
	 &[.2cm]
\\
& \bCF^-(\Sigma, \as_2,\bs)
	\arrow[r,swap, "\ul{f}_{\a_2\to \bar{\bar{\a}}_3}^{\b;\hat{z}}"]
&\buCF^-(\Sigma,\bar{\bar{\as}}_3,\bs,w,\hat{z})
	\arrow[dd, "T^n\cdot \theta_z"]
\\[1cm]
\bCF^-(\Sigma,\as_1,\bs)
	\arrow[dr, "f_{\a_1\to \a_2}^{\b}", sloped, swap]
	\arrow[dddrr,dotted, "\ul{\frP}_{\a_1\to \a_3^{\twind}}^{\b}",sloped, swap, pos=.4]
	\arrow[dd,"\id"]
	\arrow[rrd,dashed, 
		"\underline{\frh}^{\b}_{\a_1\to \a_2\to \bar{\bar{\a}}_3}",sloped, pos=.4] && 
\\
& 
\bCF^-(\Sigma, \as_2,\bs)
	\arrow[r,swap, "\ul{\frf}_{\a_2\to \bar{\bar{\a}}_3}^{\b}"]
	\arrow[ddr,dashed, "\ul{\frH}_{\a_2\to \a_3^{\twind}}^{\b}",sloped]
	\arrow[from=uu, crossing over,swap, "\id"]	
	&
\bCF^-(\Sigma,\bar{\bar{\as}}_3,\bs,z)\otimes \bF[\Z_m]
	\arrow[dd, "f_{\bar{\bar{\a}}_3\to \a_3^{\twind}}^{\b}\otimes \id"]
 \\[1cm]
\bCF^-(\Sigma,\as_1,\bs)
	\arrow[dr, "f_{\a_1\to \a_2}^{\b}", sloped, swap]
	\arrow[dddrr,dotted, sloped,swap, pos=.4, "\ul{\frP}_{\a_1\to \a_3'}^{\b}"]
	\arrow[dd,"\id"]
	\arrow[rrd,dashed, "\underline{\frH}^{\b}_{\a_1\to \a_3^{\twind}}",sloped, pos=.4] && 
\\
& \bCF^-(\Sigma, \as_2,\bs)
	\arrow[r,swap, "\ul{\frf}_{\a_2\to \a^{\twind}_3}"]
	\arrow[ddr,dashed, " \ul{\frH}_{\a_2\to \a_3'}^{\b}",sloped]
	\arrow[from=uu, crossing over,swap, "\id"]
&\bCF^-(\Sigma,\as_3^{\twind},\bs,z)\otimes \bF[\Z/m]
	\arrow[dd, "f_{\a_3^{\twind}\to \a_3'}^{\b}\otimes \id"]
\\[1cm]
\bCF^-(\Sigma,\as_1,\bs)
	\arrow[dr, "f_{\a_1\to \a_2}^{\b}", sloped, swap]
	\arrow[dd,"\id"]
	\arrow[rrd,dashed, "\underline{\frh}^{\b}_{\a_1\to \a_2\to \a'_3}",sloped,pos=.4] && 
\\
& \bCF^-(\Sigma, \as_2,\bs)
	\arrow[r,swap, "\ul{\frf}_{\a_2\to \a_3'}^{\b}"]
	\arrow[from=uu, crossing over,swap, "\id"]
&\bCF^-(\Sigma,\as_3',\bs,z)\otimes \bF[\Z_m]
	\arrow[dd, "\phi_*\otimes \id"]
\\[1cm]
\bCF^-(\Sigma,\as_1,\bs)
	\arrow[dr, "f_{\a_1\to \a_2}^{\b}", sloped, swap]
	\arrow[dd,"\id"]
	\arrow[rrd,dashed, "\underline{\frh}^{\b}_{\a_1\to \a_2\to \a_3}",sloped, pos=.4] && 
\\
& \bCF^-(\Sigma, \as_2,\bs)
	\arrow[r,swap, "\ul{\frf}_{\a_2\to \a_3}^{\b}"]
	\arrow[from=uu, crossing over,swap, "\id"]
&\bCF^-(\Sigma,\as_3,\bs,w)\otimes \bF[\Z_m]
	\arrow[dd, "\theta_w^{-1}"]
\\[1cm]
\bCF^-(\Sigma,\as_1,\bs)
	\arrow[drr,dashed,sloped,"\underline{h}^{\b;\hat{w}}_{\a_1\to \a_2\to \a_3}",pos=.4]
	\arrow[dr, "f_{\a_1\to \a_2}^{\b}",sloped, swap]&& 
	\\
& \bCF^-(\Sigma,\as_2,\bs)
	\arrow[from=uu, crossing over,swap, "\id"]
	\arrow[r,swap,"\ul{f}_{\a_2\to \a_3}^{\b;\hat{w}}"]
&\buCF^-(\Sigma,\as_3,\bs,\hat{w},z)
\end{tikzcd}
\]
\caption{The hyperbox of chain complexes whose compression is $\cC_{\frF}$. We write $\cC_{\frF}^{(1)},\dots, \cC_{\frF}^{(5)}$ for these hypercubes, ordered top-to-bottom.}\label{fig:115}
\end{figure}

 We define the hypercube $\cC_{\frF}$ as the compression of a hyperbox of size $(1,1,5)$, outlined in Figure~\ref{fig:115}.  We explain presently the maps and complexes therein, and then will prove that the diagram satisfies the hyperbox relations.

Let $W_{n+m}$ denote the 2-handle cobordism from $Y$ to $Y_{n+m}(K)$, and let $W_{n+m}'$ denote the cobordism from $Y_{n+m}(K)$ to $Y$, obtained by reversing the orientation. Let $S_{n+m}\subset W_{n+m}'$ denote the induced link cobordism from $(Y_{n+m}(K),U)$ to $(Y,K)$. We orient $S_{n+m}$ so that the chosen orientation of $K\subset Y$ coincides with the boundary orientation of $\d S_{n+m}$. Write $\widehat{S}_{n+m}\subset W_{n+m}'$ for the closed surface obtained by capping $S_{n+m}$ with Seifert surfaces on both ends.

We now describe the maps appearing in the hyperbox in Figure~\ref{fig:115}. The maps $\ul{f}_{\a_2\to \a_3}^{\b;\hat{w}}$ and $\ul{f}_{\a_2\to \bar{\bar{\a}}_3}^{\b;\hat{z}}$ are the ones which count triangles weighted by $U^{n_{w}(\phi)} T^{n_z(\psi)-n_w(\psi)}$, and $U^{n_z(\psi)} T^{n_z(\psi)-n_w(\psi)}$, respectively. These are the maps we have already encountered in the exact triangle.

We view the hyperbox in Figure~\ref{fig:115} as being obtained by stacking 5 hypercubes $\cC_{\frF}^{(1)},\dots, \cC_{\frF}^{(5)}$, ordered from top to bottom. We now explain in more detail the maps which appear in $\cC_{\frF}^{(2)}$ and $\cC_{\frF}^{(3)}$. The maps appearing in $\cC_{\frF}^{(2)}$ are easiest to define if we define them as a pairing of slightly generalized hypercubes of alpha and beta attaching curves, which are equipped with singular twisting chains. We  write down the following non-admissible diagram of alpha attaching curves:
\begin{equation}
\cL_{\a}^{(2)}:=
\begin{tikzcd}[column sep=2cm, labels=description]
\as_1
	\arrow[dddrr,dotted,sloped,pos=.4,swap, "\omega_{\a_3^{\twind}, \a_1}"]
	\arrow[dd, "\Theta_{\a_1,\a_1}^+"]
	\arrow[dr, "\Theta_{\a_2,\a_1}^{\can}",sloped]&&\\
&
\as_2
	\arrow[ddr,dashed, "\lambda_{\a_3^{\twind}, \a_2}",sloped]
	\arrow[r, "\Theta_{\bar{\bar{\a}}_3,\a_2}^{+}"]
&\bar{\bar{\as}}_3
	\arrow[dd, "\Theta_{\a_3^{\twind},\bar{\bar{\a}}_3}^+"]
	\\[1.5cm]
\as_1
	\arrow[drr,dashed, pos=.4,sloped, "\lambda_{\a_3^{\twind}, \a_1}"]
	\arrow[dr, "\Theta_{\a_2,\a_1}^{\can}",swap,sloped]&&\\
&\as_2
	\arrow[r, swap, "\Theta_{\a_3^{\twind},\a_2}^+"]
	\arrow[from=uu, crossing over, "\Theta_{\a_2,\a_2}^+"]
&\as_3^{\twind}.
\end{tikzcd}
\label{eq:L-a-2-fake-diagram}
\end{equation}

Equation~\eqref{eq:L-a-2-fake-diagram} is not a genuine diagram of attaching curves, since $\as_1$ and $\as_2$ are both repeated. In particular, the classes $\Theta_{\a_1,\a_1}^+$ and $\Theta_{\a_2,\a_2}^+$ are not genuine elements of a chain complex. Nonetheless, we may still define polygon counting maps which have them as inputs, by adopting the \emph{strict unality} convention that a triangle map with input $\Theta_{\a_i,\a_i}^+$ is the identity map, and a quadrilateral or pentagon map with $\Theta_{\a_i,\a_i}^+$ as input vanishes. (Alternatively, one could just replace the second copies of $\as_1$ and $\as_2$ with small Hamiltonian translates).

We now wish to promote $\cL_\a^2$ to a hypercube of attaching curves with a cube of singular twisting chains, as follows. We consider the diagram $(\Sigma, \as_3,\as_2,\as_1,w,z)$. Each set of attaching curves determines a 1-chain $\g_3$, $\g_2$ or $\g_1$, as follows. Pick an embedded arc on $\Sigma$ which avoids $\as_3$, $\as_2$ or $\as_1$ and connects $w$ to $z$. Then $\g_i$ is obtained by pushing the interior of this arc into the handlebody $U_{\a_3}$, $U_{\a_2}$ or $U_{\a_1}$. For $\as_1$ and $\as_2$, we pick the canonical short path on $\Sigma$. We pick the 2-chains $C_{\a_2,\a_1}$ to be a disk in the neighborhood of the canonical short path. The choice of $C_{\a_3,\a_2}$ and $C_{\a_3,\a_1}$ is not important, though we assume they are supported in the summands of $Y_{\a_3,\a_2}$ and $Y_{\a_2,\a_1}$ corresponding to the special genus 1 region of the Heegaard triple.  Additionally, we construct 2-chains $C_{\a_1,\a_1}$, $C_{\a_2,\a_2}$ and $C_{\a_3,\a_3}$, as follows. Since the arc $\g_i$ is obtained by pushing an arc off of $\Sigma$, we may push this arc off in both directions to obtain a disk in $Y_{\a_i,\a_i}$.

To equip $\cL_{\a}^{2}$ with a cube of singular twisting chains, we use the arcs $\g_i$ whenever a set of attaching curves is related to $\as_i$ by a sequence of handleslides and isotopies. For example, we decorate $\bar{\bar{\as}}_3$ and $\as_3^{\wind}$ both by $\g_3$. Similarly, if $\ds$ and $\xis$ are a pair of attaching curves in the cube, we use the 2-chains $C_{\a_i,\a_j}$ if $\ds$ (resp. $\xis$) is related to $\as_i$ (resp. $\as_j$) by a sequence of handleslides and isotopies. 

We now pick the length 1 chains in Figure~\ref{eq:L-a-2-fake-diagram} arbitrarily, so that they are homogeneously grades cycles and represent the top degree or canonical elements of homology.

\begin{lem} The cycles in Figure~\ref{eq:L-a-2-fake-diagram} satisfy
\[
\begin{split}
\ul{\frf}_{\bar{\bar{\a}}_3,\a_2,\a_1}\left(\Theta_{\bar{\bar{\a}}_3,\a_2}^+, \Theta_{\a_2,\a_1}^{\can}\right)&=0,\\
\left[f_{\a_3^{\twind},\bar{\bar{\a}}_3,\a_2}\left(\Theta_{\a_3^{\twind}, \bar{\bar{\a}}_3}^+, \Theta_{\bar{\bar{\a}}_3, \a_2}^+\right)\right]&=\left[\Theta_{\a_3^{\twind}, \a_2}^+\right],\\
\left[\ul{\frf}_{\a_3^{\twind}, \a_2,\a_1}\left(\Theta_{\a_3^\twind,\a_2}^+, \Theta_{\a_2,\a_1}^{\can}\right)\right]&=0.
\end{split}
\]
\end{lem}
\begin{proof} The second equality follows since $\as_3^\twind$ are isotopic to $\bar{\bar{\as}}_3$, so the corresponding transition map preserves the top degree generator.
We turn our attention to the first equation. Equation~\eqref{lem:model-computations-m-surgery} implies that triangle classes appear in pairs which have the same $U$ power, and represent conjugate $\Spin^c$ structures. It remains to see that the $T$-powers coincide, for each pair. The surface $\widehat{S}_{\bar{\bar{\a}}_3,\a_2,\a_1}$ is $\widehat{S}_m$, a sphere of self-interestion $-m$. Using the notation of Lemma~\ref{lem:model-computations-m-surgery}, the difference of the $T$-powers representing $\frs_k^+$ and $\frs^-_k$ is $\pm \langle c_1(\frs_k^+), \widehat{S}_{m}\rangle=\pm (2k+1) m$, which is divisible by $m$.  The third equation follows from the first, together with the well-definedness of the cobordism maps.
\end{proof}

We pick $\lambda_{\a_3^{\twind}, \a_2}$ and $\lambda_{\a_3^{\twind}, \a_1}$ so that
\begin{equation}
\begin{split}
\d (\lambda_{\a_3^{\twind}, \a_2})&= f_{\a_3^{\twind}, \bar{\bar{\a}}_3, \a_2}\left(\Theta^+_{\a_3^{\twind}, \bar{\bar{\a}}_3}, \Theta_{\bar{\bar{\a}}_3, \a_2}^+\right)+\Theta_{\a_3^{\twind},\a_2}^+,\\
\d (\lambda_{\a_3^{\twind}, \a_1})&=\ul{\frf}_{\a_3^{\twind},\a_2,\a_1}\left(\Theta_{\a_3^{\twind},\a_2}^+, \Theta_{\a_2,\a_1}^{\can}\right).
\end{split}
\label{eq:eta-classes-C-Flip-2}
\end{equation}
Furthermore, we assume that $\lambda_{\a_3^{\twind},\a_2}$ is of homogenous grading one higher than  $\Theta_{\a_3^{\twind},\a_2}^+$. This is not possible for $\lambda_{\a_3^{\twind}, \a_1}$, since $\ul{\frf}_{\a_3^{\twind}, \a_2,\a_1}$ sums over infinitely many $\Spin^c$ structures. Instead, we assume that $\lambda_{\a_3^{\twind}, \a_1}$ is decomposed as a sum
\[
\lambda_{\a_3^{\twind}, \a_1}:=\sum_{k=0}^{\infty}\lambda_{\a_3^{\twind}, \a_1;\frs_k^+},
\]
where $\frs_k^+$ is the $\Spin^c$ structure on $X_{\a_3^{\twind}, \a_2,\a_1}$ defined in Section~\ref{sec:surgery-m-framed-unknot-1}. We will think of $\lambda_{\a_3^{\twind}, \a_1; \frs_k^-}$ as also being defined, but being 0. (The asymmetry is somewhat arbitrary, but makes our discussion of gradings in Section~\ref{sec:involutive-cone} simpler). Write 
\[
\fru_k:=\{\frs_k^+,\frs_k^-\}\subset \Spin^c(X_{\a_3^\twind, \a_2,\a_1}).
\]

 We assume $\lambda_{\a_3^{\twind}, \a_1;\frs_k^+}$ is chosen so that
\begin{equation}
\d\left( \lambda_{\a_3^{\twind}, \a_1;\frs_k^+}\right)=\ul{\frf}_{\a_3^{\twind},\a_2,\a_1; \fru_k}\left(\Theta_{\a_3^{\twind},\a_2}^+, \Theta_{\a_2,\a_1}^{\can}\right), \label{eq:def-refined-eta}
\end{equation}
and also so that $\lambda_{\a_3^{\twind}, \a_1;\frs_k^+}$ is of homogeneous grading one greater than the expected grading of the right hand side of ~\eqref{eq:def-refined-eta}.

Next, we  construct the length 3 cycle $\omega_{\a_3^{\twind}, \a_1}$. Note that there is a restriction map
\[
r\colon \Spin^c\left(X_{\a_3^{\twind}, \bar{\bar{\a}}_3,\a_2,\a_1}\right)\to   \Spin^c\left(X_{\a_3^{\twind}, \a_2,a_1}\right),
\]
which is an isomorphism on the set of $\Spin^c$ structures which have torsion restriction to the boundary. Abusing notation slightly, write $\fru_k$ also for $r(\fru_k)$.

\begin{lem} For each $k\ge 0$, the chain
\[
C_{\a_3^{\twind},\a_1;\fru_k}:=\ul{\frh}_{\a_3^{\twind}, \bar{\bar{\a}}_3,\a_2,\a_1;\fru_k}\left(\Theta_{\a_3^{\twind}, \bar{\bar{\a}}_3}^+, \Theta_{\bar{\bar{\a}}_3,\a_2}^+,\Theta_{\a_2,\a_1}^{\can}\right)+\ul{\frf}_{\a_3^{\twind}, \a_2,\a_1;\fru_k}\left(\lambda_{\a_3^{\twind}, \a_2}, \Theta_{\a_2,\a_1}^{\can}\right)+\lambda_{\a_3^{\twind}, \a_1;\frs_k^+}
\]
is a boundary in $\bCF^-(\as_1,\as_3^{\twind})\otimes \bF[\Z/m]$.
\end{lem}
\begin{proof} The proof is similar to the proof of Lemma~\ref{lem:C-d-d3-boundary}. An easy computation using associativity shows that $C_{\a_3^{\twind}, \a_1;\fru_k}$ is a cycle. Note also that since we use almost complex structures which are maximally pinched along the connected sum neck, it is sufficient to consider the case then $g(\Sigma)=1$, by Propositions~\ref{prop:stabilize-triangles-general} and \ref{prop:stabilize-quadrilateral}. In the genus 1 case, each chain $C_{\a_3^{\twind}, \a_1; \fru_k}$ is of homogeneous grading $-mk(k+1)+1$ in a chain complex whose homology is $\bF\llsquare U\rrsquare\otimes \bF[\Z_m]$, which is supported in even gradings. Hence $C_{\a_3^{\twind}, \a_1;\fru_k}$ is a boundary.
\end{proof}

For each $k$, we let $\omega_{\a_3^{\twind}, \a_1;\frs_k^+}$ be a chain of homogeneous grading one higher than $C_{\a_3^{\twind}, \a_1;\fru_k}$, such that
\[
\d \left(\omega_{\a_3^{\twind}, \a_1;\frs_k^+}\right)=C_{\a_3^{\twind}, \a_1;\fru_k}.
\]
Finally, we define
\begin{equation}
\omega_{\a_3^{\twind}, \a_1}:=\sum_{k=0}^\infty \omega_{\a_3^{\twind}, \a_1;\frs_k^+}.\label{eq:omega-C-flip-2}
\end{equation}

Finally, we define the maps in $\cC_{\frF}^{(2)}$ by pairing $\cL_{\a}^{(2)}$ with the 0-dimensional hypercube of beta attaching curves consisting of $\bs$ (using the strict unality convention to adapt the construction from Section~\ref{eq:pairing-hypercube}). The pairing also requires a choice of 2-chains $C_{\a_1,\b}$, $C_{\a_2,\b}$ and $C_{\a_3,\b}$. We pick $C_{\a_1,\b}$ and $C_{\a_2,\b}$ to be the disks traced out by push-offs of the canonical short paths on the diagram. We pick $C_{\a_3,\b}$ to be a Seifert surface for $K$.

 Concretely,
\begin{equation}
\begin{split}
\ul{\frH}_{\a_1\to \a_3^\twind}^{\b}&:=\ul{\frh}_{\a_1\to \a_2\to \a_3^{\twind}}^{\b}+\ul{\frf}_{\a_1\to \a_3^\twind}^\b\\
\ul{\frH}_{\a_2\to \a_3^{\twind}}^{\b}&:=\ul{\frh}_{\a_2\to \bar{\bar{\a}}_3\to \a_3^{\twind}}^\b+\ul{\frf}_{\a_2\to \a_3^{\twind}}^{\b},\\
\ul{\frP}_{\a_1\to \a_3^{\twind}}^{\b}&:=\ul{\frp}_{\a_1\to \a_2\to \bar{\bar{\a}}_3\to \a_3^{\twind}}^{\b}+\ul{\frh}_{\a_1\to \a_2\to \a_3^{\twind}}^{\b}+\ul{\frf}_{\a_1\to \a_3^{\twind}}^{\b}.
\end{split}
\label{eq:C-F-2-length-2-and-3-maps}
\end{equation}
In the first line of~\eqref{eq:C-F-2-length-2-and-3-maps}, $\ul{\frh}_{\a_1\to \a_2\to \a_3^{\twind}}^\b$ counts rectangles weighted by $U^{n_w(\psi)} T^{\frm(w,\widehat{S}_{n+m},\psi)}$, with special inputs $\Theta_{\a_2,\a_1}^{\can}$ and $\Theta_{\a_3^{\twind}, \a_3}^+$. The map labeled $\ul{\frf}_{\a_1\to \a_3^{\twind}}^{\b}$ in the first line counts triangles with special input $\lambda_{\a_3^{\twind},\a_1}$, weighted by $U^{n_w(\psi)} T^{\frm(w,\widehat{S}_{n+m},\psi)}$ where $\widehat{S}_{n+m}$ denotes the connected sum of $\widehat{S}_n\subset W'_n(K)$ and $\widehat{S}_m\subset D(-m,1)$.
We note that there is a minor ambiguity in our notation, because the map labeled $\ul{\frf}^{\b}_{\a_1\to \a_3^{\twind}}$ on the bottom line counts triangles with special input $\omega_{\a_3^{\twind},\a_1}$, and the same $T$ and $U$ weighting as the top line.

The maps for $\cC_{\frF}^{(3)}$ are constructed analogously.

\begin{prop}
  The diagram shown in Figure~\ref{fig:115} forms a 3-dimensional hyperbox of chain complexes of size $(1,1,5)$.
\end{prop}
\begin{proof} It is sufficient to show that each of $\cC_{\frF}^{(i)}$ is a 3-dimensional hypercube, for $i\in \{1,2,3,4,5\}$.

Firstly, the diagram $\cC_{\frF}^{(4)}$ is obviously a hypercube. Next, the diagrams $\cC_{\frF}^{(2)}$ and $\cC_{\frF}^{(3)}$ are hypercubes because they are obtained by pairing two hypercubes of attaching curves with cubes of singular twisting chains, as in Section~\ref{sec:twisted-hypercubes-2}.

Next, we consider $\cC_{\frF}^{(1)}$ and $\cC_{\frF}^{(5)}$. We claim
\begin{equation}
T^n\cdot \theta_z\circ \underline{f}_{\a_2\to \bar{\bar{\a}}_3}^{\b;\hat{z}}=\ul{\frf}_{\a_2\to \bar{\bar{\a}}_3}^{\b}\quad \text{and} \quad 
\theta_{w}^{-1}\circ \ul{\frf}_{\a_2\to \a_3}^{\b}=\underline{f}_{\a_2\to \a_3}^{\b;\hat{w}}.\label{eq:change-of-twisting-triangles}
\end{equation}
Equation~\eqref{eq:change-of-twisting-triangles} follows similarly to \cite{OSIntegerSurgeries}*{pg. 18}, as we now describe. For the equation involving $\theta_z$, it suffices to show that if $\psi\in \pi_2(\Theta_{\bar{\bar{\a}}_3,\a_2}^+, \xs,\ys)$ is a class of triangles, then
\begin{equation}
(n_z-n_w)(\psi)+A_{w,z}(\ys)=\frm(z,\widehat{S}_{n+m},\psi)-n. \label{eq:alexander-grading-change}
\end{equation}
Equation~\eqref{eq:alexander-grading-change} can be verified by interpreting the triple $(\Sigma,\bar{\bar{\as}}_3,\as_2,\bs,w,z)$ as determining the link cobordism $(W_{n+m}',S_{n+m})$ from $(Y_{n+m}(K),U)$ to $(Y,K)$. Noting that $A_{w,z}(\xs)=0$ for all $\xs\in \bT_{\a_2}\cap \bT_{\b}$, since $(\Sigma,\as_2,\bs,w,z)$ represents an unknot in $Y_{n+m}(K)$, the left hand side of~\eqref{eq:alexander-grading-change} is the Alexander grading change of the link cobordism map, which coincides with $\frm(w,\widehat{S}_{n+m},\psi)$ by \cite{ZemAbsoluteGradings}*{Theorem~1.4}. On the other hand,
\begin{equation}
\frm(w,\widehat{S}_{n+m},\psi)-\frm(z,\widehat{S}_{n+m},\psi)=\frac{\left\langle c_1(\frs_w(\psi))-c_1(\frs_z(\psi)),[\widehat{S}_{n+m}]\right\rangle}{2}=-(n+m).\label{eq:diff:mw-mz}
\end{equation}
The first equality of ~\eqref{eq:diff:mw-mz} is from the definition, while the second follows from the relation $\frs_{w}(\psi)-\frs_z(\psi)=\PD [\widehat{S}_{n+m}]$; see \cite{OSIntegerSurgeries}*{pg. 11} or more generally \cite{ZemAbsoluteGradings}*{Lemma~3.9}.

 The equality involving $\theta_w$ in~\eqref{eq:change-of-twisting-triangles} is proven by the same reasoning. Similarly, 
\begin{equation}
\begin{split}
T^n\cdot\theta_z\circ\underline{h}_{\a_1\to \a_2\to \bar{\bar{\a}}_3}^{\b;\hat{z}}&=\underline{\frh}_{\a_1\to \a_2\to \bar{\bar{\a}}_3}^{\b}\quad \text{and} \\
\theta_w^{-1} \circ \underline{\frh}_{\a_1\to \a_2\to \a_3}^{\b}&=\underline{h}_{\a_1\to \a_2\to \a_3}^{\b;\hat{w}}. \label{eq:change-of-twisting-squares}
\end{split}
\end{equation}
Equations~\eqref{eq:change-of-twisting-triangles} and~\eqref{eq:change-of-twisting-squares} imply the length 2 and 3 hypercube relations for $\cC_{\frF}^{(1)}$ and $\cC_{\frF}^{(5)}$. The length 0 and 1 relations are easy to verify, so we conclude that $\cC_{\frF}^{(1)}$ and $\cC_{\frF}^{(5)}$ are hypercubes. The proof is complete.
 \end{proof}

\section{Constructing the null-homotopy \texorpdfstring{$H_E$}{H-E}}
\label{sec:N-null-homotopy}

In Section~\ref{sec:constructing-C-cen}, we showed that to finish the construction of the central hypercube, it is sufficient to
construct a map 
\[
H_E\colon \buCF^-(\Sigma\# \bar{\Sigma}, \gs, \ds')\to \buCF^-(\Sigma_0\# \bar{\Sigma},\as\cup \bar{\bs}, \Ds_3)
\]
 satisfying
\[
[\d, H_E]=F_3 h_1+F_3h_2+F_3 h_3.
\]
We recall that the domain and codomain of $H_E$ are both complexes for $\buCF^-(Y)$. In this section, we describe the construction of the map $H_E$.

\subsection{The endomorphism \texorpdfstring{$E$}{E}}

\begin{figure}[H]
	\centering
\begingroup%
  \makeatletter%
  \providecommand\color[2][]{%
    \errmessage{(Inkscape) Color is used for the text in Inkscape, but the package 'color.sty' is not loaded}%
    \renewcommand\color[2][]{}%
  }%
  \providecommand\transparent[1]{%
    \errmessage{(Inkscape) Transparency is used (non-zero) for the text in Inkscape, but the package 'transparent.sty' is not loaded}%
    \renewcommand\transparent[1]{}%
  }%
  \providecommand\rotatebox[2]{#2}%
  \newcommand*\fsize{\dimexpr\f@size pt\relax}%
  \newcommand*\lineheight[1]{\fontsize{\fsize}{#1\fsize}\selectfont}%
  \ifx\svgwidth\undefined%
    \setlength{\unitlength}{159.73141703bp}%
    \ifx\svgscale\undefined%
      \relax%
    \else%
      \setlength{\unitlength}{\unitlength * \real{\svgscale}}%
    \fi%
  \else%
    \setlength{\unitlength}{\svgwidth}%
  \fi%
  \global\let\svgwidth\undefined%
  \global\let\svgscale\undefined%
  \makeatother%
  \begin{picture}(1,1.06894878)%
    \lineheight{1}%
    \setlength\tabcolsep{0pt}%
    \put(0,0){\includegraphics[width=\unitlength,page=1]{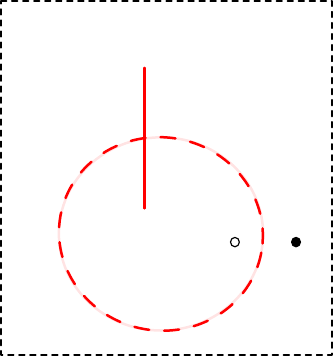}}%
    \put(0.34016772,-0.04109036){\color[rgb]{0,0,0}\makebox(0,0)[lt]{\begin{minipage}{8.98584476\unitlength}\end{minipage}}}%
    \put(0.89108923,0.37854488){\color[rgb]{0,0,0}\makebox(0,0)[t]{\lineheight{1.25}\smash{\begin{tabular}[t]{c}$w_0$\end{tabular}}}}%
    \put(0.70560751,0.37854488){\color[rgb]{0,0,0}\makebox(0,0)[t]{\lineheight{1.25}\smash{\begin{tabular}[t]{c}$z_0$\end{tabular}}}}%
    \put(0,0){\includegraphics[width=\unitlength,page=2]{fig171.pdf}}%
    \put(0.22086925,0.11560215){\color[rgb]{0,0,0}\makebox(0,0)[rt]{\lineheight{1.25}\smash{\begin{tabular}[t]{r}$c_3'$\end{tabular}}}}%
    \put(0.71336586,0.58801193){\color[rgb]{0,0,0}\makebox(0,0)[lt]{\lineheight{1.25}\smash{\begin{tabular}[t]{l}$c_3$\end{tabular}}}}%
    \put(0.39864029,0.25282862){\color[rgb]{0,0,0}\makebox(0,0)[lt]{\lineheight{1.25}\smash{\begin{tabular}[t]{l}$c_2$\end{tabular}}}}%
    \put(0.42416327,0.81283261){\color[rgb]{0,0,0}\makebox(0,0)[rt]{\lineheight{1.25}\smash{\begin{tabular}[t]{r}$c_1$\end{tabular}}}}%
    \put(0,0){\includegraphics[width=\unitlength,page=3]{fig171.pdf}}%
  \end{picture}%
\endgroup%

	\caption{The Heegaard quadruple $\cQ_1^0$.}\label{fig:171}
\end{figure}

In this section, we describe a more general way to organize the maps $h_1$, $h_2$ and $h_3$ from the central hypercube.  We consider genus 1, doubly pointed quadruples 
\[
\cQ_1^0=(\bT^2,c_3,c_2,c_1,c_3',w_0,z_0), \quad \cQ_2^0=(\bT^2,c_2,c_1,c_3',c_2',w_0,z_0)
\]
\[
\quad \text{and} \quad \cQ_3^0=(\bT^2,c_1,c_3',c_2',c_1',w_0,z_0).
\]
The quadruple $\cQ_1^0$ is shown in Figure~\ref{fig:171}, and $\cQ_2^0$ and $\cQ_3^0$ are obtained similarly.

The doubly pointed quadruple $\cQ_1^0$ determines a properly embedded surface 
$ S_{c_3,c_2,c_1,c_3'}\subset X_{c_3,c_2,c_1,c_3'}$, as described in Section~\ref{sec:twisted-hypercubes-2}.  The 4-manifold $X_{c_3,c_2,c_1,c_3'}$ has four boundary components.   Two of these boundary components are $S^3$, one is $S^1\times S^2$, and the other is $L(m,1)$. The surface $S_{c_3,c_2,c_1,c_3'}$ intersects each component of $\d X_{c_3,c_2,c_1,c_3'}$ in an unknot. Upon attaching 3- and 4-handles to $\d X_{c_3,c_2,c_1,c_3'}$, and filling in $\d S_{c_3,c_2,c_1,c_3'}$ with disks, we obtain a sphere of self-intersection $-m$ in the disk bundle $D(-m,1)$. The quadruples $\cQ_2^0$ and $\cQ_3^0$ determine surfaces and 4-manifolds with the same description.

Suppose that $(\Sigma,\as,\bs,w,z)$ is a
 doubly pointed diagram for $(Y,K)$. We form the connected sum $\Sigma\# \bT^2$ near the points $w_0\in \bT^2$ and $z\in \Sigma$. 
 We define the following sets of attaching curves on $\Sigma\# \bT^2$:
\[
\as_1=\as\cup \{c_1\}, \quad \as_2=\as'\cup \{c_2\},\quad \as_3=\as''\cup \{c_3\},
\]
\[
 \quad \bs_3=\bs\cup \{c_3'\},\quad \bs_2:=\bs'\cup \{c_2'\},\quad \text{and} \quad \bs_1:=\bs''\cup \{c_1'\}.
\]
Here $\as'$ and $\as''$ are small Hamiltonian translates of $\as$. Similarly $\bs'$ and $\bs''$ are small Hamiltonian translates of $\bs$, and $c_i'$ denotes a translate of $c_i$.

We define the quadruples
\[
\begin{split}
\cQ_1&:=(\Sigma\# \bT^2, \as_3,\as_2,\as_1,\bs_3,w,z_0)\\
\cQ_2&:=(\Sigma\# \bT^2, \as_2,\as_1,\bs_3,\bs_2,w,z_0)\\
\cQ_3&:=(\Sigma\# \bT^2, \as_1,\bs_3,\bs_2,\bs_1,w,z_0).
\end{split}
\]
We view $\cQ_1$ as being obtained by taking the connected sum of $\cQ_1^0$ with the quadruple $(\Sigma,\as'',\as',\as,\bs,w,z)$ at points $p_0$ and $p$, very near $w_0$ and $z$, respectively. Then we delete  the basepoints $w_0$ and $z$.

\begin{figure}[h]
	\centering
	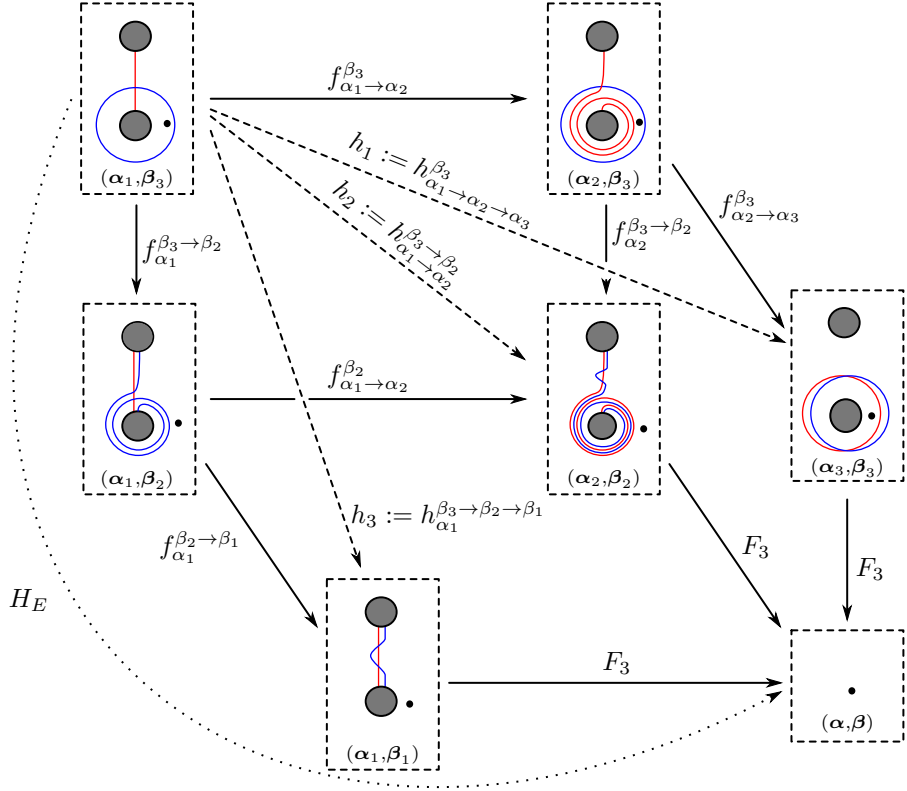
	\caption{A null-homotopy $H_E$ of $\exc$ gives a length 3 arrow which makes this diagram a hypercube of chain complexes.}\label{fig:95}
\end{figure}

In this general setting, we define maps $h_1$, $h_2$ and $h_3$ as rectangle counts
\begin{equation}
h_1:=h_{\a_1\to \a_2\to \a_3}^{\b_3},\quad h_2:=h_{\a_1\to \a_2}^{\b_3\to \b_2},\quad \text{and}\quad h_3:=h_{\a_1}^{\b_3\to \b_2\to \b_1},\label{eq:h1h2h3-slightly-more-general}
\end{equation}
where we weight a rectangle class by $U^{n_{w}(\psi)}T^{(n_{z_0}-n_{w})(\psi)}$. Here, we also use almost complex structures which are maximally pinched along a circle bounding the genus 1 region shown in Figure~\ref{fig:95}.  The maps labeled $h_1$, $h_2$ and $h_3$ in Section~\ref{sec:constructing-C-cen} coincide with the maps defined in ~\eqref{eq:h1h2h3-slightly-more-general}, for a particular choice of Heegaard diagram $(\Sigma,\as,\bs,w,z)$.

The maps $h_1$, $h_2$ and $h_3$ fit into the diagram of chain complexes shown in Figure~\ref{fig:95}.  In fact, this diagram satisfies all of the hypercube relations except for the length 3 relation. (To see the length 2 relations along the bottom and right faces, we use Proposition~\ref{prop:stabilize-triangles-general} to destabilize the triples, and Proposition~\ref{prop:nearest-point-triangles} to identify the resulting triangle maps with nearest point maps).

 We now define
\[
\exc:=F_3\circ h_1+F_3\circ h_2+F_3\circ h_3,
\]
where $F_3$ denotes a 3-handle map.  In light of the discussion in Section~\ref{sec:constructing-C-cen}, to finish the construction of $\cC_{\cen}$, it is sufficient to show the following:

\begin{prop} \label{prop:Esimeq0}
The map $\exc$ is chain homotopic to zero. Equivalently, there is a length 3 map $H_E$ which makes the diagram in Figure~\ref{fig:95} into a hypercube of chain complexes.
\end{prop}

 In the remainder of this section, we construct a null-homotopy of $\exc$, via a holomorphic deformation.

\subsection{Vector fields for deformations}

\label{sec:vector-fields}

In this section, we define several vector fields on $[0,1]\times \R$, $\Delta$ and $\Box$, which we will use in our proof of Proposition~\ref{prop:Esimeq0}. We will use these vector fields to deform the quadrilateral maps appearing in the definition of the map $E$.

 We begin by defining two vector fields on $[0,1]\times \R$, for which we write $v_Y$ and $v_H$. Let $f\colon \R\to \R$ be a smooth function, such that $f(0)=0$, $f(1)=0$ and $f(s)>0$ for $s\in (0,1)$. 
 We define $v_H$ to be
\begin{equation}
(v_H)_{(s,t)}:=f(s)\cdot \d/\d s. \label{eq:def-vector-field-v}
\end{equation}  

Next, we construct a vector field $v_Y$ on $[0,1]\times \R$ which satisies the following:
\begin{enumerate}[label=($v_Y$-\arabic*), ref=$v_Y$-\arabic*, leftmargin=*,widest=IIII]
\item $v_Y$ is invariant under the $\R$ action.
\item $v_Y$ is 0 on $\{0,1\}\times \R$, and is non-zero on $(0,1)\times \R$.
\item  $v_Y$ extends smoothly to a vector field on $\R\times \R$.
\item There is an $\R$-invariant neighborhood of the line $s=1/2$ such that
\begin{equation}
v_Y= \d/\d t+(s-1/2)\cdot \d/\d s. \label{eq:vector-fields-v-y-technical}
\end{equation}
\item The closure of each flowline of $v_Y$, except for $s=1/2$, has one transverse intersection point with $\{0,1\}\times \R$.
\end{enumerate}

\begin{figure}[ht!]
	\centering
	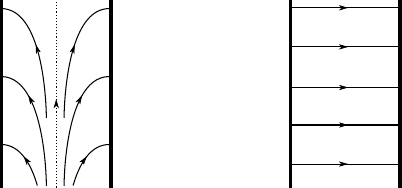
	\caption{The $\R$-invariant vector fields $v_Y$ (left) and $v_H$ (right) on $[0,1]\times \R$.}\label{fig:143}
\end{figure}

 Equation~\eqref{eq:vector-fields-v-y-technical} is useful in controlling the level structure of degenerations of matched combs. The two vector fields $v_H$ and $v_Y$ are schematically shown in Figure~\ref{fig:143}.  

\begin{figure}[ht!]
	\centering
	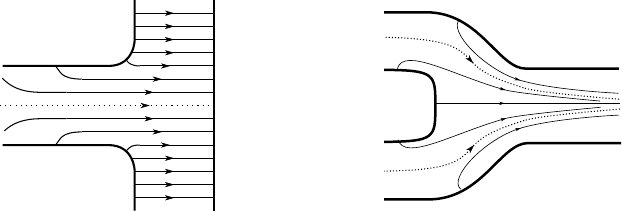
	\caption{ The vector fields $v_{\Delta,1}$ (left) and $v_{\Delta,2}$ (right).}\label{fig:140}
\end{figure}

We now discuss two important vector fields for deforming the holomorphic triangle maps. We call these $v_{\Delta,1}$ and $v_{\Delta,2}$. These are shown in Figure~\ref{fig:140}. We assume that $v_{\Delta,1}$ satisfies the following:
\begin{enumerate}[label=($v_{\Delta,1}$-\arabic*), ref=$v_{\Delta,1}$-\arabic*, leftmargin=*,widest=IIII]
\item  $v_{\Delta,1}$ vanishes on $\d \Delta$ and extends to a smooth vector field on $\C$.
\item $v_{\Delta,1}$ is non-vanishing on $\Int (\Delta)$.
\item  $v_{\Delta,1}$ is $\R$-invariant on each of the cylindrical ends of $\Delta$.
\item $v_{\Delta,1}$ coincides with $v_H$ on two cylindrical ends, and $- v_Y$ on the remaining end.
\item The closure of all but one flowline of $v_{\Delta,1}$ has two transverse intersection points with $\d \Delta$. The remaining flowline has just one intersection.
\end{enumerate}

We construct the vector field $v_{\Delta,2}$ shown in Figure~\ref{fig:140}, except we require $v_{\Delta,2}$ to coincide with $v_Y$ in all three ends of $\Delta$. Furthermore, all but two flowlines have closures with exactly one transverse intersection with $\d \Delta$, while two flowlines have closures which are bounded away from $\d \Delta$.

We now describe three families of vector fields on $\Box$. To meaningfully deform the quadrilateral maps, we need to describe pairs consisting of families of vector fields and almost complex structures on $\Box$. We think of these families as being associated to the Heegaard quadruples $\cQ_1$, $\cQ_2$ and $\cQ_3$.

\begin{figure}[ht!]
	\centering
	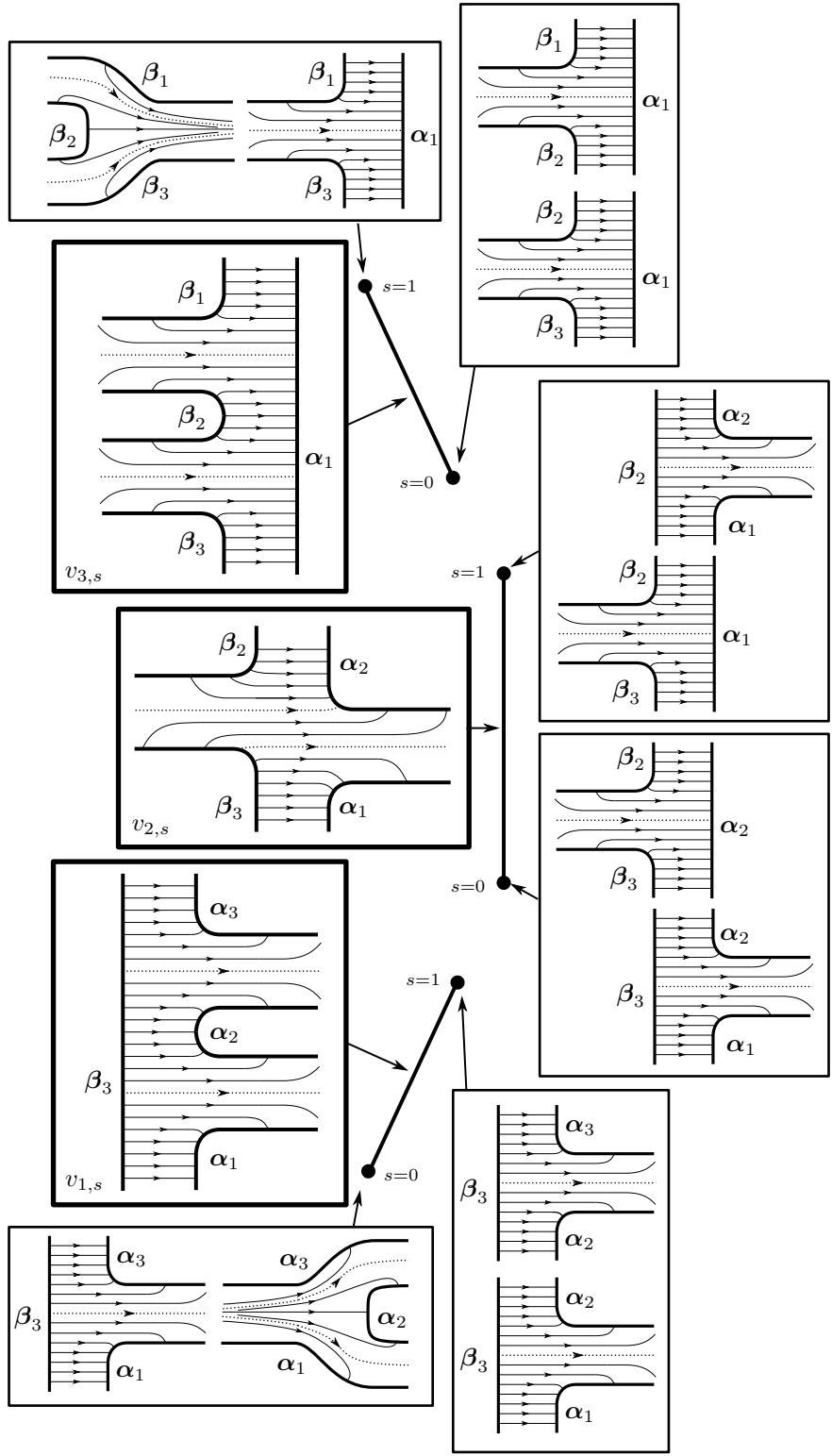
	\caption{The 1-parameter families $(v_1,j_1)$, $(v_2,j_2)$ and $(v_3,j_3)$.}\label{fig:145}
\end{figure}

For $i\in \{1,2,3\}$,  we pick smooth families of vector fields $v_i:=(v_{i,s})_{s\in (0,1)}$ and almost complex structures $j_i=(j_{i,s})_{s\in (0,1)}$, on $\Box$,  which have the configuration shown in Figure~\ref{fig:145}. We assume $v_{1}$ and $j_1$ satisfy the following:
 \begin{enumerate}[label=($v_1$-\arabic*), ref=$v_1$-\arabic*, leftmargin=*,widest=IIII]
\item\label{v-1-1} In each cylindrical end of $\Box$, $v_{1,s}$ and $j_{1,s}$ are constant in $s$ and $v_{1,s}$ coincides with $v_H$ or $v_Y$.
\item\label{v-1-2} Near $s=0$, the vector field $v_{1,s}$ is obtained by inserting copies of $([0,1]\times (-R,R),v_Y)$, for increasing $R$, into the $(\as_1,\as_3)$-subregion. Similarly, $j_{1,s}$ is obtained by inserting copies of $[0,1]\times (-R,R)$, equipped with the standard almost complex structure. The limiting vector field consists of one copy of $v_{\Delta,1}$ and one copy of $v_{\Delta,2}$.
\item\label{v-1-3} Near $s=1$, $v_{1,s}$ is obtained by inserting copies of $([0,1]\times (-R,R),v_H)$, for increasing $R$, into the $(\as_2,\bs_3)$-subregion. The analogous statement holds for $j_1$. The limiting vector field consists of two copies of $\pm v_{\Delta,1}$.
\item\label{v-1-4} $v_{1,s}|_{\d \Box}=0$ for each $s$.
 \end{enumerate}
 
We assume that the families $(v_2,j_2)$ and $(v_3,j_3)$ satisfy the analogous conditions. 
 
\begin{rem}
Families  $v_1$ and $j_1$ satisfying~\eqref{v-1-1}--\eqref{v-1-4} may be constructed as follows. We may identify the moduli space of complex rectangles with the space of configurations of two marked points along $\{1\}\times \R\subset [0,1]\times \R$ (with the constant family of almost complex structures $j_1$). The family $v_{1}$ may be constructed by rescaling the vector field $\d/ \d s$ on $[0,1]\times \R$, similar to ~\eqref{eq:def-vector-field-v}.
\end{rem}

\begin{rem}
 Many of the vector fields have one or two non-generic flowlines. For example, $v_Y$ has a single flowline ($s=1/2$) whose closure is bounded away from $\{0,1\}\times \R$. Similarly $v_{\Delta,1}$ has one flowline whose closure intersects $\d \Delta$ only once. We refer to these flowlines, and their natural analogs on $v_{\Delta,2}$ and $v_i$, as the \emph{special lines}.
\end{rem}

\subsection{Local degree theory}
\label{sec:gluing-generalities}

In this section, we frequently obtain algebraic relations by counting the ends of 1-dimensional matched moduli spaces. A natural strategy for such an argument is to show that each point in the boundary of the 1-dimensional moduli space in question has a neighborhood modeled on $[0,\infty)$. Proving such a claim usually requires defining a smooth structure on the compactification of the moduli space of unmatched curves, and proving that the relevant evaluation maps extend smoothly. Instead we opt to use degree theory argument of \cite{LOTBordered}, which is sufficient for the purpose of determining algebraic relations. In this section, we collect several background results which we will subsequently use.

 We refer the reader to \cite{LOTBordered}*{Definition~5.36} for the definition of a stratified space which we use. We now recall \cite{LOTBordered}*{Definition~5.37}:
 
 \begin{define}
 \label{def:locally-odd-degree} Let $X$ be a stratified space so that the top stratum is a smooth $m$-manifold, and suppose that $f\colon X\to [0,\infty)^m$ is a stratified map so that $f^{-1}((0,\infty)^m)$ is the top stratum of $X$. Let $p\in [0,\infty)^m$ and assume $f$ is proper near $p$. We say that $f$ is \emph{odd degree near $p$} if there is an open neighborhood $U$ of $p$ such that for any regular value $p'\in U\cap (0,\infty)^m$, the set $f^{-1}(p')$ consists of an odd number of points.
 \end{define}
 
 In practice, we will say a broken holomorphic curve appearing in the boundary of a 1-dimensional moduli space \emph{makes odd algebraic contribution}, if the closure of the 1-dimensional moduli space admits a stratified map to $[0,\infty)$, which is locally proper and odd degree near 0, such that the preimage of $0$ is the broken curve.

The following lemma is a helpful criterion for constructing stratified spaces as fibered products which have an evaluation map to $[0,\infty)$ which is odd degree near $0$:
  
  \begin{lem}\label{lem:linking}
  Suppose that $X^j$ and $Y^k$ are topological manifolds without boundary, and suppose  $X\times [0,\infty)$ and $Y\times [0,\infty)$ are given the structure of stratified spaces (i.e. given a choice of smooth structure on $\{0\}\times X$, $(0,\infty)\times X$, $\{0\}\times Y$ and $(0,\infty)\times Y$). Suppose that 
  \[
f\colon [0,\infty) \times X\to [0,\infty)\times \R^{j+k}\quad \text{and} \quad g\colon [0,\infty)\times Y\to [0,\infty)\times \R^{j+k}   
  \]
  are proper, stratified maps (i.e. continuous maps which are smooth on the interior of each strata). Furthermore, suppose that $f^{-1}(0\times \R^{j+k})=0\times X$ and similarly for $g$. Finally, suppose that $f|_{0\times X}$ and $g|_{0\times Y}$ have transverse intersection at a single point $\ve{0}=0\times \ve{0}\in 0\times \R^{j+k}\subset [0,\infty)\times \R^{j+k}$, and that $f$ and $g$ are transverse on the top strata. Then projection onto $[0,\infty)$ of the fibered product of $f$ and $g$ is odd degree near $\ve{0}$.
 \end{lem}

 \begin{proof} Let $\pi_1\colon [0,\infty)\times \R^{j+k}\to [0,\infty)$ denote projection onto $[0,\infty)$, and let $\pi_2$ denote projection onto $\R^{j+k}$. We first claim that the map 
 \[
E:= (\pi_1\circ f,\pi_1\circ g,\pi_2\circ f-\pi_2\circ g)\colon [0,\infty)^2\times X \times Y \to [0,\infty)^2\times \R^{j+k}
 \]
  is odd degree near $\ve{0}$.  To see this, note that $E|_{0\times 0\times X\times Y}$ is proper near $\ve{0}$ and also odd degree near $\ve{0}$, since $f|_{0\times X}$ and $g|_{0\times Y}$ are transverse, and hence $E|_{[0,\infty)\times 0\times X\times Y}$ is also odd degree near 0 by \cite{LOTBordered}*{Lemma~5.38~(2)}. The proof of \cite{LOTBordered}*{Lemma~5.38~(2)} implies that $E$ is odd degree near $\ve{0}$ on $[0,\infty)^2\times X\times Y$. Transversality on the top stratum to $\Delta\times \R^{j+k}$, where $\Delta\subset [0,\infty)^2$ is the diagonal, guarantees that the fiber over $\Delta\times \ve{0}\subset [0,\infty)^2\times \R^{j+k}$, with its projection to $[0,\infty)$, is proper and odd-degree near zero.  But that is equivalent to the projection onto $[0,\infty)$ of the fibered product of $f$ and $g$ being odd-degree near zero.   
\end{proof}

\subsection{The connected sum formula by deforming the matching}
\label{sec:time-dilation-complexes}

As a warmup to the holomorphic degeneration we use to analyze $E$, we describe a deformation of the matched moduli spaces in Section~\ref{sec:matched-disks} which gives the connected sum formula for $\CF^-$. The deformation
is inspired by the time-dilation degeneration described in \cite{LOTBordered}, as well as
Ozsv\'{a}th and Szab\'{o}'s bordered version of knot Floer homology \cite{OSBorderedKauffmanStates} \cite{OSMatchings} \cite{OSBorderedHFK},
and the forthcoming minus version of bordered Floer homology \cite{LOT-minus-bordered}. Our argument deforms the matching in the $[0,1]$-direction.

We are interested in the following two connected sum operations:
\begin{enumerate}[label=($\# $-\arabic*), ref=$\#$-\arabic*, leftmargin=*,widest=IIII]
\item\label{CS-1} Suppose $\cH_1=(\Sigma_1,\as_1,\bs_1,w_1)$ and $\cH_2=(\Sigma_2,\as_2,\bs_2,w_2)$ are two singly pointed Heegaard diagrams. Let $p_1$ and $p_2$  be two points on $\Sigma_1$ and $\Sigma_2$, in the complement of the alpha and beta curves, and suppose that $p_2$ is immediately adjacent to $w_2$. Let $\cH_1\# \cH_2$ denote their connected sum, taken at $p_1$ and $p_2$. In forming $\cH_1\# \cH_2$, we delete $w_2$, and leave $w_1$.
\item\label{CS-2}  There is a version of the previous operation for knots, as follows. Suppose $\cH_1=(\Sigma_1,\as_1,\bs_1,w_1,z_1)$ and $\cH_2=(\Sigma_2,\as_2,\bs_2,w_2,z_2)$ are two doubly based knot diagrams. We take the connected sum $\cH_1\# \cH_2$ at $z_1$ and $w_2$. We delete $z_1$ and $w_2$, but leave $w_1$ and $z_2$.
\end{enumerate}

We now describe a deformation of the matched complexes from Section~\ref{sec:matched-moduli-spaces}. If $\tau\in \R$, let $\rho_\tau\colon [0,1]\to [0,1]$ be the flow of $v_H$. Note that $\rho_\tau$ fixes $\{0,1\}$ for all $\tau$, but has no fixed points on $(0,1)$ for $\tau \neq 0$.

The map $\rho_\tau$ determines a diffeomorphism
\[
R^\tau\colon [0,1]\times \R\to [0,1]\times \R,
\]
given by
\begin{equation}
R^{\tau}(s,t)=(\rho_\tau(s),t). \label{eq:R-t-def}
\end{equation}
The diffeomorphism $R^\tau$ also gives a map
\[
R^\tau\colon ([0,1]\times \R)^k\to ([0,1]\times \R)^k,
\]
given by
\begin{equation}
R^\tau(x_1,\dots, x_k)=(R^\tau(x_1),\dots, R^\tau(x_k)), \label{eq:R-map-symmetric-product}
\end{equation}
which fixes the fat-diagonal setwise.

Similar to the matched moduli spaces from ~\eqref{eq:matched-moduli-space-disks}, if $S_1^{\qs_1}$ and $S_2^{\qs_2}$ are two $k$-pointed source curves, we define the $\tau$-matched moduli spaces as 
\begin{equation}
\begin{split}
&\cM\cM^{\tau}_{J_1\wedge J_2}(S_1^{\qs_1}, S_2^{\qs_2}, \phi_1,\phi_2)\\
:=&\left\{(u_1,u_2)\in \cM_{J_1}(S_1^{\qs_1},\phi_1)\times \cM_{J_2}(S_2^{\qs_2}, \phi_2)\middle \vert \begin{array}{c} (\pi_{\Sigma_1}\circ u_1)(q_{1,i})= p_1,\\
(\pi_{\Sigma_2}\circ u_2)(q_{2,i})=p_2,\\
R^{\tau}( (\pi_{\bD}\circ u_1)(q_{1,i}))= (\pi_{\bD}\circ u_2)(q_{2,i}), \\
 \text{ for } i=1,\dots, k
 \end{array}\right\}.
 \end{split}
 \label{eq:t-matched-moduli-space-disks}
\end{equation}

We define a $\tau$-matched differential on the complex $\CF^-(\cH_1, \cH_2)$ from Section~\ref{sec:totally-degenerated-ac-complexes} via the formula
\[
\d^{\tau}(\xs_1\times \xs_2):=\sum_{\substack{\phi_1\in \pi_2(\xs_1,\ys_1)\\
\phi_2\in \pi_2(\xs_2,\ys_2)\\
n_{p_1}(\phi_1)=n_{p_2}(\phi_2)\\
 \ind_{\emb}(\phi_1,\phi_2,M_{\#})=1}} \# (\cM\cM^{\tau}_{J_1\wedge J_2}(\phi_1,\phi_2)/\R)\cdot  U^{n_{w_1}(\phi_1)}\cdot (\ys_1\times \ys_2). 
\]
We weight the counts by $n_{w_1}(\phi_1)$ since we delete $w_2$ in the operation~\eqref{CS-1}. Note that if we instead used the construction in \eqref{CS-2}, we would use $w_1$ and $z_2$ to form the weights.

\begin{lem}\label{lem:t-matched-complexes}
\begin{enumerate}
\item[]
\item The endomorphism $\d^{\tau}$ is a differential.
\item If $\tau_1,\tau_2\in \R$, then $(\CF^-(\cH_1, \cH_2),\d^{\tau_1})\simeq (\CF^-(\cH_1, \cH_2), \d^{\tau_2})$.
\end{enumerate}
\end{lem}
\begin{proof} 
 The proof that $\d^\tau$ squares to zero is not substantially different than the proof of Lemma~\ref{lem:matched-d^2=0}.

The proof of the second claim is similar to the proof of  \cite{LOTBordered}*{Proposition~9.22}, and follows by defining an appropriate continuation map, as we now describe.  Let $\theta\colon \R\to \R$ be a monotonic, smooth function such that there are $a,b\in \R$ such that $\theta(t)=\tau_1$ if $t\le a$, and $\theta(t)=\tau_2$ if $t\ge b$. The map $\theta$ determines a diffeomorphism
\[
R^\theta\colon [0,1]\times \R\to [0,1]\times \R,\quad  \text{defined by} \quad 
R^\theta(s,t)=(\rho_{\theta(t)}(s),t),
\]
where $\rho_\tau$ denotes the time $\tau$ flow of the vector field $v$ in ~\eqref{eq:def-vector-field-v}. The map $R^\theta$ induces an automorphism of $([0,1]\times \R)^k$, which fixes the fat-diagonal, setwise.

We define a transition map 
\[
\Psi^{\tau_1\to \tau_2}\colon (\CF^-(\cH_1, \cH_2),\d^{\tau_1})\to (\CF^-(\cH_1, \cH_2), \d^{\tau_2})
\] 
by counting index 0 holomorphic curves which are $\theta$-matched, i.e., counting elements of the moduli spaces
\begin{equation}
\begin{split}
&\cM\cM^{\theta}_{J_1\wedge J_2}(S_1^{\qs_1}, S_2^{\qs_2}, \phi_1,\phi_2)\\
:=&\left\{(u_1,u_2)\in \cM_{J_1}(S_1^{\qs_1},\phi_1)\times \cM_{J_2}(S_2^{\qs_2}, \phi_2)\middle \vert \begin{array}{c} (\pi_{\Sigma_1}\circ u_1)(q_{1,i})=p_1,\\
(\pi_{\Sigma_2}\circ u_2)(q_{2,i})=p_2,\\
R^\theta((\pi_{\bD}\circ u_1)(q_{2,i}))=(\pi_{\bD}\circ u_2)(q_{1,i}),\\
 \text{ for } i=1,\dots, k
 \end{array}\right\},
 \end{split}
 \label{eq:phi-matched-moduli-space-disks}
\end{equation}
when $\ind_{\emb}(\phi_1,\phi_2,M_{\#})=0$.

To see that $\Psi^{\tau_1\to \tau_2}$ is a chain map, one counts the ends of the $\phi$-matched moduli spaces for pairs $(\phi_1,\phi_2)$ with $\ind_{\emb}(\phi_1,\phi_2,M_{\#})=1$. It is easily checked, using expected dimensions, as in Proposition~\ref{prop:dimension-counts}, that for generically chosen almost complex structures, the ends of index 1 moduli spaces correspond exactly to broken curves with two stories, where one story has index 1 and is $\tau_1$ or $\tau_2$ matched, and the other story has index 0 and is $\phi$-matched. 

 To see that $\Psi^{\tau_1\to \tau_2}$ is a chain homotopy equivalence, one argues that $\Psi^{\tau_2\to \tau_1}$ is a homotopy inverse, by showing that $\Psi^{\tau_2\to \tau_1}\circ \Psi^{\tau_1\to \tau_2}$ and $\Psi^{\tau_1\to \tau_2}\circ \Psi^{\tau_2\to \tau_1}$ are chain homotopic to the identity. The chain homotopy is constructed  by counting index $-1$ curves which are matched over a 2-parameter family of matching conditions.
\end{proof}

We now show that the $\tau$-matched complex is chain homotopic to the ordinary tensor product complex  $\CF^-(\cH_1)\otimes_{\bF[U]} \CF^-(\cH_2)$, by letting $\tau\to \infty$.  Write 
\[
\d^\otimes=\d_1\otimes \id+\id\otimes \d_2
\]
for the differential on  $\CF^-(\cH_1)\otimes \CF^-(\cH_2)$.

If $S^{\qs}$ is a marked source, and $q\in \qs$, then we define \emph{stable} and \emph{unstable} evaluation maps
\[
\ev_{q}^{\st}\colon \cM_J(S^{\qs},\phi)\to \{1\}\times \R  \quad \text{and} \quad \ev_{q}^{\un}\colon \cM_J(S^{\qs},\phi)\to \{0\}\times \R,
\]
via the formulas
\[
\ev_{q}^{\st}(u):=\lim_{\tau\to + \infty} (R^\tau\circ \pi_{\bD}\circ  u)(q) \quad \text{and} \quad \ev_{q}^{\un}(u):=\lim_{\tau\to - \infty} (R^\tau\circ \pi_{\bD}\circ  u)(q)
\]
where  $R^\tau$ denotes the flow of $v_H$.

\begin{define}\label{def:infinity-matched-disk}
 An \emph{$\infty$-matched holomorphic disk} representing $(\phi_1,\phi_2)\in \pi_2(\xs_1,\ys_1)\times \pi_2(\xs_2,\ys_2)$ consists of a pair $(\cU_1,\cU_2)$ of 1-story holomorphic combs on $\cH_1$ and $\cH_2$, respectively, satisfying the following:
 \begin{enumerate}
 \item $\cU_1=(u_1,a)$, where $u_1\in \cM_{J_1}(S_1^{\qs_1},\phi_1')$ satisfies \ref{M-1}--\ref{M-7}, and $a$ is the disjoint union of $n$ boundary degenerations $a_1,\dots, a_n$, which each have Maslov index $2$. Furthermore $\phi_1=\phi_1'+n\cdot [\Sigma_1]$, and $n=n_{p_2}(\phi_2)$.
 \item $\cU_2=(b,u_2)$ where $u_2\in \cM_{J_2}(S_2^{\qs_2},\phi_2')$ satisfies \ref{M-1}--\ref{M-7}, and $b$ is the disjoint union of $m$ beta-boundary degenerations $b_1,\dots, b_m$ which each have Maslov index $2$. Furthermore, $\phi_2=\phi_2'+m\cdot [\Sigma_2]$, and $m=n_{p_1}(\phi_1)$.
 \item  The collection $\qs_1$ consists of $m$ marked points $(q_{1,1},\dots, q_{1,m})$ and $\qs_2$ consists of $n$ marked points $(q_{2,1},\dots, q_{2,n})$. Furthermore, for all $i$ and $j$,
 \[
 \begin{split}
( \pi_{\Sigma_i}\circ u_i)(q_{i,j})&=p_i,\\
(\pi_{\Sigma_2}\circ u_2)\left((\pi_{\bD}\circ u_2)^{-1}\left(\ev_{q_{1,j}}^{\un}(u_1)\right)\right)
 &=\ev^{\infty}(b_j),\quad \text{and} \\
(\pi_{\Sigma_1}\circ u_1)\left((\pi_{\bD}\circ u_1)^{-1}\left(\ev^{\st}_{q_{2,j}}(u_2)\right)\right)&=\ev^{\infty}(a_j).
\end{split}
 \]
 \end{enumerate}
\end{define}

\begin{figure}[ht!]
	\centering
	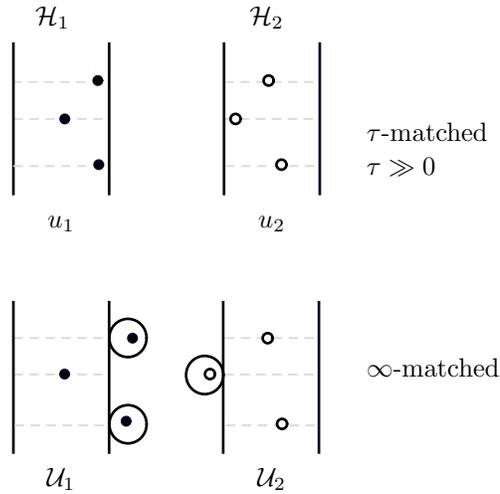
	\caption{ A schematic of a $\tau$-matched curve (top), and a nearby $\infty$-matched curve (bottom). The dots represent the image under of the marked points under the evaluation map to $[0,1]\times \R$.}\label{fig:99}
\end{figure}

A schematic of a $\tau$-matched disk and a nearby $\infty$-matched disk is shown in Figure~\ref{fig:99}. If $(\cU_1,\cU_2)$ is an $\infty$-matched disk, we define the \emph{trimming} of $(\cU_1,\cU_2)$ to be the pair $(u_1,u_2)$.

Given two classes $\phi_1$ and $\phi_2$, on $\cH_1$ and $\cH_2$, respectively, satisfying $n_{p_1}(\phi_1)=n_{p_2}(\phi_2)$, we write $\cM\cM^\infty_{J_1\wedge J_2}(\phi_1,\phi_2)$ for the moduli space of $\infty$-matched disks representing $(\phi_1,\phi_2)$.

We write $M_\otimes$ for the matching condition
\[
M_\otimes:=\left(([0,1]\times \R)^{k_1}\times ([0,1]\times \R)^{k_2}, (0,\dots, 0)\right).
\]
Note that being $M_\otimes$-matched is a vacuous condition: any two curves are $M_\otimes$-matched. Also, $M_\otimes$ depends on the integers $k_1$ and $k_2$, though we suppress this from the notation.

We define a \emph{formal expansion} map
\[
\scE\colon \pi_2(\xs_1,\ys_1)\times \pi_2(\xs_2,\ys_2)\to \pi_2(\xs_1,\ys_1)\times \pi_2(\xs_2,\ys_2),
\]
via the formula
\begin{equation}
\scE(\phi_1,\phi_2)=\left(\phi_1+n_{p_2}(\phi_2)\cdot [\Sigma_1], \phi_2+n_{p_1}(\phi_1)\cdot [\Sigma_2]\right). \label{eq:formal-expansion}
\end{equation}
If $(\phi_1, \phi_2)\in \pi_2(\xs_1,\ys_1)\times \pi_2(\xs_2,\ys_2)$, then 
\[
\ind_{\emb}(\phi_1,\phi_2,M_{\otimes})=\ind_{\emb}(\scE(\phi_1,\phi_2), M_{\#}).
\]

The tensor product differential $\d^\otimes =\d_1\otimes \id+\id\otimes \d_2$ coincides with the map obtained by counting $\# (\cM_{J_1}(\phi_1')\times \cM_{J_2}(\phi_2'))/\R$ for all pairs $(\phi_1',\phi_2')\in \pi_2(\xs_1,\ys_1)\times \pi_2(\xs_2,\ys_2)$ with 
\[
\ind_{\emb}(\phi_1',\phi_2',M_\otimes)=1.
\]
On the other hand, the count of boundary degenerations from Proposition~\ref{prop:relaxed-count} implies that if $\phi_1\in \pi_2(\xs_1,\ys_1)$ and $\phi_2\in \pi_2(\xs_2,\ys_2)$ are classes which have the same multiplicity at the connected sum points, and $\ind_{\emb}(\phi_1,\phi_2,M_\#)=1$, then
\begin{equation}
\# \left( \cM\cM^\infty_{J_1\wedge J_2}(\phi_1,\phi_2)/\R \right)\equiv\sum_{\substack{(\phi_1',\phi_2') \in \pi_2(\xs_1,\ys_1)\times \pi_2(\xs_2,\ys_2) \\ \scE(\phi_1',\phi_2')=(\phi_1,\phi_2)}} \# \left(\cM_{J_1}(\phi_1') \times \cM_{J_2}(\phi_2')\right)/\R. \label{eq:equivalence-trimmed-untrimmed}
\end{equation}

In particular, ~\eqref{eq:equivalence-trimmed-untrimmed} implies that $\d^\otimes$ coincides with the differential obtained by counting $\infty$-matched holomorphic disks representing pairs $(\phi_1',\phi_2')$ with $\ind_{\emb}(\phi_1',\phi_2',M_\#)=1$.

We define a transition map
\[
\Psi^{\tau\to \infty}\colon (\CF^-(\cH_1,\cH_2), \d^{\tau})\to (\CF^-(\cH_1)\otimes \CF^-(\cH_2),\d^\otimes),
\]
as follows. Pick a smooth function $\theta\colon \R\to \R$ such that $\theta'(t)\ge 0$ for all $t\in \R$,   $\theta(t)=\tau$ for all $t\ll 0$, and $\theta(t)\to \infty$ as $t\to \infty$. We define a moduli space $\cM\cM^\theta(\phi_1,\phi_2)$  as in ~\eqref{eq:phi-matched-moduli-space-disks}. The map $\Psi^{\tau\to \infty}$ counts elements of the $\theta$-matched moduli spaces ranging over pairs $(\phi_1,\phi_2)$ satisfying $\ind_{\emb}(\phi_1,\phi_2,M_{\#})=0$ and $n_{p_1}(\phi_1)=n_{p_2}(\phi_2)$. We may similarly define a map $\Psi^{\infty\to \tau}$ in the opposite direction.

\begin{lem}\label{lem:homotopy-equivalence-different-t}
\begin{enumerate}
\item[]
\item  The map $\Psi^{\tau\to \infty}$ is a chain map.
\item The map $\Psi^{\tau\to \infty}$ is a chain homotopy equivalence.
\end{enumerate}
\end{lem}
\begin{proof}  The first is proven by counting the ends of index 1 families of $\theta$-matched flowlines. Using transversality, the ends are constrained into the following 2-story holomorphic combs:
\begin{enumerate}
\item A $\tau$-matched holomorphic disk with $\ind_{\emb}(\phi_1,\phi_2,M_{\#})=1$, and a $\theta$-matched holomorphic disk with $\ind_{\emb}(\phi_1,\phi_2,M_{\#})=0$.
\item\label{comb:2:del-tensor} An $\infty$-matched disk, whose trimming represents a pair $(\phi_1',\phi_2')$ with $\ind_{\emb}(\phi_1',\phi_2',M_\otimes)=1$, and a $\theta$-matched holomorphic disk with $\ind_{\emb}(\phi_1,\phi_2,M_\#)=0$.
\end{enumerate}

We claim that in the latter case, the matched moduli space admits a map to $[0,\infty)$, defined on a small neighborhood of the broken curve, which is locally proper and  odd degree near 0. We sketch the argument, which is very similar to \cite{OSBorderedHFK}*{Section~7.8}. Suppose $(u_1,u_2)$ is such a limiting comb, and suppose further that $u_1$ consists of an index 1 disk with $k$ marked points, and $u_2$ consists of $k$ index 2 boundary degenerations, as well as a constant flowline. Let $c_1<c_2<\cdots<c_k$ denote the $\R$-values of the marked points of $u_1$. Let $V_0$ denote $(c_1-\veps,c_1+\veps)\times \cdots \times (c_k-\veps,c_k+\veps)$, and let $V$ denote the image of $V$ under the diagonal $\R$-action. For some small $\veps>0$, gluing gives a map
\[
\g\colon [0,\veps)^k\times V\to \bar{\cM}(\phi_2),
\]
which is a homeomorphism onto its image. Let $\bar{U}$ denote the image of $\g$. The $(c_i-\veps,c_i+\veps)$ components correspond to the height at which a boundary degeneration is glued, and the $[0,\veps)$ components are the gluing parameters. There is an evaluation map $\bar{\ev}_2\colon \bar{U}\to [0,1)^k\times \R^k$. Similarly there is an evaluation map $\ev_1\colon \cM(\phi_1)\to (0,1)^k\times \R^k$. The $\theta$-matched moduli space is the fibered product of $\cM(\phi_1)$ and $\cM(\phi_2)$ under $\ev_1$ and $R^\theta\circ \ev_2$.

If $q_i$ is an arbitrarily chosen marked point of $u_1$, then the map $\ev_\R\colon \cM(\phi_1)\to \R$ given by $u\mapsto \pi_{\R}(u(q_i))$ restricts to a diffeomorphism on a neighborhood of $u_1$. Furthermore, we may view the map as extending to a continuous map $\bar{\ev}_{\R}\colon \bar{\cM}(\phi_1)\to [-\infty,\infty]$. We claim that on the $\theta$-matched moduli spaces, the map $\bar{\ev}_{\R}$ is locally proper and odd degree near $\infty$. To see this, we note that a small modification of \cite{LOTBordered}*{Lemma~5.38 (2)} implies that the map $\bar{\ev}_2$ is proper and odd degree near $\ve{0}\times c_1\times \cdots \times c_k$. In particular, if $t\in \R$ is fixed, and $u_t:=\ev_\R^{-1}(t)\in \cM(\phi_1)$, and $R^{-\theta}\ev_1(u_t)$ is a regular value of $\ev_2$, then there are an odd number of curves in the $\theta$-matched moduli space which have $\ev_{\R}=t$. It follows that $\ev_{\R}$ is locally odd degree near $\infty$. In particular, combs satisfying~\eqref{comb:2:del-tensor} make odd algebraic contribution, so  
\[
\Psi^{\tau\to \infty}\circ \d^{\tau}+\d^{\otimes} \circ \Psi^{\tau\to \infty}=0,
\]
where we are identifying $\d^{\otimes}$ with the map which counts $\infty$-matched holomorphic disks, using~\eqref{eq:equivalence-trimmed-untrimmed}.

The second claim follows from the same reasoning as the corresponding claim of Lemma~\ref{lem:t-matched-complexes}. 
\end{proof}

\subsection{Deformed stabilizations}
\label{sec:deformed-stabilizations-disks}
In our construction of $H_E$, we will encounter some additional variations of the matching deformation considered in Section~\ref{sec:time-dilation-complexes}.  Suppose $\cH_1$ and $\cH_2$ are two Heegaard diagrams with chosen connected sum points, and $v$ is an $\R$-invariant vector field on $[0,1]\times \R$ which vanishes on $\{0,1\}\times \R$.  We may define a $(v,\tau)$-deformed  differential $\d^{v,\tau}_{\cH_1\wedge \cH_2}$, by counting holomorphic disks which are matched by the time $\tau$ flow of $v$.

We will need the following fact:

\begin{lem}\label{lem:Y-matched-differential}
 Suppose that $\cH=(\Sigma,\as,\bs,w)$ is a Heegaard diagram with a special point $p\in \Sigma\setminus (\as\cup \bs)$, and $v$ is a vector field on $[0,1]\times \R$, as above.
 \begin{enumerate}
 \item Suppose $\cH_0=(\bT^2,\a,\b)$ is a standard diagram for $S^3$ or $L(m,1)$ and  $x\in \a\cap \b$. If $\tau\in \R$ is generic, then
 \[
\d_{\cH\wedge \cH_0}^{v,\tau}(\xs\times x)=\d_{\cH}(\xs)\otimes x.
 \]
 Furthermore, if $\phi$ and $\phi_0$ are a pair of classes of disks such that $\ind_{\emb}(\phi,\phi_0,M_\#)=0$, and
 \[
\cM\cM^{v,\tau\in [\tau_1,\tau_2]}(\phi,\phi_0)\neq\emptyset,
 \]
 then $\phi$ and $\phi_0$ represent constant classes.
 \item Suppose $\cH_0=(\bT^2,\a,\b)$ is a standard diagram for $S^1\times S^2$. Write $\{\theta^+,\theta^-\}=\a\cap \b$. Then, for generic $\tau \in \R$, 
 \[
 \begin{split}
\d_{\cH\wedge \cH_0}^{v,\tau}(\xs\times \theta^+)&=\d_{\cH}(\xs)\otimes \theta^+, \quad \text{and}\\
\d_{\cH\wedge \cH_0}^{v,\tau}(\xs\times \theta^-)&=\d_{\cH}(\xs)\otimes \theta^-+\sum_{\ys\in \bT_{\a}\cap \bT_{\b}} C_{\xs,\ys}\cdot  \ys\times \theta^+,
\end{split}
 \]
 for some $C_{\xs,\ys}\in \bF[U]$. Furthermore, suppose $\phi$ and $\phi_0$ are classes of disks, such that $\ind_{\emb}(\phi,\phi_0, M_{\#})=0$. Write $\phi_0\in \pi_2(\theta_1,\theta_2)$, for $\theta_1,\theta_2\in \a\cap \b$. If $\gr(\theta_1,\theta_2)\ge 0$, and
 \[
\cM\cM^{v,\tau\in [\tau_1,\tau_2]}(\phi,\phi_0)\neq \emptyset,
 \]
 then $\phi$ and $\phi_0$ both represent  constant classes.
 \end{enumerate}
\end{lem}
\begin{proof} The computation of the differential follows from the same argument as the ordinary proof of stabilization invariance \cite{OSDisks}*{Section~10}. See Lemma~\ref{lem:lens-space-stabilization-disks} for details in our present context.

We now consider the claim about index 0 classes. We focus on the claim when $\cH_0$ is a diagram for $S^1\times S^2$, since the other cases follow a similar argument. If $\phi_0\in \pi_2(\theta_1,\theta_2)$ is a class on $\cH_0$, then 
\[
\mu(\phi_0)=2n_{p_0}(\phi_0)+\gr(\theta_1,\theta_2).
\]
 Hence, if $\phi$ is a class on $\cH$, then 
\[
\ind_{\emb}(\phi, \phi_0,M_{\#})=\mu(\phi)+\mu(\phi_0)-2n_{p_0}(\phi_0)=\mu(\phi)+\gr(\theta_1,\theta_2).
\]
If $\cM\cM^{v,\tau\in [\tau_1,\tau_2]}(\phi,\phi_0)\neq \emptyset$, then $\mu(\phi)\ge 0$ by transversality. Hence, if $\gr(\theta_1,\theta_2)\ge 0$ we must have $\mu(\phi)=0$. By transversality, $\phi$ must be the constant class. It is straightforward to see that this also implies that $\phi_0$ is a constant class.
\end{proof}

\subsection{Deformed stabilizations and holomorphic triangles}
\label{sec:deformed-stabilizations}

We need several results about deformed stabilizations and holomorphic triangles, similar to those stated in Section~\ref{sec:deformed-stabilizations-disks} for disks.

We now consider deforming the triangle maps using a smooth vector field $v$ on $\Delta$ which vanishes on $\d \Delta$ and is $\R$-invariant on the cylindrical ends of $\Delta$. (We only have need to consider the vector fields $v_{\Delta,1}$ and $v_{\Delta,2}$ from Figure~\ref{fig:140}). Suppose that $\cT=(\Sigma,\as,\bs,\gs,w)$ is a Heegaard triple with a distinguished point $p\in \Sigma\setminus (\as\cup \bs\cup \gs)$, and suppose that $\cT_0=(\bT^2,\a,\b,\a',p_0)$ is a genus 1 lens space stabilizing triple, as in Section~\ref{sec:lens-space-stabilizations}. If $\psi$ and $\psi_0$ are two classes of triangles on $\cT$ and $\cT_0$, then we may define their $\tau$-matched moduli space $\cM\cM^\tau(S^{\qs}, S_0^{\qs_0},\psi,\psi_0)$ by adapting ~\eqref{eq:t-matched-moduli-space-disks}.

We define a matched triangle map $f^\tau_{\cT \wedge \cT_0}$ by counting $\tau$-matched triangles representing pairs $(\psi_1,\psi_2)$ with $\ind_{\emb}(\psi_1,\psi_2,M_\#)=0$. Given  $\tau_1,\tau_2\in \R$ with $\tau_1<\tau_2$, we may also count representatives of pairs $(\psi,\psi_0)$ with $\ind_{\emb}(\psi,\psi_0)= -1$, which are $\tau$-matched for some $\tau\in [\tau_1,\tau_2]$. We write $f_{\cT\wedge \cT_0}^{\tau_1\to \tau_2}$ for this map.

 We have the following analog of Lemma~\ref{lem:lens-stabilization-triangles}:

\begin{lem}\label{lem:t-matched-stabilizations} Suppose that $\cT=(\Sigma,\as,\bs,\gs,w)$ is a Heegaard triple with a distinguished point $p\in \Sigma\setminus (\as\cup \bs\cup \gs)$, and $\cT_0=(\bT^2,\a,\b,\a',p_0)$ is a triple such that  $(\bT^2,\a,\b)$ is a standard diagram for $L(m,1)$, and $\a'$ is a small Hamiltonian translate of $\a$ satisfying $|\a\cap \a'|=2$. See Figure~\ref{fig:179}. 
\begin{enumerate}
\item If $\tau\in \R$ is generic $c_{\a,\b}\in \a\cap \b$ and $c_{\b,\a'}\in \b\cap \a'$, then
\[
f_{\cT\wedge \cT_0}^\tau(\xs\times c_{\a,\b}, \ys\times c_{\b,\a'})=f_{\cT}(\xs,\ys)\otimes \Theta_{\a,\a'}^-+\sum_{\zs\in \bT_{\a}\cap \bT_{\g}} C_{\xs,\ys,\zs}^\tau\cdot \zs\otimes \Theta_{\a,\a'}^+,
\]
for some $C_{\xs,\ys,\zs}^\tau\in \bF[U]$.
\item\label{lem:t-matched-stabilizations-triangles-part-2} If $\tau_1,\tau_2\in \R$ are generic, and $c_{\a,\b}\in \a\cap \b$ and $c_{\b,\a'}\in \b\cap \a'$, then
\[
f_{\cT\wedge \cT_0}^{\tau_1\to \tau_2}(\xs\times c_{\a,\b} ,\ys\times c_{\b,\a'})=
\begin{cases}
\displaystyle
\sum_{\zs\in \bT_{\a}\cap \bT_{\g}} C^{[\tau_1,\tau_2]}_{\xs,\ys,\zs} \cdot \zs\otimes \Theta^+_{\a,\a'} &\text{ if } \frs_{p_0}(c_{\a,\b})=\frs_{p_0}(c_{\b,\a'}),\\
0& \text{ otherwise},
\end{cases}
\]
for some $C_{\xs,\ys,\zs}^{[\tau_1,\tau_2]}\in \bF[U]$. 
\end{enumerate}
\end{lem}
\begin{proof} The first claim is proven no differently than Lemma~\ref{lem:lens-stabilization-triangles}, using ~\eqref{eq:matched-disks-stabilization}.

The second claim is proven as follows. Suppose $\psi$ is a class on $(\Sigma,\as,\bs,\gs)$ and $\psi_0\in \pi_2(c_{\a,\b}, c_{\b,\a'}, \Theta)$ is a class on $(\bT^2,\a,\b,\a')$, where $\Theta\in \a\cap \a'$.   Adapting ~\eqref{eq:index-connected-sum-stabilization} to our present situation, we obtain
\begin{equation}
\ind_{\emb}(\psi,\psi_0,M_\#)=\mu(\psi)+\gr(\Theta_{\a,\a'}^+,\Theta)-1.\label{eq:index-tau-matched-triangle-stabilization}
\end{equation}
Note that any holomorphic curve counted by $f_{\cT\wedge \cT_0}^{\tau_1\to \tau_2}$ contains an ordinary holomorphic triangle on $\cT$. In particular, if  $(\psi,\psi_0)$ is a pair of classes which admit a $\tau$-matched representative for some $\tau$, then $\mu(\psi)\ge 0$, by transversality for curves on $\cT$. Combined with ~\eqref{eq:index-tau-matched-triangle-stabilization}, we see that if $\Theta=\Theta_{\a,\a'}^-$, and $(\psi,\psi_0)$ has a $\tau$-matched index $-1$ representative, then $\mu(\psi)=-1$, which violates transversality. Hence we may only have representatives if $\Theta=\Theta_{\a,\a'}^+$, proving the stated formula.
\end{proof}

\subsection{Deforming $\exc$}

In this section, we describe our null-homotopy of $E$, by deforming the holomorphic quadrilateral maps using the families of vector fields $v_1,$ $v_2$ and $v_3$ on $\Box$ described in Section~\ref{sec:vector-fields}.

The time $\tau$ flow of $v_{i,s}$ gives a 2-parameter family of diffeomorphisms,
\[
R_{i,s}^{\tau}\colon \Box\to \Box,
\]
for each $i\in \{1,2,3\}$.

 Suppose $v_i$ and $j_i$ satisfy the analogs of~\eqref{v-1-1}--\eqref{v-1-4}. If $J\wedge I=(J_s\wedge I_s)_{s\in [0,1]}$ is a stratified family of almost complex structures on $(\Sigma\wedge \bT^2)\times \Box$, which covers a family $j_i$ on $\Box$, we define the $(v_i,\tau)$-matched moduli space of rectangles as follows:
\begin{equation}
\begin{split}
&\cM\cM^{v_i,\tau}_{I\wedge J}(S_0^{\qs_0},S^{\qs},\psi_0 ,\psi)\\
:=&\left\{(u,s,u_0,s_0)\in \cM_{J}(S^{\qs},\psi)\times \cM_{I}(S_0^{\qs_0}, \psi_0)\middle \vert \begin{array}{c} 
(\pi_{\Sigma}\circ u)(q_j)=p ,\\
(\pi_{\bT^2}\circ u_0)(q_{0,j})=p_0,\\
R^{\tau}_{i,s}((\pi_{\Box}\circ u_0)(q_{0,j})) = (\pi_{\Box}\circ u)(q_{j}),\\
 \text{ for } j=1,\dots, k,\\
 s=s_0
 \end{array}\right\}.
 \end{split}
 \label{eq:tau-matched-rectangles}
\end{equation}

We define maps $h_1^\tau$, $h_2^{\tau}$ and $h_3^{\tau}$ with the same domain and range as $h_1$, $h_2$ and $h_3$, by counting $(v_i,\tau)$-matched quadrilaterals on $\cQ_1$, $\cQ_2$ and $\cQ_3$. We define
\[
E^\tau:= F_3\circ h_1^\tau+F_3 \circ h_2^\tau+F_3 \circ h_3^\tau.
\]

\begin{lem}\label{lem:Etau-chain-map}
 For generic $\tau\ge 0$, the map $E^\tau$ is a chain map.
\end{lem}
\begin{proof}
 The proof follows by counting the ends of index 0, $(v_i,\tau)$-matched rectangles, representing classes $(\psi, \psi_0)$ in 
 \begin{equation}\begin{split}
 \pi_2(\Theta_{\a_3,\a_2}^+, \Theta_{\a_2,\a_1}^{\can},\xs\times x, \ys\times \theta^-), \quad \pi_2(\Theta_{\a_2,\a_1}^{\can},\xs\times x, \Theta_{\b_3,\b_2}^+, \ys\times \theta^-),\\
 \text{and} \quad \pi_2(\xs\times x, \Theta_{\b_3,\b_2}^+, \Theta_{\b_2,\b_1}^{\can}, \ys\times \theta^-).\hspace{2cm}
 \end{split} \label{eq:classes-for-E}
 \end{equation}
	Here, $\{x\}=c_1\cap c_3'$, and $\theta^-$ denotes the bottom generator of either $c_1\cap c_1'$, $c_2\cap c_2'$ or $c_3\cap c_3'$. 
	
We claim that for generic, fixed $\tau\in [0,\infty)$, the moduli spaces of index $0$, $(v_1,\tau)$-matched quadrilaterals on $\cQ_1$ have ends corresponding exactly to the following 2-story combs:
 \begin{enumerate}[label=($h_1$-\arabic*), ref=$h_1$-\arabic*, leftmargin=*,widest=IIII]
 \item\label{deg:h1-1} $s\in (0,1)$:
 \begin{itemize}
 \item An index $1$, $(v,\tau)$-matched  holomorphic disk in one of the four cylindrical ends, where $v\in \{v_H,v_Y\}$, and
 \item  a $(v_{1,s},\tau)$-matched holomorphic rectangle on $\cQ_1$ of index $-1$.
 \end{itemize}
 \item \label{deg:h1-3} $s= 0$:
 \begin{itemize}
 \item A $(v_{1,0},\tau)$-matched  triangle on $(\Sigma\wedge \bT^2,\as_3,\as_1,\bs_3)$ of index 0, and
 \item  a $(v_{1,0},\tau)$-matched triangle on $(\Sigma\wedge \bT^2,\as_3,\as_2,\as_1)$ of index 0.
 \end{itemize}
 \item\label{deg:h1-4}  $s= 1$:
 \begin{itemize}
 \item A $(v_{1,1},\tau)$-matched triangle on $(\Sigma\wedge \bT^2,\as_3,\as_2,\bs_3)$ of index 0, and
 \item a $(v_{1,1},\tau)$-matched triangle on $(\Sigma\wedge \bT^2,\as_2,\as_1,\bs_3)$ of  index 0.
 \end{itemize}
 \end{enumerate}
 
Similarly, the ends of the moduli spaces of $(v_2,\tau)$-matched quadrilaterals on $\cQ_2$ have one of the following forms: 
 \begin{enumerate}[label=($h_2$-\arabic*), ref=$h_2$-\arabic*, leftmargin=*,widest=IIII]
 \item\label{deg:h2-1}$s\in (0,1)$:
 \begin{itemize}
 \item  An index $1$, $(v,\tau)$-matched holomorphic disk in one of the 4 cylindrical ends, where $v\in \{v_H,v_Y\}$, and 
 \item a $(v_{2,s},\tau)$-matched holomorphic rectangle on $\cQ_2$ of index $-1$.
 \end{itemize}
 \item \label{deg:h2-3} $s=0$: 
 \begin{itemize}
 \item  A $(v_{2,0},\tau)$-matched triangle on $(\Sigma\wedge \bT^2,\as_2,\bs_3,\bs_2)$  of index 0, and 
  \item a $(v_{2,0},\tau)$-matched triangle $(\Sigma\wedge \bT^2,\as_2,\as_1,\bs_3)$ of index 0.
 \end{itemize}
 \item\label{deg:h2-4}  $s=1$:
 \begin{itemize}
 \item A $(v_{2,1},\tau)$-matched  triangle on $(\Sigma\wedge \bT^2,\as_2,\as_1,\bs_2)$ of index 0, and 
 \item a $\tau$-matched triangle on $(\Sigma\wedge \bT^2,\as_1,\bs_3,\bs_2)$ of  index 0.
 \end{itemize}
 \end{enumerate}
 
 Finally, the ends of the moduli spaces of $(v_3,\tau)$-matched quadrilaterals on $\cQ_3$ take one of the following forms:
  \begin{enumerate}[label=($h_3$-\arabic*), ref=$h_3$-\arabic*, leftmargin=*,widest=IIII]
 \item\label{deg:h3-1}  $s\in (0,1)$:
 \begin{itemize}
	\item  An index $1$, $(v,\tau)$-matched holomorphic disk in one of the 4 cylindrical ends, where $v\in \{v_H,v_Y\}$, and
  	\item a $(v_{3,s},\tau)$-matched holomorphic rectangle on $\cQ_3$  of index $-1$.
  \end{itemize}
 \item \label{deg:h3-3} $s= 0$:
 \begin{itemize}
 \item A $(v_{3,0},\tau)$-matched triangle on  $(\Sigma\wedge \bT^2,\as_1,\bs_3,\bs_2)$ of index 0, and
 \item  a $(v_{3,0},\tau)$-matched triangle on $(\Sigma\wedge \bT^2,\as_1,\bs_2,\bs_1)$ of index 0.
 \end{itemize}
 \item\label{deg:h3-4}  $s= 1$:
 \begin{itemize}
 \item A $(v_{3,1},\tau)$-matched triangle on $(\Sigma\wedge \bT^2,\as_1,\bs_3,\bs_1)$ of index 0, and
 \item A $(v_{3,1},\tau)$-matched triangle on  $(\Sigma\wedge \bT^2,\bs_3,\bs_2,\bs_1)$ of index 0.
 \end{itemize}
 \end{enumerate}

Furthermore, an argument similar to the one in Lemma~\ref{lem:homotopy-equivalence-different-t} implies that each broken curve labeled $(h_i$-$j)$ makes odd algebraic contribution to the count of ends of its moduli space. 

The ends of the moduli spaces in ~\eqref{eq:classes-for-E}, of type \eqref{deg:h1-1}, \eqref{deg:h2-1} and \eqref{deg:h3-1}, where additionally the index 1 disk occurs in the $(\as_1,\bs_3)$ end or the $(\as_i,\bs_i)$ end, give the $\ys$ component of
\[
(\d E^\tau+E^\tau \d)(\xs).
\]
The remaining ends of type ~\eqref{deg:h1-1},  \eqref{deg:h2-1} and \eqref{deg:h3-1} have a disk component which is $\pm (v_Y,\tau)$-matched and occurs in one of the ends labeled $(\as_2,\as_1)$, $(\as_3,\as_2)$, $(\bs_3,\bs_2)$ or $(\bs_2,\bs_1)$. The total count of these ends is zero, since $\Theta_{\a_2,\a_1}^{\can}$, $\Theta_{\a_3,\a_2}^+$, $\Theta_{\b_1,\b_2}^{\can}$ and $\Theta_{\b_2,\b_3}^+$ are cycles in the $\pm (v_Y,\tau)$-matched complexes by Lemma~\ref{lem:Y-matched-differential}.

Next, the moduli spaces \eqref{deg:h1-3} and \eqref{deg:h3-4} make total contribution zero, by Lemma~\ref{lem:model-computations-m-surgery} and a $\tau$-matched version of Proposition~\ref{prop:stabilize-triangles-general}.

We now claim that the ends labeled \eqref{deg:h1-4} cancel the ends labeled \eqref{deg:h2-3}. This follows since both ends feature a count of triangles on $(\as_2,\as_1,\bs_3)$, and a count on either $(\as_3,\as_2,\bs_3)$ or $(\as_2,\as_1,\bs_2)$. The counts of triangles on $(\Sigma\wedge \bT^2,\as_3,\as_2,\bs_3)$ and $(\Sigma\wedge \bT^2,\as_2,\as_1,\bs_2)$ coincide, since we may destabilize using Lemma~\ref{lem:lens-stabilization-triangles}, and then identify both destabilized counts with nearest point maps by Proposition~\ref{prop:nearest-point-triangles}. Hence, the  counts of the ends ~\eqref{deg:h1-4} cancel those of \eqref{deg:h2-3}, as claimed. The ends \eqref{deg:h2-4} cancel the ends of \eqref{deg:h3-3} in a similar manner.

Summing over all of the above ends, we obtain $\d E^\tau+E^\tau \d=0$, completing the proof.
\end{proof}

We now describe a candidate limiting map, $E^\infty$, which has the same domain and range as $E^\tau$. We define
\[
E^\infty:=F_3\circ h_1^\infty+F_3\circ h_2^{\infty}+F_3\circ h_3^{\infty},
\]
where $h_i^\infty$ counts \emph{aspect ratio matched}, holomorphic quadrilaterals on $\cQ_i$  of index $-1$, i.e. $E^\infty$ counts elements of the moduli spaces
\[
\cM\cM^{\infty,\ar}_{J\wedge I}(\psi,\psi_0):=\{(u,s,u_0,s_0)\in \cM_{J}(\psi)\times \cM_{I}(\psi_0):s=s_0\},
\]
ranging over pairs with
\[
\ind_{\emb}(\psi,\psi_0,M_\otimes):=\mu(\psi)+\mu(\psi_0)=-1.
\]

Note that in the definition of $\cM\cM_{J\wedge I}^{\infty,\ar}(\psi,\psi_0)$, we do not require $n_p(\psi)=n_{p_0}(\psi_0)$, and also there is no matching condition at the connected sum points.

In the map $h_i^\infty$, an aspect ratio matched pair $(u,u_0)$ representing $(\psi,\psi_0)$ is counted with a factor of
\begin{equation}
U^{n_{p_0}(\psi_0)+n_{w}(\psi)}T^{n_{p}(\psi)+n_{z_0}(\psi_0)-n_{p_0}(\psi_0)-n_{w}(\psi)}. \label{eq:weighting-trimmed-curves}
\end{equation}

The weight in~\eqref{eq:weighting-trimmed-curves} is the ordinary weight of the formal expansion $\scE(\psi,\psi_0)$, as in ~\eqref{eq:formal-expansion}. 

We will show the following:

\begin{prop}
\label{prop:deform-tau-to-infinity} If $\tau$ is generic, then $E^{\tau}\simeq E^\infty$.
\end{prop}

The proof of Proposition~\ref{prop:deform-tau-to-infinity} is technical, and is delayed until Section~\ref{sec:Et=E-infty}. The proof is in the same spirit as the proof of Lemma~\ref{lem:homotopy-equivalence-different-t}. Before embarking upon the proof of Proposition~\ref{prop:deform-tau-to-infinity}, we prove the following, which, together with Proposition~\ref{prop:deform-tau-to-infinity}, implies Proposition~\ref{prop:Esimeq0}.

\begin{lem} \label{lem:E-infty-null-homotopic}
The map $E^\infty$ satisfies $E^\infty\simeq 0.$
\end{lem}
\begin{proof} The proof follows by deforming the aspect ratio matching. For $T\in \R$, we consider the deformed moduli spaces
\[
\cM\cM^{\infty,\ar,T}_{J\wedge I}(\psi,\psi_0):=\left\{(u,s,u_0,s_0)\in \cM_{J}(\psi)\times \cM_{I}(\psi_0): \rho_T(s)=s_0\right\},
\]
where $\rho_T\colon [0,1]\to [0,1]$ denotes the time $T$ flow of the vector field from ~\eqref{eq:def-vector-field-v}.

We consider the ends of the parametrized moduli spaces
\begin{equation}
\cM\cM^{\infty, \ar, T\in [0,\infty)} (\psi,\psi_0):=\bigcup_{T\in [0,\infty)}\cM\cM^{\infty,\ar,T}_{J\wedge I}(\psi,\psi_0)\times \{T\},
\label{eq:infinity-and-aspect-ratio-matched}
\end{equation}
ranging over the classes in~\eqref{eq:classes-for-E}, with $\ind_{\emb}(\psi,\psi_0,M_\otimes)=-1$.

Firstly, there are ends at $T=0$, which contribute
\[
E^\infty(\xs).
\]

We claim that any end at $T\in (0,\infty)$ must have $s\in (0,1)$. Indeed, an end with $T\in (0,\infty)$ and $s\in \{0,1\}$ would necessarily contain four holomorphic triangles, as well as potentially a collection of holomorphic disks and boundary degenerations. By transversality, every holomorphic triangle, disk or boundary degeneration must have nonnegative Maslov index. However, the total Maslov index must be $-1$, since such a curve appears in the boundary of the moduli space of Maslov index $-1$ rectangles, leading to a contradiction.

Hence, the ends at finite $T$ are restricted to consist of holomorphic strips breaking off into the 4 cylindrical ends. Furthermore, strips breaking off into the ends labeled $(\as_3,\as_2)$, $(\as_2,\as_1)$, $(\bs_3,\bs_2)$ or $(\bs_2,\bs_1)$ cancel modulo 2, since $\Theta_{\a_3,\a_2}^+$, $\Theta_{\a_2,\a_1}^{\can}$, $\Theta_{\b_3,\b_2}^+$ and $\Theta_{\b_2,\b_1}^{\can}$ are cycles. Hence,  the ends at finite $T$ contribute
\[
(H\circ \d+\d \circ H)(\xs),
\]
where $H$ counts index $-2$ curves $(u,s,u_0,s_0)$, which satisfy $\rho_T(s)=s_0$ for some $T\in (0,\infty)$.

Finally, it remains to count the ends appearing as $T\to \infty$. By transversality, these ends consist of two holomorphic triangles, and one holomorphic rectangle. Furthermore, since there is no matching condition, the triangles must have index 0, and the rectangle index $-1$.

 More concretely, the broken curves appearing on $\cQ_1$ as $T\to \infty$ have one of the following configurations:
\begin{enumerate}[label=($h_1^{\infty,\infty}$-\arabic*), ref=$h_1^{\infty,\infty}$-\arabic*, leftmargin=*,widest=IIIIII]
\item $s\in (0,1),$ $T\to \infty$:
\begin{itemize}
\item An index 0 triangle on $(\bT^2,c_3,c_1,c_3')$,
\item an index 0 triangle on $(\bT^2,c_3,c_2,c_1)$, and
\item an index $-1$ rectangle on $(\Sigma,\as'',\as',\as,\bs)$.
\end{itemize}
\item $s=0$, $T\to \infty$:
\begin{itemize}
\item  An index $-1$ rectangle on $(\bT^2,c_3,c_2,c_1,c_3')$, 
\item an index 0 triangle on $(\Sigma,\as'',\as,\bs)$, and
\item  an index 0 triangle on $(\Sigma,\as'',\as',\as)$.
\end{itemize}
\end{enumerate}
Similarly, the broken curves arising on $\cQ_2$ have one of the following configurations:
\begin{enumerate}[label=($h_2^{\infty,\infty}$-\arabic*), ref=$h_2^{\infty,\infty}$-\arabic*, leftmargin=*,widest=IIIIII]
\item $s\in (0,1),$ $T\to \infty$:
\begin{itemize}
\item An index 0 triangle on $(\bT^2,c_1,c_3',c_2')$,
\item  an index 0 triangle on $(\bT^2,c_2,c_1,c_2')$, and
\item  an index $-1$ rectangle on $(\Sigma,\as',\as,\bs,\bs')$.
\end{itemize}
\item $s=0,$ $T\to \infty$: 
\begin{itemize}
\item An index $-1$ rectangle on $(\bT^2,c_2,c_1,c_3',c_2')$,
\item  an index 0 triangle on $(\Sigma,\as',\as,\bs)$, and
\item  an index 0 triangle on $(\Sigma,\as,\bs,\bs')$.
\end{itemize}
\end{enumerate}
Finally, the curves appearing on $\cQ_3$ consist of the following:
\begin{enumerate}[label=($h_3^{\infty,\infty}$-\arabic*), ref=$h_3^{\infty,\infty}$-\arabic*, leftmargin=*,widest=IIIIII]
\item $s\in (0,1),$ $T\to \infty$:
\begin{itemize}
\item An index 0 triangle on $(\bT^2,c_1,c_3',c_1')$,
\item  an index 0 triangle on $(\bT^2,c_3',c_2',c_1')$, and
\item an index $-1$ rectangle on $(\Sigma,\as,\bs,\bs',\bs'')$.
\end{itemize}
\item $s=0$, $T\to \infty$:
\begin{itemize}
\item  An index $-1$ rectangle on $(\bT^2,c_1,c_3',c_2',c_1')$,
\item  an index 0 triangle on $(\Sigma,\as,\bs,\bs')$, and
\item an index 0 triangle on $(\Sigma,\as,\bs',\bs'')$.
\end{itemize}
\end{enumerate}
Furthermore, we claim that each of the above broken curves labeled $(h_i^{\infty,\infty}$-$j$) makes odd algebraic contribution to the count of ends of the moduli spaces in~\eqref{eq:infinity-and-aspect-ratio-matched}. By this, we mean that in a neighborhood of a broken curve, the matched moduli spaces admit a continuous, proper map to $[0,\infty)$, with the preimage of $0$ being the broken curve, which is odd degree near $0$; see Definition~\ref{def:locally-odd-degree}.

The total count of curves satisfying each of the configurations (weighted by powers of $U$ and $T$), ranging over classes in ~\eqref{eq:classes-for-E},  vanishes modulo 2, by Lemmas~\ref{lem:genus-1-quadruple-count} and~\ref{lem:model-computations-m-surgery}.

Counting up all ends, we obtain $E^\infty=[\d,H]$, completing the proof.
\end{proof}

\subsection{Proof of Proposition~\ref{prop:deform-tau-to-infinity}}
\label{sec:Et=E-infty}

As a warmup to proving Proposition~\ref{prop:deform-tau-to-infinity}, we prove the following:

\begin{lem}\label{lem:deform-tau}
 If $\tau_1,\tau_2\in \R$ are generic, then  $E^{\tau_1}\simeq E^{\tau_2}$.
\end{lem}
\begin{proof}
We count the ends of the parametrized moduli spaces 
\[
\cM\cM^{v_i,\tau\in [\tau_1,\tau_2]}(\psi,\psi_0):=\bigcup_{\tau\in [\tau_1,\tau_2]} \cM\cM^{v_i,\tau}(\psi,\psi_0)\times \{\tau\},
\]
 ranging over pairs of classes $(\psi,\psi_0)$ as in ~\eqref{eq:classes-for-E}, which have index $-1$. These parametrized moduli spaces have ends similar to those described in Lemma~\ref{lem:Etau-chain-map}. For example, in addition to the ends at $\tau_1$ and $\tau_2$, the moduli spaces on $\cQ_1=(\Sigma \wedge \bT^2,\as_3,\as_2,\as_1,\bs_3)$ have the following ends:
 \begin{enumerate}[label=($h_1^{[\tau_1,\tau_2]}$-\arabic*), ref=$h_1^{[\tau_1,\tau_2]}$-\arabic*, leftmargin=*,widest=IIIIII]
 \item\label{deg:h1-1-family} $s\in (0,1)$, $\tau\in (\tau_1,\tau_2)$:
 \begin{itemize}
 \item An index $1$, $\tau$-matched holomorphic disk in one of the four cylindrical ends, and
 \item an index $-2$, $(v_{1,s},\tau)$-matched holomorphic quadrilateral.
 \end{itemize}
 \item\label{deg:th1-2-family}  $s\in (0,1)$, $\tau\in (\tau_1,\tau_2)$:
 \begin{itemize}
 \item An index 0, $\tau$-matched disk in the four cylindrical ends, and
 \item an index $-1$, $(v_{1,s},\tau)$-matched holomorphic quadrilateral.
 \end{itemize}
 \item \label{deg:h1-3-family} $s= 0$, $\tau\in (\tau_1,\tau_2)$:
 \begin{itemize}
 \item A $(v_{1,0},\tau)$-matched triangle on $(\Sigma\wedge \bT^2,\as_3,\as_1,\bs_3)$  of some index $k\in \{-1,0\}$, and
 \item  a $(v_{1,0},\tau)$-matched triangle on  $(\Sigma \wedge \bT^2,\as_3,\as_2,\as_1)$ of index $-1-k$.
 \end{itemize}
 \item\label{deg:h1-4-family}  $s= 1$, $\tau\in (\tau_1,\tau_2)$:
 \begin{itemize}
 \item A $(v_{1,1},\tau)$-matched triangle on $(\Sigma\wedge \bT^2,\as_3,\as_2,\bs_3)$ of index $k\in \{-1,0\}$, and 
 \item a $(v_{1,1},\tau)$-matched triangle on $(\Sigma \wedge \bT^2,\as_2,\as_1,\bs_3)$ of some index $-1-k$.
 \end{itemize}
 \end{enumerate}
Furthermore, a gluing argument similar to Lemma~\ref{lem:homotopy-equivalence-different-t} implies that each of the ends labeled $(h_1^{[\tau_1,\tau_2]}$-$j)$ makes odd algebraic contribution to the ends of its moduli space.

The only ends which did not have a counterpart in Lemma~\ref{lem:Etau-chain-map} are the ends labeled~\eqref{deg:th1-2-family} above. Note that simple index considerations do not prohibit a similar end consisting of an index $-1$ rectangle, and a finite collection of index 0 disks. Nonetheless, for generically chosen almost complex structures, the set of $\tau$ which support a $\tau$-matched, index 0 disk is discrete, and hence in each end where an index 0 disk appears, we may assume that just one appears. Lemma~\ref{lem:Y-matched-differential} implies that the count of the ends of type~\eqref{deg:th1-2-family}, for curves representing the classes shown in ~\eqref{eq:classes-for-E}, vanishes. 

Part~\eqref{lem:t-matched-stabilizations-triangles-part-2} of Lemma~\ref{lem:t-matched-stabilizations} implies that $k=0$ in the ends of type~\eqref{deg:h1-3-family} and ~\eqref{deg:h1-4-family}. Furthermore,  we claim that there are no ends with the configuration~\eqref{deg:h1-3-family}. By the previous reasoning, such an end would contain a $(v_{1,0},\tau)$-matched triangle on $(\Sigma \wedge \bT^2,\as_3,\as_2,\as_1)$ with index $-1-k$. If $\psi\# \psi_0$ is the underlying homology class, this implies that
\[
\ind(\psi\# \psi_0, M_{\#})=\mu(\psi)+\mu(\psi_0)-2 n_p(\psi)=-1.
\]
Write $\psi\in \pi_2(\Theta_{\a'',\a'}^+, \Theta_{\a',\a}^+,\zs)$. By Lemma~\ref{lem:Maslov-triangle-multi-stabilization} we have $\mu(\psi)=2n_p(\psi)+\gr(\Theta_{\a'',\a}^+,\zs)$, and hence
\[
\mu(\psi_0)+\gr(\Theta_{\a'',\a}^+,\zs)=-1.
\]
In particular, $\mu(\psi_0)\le -1$. However a $(v_{1,0},\tau)$-matched representative of $\psi\# \psi_0$ contains an ordinary representative of $\psi_0$, which would violate transversality for curves on $(\bT^2,c_3,c_2,c_1)$ with respect to a fixed almost complex structure.

 The ends of the moduli spaces $\cM\cM^{v_i,\tau\in [\tau_1,\tau_2]}(\psi,\psi_0)$ on $\cQ_2=(\Sigma\wedge \bT^2,\as_2,\as_1,\bs_3,\bs_2)$ and $\cQ_3=(\Sigma\wedge \bT^2,\as_1,\bs_3,\bs_2,\bs_1)$ have a similar description. A cancellation pattern occurs similar to how, in Lemma~\ref{lem:Etau-chain-map}, the ends \eqref{deg:h1-4}  canceled the ends \eqref{deg:h2-3}, and the ends labeled \eqref{deg:h2-4} canceled the ends labeled \eqref{deg:h3-3}. In particular, the ends labeled ~\eqref{deg:h1-4-family} cancel ends on $\cQ_2$, and another family of ends on $\cQ_2$ and $\cQ_3$ cancel. Counting all the ends, we obtain
\[
E^{\tau_1}+E^{\tau_2}=\d\circ  H^{\tau_1\to \tau_2}+H^{\tau_1\to \tau_2} \circ \d,
\] 
where $H^{\tau_1\to \tau_2}$ counts index $-2$ rectangles which are $(v_i,\tau)$-matched for some $\tau\in (\tau_1,\tau_2)$.
\end{proof}

We now move on to prove Proposition~\ref{prop:deform-tau-to-infinity}. The proof consists of carefully analyzing the ends of index $-1$ moduli spaces of rectangles which are $(v_i,\tau)$-matched for some $\tau\in [\tau_0,\infty)$. The ends appearing as $\tau\to \infty$ are more complicated than those which appeared for $\tau\in [\tau_1,\tau_2]$. Unlike in Lemma~\ref{lem:homotopy-equivalence-different-t}, which concerned holomorphic disks, the ends which appear as $\tau\to \infty$ do not correspond exactly to the map $E^\infty$.

Before moving on to consider the ends appearing as $\tau\to \infty$, we first record the ends of $(v_i,\tau)$-matched moduli spaces which are easy to understand. The ends at $\tau=\tau_0$ contribute the term
\[
E^{\tau_0}(\xs).
\]
As in Lemma~\ref{lem:Etau-chain-map}, there are ends at finite $\tau\in [\tau_0,\infty)$, similar to \eqref{deg:h1-1-family},\dots,\eqref{deg:h1-4-family}, and their analogs on the quadruples $\cQ_2$ and $\cQ_3$. The total count of these ends gives the term
\[
(\d\circ H^{\tau_0\to \infty}+H^{\tau_0\to \infty}\circ \d)(\xs),
\] 
where $H^{\tau_0\to \infty}$ counts $(v_i,\tau)$-matched, index $-2$ rectangles, for $\tau\in [\tau_0,\infty)$.

The ends of the $(v_i,\tau)$-matched moduli spaces which appear as $\tau\to \infty$ may be partitioned into two sets:
\begin{enumerate}[label=($e$-\arabic*), ref=$e$-\arabic*, leftmargin=*,widest=IIII]
\item\label{h-infty-ends-1}  Ends where $\tau\to \infty$ and $s\in (0,1)$.
\item\label{h-infty-ends-2} Ends where  $\tau\to\infty$ and $s= 0$ or $s= 1$.
\end{enumerate}

To analyze the ends labeled ~\eqref{h-infty-ends-1} and~\eqref{h-infty-ends-2}, we  introduce a new type of evaluation map. We first introduce a version of this evaluation map for disks, and then move on to an analog for triangles and quadrilaterals which is relevant to \eqref{h-infty-ends-1} and~\eqref{h-infty-ends-2}.

 Consider the vector field $v_Y$ on $[0,1]\times \R$, shown in Figure~\ref{fig:143}.  Suppose $S^{\qs}$ is marked Riemann surface, $\phi$ is a homology class of disks, and $|\qs|=k$. Write $\R^k/\R$ for $\R^k$, modulo the diagonal action of $\R$. There is an \emph{unstable asymptotic evaluation map}
\[
\ev_{\qs}^{\un}\colon \cM_J(S^{\qs},\phi)\to \R^k/\R,
\]
defined as follows.  Write $\qs=(q_1,\dots, q_k)$. Let $R^\tau$ denote the flow of $v_Y$. If $u\in\cM_J(S^{\qs},\phi)$, we consider the tuple
\begin{equation}
\lim_{\tau\to -\infty} \big( (\pi_{\R}\circ R^{\tau}\circ \pi_{\bD}\circ u)(q_1),\dots, (\pi_{\R}\circ R^{\tau}\circ \pi_{\bD}\circ u)(q_k)\big)\in \ell^k/\R. \label{eq:asymptotic-evaluation}
\end{equation}
Since $v_Y$ is $\R$-invariant, and $\ell$ is a flowline of $v_Y$, we obtain a canonical identification of $\ell$ with $\R$. We define $\ev_{\qs}^{\un}(u)$ to be the image of the tuple in \eqref{eq:asymptotic-evaluation} under the canonical identification of $\ell^k/\R$ with $\R^k/\R$. Note also the existence of the limit in ~\eqref{eq:asymptotic-evaluation} follows from ~\eqref{eq:vector-fields-v-y-technical}. 

 Similarly, if $\ell$ denotes the line $\{\tfrac{1}{2}\}\times \R$, and $\cM_J(S^{\qs}, \phi,\ell)$ denotes the set of holomorphic disks $u$ such that $(\pi_{\bD}\circ u)(\qs)\subset \ell$,  there is a \emph{stable asymptotic evaluation map}
\[
\ev_{\qs}^{\st}\colon \cM_J(S^{\qs},\phi,\ell)\to \R^k/\R,
\]
defined similarly. (In fact, in our present setting of holomorphic disks, $\ev_{\qs}^{\st}$ is induced by the ordinary evaluation $(\pi_{\R}\circ u)(q_i)$, though this will not be the case for triangles or rectangles).

We now consider analogous evaluation maps for holomorphic triangles and rectangles.  Suppose that $S^{\qs}$ is a marked Riemann surface, and $u\colon S\to \Sigma\times \Delta$ is a holomorphic triangle for the Heegaard triple $(\Sigma,\gs_1,\gs_2,\gs_3,w)$. Let $v$ denote one of the limiting vector fields on $\Delta$  shown in Figure~\ref{fig:145}. (The cyclic identification of $(\gs_1,\gs_2,\gs_3)$ with $(\as,\bs,\gs)$ is unimportant for our present purposes). If $v$ coincides with $\pm v_Y$ on the $(\gs_i,\gs_{i+1})$-cylindrical end of $\Delta$, then we write write $\ell_{\g_i,\g_{i+1}}$ for the line $\{\tfrac{1}{2}\}\times [0,\infty)\subset [0,1]\times [0,\infty)$ in the $(\gs_i,\gs_{i+1}
)$ cylindrical end of $\Delta$. Note that $\ell_{\g_i,\g_{i+1}}$ is a flowline of the vector field $v$. We refer to $\ell_{\g_i,\g_{i+1}}$ as the \emph{special line}. We write $\cW^s_v(\ell_{\g_i,\g_{i+1}})$ and (resp. $\cW^u_v(\ell_{\g_i,\g_{i+1}})$) for \emph{stable} (resp. \emph{unstable}) sets of $\ell_{\g_i,\g_{i+1}}$, by which we mean the set of points of $\Delta$ whose positive (resp. negative) flow under $v$ asymptotically approaches the line $\ell_{\g_i,\g_{i+1}}$. We may similarly define stable and unstable sets associated to other subsets of $\Delta$.

If $u$ is a holomorphic triangle with marked source $S^{\qs}$, and $\qs_0\subset \qs$ is a subset of the marked points such that
\begin{equation}
(\pi_{\Delta}\circ u)(\qs_0)\subset \cW^{s}_v(\ell_{\g_i,\g_{i+1}}), \label{eq:marked-points-in-stable-subset}
\end{equation} 
then there is a well defined element
\[
\ev^{\st}_{\qs_0}(u)\in \R^{|\qs_0|}/\R,
\]
defined similarly to ~\eqref{eq:asymptotic-evaluation}. The map $\ev_{\qs_0}^{\st}$ determines a continuous map on the subspace of $\cM_J(S^{\qs}, \phi)$ consisting of curves $u$ satisfying~\eqref{eq:marked-points-in-stable-subset}. The above terminology adapts easily for holomorphic quadrilaterals.

With this notation in place, we can now describe the ends of the moduli spaces $\cM\cM^{v_i,\tau\in [\tau_0,\infty)}(\psi,\psi_0)$. We have the following analog of Definition~\ref{def:infinity-matched-disk}:

\begin{define} \label{def:simple-infty-matched-rectangle} 
\begin{enumerate}[label=($h_1^\infty$-$e$-\arabic*), ref=$h_1^\infty$-$e$-\arabic*, leftmargin=*, widest=IIIIII]
\item[]
\item\label{deg:infinity-matched-rectangle-1}  We say a pair $\cU$ and $\cU_0$ of broken $J_s\wedge I_s$-holomorphic quadrilaterals on $\cQ_1$ is a \emph{simple, $(v_1,\infty)$-matched quadrilateral}, if it satisfies the following:
\begin{enumerate}[label=(\alph*), ref=\alph*, leftmargin=*]
\item$\cU$ consists of exactly one $J_s$-holomorphic rectangle $u\in \cM_{J_s}(S^{\qs},\psi')$ on $(\Sigma,\as'',\as',\as,\bs)$, as well as a collection of index 2 boundary degenerations $b_1,\dots, b_n$. Furthermore $\psi=\psi'+[b_1]+\dots +[b_n]$.
\item $\cU_0$ consists of exactly one $I_s$-holomorphic rectangle $u_0\in \cM_{I_s}(S_0^{\qs_0},\psi_0')$ on $(\bT^2,c_3,c_2,c_1,c_3')$, as well as a collection of index 2 boundary degenerations $d_1,\dots, d_m$. Furthermore, $\psi_0=\psi_0'+[d_1]+\cdots [d_m].$
\item $n_{p}(\psi')=m=|\qs|$ and $n_{p_0}(\psi_0')=n=|\qs_0|$. Furthermore,
\[
(\pi_{\Box}\circ u)(\qs)\subset \cW^s_{v_1}( e_{\a''}\cup e_{\a'}\cup e_{\a})\quad \text{and} \quad (\pi_{\Box}\circ u_0)(\qs_0)\subset \cW^u_{v_1}(e_{c_3}).
\]
\item \label{eq:infinity-matching-quadrilateral} The following matching conditions are satisfied:
\[
 \begin{split}
( \pi_{\Sigma}\circ u)(q_{j})&=p,\\
(\pi_{\bT^2}\circ u_0)(q_{0,j})&=p_0,\\
(\pi_{\bT^2}\circ u_0)\left((\pi_{\Box}\circ u_0)^{-1}\left(\ev_{q_{j}}^{\un}(u)\right)\right)
 &=\ev^{\infty}(d_j),\quad \text{and} \\
(\pi_{\Sigma}\circ u)\left((\pi_{\Box}\circ u)^{-1}\left(\ev^{\st}_{q_{0,j}}(u_0)\right)\right)&=\ev^{\infty}(b_j).
\end{split}
\]
\end{enumerate}
\end{enumerate}
\begin{enumerate}[label=($h_2^\infty$-$e$-\arabic*), ref=$h_2^\infty$-$e$-\arabic*, leftmargin=*, widest=IIIIII]
\item\label{deg:infinity-matched-rectangle-2} A pair $(\cU,\cU_0)$ is a \emph{simple, $(v_2,\infty)$-matched quadrilateral} on $\cQ_2$ if it satisfies the analog of \eqref{deg:infinity-matched-rectangle-1} for $\cQ_2$.
\end{enumerate}
\begin{enumerate}[label=($h_3^\infty$-$e$-\arabic*), ref=$h_3^\infty$-$e$-\arabic*, leftmargin=*, widest=IIIIII]
\item\label{deg:infinity-matched-rectangle-3} A pair $(\cU,\cU_0)$ is a \emph{simple, $(v_3,\infty)$-matched quadrilateral} on $\cQ_3$ if it satisfies the analog of \eqref{deg:infinity-matched-rectangle-1} for $\cQ_3$.
\end{enumerate}
\end{define}

\begin{figure}[H]
	\centering
	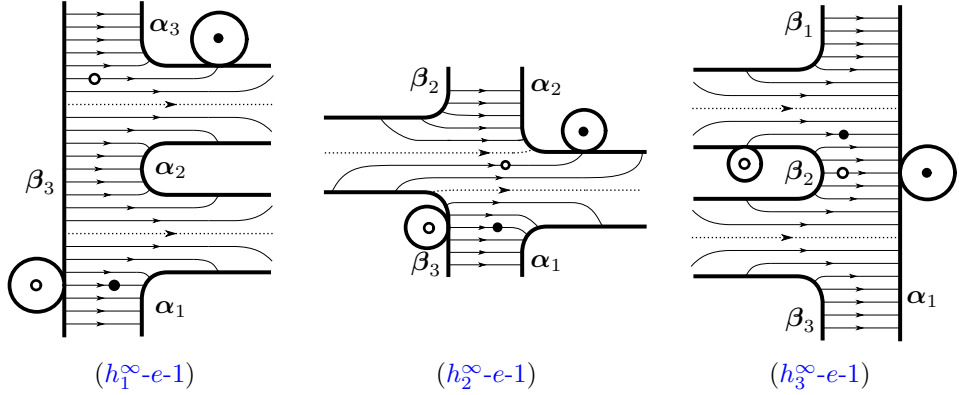
	\caption{Simple, $(v_i,\infty)$-matched quadrilaterals from Definition~\ref{def:simple-infty-matched-rectangle}. The open dots denote the images of the marked points  mapping to $p\in \Sigma$, and the solid dots denote the marked points  mapping   to $p_0\in \bT^2$.}\label{fig:147}
\end{figure}

Examples of simple, $(v_i,\infty)$-matched quadrilaterals are shown in Figure~\ref{fig:147}.  Note that \eqref{deg:infinity-matched-rectangle-1}, \eqref{deg:infinity-matched-rectangle-2} and \eqref{deg:infinity-matched-rectangle-3} are natural analogs of Definition~\ref{def:infinity-matched-disk} for rectangles.

\begin{lem}\label{lem:constrain-degenerations-finite-s}
The ends of the moduli spaces $\cM\cM^{v_i,\tau\in [\tau_0,\infty)}(\psi,\psi_0)$ which appear as $\tau\to \infty$ and $s\to s_0\in (0,1)$ (labeled~\eqref{h-infty-ends-1} above) generically correspond to simple, $(v_i,\infty)$-matched holomorphic rectangles. Furthermore, each such broken curve makes odd algebraic contribution to the ends of the moduli spaces $\cM\cM^{v_i,\tau\in [\tau_0,\infty)}(\psi,\psi_0)$.
\end{lem}
\begin{proof}
We focus on the curves appearing on $\cQ_1$. We can view the moduli space $\cM\cM^{v_i,\tau\in [\tau_0,\infty)}(\psi,\psi_0)$ as consisting of holomorphic curves mapping into $(\Sigma\sqcup \bT^2)\times \Box$, and hence limiting curves may be described in terms of holomorphic combs, via an extension of Definition~\ref{def:holomorphic-comb} to handle holomorphic polygons. Since we are restricting to the case that $s_0\in (0,1)$, the limiting collection of curves consists of a pair of holomorphic rectangles (one on each of $\bT^2\times \Box$ and $\Sigma\times \Box$), as well as a collection of holomorphic curves in the four ends, and also boundary degenerations. Let us write $\cU_1,\dots, \cU_n$ for the disk levels appearing in the $(\as_2,\as_3)$ end. Let us write $\cV_1,\dots, \cV_m$ for the disk levels appearing in the $(\as_1,\as_2)$ end. (These are not shown in Figure~\ref{fig:201}.)  For ease of notation, we group any boundary degenerations into the disk or rectangle level that they originate from.

As a first step, we will show that generically $n\le 1$ and $m\le 1$, and that there are no disk levels in the $(\as_1,\bs_3)$ or $(\as_3,\bs_3)$ ends. To simply the notation, we first show that if $m=0$, then $n\in \{0,1\}$. The same argument applies to show that more generally, $n,m\in \{0,1\}$. Similarly, it is easy to use index considerations to show that generically no curves can appear in the $(\as_1,\bs_3)$ or $(\as_3,\bs_3)$ ends (since here any disks are $v_H$-matched).

Write $\cU_0$ for the curves in the quadrilateral level (as well as any boundary degenerations branching off of this level), and let $\cU_1,\dots, \cU_n$ be the disks in the $(\as_2,\as_3)$ levels, as above. We can further partition these curves into collections $\cU^\Sigma_0,\dots, \cU_n^\Sigma$ and $\cU_0^\bT,\dots, \cU_n^{\bT}$, depending on whether they map to $\Sigma$ or $\bT^2$. See Figure~\ref{fig:201}.

 \begin{figure}[H]
	\centering
	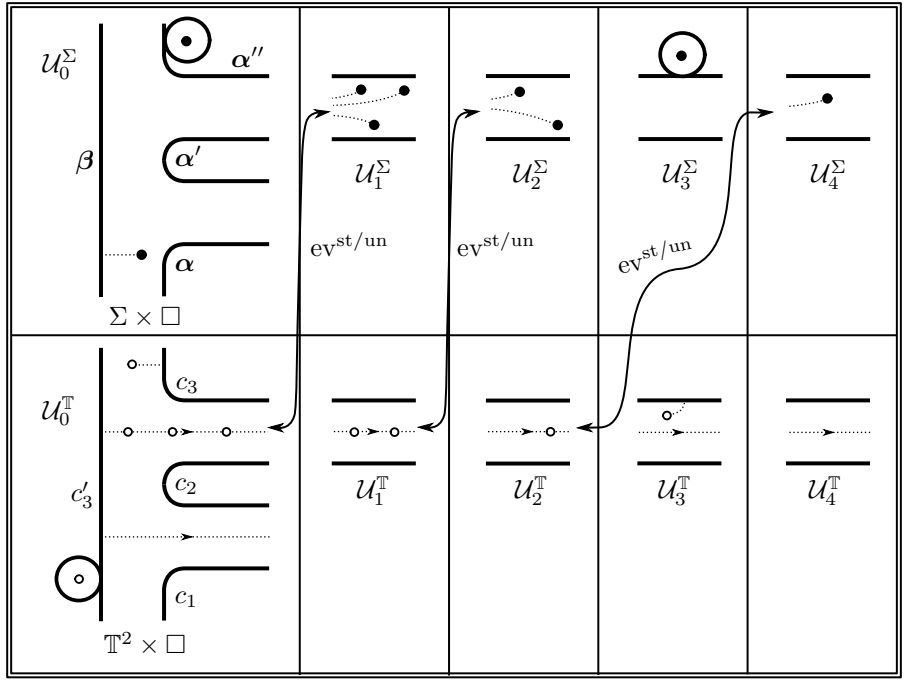
	\caption{A non-generic degeneration appearing in \eqref{h-infty-ends-1}. }\label{fig:201}
\end{figure}  

We now trim off any boundary degenerations from the collections $\cU_i^\bT$ and $\cU_i^\Sigma$ (and we call the remaining collections $\cU_i^\bT$ and $\cU_i^\Sigma$, abusing notation slightly). This leaves two types of marked points on the remaining curves. We say a marked point is \emph{coupled} if it is matched to a marked point which remains after trimming boundary degenerations. We say a marked point is \emph{decoupled} if it was matched to a marked point in a boundary degeneration.

For $i\in \{0,\dots, n\}$, let $d_i \ge 0$ denote the number coupled marked points on the curves in $\cU^{\bT}_i$. The marked points of $\cU^{\bT}_i$ must be mapped to the special line of $[0,1]\times \R$ or the special line of $\Box$ which extends into the $(\as_2,\as_3)$ end.  There is an injection $f\colon \{i: d_i>0\}\to \{1,\dots, n\}$, which sends an $i$ to the index $j$ such that the $d_i$ marked points of  $\cU_i^{\bT^2}$ are matched with marked points of $\cU_j^{\Sigma}$. Note that if $d_i>0$, then the curves of $\cU_i^{\bT}$ are matched with the curves of $\cU_{f(i)}^\Sigma$ via the asymptotic evaluation map $\ev^{\st/\un}$ applied to the coupled marked points. We write
\begin{equation}
\ve{d}:=\sum_{i=0}^n d_i. \label{eq:bold-d-sum-ds}
\end{equation}

 Write $\phi_0^\Sigma,\dots, \phi_n^\Sigma$ and $\phi_0^\bT,\dots, \phi_n^\bT$ for the homology classes of $\cU^\Sigma_0,\dots, \cU^\Sigma_n$ and $\cU_0^\bT,\dots, \cU_n^\bT$. Here, $\phi_0^\bT$ is a class of rectangles, while $\phi_i^\bT$ is a class of disks, for $i>0$, and similarly for the classes on $\Sigma$. By assumption, we are considering curves which appeared in the limit of a 1-parameter family of $(v_i,\tau)$-matched curves representing a pair $(\psi,\psi_0)$ satisfying 
 \begin{equation}
 \mu(\psi)+\mu(\psi_0)=2 n_p(\psi)-1, \quad \text{and} \quad n_{p}(\psi)=n_{p_0}(\psi_0). \label{eq:simple-infinty-matched-quad-index}
 \end{equation}
  Since we have trimmed off boundary degenerations, we obtain
\begin{equation}
\sum_{i=0}^n \left(\mu(\phi_i^\Sigma)+\mu(\phi_i^\bT)\right)=2\ve{d}-1.\label{eq:a-priori-Maslov-index-estimate}
\end{equation}
Write $n^\bT$ and $n^{\Sigma}$ for the number of non-trivial disk levels amongst the $\cU_i^\bT$ and $\cU_i^\Sigma$, respectively, ranging over $i>0$. We define 
\[
\Phi^\bT:=(\phi_0^\bT,\dots, \phi_n^\bT,d_0,\dots, d_n),
\]
an ordered tuple of classes, counts of coupled marked points. We define an ordered tuple $\Phi^{\Sigma}$ similarly.

We can consider the collection $\cU_0^{\bT},\dots, \cU_n^{\bT},\cU_0^\Sigma,\dots, \cU_n^\Sigma$ as living in their own moduli space $\cM\cM(\Phi^\Sigma,\Phi^\bT,f)$, consisting of such collections of curves which are matched under the asymptotic evaluation map. The expected dimension of this moduli space at embedded curves is given by 
\begin{equation}
\dim \cM\cM(\Phi^\Sigma,\Phi^\bT,f)=1+\sum_{i=0}^n \left( \mu(\phi_i^\Sigma)+\mu(\phi_i^\bT)-d_i-\max(0,d_i-1)\right). \label{eq:expected-dim-non-simple-rectangles}
\end{equation}
In~\eqref{eq:expected-dim-non-simple-rectangles}, $\sum_{i=0}^n (\mu(\phi_i^\Sigma)+\mu(\phi_i^{\bT}))$ is the expected dimension of the unconstrained moduli space, for a fixed almost complex structure. The $+1$ in~\eqref{eq:expected-dim-non-simple-rectangles} comes from the parameter $s$,  the summand $-d_i$ is the result of the constraint that the coupled marked points on the $\bT^2$-side must lie on the special line, and the summand of $-\max(0,d_i-1)$ given by matching condition involving the asymptotic evaluation map.

The moduli space $\cM\cM(\Phi^\Sigma,\Phi^\bT,f)$ has a free action of $\R^{n^\Sigma+n^\bT}$, where each component acts by translation on one of the $n^\Sigma+n^\bT$ non-trivial disk components of $\cM\cM(\Phi^\Sigma,\Phi^\bT,f)$. In particular, for $\cM\cM(\Phi^\Sigma,\Phi^\bT,f)$ to be non-empty, it must have dimension at least $n^\Sigma+n^\bT$. Hence~\eqref{eq:expected-dim-non-simple-rectangles} gives
\[
1+\sum_{i=0}^n \left( \mu(\phi_i^\Sigma)+\mu(\phi_i^\bT)-d_i-\max(0,d_i-1)\right)\ge n^\Sigma+n^\bT.
\]
Applying~\eqref{eq:bold-d-sum-ds} and~\eqref{eq:a-priori-Maslov-index-estimate} and rearranging, we obtain that if $\cM\cM(\Phi^\Sigma,\Phi^\bT,f)$ is non-empty, then
\begin{equation}
\ve{d}-n^\Sigma-n^\bT-\sum_{i=0}^n \max(0,d_i-1)\ge 0.\label{eq:almost-there-inequality}
\end{equation}
Since for each $i$ with $d_i>0$, the level $\cU_{f(i)}^\Sigma$ is non-trivial, we obtain that 
\[
n^\Sigma+\sum_{i=0}^n \max(0,d_i-1)= \ve{d}+n^{\Sigma}_0,
\]
where $n^{\Sigma}_0$ denotes the number of non-trivial levels of $\cU_1^\Sigma,\dots \cU_n^\Sigma$ which have no coupled marked points. Hence, from~\eqref{eq:almost-there-inequality} we obtain $-n^\bT-n^{\Sigma}_0\ge 0,$
which implies that $n^\bT=n^{\Sigma}_0=0$. It follows that $n=n^{\Sigma}\in \{0,1\}$, since all coupled marked points must be matched with the rectangle $\cU_0^\bT$, and at most one level can match this level. This completes the proof of the subclaim under the assumption that $m=0$. Furthermore, the argument easily extends to handle the case that $m$ and $n$ are both allowed to be non-zero, and implies that $n,m\in \{0,1\}$.

It follows from the above argument that, after trimming boundary degenerations, the limiting curve consists of up to four levels. There are two holomorphic quadrilaterals, $v^\Sigma$ and $v^\bT$, as well as potentially two holomorphic disks, $u_{\a'',\a'}$ and $u_{\a',\a}$ on $\bT^2\times [0,1]\times \R$, which occur in the  $(\as'',\as')$ and $(\as',\as)$ ends, respectively. The curves $u_{\a'',\a'}$ and $u_{\a',\a}$ are matched to $v^\bT$ via the asymptotic evaluation map. See Figure~\ref{fig:202}. 

 \begin{figure}[H]
	\centering
	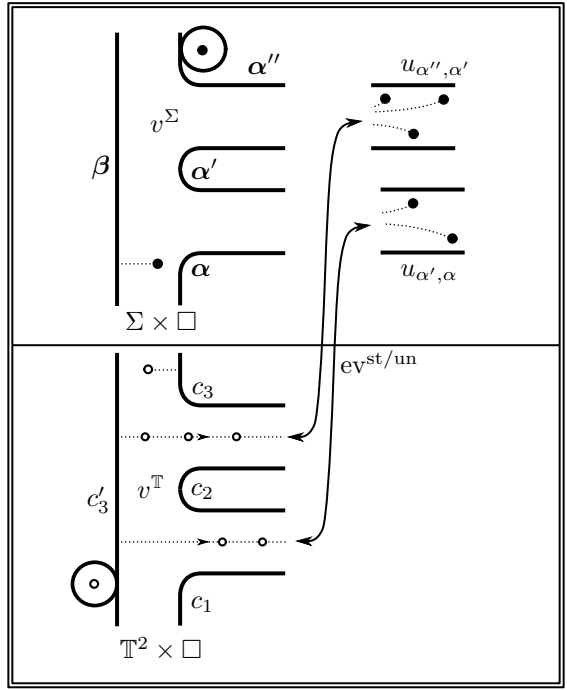
	\caption{A more generic configuration, which is still prohibited for $\cQ_1$.}\label{fig:202}
\end{figure} 

Up until now, we have not used any particular facts about the quadruple $\cQ_1$. The key observation that prohibits degenerations like those in Figure~\ref{fig:202} is that if $\phi\in \pi_2(\Theta^+,\Theta)$ is a class of disks on $(\Sigma,\as'',\as')$ or $(\Sigma,\as',\as)$, then then
\begin{equation}
\mu(\phi)=2n_{p}(\phi)+\gr(\Theta^+,\Theta)\ge 2n_p(\phi). \label{eq:maslov-index-disks-on-small-diagrams}
\end{equation}

Write $d_{\a'',\a'}$ and $d_{\a',\a}$ for the number of marked points mapped to the two special lines of $\Box$, by the holomorphic rectangle on $\Sigma\times \Box$ of a curve appearing in a limiting configuration. Since this curve was obtained in a limit of a path of curves representing a pair $(\psi,\psi_0)$ satisfying~\eqref{eq:simple-infinty-matched-quad-index}, we obtain, similarly to~\eqref{eq:a-priori-Maslov-index-estimate}, that
\begin{equation}
\mu(v^\Sigma)+\mu(v^\bT)+\mu(u_{\a'',\a'})+\mu(u_{\a',\a})=2d_{\a'',\a'}+2d_{\a',\a}-1. \label{eq:formula-maslov-index-Q-1-degen}
\end{equation}
Combining ~\eqref{eq:maslov-index-disks-on-small-diagrams} and~\eqref{eq:formula-maslov-index-Q-1-degen}, we obtain
\begin{equation}
\mu(v^\Sigma)+\mu(v^\bT)\le -1.\label{eq:almost-there-index-prohibiting-ends}
\end{equation}
However, ~\eqref{eq:almost-there-index-prohibiting-ends} implies that generically there are no representatives of the pair $(v^\Sigma,v^\bT)$, ranging over $s\in (0,1)$, if $v^{\bT}$ has one or more coupled marked points which are sent to the special line, since this adds a constraint of codimension at least 1. This completes the proof for $\cQ_1$.  The proofs for the quadruples $\cQ_2$ and $\cQ_3$ follow from similar arguments. 
\end{proof}

 The count of index 2 boundary degenerations in Proposition~\ref{prop:relaxed-count} lets us identify the count of curves of type~\eqref{deg:infinity-matched-rectangle-1}, \eqref{deg:infinity-matched-rectangle-2} and ~\eqref{deg:infinity-matched-rectangle-3} with the map $E^\infty$.

\begin{figure}[ht!]
	\centering
	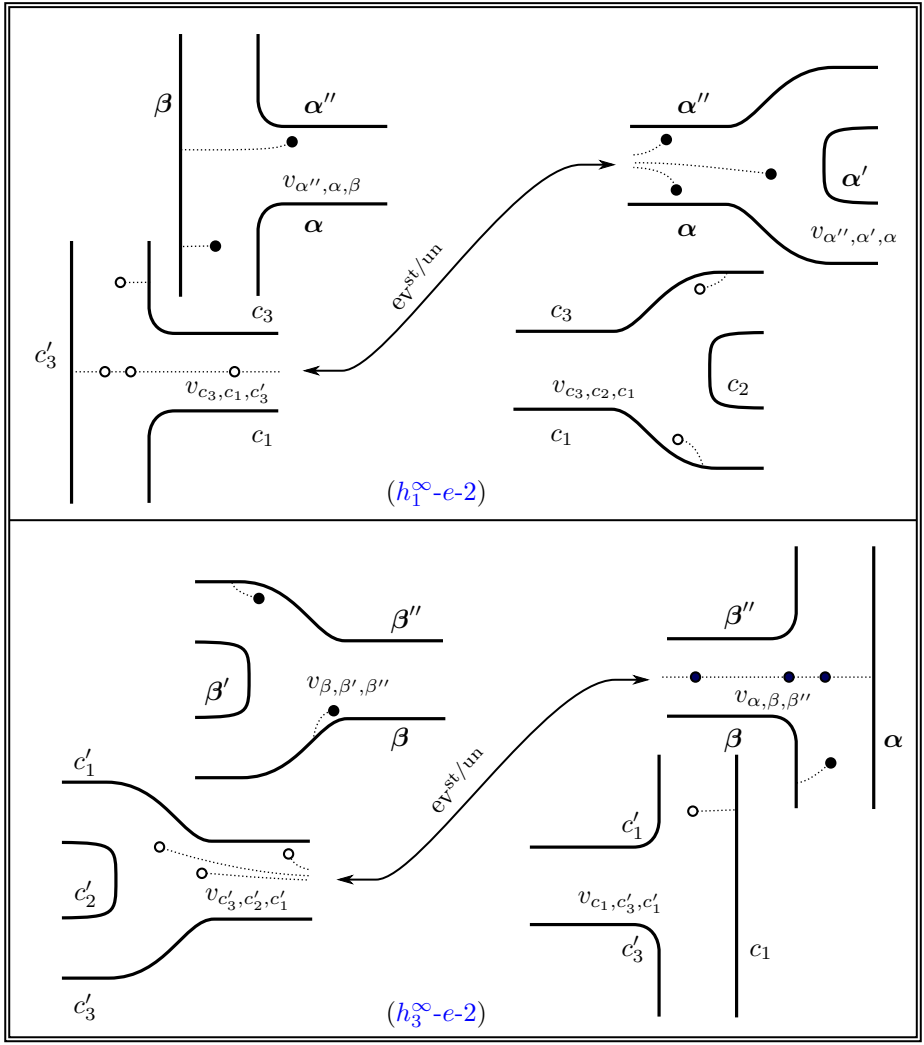
	\caption{A schematic of configurations \eqref{deg:infinity-asymptotically-matched-rectangle-h1} and~\eqref{deg:infinity-asymptotically-matched-rectangle-h3}. Open dots indicate marked points which map to $\Sigma$. Solid dots indicate marked points mapped to $\bT^2$.  Boundary degenerations have been trimmed.}\label{fig:148}
\end{figure}

We now consider the ends labeled ~\eqref{h-infty-ends-2} above, which appear as $\tau\to \infty$ and $s\to 1$ or $s\to 0$. We consider the following two types of broken curves (see Figure~\ref{fig:148}):
\begin{enumerate}[label=($h_1^\infty$-$e$-\arabic*), ref=$h_1^\infty$-$e$-\arabic*,leftmargin=*, widest=IIIIII]
\setcounter{enumi}{1}
\item\label{deg:infinity-asymptotically-matched-rectangle-h1} A pair of broken holomorphic rectangles $\cV$ and $\cV_0$ on $\cQ_1$, satisfying the following:
\begin{enumerate}[label=(\alph*), ref=\alph*]
\item $\cV$ consists of a holomorphic triangle $v_{\a'',\a,\b}$ on $(\Sigma,\as'',\as,\bs)$ and a holomorphic triangle $v_{\a'',\a',\a}$ on $(\Sigma,\as'',\as',\as)$, as well as a collection of index 2 boundary degenerations $d_1,\dots, d_n$.
\item $\cV_0$ consists of a holomorphic triangle $v_{c_3,c_1,c_3'}$ on $(\bT^2,c_3,c_1,c_3')$ and a holomorphic triangle $v_{c_3,c_2,c_1}$ on $(\bT^2,c_3,c_2,c_1)$, as well as a collection of index 2 boundary degenerations $e_1,\dots e_m$.
\item The source of the holomorphic triangle $v_{\a'',\a,\b}$ has two disjoint collections of marked points $\ve{q}^c_{\a'',\a,\b}$, $\ve{q}^{d}_{\a'',\a,\b}$, which we call the \emph{coupled} and \emph{decoupled} marked points, respectively. The other three holomorphic triangles have analogous collections of marked points.
\item We have $|\qs^d_{\a'',\a,\b}|+|\qs^d_{\a'',\a',\a}|=m$ and $|\qs^d_{c_3,c_1,c_3'}|+|\qs^d_{c_3,c_2,c_1}|=n$.
\item $v_{\a'',\a,\b}$ and $v_{c_3,c_2,c_1}$ have no coupled marked points.
\item $v_{\a'',\a',\a}$ has no decoupled marked points.
\item $|\ve{q}^c_{\a'',\a',\a}|=|\ve{q}^c_{c_3,c_1,c_3'}|$.
\item The coupled marked points of $v_{c_3,c_1,c_3'}$ map to the special line of $\Delta$.
\item The coupled marked points of $v_{\a'',\a',\a}$ and $v_{c_3,c_1,c_3'}$ are matched under the asymptotic evaluation map $\ev^{\st/\un}$, which takes values in $\R^k/\R$, where $k=|\ve{q}_{\a'',\a',\a}^c|$.
\item The decoupled marked points of  $v_{\a'',\a,\b}$ are matched to the boundary degenerations $e_1,\dots, e_m$, similar to part~\eqref{eq:infinity-matching-quadrilateral} of \eqref{deg:infinity-matched-rectangle-1}, in Definition~\ref{def:simple-infty-matched-rectangle}. Similarly the decoupled marked points of $v_{c_3,c_2,c_1}$ and $v_{c_3,c_1,c_3'}$ are matched to the boundary degenerations $d_1,\dots, d_n$.
\end{enumerate}
\end{enumerate}
\begin{enumerate}[label=($h_3^\infty$-$e$-\arabic*), ref=$h_3^\infty$-$e$-\arabic*,leftmargin=*,widest=IIIIII]
\setcounter{enumi}{1}
\item\label{deg:infinity-asymptotically-matched-rectangle-h3} A pair of broken holomorphic rectangles $\cV$ and $\cV_0$ on $\cQ_3$, satisfying the analogs of \eqref{deg:infinity-asymptotically-matched-rectangle-h1}.
\end{enumerate}

The argument from Lemma~\ref{lem:constrain-degenerations-finite-s} adapts to show the following:

\begin{lem}\label{lem:extra-ends-tau->infty,s->1}
 The ends of $\cM\cM^{\tau\in [\tau_0,\infty),v_i}_{J\wedge I}(\psi,\psi_0)$ labeled~\eqref{h-infty-ends-2} above (i.e. appearing as $\tau\to \infty$ and $s\to 0$ or $s\to 1$) correspond generically to broken curves which satisfy either~\eqref{deg:infinity-asymptotically-matched-rectangle-h1} or~\eqref{deg:infinity-asymptotically-matched-rectangle-h3}. Furthermore, any such broken curve makes odd algebraic contribution to the count of ends of its moduli space.
\end{lem}

Write $X_1^{\infty}$ and $X_3^{\infty}$ for the maps which count the elements of 0-dimensional moduli spaces of broken curves satisfying~\eqref{deg:infinity-asymptotically-matched-rectangle-h1} or ~\eqref{deg:infinity-asymptotically-matched-rectangle-h3}, respectively. Combining Lemma~\ref{lem:extra-ends-tau->infty,s->1} with the proof of Lemma~\ref{lem:deform-tau}, by counting the ends of $[\tau_0,\infty)$-matched moduli spaces, we obtain the relation
\begin{equation}
E^{\tau_0}+E^{\infty}+X_1^{\infty}+X_3^{\infty}=\d  H^{\tau_0\to \infty}+H^{\tau_0\to \infty}\d. \label{eq:Etau=Einfty-almost}
\end{equation}

We may trim off the boundary degenerations in the definitions of the maps $X_1^\infty$ and $X_3^{\infty}$, as in ~\eqref{eq:equivalence-trimmed-untrimmed}. Hence $X_1^\infty$ may be identified as the counts of quadruples of  holomorphic triangles 
\begin{equation}
(v_{\a'',\a,\b}, v_{c_3,c_1,c_3'},v_{\a'',\a',\a},v_{c_3,c_2,c_1}), \label{eq:trimmed-4-tuples-Psi-1}
\end{equation}
which are equipped with marked points and satisfy the matching conditions of \eqref{deg:infinity-asymptotically-matched-rectangle-h1}, which furthermore
represent a quadruple of triangle classes $(\psi_{\a'',\a,\b}, \psi_{c_3,c_1,c_3'}, \psi_{\a'',\a',\a}, \psi_{c_3,c_2,c_1})$ satisfying
\begin{equation}
\mu(\psi_{\a'',\a,\b}) +\mu(\psi_{c_3,c_1,c_3'})+\mu (\psi_{\a'',\a',\a})+ \mu(\psi_{c_3,c_2,c_1})=2 n_{p_0}(\psi_{\a'',\a',\a})-1. \label{eq:Maslov-index-trimmed-rectangle}
\end{equation}
The map $X_3^{\infty}$ has an analogous description. The maps $X_1^\infty$ and $X_3^{\infty}$ assign $U$ and $T$ weights to quadruples of triangles similarly to ~\eqref{eq:weighting-trimmed-curves}.

Proposition~\ref{prop:deform-tau-to-infinity} now follows from ~\eqref{eq:Etau=Einfty-almost}, together with the following lemma: 

\begin{lem}\label{lem:Psi1-Psi3}
\begin{enumerate}
\item[]
\item\label{lem:Psi1} $X_1^\infty=0$.
\item\label{lem:Psi3} $X_3^{\infty}\simeq 0.$
\end{enumerate}
\end{lem}

 Lemma~\ref{lem:Psi1-Psi3} is slightly involved, so we prove the two subclaims separately. We begin with the proof for $X_1^\infty$:

\begin{proof}[Proof of Part~\eqref{lem:Psi1} of Lemma~\ref{lem:Psi1-Psi3} ] We focus on the trimmed description of $X_1^\infty$, which counts 4-tuples of marked holomorphic triangles, as in ~\eqref{eq:trimmed-4-tuples-Psi-1}.

The map $X_1^\infty$ counts 4-tuples of holomorphic triangles, representing 4-tuples of classes
\[
(\psi_{\a'',\a,\b}, \psi_{c_3,c_1,c_3'}, \psi_{\a'',\a',\a}, \psi_{c_3,c_2,c_1}).
\]
We write $S^{\qs_{\a'',\a,\b}}$ for the marked source of the holomorphic triangle $v_{\a'',\a,\b})$. The marked points $\qs_{\a'',\a,\b}$ are partitioned into coupled and decoupled marked points, $\qs_{\a'',\a,\b}^c$ and $\qs_{\a'',\a,\b}^d$, respectively, as in ~\eqref{deg:infinity-asymptotically-matched-rectangle-h1}. The same holds for the other three holomorphic triangles, and we use the analogous notation for their marked sources.

Define
\[
N:=|\qs_{c_3,c_1,c_3'}^c|,
\]
which we note coincides with $|\qs^c_{\a'',\a',\a}|$. Note that $v_{\a'',\a,\b}$ and $v_{c_3,c_2,c_1}$ have no coupled marked points.  We recall from ~\eqref{eq:Maslov-index-trimmed-rectangle} that
\begin{equation}
\mu(\psi_{\a'',\a,\b})+\mu(\psi_{c_3,c_1,c_3'})+\mu(\psi_{\a'',\a',\a})+\mu(\psi_{c_3,c_2,c_1})=2N-1.\label{eq:sum-Maslov-indices-1}
\end{equation}

The matching condition only involves $\psi_{\a'',\a',\a}$ and $\psi_{c_3,c_1,c_3'}$, and it imposes a codimension $2N-1$ condition, by Proposition~\ref{prop:dimension-counts}. More explicitly, the restriction that $\qs_{c_3,c_1,c_3'}^c$ map to the special line of $\Delta$ imposes a codimension $N$ condition, while the matching of $v_{c_3,c_1,c_3'}$ and $v_{\a'',\a',\a}$  under the evaluation maps $\ev^{\st/\un}$, which take values in $\R^{N}/\R$,  gives an extra codimension $N-1$ condition.

 Hence, for a 4-tuple of marked sources and homology classes to have holomorphic representatives for generic almost complex structures, it is necessary that
\[
\mu(\psi_{\a'',\a,\b}) \ge 0 \quad \mu(\psi_{c_3,c_2,c_1})\ge 0 , \quad  \text{and} \quad \mu(\psi_{\a'',\a',\a})+\mu(\psi_{c_3,c_1,c_3'})\ge 2N-1.\label{eq:bound-from-asymptotic-matching-1}
\]
Equations~\eqref{eq:sum-Maslov-indices-1} and~\eqref{eq:bound-from-asymptotic-matching-1} imply that if a 4-tuple of homology classes and marked sources have matched representatives, then
\begin{equation}
\mu(\psi_{\a'',\a,\b})=\mu(\psi_{c_3,c_2,c_1})=0\quad \text{and} \quad \mu(\psi_{\a'',\a',\a})+\mu(\psi_{c_3,c_1,c_3'})= 2N-1. \label{eq:generic-maslov-index-Psi-1}
\end{equation}
In particular, the curves counted by $X_1^\infty$ on $(\bT^2,c_3,c_2,c_1)$ are exactly the Maslov index 0 triangles, with no constraints. Consequently, the equality $X_1^\infty=0$ follows immediately from Lemma~\ref{lem:model-computations-m-surgery}, the ordinary count of index 0 triangles on $(\bT^2,c_3,c_2,c_1)$.
\end{proof}

The computation for $X_3^\infty$ is slightly more involved. As in the proof of  part~\eqref{lem:Psi1} of Lemma~\ref{lem:Psi1-Psi3}, we focus on the trimmed description of $X_3^\infty$, which counts quadruples of marked holomorphic triangles
\[
(v_{c_3',c_2',c_1'},v_{c_1,c_3',c_1'},v_{\b,\b',\b''},v_{\a,\b,\b''}),
\]
such that the evaluation under $\ev^{\st/\un}$ on the coupled marked points of $v_{c_3',c_2',c_1'}$ and $v_{\a,\b,\b'}$ coincide.

For the map $X_3^\infty$, the triangles $v_{c_3',c_2',c_1'}$ have only coupled marked points, and hence the classes counted by $X_3^\infty$ satisfy
\begin{equation}
\mu(\psi_{c_1,c_3',c_1'})+\mu(\psi_{\b,\b',\b''})+\mu(\psi_{\a,\b,\b''})+\mu(\psi_{c_3',c_2',c_1'})=2n_{p_0}(\psi_{c_3',c_2',c_1'})-1.\label{eq:Psi-3-condition}
\end{equation}

The matching condition satisfied by $v_{\a,\b,\b''}$ and $v_{c_3',c_2',c_1'}$ consists of the following. Both triangles have exactly $n_{p_0}(\psi_{c_3',c_2',c_1'})$ coupled marked points. The coupled marked points of $v_{\a,\b,\b''}$ occur along the special line of $\Delta$. Furthermore,
\[
\ev^{\st}_{\qs^c_{c_3',c_2',c_1'}}(v_{c_3',c_2',c_1'})=\ev^{\un}_{\qs^c_{\a,\b,\b''}}(v_{\a,\b,\b''}).
\]
See Figure~\ref{fig:148} for a schematic.

To define our null-homotopy of $X_3^\infty$, we investigate expected dimensions further.
Let $S^{\qs}$ and $S_0^{\qs_0}$ denote two marked sources for $\psi_{\a,\b,\b''}$ and $\psi_{c_3',c_2',c_1'}$, respectively. Generically, the fibered product
\[
\cM(S^{\qs}, \psi_{\a,\b,\b''})\times_{\ev^{\st/\un}} \cM(S_0^{\qs_0}, \psi_{c_3',c_2',c_1'})
\]
will be a manifold, whose dimension near any embedded curve is given by
\begin{equation}
\begin{split}
&\dim \left(\cM(S^{\qs},\psi_{\a,\b,\b''})\times_{ \ev^{\st/\un}} \cM(S_0^{\qs_0},\psi_{c_3',c_2',c_1'})\right)\\
=& \mu(\psi_{\a,\b,\b''})+\mu(\psi_{c_3',c_2',c_1'})-2n_{p_0}(\psi_{c_3',c_2',c_1'})+1.
\end{split}\label{eq:expected-dimension-asymptotic-eval}
\end{equation}

We now deform the evaluation map. Given $r\in [1,\infty)$, we define
\[
\cM(S^{\qs}, \psi_{\a,\b,\b''})\times_{r\cdot \ev^{\st/\un}} \cM(S_0^{\qs_0}, \psi_{c_3',c_2',c_1'})
\]
to be the set of pairs of marked holomorphic triangles $(v_{c_3',c_2',c_1'}, v_{\a,\b,\b''})$ which satisfy
\begin{equation}
r\cdot \ev_{\qs^{c}_{c_3',c_2',c_1'}}^{\st}(v_{c_3',c_2',c_1'})=\ev^{\un}_{\qs^{c}_{\a,\b,\b''}}(v_{\a,\b,\b''}). \label{eq:r-dilated-triangles}
\end{equation}
In~\eqref{eq:r-dilated-triangles}, $r$ acts diagonally on $\R^k/\R$ by multiplication.

Furthermore, the parametrized moduli space 
\begin{equation}
\begin{split}
&\cM\cM^{r\cdot \ev^{\st/\un}, r\in [1,\infty)} (S^{\qs},S_0^{\qs_0}, \psi_{\a,\b,\b''}, \psi_{c_3',c_2',c_1'})\\
:=&\bigcup_{r\in [1,\infty)} \cM(S^{\qs},\psi_{\a,\b,\b''})\times_{r\cdot \ev^{\st/\un}} \cM(S_0^{\qs_0},\psi_{c_3',c_2',c_1'})\times \{r\}\label{eq:r-dilated-matched-moduli-spaces-parametrized-def}
\end{split}
\end{equation}
 will generically be a manifold, whose dimension near any embedded curve is
\begin{equation}
\begin{split}
&\dim\left(\cM\cM^{r\cdot \ev^{\st/\un}, r\in [1,\infty)} (S^{\qs},S_0^{\qs_0}, \psi_{\a,\b,\b''}, \psi_{c_3',c_2',c_1'}) \right)\\
=&\mu(\psi_{\a,\b,\b''})+\mu(\psi_{c_3',c_2',c_1'})-2n_p(\psi_{c_3',c_2',c_1'})+2.
\end{split}\label{eq:expected-dimension-asymptotic-eval-2}
\end{equation}

We define our homotopy $H_3^{\infty,[1,\infty)}$ to count quadruples of marked triangles 
\[
(v_{c_3',c_2',c_1'}, v_{c_1,c_3',c_1'}, v_{\b,\b',\b''}, v_{\a,\b,\b''}),
\] satisfying
\[
\mu(v_{c_1,c_3',c_1'})=\mu(v_{\b,\b',\b''})=0\quad \text{and} \quad \mu(v_{\a,\b,\b''})+\mu(v_{c_3',c_2',c_1'})=2n_{p_0}(v_{c_3',c_2',c_1'})-2,
\]
and such that $v_{\a,\b,\b''}$ and $v_{c_3',c_2',c_1'}$ are $r\cdot \ev^{\st/\un}$ matched for some $r\in [1,\infty)$. A quadruple of triangles is assigned $U$ and $T$ weights similarly to ~\eqref{eq:weighting-trimmed-curves}.

Part~\eqref{lem:Psi3} of Lemma~\ref{lem:Psi1-Psi3} is a consequence of the following lemma:

\begin{lem}\label{lem:H-3-null-homotopy-X-3} The maps $X_3^\infty$ and $H_3^{\infty,[1,\infty)}$ satisfy
\[
X_3^\infty=\left[\d, H_3^{\infty, [1,\infty)}\right].
\]
\end{lem}

The proof of Lemma~\ref{lem:H-3-null-homotopy-X-3} follows from a careful analysis of the ends of the moduli spaces in~\eqref{eq:r-dilated-matched-moduli-spaces-parametrized-def}, which we do in the subsequent Lemmas~\ref{lem:counting-ends-M-Psi-finite-r} and ~\ref{lem:ends-r-to-infty}. First, we introduce some notation. Suppose $\Psi$ is a tuple consisting of four homology classes of triangles $(\psi_{c_3',c_2',c_1'},\psi_{c_1,c_3',c_1'}, \psi_{\b,\b',\b''}, \psi_{\a,\b,\b''})$, as well as four marked sources; $r\ge 1$ is fixed. Write 
\[
\cM\cM^{r\cdot \ev^{\pm \infty}}(\Psi),
\]
for the moduli space of quadruples of marked triangles $(v_{c_3',c_2',c_1'}, v_{c_1,c_3',c_1'}, v_{\b,\b',\b''}, v_{\a,\b,\b''})$ such that $v_{c_3',c_2',c_1'}$ and $v_{\a,\b,\b''}$ are matched under $r\cdot \ev^{\st/\un}$. Similarly, we write
\begin{equation}
\cM\cM^{r\cdot \ev^{\st/\un}, r\in [1,\infty)}(\Psi) \label{eq:r-dilated-moduli-space-range}
\end{equation}
for the parametrized moduli space, consisting of such quadruples, which are $r$-matched for some $r\in [1,\infty)$.  

Note that the number of coupled marked points in the four homology classes of $\Psi$ is fixed. Namely, the sources of $v_{c_3',c_2',c_1'}$ and $v_{\a,\b,\b''}$ both have $n_{p_0}(c_3',c_2',c_1')$ coupled marked points. All other marked points are decoupled. We note that $\Psi$ also contains the information of a bijection between the coupled marked points of $v_{c_3',c_2',c_1'}$ and $v_{\a,\b,\b''}$.

\begin{figure}[ht!]
	\centering
\begingroup%
  \makeatletter%
  \providecommand\color[2][]{%
    \errmessage{(Inkscape) Color is used for the text in Inkscape, but the package 'color.sty' is not loaded}%
    \renewcommand\color[2][]{}%
  }%
  \providecommand\transparent[1]{%
    \errmessage{(Inkscape) Transparency is used (non-zero) for the text in Inkscape, but the package 'transparent.sty' is not loaded}%
    \renewcommand\transparent[1]{}%
  }%
  \providecommand\rotatebox[2]{#2}%
  \newcommand*\fsize{\dimexpr\f@size pt\relax}%
  \newcommand*\lineheight[1]{\fontsize{\fsize}{#1\fsize}\selectfont}%
  \ifx\svgwidth\undefined%
    \setlength{\unitlength}{144.76819794bp}%
    \ifx\svgscale\undefined%
      \relax%
    \else%
      \setlength{\unitlength}{\unitlength * \real{\svgscale}}%
    \fi%
  \else%
    \setlength{\unitlength}{\svgwidth}%
  \fi%
  \global\let\svgwidth\undefined%
  \global\let\svgscale\undefined%
  \makeatother%
  \begin{picture}(1,0.80793992)%
    \lineheight{1}%
    \setlength\tabcolsep{0pt}%
    \put(0,0){\includegraphics[width=\unitlength,page=1]{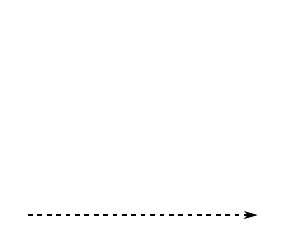}}%
    \put(0.66088054,0.0111169){\color[rgb]{0,0,0}\makebox(0,0)[lt]{\lineheight{1.25}\smash{\begin{tabular}[t]{l}$r\to \infty$\end{tabular}}}}%
    \put(0,0){\includegraphics[width=\unitlength,page=2]{fig150.pdf}}%
    \put(0.08302598,0.01307926){\color[rgb]{0,0,0}\makebox(0,0)[t]{\lineheight{1.25}\smash{\begin{tabular}[t]{c}$r=1$\end{tabular}}}}%
    \put(0,0){\includegraphics[width=\unitlength,page=3]{fig150.pdf}}%
  \end{picture}%
\endgroup%

	\caption{The union of the moduli spaces $\cM\cM^{r\cdot \ev^{\st/\un}}(\Psi)$ over different configurations $\Psi$. The $r$-regular subset is shown using thin-solid lines. The $r$-flat subset is shown with bold-solid lines.}\label{fig:150}
\end{figure}

We wish to count the ends of the moduli spaces $\cM\cM^{r\cdot \ev^{\st/\un}, r\in [1,\infty)} (\Psi)$, ranging over configurations $\Psi$ which satisfy
\begin{equation}
\mu(\psi_{c_3',c_2',c_1'})+\mu(\psi_{c_1,c_3',c_1'})+\mu(\psi_{\b,\b',\b''})+\mu(\psi_{\a,\b,\b''}) - 2n_{p_0}(\psi_{c_3',c_2',c_1'})+1=0. \label{eq:configurations-X-3-dim=1}
\end{equation}
These are exactly the configurations $\Psi$ where~\eqref{eq:r-dilated-moduli-space-range} has expected dimension 1. Using ~\eqref{eq:expected-dimension-asymptotic-eval-2}, we see that there are two ways a tuple $\Psi$ can have a 1-dimensional moduli space:
\begin{enumerate}[ref=$r$-\arabic*, label=($r$-\arabic*),leftmargin=*,widest=IIII]
\item\label{r-regular} \emph{$r$-regular}: 
\[
\mu(\psi_{c_1,c_3',c_1'})=\mu(\psi_{\b,\b',\b''})=0\quad \text{and} \quad \mu(\psi_{\a,\b,\b''})+\mu(\psi_{c_3',c_2',c_1'})=2n_{p_0}(\psi_{c_3',c_2',c_1'})-1.
\]
\item\label{r-flat} \emph{$r$-flat}: 
\[
\mu(\psi_{c_1,c_3',c_1'})+\mu(\psi_{\b,\b',\b''})=1\quad \text{and} \quad \mu(\psi_{\a,\b,\b''})+\mu(\psi_{c_3',c_2',c_1'})=2n_{p_0}(\psi_{c_3',c_2',c_1'})-2.
\]
\end{enumerate} For $r$-flat configurations $\Psi$, projection onto $r$ of $\cM\cM^{r\cdot \ev^{\st/\un}, r\in [1,\infty)}(\Psi)$ is locally constant. A schematic is shown in Figure~\ref{fig:150}.

We now enumerate the  ends of the moduli spaces $\cM\cM^{r\cdot \ev^{\st/\un}, r\in [1,\infty)}(\Psi)$ which occur at $r\in [1,\infty)$:

\begin{lem}\label{lem:counting-ends-M-Psi-finite-r} Suppose that $\Psi$ is a configuration of homology classes, decorated sources, and marked point matchings which satisfies ~\eqref{eq:configurations-X-3-dim=1}, and has incoming intersection points at a fixed $\xs\in \bT_{\a_1}\cap \bT_{\b_3}$, as well as the canonical inputs, as in~\eqref{eq:classes-for-E}. Let $k$ denote $n_{p_0}(\psi_{c_3',c_2',c_1'})$. Generically, the ends of $\cM\cM^{r\cdot \ev^{\st/\un}, r\in [1,\infty)}(\Psi)$ occurring  at finite  $r\in [1,\infty)$ correspond to the following configurations:
\begin{enumerate}[ref=$c$-\arabic*, label=($c$-\arabic*),leftmargin=*,widest=IIII]
\item\label{end:c-1-finite-r} The curves at $r=1$.
\item\label{end:c-2'-finite-r} (Disk breaking) A 5-tuple of the form $(u,v_{\b,\b',\b''}, v_{\a,\b,\b''},v_{c_3',c_2',c_1'}, v_{c_1,c_3',c_1'})$ where $u$ is a holomorphic disk of Maslov index 1, in any of the twelve cylindrical ends of the 4-manifolds corresponding to the triple diagrams $(\Sigma,\bs,\bs',\bs'')$, $(\Sigma,\as,\bs,\bs'')$, $(\bT^2, c_3',c_2',c_1')$ and $(\bT^2,c_1,c_3',c_1')$. The disk $u$ has no coupled marked points, while $v_{c_3',c_2',c_1'}$ and $v_{\a,\b,\b''}$ each have exactly $k$ coupled marked points. Furthermore, $v_{c_3',c_2',c_1'}$ has no decoupled marked points. Finally, 
\[
\mu(v_{\b,\b',\b''})=\mu(v_{c_1,c_3',c_1'})=0\quad \text{and} \quad \mu(v_{c_3',c_2',c_1'})+\mu(v_{\a,\b,\b''})=2k-2.
\]
\item\label{end:c-3'-finite-r} (Collisions of evaluation) A quadruple  $(v_{\b,\b',\b''}, v_{\a,\b,\b''}',b,v_{c_3',c_2',c_1'}, v_{c_1,c_3',c_1'})$, where $v_{\b,\b',\b''}$, $v_{c_3',c_2',c_1'}$, $v_{c_1,c_3',c_1'}$ are holomorphic triangles satisying the triangular analogs of~\ref{M-1}--\ref{M-7}. The curve $v_{\a,\b,\b''}'$, satisfies \ref{M-1}--\ref{M-7}, except that it has one additional interior puncture, which matches $p\in \Sigma$ to first order. Furthermore $(\pi_{\Delta}\circ v_{\a',\b,\b''})(q)$  is contained in the special line of $\Delta$, and also corresponds with the image of the constant bubble $b$. The asymptotic evaluations of the coupled marked points of $(v_{\a,\b,\b''}', b)$ and $v_{c_3',c_2',c_1'}$ coincide, and furthermore have exactly one repeated entry. Furthermore, the tuple $(\pi_{\Delta}\circ v_{c_3',c_2',c_1'})(\qs_{c_3',c_2',c_1'})$ consists of $k$ distinct points. Finally,
\[
\mu(v_{\b,\b',\b''})=\mu(v_{c_1,c_3',c_1'})=0\quad \text{and} \quad \mu(v_{c_3',c_2',c_1'})+\mu(v_{\a,\b,\b''})=2k-1.
\]
\end{enumerate}
\end{lem}

 \begin{figure}[ht!]
	\centering
	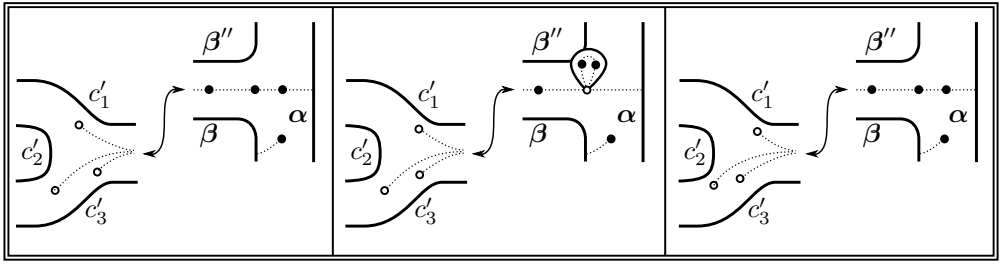
	\caption{A collision of evaluation~\eqref{end:c-3'-finite-r} (center) appearing in the boundaries of two families of $r$-matched moduli spaces (left and right).}\label{fig:207}
\end{figure}

\begin{proof} The degenerations labeled~\eqref{end:c-1-finite-r},~\eqref{end:c-2'-finite-r} and ~\eqref{end:c-3'-finite-r} occur in the correct codimension, so it suffices to prove that there are no additional ends. To prove this, we first describe which degenerations occur generically in the ends of 1-dimensional $r$-matched moduli spaces, and then use several facts about $\cQ_3$ to prohibit degenerations beyond those in the statement.

We consider a degeneration at some $r\in (1,\infty)$, consisting of a pair of combs, which together consist of two curves on $\Sigma\times \Delta$, two curves on $\bT^2\times \Delta$, as well as a collection of boundary degenerations and curves in the cylindrical ends. To analyze the limiting curves, it is helpful to combine the level structures in the $(c_3',c_2')$ and $(c_2',c_1')$ ends into a single level structure. We view limiting curves in this end as mapping into $\bT^2\times ([0,1]\times \R\sqcup [0,1]\times \R)$.

We may assume that a comb is the limit of a $r$-regular family, since $r$-flat families involve only degenerations of the curves $v_{\a,\b,\b''}$ and $v_{c_3',c_2',c_1'}$. For an $r$-flat family, one of $v_{\a,\b,\b''}$ and $v_{c_3',c_2',c_1'}$ has Maslov index 1, while the other has index 0. Hence, all degenerations correspond to an index 1 disk breaking off, which contains no coupled marked points; this is a special case of~\eqref{end:c-2'-finite-r}. Hence, it is sufficient to consider limiting combs involving only degenerations of $v_{c_3',c_2',c_1'}$ and $v_{\a,\b,\b''}$.

Next, we note that any holomorphic disk appearing in the $(\as,\bs)$ or $(\as,\bs'')$ end  has no coupled marked points, by definition of the vector field $v_3$, and hence must have Maslov index  1 (as the remaining components would otherwise violate transversality). It is easy to see that the resulting configuration satisfies~\eqref{end:c-2'-finite-r}.

 We now consider degenerations which involve non-trivial holomorphic disks breaking off in the $(\bs,\bs'')$, $(c_3',c_2')$, $(c_2',c_1')$ or $(c_3',c_1')$ ends, or non-trivial boundary degenerations. Suppose that $\cU$ and $\cV$ are the combs appearing in such a degeneration, and let $\cU$ be on the $\bT^2$-side, and $\cV$ be on the $\Sigma$-side. As described above, we may assume that this degeneration was the limit of an $r$-regular family. For notational convenience, we trim  $v_{c_1,c_3',c_1'}$ and $v_{\b,\b',\b''}$ (which are index 0 triangles) from $\cU$ and $\cV$. Write $\cU_1,\dots, \cU_n$ for the levels of $\cU$, and $\cV_1,\dots, \cV_m$ for the levels of $\cV$. See Figure~\ref{fig:204}.

 \begin{figure}[ht!]
	\centering
	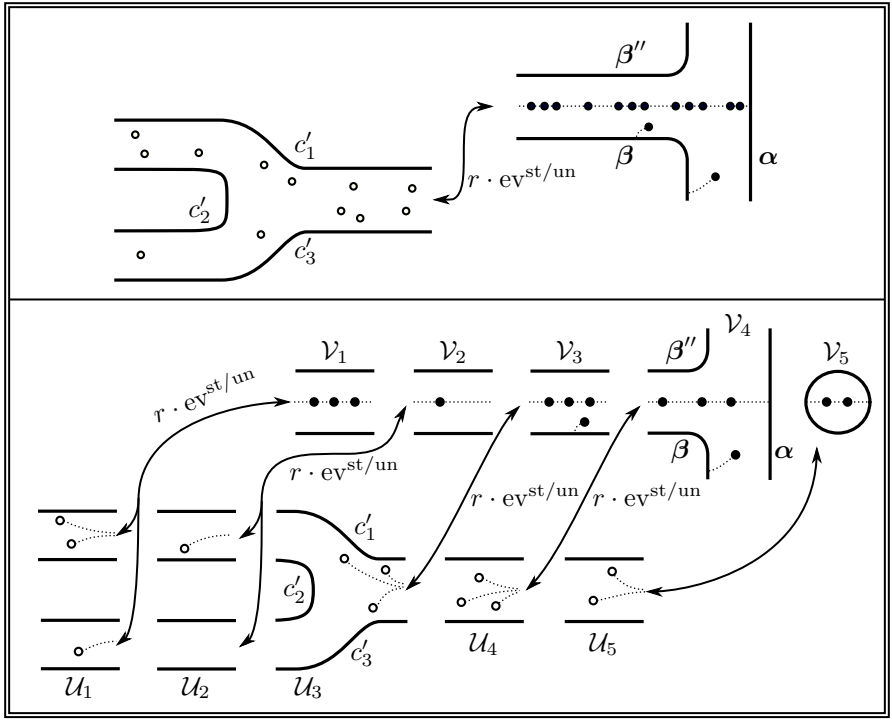
	\caption{An $\ev^{\st/\un}$-matched pair $(v_{\a,\b,\b''},v_{c_3',c_2',c_1'})$ (top) and a nearby comb (bottom). The pictured degeneration occurs in high codimension.}\label{fig:204}
\end{figure}  

 Let $d_i$ denote the number of coupled marked points of $\cU_i$, and let 
 \[
 \ve{d}=\sum_{i=1}^n d_i.
 \]
  Define 
 \[
 \bI:=\{i : 1\le i\le m, \quad d_i>0\}.
 \]
 The matching of marked points descends to an injection
\[
f\colon \bI\to \{1,\dots, n\}.
\]
When neither $\cU_i$ nor $\cV_{f(i)}$ is a boundary degeneration, the curves $\cU_i$ and $\cV_{f(i)}$ are matched under $r\cdot \ev^{\st/\un}$. (When one or both of $\cU_i$ or $\cV_{f(i)}$ is a boundary degeneration, there is still a matching condition which is satisfied, though it is more subtle and not important for our purposes).

Write $\phi_1^\bT,\dots, \phi_n^\bT$ for the homology classes of $\cU_1,\dots, \cU_n$, and $\phi_1^\Sigma,\dots, \phi_m^\Sigma$ for the homology classes of $\cV_1,\dots, \cV_m$. For each $i\in \bI$, we define
\[
\kappa_i:=\mu(\phi_i^\bT)+\mu(\phi_{f(i)}^\Sigma)-2d_i.
\]
Since $\cU$ and $\cV$ are the limit of a family of curves representing a pair $(\psi_{\a,\b,\b''}, \psi_{c_3',c_2',c_1'})$ with
\[
\mu(\psi_{\a,\b,\b''})+\mu(\psi_{c_3',c_2',c_1'})-2 \ve{d}=-1,
\]
we obtain the formula
\begin{equation}
\sum_{i\in \bI} \kappa_i +\sum_{i\in \{1,\dots, n\}\setminus \bI} \mu(\phi_i^{\bT})+\sum_{j\in \{1,\dots, m\}\setminus f(\bI)} \mu(\phi_j^\Sigma)=-1.  \label{eq:indices-sum-to--1}
\end{equation}
By transversality, if $i\in \{1,\dots, n\}\setminus \bI$ and $\phi_i^{\bT}$ is a non-trivial homology class of $g_i^\bT$-gons, then 
\begin{equation}
\mu(\phi_i^{\bT})\ge 3-g_i^{\bT}. \label{eq:lower-bound-decoupled-classes}
\end{equation}
An analogous lower bound holds for each $j\in \{1,\dots, m\}\setminus f(\bI)$.

We now claim that  for each $i\in \bI$, if $\cU_i$ is a holomorphic $g_i^{\bT}$-gon and $\cV_{f(i)}$ is a $g_i^\Sigma$-gon, then
\begin{equation}
\kappa_i\ge \begin{cases} 5-g_i^{\bT}-g_i^{\Sigma} -\delta(i)&  \text{if}\quad g_i^{\bT}>1 \quad  \text{and} \quad g_i^\Sigma>1, \\
3-g_i^{\Sigma}+d_i  & \text{if}\quad  g_i^{\bT}=1 \quad \text{and} \quad g_i^{\Sigma}>1,\\
3-g_i^{\bT}& \text{if} \quad g^{\bT}_i>1  \quad \text{and} \quad g^{\Sigma}_i=1,\\
2d_i& \text{if} \quad g_i^\Sigma=g_i^{\bT}=1.
\end{cases}
\label{eq:kappa-i-apriori-bound}
\end{equation}
for some function $\delta\colon \bI\to \{0,1\}$, satisfying
\[
\sum_{i\in \bI}\delta(i)\le 1.
\]
We note that the bound in~\eqref{eq:kappa-i-apriori-bound} is not dependent on the quadruple $\cQ_3$ in any way.

We now prove~\eqref{eq:kappa-i-apriori-bound}. Suppose first that $g_i^{\bT},g_i^{\Sigma}>1$. If $r$ is fixed, then  the moduli space of $r\cdot \ev^{\st/\un}$-matched pairs of holomorphic curves representing the classes $\phi_i^{\bT}$ and $\phi_{f(i)}^{\Sigma}$ has expected dimension $\kappa_i+1$. The expected dimension ranging over $r\in (1,\infty)$ is greater by 1.  If an $r\cdot \ev^{\st/\un}$-matched moduli space is non-empty for some $r$, then it must have dimension at least $(3-g^{\Sigma}_{i})+(3-g_i^{\bT})$, since the map $\ev^{\st/\un}$ is invariant under the translation action on the moduli spaces of disks. 

Note that if $\kappa_i+2= (3-g_i^\Sigma)+(3-g_i^\bT)$, such configurations do not appear for generic $r$.  It follows that there is at most one $i\in \bI$ such that $(\cU_i,\cV_{f(i)})$ contains no boundary degenerations and the $r$-matched moduli space generated by the homology classes of $\cU_i$ and $\cV_{f(i)}$, parametrized by $r\in (1,\infty)$, has expected dimension $6-g_i^{\bT}-g_i^\Sigma$. This pair corresponds to the $i\in \bI$ with $\delta(i)>0$. The remaining $i\in \bI$ such that $(\cU_i,\cV_{f(i)})$ has no boundary degenerations must live in $r$-matched moduli spaces (with $r$ ranging over $(1,\infty)$) which are 1-dimensional. Since $\kappa_i+1$ is the expected dimension, this immediately gives the first line of~\eqref{eq:kappa-i-apriori-bound}.

Suppose now that $g_i^{\bT}=1$ while $g_i^{\Sigma}>1$, so that $\cU_i$ is a boundary degeneration. By definition,  $\kappa_i=\mu(\phi_{f(i)}^\Sigma)$. On the other hand, the curve $\cV_{f(i)}$ contains a disk or triangle where the $d_i$ coupled marked points are mapped to the special line. The subspace of $\cM(\phi_i^{\Sigma})$ which satisfies this constraint has expected dimension $\mu(\phi_{f(i)}^\Sigma)-d_i.$ Hence, for a representative to exist, we must have $\kappa_i=\mu(\phi_i^\Sigma)\ge (3-g_i^\Sigma)+d_i,$ as claimed.  This gives the second line of~\eqref{eq:kappa-i-apriori-bound}.

Suppose now that $g_i^\Sigma=1$. In this case, $\kappa_i=\mu(\phi_i^{\bT})$. Since $\phi_i^\bT$ has a holomorphic representative, we know that $\mu(\phi_{f(i)}^{\bT})\ge 3-g_i^\bT$. This gives the third line of~\eqref{eq:kappa-i-apriori-bound}. If $g_i^\Sigma=g_i^{\bT}=1$, then the fourth line of ~\eqref{eq:kappa-i-apriori-bound} is in fact an equality.

In the setting of a general Heegaard quadruple, ~\eqref{eq:kappa-i-apriori-bound} implies that the generic degenerations of 1-dimensional $r\cdot \ev^{\st/\un}$-matched moduli spaces which are parametrized over $r\in (1,\infty)$, and also involve a non-trivial disk or boundary degeneration, must consist of one of the following  (see Figure~\ref{fig:205}):
\begin{enumerate}[ref=$g$-\arabic*, label=($g$-\arabic*),leftmargin=*,widest=IIII]
\item\label{g-c-1} A pair of holomorphic triangles which together have $\kappa_i=-2$, as well as a Maslov index 1 disk in the $(c_3',c_1')$-end, whose coupled marked points are matched to the coupled marked points on a $(\Sigma,\as)$ boundary degeneration. The disk and boundary degeneration together have $\kappa_i=1$.
\item\label{g-c-2} A pair of holomorphic triangles which together have $\kappa_i=-1$, as well as a pair of disks with a non-zero number of marked points, which together have $\kappa_i=0$.  One disk occurs in the $(c_3',c_2')$-$(c_2',c_1')$-end (recall we have collapsed the level structure over these two ends), while the other occurs in the $(\bs,\bs'')$-end.
\item\label{g-c-3} A pair of holomorphic triangles which together have $\kappa_i=-2$, as well as an index 1 disk in any end, which contains no coupled marked points.
\end{enumerate}
We remark that in~\eqref{g-c-1}, we have not ruled out the boundary degeneration from containing further ghost curve levels (e.g. resulting from the collision of two marked points). 

 \begin{figure}[ht]
	\centering
	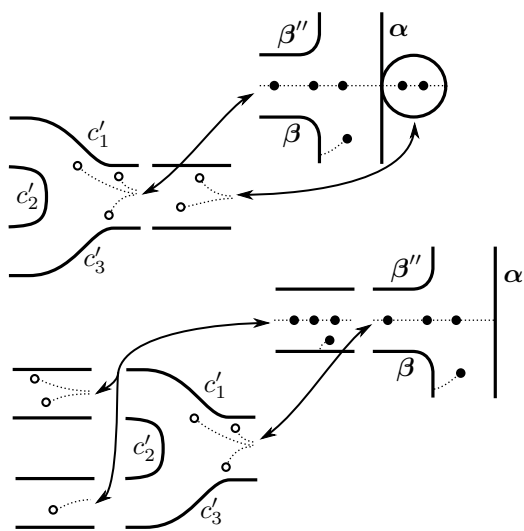
	\caption{Codimension 1 degenerations in Lemma~\ref{lem:counting-ends-M-Psi-finite-r} that are prohibited for $\cQ_3$.}\label{fig:205}
\end{figure}  

We claim that the degenerations~\eqref{g-c-1} and~\eqref{g-c-2} are non-generic for $\cQ_3$. The necessary extra fact is that any class of disks $\psi$ on one of the three ends of $(\bT^2,c_3',c_2',c_1')$ satisfies
\begin{equation}
\mu(\psi)=2n_{p_0}(\psi).
\label{eq:index-formula-disk-c-1-c-3}
\end{equation}

We consider first~\eqref{g-c-1}. Suppose a limiting comb had a matched level $(\cU_i,\cV_{f(i)})$ such that $\cU_i$ is a holomorphic disk, and $\cV_{f(i)}$ is a boundary degeneration. From~\eqref{eq:index-formula-disk-c-1-c-3} and the definition of $\kappa_i$, we obtain
\[
\kappa_i=2 d_i,
\]
which can never be 1.

 Next, we consider~\eqref{g-c-2}. Suppose a limiting comb had a matched level $(\cU_i,\cV_{f(i)})$ consisting of a holomorphic disk in the $(c_3',c_2')$-$(c_2',c_1')$-end, matched to a disk in the $(\bs,\bs'')$-end. Furthermore, suppose that $\kappa_i=0$. Combining the definition of $\kappa_i$ with~\eqref{eq:index-formula-disk-c-1-c-3}, we obtain $\kappa_i=\mu(\phi_{f(i)}^\Sigma)$, which must be at least 1 for $\phi_{f(i)}^\Sigma$ to have a holomorphic representative. We obtain a contradiction, so ends of the form~\eqref{g-c-2} are prohibited. Hence, ~\eqref{g-c-3} is the only generic degeneration for $\cQ_3$ which involves non-trivial disk or boundary degeneration levels. We note~\eqref{g-c-3} corresponds to~\eqref{end:c-2'-finite-r}.

We now consider limiting combs which do not involve a holomorphic disk or boundary degeneration breaking off. Such a limiting comb must contain a ghost curve (i.e. a curve which has constant value), or the formation of a double point. Limiting combs containing a ghost curve and double points on the $\bT$-side occur in codimension at least 2, and hence do not feature in codimension 1 degenerations featuring in the statement. This is proven as follows. A standard argument \cite{LOTBordered}*{Lemma~5.57} implies that double points and ghost curves which do not contain a marked point occur in codimension at least 2. Ghost curves containing one coupled marked point are prohibited by identical reasoning. Ghost curves containing more than one coupled marked point are prohibited by combining \cite{LOTBordered}*{Lemma~5.57} with the expected dimensions of moduli spaces with tangency constraints using Proposition~\ref{prop:dimension-counts}. 
 Double points on the $\Sigma$-side are prohibited by similar reasoning.

Ghost degenerations on the $\Sigma$-side are slightly more subtle, since the coupled marked points are already matched to the special line. Ghost curves containing 0 or 1 marked points are prohibited by identical reasoning to the $\bT$-side. More generally, suppose $\cU_i$ and $\cV_{f(i)}$ are a matched level, as above, and $\cV_{f(i)}$ contains a ghost curve with $\ell>1$ marked points. Suppose further that $\cV_{f(i)}$ is not a boundary degeneration. We can trim off the ghost curve, and fill in the remaining puncture with a marked point, to obtain a curve with a marked point which matches $p\in \Sigma$ to order $\ell-1$, and also matches the special line. An easy computation using~\eqref{eq:kappa-i-apriori-bound} shows that we may increase the lower bound of $\kappa_i$ appearing in~\eqref{eq:kappa-i-apriori-bound} by $\ell-1$. There is exactly one configuration which satisfies this modified~\eqref{eq:kappa-i-apriori-bound}, which is the case that $\cU_i$ and $\cV_{f(i)}$ are both triangles, $\ell=2$, and there are no more levels.   We restricted above to the case where the topological source of the ghost was a sphere.  Ghost curves with more topologically complicated sources can be seen to live in codimension at least 2, by adapting~\cite{LOTBordered}*{Lemma~5.57}.
\end{proof}

The ends labeled~\eqref{end:c-1-finite-r} correspond to the expression
\[
X_3^\infty(\xs).
\]
Amongst the ends labeled ~\eqref{end:c-2'-finite-r}, disks breaking into the ends $(\as,\bs)$ and $(\as,\bs'')$ correspond exactly to
\[
\left(\d H_3^{\infty, [1,\infty)}+H_3^{\infty, [1,\infty)} \d\right)(\xs).
\]
We claim that the remaining ends occuring at finite $r\in (1,\infty)$ make vanishing total contribution:

\begin{lem}
 The total weighted count of the broken curves  labeled ~\eqref{end:c-2'-finite-r}, not occuring in the $(\as,\bs)$ or $(\as,\bs'')$-cylindrical ends, as well as the broken curves labeled~\eqref{end:c-3'-finite-r}, vanishes. 
\end{lem}
\begin{proof}

 The total count of broken curves with the configuration~\eqref{end:c-2'-finite-r} where the disk occurs in one of the cylindrical ends labeled $(\bs,\bs')$, $(\bs',\bs'')$, $(c_3',c_2')$, $(c_2',c_1')$ or $(c_1,c_3')$ or $(c_1,c_1')$ vanishes since the special inputs (i.e. $\Theta_{c_1,c_3'}$, $\Theta_{\b,\b'}^+$, etc.) are cycles. 
 
 The remaining configurations from~\eqref{end:c-2'-finite-r} corresponds to a Maslov index 1 holomorphic disk, with no coupled marked points, breaking off in one of the $(\bs,\bs'')$ cylindrical ends. These degenerations occur in canceling pairs: if $\cU=(v_{\a,\b,\b''}, v_{\b,\b',\b''}, v_{c_3',c_2',c_1'}, v_{c_1,c_3',c_1'},u)$ is such a limiting comb, then we may glue $u$ to either $v_{\b,\b',\b''}$ or $v_{\a,\b,\b''}$ to obtain a 1-dimensional moduli space with $\cU$ in its boundary.   (Note that such ends correspond to transitions between the $r$-flat and $r$-regular  moduli spaces).

We now consider ~\eqref{end:c-3'-finite-r}. Suppose the punctures $q_1$ and $q_2$ bubble off in a once punctured sphere, denoted $T_{0}^{q_1,q_2}$.  We write $v_{\a,\b,\b''}'$ for the holomorphic triangle on $(\Sigma,\as,\bs,\bs'')$ obtained by trimming off the bubble, and we write $\qs$ for the set of coupled marked points on this level (including $q_1$ and $q_2$). We write $v_{c_3',c_2',c_1'}$ for the holomorphic triangle on $(\bT^2,c_3',c_2',c_1')$, and $\qs_0$ for the coupled marked points if $v_{c_3',c_2',c_1'}$. There are two marked points, $q_1',q_2'\in \qs_0$, which are matched with $q_1$ and $q_2$, respectively. In the formation of the limiting comb, we may rescale the $\Sigma$-component near $p$, and also rescale the $\Delta$ component, to obtain a limiting curve $v_0\colon T_0^{q_1,q_2}\to \C\times \C$. We view $S^2\times S^2$ as the compactification of $\C\times \C$. Note that we can complete $v_0$ to a map from $S^2$ into $S^2\times S^2$.

Let us write $\pi_1\colon \C\times \C\to \C$ for projection onto the first factor, and $\pi_2$ for projection onto the second.

The following two properties of $v_0$ are immediate:
\begin{enumerate}[ref=$v_0$-\arabic*, label=($v_0$-\arabic*), leftmargin=*, widest=IIII]
\item\label{property:bubble-1} $v_0$ maps $q_1$ and $q_2$ to $\{0\}\times \R$ (here we view $\R\subset \C$ as the image of the special line in $\Delta$ after rescaling).
\item\label{property:bubble-2} $v_0$ maps the puncture of $T_0$ to $(\infty,\infty)\in S^2\times S^2$.
\end{enumerate}
We claim that generically $v_0$ will also have the following properties:
\begin{enumerate}[ref=$v_0$-\arabic*, label=($v_0$-\arabic*), leftmargin=*, widest=IIII]
\setcounter{enumi}{2}
\item\label{property:bubble-3} The map $\pi_2\circ v_0$ has degree 1.
\item\label{property:bubble-5} If we identify $(T_0,q_0,q_1)$ with $(\C,-1,1)$, then the map $\pi_1\circ v_0$ takes the form
\[
(\pi_1\circ v_0)(z)=\alpha (z+1)(z-1)
\]
for some $\alpha\in \C^\times$.
\end{enumerate}
Property~\eqref{property:bubble-5} follows from \eqref{property:bubble-3}, as well as~\eqref{property:bubble-1}, ~\eqref{property:bubble-2} and basic complex analysis. Property~\eqref{property:bubble-3} is proven as follows.  We can view the source of $v_{\a,\b,\b''}'$ as having a marked point $q$, which is mapped to a point $(p,d)\in \Sigma\times \Delta$, which is the same point where the ghost bubble maps. If the degree of $\pi_2\circ v_0$ was greater than 1, then  $v_{\a,\b,\b''}'$ would have the property that $(\pi_\Sigma\circ v_{\a,\b,\b''}')(q)=p$ to order 1,  $(\pi_\Delta\circ v_{\a,\b,\b''}')(q)$ is in the special line, and $\pi_{\Delta}\circ v_{\a,\b,\b''}'$ has vanishing derivative. Using the dimension counts of Proposition~\ref{prop:dimension-counts}, we see that the expected dimension of the moduli space of such curves is too small (it drops by 3 instead of by 1), so does not appear in  codimension 1 degenerations.

 Let us write $\cU'$ for the four triangular components of the limit (with the ghost trimmed off). We view $(\cU',v_0)$ as being in the compactification of the moduli space  $\cM\cM^{r\cdot \ev^{\st/\un}, r\in [1,\infty)}(\Phi)$, for some collection of homology classes with marked sources and matching data, which we denote by $\Phi$.   There is exactly one other configuration $\Phi'$ which contains a canceling broken curve of the form $(\cU',v_0')$. The curve $v_0'$ is the same as $v_0$, except with the matching of $(q_1,q_2)$ and $(q_1',q_2')$ switched.
 \end{proof}

 It remains to show that each of the broken curves described above make odd algebraic contribution to the count of the ends of the matched moduli spaces.

\begin{lem}\label{lem:gluing-bubbles-1}
 Near a broken curve labeled \eqref{end:c-2'-finite-r} or~\eqref{end:c-3'-finite-r}, there is a neighborhood in the compactification of the matched moduli spaces which admits a proper continuous map to $[0,\infty)$ which is locally odd degree and proper, such that the preimage of $0$ is the broken curve. 
\end{lem} 

 \begin{proof}
 We focus on the ends labeled~\eqref{end:c-3'-finite-r}, since the ends labeled~\eqref{end:c-2'-finite-r} can be handled by a similar argument. Let us write $(\cU',v_0)$ for a broken curve satisfying~\eqref{end:c-3'-finite-r}.

Since the source of $v_0$ is a once punctured sphere, the map $v_0$ achieves transversality for any split almost complex structure; see \cite{McDuffSalamonSymplectic}*{Lemma~3.3.1}. (Compare \cite{LOTBordered}*{Proposition~5.16}).

 We first establish some notation. First, trim off the curves $v_{\b,\b',\b''}$ and $v_{c_1,c_3',c_1'}$, since they are index 0 and have no constraint. Write $v_{\a,\b,\b''}'$ for the limiting triangle on $(\Sigma,\as,\bs,\bs'')$ with the bubble trimmed off. By filling in the order two Reeb orbit with a marked point, we can view $v_{\a,\b,\b''}'$ as being an element of the moduli space of holomorphic triangles representing $\psi_{\a,\b,\b''}$, which have $n-1$ marked points, $n-2$ of which are order 1 at $p$, and $1$ of which is order two, but which have no constraint with regards to the special line. We write $W_0$ for a small neighborhoof of $v_{\a,\b,\b''}'$ in this moduli space.

 We write $\cN_2$ for the moduli space of holomorphic maps $v_0\colon T_0\to S^2\times S^2$ with 2 marked points, and one puncture, which satisfy the following conditions. We assume that $\pi_1\circ v_0$ is asymptotic to $\infty$ with order 2 at the puncture of $T_0$, and $(\pi_1\circ v_0)(q_i)=0$ for $i\in \{1,2\}$. Finally, we assume that $\pi_2\circ v_0$ is degree 1. 
 
 We write $G$ for the group of transformations of $\C$ given by $\{z\mapsto a z+b: a>0, b\in \C\}$. We write $\tilde{\cN}_2$ for the quotient of $\cN_2$ by $G$ (where $G$ acts on the second factor of a map). Note that complex analysis identifies $\tilde{\cN}_2$ with $S^1$.

Standard gluing results imply that the above neighborhood $W_0$ of $v_{\a,\b,\b''}'$ may be chosen so that there is a map
\[
\g\colon W_0\times \tilde{\cN}_2\times [0,\veps)\to \bar{W},
\] 
which is a homeomorphism onto its image.  Here, $\bar{W}$ denotes a neighborhood of $(v_{\a,\b,\b''}',v_0)$ in the compactification of holomorphic triangles which have $n$ marked points which each match the point $p\in \Sigma$ to order 1. We write $W$ for the interior of $\bar{W}$. The factor $[0,\veps)$ is the gluing parameter.

 Write $U$ for a neighborhood of $v_{c_3',c_2',c_1'}$ inside the moduli space of triangles representing $\psi_{c_3',c_2',c_1'}$ which have $n$ marked points which match $p_0\in \bT^2$ to order 1.
 
Let $r_0$ denote the value of $r$ where the bubble forms. Note that we can view the evaluation map $\ev^{\st}$ on $U$ as taking values in $\Delta^n/\R$, where $\R$ acts diagonally by the flow of $v_{3,1}$. We can pick a diffeomorphism of $\Delta$ with $\C$, such that the time $t$ flow of $v_{3,1}$ coincides with the map $z\mapsto z+t$. Hence, we will view the evaluation map as taking values in $\C^n/\R$, and we assume the special line is the real axis. Similarly, $r\cdot \ev^{\st}$ can be thought of as having values in the same space, but taking values  in $\R^n/\R\subset \C^n/\R$.

Define 
\[
\ev_{q_1,q_2}:=\re(\ev_{q_2})-\re(\ev_{q_1}),
\]
and write $\bar{\ev}_{q_1,q_2}$ for the continuous extension to the compactification. Note that near $(\cU',v_0)$, $\bar{\ev}_{q_1,q_2}$ is either nonnegative or nonpositive (depending on the ordering of the marked points along the special line), and the only zero is $(\cU',v_0)$. For concreteness, assume $\bar{\ev}_{q_1,q_2}\ge 0$. We will show that the map $\bar{\ev}_{q_1,q_2}$ is locally degree 1 near $(\cU',v_0)$, which will prove the main claim.

Consider the diagram
\begin{equation}
\begin{tikzcd}[row sep=.5cm]
\R \times U
	\ar[dr, "\ev_1",swap]
&&
\R\times \bar{W}
	\ar[dl, "\ev_2"]
&
 \R\times W_0 \times \tilde{\cN}_2\times  [0,\veps)
	\ar[l, "\id\times \g", "\sim"'] 
 \\
& \R\times \C^n/ \R
\end{tikzcd}
\label{eq:gluing-diagram-1}
\end{equation}
where
\[
\ev_1(r,u)=(r, r\cdot \ev^{\st}(u))\quad \text{and} \quad\ev_2(r,v)=(r,\bar{\ev}(v)).
\]

Near $(\cU',v_0)$, the moduli space $\cM\cM^{r\cdot \ev^{\st/\un}, r\in [1,\infty)}(\Phi)$ is the fibered product of $\ev_1$ and $\ev_2$.  We define $\cI_U^{\ge 0}$ (resp. $\cI_U^{>0}$ or $\cI_U^{0}$) to be subset of $\R\times U$ where $\ev_{q_1,q_2}\ge 0$ (resp. $\ev_{q_1,q_2}>0$ or $\ev_{q_1,q_2}=0$).
For generically chosen almost complex structures, the following two statements follow from an adaptation of Proposition~\ref{prop:dimension-counts}:
\begin{itemize}
\item The maps $\ev_1$ and $\ev_2$ are transverse, when restricted to $\cI^{>0}_{U}$ and $\R\times W$, respectively.
\item The maps $\ev_1$ and $\ev_2$ are transverse, when restricted to $\cI^0_U$ and $\R\times W_0$ are transverse, and have $\cU'$ as an intersection of their images.
\end{itemize}
Finally, we note that $\cM\cM^{r\cdot \ev^{\st/\un} r\in [1,\infty)}(\Phi)$ is identified with the intersection of the images of  $\cI_U^{\ge 0}$ and $\R\times \bar{W}$.
Lemma~\ref{lem:linking} can be used to show that the map $\bar{\ev}_{q_1,q_2}$ is locally proper and odd degree near $0$ on $\cM\cM^{r\cdot \ev^{\st/\un} r\in [1,\infty)}(\Phi)$, completing the proof. 
\end{proof}

It remains to count the ends which appear as $r\to \infty$.

\begin{lem}\label{lem:ends-r-to-infty} The ends appearing as $r\to\infty$ consist exactly of 4-tuples of curves 
\[
(v_{c_3',c_2',c_1'}, v_{c_1,c_3',c_1'}, v_{\b,\b',\b''}, v_{\a,\b,\b''}),
\] as follows.
\begin{enumerate}
\item The curves $v_{c_3',c_2',c_1'}$, $v_{c_1,c_3',c_1'}$, $v_{\b,\b',\b''}$ all have Maslov index 0, and have no marked point constraints.
\item If $k=n_{p_0}(v_{c_3',c_2',c_1'})>0$, then $v_{\a,\b,\b''}$ has Maslov index $2k-1$. Furthermore, $v_{\a,\b,\b''}$ has exactly one coupled marked point (and any number of decoupled marked points). This marked point projects to the special line, and matches $p$ to order $k-1$. If $k=0$, then $v_{\a,\b,\b''}$ has Maslov index 0 and no constraint.
\end{enumerate}
 Furthermore,
\begin{enumerate}[ref=$g$-\arabic*, label=($g$-\arabic*), leftmargin=*, widest=IIII]
\item \label{r->infty:gluing-1}   The total count of such broken curves, weighted by their appropriate powers of $U$ and $T$, is 0.
\item \label{r->infty:gluing-2}
Each of the above broken curves has a neighborhood in the compactification of $r\cdot \ev^{\st/\un}$-matched curves which admits a map to $[0,1)$ which is degree 1 near 0.
\end{enumerate}
\end{lem}
\begin{proof} Similar to Lemma~\ref{lem:counting-ends-M-Psi-finite-r}, we first give a general description, and then specialize to the quadruple $\cQ_3$. The argument therein adapts without substantial complication to show that the limiting curves as $r\to \infty$ consist of the following configurations:
\begin{itemize}
\item An index $0$ triangle on $(\Sigma,\bs,\bs',\bs'')$, with only decoupled marked points.
\item An index $0$ triangle on $(\bT^2,c_1,c_3',c_1')$, with only decoupled marked points.
\item An index $0$ triangle on $(\bT^2,c_3',c_2',c_1')$, with only coupled marked points.
\item A collection of Maslov index 1 disks in the $(c_1',c_3')$ and $(c_3',c_2')$-$(c_2',c_1')$-cylindrical ends of $(\bT^2,c_3',c_2',c_1')$. (An index 1 disk in the $(c_3',c_2')$-$(c_2',c_1')$ consists of an index 1 disk in either of the $(c_3',c_2')$ or $(c_2',c_1')$ ends, and a constant disk in the other). All marked points here are coupled. Suppose there are $n\ge 1$ levels here, and write $k_1,\dots,k_n$ for the number of marked points in each level.
\item A broken holomorphic triangle on $(\Sigma,\as,\bs,\bs'')$ of Maslov index $2(k_1+\cdots +k_n)-1$, with $n$ interior punctures, which each project to the special line. Furthermore, the $i^{\text{th}}$ puncture matches the point $p\in \Sigma$ to order $k_i$. Additionally, there is one ghost component for each $i\in \{1,\dots, n\}$ such that $k_i>1$. This ghost component consists of a sphere with one puncture, and $k_i$ marked points.
\end{itemize}
See Figure~\ref{fig:206} for a schematic. We note that in the limit, we can rescale any ghost bubble to obtain a map of a once punctured sphere into $\C\times \C$, as we did in Lemma~\ref{lem:counting-ends-M-Psi-finite-r}. For the component corresponding to $i\in \{1,\dots, n\}$, the $k_i$ marked points of this rescaled curve are mapped to $\{0\}\times \R\subset \C\times \C$. If we identify the source with $\C$, we may identify the marked points with numbers $\alpha_1,\dots, \alpha_{k_i}\in \R$. Furthermore $[(\alpha_{1},\dots, \alpha_{k_i})]\in \R^{k_i}/\R$ coincides up to an overall scalar factor with the asymptotic evaluation map, applied to the corresponding holomorphic triangle on $(\bT,c_3',c_2',c_1')$ or holomorphic disk mapping into one of its cylindrical ends. Gluing gives an identification of a neighborhood of such a configuration with $(0,1]$.

In the present case, that Maslov index formula from~\eqref{eq:index-formula-disk-c-1-c-3} implies immediately that $n=1$, since there are no holomorphic disks in the ends of $(\bT^2,c_3',c_2',c_1')$ which have Maslov index 1.

Claim~\eqref{r->infty:gluing-1}, the count of such broken curves, follows from Lemma~\ref{lem:model-computations-m-surgery}, since the holomorphic triangles on $(\bT^2,c_3',c_2',c_1')$ are unconstrained. Some care is required to verify that the weights of the curves cancel, since Lemma~\ref{lem:model-computations-m-surgery} holds with the algebra weights $U^{n_{p_0}}(\psi_0) T^{n_{z_0}(\psi_0)-n_{w_0}(\psi_0)}$, and for the curves in the limit, the punctures at $p_0$ do not contribute any $U$-weight or $T$-weight. Instead, a class $\psi_0$ is given $T$-weight $T^{n_{z_0}(\psi_0)}$ and no $U$-weight. We observe that the index 0 triangles on $(\bT^2,c_3',c_2',c_1')$ may be indexed by $s\in 2\Z+1$. The triangle $u_s$ has
\[
n_{w_0}(u_s)=\frac{m(s^2-1)}{8}\quad \text{and} \quad n_{z_0}(u_s)=\frac{m((s+2)^2-1)}{8}.
\]
In particular, 
\[
n_{z_0}(u_s)\equiv n_{z_0}(u_{-s}) \pmod{m}.
\]
In particular, the two triangles $u_s$ and $u_{-s}$ have the same multiplicity at $w_0$ and contribute the same $U$ and $T$-weight. Hence, the curve quadruples appearing as $r\to \infty$ will cancel when counted with $U$ and $T$ weights.

We remark that it is in this model computatation that we use the choice of the involution described in Remark~\ref{rem:careful-twisting}. Indeed if we had taken the connected sum at $z_0\in \bT^2$ instead of $w_0$, then the curves $u_{s}$ and $u_{-s-4}$ would have the same multiplicity at the puncture $z_0$, but these curves contribute different $U$-weights since the exponent of $U$ is not taken with values modulo $m$.

We now consider claim~\eqref{r->infty:gluing-2}. The proof is very similar to the proof of Lemma~\ref{lem:gluing-bubbles-1}. We let $\cN_n$ denote the moduli space of once punctured spheres with $n$ marked points, which map to $\C\times \C$, have an order $n$ puncture mapping to $(\infty,\infty)$, whose projection to the second factor is degree 1. We further assume that each marked point projects to $0$ in the first factor. The group $G$ of transformations of $\C$ given by $\{z\mapsto az+b: a>0, b\in \C\}$ acts on $\cN_n$, and we write $\tilde{\cN}_n$ for the quotient. We let $q_1$ and $q_2$ be two of the coupled marked points of the curve on $(\Sigma,\as,\bs,\bs'')$, and let $q_1'$ and $q_2'$ be the corresponding marked points of the curve of $(\bT^2,c_3',c_2',c_1')$. Write $(v_{\a,\b,\b''}',v_0)$ for the limiting triangle and marked bubble on $(\Sigma,\as,\bs,\bs'')$. Let $W_0$ denote a neighborhood of $v_{\a,\b,\b''}'$ in the moduli space of holomorphic curves representing $\psi_{\a,\b,\b''}$ which have an order $n$ marked points which match $p\in \Sigma$ to order 1, but which do not necessarily map to the special line. Standard gluing results give a homeomorphism 
\[
\g\colon W_0\times \tilde{\cN}_n\times [0,\veps)\to \bar{W},
\]
where $\bar{W}$ is a neighborhood of $(v_{\a,\b,\b''}',v_0)$ in the compactification triangles on $(\Sigma,\as,\bs,\bs'')$, with $n$ marked points which each match $p$ to order 1.

 \begin{figure}[ht!]
	\centering
	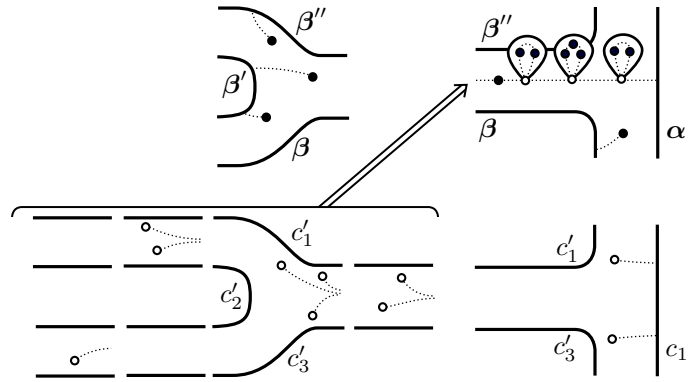
	\caption{A curve appearing as $r\to \infty$. For $\cQ_3$, the only generic configuration occurs when there are no disk levels.}\label{fig:206}
\end{figure}

We let $U$ denote a neighborhood of the limiting representative of $\psi_{c_3',c_2',c_1'}$. In the present case, $U$ consists of a single point, since $\psi_{c_3',c_2',c_1'}$ has Maslov index 0. In particular, we may assume that $\re(\ev_{q_1'}-\ev_{q_2'})$ is non-zero.
Let $V$ denote a neighborhood in $\tilde{\cN}_n$ of the bubble curve $v_0$.  We assume that on $V$, any curve satisfies $\re(\ev_{q_1}-\ev_{q_2}) \ge 0$.  Let $\bar{W}_V$ denote the image under $\gamma$ of $W_0\times V\times [0,\veps)$.

It is helpful to consider the parameter $r\in (0,\infty]$ as instead taking values in $[0,\infty)$, which we do by the change of variables $s=1/r$.
 We consider the diagram
\begin{equation}
\begin{tikzcd}[row sep=.5cm]
[0,\infty) \times U
	\ar[dr,swap, "\ev_1"]
&[-.3cm] &[-.3cm]
\R\times \bar{W}_V
	\ar[dl, "\ev_2"]
&
 \R\times W_0 \times V\times  [0,\veps)
	\ar[l, "\id\times \g", "\sim"'] 
 \\
& \R\times \C^n/G\times\R\times [0,\infty)
\end{tikzcd}
\label{eq:gluing-diagram-2}
\end{equation}
In~\eqref{eq:gluing-diagram-2},
\[
\ev_1(s,u)=(s,\ev^{\st}(u), 0, s\cdot \re(\ev_{q_1}(u)-\ev_{q_2}(u))),
\]
and
\[
\ev_2(s,v)=(s, \bar{\ev}(v), \im(\ev_{q_1'}(v)), \re(\ev_{q_1'}(v)-\ev_{q_2'}(v))).  
\]
 The images of the codimension 1 facets $\{0\}\times U$ and $\R\times \d\bar{W}_V$ have a transverse intersection point corresponding to $(v_{c_{3',c_2',c_1'}},v_{\a,\b,\b''},v_0)$, so Lemma~\ref{lem:linking} implies that the map $\re(\ev_{q_1'}-\ev_{q_2'})$ is locally proper and odd degree near 0 on the fibered product of the diagram in~\eqref{eq:gluing-diagram-2}, which one may easily show coincides with $\cM\cM^{r\cdot \ev^{\st/\un}, r\in [1,\infty)}(\Phi)$ near the limiting curve.
\end{proof}

Lemma~\ref{lem:H-3-null-homotopy-X-3} now follows from Lemmas~\ref{lem:counting-ends-M-Psi-finite-r} and~\ref{lem:ends-r-to-infty}. As a consequence, Proposition~\ref{prop:deform-tau-to-infinity} follows, and hence from Lemma~\ref{lem:E-infty-null-homotopic} we obtain Proposition~\ref{prop:Esimeq0}.

\section{The involutive  mapping cone formula and surgery exact sequence}
\label{sec:involutive-cone}

In this Section, we prove Theorem~\ref{thm:cone-intro} of the introduction.  Let 
\[
\bA=\prod_{s\in \Z} \bm{A}_s\quad \text{and} \quad \bB:=\prod_{s\in \Z} \bm{B}_s.
\]
 Ozsv\'{a}th and Szab\'{o} describe a map $D_n\colon \bA\to \bB$ so that $\bCF^-(Y_n)$ is quasi-isomorphic to $\bX_n:=\Cone(D_n\colon \bA\to \bB)$.

  It is also helpful to consider a smaller model of the mapping cone. If $b\in \N$, we define the \emph{horizontal truncation} of $\bX_n$, denoted $\bX_n\langle b\rangle $, to be the submodule
\begin{equation}
\bigoplus_{-b\le s\le b} \bm{A}_s\oplus  \bigoplus_{-b+n\le s\le b} \bm{B}_s. \label{eq:horizontal-truncation}
\end{equation}
We endow $\bX_n\langle b\rangle $ with a differential $D_n$ by viewing $\bX_n\langle b \rangle$ as a quotient of $\bX_n$ if $n>0$, and viewing $\bX_n\langle b\rangle$ as a subcomplex if $n<0$. We write  $\bA_n\langle b\rangle$ and $\bB_n\langle b \rangle$ for the two summands of $\bX_n\langle b \rangle$. 
 
 If $n>0$, then for large $b$, the quotient map
\[
\Pi\colon\bX_n\to \bX_n\langle b\rangle
\]
is a homotopy equivalence. For $n<0$, the inclusion from $\bX_n\langle b\rangle$ to $\bX_n$ is a homotopy equivalence for negative $n$. See \cite{OSIntegerSurgeries}*{Lemma~4.3}.
 
 In this section, we prove the following:

\begin{thm}\label{thm:mapping-cone} Suppose $n\in \Z\setminus \{0\}$, $Y$ is an integer homology 3-sphere, and $K\subset Y$ is a knot. The involutive Heegaard Floer homology $\bCFI^-(Y_n)$ is homotopy equivalent over the ring $\bF\llsquare U \rrsquare[Q]/Q^2$ to a complex $\bF[Q]/(Q^2)\otimes_\bF (\bA\oplus \bB)$  with differential specified by some $\bF\llsquare U \rrsquare$-linear maps as below:
\[
\bX \bI_n= 
 \begin{tikzcd}[column sep=2cm, row sep=2cm, labels=description]
\bA
	\ar[d, "Q\cdot (\id +\iota_{\bA})",swap]
	\ar[dr, dashed, "Q\cdot H"]
	\ar[r, "D_n"]
& \bB
	\ar[d, "Q\cdot (\id+\iota_{\bB})"]\\
Q\cdot \bA
	\ar[r, "D_n"]
& Q\cdot \bB.
\end{tikzcd}
\]
 In fact, the $\iota$-complex $(\bCF^-(Y_n), \iota)$ is $\iota$-homotopy equivalent to the complex $\bX_n:=\Cone\left(D_n\colon \bA\to \bB\right)$, with involution $\iota_{\bX}:=\iota_{\bA}+H+\iota_{\bB}$. Furthermore, the following hold:
\begin{enumerate}
\item \label{cone:1} $D_n=v+h$, where $v$ sends $\bm{A}_s$ to $\bm{B}_s$, and $h$ sends $\bm{A}_s$ to $\bm{B}_{s+n}$. 
\item\label{cone:2} $\iota_{\bA}$ sends $\bm{A}_s$ to $\bm{A}_{-s}$. Furthermore, $\iota_{\bA}$ coincides with $x\mapsto U^{s} \cdot  \iota_K(x)$, where $\iota_K$ is the knot involution.
\item\label{cone:4} $H=k+j$, where $k$ maps $\bm{A}_s$ to $\bm{B}_{-s+n}$ and $j$ maps $\bm{A}_s$ to $\bm{B}_{-s}$.
\end{enumerate}

The analogous claims hold for the truncated complexes $\bXI_n\langle b\rangle$, provided $b$ is sufficiently large.
\end{thm}

\begin{rem}
It is an easy consequence of Theorem~\ref{thm:mapping-cone} that the uncompleted complex $\CFI^-(Y_n(K))$ is $\iota$-homotopy equivalent to the truncated complex $\XI_n\langle b\rangle$, taken with coefficients in $\bF[U,Q]/Q^2$.  This holds since every free, finitely generated,  relatively $\Z$-graded chain complex $C$ over $\bF\llsquare U\rrsquare$ is chain isomorphic to $C'\otimes_{\bF[U]}\bF\llsquare U\rrsquare$, for some free, finitely generated, relatively $\Z$-graded chain complex $C'$ over $\bF[U]$. Furthermore, the chain complex $C'$ is unique up to chain isomorphism. Indeed, if $C$ and $D$ are free, finitely generated, relatively $\Z$ graded chain complexes over $\bF\llsquare U\rrsquare$, then $C$ and $D$ are homotopy equivalent over $\bF\llsquare U\rrsquare$ if and only if the corresponding complexes $C'$ and $D'$ are homotopy equivalent over $\bF[U]$.
\end{rem}

The complex $\bX \bI_n$ splits over conjugacy classes of $\Spin^c$ structures, as we now describe. Suppose $\varpi=\{\frs,\bar{\frs}\}$ is a conjugacy class of $\Spin^c$ structures on $Y_n$. Using the standard identification $\Spin^c(Y_n)\iso \Z_n$, we write $\varpi_i:=\{\frs_i,\bar{\frs}_i\}$, which consists of either one or two elements. We define
\[
\bA_{\{i,-i\}}:=\prod_{\substack{s\in \Z\\ s\equiv \pm i \mod n}} \bm{A}_s,
\]
and we define $\bB_{\{i,-i\}}$ similarly. The complex $\XI_n$ splits as
\[
\XI_n=\bigoplus_{\{i,-i\}\in \Z_n/\text{conj} } \bX\bI_{n;\{i,-i\}},
\]
where $\bX\bI_{n;\{i,-i\}}$ has underlying group $\bF[Q]/Q^2\otimes (\bA_{\{i,-i\}}\oplus \bB_{\{i,-i\}})$.  Furthermore,  
\[
\bX\bI_{n;\{i,-i\}}\simeq \bCFI^-(Y_n, \varpi_i).
\]

\subsection{Truncations}

Following \cite{OSIntegerSurgeries} and \cite{MOIntegerSurgery}, we define several truncation operations.

\begin{define}
Suppose $(C,\d)$ is a finitely generated, free complex over $\bF[U]$ or $\bF\llsquare U \rrsquare$ and $\delta\in \N$. The \emph{vertical trunction} of $(C,\d)$ is the complex $C^\delta:=C\otimes_{\bF[U]} \bF[U]/U^\delta$, with differential induced by $\d$.
\end{define}

The following algebraic lemma is a variation of a standard result (see \cite{OSIntegerSurgeries}*{Section~4.1}, \cite{MOIntegerSurgery}*{Section~10.2}):
\begin{lem} \label{lem:truncation-algebra} Suppose $\cR$ is $\bF[U]$ or $\bF\llsquare U \rrsquare$. Suppose that $(C_1,\d_1)$ and $(C_2,\d_2)$ are two chain complexes over $\cR$ or $\cR[Q]/Q^2$, which are homotopy equivalent to finitely generated, free chain complexes over $\cR$ or $\cR[Q]/Q^2$, which are relatively $\Z$-graded. Then $(C_1,\d_1)$ and $(C_2,\d_2)$ are homotopy equivalent if and only if  $(C_1^\delta,\d_1^{\delta_1})$ and $(C_2^{\delta},\d_2^{\delta_2})$ are homotopy equivalent for all $\delta$.
\end{lem}

We have also encountered the \emph{horizontal truncation} operation in~\eqref{eq:horizontal-truncation}, denoted $\bX_n\langle b\rangle$. We will use the following result:

\begin{lem}\label{lem:homotopy-equivalence-quasi}
If $\scC$ is an $\iota_K$-complex, equipped with a choice of flip-map, then $\bX_n(\scC)\langle b\rangle$ and $\bX_n(\scC)$ are homotopy equivalent over the ring $\bF\llsquare U\rrsquare$, for sufficiently large $b$.
\end{lem}

Ozsv\'{a}th and Szab\'{o} prove that $\bX_n(\scC)\langle b\rangle$ and $\bX_n(\scC)$ are quasi-isomorphic, by exhibiting a quasi-isomorphism from $\bX_n(\scC)$ to $\bX_n(\scC)\langle b\rangle$ if $n>0$ (projection), and a quasi-isomorphism from $\bX_n(\scC)\langle b\rangle $ to $\bX_n(\scC)$ if $n<0$ (inclusion) \cite{OSIntegerSurgeries}*{Lemma~4.3}. For bounded, free complexes, being quasi-isomorphic is equivalent to being homotopy equivalent (cf. \cite{Weibel}*{Lemma~10.4.6}). However, $\bX_n(\scC)$ is not free (it is an infinite direct product), and the set of gradings of the generators is not bounded. Nonetheless, it is not too hard to construct homotopy inverses by hand, using the fact that $v_s$ and (resp. $h_s$) admits a homotopy inverse if $s\gg 0$ (resp. $s\ll 0$). In slightly more detail, if $n>0$, one first constructs a map from $\bX_n(\scC)\langle b\rangle$ to $\bX_n(\scC)\langle b+1\rangle$, which is a homotopy inverse to projection, and which is furthermore a genuine one-sided inverse to projection (such a map may be chosen to have a similar shape to Figure~\ref{fig:local-map-F0-schematic}). By composing, we obtain a map $\bX_n(\scC)\langle b\rangle\to \bX_n(\scC)\langle b+k\rangle$, for any $k$. This construction naturally gives a limiting map, from $\bX_n(\scC)\langle b\rangle$ to $\bX_n(\scC)$, which is a homotopy inverse to projection. We leave the details to the reader.

\subsection{Filling hypercubes}

A key lemma for passing from the involutive hypercube to the involutive mapping cone formula is the following algebraic lemma:

\begin{lem}\label{lem:cone-filling-lemma} Suppose that we have a diagram of chain complexes 
\begin{equation}\label{eq:hcoherent-cube}
\begin{tikzcd}[column sep={1.5cm,between origins},row sep=.8cm,labels=description]
C^{000}
	\ar[dr, "D_{000}^{100}", pos=.6]
	\ar[ddd]
	\ar[dddrr,dashed]
	\ar[rr]
&[.7 cm]\,
&C^{010}
	\ar[rd]
	\ar[ddd]
	\ar[ddddr,dashed]
	&[.7 cm]\,
\\
&[.7 cm] C^{100}
	\ar[rr, crossing over]
&\,
&[.7 cm]C^{110}
	\ar[from=ulll, dashed, crossing over]
	\ar[ddd]
\\
\\
C^{001}
	\ar[dr]
	\ar[rr]
	\ar[drrr,dashed]	
	&[.7 cm]\,&
C^{011}
	\ar[dr]&[.9 cm]\,\\
& [.7 cm] C^{101}
	\ar[from=uuuul, dashed, crossing over]
	\ar[from=uuu,crossing over]\ar[rr]&\,
&[.9 cm] C^{111}
\end{tikzcd}
\end{equation}

\noindent consisting of the maps shown and the internal differentials of the complexes, with maps denoted as in Definition~\ref{def:hyperbox-chain-complexes}. Suppose that all faces, except for the front face, satisfy the 2-dimensional hypercube relation. Further, assume that
\[
D^{100}_{000}\colon C^{000}\to C^{100}
\]
is a chain homotopy equivalence. Then there exist maps
\[D^{011}_{100}: C^{100 } \rightarrow C^{111}\]
\[D^{111}_{000}: C^{000} \rightarrow C^{111}\]
\noindent which make $(C,D)$ into a three dimensional hypercube. In particular, if $G_{\mathrm{back}}$, $G_{\mathrm{bottom}}$, and $G_{\mathrm{top}}$ are the filtered maps
\[
G_{\mathrm{back}}\colon \Cone \left( C^{000}  \to
C^{010}\right)\to \Cone\left(C^{001} \to C^{011} \right)
\]
\[
G_{\mathrm{bottom}}\colon \Cone \left( C^{001}  \to
C^{011}\right)\to \Cone\left(C^{101} \to C^{111} \right)
\]
\[
G_{\mathrm{top}}\colon \Cone \left( C^{000}  \to
C^{010}\right)\to \Cone\left(C^{100} \to C^{110} \right)
\]
\noindent and $G_{\mathrm{front}}$ is the induced filtered map
\[
G_{\mathrm{front}}\colon \Cone \left( C^{100}  \to
C^{110}\right)\to \Cone\left(C^{101} \to C^{111} \right),
\]
\noindent then the compositions $G_{\mathrm{front}} \circ G_{\mathrm{top}}$ and $G_{\mathrm{bottom}} \circ G_{\mathrm{back}}$ are chain homotopy equivalent. \end{lem}

\begin{proof} Let $E: C^{100}\rightarrow C^{000}$ be a homotopy inverse to $D_{000}^{100}$. Choose homotopies 
\[
\begin{split}
H &\colon C^{000} \rightarrow C^{000}, \quad \text{and }\\
J&\colon  C^{100} \rightarrow C^{100}
\end{split}
\]
 such that $\partial H + H\partial = \mathrm{Id} + E \circ D_{000}^{100}$ and  $\partial J + J \partial = \mathrm{Id} + D_{000}^{100} \circ E$.
Observe that the sum $f=D^{100}_{000} \circ H + J \circ D_{000}^{100}$ is a chain map, since

\begin{align*}
\partial f + f\partial &= \partial(D^{100}_{000} \circ H + J \circ D_{000}^{100}) + (D^{100}_{000} \circ H + J \circ D_{000}^{100})\partial \\
					   &= D^{100}_{000}\circ (\partial H + H \partial) + (\partial J + J \partial)\circ D^{100}_{000} \\
					   &= D^{100}_{000}(\mathrm{Id} + E \circ D_{000}^{100}) + (\mathrm{Id} + D_{000}^{100} \circ E) \circ D^{100}_{000}\\
					   &= 0
\end{align*}

\noindent Let $K$ denote the sum of all of the length two compositions of maps from $C^{000}$ to $C^{111}$ appearing in the diagram (\ref{eq:hcoherent-cube}), so that

\[
 K = D_{011}^{100}  \circ D_{000}^{011} + D_{010}^{101}\circ D_{000}^{010} + D_{001}^{110}\circ D_{000}^{001} + D_{101}^{010}\circ D_{000}^{101} + D_{110}^{001} \circ D_{000}^{110}. \]

\noindent Observe that 

\[
\partial K + K \partial = \left(D^{010}_{101} \circ D^{001}_{100} + D^{001}_{110} \circ D^{010}_{100}\right) \circ D^{100}_{000}
\]

\noindent Using this notation, we claim that the following maps complete the hypercube
\begin{align*}
D^{011}_{100} &= K \circ E + \left(D^{010}_{101} \circ D^{001}_{100} + D^{001}_{110} \circ D^{010}_{100} \right)\circ J. \\ 
D^{111}_{000} &= K \circ H + D^{011}_{100}\circ f
\end{align*}
\noindent There are two relations to check: the 2-dimensional hypercube relation on the front face and the 3-dimensional hypercube relation for the entire cube. We see that
\begin{align*}
\partial D^{011}_{100} + D^{011}_{100} \partial &= (\partial K + K \partial) \circ E + \left(D^{010}_{101} \circ D^{001}_{100} + D^{001}_{110} \circ D^{010}_{100} \right)\circ (\partial J + J \partial) \\
&= \left(D^{010}_{101} \circ D^{001}_{100} + D^{001}_{110} \circ D^{010}_{100} \right)\circ D_{000}^{100}\circ E \\
& \qquad+ \left(D^{010}_{101} \circ D^{001}_{100} + D^{001}_{110} \circ D^{010}_{100} \right)\circ (\mathrm{Id} + D_{000}^{100} \circ E)\\
&=D^{010}_{101} \circ D^{001}_{100} + D^{001}_{110} \circ D^{010}_{100}
\end{align*}
\noindent which confirms the length 2 hypercube relation on the front face. We also have that
\begin{align*}
\partial D^{111}_{000} + D^{111}_{000} \partial &= (\partial K + K \partial)\circ H + K \circ (\partial H + H \partial) + \left(\partial D^{011}_{100} + D^{011}_{100} \partial\right) \circ f \\
						&= \left(\left(D^{010}_{101} \circ D^{001}_{100} + D^{001}_{110} \circ D^{010}_{100}\right) \circ D^{100}_{000}\right) \circ H \\
						&\qquad+ K \circ \left(\mathrm{Id} + E \circ D_{000}^{100}\right) + \left(D^{010}_{101} \circ D^{001}_{100} + D^{001}_{110} \circ D^{010}_{100}\right)\circ f\\
						&= K + K \circ E \circ D^{100}_{000} + \left(D^{010}_{101} \circ D^{001}_{100} + D^{001}_{110} \circ D^{010}_{100}\right) \circ (D^{100}_{000} \circ H +f )\\
						&= K + K \circ E \circ D^{100}_{000} + \left(D^{010}_{101} \circ D^{001}_{100} + D^{001}_{110} \circ D^{010}_{100}\right) \circ (J \circ D^{100}_{000})\\
						&= K + D^{011}_{100}\circ D^{100}_{000}
\end{align*}
\noindent This is exactly the length 3 hypercube relation.

\end{proof}

\begin{rem}\label{rem:cube-filling-lemma-iso}
If the map $D_{000}^{100}$ is a chain isomorphism, then the maps $D_{100}^{011}$  and $D_{000}^{111}$ have a particularly simple form. In this case, $D_{100}^{011}=K\circ (D_{000}^{100})^{-1}$ and $D_{000}^{111}\equiv 0$.
\end{rem}

\begin{rem}
	Lemma~\ref{lem:cone-filling-lemma} can also be viewed as a consequence of the fact that $\bF[U]$-chain complexes form an infinity category in the sense of Lurie \cite{LurieHTT}.  The defining property of an infinity category $\mathcal{C}$ is that diagrams on ``inner horns" $\Lambda^n_i$ for $0<i<n$ in $\mathcal{C}$ extend to diagrams on the standard simplex $\Delta^n$.  There is a precise sense in which cube (or other diagrammatic) extension problems are related to extension problems along simplices: this is the idea of anodyne extensions.  From this perspective, the diagram in Lemma \ref{lem:cone-filling-lemma} is not analogous to an inner horn, but to a horn $\Lambda^3_0\subset \Delta^3$.  However, an early result in the theory \cite{Joyal} shows that if the ``first arrow" in a horn $\Lambda^n_0$ is an equivalence, then the horn extends.  The proof of Lemma \ref{lem:cone-filling-lemma} can be viewed as writing down the map explicitly that Joyal constructs (inductively), slightly modified to treat cubes in place of simplices.
\end{rem}

\subsection{Gradings}

In this section, we describe how the maps in the involutive surgery hypercube decompose with respect to $\Spin^c$ structures, and their associated gradings. Consider the vertical truncation of the main hypercube:
\begin{figure}[H]
\begin{tikzcd}[column sep=2cm, row sep=2cm,labels=description]
\CF^\delta(\as_1,\bs)
	\ar[drr, dashed, "(h_{\a_1\to \a_2\to \a_3}^{\b})^{\delta}",sloped,]
	\ar[dr,"(f_{\a_1\to \a_2}^{\b})^{\delta}",sloped]
	\ar[dd, "\iota_1^\delta"]
	\ar[dddrr,dotted, sloped, swap,pos=.55,"(P_{\a_1\to \a_3}^{\b})^{\delta}"]
\\[-.7cm]
&[-.5cm] \CF^\delta(\as_2,\bs)
	\ar[r,"(f_{\a_2\to \a_3}^{\b})^{\delta}",swap]
	\ar[ddr,dashed, sloped, "(H_{\a_2\to \a_3}^{\b})^{\delta}"]
&[.8cm] \uCF^\delta(\as_3,\bs)
	\ar[dd, "\iota_3^\delta"]
\\[-.5cm]
\CF^\delta(\as_1,\bs)
	\ar[dr,"(f_{\a_1\to \a_2}^{\b})^{\delta}",swap,sloped]
	\ar[drr, dashed, "(h_{\a_1\to \a_2\to \a_3}^{\b})^{\delta}",sloped,pos=.6]
\\[-.7cm]
& \CF^\delta(\as_2,\bs)
	\ar[r,"(f_{\a_2\to \a_3}^{\b})^{\delta}",swap]
	\ar[from=uu, crossing over, "\iota_2^{\delta}"]
	\ar[from=uuul,crossing over, dashed, "(H_{\a_1\to \a_2}^{\b})^{\delta}",sloped]
&
\uCF^\delta(\as_3,\bs)
\end{tikzcd}
\caption{Vertically truncating the main hypercube.}
\label{fig:main-hypercube-relabeled}
\end{figure}
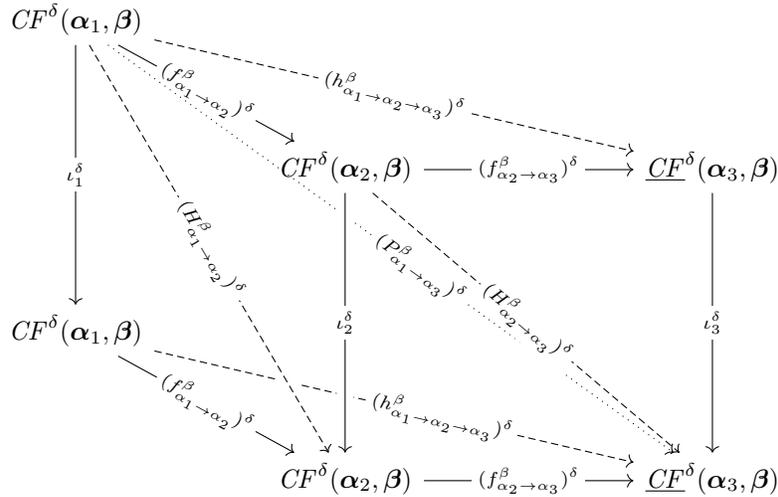

We give the complexes in Figure~\ref{fig:main-hypercube-relabeled} the  absolute grading induced by $\CF^-(Y)$, $\CF^-(Y_n)$ and $\CF^-(Y_{n+m})$. For $\uCF^{\delta}(\as_3,\bs)$, we give $T$ grading $0$. (Ozsv\'{a}th and Szab\'{o} refer to these as the \emph{old} gradings in \cite{OSIntegerSurgeries}).

\begin{figure}[H]
	\centering
\begingroup%
  \makeatletter%
  \providecommand\color[2][]{%
    \errmessage{(Inkscape) Color is used for the text in Inkscape, but the package 'color.sty' is not loaded}%
    \renewcommand\color[2][]{}%
  }%
  \providecommand\transparent[1]{%
    \errmessage{(Inkscape) Transparency is used (non-zero) for the text in Inkscape, but the package 'transparent.sty' is not loaded}%
    \renewcommand\transparent[1]{}%
  }%
  \providecommand\rotatebox[2]{#2}%
  \newcommand*\fsize{\dimexpr\f@size pt\relax}%
  \newcommand*\lineheight[1]{\fontsize{\fsize}{#1\fsize}\selectfont}%
  \ifx\svgwidth\undefined%
    \setlength{\unitlength}{122.01727454bp}%
    \ifx\svgscale\undefined%
      \relax%
    \else%
      \setlength{\unitlength}{\unitlength * \real{\svgscale}}%
    \fi%
  \else%
    \setlength{\unitlength}{\svgwidth}%
  \fi%
  \global\let\svgwidth\undefined%
  \global\let\svgscale\undefined%
  \makeatother%
  \begin{picture}(1,0.90688435)%
    \lineheight{1}%
    \setlength\tabcolsep{0pt}%
    \put(0,0){\includegraphics[width=\unitlength,page=1]{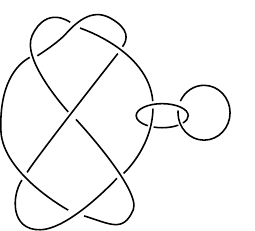}}%
    \put(0.92867301,0.43382859){\color[rgb]{0,0,0}\makebox(0,0)[lt]{\lineheight{1.25}\smash{\begin{tabular}[t]{l}$m$\end{tabular}}}}%
    \put(0.61668249,0.33261363){\color[rgb]{0,0,0}\makebox(0,0)[lt]{\lineheight{1.25}\smash{\begin{tabular}[t]{l}$0$\end{tabular}}}}%
    \put(0.44023213,0.88165739){\color[rgb]{0,0,0}\makebox(0,0)[lt]{\lineheight{1.25}\smash{\begin{tabular}[t]{l}$K$\end{tabular}}}}%
  \end{picture}%
\endgroup%

	\caption{The cobordism $W_m\colon Y_n(K)\# L(m,1)\to Y_{n+m}(K)$ is obtained by attaching a 2-handle along the 0-framed unknot.}\label{fig:193}
\end{figure}
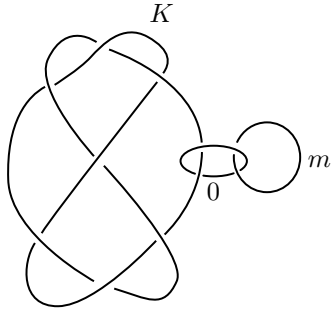

We write
\[
\begin{split}
W_n&\colon Y\to Y_n(K),\\
W_{n,m+n}& \colon Y_n(K)\# L(m,1)\to Y_{n+m}(K), \quad \text{and}\\
W_{n+m}&\colon Y\to Y_{n+m}(K)
\end{split}
\]
for the standard 2-handle cobordisms. We write $W_n'\colon Y_n(K)\to Y$ (and so forth) for the cobordism obtained by turning around and reversing the orientation of $W_n$. We note that $W_{n,m+n}$ is obtained by attaching a 0-framed 2-handle to $Y_n(K)\# L(m,1)$, which clasps $K$ and an $m$-framed unknot corresponding to the $L(m,1)$-summand. See Figure~\ref{fig:193}.

 A straightforward computation (see, e.g. \cite{MOIntegerSurgery}*{Lemma~11.1}) shows that $H_2(W_{n,m+n};\Z)$ is free, of rank 1, and is generated by a surface $\Sigma$ with
 \begin{equation}
\Sigma^2=-\frac{mn(m+n)}{\gcd(m,n)^2}. \label{eq:self-intersection-number}
 \end{equation}

Note that $W_n$ coincides with $X_{\a_3,\a_1,\b}$, after filling in $Y_{\a_3,\a_1}$ with 3-handles and a 4-handle. Similarly, $W_{n,n+m}$ is obtained from $X_{\a_2,\a_1,\b}$ by filling in the $S^1\times S^2$ summands of $Y_{\a_2,\a_1}$, and removing a neighborhood of an arc connecting $Y_{\a_1,\b}$ and $Y_{\a_2,\a_1}$. Finally, the 4-manifold $X_{\a_3,\a_2,\a_1,\b}$ may be identified with $W_n'\# D(-m,1)$, after filling in $S^1\times S^2$ summands.

The maps in the main hypercube are not homogeneously graded, but instead decompose as sums of homogeneously graded maps, indexed by $\Spin^c$ structures. The following lemma describes how the maps in the main hypercube decompose into homogeneously graded maps:

\begin{lem}\phantom{0} \label{lem:gradings}
\begin{enumerate}
\item[]
\item The maps $f_{\a_1\to \a_2}^{\b}$ and $H_{\a_1\to \a_2}^{\b}$ decompose as
\[
f_{\a_1\to \a_2}^{\b}=\sum_{\substack{\frs\in \Spin^c(W_{n,m+n})\\ \frs|_{L(m,1)}=\frs_0}} f_{\a_1\to \a_2;\frs}^{\b}\quad \text{and} \quad H_{\a_1\to \a_2}^{\b}=\sum_{\substack{\frs\in \Spin^c(W_{n,m+n})\\ \frs|_{L(m,1)}=\frs_0}} H_{\a_1\to \a_2;\frs}^{\b}.
\]
Furthermore, $f_{\a_1\to \a_2;\frs}^{\b}$ and $H_{\a_1\to \a_2;\frs}^{\b}$ are of homogeneous degrees
\[
\frac{c_1^2(\frs)+m}{4} \quad\text{and} \quad \frac{c_1^2(\frs)+m}{4}+1,
\]
respectively.
\item The maps $f_{\a_2\to \a_3}^{\b}$ and $H_{\a_2\to \a_3}^{\b}$ decompose as
\[
f_{\a_2\to \a_3}^{\b}=\sum_{\frs\in \Spin^c(W_{n+m}')} f_{\a_2\to \a_3;\frs}^{\b}\quad \text{and} \quad H_{\a_2\to \a_3}^{\b}=\sum_{\frs\in \Spin^c(W_{n+m}')} H_{\a_2\to \a_3;\frs}^{\b}.
\]
 Furthermore, $f_{\a_2\to \a_3;\frs}^{\b}$ and $H_{\a_2\to \a_3;\frs}^{\b}$ have homogeneous degree
\[
\frac{c_1^2(\frs)+1}{4}\quad \text{and}\quad \frac{c_1^2(\frs)+1}{4}+1,
\]
respectively.
\item The maps $h_{\a_1\to \a_2\to \a_3}^{\b}$ and $P_{\a_1\to \a_3}^{\b}$ decompose as
\[
\begin{split}
h_{\a_1\to \a_2\to \a_3}^{\b}&=\sum_{\substack{\frs \in \Spin^c(W_n'\# D(-m,1))\\ \frs|_{L(m,1)}=\frs_0} } h_{\a_1\to \a_2\to \a_3;\frs}^{\b} \quad \text{and} 
\\
P_{\a_1\to \a_2}^{\b} &=\sum_{\substack{\frs \in \Spin^c(W_n'\# D(-m,1))\\ \frs|_{L(m,1)}=\frs_0}} P_{\a_1\to\a_3;\frs}^{\b}.
\end{split}
\]
 Furthermore, $h_{\a_1\to \a_2\to \a_3;\frs}^{\b}$ and $P_{\a_1\to \a_3;\frs}^{\b}$ are of homogeneous degree
\[
\frac{c_1^2(\frs|_{W_{n,n+m}})+c_1^2(\frs|_{W_{n+m}'})+m+1}{4}+1 \quad \text{and} \quad \frac{c_1^2(\frs|_{W_{n,n+m}})+c_1^2(\frs|_{W_{n+m}'})+m+1}{4}+2,
\]
respectively.
\end{enumerate}
\end{lem}
\begin{proof}
For the maps along the top and bottom of the hypercube, the claim is immediate. The proof is slightly more subtle for the maps $H_{\a_1\to \a_2}^{\b}$, $H_{\a_2\to \a_3}^{\b}$ and $P_{\a_1\to \a_3}^{\b}$, as they are not defined via a single count of holomorphic polygons. These maps are defined via the compression operation for hyperboxes, and hence involve many different maps. 

We focus our attention first on the map $H_{\a_2\to \a_3}^{\b}$. This map is the diagonal map in the compression of a hyperbox of size $(1,d)$ for some large $d\ge 0$. The compression operation for such hyperboxes is illustrated in Figure~\ref{fig:80}. Each of the vertical maps represents either a change of diagram map, the cobordism map for a collection of 1-handles which are canceled in a later step by a collection of 2-handles, or the cobordism map for surgery on 0-framed unlink which is canceled by a collection of 3-handles in a later step. Each summand of $H_{\a_1\to \a_2}^{\b}$ contains exactly one length 2 map from a constituent hypercube. Each of these length 2 maps counts holomorphic triangles or rectangles on a triple or quadruple which topologically represents the 2-handle cobordism $W_{n+m}'\colon Y_{m+n}\to Y$, or a small modification of this cobordism. For example, the corresponding length 2 map in $\cC_{\cen}^{(3)}$ counts holomorphic triangles and rectangles on diagrams which represent a cobordism from $L(m,1)\# Y_{m+n}$ to $S^1\times S^2\# Y$, obtained by surgery on the dual of $K$ in $Y_{m+n}$, as well as a knot in $L(m,1)$. By direct inspection,  the maps decompose over $\Spin^c$ structures on $W_{m+n}'$ into a sum of maps which have the stated grading. The claim for $H_{\a_1\to \a_2}^{\b}$ is not substantially different.

The argument for the length 3 map $P_{\a_1\to \a_3}^{\b}$, is similar, but is slightly more involved. The compression operation for hyperboxes of size $(1,1,d)$ is illustrated in Figure~\ref{fig:79}. The summands of $P_{\a_1\to \a_3}^{\b}$ each consist either of a composition of numerous length 1 maps and exactly two length 2 maps, or numerous length 1 maps and exactly one length 3 map. For the summands of $P_{\a_1\to \a_3}^{\b}$ which contain two length 2 maps, the argument is not substantially different than the above argument for the length 2 maps $H_{\a_1\to \a_2}^{\b}$ and $H_{\a_2\to \a_3}^{\b}$. For the summands containing a length 3 map, the situation is similar, but in principle there are more checks to be done, as we now describe. Most of the length 3 maps which contribute to $P_{\a_1\to \a_3}^{\b}$ count pentagons on a Heegaard pentuple which represents surgery on the dual knot in $Y_n$, as well as a $-m$-framed unknot which is unlinked from $K$ (the associated 2-handle cobordism is $W_n'\# D(-m,1)$). There are also additional triangle and quadrilateral counts which arise from the higher compositions in a pairing of hypercubes of Lagrangians. For example, in the construction of $\cC_{\cen}^{(3)}$ in Section~\ref{sec:hypercube-Cc3},  the length 3 map has summands which count rectangles with one input labeled $\lambda$ and one input labeled $\Theta$, and also a summand corresponding to counting triangles with the special input $\omega_{\dt,\dt_3}$. There is a small subtlety in that the triangle map with special input $\omega_{\dt,\dt_3}$ counts triangles on a triple representing the cobordism $W'_n$. However we constructed the chain $\omega_{\dt,\dt_3}$ to decompose over $\Spin^c$ structures on $D(-m,1)$ (see~\eqref{eq:omega-d-d3-def}), and hence the corresponding triangle map decomposes over $\Spin^c(W'_n\# D(-m,1))$, as claimed. Finally, as described in Section~\ref{sec:constructing-C-cen}, the length 3 map of $\cC_{\cen}$ has an additional contribution of $H_E f_{\g}^{\dt\to \dt'}$. The map $f_{\g}^{\dt\to \dt'}$ is the 2-handle map from $Y_n$ to $Y$, while $H_E$ counts more complicated holomorphic curves on quadruples which themselves represent $([0,1]\times Y)\# D(-m,1)$ (followed by a cobordism for $0$-surgery on an unknot). Hence, it is straightforward to verify that $H_E$ decomposes over $\Spin^c(D(-m,1))$ into summands of the expected grading.
\end{proof}

Ozsv\'{a}th and Szab\'{o} additionally proved several strong restrictions on some of the vertically truncated maps appearing in Figure~\eqref{fig:main-hypercube-relabeled}, when $m$ is sufficiently large. To state their results, and how they adapt to our present situation,  we first recall some background and notation. There is a standard identification $\Spin^c(Y_{m+n}(K))\iso \Z/(m+n)$, as follows. Let $F$ be a Seifert surface in $Y$ for $K$, and let $\hat{F}\subset W_{m+n}'(K)$ be the surface which is capped off by the core of the 2-handle. If $i\in \Z/(m+n)$, we define $\frs_i\in \Spin^c(Y_{m+n}(K))$ to be the $\Spin^c$ structure which admits an extension $\frt_i$ over $W_{m+n}'(K)$ satisfying
\[
\langle c_1(\frt_i), [\hat{F}] \rangle\equiv n+m+2i\mod {2(n+m)}.
\]

We write $\frx_s$ and $\fry_s$ for the two $\Spin^c$ structures on $W_{n+m}'(K)$ which satisfy
\begin{equation}
\langle c_1(\frx_s),[\hat{F}]\rangle =2s-m-n\quad \text{and} \quad \langle c_1(\fry_s),[\hat{F}]\rangle=2s+m+n, \label{eq:characterize-x-s-y-s}
\end{equation}
respectively. If $-(m+n)/2\le  s< (m+n)/2$, then $\frx_s$ and $\fry_s$ are the $\Spin^c$ structures with maximal and second to maximal $c_1^2$. (See \cite{OSIntegerSurgeries}*{Lemma~4.4}).

\begin{prop}\label{prop:gradings+spinc-structures}
 Suppose $\delta$ is fixed.
 \begin{enumerate}
 \item\label{prop:gradings-1} If $m$ is sufficiently large, then for each $\frt\in \Spin^c(Y_{n+m})$, there is at most one $\Spin^c$ structure $\frs$ on $W_{n,m+n}$ which restricts to $\frt$ on $Y_{n+m}$, for which the map $(f_{\a_1\to \a_2;\frs}^{\b})^{\delta}$ is non-trivial. This is also the only $\Spin^c$ structure for which $(H_{\a_1\to \a_2;\frs}^{\b})^{\delta}$ may be non-trivial. 
 \item\label{prop:gradings-2} There is an integer $b>0$ with the following property. If $m$ is sufficiently large, and $-(m+n)/2\le s< (m+n)/2$, then the maps $(f_{\a_2\to \a_3;\frs})^{\delta}$ are non-vanishing on $\bCF^-(\as_2,\bs,[s])$ only if $\frs=\frx_s$ or $\frs=\fry_s$. Furthermore, if $s>b$, then $(f_{\a_2\to \a_3;\frs}^{\b})^{\delta}$ is non-vanishing only if $\frs=\frx_s$. If $s<-b$, then $(f_{\a_2\to \a_3;\frs}^{\b})^\delta$ is non-vanishing only if $\frs=\fry_s$. The same statements hold for the maps $(H_{\a_2\to \a_3;\frs}^{\b})^{\delta}$.
 \item \label{prop:gradings-3} If $m$ is sufficiently large, then $(h_{\a_1\to \a_2\to \a_3;\frs}^{\b})^{\delta}$ is non-vanishing only if $\frs|_{W_{n,n+m}}$ and $\frs|_{W_{n+m}'}$ are contained on the list of $\Spin^c$ structures described in parts ~\eqref{prop:gradings-1} and ~\eqref{prop:gradings-2}. The same holds for the maps $(P_{\a_1\to \a_3;\frs}^{\b})^{\delta}$.
 \end{enumerate} 
\end{prop}

Proposition~\ref{prop:gradings+spinc-structures} is proven by Ozsv\'{a}th and Szab\'{o} for the maps  $(f_{\a_1\to \a_2;\frs}^{\b})^{\delta}$, $(f_{\a_2\to \a_3;\frs}^{\b})^{\delta}$, and $(h_{\a_1\to \a_2\to \a_3;\frs}^{\b})^{\delta}$ when $n>0$. See specifically \cite{OSIntegerSurgeries}*{Lemma~4.4, Proposition~4.6, and Lemma 4.8}. Note that Ozsv\'{a}th and Szab\'{o} focus on the case that $m=nk$, for $k\in \N$, which simplifies equation~\eqref{eq:self-intersection-number}, though their proof extends to the more general case. They argued by showing that for sufficiently large $m$, the maps $(f_{\a_1\to \a_2;\frs}^{\b})^{\delta}$, $(f_{\a_2\to \a_3;\frs}^{\b})^{\delta}$ and $(h_{\a_1\to \a_2\to \a_3;\frs}^{\b})^\delta$ land outside of the range of gradings supported by the vertically truncated complexes appearing in Figure~\ref{fig:main-hypercube-relabeled}, except for the $\Spin^c$ structures described in the statement of Proposition~\ref{prop:gradings+spinc-structures}. Using the grading change formulas described in Lemma~\ref{lem:gradings}, it is easy to check that Ozsv\'{a}th and Szab\'{o}'s argument implies the same statement for the maps $H_{\a_1\to \a_2;\frs}^{\b}$, $H_{\a_2\to \a_3;\frs}^{\b}$ and $P_{\a_1\to \a_3;\frs}^{\b}$. 

For negative $n$, Ozsv\'{a}th and Szab\'{o}'s proof does not immediately apply (the main challenge being part~\eqref{prop:gradings-3}; note that for negative surgeries they instead built a map from $\Cone(D_n)$ to $\bCF^-(Y_n)$ whose grading properties are easier to analyze). We delay our proof for negative $n$ until later, where we cover it in greater generality in Section~\ref{sec:rational}. See specifically Lemmas~\ref{lem:len-1-homogeneously-graded-pt1}, ~\ref{lem:len-1-homogeneously-graded-pt2}, \ref{lem:len-1-homogeneously-graded-pt3} and \ref{lem:homogeneous-grading-lemma}.

\subsection{Proof of Theorem~\ref{thm:mapping-cone}}
\label{sec:truncated}

In this section, we prove the mapping cone formula, as stated in Theorem~\ref{thm:mapping-cone}.

\begin{proof}[Proof of Theorem~\ref{thm:mapping-cone}]

 Proposition~\ref{prop:gradings+spinc-structures} implies that if $\delta>0$ is fixed, and $m$ is sufficiently large, then the maps $(f_{\a_2\to \a_3;\frs}^{\b})^{\delta}$ and $(H_{\a_2\to \a_3;\frs}^{\b})^{\delta}$ will be non-zero only if $\frs$ is one of $\frx_s$ or $\fry_s$ for some $-(m+n)/2< s\le  (m+n)/2$.  We consider the hypercube in Figure~\ref{fig:relate-3-manifold-knot-Floer-step-1}, which relates the 3-manifold complexes to the knot Floer complexes.
\begin{figure}[H]
\begin{tikzcd}[column sep={3.5cm,between origins},
row sep=.3cm,labels=description]
\CF^\delta(\as_2,\bs)
	\ar[dd, swap,"\iota_2^\delta"]
	\ar[dr, "\Gamma^\delta"]
	\ar[rr, "(f_{\a_2\to \a_3}^{\b})^\delta"]
	\ar[ddrr,dashed, sloped, "(H_{\a_2\to \a_3}^{\b})^\delta", pos=.4]
&&[-1.5cm]
\uCF^\delta(\as_3,\bs)
	\ar[dd, "\iota_3^\delta"]
	\ar[dr, "\theta_w^\delta"]
&
\\
&\bA^\delta\langle b\rangle
	\ar[rr,crossing over, "D_n^\delta"]
&&
\bB^\delta\langle b \rangle
	\ar[dd, "\iota_{\bB}^\delta"]
	\ar[from=ulll, dashed,crossing over, "L^\delta"]
\\[2cm]
\CF^\delta(\as_2,\bs)
	\ar[rr, "(f_{\a_2\to \a_3}^{\b})^\delta", pos=.25]
	\ar[dr, "\Gamma^\delta"]
	\ar[drrr,dashed, "L^\delta"]
	&
&\uCF^\delta(\as_3,\bs) 
	\ar[dr, "\theta_w^\delta"]	
&
\\
&
\bA^{\delta}\langle b \rangle
	\ar[rr, "D_n^\delta"]
	\ar[from =uu, crossing over, "\iota_{\bA}^\delta"]
	\ar[from=uuul,dashed,crossing over,"J^\delta"]
	&&
\bB^{\delta}\langle b \rangle
	\ar[from=uull,dashed, crossing over, "H_m^\delta"]
\end{tikzcd}
\caption{A hypercube which relates the 3-manifold complexes with the knot Floer complexes. There is a length 3 map, which is not shown.}
\label{fig:relate-3-manifold-knot-Floer-step-1}
\end{figure}
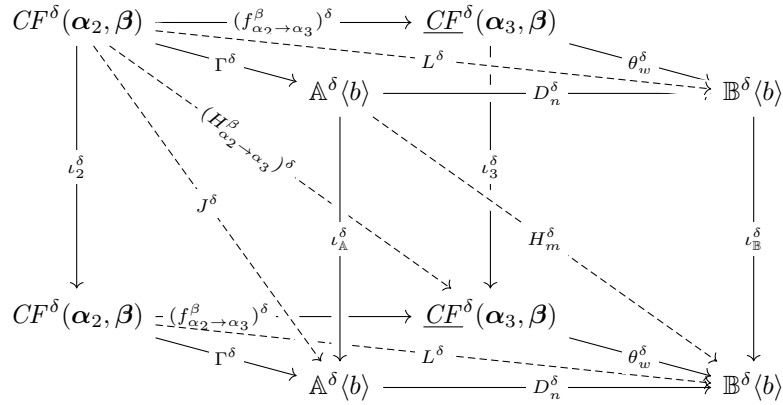

We now explain the terms appearing in Figure~\ref{fig:relate-3-manifold-knot-Floer-step-1}. The truncation $\bA\langle b\rangle$ is the direct sum of the complexes $A_s$, ranging over $-b\le s\le b$. Similarly $\bB\langle b\rangle$ is the direct sum of complexes $B_s$ ranging over $-b+n\le s\le b$. We are assuming here that 
\begin{equation}
0<b<(m+n)/2.
\label{eq:b<m+n/2}
\end{equation}

 The map $\Gamma$ is the direct sum of the maps
\[
\Gamma_s\colon \CF^-(\as_2,\bs,[s])\to A_s,
\]
described by Ozsv\'{a}th and Szab\'{o} \cite{OSKnots}*{Theorem~4.4} \cite{OSIntegerSurgeries}*{Theorem~2.3}, ranging over $-b\le s\le b$.  The map $\Gamma_s$ counts holomorphic triangles on $(\Sigma,\as_3,\as_2,\bs,w,z)$ which satisfy $\frs_{w}(\psi)=\frx_s$.

 Note that if $m+n$ is even and positive, then $Y_{n+m}(K)$ has two self-conjugate $\Spin^c$ structures, corresponding to $s=0$ and $s=(m+n)/2$. The condition in~\eqref{eq:b<m+n/2}
ensures that the one corresponding to $s=(m+n)/2$ is sent to zero by $ \Gamma^\delta$. Similarly, in the map $\theta_w^\delta$, we are implicitly projecting $\uCF^\delta(\as_3,\bs,[s])$ to zero unless $-b+n\le s\le b$.

We recall from \cite{HMInvolutive}*{Section~6.6} that there is a map $J_0$ making the following diagram a hypercube:
\begin{equation}
\begin{tikzcd}[row sep=1.5 cm, column sep =1.5 cm]
\CF^-(\as_2,\bs,[0])
	\ar[r, "\Gamma_{0}"]
	\ar[d, "\iota_2"]
	\ar[dr, dashed, "J_0"]	
&A_0
	\ar[d,"\iota_{\bA,0}"]\\
\CF^-(\as_2,\bs,[0])
	\ar[r,"\Gamma_0"]& 
A_0,
\end{tikzcd}
\end{equation}
where $\iota_{\bA,0}$ is the map induced by the knot involution $\iota_K$. Applying the same argument, but to other $\Spin^c$ structures, gives a map $J_s$ which makes the following diagram a hypercube of chain complexes:
\begin{equation}
\begin{tikzcd}[row sep=1.5 cm, column sep =1.5 cm]
\CF^-(\as_2,\bs,[s])
	\ar[r, "\Gamma_{s}"]
	\ar[d, "\iota_2"]
	\ar[dr,dashed, "J_s"]	
&A_s
	\ar[d,"\iota_{\bA,s}"]\\
\CF^-(\as_2,\bs,[-s])
	\ar[r,"\Gamma_{-s}"]& 
A_{-s},
\end{tikzcd}
\end{equation}
where $\iota_{\bA,s}\colon A_s\to A_{-s}$ denotes the map 
\[
x\mapsto U^{s}\cdot  \iota_K(x).
\]
We define $J$ to be the direct sum of the maps $J_s$, and define the map $J^\delta$ appearing along the left side of Figure~\ref{fig:relate-3-manifold-knot-Floer-step-1} to be the truncation, summed only over $s$ satisfying $-b\le s\le b$.

In Figure~\ref{fig:relate-3-manifold-knot-Floer-step-1}, we define $D_n$, $L$ and $H$ as follows. Let $I_s$ be a homotopy inverse of $\Gamma_s$. We defined $D_n$ on $A_s$ via the formula
\begin{equation}
D_{n,s}:=\theta_w \circ(f_{\a_2\to \a_3;\frx_s}^{\b}+f_{\a_2\to \a_3;\fry_s}^{\b})\circ I_s.
\label{eq:definition-D-n-cobordism-maps}
\end{equation}
To construct $L$, we pick a homogeneous null-homotopy $h$ of $\id+I\circ \Gamma$, and set
\[
L=\theta_w\circ f_{\a_2\to \a_3}^{\b}\circ h.
\]
Along the right face, we define $\iota_{\bB}$ to be $\theta_w\circ \iota_3\circ \theta_w^{-1}$, so that the right face commutes as well.

We define the map
\[
H_m\colon \bA_n\langle m\rangle \to \bB_n\langle m\rangle
\]
appearing along the front face of Figure~\ref{fig:relate-3-manifold-knot-Floer-step-1} to be similar to the one constructed in Lemma~\ref{lem:cone-filling-lemma}. The only caveat here is that the map labeled $\Gamma^\delta$ is not in general a homotopy equivalence, because of the projection maps. However, Lemma~\ref{lem:cone-filling-lemma} may be applied to a restricted version of the hypercube, shown in Figure~\ref{fig:relate-3-manifold-knot-Floer-step-2}.  Lemma~\ref{lem:cone-filling-lemma} also gives a length 3 map, which turns the diagram into a hypercube, though we have omitted this from the figure.

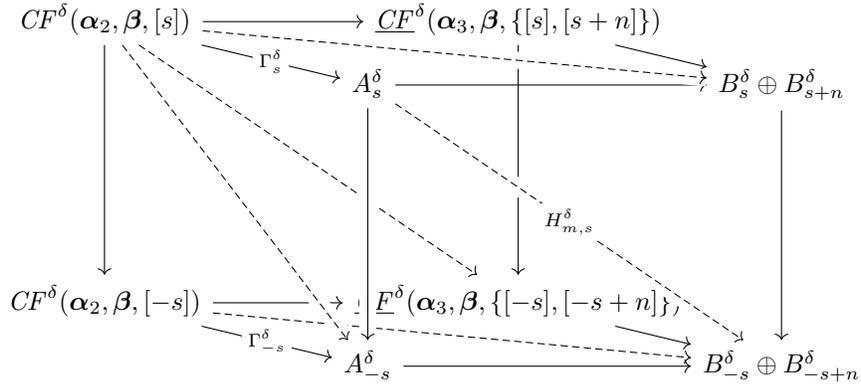
\begin{figure}[H]
\begin{tikzcd}[column sep={3.5cm,between origins},row sep=.2cm,labels=description]
\CF^\delta(\as_2,\bs,[s])
	\ar[dd]
	\ar[dr, "\Gamma_s^\delta"]
	\ar[rr]
	\ar[ddrr,dashed,crossing over]
&&[-1.5cm]
\uCF^\delta(\as_3,\bs,\{[s],[s+n]\})
	\ar[dd]
	\ar[dr]
&
\\
&A_s^\delta
	\ar[rr]
&&
B_s^\delta \oplus B_{s+n}^\delta
	\ar[dd]
	\ar[from=ulll,dashed,crossing over]
\\[2cm]
\CF^\delta(\as_2,\bs,[-s])
	\ar[rr]
	\ar[dr, "\Gamma_{-s}^\delta"]
	\ar[drrr,dashed]
	&
&\uCF^\delta(\as_3,\bs,\{[-s],[-s+n]\}) 
	\ar[dr]	
&
\\
&
A_{-s}^\delta
	\ar[rr,]
	\ar[from =uu, crossing over]
	\ar[from=uuul,dashed,crossing over]
	&&
B_{-s}^\delta\oplus B_{-s+n}^\delta
	\ar[from=uull,dashed, crossing over, "H_{m,s}^\delta"]
\end{tikzcd}
\caption{A subdiagram of Figure~\ref{fig:relate-3-manifold-knot-Floer-step-1}.}
\label{fig:relate-3-manifold-knot-Floer-step-2}
\end{figure}
We claim that the resulting diagram in Figure~\ref{fig:relate-3-manifold-knot-Floer-step-1} is a hypercube of chain complexes. This mostly follows immediately from the hypercube relations in Figure~\ref{fig:relate-3-manifold-knot-Floer-step-2}. For positive surgeries, the hypercube relations are actually satisfied on the full minus version (i.e. with no vertical truncations), as long as we restrict each map on the back face to the $\Spin^c$ structures of the form $\frx_s$ and $\fry_s$. For negative $n$, the hypercube relations will only hold on the vertically truncated complexes. The possibility that requires checking is if $s>b$ but $s+n\le b$. In this case, the hypercube relations follow from the fact that $f_{2,\fry_{s}}^{\delta}$ vanishes if $s$ is sufficiently large, relative to $\delta$, by Proposition~\ref{prop:gradings+spinc-structures}.

 Note that claims~\eqref{cone:1} and~\eqref{cone:4} hold for $D_n$ and $H_m$, since by construction we are only using the contributions from the $\Spin^c$ structures $\frx_s$ and $\fry_s$.

We stack and compress the hypercubes shown in Figures~\ref{fig:main-hypercube-relabeled}, and~\ref{fig:relate-3-manifold-knot-Floer-step-1}. The maps along the top and bottom faces of the resulting hypercube coincide. Hence, we may add the identity to each vertical map, and the hypercube relations are preserved. We obtain a hypercube of the form shown in  Figure~\ref{fig:truncated-but-sufficient-hypercube}. 
\begin{figure}[H]
\begin{tikzcd}[column sep=3cm, row sep=.7cm, labels=description]
\CF^\delta(\as_1,\bs)
	\ar[drr, dashed,sloped,]
	\ar[dr,sloped]
	\ar[dd, "\id+\iota_1^\delta"]
	\ar[dddrr,dotted, sloped, swap]
&[-2cm]
&[.3cm]
\\[-.2cm]
& \bA^\delta_n\langle b \rangle
	\ar[r,"D_n^\delta",swap]
	\ar[ddr,dashed,  "H^\delta_m"]
& \bB^\delta_n\langle b\rangle
	\ar[dd, "\id+\iota_{\bB}^\delta"]
\\[1.3cm]
\CF^\delta(\as_1,\bs)
	\ar[dr,swap,sloped]
	\ar[drr, dashed,sloped, near end]
\\[-.2cm]
& \bA^{\delta}_n\langle b\rangle
	\ar[r,"D_n^\delta",swap]
	\ar[from=uu, crossing over, "\id+\iota_{\bA}^\delta"]
	\ar[from=uuul,crossing over, dashed]
&
\bB^{\delta}_n\langle b\rangle 
\end{tikzcd}
\caption{A hypercube relating $\CFI^\delta(Y_n)$ and a truncated mapping cone.}
\label{fig:truncated-but-sufficient-hypercube}
\end{figure}
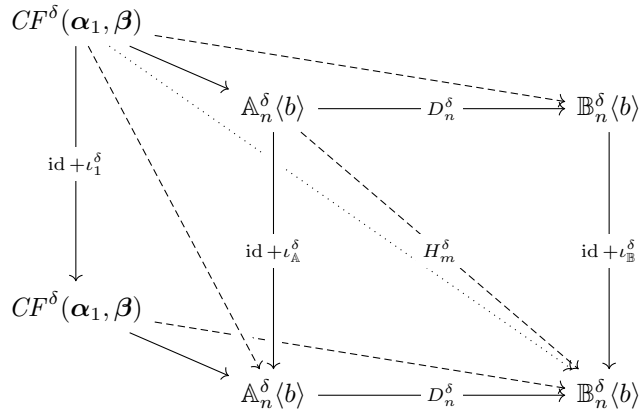

The diagram in Figure~\ref{fig:truncated-but-sufficient-hypercube} determines a $\bF[U,Q]/Q^2$ equivariant map from $\CFI^\delta(Y_n)$ to the chain complex associated to the 2-dimensional hypercube on the front side of Figure~\ref{fig:truncated-but-sufficient-hypercube}. The maps along the top and bottom faces coincide with Ozsv\'{a}th and Szab\'{o}'s homotopy equivalence between $\CF^\delta(\as_1,\bs)$ and $\Cone(D_n^\delta \colon \bA_n^\delta\langle b\rangle \to \bA_n^\delta\langle b\rangle).$

We now claim that we can obtain from the diagram in Figure~\ref{fig:truncated-but-sufficient-hypercube} both finite and infinite models of the involutive mapping cone, with maps which satisfy the conditions outlined in Theorem~\ref{thm:mapping-cone}.

The required $m$, necessary for a hypercube as in Figure~\ref{fig:truncated-but-sufficient-hypercube} to exist, depends on $\delta$. We let $\delta_i$ be a sequence approaching $+\infty$, and we $m_i$ be a sequence of integers, such that $m_i$ is sufficiently large for the hypercube in Figure~\ref{fig:truncated-but-sufficient-hypercube} to exist. For convenience, take $b=\lfloor \frac{m+n-1}{2}\rfloor$, which ensures~\eqref{eq:b<m+n/2} is satisfied.

 Since we are working over $\bF=\Z/2$, the set of elements of $A_s$ and $B_s$ in each grading is finite. Consequently, by picking a subsequence, we may assume that $H_{m_i}(x)$ is eventually constant, for each $x\in A_s$. Consequently, we obtain a limiting map, which we call $\tilde{H}$, which we extend $\bF\llsquare U\rrsquare$-equivariantly. The map $\tilde{H}$ clearly makes the untruncated analog of the front of the hypercube in Figure~\ref{fig:truncated-but-sufficient-hypercube} into a 2-dimensional hypercube of chain complexes. Analogously, we may assume that the maps labeled $D_n$ are eventually constant on each $x\in A_s$, for each $s$. Hence, we obtain a limiting map, for which we write $\tilde{D}_n$.

 Let us write $\tilde{\bXI}_n$ for the infinite involutive mapping cone, constructed with the maps $\tilde{D}_n$ and $\tilde{H}$. There are also several horizontal truncations to consider, namely the one appearing in Figure~\ref{fig:truncated-but-sufficient-hypercube}, constructed for $m_i$, for which we write $\bXI_{n}(m_i)$ (it is easy to check that the hypercube relations are satisfied without vertically truncating). We may also horizontally truncate $\tilde{\bXI}_n$. We write $\tilde{\bXI}_n\langle m\rangle$ for this complex.
 
 We make the following subclaims:
 \begin{enumerate}[ref=$t$-\arabic*, label=($t$-\arabic*),leftmargin=*]
 \item\label{claim-truncations-1} If $m$ is large, then $\tilde{\bXI}_{n} $ and $\tilde{\XI}_{n}\langle m\rangle$ are homotopy equivalent over $\bF\llsquare U\rrsquare [Q]/Q^2$.
 \item\label{claim-truncations-2} If $m_i$ is large, then $\tilde{\XI}_n\langle m_i\rangle$ and $\XI_{n}(m_i)$ are homotopy equivalent over $\bF\llsquare U\rrsquare[Q]/Q^2$.
 \end{enumerate}

 The proof of claim~\eqref{claim-truncations-1} is as follows. Let us assume $n>0$, for the sake of demonstration. In this case, projection from $\tilde{\XI}_{n}$ to $\tilde{\XI}_n \langle m\rangle$ is an $\bF\llsquare U\rrsquare [Q]/Q^2$-equivariant chain map. Furthermore, the map obtained by setting $Q=0$ is a homotopy equivalence, by Lemma~\ref{lem:homotopy-equivalence-quasi}. Hence, Lemma~\ref{lem:homotopy-equivalence-iota-complexes} implies that the involutive version is an $\bF\llsquare U\rrsquare [Q]/Q^2$-equivariant homotopy equivalence. The proof for negative coefficients is essentially the same, except one instead considers the inclusion of $\tilde{\XI}_n\langle m\rangle$ into $\tilde{\XI}_n$, which is a homotopy equivalence for large $m$.
  
The proof of claim~\eqref{claim-truncations-2} is similar. Fix an integer $b$.  For sufficiently large $m_i$, the complex $\XI_{n}(m_i)\langle b\rangle$ is equal to $\tilde{\XI}_{n}\langle b\rangle$ (since $\tilde{\XI}_n$ is obtained as a limit). As before, let us focus on the case that $n$ is positive. Hence if $b$ is sufficiently large, and $m_i$ is chosen sufficiently large, relative to $b$, then the projection for $\tilde{\XI}_n$ onto $\tilde{\XI}_n\langle b\rangle$ is an $\bF\llsquare U\rrsquare[Q]/Q^2$-equivariant homotopy equivalence,  and projection of $\XI_{n}(m_i)$ onto $\tilde{\XI}_n\langle b\rangle$ is also an $\bF\llsquare U\rrsquare[Q]/Q^2$-equivariant homotopy equivalence. Claim~\eqref{claim-truncations-2} now follows.

The diagram in Figure~\ref{fig:truncated-but-sufficient-hypercube} induces an $\bF[U,Q]/Q^2$-equivariant chain map between $\CFI^{\delta_i}(Y_n(K))$ and $\bXI^{\delta_i}_{n}(m_i)$, for each $i$, which we see is an $\bF[U,Q]/Q^2$-homotopy equivalence, by combining the fact that the $Q=0$ specializing is an $\bF[U]$-equivariant homotopy equivalence with Lemma~\ref{lem:homotopy-equivalence-iota-complexes}. By claims~\eqref{claim-truncations-1} and ~\eqref{claim-truncations-2}, we obtain also that $\CFI^\delta(Y_n(K))$ and $\tilde{\XI}^\delta_n$ are homotopy equivalent over $\bF[U,Q]/Q^2$, for all $\delta$. Lemma~\ref{lem:truncation-algebra} now implies the claim for the complexes without the vertical truncations.
 \end{proof}

\subsection{Proof of Theorems~\ref{thm:exact-sequence-1} and \ref{thm:exact-sequence-2}}

We now prove the exact sequences described in the introduction:
\begin{proof}[Proof of Theorem~\ref{thm:exact-sequence-2}]

  The front face of the involutive surgery hypercube of Theorem~\ref{thm:main-hypercube} determines an $\bF\llsquare U\rrsquare [Q]/Q^2$-equivariant chain map 
  \[
  F_2\colon \bCFI^-(Y_{n+m}(K))\to \underline{\bCFI}^-(Y).
  \]
  Furthermore, the remaining maps in the hypercube determine an $\bF\llsquare U\rrsquare[Q]/Q^2$-equivariant chain map
  \[
  F_1\colon \bCFI^-(Y_n(K))\to \Cone(F_2).
  \] 
  Furthermore, $F_1$ is a homotopy equivalence of $\bF\llsquare U\rrsquare [Q]/Q^2$-modules by Lemma~\ref{lem:homotopy-equivalence-iota-complexes}, since Ozsv\'{a}th and Szab\'{o} proved that the map obtained by setting $Q=0$ is a homotopy equivalence of $\bF\llsquare U\rrsquare$-modules \cite{OSIntegerSurgeries}*{Theorem~3.1}. The long exact sequence of a mapping cone now gives the exact sequence
\[
\cdots \to \bHFI^-(Y_n(K))\to \bHFI^-(Y_{n+m}(K))\to \underline{\bHFI}^-(Y)\to \cdots.
\]
  
It remains only to prove the description of $\underline{\bHFI}^-(Y)$ stated in~\eqref{eq:decomposition-twisted-HFI-intro}. The map $\theta_w$ gives a chain isomorphism $\bCFI^-(Y)\iso \bCF^-(Y)\otimes \bF[\Z/m].$ The map $\underline{\iota}_3$ sends $\bCF^-(Y)\otimes T^s$ to $\bCF^-(Y)\otimes T^{-s+n}$, but otherwise coincides with the involution on $\bCF^-(Y)$ by Proposition~\ref{prop:doubling-transition-formula} and Lemma~\ref{lem:columns=change-of-diagrams}. It is also important to check the path of the basepoint $w$ in the corresponding loop of Heegaard diagrams. We observe that our model of $\iota_K$ moves $w$ to $\bar z$ negatively along $K$, while the flip map moves $z$ to $w$ in the positive direction along $K$, so the overall path of $w$ is null-homotopic.

 Hence, 
\[
\underline{\bHFI}^-(Y)\iso \bigoplus_{i=1}^{N} \bHFI^-(Y)\oplus \bigoplus_{i=1}^M (\bHF^-(Y)[-1]\oplus \bHF^-(Y)),
\]
where $N$ is the number of subsets $\{s,-s+n\}\subset \Z/m$ which consist of one element, and $M$ is the number of subsets $\{s,-s+n\}$ which consist of two elements. Note that $\{s,-s+n\}$ consists of one element if and only if $2s\equiv n\mod m$. If $m$ is odd, then the map $\times 2\colon \Z/m\to \Z/m$ is a bijection. If $m$ is even, then $\times 2$ has image exactly equal to the even numbers, and the fiber of each even number has cardinality 2. Hence, if $m$ is odd, then $N=1$ and $M=(m-1)/2$. If $m$ is even, and $n$ is odd, then $N=0$ and $M=m/2$. If $m$ and $n$ are even, then $N=2$ and $M=(m-2)/2$. 
\end{proof}

\begin{proof}[Proof of Theorem~\ref{thm:exact-sequence-1}]
 The proof is very similar to the proof of Theorem~\ref{thm:exact-sequence-2}. We build a hypercube similar to the hypercube in Theorem~\ref{thm:main-hypercube}, but which has $m=1$. The construction is essentially the same as Figure~\ref{fig:large-hyperbox}, except we only keep track of the basepoint labeled with a hat, and we use no twisted coefficients. With this hypercube in hand, the proof follows no differently than the proof of Theorem~\ref{thm:exact-sequence-2}.
\end{proof}

\section{The strong cone formula for L-spaces}
\label{sec:strong-cone-formula}

In this section, we prove the strong version of the mapping cone formula, for knots in L-spaces:

\begin{thm}\label{thm:strong-cone-formula-proven}
Suppose $n\neq 0$ and $K$ is a knot in an L-space integer homology 3-sphere. Then $\bCFI^-(Y_n(K))$ is homotopy equivalent over $\bF\llsquare U\rrsquare [Q]/Q^2$ to a complex of the form
\[
\bXI_n(Y,K):= \begin{tikzcd}[column sep=2cm, row sep=2cm, labels=description]
\bA
	\ar[d, "Q\cdot (\id +\iota_{\bA})",swap]
	\ar[dr, dashed, "Q\cdot H\vopp"]
	\ar[r, "D_n"]
& \bB
	\ar[d, "Q\cdot (\id+\iota_{\bB})"]\\
Q\cdot \bA
	\ar[r, "D_n"]
& Q\cdot \bB,
\end{tikzcd}
\]
for some $\bF\llsquare U\rrsquare$-equivariant maps $D_n,$ $\iota_{\bA}$, $\iota_{\bB}$ and $H$. Furthermore, the following hold:
\begin{enumerate}
\item $H$ sends $\Bopp_s$ to $B_{-s}$.
\item $D_n=v+\frF \vopp$, where $v$ denotes the inclusion of $A_s$ into $B_s$, and $\vopp$ denotes inclusion of $A_s$ into $\Bopp_s$.
\item $\iota_{\bB}=\frF \iota_K$.
\item $\iota_{\bA}$ coincides with the conjugation of $\iota_K$ by the shift isomorphism in~\eqref{eq:shift-map-def}.
\end{enumerate}
In particular, $(\bCF^-(Y_n(K)),\iota)$ is $\iota$-homotopy equivalent to $\XI_n^{\alg}(\CFK^\infty(Y,K),\iota_K)$.
\end{thm}

Interestingly, our argument does not seem to work unless $Y$ is an L-space. See Remark~\ref{rem:natural=challenge}.

\subsection{Overview of the proof of Theorem~\ref{thm:strong-cone-formula-proven}}

We begin by outlining the proof of Theorem~\ref{thm:strong-cone-formula-proven}. We describe more details in the subsequent sections. We begin by considering a vertical truncation of the involutive surgery hypercube involving $Y_n(K)$, $Y_{n+m}(K)$ and $Y$, as shown in Figure~\ref{fig:main-hypercube-relabeled}. It follows from Lemma~\ref{lem:gradings} that for sufficiently large $m$, the only $\Spin^c$ structures which contribute to the front face are $\frx_s$ or $\fry_s$, ranging over $-(m+n)/2\le s\le (m+n)/2$. Here $\frx_s$ and $\fry_s$ are characterized by~\eqref{eq:characterize-x-s-y-s}.

\begin{figure}[h]
\begin{tikzcd}[column sep={2.5cm,between origins},
row sep=.6cm,labels=description
]
\CF^-(Y_{n+m},[s])
	\ar[dd, swap,"\iota_2"]
	\ar[dr,  "\Gamma_s"]
	\ar[rr, " f_{2;\frx_s}"]
	\ar[ddrr,dashed,near end,"H_{2;\frx_s}"]
&&[1.2cm]
\uCF^-(Y,[s])
	\ar[dd, "\iota_3"]
	\ar[dr,"\theta_w"]
&
\\
&A_s(Y,K)
	\ar[rr,line width=2mm,dash,color=white,opacity=.7]
&&
B_s(Y,K)
	\ar[from=ll,crossing over, "v_s"]
	\ar[dd, "\iota_{\bB}"]
\\[1.8cm]
\CF^-(Y_{n+m},[-s])
	\ar[rr, "f_{2;\fry_{-s}}", pos=.6]
	\ar[dr,"\Gamma_{-s}"]
	\ar[drrr,dashed, "L_{-s}"]
&&\uCF^-(Y,[-s+n])
	\ar[dr, " \theta_w"]	
&
\\
&
A_{-s}(Y,K)
	\ar[rr, "\frF_{z\to w} \vopp_{-s}"]
	\ar[from =uu, "\iota_{\bA}",crossing over]
	\ar[from=uuul,dashed,"J_s",crossing over]
	&&
B_{-s+n}(Y,K)
\end{tikzcd}
\begin{tikzcd}[column sep={2.6cm,between origins},
row sep=.6cm,labels=description,
]
\CF^-(Y_{n+m},[s])
	\ar[dd, swap,"\iota_2"]
	\ar[dr,  "\Gamma_s"]
	\ar[rr, " f_{2;\fry_s}"]
	\ar[ddrr,dashed,near end, "H_{2;\fry_s}"]
&&[.7cm]
\uCF^-(Y,[s+n])
	\ar[dd, "\iota_3"]
	\ar[dr,"\theta_w"]
&
\\
&A_s(Y,K)
&&
B_{s+n}(Y,K)
	\ar[dd, "\iota_{\bB}"]
	\ar[from=ll, "\frF_{z\to w}\vopp_s"]
	\ar[from=ulll,dashed, "L_{s}",crossing over]
\\[1.8cm]
\CF^-(Y_{n+m},[-s])
	\ar[rr, "f_{2;\frx_{-s}}", pos=.6]
	\ar[dr,"\Gamma_{-s}"]
&&\uCF^-(Y,[-s])
	\ar[dr, " \theta_w"]	
&
\\
&
A_{-s}(Y,K)
	\ar[rr, "v_{-s}"]
	\ar[from =uu,"\iota_{\bA}",crossing over]
	\ar[from=uuul,dashed, "J_s",crossing over]
	&&
B_{-s}(Y,K)
	\ar[from=uull,line width=2mm,dash,color=white,opacity=.7]
	\ar[from=uull, dashed, near start,"H_s\vopp_s"]
\end{tikzcd}
\caption{The hypercubes $\cK_{\frx_s}$ (top) and $\cK_{\fry_s}$ (bottom). The length 3 arrows are not drawn. The maps labeled $L_s$ coincide in the two diagrams, and similarly for the maps $J_s$ and $\Gamma_s$. Both cubes have a length 3 map which is not shown.}
\label{eq:cube-frx-s-and-fry-s}
\end{figure}
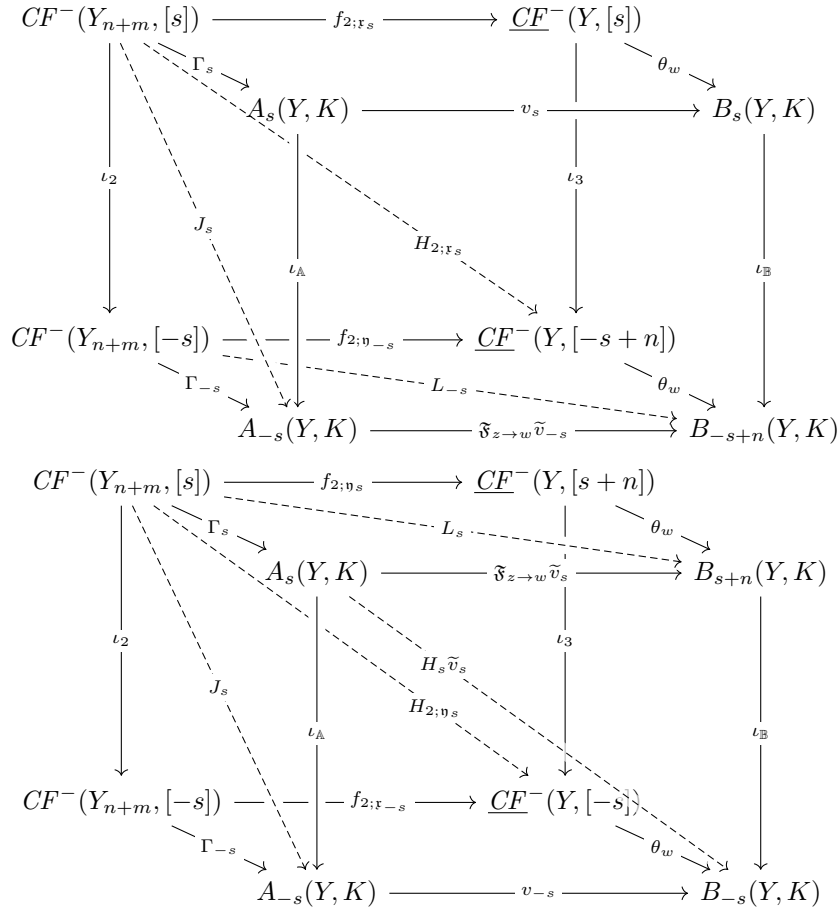

The key technical step of our proof of Theorem~\ref{thm:strong-cone-formula-proven} is to construct two 3-dimensional hypercubes, $\cK_{\frx_s}$ and $\cK_{\fry_s}$, for each $s\in \Z$. After restricting $\Spin^c$ structures and vertically truncating, the back faces of the hypercubes $\cK_{\frx_s}$ and $\cK_{\fry_s}$ match the front face of the main involutive surgery hypercube in Figure~\ref{fig:main-hypercube-relabeled}.  The hypercubes $\cK_{\frx_s}$ and $\cK_{\fry_s}$ are shown schematically in Figure~\ref{eq:cube-frx-s-and-fry-s}. The maps labeled $L_s$ coincide in the two diagrams, and similarly for $J_s$ and $\Gamma_s$. Therein, the indexing of the maps is consistent with notation in Theorem~\ref{thm:main-hypercube} of the introduction, so $f_{2;\frx_s}$ denotes $f_{\a_2\to \a_3;\frx_s}^{\b}$, and similarly for the other maps.

We now prove Theorem~\ref{thm:strong-cone-formula-proven}, assuming the existence of the hypercubes $\cK_{\frx_s}$ and $\cK_{\fry_s}$:
\begin{proof}[Proof of Theorem~\ref{thm:strong-cone-formula-proven}]

We construct hypercubes $\cK_{\frX}\langle b\rangle $ by taking the direct sum of the hypercubes $\cK_{\frx_s}$ ranging over $-b\le s\le b$. We define $\cK_{\frY}\langle b\rangle$ similarly by taking the direct sum of the $\cK_{\fry_s}$ for $-b\le s\le b$. We now combine $\cK_{\frX}\langle b\rangle$ and $\cK_{\frY}\langle b\rangle$ to construct a new cube, $\cK_{\{\frX,\frY\}}\langle b\rangle$, the vertical truncation of which is shown in  Figure~\ref{fig:C-x-y-full}. Here, $b$ is an integer satisfying
\[
0\ll b<(m+n)/2.
\]
The maps along the left face of $\cK_{\{\frX,\frY\}}\langle b\rangle$ coincide with the ones in $\cK_{\frX}\langle b\rangle$ (which themselves coincide with those of $\cK_{\frY}\langle b\rangle$). The maps along the top, front, bottom and back faces of $\cK_{\{\frX,\frY\}}^\delta\langle b\rangle$, except for the purely vertical maps (i.e. the involutions), are defined to be the sum of the maps from $\cK_{\frX}\langle b\rangle$ and $\cK_{\frY}\langle b\rangle$.  The map along the right face of $\cK_{\{\frX,\frY\}}^\delta\langle b\rangle$ is zero. The construction of $\cK_{\{\frX,\frY\}}\langle b\rangle$ is essentially the same as the operation of \emph{gluing} hypercubes, described by Manolescu and Ozsv\'{a}th \cite{MOIntegerSurgery}*{Proposition~12.10}.

By a slight abuse of notation, we  define $\Gamma$, $L$ and $J$ in Figure~\ref{fig:C-x-y-full} to project to 0 any summands corresponding to  $[s]\in \Spin^c(Y_{n+m})$ for $s$ satisfying $|s|>b$. Similarly, in the map labeled $\theta_w$ we are also implicitly post-composing with the map which projects $B_s$ to $0$ unless $-b+n\le s\le b$.

\begin{figure}[h]
\begin{tikzcd}[column sep={3.5cm,between origins},row sep=.4cm,labels=description]
\CF^\delta(Y_{n+m})
	\ar[dd, "\iota_2^\delta"]
	\ar[dr, "\Gamma^\delta"]
	\ar[rr,"f_2^\delta"]
	\ar[ddrr,dashed, "H_2^\delta"]
&&[-1.5cm]
\uCF^\delta(Y)
	\ar[dd, "\iota_3^\delta"]
	\ar[dr, "\theta_w^\delta"]
&
\\
&\bA^\delta\langle b\rangle
	\ar[rr, "D_n^\delta"]
&&
\bB^\delta\langle b\rangle
	\ar[dd, "\iota_{\bB}^\delta"]
	\ar[from=ulll,dashed,crossing over, "L^\delta"]
	\ar[from=ll,crossing over, "D_n^\delta"]
\\[2cm]
\CF^\delta(Y_{n+m})
	\ar[rr, "f_2^\delta", pos=.37]
	\ar[dr, "\Gamma^\delta"]
	\ar[drrr,dashed, "L^\delta"]
	&
&\uCF^\delta(Y) 
	\ar[dr, "\theta_w^\delta"]	
&
\\
&
\bA^\delta \langle b\rangle
	\ar[rr, "D_n^\delta"]
	\ar[from =uu, crossing over, "\iota_{\bA}^\delta"]
	\ar[from=uuul,dashed,crossing over, "J^\delta"]
	&&
\bB^\delta \langle b\rangle
	\ar[from=uull,dashed, crossing over, "H^\delta\vopp^\delta"]
\end{tikzcd}
\caption{The hypercube $\cK_{\{\frX,\frY\}}^\delta\langle b\rangle$. There is a length 3 map which is not shown.}
\label{fig:C-x-y-full}
\end{figure}

For $s$ satisfying $-b\le s\le b$, the map $\Gamma_s$ coincides with the large surgery map, and hence is a homotopy equivalence. Hence, the same argument as at the end of the proof of Theorem~\ref{thm:mapping-cone} (with $\cK_{\{\frX,\frY\}}^\delta\langle b\rangle$ replacing the cube in Figure~\ref{fig:relate-3-manifold-knot-Floer-step-2}) shows that $(\bCF^-(Y_n(K)),\iota)$ is $\iota$-homotopy equivalent to $\XI^{\alg}_n(\CFK^\infty(K),\iota_K)$, completing the proof.
\end{proof}

We now outline the construction of the hypercubes $\cK_{\frx_s}$ and $\cK_{\fry_s}$, which will be described in more detail in the subsequent sections. The main involutive surgery hypercube of Theorem~\ref{thm:main-hypercube} is constructed as the compression of the hyperbox of size $(1,1,7)$ shown in Figure~\ref{fig:large-hyperbox}. The hyperbox therein is constructed by stacking the seven hypercubes
\begin{equation}
\cC_{\oneh},\quad
\cC_{\twoh}^1, \quad
\cC_{\cen}, \quad
\cC_{\twoh}^2,\quad
\cC_{\thrh},\quad
\cC_{\eta} \quad \text{and}\quad 
\cC_{\frF}.
\label{eq:C-hypercubes-overview}
\end{equation}

The hypercube $\cK_{\frx_s}$ is built as the compression of a hyperbox of size $(1,1,7)$,  composed of seven 3-dimensional hypercubes
\begin{equation}
\cK_{\oneh;\frx_s},\quad
\cK_{\twoh;\frx_s}^1, \quad
\cK_{\cen;\frx_s}, \quad
\cK_{\twoh;\frx_s}^2,\quad
\cK_{\thrh;\frx_s},\quad
\cK_{\eta;\frx_s} \quad \text{and}\quad 
\cK_{\frF;\frx_s}.
\label{eq:K-x-s-hypercubes-overview}
\end{equation}
The back face of each of the hypercubes in~\eqref{eq:K-x-s-hypercubes-overview} matches a $\Spin^c$ restriction of the front face of the corresponding hypercube in~\eqref{eq:C-hypercubes-overview}. The hypercube $\cK_{\fry_s}$ is constructed similarly as the compression of a hyperbox of size $(1,1,7)$, consisting of seven hypercubes
\begin{equation}
\cK_{\oneh;\fry_s},\quad
\cK_{\twoh;\fry_s}^1, \quad
\cK_{\cen;\fry_s}, \quad
\cK_{\twoh;\fry_s}^2,\quad
\cK_{\thrh;\fry_s},\quad
\cK_{\eta;\fry_s} \quad \text{and}\quad 
\cK_{\frF;\fry_s}.
\label{eq:K-y-s-hypercubes-overview}
\end{equation}

We now introduce some additional notation which will be helpful in this section. Recall that the description of twisted Floer complexes from Section~\ref{sec:twisted-complexes-def} involved a choice of orientation $\fro$ for $K$. In this section, we fix an orientation $\fro$ on $K$, which never changes in the course of the argument. We will need to frequently change the special basepoint. To simplify the notation, if $\cH=(\Sigma,\as,\bs,w,z)$ is a doubly pointed diagram for $(Y,K)$, we write $\cH_w$ for $\cH$, equipped with the enhanced orientation $(\fro,w)$. Similarly, we write $\cH_z$ for $\cH$ equipped with the enhanced orientation $(\fro,z)$.

\subsection{The hypercube $\cK_{\frx_s}$}

 We now describe the construction of $\cK_{\frx_s}$ in detail. This amounts to constructing the seven hypercubes shown in ~\eqref{eq:K-x-s-hypercubes-overview}.
 
  First, suppose that $\cC_\circ$ is one of the first five hypercubes in~\eqref{eq:C-hypercubes-overview}. We will describe the construction of the corresponding hypercube in~\eqref{eq:K-x-s-hypercubes-overview}, which we will denote by $\cK_{\circ;\frx_s}$. Let $\cH^1_0$, $\cH^2_0$, and $\cH^3_0$ denote the diagrams appearing along the top face of $\cC_\circ$, and let $\cH^1_1$, $\cH^2_1$ and $\cH^3_1$ denote the diagrams appearing along the bottom face. We will abuse notation slightly and  write $w$ and $z$ for the basepoints on $\cH_0^3$ and $\cH_1^3$, even though they trace a half-twist along $K$ in the doubling model of the involution; see Sections~\ref{sec:doubled-knots} and~\ref{sec:involution-doubling}.
  
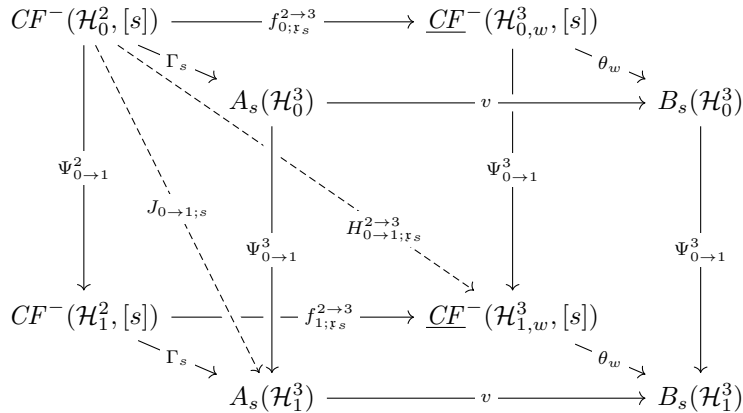
\begin{figure}[h]
\begin{tikzcd}[column sep={2.5cm,between origins},row sep=.4cm,labels=description]
\CF^-(\cH^2_0,[s])
	\ar[dd, swap,"\Psi_{0\to 1}^{2}"]
	\ar[dr,  "\Gamma_s"]
	\ar[rr, " f^{2\to 3}_{0;\frx_s}"]
	\ar[ddrr,dashed,near end, "H^{2\to 3}_{0\to 1;\frx_s}"]
&&[.7cm]
\uCF^-(\cH_{0,w}^3,[s])
	\ar[dd, "\Psi^{3}_{0\to 1}"]
	\ar[dr,"\theta_w"]
&
\\
&A_s(\cH^3_{0})
	\ar[rr,crossing over,  "v"]
&&
B_{s}(\cH_{0}^3)
	\ar[dd, "\Psi^{3}_{0\to 1}"]
\\[1.8cm]
\CF^-(\cH^2_{1},[s])
	\ar[rr, "f^{2\to 3}_{1;\frx_s}", pos=.65]
	\ar[dr,"\Gamma_{s}"]
&&\uCF^-(\cH^3_{1,w},[s])
	\ar[dr, " \theta_w"]	
&
\\
&
A_{s}(\cH^3_{1})
	\ar[rr, "v"]
	\ar[from =uu, crossing over,"\Psi^{3}_{0\to 1}"]
	\ar[from=uuul,dashed, crossing over, "J_{0\to 1;s}"]
	&&
B_{s}(\cH^3_{1})
\end{tikzcd}
\caption{A schematic of the hypercube $\cK_{\circ,\frx_s}$, which is one of the first five hypercubes of~\eqref{eq:K-x-s-hypercubes-overview}. There is no length 3 map.}
\label{eq:C-frx-s-i->i+1}
\end{figure}

   The cube $\cK_{\circ;\frx_s}$ is illustrated schematically in Figure~\ref{eq:C-frx-s-i->i+1}. The cube $\cK_{\circ;\frx_s}$ has no length 3 map. The maps $\Psi_{0\to 1}^2$ and $\Psi_{0\to 1}^3$ appearing along the back face are the same ones which appear along the front face of $\cC_{\circ}$. The maps labeled $\Psi_{0\to 1}^3$ appearing along the front face are defined using the same sequence of Heegaard diagrams as for the map $\Psi_{0\to 1}^3$ along the back face, as well as the same holomorphic curve counts, but instead using the weights $\knotU^{n_{w}(\psi)} \knotV^{n_{z}(\psi)}$. Since the map constructed in this manner preserves the Alexander grading (by explicit examination of the construction of each hypercube), we conclude that there is an induced map, also denoted $\Psi_{0\to 1}^3$, which sends $A_s(\cH_0^3)$ to $A_s(\cH_1^3)$, and similarly for the $B_s$ complexes.
   
  The map $\Gamma_s$ is the composition of the the shift isomorphism between $\cA_s$ and $A_s$, together with the map from $\CF^-(\cH_0^2,[s])$ to $\cA_s(\cH_0^2)$, which  counts holomorphic triangles representing classes which satisfy $\frs_{w}(\psi)=\frx_s$, weighted by a factor of $\knotU^{n_w(\psi)}\knotV^{n_z(\psi)}$. This map counts the same triangles as $f_{0;\frx_s}^{2\to 3}$.  Similarly, we define $J_{0\to 1;s}$ by composing the shift map from $\cA_s$ to $A_s$ with the analogs on $\cA_s$ of each constituent map of $H_{0\to 1;\frx_s}^{2\to 3}$. For example, if  $H_{0\to 1;\frx_s}^{2\to 3}$ had a factor or summand which counted holomorphic rectangles representing classes satisfying $\frs_{w}(\psi)=\frx_s$, then the map $J_{0\to 1;s}$ would also have a constituent map which counted rectangles weighted by $\knotU^{n_w(\psi)} \knotV^{n_z(\psi)}$. 
 
  To prove that $\Gamma_s$ and $J_{0\to 1;s}$ have codomain in $A_s$, as claimed, one needs to check that the maps $\Gamma_s$ and $J_{0\to 1;s}$ change Alexander grading by $s$ (viewing $\CF^-(\cH_0^2,[s])$ and $\CF^-(\cH_1^2,[s])$ as being concentrated in Alexander grading 0). This follows from \cite{ZemAbsoluteGradings}*{Theorem~1.4} (compare \cite{OSRationalSurgeries}*{Proposition~2.2}). 
  
  Next, we verify the hypercube relations for $\cK_{\circ;\frx_s}$. The length 2 relations along the left, back, right and front faces are straightforward. The length 2
relation along the top and bottom faces, as well as the length 3 relation, follow from the fact that $\Gamma_s$ and $J_{0\to 1;s}$ are both given as a map which changes Alexander grading by $s$, followed by the shift isomorphism, as well as the fact that they both  count the same holomorphic curves as $f_{0;\frx_s}^{2\to 3}$ or $f_{1;\frx_s}^{2\to 3}$, and $H_{0\to 1;\frx_s}^{2\to 3}$, respectively.

We now describe the hypercube $\cK_{\eta;\frx_s}$, which is shown in Figure~\ref{eq:C-frx-eta}. Therein, we have added some decoration to the map $\Gamma_s$ to disambiguate terms.  We write $\Gamma_{\frs;w}$ for the map which counts holomorphic triangles satisfying $\frs_w(\psi)=\frs$,  which are weighted by a factor of $\knotU^{n_w(\psi)} \knotV^{n_z(\psi)}$, composed with the shift isomorphism between $\cA_s$ and $A_s$. The map $\Gamma_s$ in Figure~\ref{eq:C-frx-s-i->i+1} coincides with $\Gamma_{\frx_s;w}$. The hypercube relations therein are easily checked. 
 
\begin{figure}[h]
\begin{tikzcd}[column sep={2.5cm,between origins},row sep=.4cm,labels=description]
\CF^-(\bar{\cH}^2,[s])
	\ar[dd, swap,"\eta_2"]
	\ar[dr,  "\Gamma_{\frx_s;z}"]
	\ar[rr, " f_{\frx_s;z}^{2\to 3}",pos=.6]
&&[.7cm]
\uCF^-(\bar{\cH}^3_{z},[s])
	\ar[dd, "\ul{\eta}_3"]
	\ar[dr,"\theta_z"]
&
\\
&A_s(\bar{\cH}^3)
	\ar[rr,crossing over,  "v"]
&&
B_{s}(\bar{\cH}^3)
	\ar[dd, "U^s \eta_K"]
\\[1.8cm]
\CF^-(\cH^2,[-s])
	\ar[rr, "f^{2\to 3}_{\fry_{-s};z}", pos=.6]
	\ar[dr,"\Gamma_{\bar{\frx}_s;z}"]
&&\uCF^-(\cH^3_{z},[-s])
	\ar[dr, " \theta_z"]	
&
\\
&
A_{-s}(\cH^3)
	\ar[rr, "\vopp"]
	\ar[from =uu, crossing over,"U^s\eta_K"]
	&&
\Bopp_{-s}(\cH^3)
\end{tikzcd}
\caption{The hypercube $\cK_{\eta;\frx_s}$. There are no length 2 or 3 maps.}
\label{eq:C-frx-eta}
\end{figure}

We now describe a helpful conjugation relation for the maps labeled by $\Gamma$. We claim that for all $s\in \Z$, 
\begin{equation}
\Gamma_{\frx_s;w}=\Gamma_{\fry_{s};z}\qquad \text{and} \qquad \bar{\frx}_s=\fry_{-s},\label{eq:Gamma-conjugation}
\end{equation}
The second part of~\eqref{eq:Gamma-conjugation} is immediate from the definition of $\frx_s$ and $\fry_s$ in~\eqref{eq:characterize-x-s-y-s}. The first part follows from the relation
\begin{equation}
\frs_w(\psi)-\frs_z(\psi)=\PD[\widehat{S}],\label{eq:strong-cone-change-spinc}
\end{equation}
where $\widehat{S}$ is the class in $H_2(W_{n+m}'(K))$ obtained by capping the core of the 2-handle. Equation~\eqref{eq:strong-cone-change-spinc} follows from ~\cite{ZemAbsoluteGradings}*{Lemma~3.9}. Hence $\Gamma_{\frx_s;w}$ and $\Gamma_{\fry_s;z}$ count the same triangles, with the same powers of $\knotU$ and $\knotV$. (Note that the choice of basepoint in the map $\Gamma_{\fry_s;z}$ affects only which classes are counted, via the choice of $\Spin^c$ structure, but not the $\knotU$ or $\knotV$ weight). Equation~\eqref{eq:Gamma-conjugation} follows.

\begin{figure}[h]
\begin{tikzcd}[
column sep={2.5cm,between origins},
row sep=.5cm,labels=description,
fill opacity=.7, text opacity=1,
execute at end picture={
\foreach \Nombre in  {A,B,C,D}
  {\coordinate (\Nombre) at (\Nombre.center);}
\fill[opacity=0.1] 
  (A) -- (B) -- (C) -- (D) -- cycle;
}]
\CF^-(\cH^2, [-s])
	\ar[dd, "\id"]
	\ar[dr,  "\id"]
	\ar[rr, "f_{\fry_{-s};z}^{2\to 3}"]
	\ar[ddrr,dashed, pos=.7, "H^{2\to 3}_{z\to w}"]
&&
\uCF^-(\cH^3_{z},[-s])
	\ar[dd, "\frF_{z\to w}"]
	\ar[dr,"\id"]
	\ar[dddr,dashed, "H^{\frF}_{z\to w\to z}"]
	\ar[rr, "\id"]
&&
\uCF^-(\cH^3_{z},[-s])
	\ar[dd, "\frF_{z\to w}"]
	\ar[dr,"\id"]
&\,
\\
&
|[alias=A]|\CF^-(\cH^2,[-s])
	\ar[rr,crossing over, "f^{2\to 3}_{\fry_{-s};z}"]
&&
\uCF^-(\cH_{z}^3,[-s])
	\ar[rr, "\id"]
&&
|[alias=B]|\uCF^-(\cH^3_{z},[-s])
	\ar[dd, "\frF_{z\to w}"]
\\[1.8cm]
\CF^-(\cH^2, [-s]) 
	\ar[rr, "f_{\fry_{-s};w}^{2\to 3}",pos=.25]
	\ar[dr,"\id"]
	\ar[drrr,dashed,pos=.6, "H^{2\to 3}_{w\to z}"]
&&
\uCF^-(\cH^3_{w},[-s+n])
	\ar[dr, "\frF_{w\to z}"]	
	\ar[rr, "\id",pos=.7]
	\ar[drrr,dashed, "H^{\frF}_{w\to z\to w}"]
&&
\uCF^-(\cH^3_{w},[-s+n])
	\ar[dr, "\id"]
\\
&
|[alias=D]|\CF^-(\cH^2,[-s])
	\ar[rr, "f_{\fry_{-s};z}^{2\to 3}"]
	\ar[from =uu, crossing over,swap, "\id"]
&&
\uCF^-(\cH^3_{z},[-s])
	\ar[from=uu,crossing over, "\id"]
	\ar[rr, "\frF_{z\to w} "]
&& |[alias=C]|\uCF^-(\cH^3_{w},[-s+n])
\end{tikzcd}
\begin{tikzcd}
[column sep={2.3cm,between origins},row sep=.4cm,labels=description,
fill opacity=.7, text opacity=1,
execute at end picture={
\foreach \Nombre in  {A,B,C,D}
  {\coordinate (\Nombre) at (\Nombre.center);}
\fill[opacity=0.1] 
  (A) -- (B) -- (C) -- (D) -- cycle;
}]
|[alias=A]|\CF^-(\cH^2,[-s])
	\ar[dd, swap,"\id"]
	\ar[dr,  "\Gamma_{\frx_{-s};w}"]
	\ar[rr, " f^{2\to 3}_{\fry_{-s};z}"]
&&[.7cm]
\uCF^-(\cH^3_{z},[-s])
	\ar[dd, "\id"]
	\ar[dr,"\theta_z"]
	\ar[rr, "\id"]
&&
|[alias=B]|\uCF^-(\cH^3_{z},[-s])
	\ar[dd, "\frF_{z\to w}"]
	\ar[dr,"\theta_z"]
&\,
\\
&A_{-s}(\cH^{3})
	\ar[rr,crossing over,  "\vopp"]
&&
\Bopp_{-s}(\cH^{3})
	\ar[dd, "\id"]
	\ar[rr,crossing over,  "\id"]
&&
\Bopp_{-s}(\cH^3)
	\ar[dd, "\frF_{z\to w}"]
\\[1.8cm]
|[alias=D]|\CF^-(\cH_2,[s])
	\ar[rr, "f^{2\to 3}_{\fry_{-s};z}", near end]
	\ar[dr,"\Gamma_{\frx_{-s};w}"]
&&\uCF^-(\cH^3_{z},[-s])
	\ar[dr, " \theta_z"]	
	\ar[rr, "\frF_{z\to w}"]
&&
|[alias=C]|\uCF^-(\cH^3_{w},[-s+n])
	\ar[dr, "\theta_w"]
&
\\
&
A_{-s}(\cH^3)
	\ar[rr, "\vopp"]
	\ar[from =uu, crossing over,"\id"]
	&&
\Bopp_{-s}(\cH^3)
	\ar[rr,"\frF_{z\to w}"]
&&
B_{-s+n}(\cH^3)
\end{tikzcd}
\caption{The hyperbox whose compression is $\cK_{\frF;\frx_s}$. The top two cubes each have a length 3 map, which is not shown, while the bottom two have no length 3 map. We stack the top and bottom hyperboxes along the gray faces.  We think of this hyperbox as being made of four hypercubes, $\cK_{\frF;\frx_s}^{1}, \dots \cK_{\frF;\frx_s}^{4}$, which are the top-left, top-right, bottom-left, and bottom-right hypercubes, respectively.}
\label{eq:C-frx-s-flip-full}
\end{figure}

\begin{figure}[h]
\begin{tikzcd}
[column sep={1.5cm,between origins},row sep=.8cm,labels=description,fill opacity=.7, text opacity=1]
(\as_2,\bs)
	\ar[dr,"\id"]
	\ar[dd,"\id"]
	\ar[rrrrr, "f_{\a_2\to \a_3}^{\b}"]
&&&&[-1cm]&
(\as_3,\bs,z)
	\ar[dr,"\id"]
	\ar[dd, "\phi_*"]
&&&
\\
&
(\as_2,\bs)
	\ar[dr,"\id"]
	\ar[dd,"\id"]
&&&&&
(\as_3,\bs,z)
	\ar[dd, "\id"]
	\ar[dr,"\id"]
\\
(\as_2,\bs)
	\ar[dr,"\id"]
	\ar[dd,"\id"]
&&
(\as_2,\bs)
	\ar[dd,"\id"]
	\ar[dr,"\id"]
&&&
(\phi(\as_3),\bs,w)
	\ar[dr, "\phi^{-1}_*"]
	\ar[dd]
	\ar[dddr,dashed]
&&
(\as_3,\bs,z)
	\ar[dr,"\id"]
	\ar[dd,"\id"]
\\
&(\as_2,\bs)
	\ar[dr,"\id"]
	\ar[dd,"\id"]
&&
(\as_2,\bs)
	\ar[dd,"\id"]
&&&
(\as_3,\bs,z)
	\ar[dr,"\id"]
	\ar[dd,"\phi_*"]
&&
(\as_3,\bs,z)
	\ar[dd,"\id"]
\\
(\as_2,\bs)
	\ar[dr,"\id"]
	\ar[dd,"\id"]
&&
(\as_2,\bs)
	\ar[dd,"\id"]
	\ar[dr,"\id"]
&&&
(\as_3^{\wind},\bs,w)
	\ar[dr]
	\ar[dd]
	\ar[dddr,dashed]
&&
(\as_3,\bs,z)
	\ar[dr,"\id"]
	\ar[dd,"\id"]
\\
&(\as_2,\bs)
	\ar[dr,"\id"]
	\ar[dd,"\id"]
&&
(\as_2,\bs)
	\ar[dd,"\id"]
&&&
(\phi(\as_3),\bs,w)
	\ar[dd]
	\ar[dr,"\phi^{-1}_*"]
	\ar[dddr,dashed]
&&(\as_3,\bs,z)
	\ar[dd,"\id"]
\\
(\as_2,\bs)
	\ar[dr,"\id"]
&&(\as_2,\bs)
	\ar[dd,"\id"]
	\ar[dr,"\id"]
&&&
(\as_3,\bs,w)
	\ar[dr]
&&
(\as_3,\bs,z)
	\ar[dr, "\id"]
	\ar[dd,"\phi"]
\\
&(\as_2,\bs)
	\ar[dr,"\id"]
&&
(\as_2,\bs)
	\ar[dd,"\id"]
&&&
(\as_3^{\wind},\bs,w)
	\ar[dr]
&&
(\as_3,\bs,z)
	\ar[dd,"\id"]
\\
&&(\as_2,\bs)
	\ar[dr,"\id"]
&&&&&
(\phi(\as_3),\bs,w)
	\ar[dr,"\phi^{-1}"]
\\
&&&
(\as_2,\bs)
	\ar[rrrrr, "f_{\a_2\to \a_3}^{\b}"]
&&&&&
(\as_3,\bs,z)
\end{tikzcd}
\caption{The hyperbox whose compression is $\cK^{1}_{\frF;\frx_s}$. The solid arrows count holomorphic triangles, while the dashed arrows count holomorphic quadrilaterals and holomorphic triangles.  Many arrows from the left-side to the right-side are not drawn.}
\label{eq:C-frx-s-flip-1}
\end{figure}
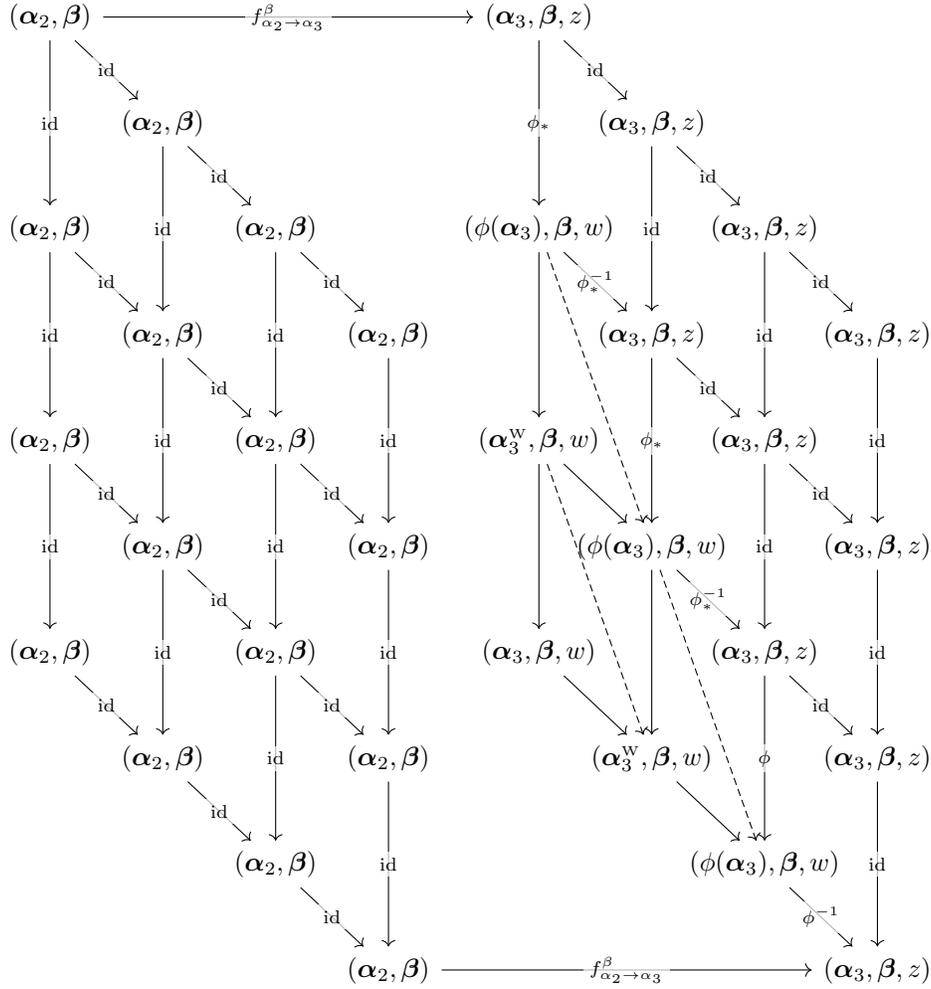

We now construct $\cK_{\frF;\frx_s}$ by stacking and compressing four hypercubes, $\cK^{1}_{\frF;\frx_s},\dots, \cK^{4}_{\frF;\frx_s}$, as shown in Figure~\ref{eq:C-frx-s-flip-full}. We begin by describing the hypercube $\cK^{1}_{\frF;\frx_s}$, which is the top left hypercube in Figure~\ref{eq:C-frx-s-flip-full}. Let $\phi\colon \Sigma\to \Sigma$ denote the diffeomorphism which moves $z$ to $w$ along the canonical short path in the Heegaard diagram. We construct $\cK^{1}_{\frF;\frx_s}$ as the compression of a hyperbox of size $(1,3,3)$. This hyperbox is shown schematically in Figure~\ref{eq:C-frx-s-flip-1}. Therein, we assume that $\as_3^{\twind}$ is a sufficiently wound copy of $\as_3$, which we additionally assume is wound so as not to intersect the boundary of the surgery region. We view this hyperbox as consisting of nine hypercubes. Each of these hypercubes may be constructed by pairing a hyperbox of alpha attaching circles with the 0-dimensional hypercube consisting of $\bs$. We also compose some of the maps appearing in these hypercubes with one of the tautological maps $\phi_*$ or $\phi^{-1}_*$. When a single attaching curve appears at different locations in the hyperbox, one would need to replace certain copies with small Hamiltonian translates. Furthermore, we are implicitly using Proposition~\ref{prop:nearest-point-triangles} to identity the maps labeled by $\id$ with the corresponding small triangle maps.

We now construct the hypercube $\cK^{2}_{\frF;\frx_s}$, which is the top right hypercube shown in Figure~\ref{eq:C-frx-s-flip-full}. The length 2 map along the left side of $\cK^{2}_{\frF;\frx_s}$ is constructed in $\cK^{1}_{\frF;\frx_s}$. The map along the bottom is constructed similarly. A length 3 map may be defined by noting that every complex appearing at a vertex of $\cK^{2}_{\frF;\frx_s}$ is homotopy equivalent to $\bF\llsquare U\rrsquare$, since $Y$ is an L-space. The expression 
\[
M:=\frF_{z\to w} H^{\frF}_{z\to w\to z}+H^{\frF}_{w\to z\to w}\frF_{z\to w}
\]
is an $\bF[U]$-equivariant chain map, which is $+1$ graded and hence is chain homotopic to zero. We let the length 3 map in $\cK^{2}_{\frF;\frx_s}$ be any null-homotopy of $M$ which is homogeneously graded. (Alternatively, one could define $\cK^{2}_{\frF;\frx_s}$ by building a hyperbox of size $(3,3,3)$, similar to the one used to build $\cK^{1}_{\frF;\frx_s}$, but with two more layers).

Finally, we consider $\cK_{\frF;\frx_s}^{3}$ and $\cK_{\frF;\frx_s}^{4}$, the bottom two hypercubes in Figure~\ref{eq:C-frx-s-flip-full}. We note that the hypercube relations for these cubes are automatic, since all faces commute on the nose, so no length 2 or 3 arrows are necessary.

\subsection{The hypercubes \texorpdfstring{$\cK_{\oneh;\fry_s}$, $\cK_{\twoh;\fry_s}^1$, $\cK_{\twoh;\fry_s}^2$ and $\cK_{\thrh;\fry_s}$}{K-1h-ys, K-2h-1-ys, K-2h-2-ys and K-3h-ys}}
\label{sec:more-flips1}

\begin{figure}[p]
\begin{tikzcd}[column sep={2.5cm,between origins},
row sep=.55cm,labels=description,
fill opacity=.7, text opacity=1
]
\CF^-(\cH^2_0,[s])
	\ar[dd, swap,"\Psi_{0\to 1}^{2}"]
	\ar[dr,  "\Gamma_s"]
	\ar[rr, " f^{2\to 3}_{0;\fry_s;w}"]
	\ar[ddrr,dashed,near end, "H^{2\to 3}_{0\to 1;\fry_s;w}",sloped]
&&[.7cm]
\uCF^-(\cH^3_{0,w},[s+n])
	\ar[dd, "\Psi^{3}_{0\to 1}"]
	\ar[dr,"\theta_w"]
&
\\
&A_s(\cH^3_{0})
	\ar[rr,line width=2mm,dash,color=white,opacity=.7]
	\ar[rr,near start,  "\frF_{z\to w}\vopp"]
&&
B_{s+n}(\cH^3_{0})
	\ar[dd, "\Psi^{3}_{0\to 1}"]
	\ar[from=ulll,line width=2mm,dash,color=white,opacity=.7]
	\ar[from=ulll,dashed, "H^{2\to 3}_{0}"]
\\[1.8cm]
\CF^-(\cH^2_{1},[s])
	\ar[rr, "f^{2\to 3}_{1;\fry_s;w}", near end]
	\ar[dr,"\Gamma_{s}"]
	\ar[drrr,dashed, "H^{2\to 3}_{1}"]
&&\uCF^-(\cH^3_{1,w},[s+n])
	\ar[dr, " \theta_w"]	
&
\\
&
A_{s}(\cH^3_{1})
	\ar[rr, "\frF_{z\to w}\vopp"]
	\ar[from=uu,line width=2mm,dash,color=white,opacity=.7]
	\ar[from =uu, "\Psi^{3}_{0\to 1}"]
	\ar[from=uuul,line width=2mm,dash,color=white,opacity=.7]
	\ar[from=uuul,dashed, "J_{0\to 1;s}"]
	&&
B_{s+n}(\cH^3_{1})
	\ar[from=uull,line width=2mm,dash,color=white,opacity=.7]
	\ar[from=uull, dashed, near start," H^{3}_{0\to 1}\vopp", sloped]
\end{tikzcd}
\begin{tikzcd}[
column sep={2.5cm,between origins},
row sep=.55cm,labels=description,fill opacity=.7, text opacity=1,
execute at end picture={
\foreach \Nombre in  {A,B,C,D}
  {\coordinate (\Nombre) at (\Nombre.center);}
\fill[opacity=0.1] 
  (A) -- (B) -- (C) -- (D) -- cycle;
}
]
\CF^-(\cH^2_0, [s])
	\ar[dd, "\Psi_{0\to 1}^2"]
	\ar[dr,  "\id"]
	\ar[rr, "f^{2\to 3}_{0;\fry_{s};w}"]
	\ar[ddrr,dashed, pos=.7, "H^{2\to 3}_{0\to 1;\fry_s;w}",sloped]
&&
\uCF^-(\cH^3_{0,w},[s+n])
	\ar[dd, "\Psi_{0\to 1}^3"]
	\ar[dr,"\frF_{w\to z}"]
	\ar[dddr,dashed, "H_{0\to 1;w\to z}^3",sloped]
	\ar[rr, "\id"]
&&
\uCF^-(\cH^3_{0,w},[s+n])
	\ar[dd, "\Psi_{0\to 1}^3"]
	\ar[dr,"\id"]
&\,
\\
&
|[alias=A]|\CF^-(\cH^2_0,[s])
	\ar[rr,line width=2mm,dash,color=white,opacity=.7]
	\ar[rr, "f^{2\to 3}_{0;\frx_{s};z}"]
&&
\uCF^-(\cH_{0,z}^3,[s])
	\ar[from=ulll,line width=2mm,dash,color=white,opacity=.7]
	\ar[from=ulll,dashed, "H_{0;w\to z}^{2\to 3}"]
&& |[alias=B]|\uCF^-(\cH^3_{0,w},[s+n])
	\ar[dd, "\Psi_{0\to 1}^3"]
	\ar[from=ulll,line width=2mm,dash,color=white,opacity=.7]
	\ar[from=ulll,dashed, "H^3_{w\to z\to w}"]
	\ar[from=ll,line width=2mm,dash,color=white,opacity=.7]
	\ar[from=ll, "\frF_{z\to w}"]
\\[1.8cm]
\CF^-(\cH^2_{1}, [s]) 
	\ar[rr, "f^{2\to 3}_{1;\fry_{s};w}",pos=.3]
	\ar[dr,"\id"]
	\ar[drrr,dashed,pos=.6, "H^{2\to 3}_{1;w\to z}"]
&&
\uCF^-(\cH^3_{1,w},[s+n])
	\ar[dr, "\frF_{w\to z}"]
	\ar[rr,"\id", pos=.3]
	\ar[drrr,dashed, "H^3_{w\to z\to w}"]
&&
\uCF^-(\cH^3_{1,w},[s+n])
	\ar[dr, "\id"]
 &\,
\\
&
|[alias=D]|\CF^-(\cH^2_{1},[s])
	\ar[rr, "f^{2\to 3}_{1;\frx_{s};z}"]
	\ar[from=uu,line width=2mm,dash,color=white,opacity=.7]
	\ar[from =uu,  "\Psi_{0\to 1}^2"]
&&
\uCF^-(\cH^3_{1,z},[s])
	\ar[from=uu,line width=2mm,dash,color=white,opacity=.7]
	\ar[from=uu, "\Psi_{0\to 1}^3"]
	\ar[rr,"\frF_{z\to w}"]
	\ar[from=uull,line width=2mm,dash,color=white,opacity=.7]
	\ar[from=uull,dashed, sloped, "H^{2\to 3}_{0\to 1;\frx_s;z}"]
&&|[alias=C]| \uCF^-(\cH^3_{1,w},[s+n])
	\ar[from=uull,line width=2mm,dash,color=white,opacity=.7]
	\ar[from=uull,dashed, sloped,"H_{0\to 1;z\to w}^3", pos=.4]
\end{tikzcd}
\begin{tikzcd}[
column sep={2.5cm,between origins},
row sep=.55cm,labels=description,fill opacity=.7, text opacity=1,
execute at end picture={
\foreach \Nombre in  {A,B,C,D}
  {\coordinate (\Nombre) at (\Nombre.center);}
\fill[opacity=0.1] 
  (A) -- (B) -- (C) -- (D) -- cycle;
}]
|[alias=A]|\CF^-(\cH^2_0, [s])
	\ar[dd, "\Psi_{0\to 1}^2"]
	\ar[dr,  "\Gamma_s"]
	\ar[rr, "f^{2\to 3}_{0;\frx_s;z}"]
	\ar[ddrr,dashed,sloped, "H_{0\to 1;\frx_s;z}^{2\to 3}", pos=.7]
&&
\uCF^-(\cH^3_{0,z},[s])
	\ar[dd, "\Psi_{0\to 1}^3"]
	\ar[dr,"\theta_z"]
 	\ar[ddrr,dashed,sloped, "H_{0\to 1;z\to w}^{3}", pos=.7]
	\ar[rr, "\frF_{z\to w}"]
&
&|[alias=B]|\uCF^-(\cH^3_{w},[s+n])
	\ar[dd, "\Psi_{0\to 1}^3"]
	\ar[dr,"\theta_w"]
&\,
\\
&A_{s}(\cH^3_{0})
	\ar[rr,line width=2mm,dash,color=white,opacity=.7]
	\ar[rr, "\vopp"]
&&
\Bopp_{s}(\cH^3_{0})
	\ar[rr, "\frF_{z\to w}"]
&& B_{s+n}(\cH^3_{0})
	\ar[dd, "\Psi_{0\to 1}^3"]
\\[1.8cm]
|[alias=D]|
\CF^-(\cH^2_{1}, [s]) 
	\ar[rr, "f^{2\to 3}_{1;\frx_s;z}",pos=.76]
	\ar[dr,"\Gamma_s"]
	&
&\uCF^-(\cH^3_{1,z},[s])
	\ar[dr, "\theta_z"]	
	\ar[rr, "\frF_{z\to w}"]
&
&|[alias=C]|\uCF^-(\cH^3_{1,w},[s+n])
	\ar[dr, "\theta_w"]
\\
&
A_{s}(\cH^3_{1})
	\ar[rr, "\vopp"]
	\ar[from =uu,line width=2mm,dash,color=white,opacity=.7]
	\ar[from =uu,  "\Psi_{0\to 1}^3"]
	\ar[from=uuul,line width=2mm, dash, color=white, opacity=.7]
	\ar[from=uuul, "J_{0\to 1;s}",dashed]
	&&
\Bopp_{s}(\cH^3_{1})
\ar[rr, "\frF_{z\to w}"]
	\ar[from=uu,line width=2mm,dash,color=white,opacity=.7]
	\ar[from=uu,"\Psi_{0\to 1}^3"]
&& B_{s+n}(\cH^3_{1})
	\ar[from=uull,line width=2mm,dash,color=white,opacity=.7]
	\ar[from=uull,dashed, "H^3_{0\to 1}",sloped, pos=.36]
\end{tikzcd}
\caption{ A schematic of the construction of a hypercube $\cK_{\circ;\fry_s}\in \{
\cK_{\oneh;\fry_s},\cK_{\twoh;\fry_s}^1,\cK_{\twoh;\fry_s}^2,\cK_{\thrh;\fry_s}\}$. We construct $\cK_{\circ;\fry_s}$ (top) as the compression of the hyperbox obtained by stacking and compressing the middle and bottom hyperboxes along their gray faces. We write $\cK_{\circ;\fry_s}^{1}$, $\cK_{\circ;\fry_s}^{2}$, $\cK_{\circ;\fry_s}^{3}$ and $\cK_{\circ;\fry_s}^{4}$ for the  
middle-left, middle-right, bottom-left and bottom-right hypercubes, respectively.}
\label{eq:C-fry-s-Psi-i-full}
\end{figure}

In this section, we describe the construction of the hypercubes $\cK_{\oneh;\fry_s}$,
 $\cK_{\twoh;\fry_s}^1,$  $\cK_{\twoh;\fry_s}^2,$ and $\cK_{\thrh;\fry_s}$. Let $\cK_{\circ;\fry_s}$ denote one of these hypercubes. We will construct $\cK_{\circ;\fry_s}$ by stacking and compressing four hypercubes
 $\cK_{\circ;\fry_s}^{i}$ for $i\in \{1,2,3, 4\}$. The hypercubes $\cK_{\circ;\fry_s}$ and $\cK_{\circ;\fry_s}^{i}$ are schematically illustrated in Figure~\ref{eq:C-fry-s-Psi-i-full}.

In the doubling model of the involution, each of the factors of the involution moves the basepoints slightly, resulting in a cumulative half-twist along $K$ (see Section~\ref{sec:involution-doubling}). 
As a consequence, the basepoints on the top faces of some of the hypercubes in~\eqref{eq:C-hypercubes-overview} differ from the basepoints on the bottom face. In this section, we abuse notation slightly, and denote every basepoint by $w$ or $z$, with the convention that $K$ intersects the Heegaard surface positively (resp. negatively) at $z$ (resp. $w$).

Let $\cH^1_0,$ $\cH^2_0$ and $\cH^3_0$ denote the Heegaard diagrams appearing on the top face of $\cC_{\circ;\fry_s}$, and let $\cH_1^1$, $\cH_1^2$ and $\cH_1^3$ denote those appearing along the bottom face. For $i\in \{0,1\}$, let $\cT_i$ denote the Heegaard triple used to compute the maps from $\bCF^-(\cH_i^2)$ to $\buCF^-(\cH_i^3)$.
We have the following cases, depending on the hypercube $\cK_{\circ,\fry_s}$:

\begin{enumerate}[ref=$K$-\arabic*, label=($K$-\arabic*),  leftmargin=*, widest=IIII]
\item\label{case:1-triple-fry-s}  $\cT_{1}$ is obtained from $\cT_0$ by a stabilization or destabilization with index tuple $(0,0,0)$ (i.e. $\cT_{1}$ is obtained by attaching or removing a collection of 1-handle stabilizations from $\cT_{0}$). This configuration appears in $\cK_{\oneh;\fry_s}$ and $\cK_{\thrh;\fry_s}$.
\item\label{case:2-triple-fry-s}  $\cT_0=(\Sigma,\as^3,\as^2,\bs_0,w,z)$ and $\cT_{1}=(\Sigma,\as^3,\as^2,\bs_1,w,z)$. Furthermore, for $j\in \{2,3\},$ the triple $(\Sigma,\as^j,\bs_0,\bs_1,w,z)$ represents surgery on an $\ell$-component link in 
\[
(Y\setminus N(K) )\# (S^1\times S^2)^{\# \ell}
\]
 which cancels the $S^1\times S^2$-summands. 
 This configuration corresponds to $\cK_{\twoh;\fry_s}^1$.
\item\label{case:3-triple-fry-s} $\cT_0=(\Sigma,\as, \bs^2_0,\bs^3_0,w,z)$ and $\cT_{1}=(\Sigma,\as,\bs^2_1,\bs^3_1,w,z)$. Furthermore the triples $(\Sigma,\as,\bs^2_0,\bs^2_1,w)$ and $(\Sigma,\as,\bs^3_0,\bs^3_1,w,z)$ both represent surgery on a 0-framed unlink in $Y\setminus N(K)$. 
This configuration corresponds to $\cK_{\twoh;\fry_s}^{2}$.
\end{enumerate} 
Note that the notation of the alpha and beta curves in configurations~\eqref{case:1-triple-fry-s}, ~\eqref{case:2-triple-fry-s} and~\eqref{case:3-triple-fry-s} does not coincide with the notation of Figure~\ref{fig:large-hyperbox}. We have condensed it for convenience.

We describe the construction of $\cK_{\twoh;\fry_s}^{1}$ (i.e. case~\eqref{case:2-triple-fry-s}) since the other two configurations are handled using similar reasoning. The hypercube $\cK_{\twoh;\fry_s}^{1,1}$ is depicted schematically in Figure~\ref{eq:hypercube-naturality-flip-map}. We think of this cube as relating 2-dimensional faces of the flip-map hypercubes for the two triples $(\Sigma,\as^3,\as^2,\bs_0,w,z)$ and $(\Sigma,\as^3,\as^2,\bs_1,w,z)$, which appear as the top and bottom faces of Figure~\ref{eq:hypercube-naturality-flip-map}.

\begin{figure}[h]
\begin{tikzcd}[
column sep={3.6cm,between origins},
row sep=.8cm,labels=description]
\CF^-(\as^2,\bs_0,[s])
	\ar[dd, swap,"f_{\a^2}^{\b_0\to \b_1}"]
	\ar[dr,  "\id"]
	\ar[rr, " f_{\a^2\to \a^3;\fry_s;w}^{\b_0}"]
	\ar[ddrr,dashed, sloped, "H_{\a^2\to \a^3;\fry_s,w}^{\b_0\to \b_1}",pos=.3]
&&[-1.5cm]
\uCF^-(\as^3,\bs_0,\hat{w},z,[s+n])
	\ar[dd, "f_{\a^3}^{\b_0\to \b_1}"]
	\ar[dr,"\frF_{w\to z}"]
	\ar[dddr,dashed, "H_{\a^3;w\to z}^{\b_0\to \b_1}", sloped]
&\,
\\
&\CF^-(\as^2,\bs_0,[s])
	\ar[rr,crossing over, "f_{\a^2\to \a^3;\fry_s;z}^{\b_0}"]
&&
\uCF^-(\as^3,\bs_0,w,\hat{z},[s])
	\ar[dd, "f_{\a^3}^{\b_0\to \b_1}"]
	\ar[from=ulll,dashed,crossing over, "H_{\a^2\to \a^3;\fry_s;w\to z}^{\b_0}",sloped]
\\[1.8cm]
\CF^-(\as^2,\bs_1,[s])
	\ar[rr, "f_{\a^2\to \a^3;\fry_s;w}^{\b_1}"]
	\ar[dr,"\id"]
	\ar[drrr,dashed, "H_{\a^2\to \a^3;\fry_s;w\to z}^{\b_1}",sloped, pos=.65]
&&\uCF^-(\as^3,\bs_1,\hat{w},z,[s+n])
	\ar[dr, "\frF_{w\to z}"]	
&
\\
&
\CF^-(\as^2,\bs_1,[s])
	\ar[rr, "f_{\a^2\to \a^3;\fry_s;z}^{\b_1}"]
	\ar[from =uu, crossing over,"f_{\a^2}^{\b_0\to \b_1}"]
	&&
\uCF^-(\as^3,\bs_1,w,\hat{z},[s])
	\ar[from=uull,crossing over, dashed, sloped, "H_{\a^2\to \a^3;\fry_s;z}^{\b_0\to \b_1}"]
\end{tikzcd}
\caption{The hypercube $\cK^{1,1}_{\twoh;\fry_s}$ (i.e. case~\eqref{case:2-triple-fry-s}). There is a length 3 arrow which is not shown.}
\label{eq:hypercube-naturality-flip-map}
\end{figure}

 We begin by describing the maps appearing on the back and front faces of $\cK_{\twoh;\fry_s}^{1,1}$.   We pick cycles $\Theta_{\a^3,\a^2}^+$ and $\Theta_{\b_0,\b_1}^+$, representing the top degree generators. The length 1 maps appearing along the back face are the holomorphic triangle maps with inputs $\Theta_{\a^3,\a^2}^+$ and $\Theta_{\b_0,\b_1}^+$. In the case of $f_{\a^2\to \a^3;\fry_s;w}^{\b_0}$, we further restrict to triangles where $\frs_w(\psi)=\fry_s$. The length 2 map $H_{\a^2\to \a^3;\fry_s;w}^{\b_0\to \b_1}$ similarly counts rectangles with the two special inputs $\Theta_{\a^3,\a^2}^+$ and $\Theta_{\b_0,\b_1}^+$, again with the restriction that $\frs_w(\psi)=\fry_s$. Along the back face, we use the weights $U^{n_w(\psi)} T^{n_z(\psi)-n_w(\psi)}$. The maps along the front face are defined similarly, but instead use $z$ as the special basepoint, and are weighted by $U^{n_z(\psi)} T^{n_z(\psi)-n_w(\psi)}$.

  The remaining maps in the hypercube $\cK_{\twoh;\fry_s}^{1,1}$ are constructed by stacking and compressing 5 hypercubes, mirroring the construction of the flip-map hypercube in Section~\ref{sec:flip-hypercube}. In analogy to the construction therein, the five hypercubes are  as follows:
  \begin{enumerate}[ref=$K$-$2$-\arabic*, label=($K$-2-\arabic*),  leftmargin=*, widest=IIIIII]
  \item \label{cube:Psi-2-1} One hypercube involving the map $\theta_w$ and $T^{-n}$.
  \item \label{cube:Psi-2-2} One hypercube involving the tautological map for a diffeomorphism $\phi$ which moves $w$ to $z$.
  \item \label{cube:Psi-2-3} Two hypercubes for changing $\phi(\as^3)$ to $\as^3$.
  \item \label{cube:Psi-2-4} One hypercube involving the map $\theta_z$.
  \end{enumerate}
The hypercube labeled \eqref{cube:Psi-2-1}, involving $\theta_w$, is shown in Figure~\ref{eq:hypercube-naturality-flip-map-1}, and the hypercube ~\eqref{cube:Psi-2-4} involving $\theta_z$ is similar, though with no weight by a power of $T$ on the $\theta_z$ arrows.
\begin{figure}[h]
\begin{tikzcd}[
column sep={3.7cm,between origins},
row sep=.8cm,labels=description]
\CF^-(\as^2,\bs_0,[s])
	\ar[dd, swap,"f_{\a^2}^{\b_0\to \b_1}"]
	\ar[dr,  "\id"]
	\ar[rr, " f_{\a^2\to \a^3;\fry_s;w}^{\b_0}"]
	\ar[ddrr,dashed, "H_{\a^2\to \a^3;\fry_s;w}^{\b_0\to \b_1}", sloped]
&&[-1.5cm]
\uCF^-(\as^3,\bs_0,\hat{w},z,[s+n])
	\ar[dd, "f_{\a^3}^{\b_0\to \b_1}"]
	\ar[dr,"T^{-n}\cdot\theta_w"]
&
\\
&\CF^-(\as^2,\bs_0,[s])
	\ar[rr,crossing over,  "f_{\a^2\to \a^3;\fry_s}^{\b_0}\otimes T^s"]
&&
\CF^-(\as^3,\bs_0,w)\otimes T^s
	\ar[dd, "f_{\a^3}^{\b_0\to \b_1}"]
\\[1.8cm]
\CF^-(\as^2,\bs_1,[s])
	\ar[rr, "f_{\a^2\to \a^3;\fry_s;w}^{\b_1}"]
	\ar[dr,"\id"]
&&\uCF^-(\as^3,\bs_1,\hat{w},z,[s+n])
	\ar[dr, "T^{-n}\cdot \theta_w"]	
&
\\
&
\CF^-(\as^2,\bs_1,[s])
	\ar[rr, "f_{\a^2\to \a^3;\fry_s}^{\b_1}\otimes T^s"]
	\ar[from =uu, crossing over,"f_{\a^2}^{\b_0\to \b_1}"]
	&&
\CF^-(\as^3,\bs_1,w)\otimes T^s
	\ar[from=uull,crossing over, dashed,sloped, "H_{\a^2\to \a^3;\fry_s}^{\b_0\to \b_1}\otimes T^s"]
\end{tikzcd}
\caption{The hypercube labeled~\eqref{cube:Psi-2-1} in the construction of $\cK_{\twoh;\fry_s}^{1,1}$.}
\label{eq:hypercube-naturality-flip-map-1}
\end{figure}

The hypercube involving the tautological map $\phi_*$, labeled~\eqref{cube:Psi-2-2},  is straightforward to construct, and has the map $\phi_*$ in the direction which comes out of the page, and no length 2 maps except on the front and back faces.

We now describe the two hypercubes which move $\phi(\as^3)$ to $\as^3$, and are labeled~\eqref{cube:Psi-2-3} above. We construct these two hypercubes by pairing  a hyperbox of alpha attaching curves $\cL_{\a}$ with a hyperbox of beta attaching curves, as shown in~\eqref{eq:hypercubes-attaching-curves-Psi-2-3}. Therein, $\as^{2}{}'$ denotes a small isotopy of $\as^2$. The curves $\as^{3}_{\twind}$ are obtained by winding $\as^3$ to achieve admissibility, in a way which does not intersect the boundary of the 1-handle region.

\begin{equation}
\cL_{\a}:= 
\begin{tikzcd}[column sep=2cm, row sep=2cm,labels=description]
\as^2
	\ar[r, "\Theta_{\phi(\a^3), \a^2}^+"] 
	\ar[d, "\Theta_{\a^2{}',\a^2}^+"]
	\ar[dr,dashed, "\lambda_{\a^{3}_{\twind}, \a^2}"]
&
\phi(\as^3)
	\ar[d, "\Theta^+_{\a^{3}_{\twind},\phi(\a^3)}"]
\\
\as^2{}'
	\ar[d, "\Theta_{\a^2,\a^2{}'}^+"]
	\ar[r, "\Theta_{\a^{3}_{\twind}, \a^2{}'}^+"]
	\ar[dr, dashed, "\lambda_{\a^3,\a^2{}'}"]
&
\as^{3}_{\twind}
	\ar[d, "\Theta_{\a^3,\a^{3}_{\twind}}^+"]
\\
\as^2
	\ar[r, "\Theta_{\a^3,\a^2}^+"]
&
\as^3
\end{tikzcd}
\qquad \qquad
\cL_{\b}:=\begin{tikzcd}[column sep=2cm,labels=description]
\bs_0
\ar[r, "\Theta_{\b_0,\b_1}^+"]
&\bs_1
 \end{tikzcd}
\label{eq:hypercubes-attaching-curves-Psi-2-3}
\end{equation}

Next, we discuss the hypercube $\cK_{\twoh;\fry_s}^{1,2}$, which is schematically  shown in Figure~\ref{eq:C-fry-s-Psi-i-2-concrete}. The hypercube $\cK_{\twoh;\fry_s}^{1,2}$ is defined as the compression of a hypercube of size $(5,5,1)$. The top face of this hypercube is shown in Figure~\ref{fig:K-y-s-Phi-2-top-face}. The bottom face is similar, but has $\bs_1$ instead of $\bs_0$. We leave the remaining details of the construction of $\cK_{\twoh;\fry_s}^{1,2}$ to the reader, since the construction is similar to the description of $\cK_{\frF;\frx_s}^{1}$ in Figure~\ref{eq:C-frx-s-flip-1}. Similarly, we leave the construction of $\cK_{\twoh;\fry_s}^{1,3}$ and $\cK_{\twoh;\fry_s}^{1,4}$ to the reader, since they are similar to $\cK_{\frF;\frx_s}^{3}$ and $\cK_{\frF;\frx_s}^{4}$, respectively.

\begin{figure}[h]
\begin{tikzcd}[
column sep={3cm,between origins},
row sep=.5cm,labels=description]
\uCF^-(\as^3,\bs_0,\hat{w},z,[s+n])
	\ar[dd, "f_{\a^3}^{\b_0\to \b_1}"]
	\ar[dr,"\frF_{w\to z}"]
	\ar[rr, "\id"]
&&
\uCF^-(\as^3,\bs_0,\hat{w},z,[s+n])
	\ar[dd, "f_{\a^3}^{\b_0\to \b_1}"]
	\ar[dr,"\id"]
&\,
\\
&
\uCF^-(\as^3,\bs_0,w,\hat{z},[s])
	\ar[rr, "\frF_{z\to w}"]
&& \uCF^-(\as^3,\bs_0,\hat{w},z,[s+n])
	\ar[dd, "f_{\a^3}^{\b_0\to \b_1}"]
	\ar[from=ulll, crossing over, dashed, sloped, "H^{\frF}_{w\to z\to w}"]
\\[1.8cm]
\uCF^-(\as^3,\bs_1,\hat{w},z,[s+n])
	\ar[dr, "\frF_{w\to z}"]	
	\ar[rr, "\id"]
	\ar[drrr,dashed]
&&
\uCF^-(\as^3,\bs_1,\hat{w},z,[s+n])
	\ar[dr, "\id"]
\\
&
\uCF^-(\as^3,\bs_1,w,\hat{z},[s])
	\ar[rr, "\frF_{z\to w} "]
	\ar[from=uu,crossing over, "f_{\a^3}^{\b_0\to \b_1}"]
	\ar[from=uuul, crossing over, dashed]
&& \uCF^-(\as^3,\bs_1,\hat{w},z,[s+n])
	\ar[from=uull,line width=2mm,dash,color=white,opacity=.7]
	\ar[from=uull,dashed]
\end{tikzcd}
\caption{The hypercube $\cK^{1,2}_{\twoh;\fry_s}$ from case~\eqref{case:2-triple-fry-s}. There is a length 3 map, which is not shown.}
\label{eq:C-fry-s-Psi-i-2-concrete}
\end{figure}

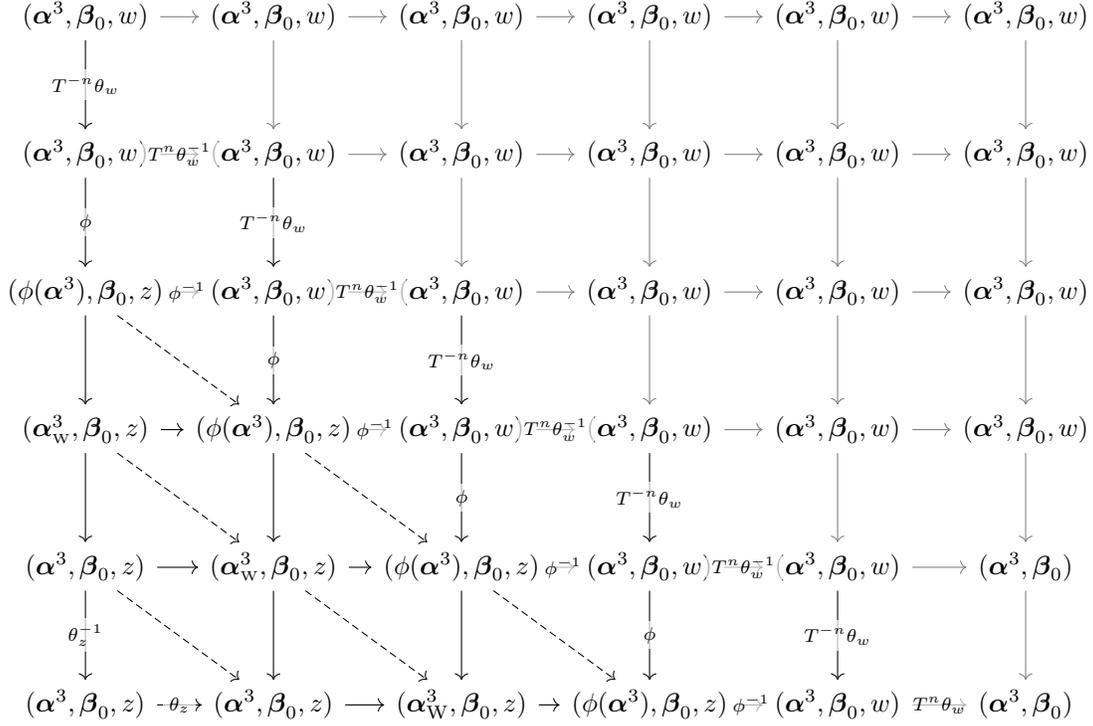
\begin{figure}[h]
\begin{tikzcd}
[column sep={2.5cm,between origins},row sep=1.2cm,labels=description,fill opacity=.7, text opacity=1]
(\as^3,\bs_0,w)
	\ar[r, color=gray]
	\ar[d, "T^{-n} \theta_w"]
&
(\as^3,\bs_0,w)
	\ar[r, color=gray]
	\ar[d, color=gray] 
&
(\as^3,\bs_0,w)
	\ar[d,color=gray]
	\ar[r,color=gray]
&
(\as^3,\bs_0,w)
	\ar[d,color=gray]
	\ar[r,color=gray]
&
(\as^3,\bs_0,w)
	\ar[r,color=gray]
	\ar[d,color=gray]
&
(\as^3,\bs_0,w)
	\ar[d,color=gray]
\\
(\as^3,\bs_0,w)
	\ar[r, "T^n\theta_w^{-1}"]
	\ar[d, "\phi"]
&
(\as^3,\bs_0,w)
	\ar[r, color=gray]
	\ar[d, "T^{-n} \theta_w"] 
&
(\as^3,\bs_0,w)
	\ar[d,color=gray]
	\ar[r,color=gray]
&
(\as^3,\bs_0,w)
	\ar[d,color=gray]
	\ar[r,color=gray]
&
(\as^3,\bs_0,w)
	\ar[d,color=gray]
	\ar[r,color=gray]
&
(\as^3,\bs_0,w)
	\ar[d,color=gray]
\\
(\phi(\as^3),\bs_0,z)
	\ar[r, "\phi^{-1}"]
	\ar[d]
	\ar[dr,dashed]
&
(\as^3,\bs_0,w)
	\ar[r, "T^n \theta_w^{-1}"]
	\ar[d,"\phi"]
&
(\as^3,\bs_0,w)
	\ar[r,color=gray]
	\ar[d,"T^{-n} \theta_w"]
&
(\as^3,\bs_0,w)
	\ar[d,color=gray]
	\ar[r,color=gray]
&
(\as^3,\bs_0,w)
	\ar[r,color=gray]
	\ar[d,color=gray]
&
(\as^3,\bs_0,w)
	\ar[d,color=gray]
\\
(\as^3_{\wind},\bs_0,z)
	\ar[r]
	\ar[d]
	\ar[dr,dashed]
&
(\phi(\as^3),\bs_0,z)
	\ar[r, "\phi^{-1}"]
	\ar[d]
	\ar[dr,dashed]
&
(\as^3,\bs_0,w)
	\ar[r, "T^n \theta^{-1}_w"]
	\ar[d,"\phi"]
&
(\as^3,\bs_0,w)
	\ar[r,color=gray]
	\ar[d,"T^{-n}\theta_w"]
&
(\as^3,\bs_0,w)
	\ar[r,color=gray]
	\ar[d,color=gray]
&
(\as^3,\bs_0,w)
	\ar[d,color=gray]
\\
(\as^3,\bs_0,z)
	\ar[r]
	\ar[d, "\theta_z^{-1}"]
	\ar[dr,dashed]
&
(\as^3_{\wind},\bs_0,z)
	\ar[r]
	\ar[d]
	\ar[dr,dashed]
&
(\phi(\as^3),\bs_0,z)
	\ar[r, "\phi^{-1}"]
	\ar[d]
	\ar[dr,dashed]
&
(\as^3,\bs_0,w)
	\ar[r, "T^n\theta_w^{-1}"]
	\ar[d, "\phi"]
&
(\as^3,\bs_0,w)
	\ar[d,  "T^{-n} \theta_w"]
	\ar[r,color=gray]
&
(\as^3,\bs_0)
	\ar[d,color=gray]
\\
(\as^3,\bs_0,z)
	\ar[r, "\theta_{z}"]
&
(\as^3,\bs_0,z)
	\ar[r]
&
(\as^3_{\twind},\bs_0,z)
	\ar[r]
&
(\phi(\as^3),\bs_0,z)
	\ar[r, "\phi^{-1}"]
&
(\as^3,\bs_0,w)
	\ar[r, "T^n \theta_w"]
&
(\as^3,\bs_0)
\end{tikzcd}
\caption{The top face of the hyperbox used to construct $\cK_{\twoh;\fry_s}^{1;2}$. The bottom face has $\bs_1$ in place of $\bs_0$. A gray, unmarked arrow indicates the identity map.}
\label{fig:K-y-s-Phi-2-top-face}
\end{figure}

\subsection{The hypercube \texorpdfstring{$\cK_{\cen;\fry_s}^{(1)}$}{K-cen-y-s-(1)}}

We now construct the hypercube $\cK_{\cen;\fry_s}^{(1)}$. The back face of $\cK_{\cen;\fry_s}^{(1)}$ will match a $\Spin^c$ restriction the front face of the hypercube $\cC_{\cen}^{(1)}$ from Section~\ref{sec:central-hypercube}.  The hypercube $\cK_{\cen;\fry_s}^{(1)}$ is constructed as a compression of four hypercubes $\cK_{\cen;\fry_s}^{(1),j}$, for $j\in \{1,2,3,4\}$, which fit into a similar configuration as the hypercubes $\cK_{i;\fry_s}^{\Psi,j}$ in Figure~\ref{eq:C-fry-s-Psi-i-full}. As in the previous section, we write $w$ and $z$ for the basepoints on the diagrams in the central hypercube, even though they have moved in the doubling operation from their original positions on $\cH^3$ (see Section~\ref{sec:doubled-knots}).

\begin{figure}[h]
\begin{tikzcd}[
column sep={3.5cm,between origins},
row sep=.7cm,labels=description,
fill opacity=.7, text opacity=1,
execute at end picture={
\foreach \Nombre in  {A,B,C,D}
  {\coordinate (\Nombre) at (\Nombre.center);}
\fill[opacity=0.1] 
  (A) -- (B) -- (C) -- (D) -- cycle;
}
]
\CF^-(\gs',\ds,[s])
	\ar[dd, swap,"\id"]
	\ar[dr,  "\id"]
	\ar[rr, "f_{\g'\to \g'';\fry_s;w}^{\dt}"]
	\ar[ddrr,dashed,sloped,near end]
&&[-1.5cm]
|[alias=A]|\uCF^-(\gs'',\ds,\hat{w},z,[s+n])
	\ar[dd, "\sigma F_3 f_{\g''}^{\dt\to \dt'}"]
	\ar[dr,"\frF_{w\to z}"]
	\ar[dddr,dashed]
&
\\
&\CF^-(\gs',\ds,[s])
	\ar[rr,"f_{\g'\to \g'';\fry_s;z}^{\dt}",crossing over]
&&
|[alias=B]|\uCF^-(\gs'',\ds,w,\hat{z},[s])
	\ar[dd, "\sigma F_3 f_{\g''}^{\dt\to \dt'}"]
	\ar[from=ulll,dashed,crossing over]
\\[1.8cm]
\CF^-(\gs',\ds,[s])
	\ar[rr, "f_{\g';\fry_s;w}^{\dt\to \dt'}"]
	\ar[dr,"\id"]
	\ar[drrr,dashed]
&&|[alias=D]|\uCF^-(\gs',\ds', \hat{w},z,[s+n])
	\ar[dr, "\frF_{w\to z}"]	
&
\\
&
\CF^-(\gs',\ds,[s])
	\ar[rr, "f_{\g';\fry_s;z}^{\dt\to \dt'}"]
	\ar[from =uu, crossing over,"\id"]
	&&
|[alias=C]|\uCF^-(\gs',\ds',w,\hat{z},[s])
	\ar[from=uull,crossing over, dashed]
\end{tikzcd}
\begin{tikzcd}[
column sep={3.5cm,between origins},
row sep=.7cm,labels=description,
fill opacity=.7, text opacity=1,
execute at end picture={
\foreach \Nombre in  {A,B,C,D}
  {\coordinate (\Nombre) at (\Nombre.center);}
\fill[opacity=0.1] 
  (A) -- (B) -- (C) -- (D) -- cycle;
}
]
|[alias=A]|\uCF^-(\gs'',\ds, \hat{w},z,[s+n])
	\ar[dd, swap,"\sigma F_3 f_{\g''}^{\dt\to \dt'}"]
	\ar[dr,  "\frF_{w\to z}"]
	\ar[rr, "\id"]
	\ar[ddrr,dashed,sloped,near end]
&&[-1.5cm]
\uCF^-(\gs'',\ds,\hat{w},z,[s+n])
	\ar[dd, "\sigma F_3 f_{\g''}^{\dt\to \dt'}"]
	\ar[dr,"\id"]
&
\\
&|[alias=B]|\uCF^-(\gs'',\ds, w,\hat{z},[s])
	\ar[rr,"\frF_{z\to w}",crossing over]
&&
\uCF^-(\gs'',\ds,w,\hat{z},[s+n])
	\ar[dd, "\sigma F_3 f_{\g''}^{\dt\to \dt'}"]
	\ar[from=ulll,dashed,crossing over]
\\[1.8cm]
|[alias=D]|\uCF^-(\gs',\ds',\hat{w},z,[s+n])
	\ar[rr, "\id"]
	\ar[dr,"\frF_{w\to z}"]
	\ar[drrr,dashed]
&&\uCF^-(\gs',\ds', \hat{w},z,[s+n])
	\ar[dr, "\id"]	
&
\\
&
|[alias=C]|\uCF^-(\gs',\ds',w,\hat{z},[s])
	\ar[rr, "\frF_{z\to w}"]
	\ar[from =uu, crossing over,"\sigma F_3 f_{\g''}^{\dt\to \dt'}"]
	\ar[from=uuul,dashed, crossing over]
	&&
\uCF^-(\gs',\ds',\hat{w},z,[s+n])
	\ar[from=uull,crossing over, dashed]
\end{tikzcd}
\\
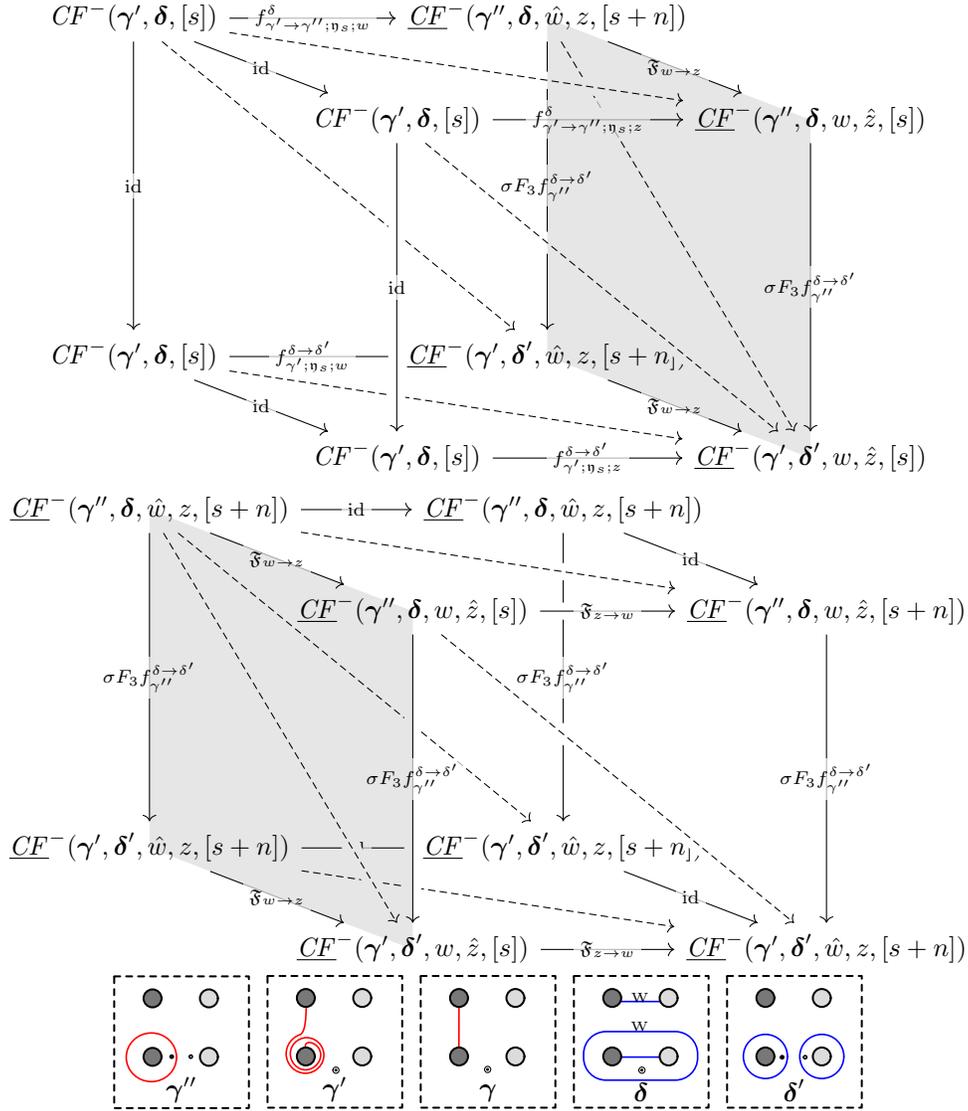
\caption{The hypercubes $\cK_{\cen;\fry_s}^{(1),1}$ (top) and $\cK_{\cen;\fry_s}^{(1),2}$ (middle). Both cubes have a length 3 arrow which is not drawn. The cubes are stacked along the gray faces. On the bottom are the attaching curves in the surgery region.}
\label{eq:hypercube-naturality-flip-map-C1-cen}
\end{figure}

We focus on the construction of $\cK_{\cen;\fry_s}^{(1),1}$, since $\cK_{\cen;\fry_s}^{(1),2}$ is simpler, and furthermore the hypercubes $\cK_{\cen;\fry_s}^{(1),3}$ and $\cK_{\cen;\fry_s}^{(1),4}$ involve only the maps which appeared on $\cK_{\cen;\fry_s}^{(1),1}$ and $\cK_{\cen;\fry_s}^{(1),2}$, as well as canonical isomorphisms and inclusions. The hypercubes $\cK_{\cen,\fry_s}^{(1),1}$ and $\cK_{\cen;\fry_s}^{(1),2}$ take the schematic form shown in  Figure~\ref{eq:hypercube-naturality-flip-map-C1-cen}. 

We first describe the construction of $\cK_{\cen;\fry_s}^{(1),1}$. The length 1 and 2 maps along the back face are constructed the same as in  $\cC_{\cen}^{(1)}$; see Figure~\ref{fig:176}. The length 1 and 2 maps along the front face are similar, except they use $z$ as the special basepoint, instead of $w$.


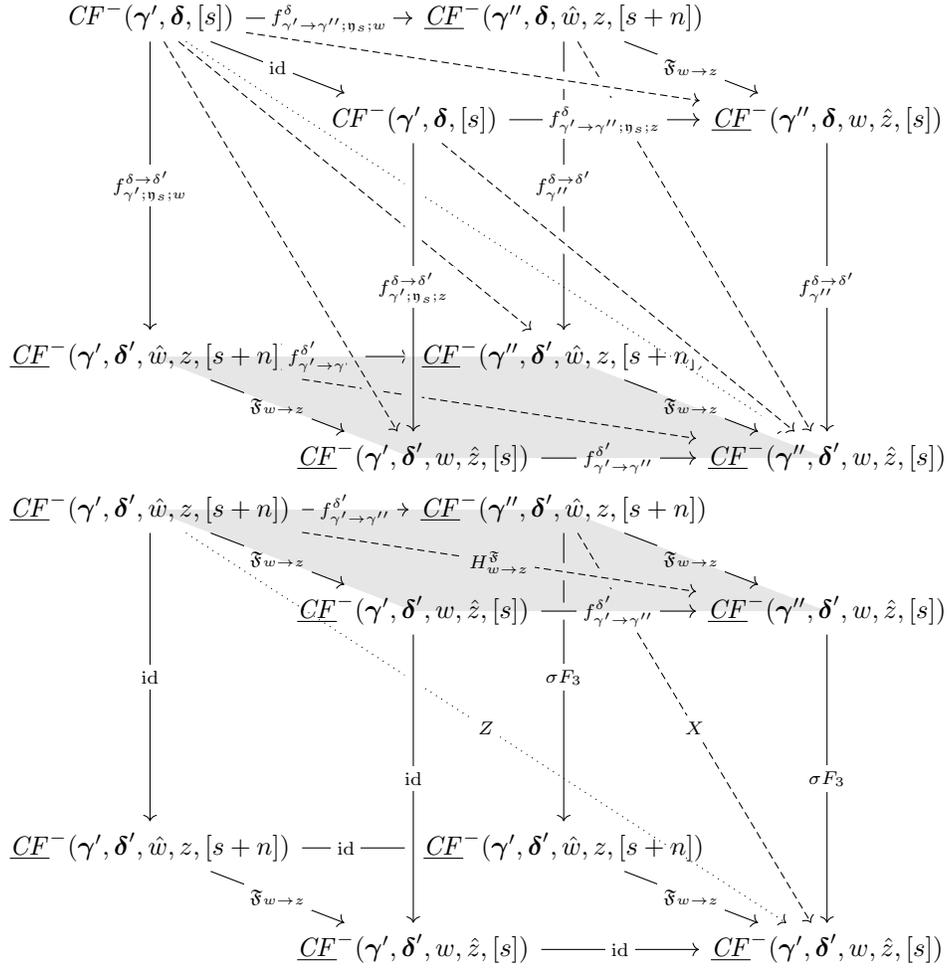
\begin{figure}[h]
\begin{tikzcd}[
column sep={3.5cm,between origins},
row sep=.7cm,labels=description,
execute at end picture={
\foreach \Nombre in  {A,B,C,D}
  {\coordinate (\Nombre) at (\Nombre.center);}
\fill[opacity=0.1] 
  (A) -- (B) -- (C) -- (D) -- cycle;
}]
\CF^-(\gs',\ds,[s])
	\ar[dd, "f_{\g';\fry_s;w}^{\dt\to \dt'}"]
	\ar[dr,  "\id"]
	\ar[rr, "f_{\g'\to \g'';\fry_s;w}^{\dt}"]
	\ar[ddrr,dashed,near end]
	\ar[dddrrr,dotted]
&&[-1.5cm]
\uCF^-(\gs'',\ds,\hat{w},z,[s+n])
	\ar[dd, "f_{\g''}^{\dt\to \dt'}"]
	\ar[dr,"\frF_{w\to z}"]
	\ar[dddr,dashed]
&
\\
&\CF^-(\gs',\ds,[s])
	\ar[rr,crossing over, "f_{\g'\to \g'';\fry_s;z}^{\dt}"]
&&
\uCF^-(\gs'',\ds,w,\hat{z},[s])
	\ar[dd, "f_{\g''}^{\dt\to \dt'}"]
	\ar[from=ulll,dashed,crossing over]
\\[1.8cm]
|[alias=A]|\uCF^-(\gs',\ds',\hat{w},z,[s+n])
	\ar[rr, "f_{\g'\to \g''}^{\dt'}", pos=.2]
	\ar[dr,"\frF_{w\to z}"]
	\ar[drrr,dashed]
&&
|[alias=B]|\uCF^-(\gs'',\ds', \hat{w},z,[s+n])
	\ar[dr, "\frF_{w\to z}"]	
&
\\
&
|[alias=D]|\uCF^-(\gs',\ds',w,\hat{z},[s])
	\ar[rr, "f_{\g'\to \g''}^{\dt'}"]
	\ar[from =uu, crossing over, "f_{\g';\fry_s;z}^{\dt\to \dt'}"]
	\ar[from=uuul,dashed, crossing over]
	&&
|[alias=C]|\uCF^-(\gs'',\ds',w,\hat{z},[s])
	\ar[from=uull,crossing over, dashed]
\end{tikzcd}

\begin{tikzcd}[
column sep={3.5cm,between origins},
row sep=.7cm,labels=description,
execute at end picture={
\foreach \Nombre in  {A,B,C,D}
  {\coordinate (\Nombre) at (\Nombre.center);}
\fill[opacity=0.1] 
  (A) -- (B) -- (C) -- (D) -- cycle;
}]
|[alias=A]|\uCF^-(\gs',\ds',\hat{w},z,[s+n])
	\ar[dd, "\id"]
	\ar[dr,  "\frF_{w\to z}"]
	\ar[rr, "f_{\g'\to \g''}^{\dt'}"]
	\ar[dddrrr,dotted,"Z"]
&&[-1.5cm]
|[alias=B]|\uCF^-(\gs'',\ds',\hat{w},z,[s+n])
	\ar[dd, "\sigma F_3"]
	\ar[dr,"\frF_{w\to z}"]
	\ar[dddr,dashed, "X"]
&
\\
&
|[alias=D]|\uCF^-(\gs',\ds', w,\hat{z},[s])
	\ar[rr,crossing over, "f_{\g'\to \g''}^{\dt'}"]
&&
|[alias=C]|\uCF^-(\gs'',\ds',w,\hat{z},[s])
	\ar[dd, "\sigma F_3"]
	\ar[from=ulll,dashed,crossing over, "H^{\frF}_{w\to z}"]
\\[1.8cm]
\uCF^-(\gs',\ds',\hat{w},z,[s+n])
	\ar[rr, "\id", pos=.4]
	\ar[dr,"\frF_{w\to z}"]
&&\uCF^-(\gs',\ds', \hat{w},z,[s+n])
	\ar[dr, "\frF_{w\to z}"]	
&
\\
&
\uCF^-(\gs',\ds',w,\hat{z},[s])
	\ar[rr, "\id"]
	\ar[from =uu, crossing over, "\id"]
	&&
\uCF^-(\gs',\ds',w,\hat{z},[s])
\end{tikzcd}
\caption{The hypercube $\cY_r^{(1)}$ is the compression of these two hypercubes, stacked vertically. }
\label{eq:hypercube-naturality-flip-map-C1-cen-2}
\end{figure}

\begin{figure}[h]
\begin{tikzcd}[
column sep={3.5cm,between origins},
row sep=.8cm,labels=description]
\CF^-(\gs',\ds,[s])
	\ar[dd, "\id"]
	\ar[dr, "\id"]
	\ar[rr, "\id"]
&&[-1.5cm]
\CF^-(\gs',\ds,[s])
	\ar[dd, "f_{\g';\fry_s;w}^{\dt\to \dt'}"]
	\ar[dr,"\id"]
	\ar[dddr,dashed]
&
\\
&\CF^-(\gs',\ds,[s])
	\ar[rr,crossing over, "\id"]
&&
\CF^-(\gs',\ds,[s])
	\ar[dd, "f_{\g';\fry_s;z}^{\dt\to \dt'}"]
\\[1.8cm]
\uCF^-(\gs',\ds,[s])
	\ar[rr, "f_{\g';\fry_s;w}^{\dt\to \dt'}", pos=.4]
	\ar[dr,"\id"]
	\ar[drrr,dashed]
&&\uCF^-(\gs',\ds', \hat{w},z,[s+n])
	\ar[dr, "\frF_{w\to z}"]	
&
\\
&
\CF^-(\gs',\ds,[s])
	\ar[rr, "f_{\g';\fry_s;z}^{\dt\to\dt'}"]
	\ar[from =uu, crossing over, "\id"]
	&&
\uCF^-(\gs',\ds',w,\hat{z},[s])
\end{tikzcd}
\caption{The hypercube $\cY_l^{(1)}$. The length 2 maps along the bottom and right side coincide, and there is no length 3 map.}
\label{fig:hypercube-K-l}
\end{figure}
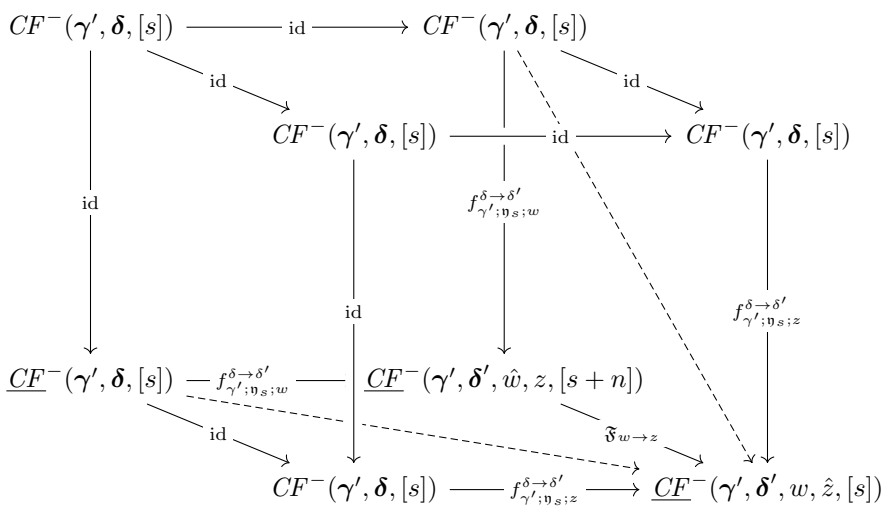

The remaining maps  in $\cK^{ (1),1}_{\cen;\fry_s}$ are obtained by constructing $\cK_{\cen;\fry_s}^{(1),1}$ by stacking and compressing two hypercubes $\cY_{l}^{(1)}$ and $\cY_{r}^{(1)}$, which we stack and compress along the left face of $\cY_{r}^{(1)}$ and the right face of $\cY_{l}^{(1)}$. The hypercube $\cY_r^{(1)}$ is shown in Figure~\ref{eq:hypercube-naturality-flip-map-C1-cen-2}. The hypercube $\cY_l^{(1)}$ is shown in Figure~\ref{fig:hypercube-K-l}.

We now describe the top hypercube of Figure~\ref{eq:hypercube-naturality-flip-map-C1-cen-2}, which features in the construction of $\cY_{r}^{(1)}$. This hypercube is constructed similarly to the flip-map hypercube. There is a new subtlety in that the canonical short path between $w$ and $z$ now intersects both an alpha curve and a beta curve. Hence, we define $\cY_{r}^{(1)}$ as the compression of a hyperbox of size $(1,7,1)$. The first constituent hypercube corresponds to $\theta_w$, the second corresponds to the canonical diffeomorphism moving $w$ to $z$. Next, there are four hypercubes: two for changing the alpha curves,and two for changing the beta curves. Finally, there is a hypercube involving $\theta_w$.

Next, we describe the bottom hypercube of Figure~\ref{eq:hypercube-naturality-flip-map-C1-cen-2}. All of the maps except $X$ and $Z$ have already been constructed. Our construction of $X$ and $Z$ requires $Y$ to be an L-space, as follows. Firstly, we pick $X$ to be any homogeneously graded map which is a homotopy between $\sigma F_3 \frF_{w\to z}$ and $\frF_{w\to z} \sigma F_3$. Such a map $X$ exists  by the following reasoning. Let $V=\bF\oplus \bF[1]$ (i.e. two copies of $\bF$, with the grading of one copy shifted down by 1). For all $s$, there are homotopy equivalences
\begin{equation}
\uCF^-(\gs',\ds',\hat{w},z,[s])\simeq \bF[U]\quad \text{and} \quad \uCF^-(\gs'',\ds',\hat{w},z,[s])\simeq \bF[U]\otimes V,
\label{eq:L-space-simplifies-situation}
\end{equation}
and similarly when $z$ is the special basepoint.
We note that $\Hom_{\bF[U]}(\bF[U]\otimes V, \bF[U])$ has rank 1 over $\bF$ in non-positive gradings (and is $(0)$ in positive gradings). Hence, if two $\bF[U]$-equivariant chain maps from $\uCF^-(\gs'',\ds',\hat{w},z,[s])$ to $\uCF^-(\gs'',\ds',\hat{w},z,[s'])$ have the same homogeneous grading, and induce non-zero elements of $\Hom_{\bF[U]}(\bF[U]\otimes V, \bF[U])$, then they are chain homotopic. This gives the map $X$.

We now construct $Z$. Let $M$ denote the sum of compositions of length 1 and 2 maps appearing in the length 3 relation of the cube:
\[
M:= X f_{\g'\to \g''}^{\dt'}+\sigma F_3 H^{\frF}_{w\to z}.
\]
The map $M$ maps $\uCF^-(\gs',\ds',\hat{w},z,[s+n])$ to $\uCF^-(\gs',\ds',w,\hat{z}, [s])$.  Equation~\eqref{eq:L-space-simplifies-situation} implies that these chain complexes, which are clearly homotopy equivalent as absolutely graded chain complexes via the map $\frF_{w\to z}$, are also both homotopy equivalent to $\bF[U]$. Since $M$ is a chain map and has grading $+1$, it must be chain homotopic to 0. We let $Z$ by any homogeneously graded map satisfying
\[
[\d, Z]=M.
\]

Finally, we consider the hypercube $\cY_l^{(1)}$, which is shown in Figure~\ref{fig:hypercube-K-l}. In this hypercube, we use the same length 2 map along the right face as along the bottom face. This map was constructed in $\cY_r^{(2)}$.

\begin{figure}[h]
\begingroup%
  \makeatletter%
  \providecommand\color[2][]{%
    \errmessage{(Inkscape) Color is used for the text in Inkscape, but the package 'color.sty' is not loaded}%
    \renewcommand\color[2][]{}%
  }%
  \providecommand\transparent[1]{%
    \errmessage{(Inkscape) Transparency is used (non-zero) for the text in Inkscape, but the package 'transparent.sty' is not loaded}%
    \renewcommand\transparent[1]{}%
  }%
  \providecommand\rotatebox[2]{#2}%
  \newcommand*\fsize{\dimexpr\f@size pt\relax}%
  \newcommand*\lineheight[1]{\fontsize{\fsize}{#1\fsize}\selectfont}%
  \ifx\svgwidth\undefined%
    \setlength{\unitlength}{219.06107713bp}%
    \ifx\svgscale\undefined%
      \relax%
    \else%
      \setlength{\unitlength}{\unitlength * \real{\svgscale}}%
    \fi%
  \else%
    \setlength{\unitlength}{\svgwidth}%
  \fi%
  \global\let\svgwidth\undefined%
  \global\let\svgscale\undefined%
  \makeatother%
  \begin{picture}(1,0.27253986)%
    \lineheight{1}%
    \setlength\tabcolsep{0pt}%
    \put(0,0){\includegraphics[width=\unitlength,page=1]{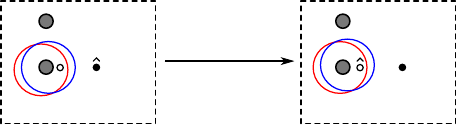}}%
    \put(0.1645526,0.02127871){\color[rgb]{0,0,0}\makebox(0,0)[t]{\lineheight{1.25}\smash{\begin{tabular}[t]{c}$(\gs'',\ds',\hat{w},z)$\end{tabular}}}}%
    \put(0.82850194,0.02225233){\color[rgb]{0,0,0}\makebox(0,0)[t]{\lineheight{1.25}\smash{\begin{tabular}[t]{c}$(\gs'',\ds',w,\hat{z})$\end{tabular}}}}%
    \put(0.49408752,0.16751976){\color[rgb]{0,0,0}\makebox(0,0)[t]{\lineheight{1.25}\smash{\begin{tabular}[t]{c}$\frF_{w\to z}$\end{tabular}}}}%
    \put(0,0){\includegraphics[width=\unitlength,page=2]{fig199.pdf}}%
  \end{picture}%
\endgroup%

\caption{Moving $w$ to $z$ through the 3-handle region, in the hypercube $\cY_r^{(1)}$. See Remark~\ref{rem:natural=challenge}.}
\label{fig:natural=challenge}
\end{figure}

\begin{rem}\label{rem:natural=challenge}
 The construction of the maps $X$ and $Z$ in $\cY_r^{(1)}$, shown in Figure~\ref{eq:hypercube-naturality-flip-map-C1-cen-2}, is unusual in that we used the fact that $Y$ was an L-space to algebraically determine the existence of two maps. It does not seem possible to construct these maps by compressing a hypercube analogously to our construction of the flip-map hypercube. The challenge is that the map $\frF_{w\to z}$ on the top right of the bottom cube in Figure~\ref{eq:hypercube-naturality-flip-map-C1-cen-2} involves moving $w$ to $z$ along a path which goes through a 3-handle region. See Figure~\ref{fig:natural=challenge}. The natural algebraic relation involving the 3-handle map does not follow from the framework of stabilizations in Section~\ref{sec:stabilization}.
\end{rem}

\subsection{The hypercube \texorpdfstring{$\cK_{\cen;\fry_s}^{(2)}$}{K-cen-y-s-(2)}}

\begin{figure}[ht!]
\[
\cK_{\cen;\fry_s}^{(2),1}=\hspace{-5mm}
\begin{tikzcd}[
column sep={3.5cm,between origins},
row sep=.5cm,labels=description,
execute at end picture={
\foreach \Nombre in  {A,B,C,D}
  {\coordinate (\Nombre) at (\Nombre.center);}
\fill[opacity=0.1] 
  (A) -- (B) -- (C) -- (D) -- cycle;
}]
\CF^-(\gs',\ds,[s])
	\ar[dd, swap,"\sigma_0 F_3 f_{\g'}^{\dt\to \tilde{\dt}'}"]
	\ar[dr,  "\id",swap]
	\ar[rr, pos=.67, "f_{\g';\fry_s;w}^{\dt\to \dt'}"]
	\ar[ddrr,dashed,sloped,near end]
&&[-1.5cm]
|[alias=A]|\uCF^-(\gs',\ds',\hat{w},z,[s+n])
	\ar[dd, "\sigma_0 F_3 f_{\g'}^{\dt'\to \dt''}", pos=.4]
	\ar[dr,"\frF_{w\to z}"]
	\ar[dddr,dashed]
&
\\
&\CF^-(\gs',\ds,[s])
	\ar[rr,crossing over, "f_{\g';\fry_s;z}^{\dt\to \dt'}"]
&&
|[alias=B]|\uCF^-(\gs',\ds',w,\hat{z}, [s])
	\ar[dd, "\sigma_0 F_3 f_{\g'}^{\dt'\to \dt''}"]
	\ar[from=ulll,dashed,crossing over, "H^{\frF}_{w\to z}"]
\\[1.8cm]
\CF^-(\gs,\tilde{\ds}',[s])
	\ar[rr, "f_{\g;\fry_s;w}^{\tilde{\dt}'\to \dt''}"]
	\ar[dr,"\id"]
	\ar[drrr,dashed, "H^{\frF}_{w\to z}"]
&&
|[alias=D]|\uCF^-(\gs,\ds'', \hat{w},z,[s+n])
	\ar[dr, "\frF_{w\to z}"]	
&
\\
&
\CF^-(\gs,\tilde{\ds}', [s])
	\ar[rr, "f_{\g';\fry_s;w}^{\tilde{\dt}'\to \dt''}"]
	\ar[from =uu, crossing over, "\sigma_0 F_3 f_{\g'}^{\dt\to \tilde{\dt}'}",swap]
	&&
|[alias=C]|\uCF^-(\gs,\ds'',w,\hat{z},[s])
	\ar[from=uull,crossing over, dashed]
\end{tikzcd}
\]
\[
\cK_{\cen;\fry_s}^{(2),2}=\hspace{-10mm}
\begin{tikzcd}[
column sep={3.5cm,between origins},
row sep=.5cm,labels=description,
execute at end picture={
\foreach \Nombre in  {A,B,C,D}
  {\coordinate (\Nombre) at (\Nombre.center);}
\fill[opacity=0.1] 
  (A) -- (B) -- (C) -- (D) -- cycle;
}]
|[alias=A]|\CF^-(\gs',\ds',\hat{w},z, [s+n])
	\ar[dd, swap,"\sigma_0 F_3 f_{\g'}^{\dt'\to \dt''}"]
	\ar[dr,  "\frF_{w\to z}",swap]
	\ar[rr, "\id"]
	\ar[ddrr,dashed]
&&[-1.5cm]
\uCF^-(\gs',\ds',\hat{w},z,[s+n])
	\ar[dd, "\sigma_0 F_3 f_{\g'}^{\dt'\to \dt''}", pos=.4]
	\ar[dr,"\id"]
&
\\
&
|[alias=B]|\CF^-(\gs',\ds', w,\hat{z}, [s])
	\ar[rr,crossing over, "\frF_{z\to w}"]
&&
\uCF^-(\gs',\ds',\hat{w},z,[s+n])
	\ar[dd, "\sigma_0 F_3 f_{\g'}^{\dt'\to \dt''}"]
	\ar[from=ulll,dashed,crossing over, "H^{\frF}_{w\to z\to w}"]
\\[1.8cm]
|[alias=D]|\CF^-(\gs,\ds'', \hat{w},z,[s+n])
	\ar[rr, "\id",pos=.2]
	\ar[dr,"\frF_{w\to z}"]
	\ar[drrr,dashed, "H^{\frF}_{w\to z\to w}"]
&&\uCF^-(\gs,\ds'', \hat{w},z,[s+n])
	\ar[dr, "\id"]	
&
\\
&
|[alias=C]|\CF^-(\gs,\ds'', w, \hat{z},[s])
	\ar[rr,"\frF_{z\to w}"]
	\ar[from =uu, crossing over, "\sigma_0 F_3 f_{\g'}^{\dt'\to \dt''}",swap]
	\ar[from=uuul,dashed, crossing over]
	&&
\uCF^-(\gs,\ds'',\hat{w},z,[s+n])
	\ar[from=uull,crossing over, dashed]
\end{tikzcd}
\]
\\
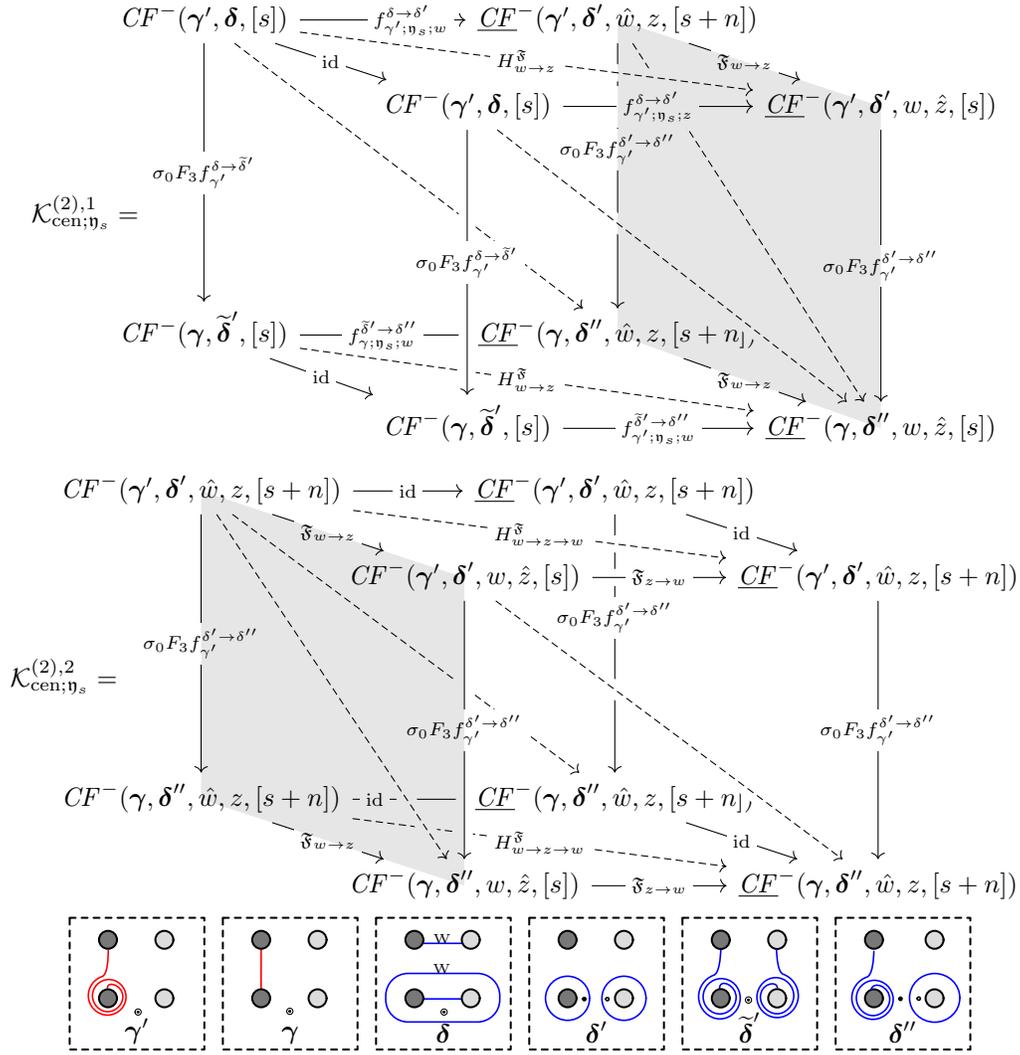
\caption{The hypercubes $\cK_{\cen;\fry_s}^{(2),1}$  and $\cK_{\cen;\fry_s}^{(2),2}$. Both have a length 3 arrow (not shown).}
\label{fig:hypercube-naturality-flip-map-C2-cen}
\end{figure}

We construct $\cK_{\cen;\fry_s}^{(2)}$ by stacking and compressing 4 hypercubes, $\cK_{\cen;\fry_s}^{(2),j}$, for $j\in \{1,2,3,4\}$, similarly to Figure~\ref{eq:C-fry-s-Psi-i-full}.  As in the previous section, we focus on the hypercubes $\cK_{\cen;\fry_s}^{(2),1}$ and $\cK_{\cen;\fry_s}^{ (2),2}$. We illustrate these two cubes schematically in Figure~\ref{fig:hypercube-naturality-flip-map-C2-cen}.

\begin{figure}[h]
\[
\cY_0^{(2)}=\hspace{-5mm}
\begin{tikzcd}[
column sep={3.5cm,between origins},
row sep=.5cm,labels=description,
execute at end picture={
\foreach \Nombre in  {A,B,C,D}
  {\coordinate (\Nombre) at (\Nombre.center);}
\fill[opacity=0.1] 
  (A) -- (B) -- (C) -- (D) -- cycle;
}]
\CF^-(\gs',\ds,[s])
	\ar[dd, "f_{\g'}^{\dt\to \tilde{\dt}'}"]
	\ar[dr,  "\id",swap]
	\ar[rr, "f_{\g';\fry_s;w}^{\dt\to \dt'}", pos=.6]
	\ar[ddrr,dashed,near end]
&&[-1.5cm]
\uCF^-(\gs',\ds',\hat{w},z,[s+n])
	\ar[dd, "f_{\g'}^{\dt'\to \dt''}"]
	\ar[dr,"\frF_{w\to z}"]
	\ar[dddr,dashed]
&
\\
&\CF^-(\gs',\ds,[s])
	\ar[rr,crossing over,"f_{\g';\fry_s;z}^{\dt\to \dt'}"]
&&
\uCF^-(\gs',\ds',w,\hat{z},[s])
	\ar[dd, "f_{\g'}^{\dt'\to \dt''}"]
	\ar[from=ulll,dashed,crossing over, "H^{\frF}_{w\to z}"]
\\[1.8cm]
|[alias=A]|\CF^-(\gs',\tilde{\ds}',[s])
	\ar[rr, "f_{\g';\fry_s;w}^{\tilde{\dt}'\to \dt''}"]
	\ar[dr,"\id"]
	\ar[drrr,dashed, "H^{\frF}_{w\to z}"]
&&|[alias=B]|\uCF^-(\gs',\ds'', \hat{w},z,[s+n])
	\ar[dr, "\frF_{w\to z}"]	
&
\\
&
|[alias=D]|\CF^-(\gs',\tilde{\ds}',[s])
	\ar[rr, "f_{\g';\fry_s;z}^{\tilde{\dt}'\to \dt''}"]
	\ar[from =uu, crossing over, "f_{\g'}^{\dt\to \tilde{\dt}'}"]
	&&
|[alias=C]|\uCF^-(\gs',\ds'',w,\hat{z},[s])
	\ar[from=uull,crossing over, dashed]
\end{tikzcd}
\]
\[
\cY_1^{(2)}=\hspace{-5mm}
\begin{tikzcd}[
column sep={3.5cm,between origins},
row sep=.5cm,labels=description,
execute at end picture={
\foreach \Nombre in  {A,B,C,D}
  {\coordinate (\Nombre) at (\Nombre.center);}
\fill[opacity=0.1] 
  (A) -- (B) -- (C) -- (D) -- cycle;
}]
|[alias=A]|\CF^-(\gs',\tilde{\ds}',[s])
	\ar[dd, swap, "\sigma_0 F_3"]
	\ar[dr,  "\id",swap]
	\ar[rr, "f_{\g';\fry_s;w}^{\tilde{\dt}'\to \dt''}"]
&&[-1.5cm]
|[alias=B]|\uCF^-(\gs',\ds'',\hat{w},z,[s+n])
	\ar[dd, "\sigma_0 F_3"]
	\ar[dr,"\frF_{w\to z}"]
&
\\
&|[alias=D]|\CF^-(\gs',\tilde{\ds}',[s])
	\ar[rr,crossing over, "f_{\g';\fry_s;z}^{\tilde{\dt}'\to \dt''}"]
&&
|[alias=C]|\uCF^-(\gs',\ds'',w,\hat{z},[s])
	\ar[dd, "\sigma_0 F_3"]
	\ar[from=ulll,dashed,crossing over, "H^{\frF}_{w\to z}"]
\\[1.8cm]
\CF^-(\gs,\tilde{\ds}',[s])
	\ar[rr, "f_{\g;\fry_s;w}^{\tilde{\dt}'\to \dt''}"]
	\ar[dr,"\id"]
	\ar[drrr,dashed, "H^{\frF}_{w\to z}"]
&&\uCF^-(\gs,\ds'', \hat{w},z,[s+n])
	\ar[dr, "\frF_{w\to z}"]	
&
\\
&
\CF^-(\gs,\tilde{\ds}',[s])
	\ar[rr, "f_{\g;\fry_s;z}^{\tilde{\dt}'\to \dt''}"]
	\ar[from =uu, crossing over, "\sigma_0 F_3"]
	&&
\uCF^-(\gs,\ds'',w,\hat{z},[s])
\end{tikzcd}
\]
\caption{A hyperbox whose compression is  $\cK_{\cen;\fry_s}^{(2),1}$. The hypercube $\cY_0^{(2)}$ has a length 3 arrow (not shown), whereas $\cY_1^{(2)}$ does not.}
\label{fig:C-fry-s-Psi-i-cen-2}
\end{figure}
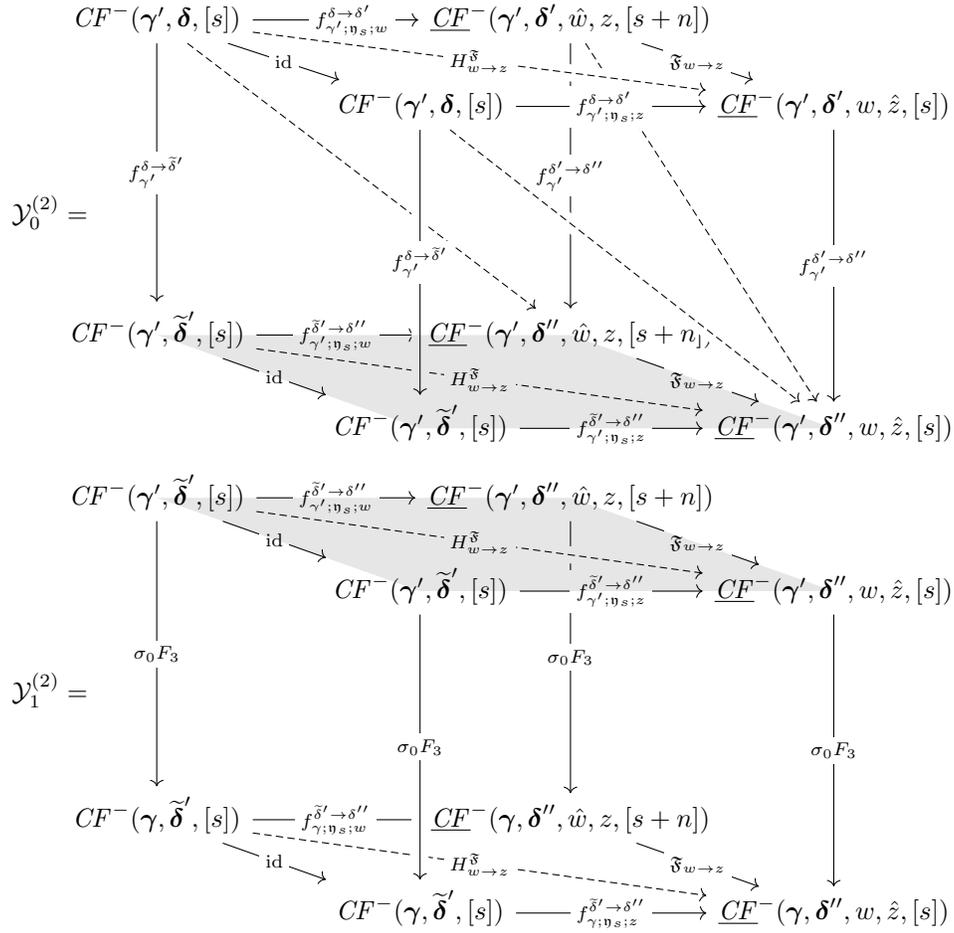

The hypercube $\cK_{\cen;\fry_s}^{(2),1}$ is constructed as the compression of the hyperbox shown in Figure~\ref{fig:C-fry-s-Psi-i-cen-2}. We write $\cY_0^{(2)}$ for the top hypercube, therein, and $\cY_1^{(2)}$ for the bottom.  The hypercube $\cY_0^{(2)}$ is defined as a compression of a hypercube of size $(1,5,1)$, similar to the flip-map hypercube. In the hyperbox for $\cY_0^{(2)}$, there is one hypercube built using $\theta_w$, one using the tautological diffeomorphism map, one hypercube for winding $\ds''$, one cube for unwinding $\ds''$, and a final hypercube involving $\theta_z$. 

The  hypercube $\cY_1^{(2)}$ is defined by stacking and compressing several hypercubes for destabilization, as we now describe. Since the hypercube $\cY_0^{(2)}$ is the compression of hyperbox of size $(1,5,1)$, the bottom face of $\cY_0^{(2)}$ is the compression of a hyperbox of size $(1,5)$.  Each of the 3-dimensional hypercubes comprising $\cY_0^{(2)}$ has along its bottom face only complexes which have a special genus 1 region representing $S^1\times S^2$, $S^3$, or $L(m,1)$. Furthermore, since the canonical path connecting $w$ and $z$ does not cross any of the curves in the special genus 1 region, we may assume that the special genus 1 region is standard for each complex appearing in the construction of the bottom face of $\cY_0^{(2)}$. For each cube used to build $\cY^{(2)}_0$, we build a corresponding hypercube of destabilization in $\cY^{(2)}_1$. The maps along the top face of each of these cubes for $\cY^{(2)}_1$ coincide with the maps appearing along the bottom of the  corresponding hypercube used for $\cY^{(2)}_0$. The maps along the bottom are given by the same set of moves used for $\cY^{(2)}_0$, but applied to the diagram obtained by removing the 3-handle region and adding a lens-space stabilization. The length 1 vertical maps are all $\sigma_0 F_3$. Propositions~\ref{prop:stabilize-triangles-general} and~\ref{prop:stabilize-quadrilateral} imply that the resulting diagram satisfies the hyperbox relations. Finally $\cY_1^{(2)}$ is defined as the compression of this hyperbox of size $(1,5,1)$.

The construction of $\cK_{\cen;\fry_s}^{(2),2}$, shown on the bottom of Figure~\ref{fig:hypercube-naturality-flip-map-C2-cen}, follows from a similar construction. The hypercubes $\cK_{\cen;\fry_s}^{(2),3}$ and $\cK_{\cen;\fry_s}^{(2),4}$ are not substantially different than the corresponding hypercubes of $\cK_{\twoh;\fry_s}^{1, 3}$ and $\cK_{\twoh;\fry_s}^{1,4}$, shown in Figure~\ref{eq:C-fry-s-Psi-i-full}.

\subsection{The hypercube \texorpdfstring{$\cK_{\cen;\fry_s}^{(3)}$}{K-cen-y-s-(3)}}

\begin{figure}[h]
\begin{tikzcd}[
column sep={3.8cm,between origins},
row sep=.5cm,labels=description,
execute at end picture={
\foreach \Nombre in  {A,B,C,D}
  {\coordinate (\Nombre) at (\Nombre.center);}
\fill[opacity=0.1] 
  (A) -- (B) -- (C) -- (D) -- cycle;
}
]
\CF^-(\gs,\tilde{\ds}',[s])
	\ar[dd, swap, "f_{\g}^{\tilde{\dt}'\to \dt_2}"]
	\ar[dr,  "\id",swap]
	\ar[rr, "f_{\g;\fry_s;w}^{\tilde{\dt}'\to \dt''}", pos=.6]
	\ar[ddrr,dashed,near end]
&&[-1.5cm]
\uCF^-(\gs,\ds'',\hat{w},z,[s+n])
	\ar[dd, "f_{\g}^{\dt''\to \dt_3}"]
	\ar[dr,"\frF_{w\to z}"]
	\ar[dddr,dashed]
&
\\
&\CF^-(\gs,\tilde{\ds}',[s])
	\ar[rr,"f_{\g;\fry_s;z}^{\tilde{\dt}'\to \dt''}"]
&&
\uCF^-(\gs,\ds'',w,\hat{z},[s])
	\ar[dd, "f_{\g}^{\dt''\to \dt_3}"]
	\ar[from=ulll,dashed,crossing over, "H^{\frF}_{w\to z}"]
\\[1.8cm]
|[alias=A]|\CF^-(\gs,\ds_2,[s])
	\ar[rr, "f_{\g;\fry_s;w}^{\dt_2\to \dt_3}"]
	\ar[dr,"\id"]
	\ar[drrr,dashed, "H^{\frF}_{w\to z}"]
&&|[alias=B]|\uCF^-(\gs,\ds_3, \hat{w},z,[s+n])
	\ar[dr, "\frF_{w\to z}"]	
&
\\
&
|[alias=D]|\CF^-(\gs,\ds_2,[s])
	\ar[rr, "f_{\g;\fry_s;z}^{\dt_2\to \dt_3}"]
	\ar[from =uu, crossing over,"f_{\g}^{\tilde{\dt}'\to \dt_2}"]
	&&
|[alias=C]|\uCF^-(\gs,\ds_3,w,\hat{z},[s])
	\ar[from=uull,crossing over, dashed]
\end{tikzcd}
\begin{tikzcd}[
column sep={3.8cm,between origins},
row sep=.5cm,labels=description,
execute at end picture={
\foreach \Nombre in  {A,B,C,D}
  {\coordinate (\Nombre) at (\Nombre.center);}
\fill[opacity=0.1] 
  (A) -- (B) -- (C) -- (D) -- cycle;
}
]
|[alias=A]|\CF^-(\gs,\ds_2,[s])
	\ar[dd, "F_3"]
	\ar[dr,  "\id",swap]
	\ar[rr, "f_{\g;\fry_s;w}^{\dt_2\to \dt_3}"]
&&[-1.5cm]
|[alias=B]|\uCF^-(\gs,\ds_3,\hat{w},z,[s+n])
	\ar[dd, "F_3"]
	\ar[dr,"\frF_{w\to z}"]
&
\\
&|[alias=D]|\CF^-(\gs,\ds_2,[s])
	\ar[rr,crossing over, "f_{\g;\fry_s;z}^{\dt_2\to \dt_3}"]
&&
|[alias=C]|\uCF^-(\gs,\ds_3,w,\hat{z},[s+n])
	\ar[dd, "F_3"]
	\ar[from=ulll,dashed,crossing over, "H^{\frF}_{w\to z}"]
\\[1.8cm]
\CF^-(\as\cup \bar{\bs},\Ds_2,[s])
	\ar[rr, "f_{\a\cup \bar{\b};\fry_s;w}^{\Dt_2\to \Dt_3}"]
	\ar[dr,"\id"]
	\ar[drrr,dashed, "H^{\frF}_{w\to z}"]
&&\uCF^-(\as\cup \bar{\bs},\Ds_3, \hat{w},z,[s+n])
	\ar[dr, "\frF_{w\to z}"]	
&
\\
&
\CF^-(\as\cup \bar{\bs},\Ds_2,[s])
	\ar[rr, "f_{\a\cup \bar{\b};\fry_s;z}^{\Dt_2\to \Dt_3}"]
	\ar[from =uu, crossing over,"F_3"]
	&&
\uCF^-(\as\cup \bar{\bs},\Ds_3,w,\hat{z},[s])
\end{tikzcd}
\begingroup%
  \makeatletter%
  \providecommand\color[2][]{%
    \errmessage{(Inkscape) Color is used for the text in Inkscape, but the package 'color.sty' is not loaded}%
    \renewcommand\color[2][]{}%
  }%
  \providecommand\transparent[1]{%
    \errmessage{(Inkscape) Transparency is used (non-zero) for the text in Inkscape, but the package 'transparent.sty' is not loaded}%
    \renewcommand\transparent[1]{}%
  }%
  \providecommand\rotatebox[2]{#2}%
  \newcommand*\fsize{\dimexpr\f@size pt\relax}%
  \newcommand*\lineheight[1]{\fontsize{\fsize}{#1\fsize}\selectfont}%
  \ifx\svgwidth\undefined%
    \setlength{\unitlength}{282.30216095bp}%
    \ifx\svgscale\undefined%
      \relax%
    \else%
      \setlength{\unitlength}{\unitlength * \real{\svgscale}}%
    \fi%
  \else%
    \setlength{\unitlength}{\svgwidth}%
  \fi%
  \global\let\svgwidth\undefined%
  \global\let\svgscale\undefined%
  \makeatother%
  \begin{picture}(1,0.17864335)%
    \lineheight{1}%
    \setlength\tabcolsep{0pt}%
    \put(0,0){\includegraphics[width=\unitlength,page=1]{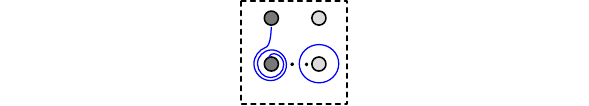}}%
    \put(0.49972577,0.01006762){\color[rgb]{0,0,0}\makebox(0,0)[t]{\lineheight{1.25}\smash{\begin{tabular}[t]{c}$\ds''$\end{tabular}}}}%
    \put(0,0){\includegraphics[width=\unitlength,page=2]{fig198.pdf}}%
    \put(0.29557626,0.01070444){\color[rgb]{0,0,0}\makebox(0,0)[t]{\lineheight{1.25}\smash{\begin{tabular}[t]{c}$\tilde{\ds}'$\end{tabular}}}}%
    \put(0,0){\includegraphics[width=\unitlength,page=3]{fig198.pdf}}%
    \put(0.09121015,0.01634918){\color[rgb]{0,0,0}\makebox(0,0)[t]{\lineheight{1.25}\smash{\begin{tabular}[t]{c}$\gs$\end{tabular}}}}%
    \put(0,0){\includegraphics[width=\unitlength,page=4]{fig198.pdf}}%
    \put(0.90840156,0.01006769){\color[rgb]{0,0,0}\makebox(0,0)[t]{\lineheight{1.25}\smash{\begin{tabular}[t]{c}$\ds_3$\end{tabular}}}}%
    \put(0,0){\includegraphics[width=\unitlength,page=5]{fig198.pdf}}%
    \put(0.70425215,0.01070452){\color[rgb]{0,0,0}\makebox(0,0)[t]{\lineheight{1.25}\smash{\begin{tabular}[t]{c}$\ds_2$\end{tabular}}}}%
    \put(0,0){\includegraphics[width=\unitlength,page=6]{fig198.pdf}}%
  \end{picture}%
\endgroup%

\caption{The two hypercubes used to build $\cK^{(3),1}_{\cen;\fry_s}$. The cubes are stacked along the gray faces.}
\label{eq:hypercube-naturality-flip-map-C3-cen-1}
\end{figure}

Similar to $\cK_{\cen;\fry_s}^{(1)}$ and $\cK_{\cen;\fry_s}^{(2)}$, the hypercube $\cK_{\cen;\fry_s}^{(3)}$ is defined by stacking and compressing four hypercubes $\cK_{\cen;\fry_s}^{(3),j}$, for $j\in \{1,2,3,4\}$. The construction is similar to $\cK_{\cen;\fry_s}^{(2)}$, so we focus on $\cK_{\cen;\fry_s}^{(3),1}$, and leave the other three hypercubes to the reader.  Recall that $\cC_{\cen}^{(3)}$ is defined by stacking $\cC_{\aux,\Pi^{\can}}^{(3)}$ and a hypercube for destabilization $\cC_{\cen,\thrh}^{(3)}$. The hypercube $\cK_{\cen;\fry_s}^{(3),1}$ is the compression of the hyperbox shown in Figure~\ref{eq:hypercube-naturality-flip-map-C3-cen-1}.

The hypercubes in Figure~\ref{eq:hypercube-naturality-flip-map-C3-cen-1} are each constructed as the compression of a hypercube of size $(1,4,1)$, similar to the ones appearing in the construction of the cube $\cK^{(2),1}_{\cen;\fry_s}$. The construction of the hypercube $\cK_{\cen;\fry_s}^{(3),2}$ follows from the same line of reasoning used to construct the hypercube $\cK_{\cen;\fry_s}^{(2),2}$. Hypercubes $\cK_{\cen;\fry_s}^{(3),3}$ and $\cK_{\cen;\fry_s}^{(3),4}$ are constructed similarly.

\subsection{The hypercubes \texorpdfstring{$\cK_{\eta;\fry_s}$}{K-y-s-eta} and \texorpdfstring{$\cK_{\frF;\fry_s}$}{K-y-s-F}}

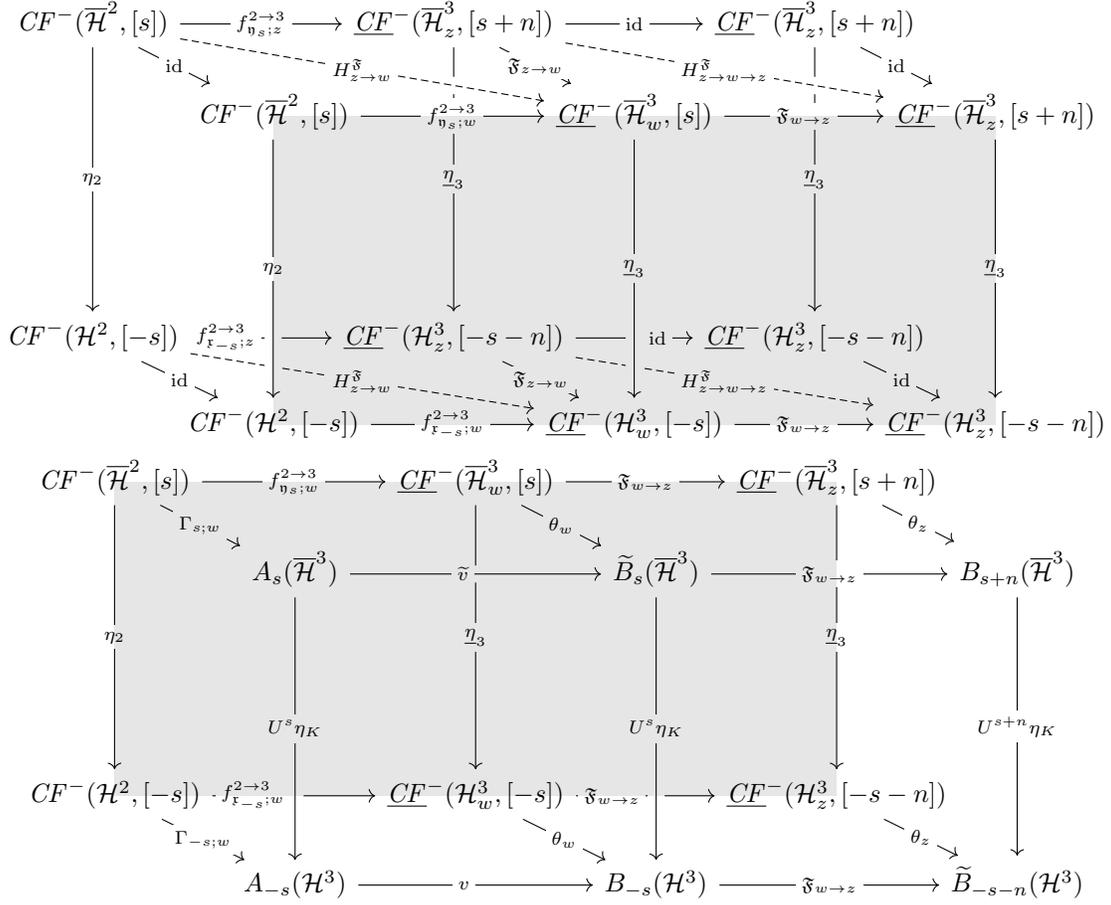
\begin{figure}[h]
\begin{tikzcd}[
column sep={2.4cm,between origins},
row sep=.5cm,labels=description,
execute at end picture={
\foreach \Nombre in  {A,B,C,D}
  {\coordinate (\Nombre) at (\Nombre.center);}
\fill[opacity=0.1] 
  (A) -- (B) -- (C) -- (D) -- cycle;
}
]
\CF^-(\bar{\cH}^2, [s])
	\ar[dd, "\eta_2"]
	\ar[dr,  "\id"]
	\ar[rr, "f^{2\to 3}_{\fry_s;z}"]
&&
\uCF^-(\bar{\cH}^3_{z},[s+n])
	\ar[dd, "\ul{\eta}_3"]
	\ar[dr,"\frF_{z\to w}"]
 \ar[rr, "\id"]
&
&\uCF^-(\bar{\cH}^3_{z},[s+n])
\ar[dd, "\ul{\eta}_3"]
\ar[dr,"\id"]
&\,
\\
&|[alias=A]|
\CF^-(\bar{\cH}^2,[s])
	\ar[rr,crossing over, "f^{2\to 3}_{\fry_s;w}"]
&&
\uCF^-(\bar{\cH}^3_{w},[s])
	\ar[rr, "\frF_{w\to z}"]
	\ar[from=ulll, crossing over, dashed, "H^{\frF}_{z\to w}"]
&& |[alias=B]|\uCF^-(\bar{\cH}^3_{z},[s+n])
	\ar[dd, "\ul{\eta}_3"]
	\ar[from=ulll, crossing over,dashed, "H^{\frF}_{z\to w\to z}"]
\\[1.8cm]
\CF^-(\cH^2, [-s]) 
	\ar[rr, "f^{2\to 3}_{\frx_{-s};z}",near start]
	\ar[dr,"\id"]
	\ar[drrr,dashed,"H^{\frF}_{z\to w}"]
	&
&\uCF^-(\cH^3_{z},[-s-n])
	\ar[dr, "\frF_{z\to w}"]	
	\ar[rr,pos=.7, "\id"]
	\ar[drrr,dashed, "H^{\frF}_{z\to w\to z}"]
&
&\uCF^-(\cH^3_{z},[-s-n])
	\ar[dr, "\id"]
\\
&
|[alias=D]|\CF^-(\cH^2,[-s])
	\ar[rr, "f^{2\to 3}_{\frx_{-s};w}"]
	\ar[from =uu, crossing over,swap, "\eta_2"]
	&&
\uCF^-(\cH_{w}^3,[-s])
	\ar[rr, "\frF_{w\to z}"]
	\ar[from=uu,crossing over, "\ul{\eta}_3"]
&&|[alias=C]| \uCF^-(\cH_{z}^3,[-s-n])
\end{tikzcd}
\begin{tikzcd}[
column sep={2.4cm,between origins},
row sep=.5cm,labels=description,
execute at end picture={
\foreach \Nombre in  {A,B,C,D}
  {\coordinate (\Nombre) at (\Nombre.center);}
\fill[opacity=0.1] 
  (A) -- (B) -- (C) -- (D) -- cycle;
}
]
|[alias=A]|\CF^-(\bar{\cH}^2, [s])
	\ar[dd, "\eta_2"]
	\ar[dr,  "\Gamma_{s;w}"]
	\ar[rr, "f^{2\to 3}_{\fry_{s};w}"]
&&
\uCF^-(\bar{\cH}^3_{w},[s])
	\ar[dd, "\ul{\eta}_3"]
	\ar[dr,"\theta_w"]
 \ar[rr, "\frF_{w\to z}"]
&
&|[alias=B]|\uCF^-(\bar{\cH}^3_{z},[s+n])
	\ar[dd, "\ul{\eta}_3"]
	\ar[dr,"\theta_z"]
&\,
\\
&A_s(\bar{\cH}^3)
	\ar[rr,crossing over,pos=.45, "\vopp"]
&&
\Bopp_{s}(\bar{\cH}^3)
	\ar[rr, "\frF_{w\to z}"]
&& 
B_{s+n}(\bar{\cH}^3)
	\ar[dd, "U^{s+n}\eta_K"]
\\[1.8cm]
|[alias=D]|\CF^-(\cH^2, [-s]) 
	\ar[rr, "f^{2\to 3}_{\frx_{-s};w}",near start]
	\ar[dr,"\Gamma_{-s;w}"]
	&
&\uCF^-(\cH^3_{w},[-s])
	\ar[dr, "\theta_w"]	
	\ar[rr, "\frF_{w\to z}",pos=.26]
&
&|[alias=C]|\uCF^-(\cH^3_{z},[-s-n])
	\ar[dr, "\theta_z"]
\\
&
A_{-s}(\cH^3)
	\ar[rr, pos=.45,"v"]
	\ar[from =uu, crossing over,swap, "U^s \eta_K"]
	&&
B_{-s}(\cH^3)
	\ar[rr, "\frF_{w\to z}"]
	\ar[from=uu,crossing over, "U^{s}\eta_K"]
&& \Bopp_{-s-n}(\cH^3)
\end{tikzcd}
\caption{The hypercube $\cK_{\eta;\fry_s}$ is obtained by compressing the above hypercube. The gray faces indicate how the hyperboxes are stacked.}
\label{eq:C-fry-s-eta}
\end{figure}

We now construct the hypercube $\cK_{\eta;\fry_s}$. We define it to be the compression of the hyperbox of size $(2,2,1)$ shown in Figure~\ref{eq:C-fry-s-eta}.

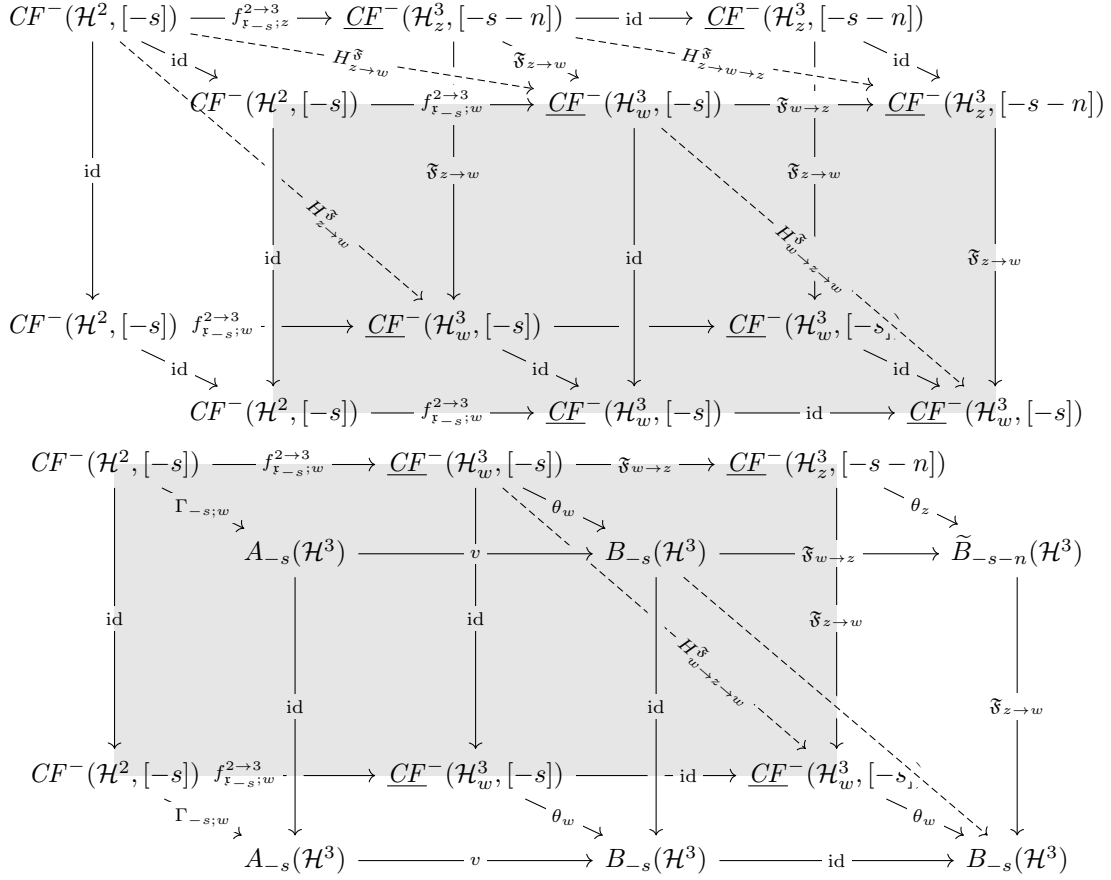
\begin{figure}[h]
\begin{tikzcd}[
column sep={2.4cm,between origins},
row sep=.5cm,labels=description,
execute at end picture={
\foreach \Nombre in  {A,B,C,D}
  {\coordinate (\Nombre) at (\Nombre.center);}
\fill[opacity=0.1] 
  (A) -- (B) -- (C) -- (D) -- cycle;
}
]
\CF^-(\cH^2, [-s])
	\ar[dd, "\id"]
	\ar[dr,  "\id"]
	\ar[rr, "f^{2\to 3}_{\frx_{-s};z}"]
 	\ar[ddrr,dashed,"H^{\frF}_{z\to w}",sloped,pos=.7]
&&
\uCF^-(\cH^3_{z},[-s-n])
	\ar[dd, "\frF_{z\to w}"]
	\ar[dr,"\frF_{z\to w}"]
 \ar[rr, "\id"]
&
&\uCF^-(\cH^3_{z},[-s-n])
	\ar[dd, "\frF_{z\to w}"]
	\ar[dr,"\id"]
&\,
\\
&|[alias=A]|\CF^-(\cH^2,[-s])
	\ar[rr,crossing over, "f^{2\to 3}_{\frx_{-s};w}"]
&&
\uCF^-(\cH^3_{w},[-s])
	\ar[rr, "\frF_{w\to z}"]
	\ar[from=ulll, crossing over, dashed, "H^{\frF}_{z\to w}",sloped]
&& |[alias=B]|\uCF^-(\cH^3_{z},[-s-n])
	\ar[dd, "\frF_{ z\to w}"]
	\ar[from=ulll, crossing over,dashed, "H^{\frF}_{z\to w\to z}",sloped]
\\[1.8cm]
\CF^-(\cH^2, [-s]) 
	\ar[rr, "f^{2\to 3}_{\frx_{-s};w}",pos=.2]
	\ar[dr,"\id"]
	&
&\uCF^-(\cH^3_{w},[-s])
	\ar[dr, "\id"]	
	\ar[rr, "\id"]
&
&\uCF^-(\cH^3_{w},[-s])
	\ar[dr, "\id"]
\\
&
|[alias=D]|\CF^-(\cH^2,[-s])
	\ar[rr, "f^{2\to 3}_{\frx_{-s};w}"]
	\ar[from =uu, crossing over,swap, "\id"]
	&&
\uCF^-(\cH_{w}^3,[-s])
	\ar[rr, "\id"]
	\ar[from=uu,crossing over, "\id"]
&&|[alias=C]| \uCF^-(\cH_{w}^3,[-s])
	\ar[from=uull,dashed, crossing over,sloped, "H^{\frF}_{w\to z\to w}"]
\end{tikzcd}
\begin{tikzcd}[
column sep={2.4cm,between origins},
row sep=.5cm,labels=description,
execute at end picture={
\foreach \Nombre in  {A,B,C,D}
  {\coordinate (\Nombre) at (\Nombre.center);}
\fill[opacity=0.1] 
  (A) -- (B) -- (C) -- (D) -- cycle;
}
]
|[alias=A]|\CF^-(\cH^2, [-s])
	\ar[dd, "\id"]
	\ar[dr,  "\Gamma_{-s;w}"]
	\ar[rr, "f^{2\to 3}_{\frx_{-s};w}"]
&&
\uCF^-(\cH^3_{w},[-s])
	\ar[dd, "\id"]
	\ar[dr,"\theta_w"]
 	\ar[ddrr,dashed, "H^{\frF}_{w\to z\to w}", pos=.7,sloped]
	\ar[rr, "\frF_{w\to z}"]
&
&|[alias=B]|\uCF^-(\cH^3_{z},[-s-n])
	\ar[dd, "\frF_{z\to w}"]
	\ar[dr,"\theta_z"]
&\,
\\
&A_{-s}(\cH^3)
	\ar[rr,crossing over, "v"]
&&
B_{-s}(\cH^3)
	\ar[rr, "\frF_{w\to z}"]
&& \Bopp_{-s-n}(\cH^3)
	\ar[dd, "\frF_{z\to w}"]
\\[1.8cm]
|[alias=D]|\CF^-(\cH^2, [-s]) 
	\ar[rr, "f^{2\to 3}_{\frx_{-s};w}",pos=.2]
	\ar[dr,"\Gamma_{-s;w}"]
	&
&\uCF^-(\cH^3_{w},[-s])
	\ar[dr, "\theta_w"]	
	\ar[rr, "\id",pos=.7]
&
&|[alias=C]|\uCF^-(\cH^3_{w},[-s])
	\ar[dr, "\theta_w"]
\\
&
A_{-s}(\cH^3)
	\ar[rr, "v"]
	\ar[from =uu, crossing over,swap, "\id"]
	&&
B_{-s}(\cH^3)
	\ar[rr, "\id"]
	\ar[from=uu,crossing over, "\id"]
&& B_{-s}(\cH^3)
	\ar[from=uull,dashed, crossing over]
\end{tikzcd}
\caption{The hypercube $\cK_{\frF;\fry_s}$ is obtained by compressing this hyperbox. The gray faces are stacked. The hypercube $\cK_{\frF;\fry_s}^{2}$ has a length 3 map, which is not shown, while the other three hypercubes do not.}
\label{eq:C-fry-s-flip}
\end{figure}

Next, we consider the hypercube $\cK_{\frF;\fry_s}$. It is obtained by compressing the hyperbox described schematically in Figure~\ref{eq:C-fry-s-flip}.

We now explain the maps appearing in Figure~\ref{eq:C-fry-s-flip}. The length 2 maps in the top left cube are constructed in the flip-map hypercube. The map along the top and the map along the back are equal. The hypercube relations are clearly satisfied. The length 2 and 3 maps in the top right cube are constructed in the same manner as those appearing in $\cK_{\frF;\frx_s}^{2}$. See Figure~\ref{eq:C-frx-s-flip-full}. The bottom left and right cubes of Figure~\ref{eq:C-fry-s-flip} are tautological.

\section{Rational surgeries and 0-surgeries}
\label{sec:rational}
In this section, we prove a version of the involutive surgeries formula for Morse surgeries on rationally null-homologous knots, and as an application, for rational surgeries on knots in $S^3$. We begin with some basic topology, then state the involutive mapping cone formula in these cases and describe how the proofs differ from the arguments for null-homologous knots. In Section \ref{subsec:0-surgery}, we describe a version of the mapping cone formula for 0-surgeries.

\subsection{Rationally null-homologous knots}
\label{sec:rationally-null-homologous}

 In this section, we describe a few basic topological notions related to rationally null-homologous knots.

Suppose $Y$ is a rational homology 3-sphere and $K\subset Y$ is a knot, which has order $d>0$ in $H_1(Y;\Z)$. 
A \emph{rational Seifert surface} for $K$ is an integral 2-chain $F$ satisfying 
\[
\d F=d\cdot [K].
\]

We now describe $H_1(Y\setminus K;\Z)$. There is a Meyer-Vietoris exact sequence
\begin{equation}
0\to \Z\xrightarrow{\sigma} H_1(Y\setminus K)\to H_1(Y)\to 0, \label{eq:exact-sequence-homology}
\end{equation}
where $\sigma(n)=n\cdot \mu$ and $\mu$ is a meridian. Using a rational Seifert surface $F$, we obtain a map 
\[
\ell\colon H_1(Y\setminus K)\to \Q
\]
given by $\ell(x)=\# (x\cap F)/d$. Note that $\ell\circ \sigma$ gives the standard inclusion of $\Z$ into $\Q$. We obtain an embedding of groups
\begin{equation}
(\ell,i)\colon H_1(Y\setminus K)\to \Q\oplus H_1(Y), \label{eq:H1Y-K-embedding}
\end{equation}
where $i$ is the inclusion of $H_1(Y\setminus K)$ into $H_1(Y)$. We can explicitly describe the image of $(\ell,i)$. If $h\in H_1(Y)$, let $\hat{h}\in H_1(Y\setminus K)$ be any lift.
We define
\begin{equation}
\hat{A}(h):=[\#(\hat{h}\cap F)/d]\in \Q/\Z. \label{eq:A-of-homology-class}
\end{equation}
The quantity $\hat{A}(h)$ is well defined in $\Q/\Z$, since $\hat{h}$ is well-defined up to addition of an integral multiple of $\mu$, and $\ell(x+\mu)=\ell(x)+1$. The image of $(\ell,i)$ consists of all pairs $(t,h)$ where $t\in \hat{A}(h)$ viewed as a subset of $\Q$.

Next, we discuss surgeries on rationally null-homologous knots. We say a framing is \emph{Morse} if it is induced by a simple closed curve on $\d N(K)$ which intersects the meridian exactly once. For null-homologous knots in an integer homology sphere, there is a canonical identification of the set of Morse framings with $\Z$. For rationally null-homologous knots, the set of Morse framings is canonically identified with $\hat{A}([K])$, viewed as a subset of $\Q$.
 The correspondence is given by
\begin{equation}
\lambda \mapsto \# (\lambda\cap F)/d. \label{eq:Morse-framing-rational}
\end{equation}

Suppose $\lambda$ is a Morse framing. Write $W_{\lambda}(K)\colon Y\to Y_{\lambda}(K)$ for the corresponding 2-handle cobordism, and write $W_{\lambda}'(K)\colon Y_{\lambda}(K)\to Y$ for the 2-handle cobordism with the opposite orientation. Applying Mayer-Vietoris to the attaching region in $Y$, we obtain an exact sequence
\begin{equation}
0\to \Z\to H^2(W_{\lambda}'(K))\to H^2(Y)\to 0, \label{eq:exact-sequence-H^2-Y}
\end{equation}
where $\Z$ acts by the Poincar\'{e} dual of the co-core of the 2-handle.
Similar to~\eqref{eq:H1Y-K-embedding}, we obtain an embedding
\begin{equation}
H^2(W_{\lambda}'(K))\to \Q\oplus H^2(Y), \label{eq:H2W-lambda-isomorphism}
\end{equation}
by combining the restriction from $H^2(W_{\lambda}'(Y))\to H^2(Y)$ with the map $H^2(W_{\lambda}'(K))\to \Q$ given by 
\[
\eta\mapsto \langle \eta, \hat{\Sigma}_d \rangle /d.
\]
 Here, $\hat{\Sigma}_d$ is the surface obtained by taking $d$ copies of the core $\Sigma$ of the 2-handle and capping using a rational Seifert surface.

Since $W_{\lambda}'(K)$ is also obtained by attaching a 2-handle to $Y_{\lambda}(K)$, we analogously obtain an exact sequence
\begin{equation}
0\to \Z\to H^2(W_{\lambda}'(K))\to H^2(Y_{\lambda}(K))\to 0,\label{eq:exact-sequence-W-lambda'-2}
\end{equation}
where $\Z$ now acts by the Poincar\'{e} dual of the core of the 2-handle  which is attached to $K$. Combining~\eqref{eq:H2W-lambda-isomorphism} we obtain the description 
\[
H^2(Y_{\lambda}(K))\iso H^2(W_{\lambda}'(K))/(-\lambda,\PD[K]),
\]
where we view $H^2(W_{\lambda}'(K))$ as being a subgroup $\Q\oplus H^2(Y)$.

\subsection{\texorpdfstring{$\Spin^c$}{Spin-c} structures}
\label{sec:spin^c-structures-1}

Similar to ~\eqref{eq:exact-sequence-H^2-Y}, there is an exact sequence of affine spaces
\begin{equation}
0\to \Z\to \Spin^c(W'_\lambda(K))\to \Spin^c(Y)\to 0. \label{eq:exact-sequence-spin^c-Y}
\end{equation}
The fiber of the map from $\Spin^c(W'_{\lambda}(K))$ to  $\Spin^c(Y)$ is generated by the Poincar\'{e} dual to the co-core of the 2-handle. There is a natural map $A\colon \Spin^c(W'_\lambda(K))\to \Q$, given by
\begin{equation}
A(\frs)=\frac{\langle c_1(\frs),[\hat{\Sigma}_d]\rangle -[\Sigma]\cdot [\hat{\Sigma}_d]}{2d},
\label{eq:G(s)-def}
\end{equation}
where $\Sigma$ denotes the core of the 2-handle, and $\hat{\Sigma}_d$ denotes $d[\Sigma]$, capped off using a rational Seifert surface. We obtain an embedding 
\begin{equation}
\Spin^c(W_{\lambda}'(K))\to \Q\times \Spin^c(Y), \label{eq:map-of-spin^c-structures-cobordism}
\end{equation}
where $\frs\mapsto (A(\frs), \frs|_Y)$.

We write 
\[
\bH(Y,K)\subset \Q\times \Spin^c(Y)
\]
for the image of the above embedding. We can describe $\bH(Y,K)$ more explicitly.  There is a map 
\begin{equation}
\hat{A}\colon \Spin^c(Y)\to \Q/\Z,\label{eq:A-hat-def}
\end{equation}
 obtained by lifting $\frs$ to an element $\fru\in \Spin^c(W_{\lambda}'(K))$, and setting $\hat{A}(\frs)=A(\fru)$. Since adding the Poincar\'{e} dual of the co-core of the 2-handle changes $L$ by $\pm 1$, the coset $\hat{A}(\frs)$ is well-defined. The image of the map in~\eqref{eq:map-of-spin^c-structures-cobordism} consists of all pairs $(t,\frs)$ where $t\in \hat{A}(\frs)$. We write $\fru_{t,\frs}$ for the $\Spin^c$ structure on $W_\lambda'(K)$ corresponding to the pair $(t,\frs)\in \Q\times \Spin^c(Y)$.
 
The set $\bH(Y,K)$ is independent of the framing $\lambda$. In fact, from the construction of the absolute gradings on $\cCFK^-(Y,K)$ (cf. Section~\ref{sec:HFK-rational}), it is immediate that the Alexander grading on $\cCFK^-(Y,K,\frs)$ takes values in in $t_0+\Z$, where $t_0$ is any element of $\Q$ in the fiber over $\frs$ in the exact sequence of \ref{eq:exact-sequence-spin^c-Y}. (A purely topological proof of this fact may be given by tracing through the proof of the well-definedness of the absolute gradings in \cite{OSIntersectionForms}).

Similar to~\eqref{eq:exact-sequence-spin^c-Y}, there is an exact sequence
\[
0\to \Z\to \Spin^c(W_{\lambda}'(K))\to \Spin^c(Y_\lambda(K))\to 0,
\]
where now $\Z$ acts by the Poincar\'{e} dual to the core of the 2-handle attached along $K\subset Y$. In terms of the embedding in ~\eqref{eq:map-of-spin^c-structures-cobordism}, the Poincar\'{e} dual of the core acts on $\Spin^c(Y)$ by $\PD[K]$, and on $\Q$ by $\langle  \PD[\Sigma], \hat{\Sigma}_d\rangle/d=\# (K\cap -F)=-\lambda$. Hence, $\Spin^c(Y_{\lambda}(K))$ is isomorphic to the quotient of $\bH(Y,K)$ by the action of $(-\lambda,\PD[K])$.

The first Chern class has a simple description in this framework. Namely, if $(t,\frs)\in \bH(Y,K)$ and $\fru_{t,\frs}\in \Spin^c(W'_\lambda(K))$ is the corresponding $\Spin^c$ structure, then  we can rearrange~\eqref{eq:G(s)-def} to obtain
\begin{equation}
\frac{\langle c_1(\fru_{t,\frs}), [\hat{\Sigma}_d]\rangle}{d}=2t-\lambda.
\end{equation}
Hence 
\begin{equation}
c_1^2(\fru_{t,\frs})=\frac{(2t-\lambda)^2}{\lambda}.\label{eq:c_1^2-rational-surgeries}
\end{equation}

In terms of the embedding in ~\eqref{eq:map-of-spin^c-structures-cobordism}, the conjugation action on $\bH(Y,K)$ has the following description: 
\begin{equation}
(t,\frs)\mapsto  (-t+\lambda,\bar{\frs}). \label{eq:conjugation-H}
\end{equation}

There are two families of $\Spin^c$ structures on $W_{\lambda}'(K)\colon Y_{\lambda}(K)\to Y$ which we are interested in when $\lambda$ is large.  Each element $\frt\in \Spin^c(Y_\lambda(K))\iso \bH(Y,K)/(-\lambda,\PD[K])$ has a unique representative $(t_0,\frs_0)$, where $-\lambda/2< t_0\le \lambda/2$. Furthermore, any other $\Spin^c$ structure on $W_{\lambda}'(K)$ which restricts to $\frt$ has the form $(t_0+k\lambda, \frs_0-k\PD[K])$, for an integer $k$. We define 
\[
\frx_{\frt}=(t_0,\frs_0), \quad  \text{and} \quad  \fry_{\frt}=(t_0+\lambda, \frs_0-\PD[K]).
\]

The following is a restatement of \cite{OSRationalSurgeries}*{Proposition~4.2} using our present notation:

\begin{lem}\label{lem:len-1-homogeneously-graded-pt1}  Suppose that $K$ is a rationally null-homologous knot in $Y$, and $C$ is a fixed number.
 The following hold:
\begin{enumerate}
\item\label{lem:f-2-rational-pt1}  If $\lambda$ is a sufficiently large Morse framing on $K$, then the $\Spin^c$ structures $\frx_\frt$ and $\fry_{\frt}$ realize the maximal and second maximal squares over all $\Spin^c$ structures on $W_{\lambda}'(K)$ which restrict each $\frt\in \Spin^c(Y_{\lambda}(K))$. Furthermore,
\[
c_1(\frz)^2\le \max \{c_1^2(\frx_{\frt}), c_1^2(\fry_{\frt}) \}-C,
\]
whenever $\frz|_{Y_\lambda(K)}=\frt$ and $\frz\not\in \{ \frx_\frt,\fry_{\frt}\}.$ 
\item \label{lem:f-2-rational-pt2} There is a number $b$ such that if $t_0<-b$ or $b<t_0$ and $\lambda$ is sufficiently large, then
\[
 C\le \max \{c_1^2(\frx_{\frt}), c_1^2(\fry_{\frt}) \}-\min\{c_1^2(\frx_{\frt}),c_1^2(\fry_{\frt})\},
\]
whenever $(t_0,\frs)\in \bH(Y,K)$ is the unique representative of $\frt$ with $-\lambda/2<t_0\le \lambda/2$.
\end{enumerate}
\end{lem}
\begin{proof} Suppose $\frt\in Y_{\lambda}(K)$, and let $(t_0,\frs)$ be the unique extension to $W_{\lambda}'(K)$ with $-\lambda /2<t_0\le \lambda/2$. Any other $\Spin^c$ structure on $W_{\lambda}'$ which restricts to $\frt$ may be written as $(t_0+s\lambda,\frs+s\cdot \PD[K])$, for some $s\in \Z$. Equation~\eqref{eq:c_1^2-rational-surgeries} implies that
\[
c_1^2(t_0+s\lambda, \frs_0+s \PD[K])=-\frac{(2t_0+(2s-1)\lambda)^2}{\lambda}.
\]
Ranging over all real $x$, the function $-(2t_0+x\lambda)^2$ is a quadratic which is maximized at $x=-2t_0/\lambda$, which by assumption lies in the half open interval $(-1,1]$. Ranging only over  odd integers $x$, the quadratic function $-(2t_0+x\lambda)^2$ will take on its two largest values at $x=-1$ and $x=1$, which correspond to $s=0$ and $s=1$. These two values  themselves correspond to $\frx_{\frt}$ and $\fry_{\frt}$, respectively.

For the stated inequality, suppose first that $t_0\le  0$. Then the $\Spin^c$ structures with maximal and third to maximal squares are $(t_0,\frs)$ and $(t_0-\lambda, \frs+\PD[K])$. Using~\eqref{eq:c_1^2-rational-surgeries}, the drop in $c_1^2$ between these two $\Spin^c$ structures is
\begin{equation}
-\frac{(2t_0-\lambda)^2}{\lambda}+\frac{(2t_0-3\lambda)^2}{\lambda}=8(\lambda- t_0)\ge 4 \lambda. \label{eq:q-0<0-inequality}
\end{equation}
Similarly, if $t_0\ge 0$, then the maximal and third to maximal $\Spin^c$ structures are $(t_0+\lambda, \frs-\PD[K])$ and $(t_0+2\lambda, \frs-2\PD[K])$, the difference of whose squares is
\begin{equation}
-\frac{(2t_0+\lambda)^2}{\lambda}+\frac{(2t_0+3\lambda)^2}{\lambda}=8(\lambda+t_0)\ge 4\lambda.\label{eq:q-0>0-inequality}
\end{equation}
Combining~\eqref{eq:q-0<0-inequality} and~\eqref{eq:q-0>0-inequality} easily gives the stated inequality, for sufficiently large $\lambda$.

If $t_0>0$, then the maximal square is achieved by $\fry_{\frt}$, and the second maximal is obtained by $\frx_{\frt}$.

 Hence,
\begin{equation}
\begin{split}
c_1^2(\fry_{\frt})-c_1^2(\frx_{\frt})=&\frac{(2t_0+\lambda)^2}{\lambda}-\frac{(2t_0-\lambda)^2}{\lambda}\\
=&8 t_0.
\end{split}
\label{eq:non-maximal-spinc}
\end{equation}

If $t_0<0$, then the maximal $\Spin^c$ structure is $\frx_{\frt}$, and the second maximal is $\fry_{\frt}$. Arguing as in ~\eqref{eq:non-maximal-spinc} shows that $c_1^2(\frx_{\frt})-c_1^2(\fry_{\frt})= -8t_0$. Combining these, we obtain 
\[
\max\{c_1^2(\frx_{\frt}), c_1^2(\fry_{\frt})\}-\min\{c_1^2(\frx_{\frt}), c_1^2(\fry_{\frt})\}= 8 |t_0|.
\]
Setting $b=C/8$ gives the main claim.
\end{proof}

\subsection{The cobordism \texorpdfstring{$W_{\lambda,\lambda+m}$}{W-lambda-lambda+m}}

\label{sec:cobordism-W-lambda-m}

Suppose $K$ is a knot in a rational homology 3-sphere $Y$, and let $K'$ and $K''$ denote two unknots, so that $K'$ is unlinked from $K$, but $K''$ links both, as in Figure~\ref{fig:193}. Let $L=K\cup K'\cup K''$ and let $\Lambda=(\lambda,m,0)$ denote a Morse framing on $L$, where $m$ is a positive integer. We can view the cobordism $W_{\lambda,\lambda+m}$ as going from $Y_\lambda(K)\# L(m,1)$ to $Y_{\lambda+m}(K)$, and being obtained by surgery on $K''$, after already having performed surgery on $K$ and $K'$.

Let $D'$ and $D''$ be Seifert disks of $K'$ and $K''$, and let $F$ be a rational Seifert surface of $K$. Let $d>0$ be the order of $K$ in $H_1(Y;\Z)$.

We obtain classes $\hat{F}$, $\hat{D}'$, and $\hat{D}''$  in the surgery cobordism $W_{\Lambda}(L)$ by capping  $F$, $D'$, and $D''$ with the core of 2-handles. These form a basis of $H_2(W_{\Lambda}(L);\Z)$. The intersection form is clearly
\[
Q=\begin{pmatrix} 
d^2\lambda&0&d\\
0&m &1\\
d &1  &0
\end{pmatrix}
\]
 The homology of $H_2(W_{\lambda,\lambda+m};\Z)$ is the orthogonal complement (over $\Z$) of the span of the first two basis elements. The orthogonal complement is spanned by the vector
 \[
\hat{\Sigma}=\begin{pmatrix}
m/g\\
d\lambda/g\\
-d\lambda m/g
\end{pmatrix} ,
 \]
 where $g=\gcd(m,d\lambda)$. Clearly,
 \[
 \hat{\Sigma}^2=-\frac{d\lambda m (dm+d\lambda)}{g^2}.
 \]

We now wish to consider gradings on the cobordism map, which will be helpful to prove that the vertical truncation of the main hypercube is homogeneously graded for both positive and negative framings $\lambda$.

\begin{lem}\label{lem:len-1-homogeneously-graded-pt2}
 There is a $C$, such that for all large $m$ and all $\frs\in\Spin^c(Y_{\lambda+m}(K))$, the homogeneous elements of $\CF^\delta(Y_{\lambda+m}(K),\frt)$ have gradings which lies in the interval $[-C,\frac{\lambda+m}{4}+C]$.
\end{lem}
\begin{proof}
Using the large surgeries theorem \cite{OSRationalSurgeries}*{Proposition~4.2}, it follows that there is a $C>0$ such that the homogeneous elements of $\CF^\delta(Y_{\lambda+m}(K),\frt)$ have grading which fall within $C$ of 
\[
-\max\{(c_1(\fru)^2+1)/4:\fru\in \Spin^c(W_{\lambda+m}'(K), \fru|_{Y_{\lambda+m}}=\frt\}.
\] 
The Chern classes with maximal and second to maximal square are $\frx_{\frt}$ and $\fry_{\frt}$, which have squares
\[
c_1(\frx_{\frt})^2=-\frac{(2t_0-\lambda-m)^2}{\lambda+m}\quad \text{and} \quad c_1(\fry_{\frt})^2=-\frac{(2t_0+\lambda-m)^2}{\lambda+m},
\]
where $-(\lambda+m)/2<t_0\le (\lambda+m)/2$. Assuming $m$ is sufficiently large, the maximum of these two quantities always lies in the interval $[-(\lambda+m),0]$, implying the main statement.
\end{proof}

The following is very similar to \cite{OSRationalSurgeries}*{Lemma~6.7} (the only difference is that it is stated for negative surgeries as well as positive):
\begin{lem}\label{lem:len-1-homogeneously-graded-pt3}
 Suppose $\lambda\neq 0$ and $C$ is some fixed number. For all $m$ sufficiently large, each $\Spin^c$ structure on $\d W_{\lambda,\lambda+m}$ has at most one extension $\fru$ to $W_{\lambda,\lambda+m}(K)$ such that $(c_1^2(\fru)+m)/4$ lies in the interval $[-C,(\lambda+m)/4+C]$.
\end{lem}
\begin{proof}
 The proof is naturally broken into cases, depending on the sign of $\lambda$. We consider negative $\lambda$ first. If $\lambda$ is negative, the cobordism $W_{\lambda,\lambda+m}$ is positive definite, so it is sufficient to show that for each $\Spin^c$ structure on $\d W_{\lambda,\lambda+m}(K)$, there is at most one extension $\fru$ with 
 \begin{equation}
c_1^2(\fru)\le\lambda+C. \label{eq:spinc-relation-fru-W-lambda-m}
 \end{equation}
 We can write $c_1(\fru)=\alpha\cdot \PD[\hat{\Sigma}]$ for some $\alpha\in \Q$, and ~\eqref{eq:spinc-relation-fru-W-lambda-m} implies
 \[
|\alpha|\le \sqrt{\frac{C+\lambda}{\hat{\Sigma}^2}}. 
 \]
 For sufficiently large $m$, we may assume that $|\alpha|\le 1/2$ for all $\fru$ satisfying~\eqref{eq:spinc-relation-fru-W-lambda-m}. If $\fru'$ is another $\Spin^c$ structure on $W_{\lambda,\lambda+m}$ which has the same restriction to $\d W_{\lambda,\lambda+m}$, then we may write $\fru'=\fru+j\PD[\hat{\Sigma}]$, for some $j\in \Z$. We compute that
 \[
 \begin{split}
c_1^2(\fru')-c_1^2(\fru)&=4 j\langle c_1(\fru),[\hat{\Sigma}]\rangle +4 j^2 \hat{\Sigma}^2\\
&= 4j^2 \hat{\Sigma}^2(\alpha/j+1).
\end{split}
 \]
 By letting $m$ be sufficiently large, the inequality~\eqref{eq:spinc-relation-fru-W-lambda-m} will be violated for $\fru'$.
 
 The case when $\lambda$ is positive is similar (cf. \cite{OSRationalSurgeries}*{Lemma~6.7}).
\end{proof}

\subsection{Negative surgeries}

The involutive hypercube we have constructed gives a morphism from $\CFI^\delta(Y_{\lambda}(K))$ to $\XI_\lambda^\delta$. For positive $\lambda$, the analysis performed by Ozsv\'{a}th and Szab\'{o} \cite{OSRationalSurgeries}*{Section~6}, together with the obvious analog of Lemma~\ref{lem:gradings}, show  that the map is homogeneously graded,  with respect to a relative $\Z$ valued grading on $\XI_\lambda^\delta$. For negative coefficients, our argument thus far shows only that the hypercube determines an $\iota$-homotopy equivalence, but we have not shown it to be homogeneously graded. Note that Ozsv\'{a}th and Szab\'{o} construct a different 2-dimensional hypercube to prove the statement about gradings in the case of negative surgery coefficients. Instead of building another involutive hypercube to mirror Ozsv\'{a}th and Szab\'{o}'s approach, we opt to to show that the vertical truncation of the main hypercube is
also homogeneously graded for negative surgery coefficients, when $m$ is sufficiently large.

The following is the key new lemma for handling gradings in the main hypercube for negative surgeries:

\begin{lem}\label{lem:homogeneous-grading-lemma} Suppose $C$ is fixed. For sufficiently large $m$, the following holds: if $\frt\in \Spin^c(Y_{\lambda}(K))$ and $\frs\in \Spin^c(Y)$, then for each $t\in \Q$, the set of real numbers
\[
G_{\ws}(t,\frs,\frt):=\left\{ \gr_{\ws}(\fru)\middle\vert \begin{array}{c} \fru\in \Spin^c(W_{\lambda}'(K)\# D(-m,1))\\
\fru|_{Y_{\lambda}(K)}=\frt\\
\fru|_{Y}=\frs\\
\fru|_{L(m,1)}=\frs_{\can}\\
 A(\fru)- t\in m\Z\end{array}\right\}
\]
contains no two elements which have distance less than $C$.
\end{lem}

\begin{proof} Each $\Spin^c$ structure on $W'_\lambda(K)\# D(-m,1)$ decomposes as a connected sum $\frw\# \frz$ where $\frw\in \Spin^c(W'_{\lambda}(K))$ and $\frz\in \Spin^c(D(-m,1))$. With this notation,
\[
A(\frw\# \frz)=\frac{\langle c_1(\frw),\hat{\Sigma}_d\rangle-[\Sigma]\cdot [\hat{\Sigma}_d]}{2d}+\frac{\langle c_1(\frz),S\rangle -S^2}{2},
\]
where $\hat{\Sigma}_d\subset W_{\lambda}'(K)$ is $d$ copies of the core $\Sigma$ of the 2-handle capped off with a rational Seifert surface, and $S$ is the 2-sphere in $D(-m,1)$ of self-intersection $-m$. Suppose that $\fru_0=\frw_0\# \frz_0$ satisfies $A(\frw_0\# \frz_0)\equiv s_0\pmod{m}$.   Any other $\Spin^c$ structure  with the same boundary restrictions may be written as
\[
\fru_{i,j}=(\frw_0+i\PD[\hat{\Sigma}_d])\#(\frz_0+j \PD[S]),
\]
for integers $i$ and $j$. Note that adding $\PD[S]$ changes $A$ by $m$, and adding $\hat{\Sigma}_d$ changes $A$ by $[\hat{\Sigma}_d]^2/d=-d\lambda$. Hence $A(\fru_0)-A(\fru_{i,j})\in m\Z$ if and only if $id\lambda\equiv 0\pmod{m}$. We compute that
\begin{equation}
\begin{split}
\gr_{\ws}(\fru_{i,j})-\gr_{\ws}(\fru_0)&=\frac{c_1^2(\fru_{i,j})-c_1^2(\fru_{0})}{4}\\
 &=i\langle c_1(\frw_0), [\hat{\Sigma}_d]\rangle +i^2 [\hat{\Sigma}_d]^2 +j\langle c_1(\frz_0), [S]\rangle +j^2 S^2\\
 &=i \langle c_1(\frw_0),[F]\rangle -i^2d\lambda  +j k m-j^2 m, 
 \end{split}
 \label{eq:difference-in-w-gradings}
\end{equation}
for some $k\in \Z$. Note that since $id\lambda\equiv 0\pmod{m}$, the integer $i$ is divisible by $m/\gcd(m,d\lambda)$. Hence $G_{\ws}(t,\frs,\frt)$ is a subset of a lattice of the form $\alpha+(m/\gcd(m,d\lambda))\cdot \Z)$, where $\alpha\in \Q$. The main claim follows, by picking $m$ so that $m/d\lambda>C$, since $|m/\gcd(m,d\lambda)|\ge |m/d\lambda|.$
\end{proof}

\subsection{Knot Floer homology and rationally null-homologous knots}
\label{sec:HFK-rational}
In this section, we consider the Floer homology of rationally null-homologous knots. We focus on knots in rational homology 3-spheres, to simplify the presentation. However, our discussion holds more generally if we focus on torsion $\Spin^c$ structures.

Suppose $K$ is a rationally null-homologous knot in a rational homology 3-sphere $Y$, and $\frs\in \Spin^c(Y)$. We pick a Heegaard diagram $(\Sigma,\as,\bs,w,z)$ for $(Y,K)$, and define the full knot Floer complex $\cCFK^-(Y,K,\frs)$ to be the free $\bF[\knotU,\knotV]$-module generated by $\xs\in \bT_{\a}\cap \bT_{\b}$ satisfying $\frs_{w}(\xs)=\frs$.

The module $\cCFK^-(Y,K,\frs)$ admits two gradings, $\gr_{\ws}$ and $\gr_{\zs}$, for which $\knotU$ and $\knotV$ having bigradings $(-2,0)$ and $(0,-2)$. These gradings are obtained from the absolute Maslov gradings on $\CF^-(Y,\frs)$ and $\CF^-(Y,\frs-\PD[K])$, via the canonical isomorphisms
\[
\begin{split}
\cCFK^-(Y,K,\frs)\otimes \bF[\knotU,\knotV]/(\knotV-1)&\iso \CF^-(Y,\frs)\quad \text{and} \\
 \cCFK^-(Y,K,\frs)\otimes \bF[\knotU,\knotV]/(\knotU-1)&\iso \CF^-(Y,\frs-\PD[K]).
\end{split}
\]
The difference in $\Spin^c$ structures on the right hand side arises from the relation 
\begin{equation}
\frs_{w}(\xs)-\frs_{z}(\xs)=\PD[K]. \label{eq:change-of-spinc-structure}
\end{equation}

We take as the definition of the absolute Alexander grading the following:
\begin{equation}
A:=\frac{1}{2}(\gr_{\ws}-\gr_{\zs}).
\label{eq:def-Alexander-grading}
\end{equation}
The equivalence with other formulations is verified in \cite{HeddenLevineSurgery}*{Lemma~2.10}. An alternate approach for null-homologous knots may be found in \cite{ZemAbsoluteGradings}*{Theorem 2.13 (e); Proposition~8.1}).

The Alexander grading takes values in $\Q$, but is relatively $\Z$ valued. If $s\in \Q$, 
we define $\cA_{(t,\frs)}(Y,K)$ to be the subset of $\cCFK^-(Y,K,\frs)$ in Alexander grading $t$, and we set  $\cB_{(t,\frs)}(Y,K)$ to be the subset of $\knotV^{-1}\cdot\cCFK^-(Y,K,\frs)$ in Alexander grading $t$.

The construction of an involution on knot Floer homology described in \cite{HMInvolutive} extends to rationally null-homologous knots. However, there is a new subtlety due to~\eqref{eq:change-of-spinc-structure}, in that the involution has the form
\[
\iota_K\colon \cCFK^-(Y,K,\frs)\to \cCFK^-(Y,K,\frs-\PD[K]).
\]

To streamline the presentation, we find it helpful to consider the grading changes associated to link cobordisms, similar to the presentation of \cite{ZemAbsoluteGradings}. The Alexander grading formulas therein are stated for null-homologous links. We now describe how to quickly derive the Alexander grading changes for cobordisms between rationally null-homologous knots from \cite{ZemAbsoluteGradings}. More generally, one could easily adapt the line of reasoning in \cite{ZemAbsoluteGradings} to handle cobordisms between rationally null-homologous \emph{links}, but we will not do so here. 

Suppose $(W,\Sigma)\colon (Y_1,K_1)\to (Y_2,K_2)$ is an oriented, genus 0 knot cobordism between two rationally null-homologous knots. Further, suppose $\frs\in \Spin^c(W)$ has torsion restriction to $Y_1$ and $Y_2$.
We consider the link cobordism map $F_{W,\Sigma,\frs}$ for $(W,\Sigma)$. (The cobordism maps of \cite{ZemCFLTQFT} also require a decoration on $\Sigma$; here we assume it consists of two parallel arcs running from $K_1$ to $K_2$). According to \cite{ZemAbsoluteGradings}*{Theorem~1.4}, the $\gr_{\ws}$ and $\gr_{\zs}$ grading changes are given by
\begin{equation}
\begin{split}
\gr_{\ws}(F_{W,\Sigma,\frs}(\xs))-\gr_{\ws}(\xs)&=\frac{c_1^2(\frs)-2\chi(W)-3\sigma(W)}{4},\\
\gr_{\zs}(F_{W,\Sigma,\frs}(\xs))-\gr_{\zs}(\xs)&=\frac{c_1^2(\frs-\PD[\Sigma])-2\chi(W)-3\sigma(W)}{4}
\end{split}
\end{equation}
Here, we are viewing $[\Sigma]$ as a class in $H_2(W,\d W)$.  From~\eqref{eq:def-Alexander-grading}, we obtain
\begin{equation}
A(F_{W,\Sigma,\frs}(\xs))-A(\xs)=\frac{c_1^2(\frs)-c_1^2(\frs-\PD[\Sigma])}{2}.\label{eq:grading-change-1}
\end{equation}
If $N$ is an integer such that $N\cdot [K_1]=0$ and $N\cdot [K_2]=0$, as elements of $H_1(Y_1)$ and $H_1(Y_2)$, then it is easy to rearrange ~\eqref{eq:grading-change-1} to obtain
\begin{equation}
A(F_{W,\Sigma,\frs}(\xs))-A(\xs)=\frac{\langle c_1(\frs),[\hat{\Sigma}_N]\rangle -[\Sigma]\cdot [\hat{\Sigma}_N]}{2N}, \label{eq:Alexander-gradings-more-useful}
\end{equation}
where $\hat{\Sigma}_N$ denotes a class in $H_2(W)$, which is a preimage of $N\cdot [\Sigma]$ under the map $H_2(W)\to H_2(W,\d W)$.

 Note that given the construction of the absolute gradings in \cite{OSIntersectionForms}, as well as our definition in~\eqref{eq:def-Alexander-grading}, we see that that $\hat{A}(\frs)$ is exactly the subset of $\Q$ that the Alexander grading on $\cCFK^-(Y,K,\frs)$ takes values in.

\subsection{Ozsv\'{a}th and Szab\'{o}'s mapping cone formula for rationally null-homologous knots}
\label{sec:original-mapping-cone-rational}

We now restate Ozsv\'{a}th and Szab\'{o}'s mapping cone formula using the notation from Sections~\ref{sec:rationally-null-homologous} and \ref{sec:spin^c-structures-1}. We set
\[
\frA(Y,K)=\prod_{(t,\frs)\in \bH(Y,K)} \bm{\cA}_{(t,\frs)}(Y,K)\quad \text{and} \quad \frB(Y,K)=\prod_{(t,\frs)\in \bH(Y,K)} \bm{\cB}_{(t,\frs)}(Y,K),
\]
where $\bm{\cA}_{(t,\frs)}(Y,K)$ denotes $\cA_{(t,\frs)}(Y,K)\otimes_{\bF[U]} \bF\llsquare U\rrsquare$, and similarly for $\bm{\cB}_{(t,\frs)}(Y,K)$.

The differential $D_{\lambda}$ is the sum of two maps, $v$ and $h$. The map  $v$ restricts to the inclusion of $\bm{\cA}_{(t,\frs)}(Y,K)$ into $\bm{\cB}_{(t,\frs)}(Y,K)$. Similarly,
\[
h_{(t,\frs)}=\frF_{z\to w}\circ \vopp_{(t,\frs)}, 
\]
where $\vopp_{(t,\frs)}$ is the inclusion of $\bm{\cA}_t(Y,K,\frs)$ into $\tilde{\bm{\cB}}_{t}(Y,K,\frs)$, and 
\[
\frF_{z\to w} \colon \tilde{\bm{\cB}}_{t}(Y,K,\frs)\to \bm{\cB}_{t+\lambda}(Y,K,\frs-\PD[K]),
\]
is the the canonical homotopy equivalence, as follows. Write $w$ and $z$ for the two basepoints of $K$. We consider the canonical isomorphisms
\[
\begin{split}
\theta_z&\colon \tilde{\bm{\cB}}_{(t,\frs)}(Y,K,w,z)\to \bCF^-(Y,z,\frs-\PD[K])\otimes T^t \quad \text{and}
\\ \theta_w&\colon \bm{\cB}_{(t+\lambda,\frs-\PD[K])}(Y,K,w,z)\to \bCF^-(Y,w,\frs-\PD[K])\otimes T^{t+\lambda}.
\end{split}
\]
The change in $\Spin^c$ structure appearing in the definition of $\theta_z$ is due to the relation
\[
\frs_w(\xs)-\frs_z(\xs)=\PD[K].
\]
We let $\phi$ be a diffeomorphism which moves $z$ to $w$ along a subarc of $K$.
We define
\[
\frF_{z\to w}:=\theta_w^{-1}\circ (\phi_*\otimes T^{\lambda}) \circ \theta_z.
\]

\begin{rem}
 We have defined $h_{(t,\frs)}$ to take $\bm{\cA}_{(t,\frs)}(Y,K)$ to $\bm{\cB}_{(t+\lambda,\frs-\PD[K])}(Y,K)$. The sign of $\PD[K]$ in our formula differs from Ozsv\'{a}th and Szab\'{o}'s \cite{OSRationalSurgeries}. This is due to inconsistencies in the literature regarding the string orientation of a knot associated to a doubly pointed Heegaard diagram. A discussion may be found in \cite{HeddenLevineSurgery}*{Section~2.2}. Note that our convention on string orientation is opposite to the one advocated for therein. Our reason is that when considering knot cobordisms, it is natural to define the Alexander grading using a Seifert surface whose oriented boundary is $-K$, so that one can cap off the knot cobordism to obtain a class $H_2(W;\Z)$. In \cite{HeddenLevineSurgery}, the authors define the Alexander grading using a Seifert surface whose boundary is $K$.
\end{rem}

Ozsv\'{a}th and Szab\'{o}'s mapping cone formula can be restated in our notation as saying that 
\[
\bCF^-(Y_{\lambda}(K))\simeq  \Cone(D_{\lambda}\colon \frA\to \frB).
\]
 We now sketch the proof of Ozsv\'{a}th and Szab\'{o}'s mapping cone formula, stated in terms of our notation. Understanding the proof will allow us to efficiently understand the proof of the involutive version.

We consider the Floer complexes $\bCF^-(Y_\lambda(K))$, $\bCF^-(Y_{\lambda+m}(K))$ and $\buCF^-(Y)$, as before. Here the twisted complex $\buCF^-(Y)$ has coefficients in $\Z[\Z/m]$, exactly as in the proof for null-homologous knots. The proof that $\bCF^-(Y_{\lambda}(K))$ is homotopy equivalent to the mapping cone of 
\[
f_2\colon \bCF^-(Y_{\lambda+m}(K))\to \buCF^-(Y)
\]
 is no different than the proof for null-homologous knots. Taking vertical truncations presents no additional complications.

  One now claims that there is a natural grading on $\Cone(f_2^\delta)$ with respect to which the homotopy equivalence is homogeneously graded. Ozsv\'{a}th and Szab\'{o} proved that for positive $\lambda$, the map from $\CF^\delta(Y_{\lambda}(K))$ to $\Cone(f_2^\delta)$ is homogeneously graded, and for negative surgery coefficients, the map from $\Cone(f_2^{\delta})$ to $\CF^\delta(Y_{\lambda}(K))$ is homogenously graded. When we consider the involutive version, it is helpful also to know that the map from $\CF^{\delta}(Y_{\lambda}(K))$ to $\Cone(f_2^\delta)$ is also homogeneously graded for negative $\lambda$, whenever $m$ is sufficiently large, compared to $\delta$. To this end,
 Lemma~\ref{lem:len-1-homogeneously-graded-pt1}, ~\ref{lem:len-1-homogeneously-graded-pt2} and~\ref{lem:len-1-homogeneously-graded-pt3} imply that the length 1 maps are homogeneously graded, and  Lemma~\ref{lem:homogeneous-grading-lemma}  implies that the length 2 map is homogeneously graded, when $m$ is sufficiently large compared with $\delta$. We note that the argument for the length 1 maps is not substantially different than Ozsv\'{a}th and Szab\'{o}'s argument, though the argument for the length 2 map is new to our present work.

Next, one builds a 2-dimensional hypercube which relates $\Cone(f_2^{\delta})$ with a horizontal truncation of the mapping cone complex. To this end, we first note that $\theta_w$ should be viewed as  an injection
\[
\theta_{w}\colon \buCF^-(Y)\to \bCF^-(Y)\otimes_{\bF} \bF[\Q/m\Z].
\]
The image of $\theta_w$ consists of the $\bF\llsquare U\rrsquare$-span of elements of the form $\xs\otimes T^t$ where $t\in \hat{A}(\frs_{w}(\xs))$, where $\hat{A}$ is defined in~\eqref{eq:A-hat-def}.

The hypercube one builds has the following form:
\begin{equation}
\begin{tikzcd}[row sep=1.5cm,column sep=1.5cm, labels=description]
\CF^\delta(Y_{\lambda+m}(K))
	\arrow[d, "\Gamma^\delta"]
	\arrow[dr,dashed, "L^\delta"]
	\arrow[r, "f_2^\delta"]&
\uCF^\delta(Y)
	\arrow[d, "\kappa^\delta"]\\
\frA^\delta\langle b\rangle 
	\arrow[r, "v^\delta+h^\delta"]& 
	\frB^\delta\langle b\rangle
\end{tikzcd}
\label{eq:rational-mapping-cone-hypercube}
\end{equation}
We now explain the expressions in~\eqref{eq:rational-mapping-cone-hypercube}. First, the complex $\frA\langle b\rangle$ is the direct sum over all $\cA_{(t,\frs)}(Y,K)$ where $-b\le t\le b$, and $\frB\langle b\rangle$ denotes the direct sum over all $\cB_{(t,\frs)}$ where $-b+\lambda\le t\le b$. Similar to~\eqref{eq:b<m+n/2}, we assume
\[
0<b<(\lambda+m)/2.
\]
Next, we describe $\kappa$. There is a canonical map from the image of $\theta_w$ to $\frB^\delta\langle b\rangle$, which sends $\CF^-(Y)\otimes T^\alpha$ to $\frB_{[\alpha]_m}(Y,K)$, where $[\cdot]_m\colon \Q/m\Z\to \Q$ is the map which sends an element $t$ to its unique representative $t_0$ in $\Q$ which satisfies $-m/2\le t_0< m/2$. We write 
\[
\kappa \colon \buCF^-(Y)\to \frB^-(Y,K)\langle b\rangle,
\]
for the composition of $\theta_w$ with this map, followed by projection onto $\frB^-\langle b\rangle$. Concretely
\begin{equation}
\kappa(T^j U^i \cdot \xs)=
\begin{cases}
 \knotU^i \knotV^{i+[j+A_{w,z}(\xs)]_m-A_{w,z}(\xs)} \cdot \xs,&\text{ if } - b+\lambda\le [j+A_{w,z}(\xs)]_m\le b\\
0& \text{ otherwise}.
\end{cases}
\label{eq:kappa}
\end{equation}

In~\eqref{eq:rational-mapping-cone-hypercube}, the map $\Gamma^\delta$ counts holomorphic triangles representing the $\Spin^c$ structures of the form $\frx_{t,\frs}$, where $-b\le t\le b.$ Note that, essentially by definition, the triangle-counting map, restricted to the $\Spin^c$ structure $\frx_{t,\frs}$ induces Alexander grading change $t$ (viewed as a knot cobordism map from the unknot in $Y_{\lambda+m}(K)$ to $K$ in $Y$), and hence has image in $\cA_{(t,\frs)}^\delta$.

By Lemma~\ref{lem:len-1-homogeneously-graded-pt1}, the only $\Spin^c$ structures which contribute to $f_2^\delta$ for sufficiently large $m$ are of the form $\frx_{t,\frs}$ or $\fry_{t,\frs}$, for $-(\lambda+m)/2\le t\le (\lambda+m)/2$. Furthermore, $\theta_{w}\circ f_{2;\frx_{t,\frs}}$ has image in $\CF^-(Y,\frs)\otimes T^{t}$, since $A(\frx_{t,\frs})=t$, and $\theta_{w}\circ f_{2;\fry_{t,\frs}}$ has image in $\CF^-(Y,\frs-\PD[K])\otimes T^{t+\lambda}$, since $A(\fry_{t,\frs})=t+\lambda+m$. The map $L^\delta$ appearing in~\eqref{eq:rational-mapping-cone-hypercube} is constructed similarly to the analogous map in Figure~\ref{fig:relate-3-manifold-knot-Floer-step-1}. 

Since $b<(\lambda+m)/2$, the truncations $\frA^\delta\langle b\rangle$ and $\frB^\delta \langle b\rangle$ are not themselves homotopy equivalent to $\CF^\delta(Y_{\lambda+m}(K))$ and $\uCF^\delta(Y)$, since $\Gamma^\delta$ projects $\Spin^c$ structures outside of the range $-b\le t\le b$ to 0, and a similar statement holds for $\kappa^\delta$. Nonetheless, it follows from part~\eqref{lem:f-2-rational-pt2} of Lemma~\ref{lem:len-1-homogeneously-graded-pt1} that the map determined by $\Gamma^\delta$, $L^\delta$ and $\kappa_w^\delta$ is a homotopy equivalence. Compare the proof of Theorem~\ref{thm:mapping-cone}. Letting $m$ and $\delta$ become sufficiently large, we obtain, in the same manner as mapping cone for null-homologous knots, a homotopy equivalence between $\Cone(f_2)$ and the infinite mapping cone $\Cone(D_{\lambda}\colon \frA\to \frB)$. Compare the proof of Theorem~\ref{thm:mapping-cone}.

Finally, we make one comment on gradings. Ozsv\'{a}th and Szab\'{o} proved that the map from $\CF^\delta(Y_\lambda(K))$ to $\Cone(f_2^\delta)$ is homogeneously graded when $\lambda$ is positive. When $\lambda$ is negative, they instead showed that the homotopy equivalence, which goes from $\Cone(f_2^\delta)$ to $\CF^\delta(Y_{\lambda}(K))$, is homogeneously graded. We note it follows from Lemma~\ref{lem:homogeneous-grading-lemma} that the map from $\CF^\delta(Y_{\lambda}(K))$ to $\Cone(f_2^\delta)$ is homogeneously graded, even when $\lambda$ is negative, provided $m$ is sufficiently large.

\subsection{The involutive mapping cone formula for rationally null-homologous knots}

Given our sketch of Ozsv\'{a}th and Szab\'{o}'s mapping cone formula for rationally null-homologous knots in Section~\ref{sec:original-mapping-cone-rational}, it is clear that our proof of Theorem~\ref{thm:mapping-cone} extends without any complication to prove the following:

\begin{thm}\label{thm:rationally-null-mapping-cone} Suppose $\lambda$ is a Morse framing on a rationally null-homologous knot $K$ in an rational homology 3-sphere $Y$ (we view $\lambda$ as an element of $\Q$), such that $Y_{\lambda}(K)$ is also a rational homology 3-sphere. The involutive Heegaard Floer homology $\bCFI^-(Y_\lambda(K))$ is homotopy equivalent over the ring $\bF\llsquare U \rrsquare[Q]/Q^2$ to a complex $\bF[Q]/(Q^2)\otimes_\bF (\frA\oplus \frB)$ with differential specified by some $\bF\llsquare U \rrsquare$-linear maps as below:
\[
\bX \bI_\lambda= 
 \begin{tikzcd}[column sep=2cm, row sep=2cm, labels=description]
\frA
	\ar[d, "Q\cdot (\id +\iota_{\bA})",swap]
	\ar[dr, dashed, "Q\cdot H"]
	\ar[r, "D_\lambda"]
& \frB
	\ar[d, "Q\cdot (\id+\iota_{\bB})"]\\
Q\cdot \frA
	\ar[r, "D_\lambda"]
& Q\cdot \frB.
\end{tikzcd}
\]
 In fact, the  $\iota$-complex $(\bCF^-(Y_\lambda(K)), \iota)$ is $\iota$-homotopy equivalent to the complex 
 \[
 \bX_\lambda:=\Cone\left(D_\lambda\colon \frA\to \frB\right),
 \]
  with involution $\iota_{\bX}:=\iota_{\bA}+H+\iota_{\bB}$. Furthermore, the following hold:
\begin{enumerate}
\item \label{cone:1-rational} $D_\lambda=v+h$, where $v$ sends $\bm{\cA}_{(t,\frs)}$ to $\bm{\cB}_{(t,\frs)}$, and $h$ sends $\bm{\cA}_{(t,\frs)}$ to $\bm{\cB}_{(t+\lambda,\frs-\PD[K])}$. 
\item\label{cone:2-rational} $\iota_{\bA}$ sends $\bm{\cA}_{(t,\frs)}$ to $\bm{\cA}_{(-t,\bar{\frs}+\PD[K])}$. Furthermore, $\iota_{\bA}$ coincides with the knot involution.
\item \label{cone:3-rational} $\iota_{\bB}$ sends $\bm{\cB}_{(t,\frs)}$ to $\bm{\cB}_{(-t+\lambda,\bar{\frs})}$, and coincides with the map induced by $\iota_{\bA}$, composed with $\frF_{z\to w}$.
\item\label{cone:4-rational} $H=k+j$, where $k$ maps $\bm{\cA}_{(t,\frs)}$ to $\bm{\cB}_{(-t+\lambda,\bar{\frs})}$ and $j$ maps $\bm{\cA}_{(t,\frs)}$ to $\bm{\cB}_{(-t,\bar{\frs}+\PD[K])}$.
\end{enumerate}
\end{thm}

Given an knot $K$ in a rational homology 3-sphere $Y$, which is also an L-space, the construction of Proposition~\ref{prop:canonical-cone} carries over without substantial change to give an algebraic involutive mapping cone
\[
\XI^{\alg}_{\lambda}(\cCFK^-(Y,K),\iota_K),
\]
whenever $\lambda$ is a non-zero Morse framing on $K$, which is well-defined up to $\iota$-homotopy equivalence. The $\iota$-complex $\XI^{\alg}_{\lambda}$ is characterized by the property that $h$ and $H$ factor through the inclusion map of $\cA_{(t,\frs)}(Y,K)$ into $\tilde{\cB}_{(t,\frs)}(Y,K)$. The proof of Theorem~\ref{thm:strong-cone-formula-proven} adapts with only notational changes to give the following:

\begin{thm}
\label{thm:rational-surgeries-in-L-spaces}
 Suppose $K$ is a rationally null-homologous knot in a rational homology 3-sphere  $Y$, with non-zero Morse framing. If $Y$ is an L-space, then $\XI_{\lambda}^{\alg}(\CFK^\infty(Y,K),\iota_K)$ is $\iota$-homotopy equivalent to $\bCFI^-(Y_{\lambda}(K))$. 
\end{thm}

\subsection{Rational surgeries}
\label{subsec:rational-surgeries}
We now apply Theorem~\ref{thm:rationally-null-mapping-cone} to compute rational surgeries. By using the involutive connected sum formula for knots, the proof follows from the same line of reasoning as \cite{OSRationalSurgeries}*{Theorem~1.1}. We now describe the proof, using our present notation.

\begin{figure}[ht!]
	\centering
\begingroup%
  \makeatletter%
  \providecommand\color[2][]{%
    \errmessage{(Inkscape) Color is used for the text in Inkscape, but the package 'color.sty' is not loaded}%
    \renewcommand\color[2][]{}%
  }%
  \providecommand\transparent[1]{%
    \errmessage{(Inkscape) Transparency is used (non-zero) for the text in Inkscape, but the package 'transparent.sty' is not loaded}%
    \renewcommand\transparent[1]{}%
  }%
  \providecommand\rotatebox[2]{#2}%
  \newcommand*\fsize{\dimexpr\f@size pt\relax}%
  \newcommand*\lineheight[1]{\fontsize{\fsize}{#1\fsize}\selectfont}%
  \ifx\svgwidth\undefined%
    \setlength{\unitlength}{248.97595316bp}%
    \ifx\svgscale\undefined%
      \relax%
    \else%
      \setlength{\unitlength}{\unitlength * \real{\svgscale}}%
    \fi%
  \else%
    \setlength{\unitlength}{\svgwidth}%
  \fi%
  \global\let\svgwidth\undefined%
  \global\let\svgscale\undefined%
  \makeatother%
  \begin{picture}(1,0.41201349)%
    \lineheight{1}%
    \setlength\tabcolsep{0pt}%
    \put(0,0){\includegraphics[width=\unitlength,page=1]{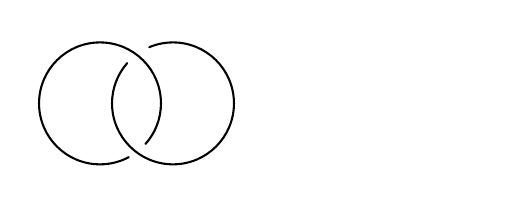}}%
    \put(0.08670666,0.2876456){\color[rgb]{0,0,0}\makebox(0,0)[rt]{\lineheight{1.25}\smash{\begin{tabular}[t]{r}$O_{q/r}$\end{tabular}}}}%
    \put(0.44306742,0.2876456){\color[rgb]{0,0,0}\makebox(0,0)[lt]{\lineheight{1.25}\smash{\begin{tabular}[t]{l}$-q/r$\end{tabular}}}}%
    \put(0,0){\includegraphics[width=\unitlength,page=2]{fig-r1.pdf}}%
    \put(0.88296481,0.21574309){\color[rgb]{1,0,0}\makebox(0,0)[lt]{\lineheight{1.25}\smash{\begin{tabular}[t]{l}$\a$\end{tabular}}}}%
    \put(0.8324006,0.30369263){\color[rgb]{0,0,1}\makebox(0,0)[lt]{\lineheight{1.25}\smash{\begin{tabular}[t]{l}$\b$\end{tabular}}}}%
    \put(0,0){\includegraphics[width=\unitlength,page=3]{fig-r1.pdf}}%
  \end{picture}%
\endgroup%

	\caption{Left: The knot $O_{q/r}$ inside of $-L(q,r)$. Right: A Heegaard diagram for $O_{5/2}$. The knot $O_{5/2}$ is shown (dotted), with the orientation from the Heegaard diagram.}\label{fig:r1}
\end{figure}

Suppose that $\lambda=p/q$ is a non-zero rational framing on $K\subset S^3$, such that $\gcd(p,q)=1$. Let $n$ denote $\lfloor p/q \rfloor$. Define $r/q=p/q-n$. We will take the connected sum of $K$ with the knot $O_{q/r}\subset -L(q,r)$. See Figure~\ref{fig:r1}. Ozsv\'{a}th and Szab\'{o} call $O_{q/r}$ a \emph{U-knot}, since $\widehat{\HFK}$ has rank 1 in each $\Spin^c$ structure.

 We write $x_0,\dots, x_{q-1}$ for the intersection points, ordered left to right (for some choice of starting points). It is easy to check that $\epsilon(x_i,x_{i+r})=\PD[O_{q/r}]$, where $O_{q/r}$ is oriented as in Figure~\ref{fig:r1}, and $\epsilon$ is defined in \cite{OSDisks}*{Section~2.4}. According to \cite{OSDisks}*{Lemma~2.19},
 $\frs_w(x_{i+r})-\frs_w(x_i)=\epsilon(x_i,x_{i+r})=\PD[O_{q/r}]$. On the other hand $\frs_{w}(x_i)-\frs_z(x_i)=\PD[O_{q/r}]$, so $\gr_{z}(x_i)=\gr_w(x_{i-r})$.

We have
\begin{equation}
A(x_i)=\frac{1}{2}(\gr_w(x_i)-\gr_z(x_i))=\frac{1}{2}(\gr_w(x_i)-\gr_w(x_{i-r})). \label{eq:difference-of-spinc-structures}
\end{equation}
According to \cite{LeeLipshitz}*{Proposition~5.3} (cf. \cite{OSIntersectionForms}*{Proposition~4.8}), there is an identification of $\Spin^c(-L(q,r))\iso \Z/q$, such that if $r\le i\le r+q-1$, then
\begin{equation}
\gr_w(x_i)-\gr_w(x_{i-r})=\frac{1}{q}(q-1-2i), \label{eq:LeeLipshitz}
\end{equation}
(with subscripts taken modulo $q$). Combining~\eqref{eq:difference-of-spinc-structures} and~\eqref{eq:LeeLipshitz}, it follows that 
\[
A(x_i)=\frac{1}{2q}(q-1-2i),
\]
for $0\le i<q$.

On the Heegaard diagram for $O_{q/r}\subset- L(q,r)$, the $\beta$ curve satisfies $[\beta]=q\cdot O_{q/r}+r\cdot [\alpha]$. This relation gives a rational Seifert surface, by taking a 2-chain $F$ which satisfies $\d F=q\cdot O_{q/r}+r\cdot\alpha-\beta$, and  capping with compressing disks for $\alpha$ and $\beta$. Furthermore, $O_{q/r}$, pushed into the alpha handlebody, has intersection number $r$ with this 2-chain, and hence $O_{q/r}$ has rational self-linking number $r/q$.  

The framing $p/q$ on $K$ corresponds to a Morse framing on $K\# O_{q/r}$, as follows. Let $n=\lfloor p/q \rfloor$ and $r/q=p/q-n$, as above. We  view $K$ as having integral framing $n=\lfloor p/q\rfloor$ and the knot $O_{q/r}$ as having the Morse framing induced by the 0-framing in Figure~\ref{eq:Alexander-gradings}, if we ignore the $-q/r$-framed unknot. These two choices induce a Morse framing on $K\# O_{q/r}\subset -L(q,r)$, which corresponds to the rational number $p/q$, under the map in ~\eqref{eq:Morse-framing-rational}.

In particular, the intersection points have Alexander gradings
\begin{equation}
\frac{-q+1}{2q}, \frac{-q+3}{2q}, \dots, \frac{q-3}{2q},\frac{q-1}{2q}. \label{eq:Alexander-gradings}
\end{equation}

Now, if $x_i$ is the intersection point in $\Spin^c$ structure $\frs_i$, then $\cA_{(t,\frs_i)}(K\# O_{q/r})\simeq \cA_{t-A(x_i)}(K)$, by the connected sum formula. Note that the projection map from $\Spin^c(W'_{\lambda}(K))\to \Q$ is an injection (i.e. all $\Spin^c$ structures on $-L(q,r)$ have different values under $\hat{A}$), so we lose no information in writing just $\cA_{t}(K\# O_{q/r})$ instead of $\cA_{(t,\frs_i)}(K\# O_{q/r})$. By the K\"{u}nneth theorem, we have
\begin{equation}
\cA_{t}(K\# O_{q/r})\iso \cA_{[t]}(S^3,K), \label{eq:Kunneth-A-Lqr}
\end{equation}
where $[t]$ denotes the closest integer to $t$. The value of $[t]$ on half-integers is not important, since we never need to apply it there. Note also that the knot involution $\iota_{K\# O_{q/r}}$ coincides with $\iota_K$ under the isomorphism in~\eqref{eq:Kunneth-A-Lqr}, by the connected sum formula for involutive knot Floer homology \cite{ZemConnectedSums} as well as the fact that $O_{q/r}$ is a U-knot (so the maps $\Phi$ and $\Psi$ thereon vanish).

 We define
\begin{equation}
\begin{split}
\bA_{p/q}(K):=\frA(-L(q,r),K\# O_{q/r})&\iso \prod_{s\in \Z} \bm{\cA}_{[\frac{q+2s+1}{2q}]}(K)\quad \text{and} \\
\bB_{p/q}(K):=\frB(-L(q,r),K\# O_{q/r})&\iso \prod_{s\in \Z} \bm{\cB}_{[\frac{q+2s+1}{2q}]}(K).
\end{split}
\label{eq:rational-surgery-chain-complex}
\end{equation}
We think of each $\bm{\cA}_{[\frac{q+2s+1}{2q}]}(K)$ as actually being $\bm{\cA}_{[\frac{q+2s+1}{2q}]}(K)\times \{q+2s+1\}$ (so as to remember index $s$), but we suppress this from the notation.

Of course, there are canonical isomorphisms $\bm{\cA}_s(K)\iso \bm{A}_s(K)$ and $\bm{\cB}_s(K)\iso \bm{B}_s(K)$, so we may translate~\eqref{eq:rational-surgery-chain-complex} into the more familiar versions of complexes, as in \cite{OSKnots}. 

Applying Theorem~\ref{thm:rationally-null-mapping-cone}, we obtain the following:

\begin{thm}\label{thm:rational-surgery-mapping-cone} Suppose that $K$ is a knot in $S^3$, and $p$ and $q$ are relatively prime integers, neither of which is 0. Then $\bCFI^-(S^3_{p/q}(K))$ is homotopy equivalent over the ring $\bF\llsquare U \rrsquare[Q]/Q^2$ to a complex $\bF[Q]/(Q^2)\otimes_\bF (\bA_{p/q}\oplus \bB_{p/q})[-1]$, where $\bA_{p/q}$ and $\bB_{p/q}$ are as in~\eqref{eq:rational-surgery-chain-complex}, with differential specified by some $\bF\llsquare U \rrsquare$-linear maps as below:
\[
\bX \bI_{p/q}= 
 \begin{tikzcd}[column sep=2cm, row sep=2cm, labels=description]
\bA_{p/q}
	\ar[d, "Q\cdot (\id +\iota_{\bA})",swap]
	\ar[dr, dashed, "Q\cdot H"]
	\ar[r, "D_{p/q}"]
& \bB_{p/q}
	\ar[d, "Q\cdot (\id+\iota_{\bB})"]\\
Q\cdot \bA_{p/q}
	\ar[r, "D_{p/q}"]
& Q\cdot \bB_{p/q},
\end{tikzcd}.
\]
Equivalently, the  $\iota$-complex $(\bCF^-(S^3_{p/q}(K)), \iota)$ is $\iota$-homotopy equivalent to the complex $\bX_{p/q}:=\Cone\left(D_{p/q}\colon \bA_{p/q}\to \bB_{p/q}\right)$, with involution $\iota_{\bX}:=\iota_{\bA}+H+\iota_{\bB}$. Furthermore, the following hold:
\begin{enumerate}
\item \label{cone:1-pq} $D_{p/q}=v+h$, where $v$ sends $\bm{A}_{[t]}(K)$ to $\bm{B}_{[t]}(K)$, and $h$ sends $\bm{A}_{[t]}(K)$ to $\bm{B}_{[t+p/q]}$. 
\item\label{cone:2-pq} $\iota_{\bA}$ sends $\bm{A}_{[t]}$ to $\bm{A}_{[-t]}$. Furthermore, $\iota_{\bA}$ is the map induced by the knot involution under the isomorphism $\cA_n(K)\iso A_n(K)$.
\item \label{cone:3-pq} $\iota_{\bB}$ sends $\bm{B}_{[t]}$ to $\bm{B}_{[-t+p/q]}$, and coincides with the map induced by $\iota_{\bA}$, composed with $\frF_{z\to w}$.
\item\label{cone:4-pq} $H=k+j$, where $k$ maps $\bm{A}_{[t]}$ to $\bm{B}_{[-t+p/q]}$ and $j$ maps $\bm{A}_{[t]}$ to $\bm{B}_{[-t]}$.
\end{enumerate}
There is an algebraic $\iota$-complex $\XI_{p/q}^{\alg}(\CFK^\infty(K),\iota_K)$, characterized by the property that the maps $H$ and $h$ factor through $\vopp$, which is well-defined up to $\iota$-homotopy equivalence. Furthermore
\[
(\bCFI^-(S^3_{p/q}(K)),\iota)\simeq \XI^{\alg}_{p/q}(\CFK^\infty(K),\iota_K).
\]
\end{thm}

We make several basic comments about the rational surgeries mapping cone. Our rational surgeries involutive mapping cone is indexed slightly differently than Ozsv\'{a}th and Szab\'{o}'s mapping cone \cite{OSRationalSurgeries}. In their notation
\[
\bA_{p/q}:=\prod_{s\in \Z} A_{\lceil s/q \rceil}\times \{s\}\quad \text{and} \quad \bB_{p/q}:=\prod_{s\in \Z} B_{\lceil s/q\rceil}\times \{s\}.
\]
That is, $\bA_{p/q}$ and $\bB_{p/q}$ each consist of $q$ copies of each $A_s$ and $B_s$. In the differential $D_{p/q}:=v+h$, the map $v$ sends $A_{\lceil s/q\rceil}\times \{s\}$ to $B_{\lceil s/q\rceil} \times \{s\}$, and $h$ sends $A_{\lceil s/q\rceil}\times \{s\}$ to $B_{\lceil s+p/q\rceil}\times \{s+p\}$.

Next, we discuss the decomposition of the two cones over $\Spin^c$ structures. We have $\Spin^c(S^3_{p/q}(K))\iso \Z/p$. In Ozsv\'{a}th and Szab\'{o}'s indexing, the cone decomposes over $\Z/p$ congruence classes of the $s$ index. 

The conjugation action on $\Spin^c$ structures is a bit more natural with respect to the description in Theorem~\ref{thm:rational-surgery-mapping-cone}. Define
\[
\bH:=\left\{ \frac{q+2s+1}{2q}:s\in \Z\right\}.
\]
There is an isomorphism $\Spin^c(S^3_{p/q}(K))\iso \bH/(p/q)$, and analogous to~\eqref{eq:conjugation-H} the conjugation action sends $x\in \bH/(p/q)$ to $-x$. Each element $(q+2s+1)/2q\in \bH/(p/q)$ has a representative in the interval $[-p/2q,p/2q]$. The representative is unique, unless $p$ and $q$ differ in parity, in which case the $\Spin^c$ structure corresponding to $-p/2q$  has two representatives in the inteveral (the other being $p/2q$).  The self-conjugate $\Spin^c$ structures are as follows:
\begin{itemize}
\item If $p$ and $q$ are both odd, then there is one self-conjugate $\Spin^c$ structure, which corresponds to $0\in \bH$.
\item If $q$ is odd and $p$ is even, then there are two self-conjugate $\Spin^c$ structures, which are the ones which correspond to $0,p/2q\in \bH$ 
\item If $q$ is even and $p$ is odd, then there is one self-conjugate $\Spin^c$ structure, which corresponds to $p/2q$.
\end{itemize}

The techniques of Section~\ref{sec:local-class-surgeries} extend without complication to compute the local class of rational surgeries on knots, as follows:

\begin{prop}\label{prop:local-class-surgeries-rational} Suppose that $K$ is a knot in $S^3$, and $p,q$ are relatively prime, positive integers.
\begin{enumerate}
\item If $p$ and $q$ are both odd, then there is only one self-conjugate $\Spin^c$ structure on $S^3_{p/q}$, corresponding to $0\in \bH$. Moreover, $\CFI^-(S^3_{p/q}(K),[0])$ is locally equivalent to $(A_0(K),\iota_K)$.
\item If $q$ is odd and $p$ is even, then there are two self-conjugate $\Spin^c$ structures, corresponding to $0$ and $p/2q\in \bH$. Then, $\CFI^-(S^3_{p/q}(K),[0])$ is locally equivalent to $(A_0(K),\iota_K)$, while $\CFI^-(S^3_{p/q}(K),[p/2q])$ is locally equivalent to
\begin{equation}
\begin{tikzcd} [column sep={1cm,between origins},labels=description] 
A_{[p/2q]}
\ar[dr, "v"]
& & A_{[p/2q]} \ar[dl,"v"]\\
& B_{[p/2q]}
\end{tikzcd}
\label{eq:local-class-rational-even}
\end{equation}
equipped with the involution which interchanges the two copies of $A_{[p/2q]}$, and fixes $B_{[p/2q]}$.
\item If $q$ is even and $p$ is odd, then there is one self-conjugate $\Spin^c$ structure, which corresponds to $p/2q\in \bH$. Moreover, $\CFI^-(S^3_{p/q}(K),[p/2q])$ is locally equivalent to the complex shown in~\eqref{eq:local-class-rational-even}.
\end{enumerate}
\end{prop}

The proof of Proposition~\ref{prop:local-class-surgeries-rational} follows the same strategy as the proof of Proposition~\ref{prop:local-equivalence-class}. Namely, that proof implies that the local class of $\XI_{p/q}^{\alg}(\scC)$ coincides with the minimal truncation of $\XI_{p/q}^{\alg}(\scC)$, which is exactly the claim of Proposition~\ref{prop:local-class-surgeries-rational}. As a corollary of Proposition~\ref{prop:local-class-surgeries-rational}, together with the proof of Proposition~\ref{prop:correction-terms-main-computation}, we obtain the following (stated as Proposition~\ref{prop:correction-terms-intro} in the introduction):

\begin{prop}\label{prop:correction-terms-rational}
 Suppose that $p,q>0$ are relatively prime integers and $K$ is a knot in $S^3$.
\begin{enumerate}
\item If $p$ and $q$ are odd, then there is one self-conjugate $\Spin^c$ structure  $[0]$ on $S^3_{p/q}(K)$, and furthermore
\[
\dl(S^3_{p/q}(K),[0])=d(L(p,q),[0])-2 \Vl_0(K), \qquad \du(S^3_{p/q}(K),[0])=d(L(p,q),[0])-2\Vu_0(K).
\]
\item If $p$ is even and $q$ is odd, there are two self-conjugate $\Spin^c$ structures, which we denote $[0]$ and $[p/2q]$. The correction terms for $[0]$ are the same as the previous case, while
\[
\dl(S^3_{p/q}(K), [p/2q])=d(S_{p/q}^3(K), [p/2q]), \qquad \du(S_{p/q}^3(K), [p/2q])=d(L(p,q), [p/2q]).
\]
\item If $q$ is even and $p$ is odd, then there is just one self-conjugate $\Spin^c$ structure $[p/2q]$, and
\[
\dl(S^3_{p/q}(K), [p/2q])=d(S_{p/q}^3(K), [p/2q]), \qquad \du(S_{p/q}^3(K), [p/2q])=d(L(p,q), [p/2q]).
\]
\end{enumerate} 
\end{prop}

\subsection{0-surgeries}
\label{subsec:0-surgery}

As a final variation on the surgery formula, we discuss 0-surgeries. Suppose $K$ is a null-homologous knot in an integer homology 3-sphere $Y$. It follows from \cite{OSIntegerSurgeries}*{Section~4.8} that $\CF^-(Y_0(K),[0])\simeq \Cone(v_0+h_0\colon A_0\to B_0)$. There is also a simple analog of Theorem~\ref{thm:mapping-cone} which holds for 0-surgeries.  The proof of Theorem~\ref{thm:mapping-cone} carries over to show that $\CFI^-(Y_0(K), [0])$ is homotopy equivalent over the ring $\bF[U,Q]/Q^2$ to a chain complex of the form
\begin{equation}
\begin{tikzcd}[column sep=2cm, row sep=2cm, labels=description]
A_0
	\ar[d, "Q\cdot (\id +\iota_{\bA})",swap]
	\ar[dr, dashed, "Q\cdot H_0"]
	\ar[r, "D_0"]
& B_0
	\ar[d, "Q\cdot (\id+\iota_{\bB})"]\\
Q\cdot A_0
	\ar[r, "D_0"]
& Q\cdot B_0
\end{tikzcd}
\label{eq:diagram-for-0-surgery}
\end{equation}
where
\[
\begin{split}
D_0&=v_0+h_0,\\
\iota_{\bA}&=\iota_K,\\
 \iota_{\bB}&=\frF_0\iota_K.
\end{split}
\]

Mirroring this, if $\scC$ is an $\iota_K$-complex of L-space type, we can construct an algebraic version of the involutive mapping cone complex $\XI_{0}^{\alg}(\scC)$, with underlying chain complex as in~\eqref{eq:diagram-for-0-surgery}, which is characterized by the additional requirement that $h_0=\frF_0 \vopp_0$ and $H_0=H\vopp$, for an $\bF[U]$-equivariant map $H\colon \Bopp_0\to B_0$. Theorem~\ref{thm:strong-cone-formula-proven} extends into this setting to show that if $K$ is a knot in $S^3$, then
 \[
 \CFI^-(S^3_0(K),[0])\simeq \XI^{\alg}_0(\CFK^\infty(K),\iota_K).
\]

\bibliographystyle{custom}
\def\MR#1{}
\bibliography{biblio}

\end{document}